\pdfoutput=1
\documentclass[11pt,reqno]{amsart}

\numberwithin{equation}{section}

\usepackage[tt=false]{libertine}

\usepackage{amssymb}
\usepackage[varbb]{newpxmath}
\usepackage{mathrsfs}

\usepackage{upgreek}

\let\savedbigtimes\bigtimes
\let\bigtimes\relax
\usepackage{mathabx} 
\let\bigtimes\savedbigtimes

\usepackage[margin=1in]{geometry}
\usepackage{enumerate}

\usepackage[usenames,dvipsnames,table]{xcolor}
\usepackage[colorlinks=true,
    linkcolor=teal!60!black,
    citecolor=PineGreen,
    urlcolor=RedViolet]{hyperref}
\hypersetup{bookmarksopen=true}

\usepackage{tikz}
\usetikzlibrary{positioning, arrows.meta, calc}

\usepackage[yyyymmdd,hhmmss]{datetime}

\usepackage{mathtools}

\usepackage{xparse}

\newtheorem{thm}{Theorem}[section] 
\newtheorem{lem}[thm]{Lemma}
\newtheorem{ppn}[thm]{Proposition}
\newtheorem{cor}[thm]{Corollary}

\theoremstyle{definition}
\newtheorem{dfn}[thm]{Definition}
\newtheorem{rmk}[thm]{Remark}
\newtheorem*{rmk*}{Remark}
\newtheorem{ass}[thm]{Assumption}
\newtheorem{xmp}[thm]{Example}

\def\beq#1\eeq{%
    \begin{equation}%
    #1%
    \end{equation}%
}

\newcommand{\R}{\mathbb{R}}
\newcommand{\bbR}{\mathbb{R}}
\newcommand{\E}{\mathbb{E}}
\newcommand{\bbE}{\mathbb{E}}

\renewcommand{\P}{\mathbb{P}}
\newcommand{\bbP}{\mathbb{P}}
\newcommand{\Q}{\mathbb{Q}}

\newcommand{\bbW}{\mathbb{W}}
\newcommand{\bbZ}{\mathbb{Z}}

\newcommand{\cL}{\mathcal{L}}
\newcommand{\cA}{\mathcal{A}}

\newcommand{\cP}{\mathcal{P}}
\newcommand{\cF}{\mathcal{F}}
\newcommand{\cH}{\mathcal{H}}
\newcommand{\cM}{\mathcal{M}}
\newcommand{\ocM}{\bar{\mathcal{M}}}

\newcommand{\cG}{\mathscr{G}}
\newcommand{\cFt}{{\mathscr{F}}}

\newcommand{\be}{\boldsymbol{e}}

\def\bsig{{\boldsymbol \sigma}}

\def\bb{{\boldsymbol{b}}}
\def\bg{{\boldsymbol{g}}}
\def\bm{{\boldsymbol{m}}}
\def\bs{{\boldsymbol{s}}}
\def\br{{\boldsymbol{r}}}
\def\bx{{\boldsymbol{x}}}
\def\bu{{\boldsymbol{u}}}
\def\bv{{\boldsymbol{v}}}
\def\bw{{\boldsymbol{w}}}
\def\bmeta{{\boldsymbol{\eta}}}
\def\bXi{{\boldsymbol \Xi}}
\def\bTe{{\boldsymbol \Theta}}
\def\by{{\boldsymbol y}}
\def\bz{{\boldsymbol z}}

\def\bA{{\boldsymbol A}}
\def\bG{{\boldsymbol G}}
\def\bJ{{\boldsymbol J}}

\newcommand{\tbg}{{\tilde \bg}}

\newcommand{\vl}{\underline{\smash{\ell}}}
\newcommand{\vX}{\underline{\smash{X}}}
\newcommand{\vY}{\underline{\smash{Y}}}
\newcommand{\tvY}{\underline{\smash{\tilde{Y}}}}
\newcommand{\hvX}{\underline{\smash{\iX}}}
\newcommand{\hvY}{\underline{\smash{\iY}}}
\newcommand{\hvZ}{\underline{\smash{\iZ}}}
\newcommand{\vZ}{\underline{\smash{Z}}}
\newcommand{\vD}{\underline{\smash{D}}}
\newcommand{\vQ}{\underline{\smash{Q}}}
\newcommand{\vtQ}{\underline{\smash{\tilde{Q}}}}
\newcommand{\barvD}{\underline{\smash{\bar{D}}}}

\newcommand{\ubl}{\ell_{\max}}

\newcommand{\st}{{\mathsf{t}}} % matrix transpose
\newcommand{\Cost}{\textup{\textsf{cost}}}
\def\budget{{\mathsf{budget}}}

\def\incr{{\mathsf{incr}}}
\def\sign{{\mathsf{sign}}}
\def\sym{{\mathsf{sy}}}
\def\sph{{\mathsf{sph}}}

\newcommand{\bP}{\textup{\textbf{P}}}
\newcommand{\bE}{\textup{\textbf{E}}}

\newcommand{\bX}{\textup{\textbf{X}}}

\newcommand{\bZ}{\textup{\textbf{Z}}}
\newcommand{\dmax}{d_{\max}}
\def\bzero{\textup{\textbf{0}}}
\def\ind{\textup{\textbf{1}}}

\DeclareMathOperator{\dist}{dist}
\DeclareMathOperator{\plim}{p-lim}

\def\ons{\textup{\textbf{ons}}}

\def\IAMP{\mathsf{IAMP}}
\def\BOGP{{\mathsf{BOGP}}}

\newcommand{\RomI}{\textup{I}}
\newcommand{\RomII}{\textup{II}}
\newcommand{\RomIII}{\textup{III}}
\newcommand{\Ising}{\textup{Is}}
\newcommand{\Isplus}{{\frozen,\Ising}}

\newcommand{\op}{\textup{op}}

\def\oJ{\bar{J}}
\def\oL{\bar{L}}

\newcommand{\plus}{{\textup{+}}}

\newcommand{\BUCKETS}{\mathscr{B}}
\newcommand{\trK}{{\bar{K}}}
\newcommand{\trstar}{{\textup{tr}*}}
\newcommand{\trunc}{\textup{tr}}

\newcommand{\frozen}{\plus}

\newcommand{\frstar}{{\frozen*}}

\newcommand{\cFtplus}{\cFt^\frozen}

\newcommand{\sep}{\tilde{\varepsilon}}
\newcommand{\sphi}{{\varphi_\textup{sp}}} 
\newcommand{\sphicirc}{\varphi_\circ} 
\newcommand{\xxx}{\underline{\smash{\mathfrak{x}}}}
\newcommand{\sB}{{\tilde{B}}}
\newcommand{\sV}{{\tilde{V}}}
\newcommand{\sU}{{\tilde{U}}}
\newcommand{\sF}{{\tilde{F}}}
\newcommand{\sbbb}{{\tilde{\bbb}}}
\newcommand{\svvv}{{\tilde{\vvv}}}
\newcommand{\suuu}{{\tilde{\uuu}}}
\newcommand{\swww}{{\tilde{\www}}}
\newcommand{\srrr}{\tilde{\rrr}}
\newcommand{\sPsi}{{\tilde{\Psi}}}
\newcommand{\sE}{{\tilde{\bE}}}
\newcommand{\sGam}{{\tilde{\Gamma}}}

\newcommand{\tep}{\bar{\epsilon}}
\newcommand{\tphi}{{\varphi_\textup{tm}}} 
\newcommand{\sss}{{\mathfrak{s}}}

\newcommand{\tbbb}{{\bar{\bbb}}}
\newcommand{\tvvv}{{\bar{\vvv}}}
\newcommand{\tuuu}{{\bar{\uuu}}}
\newcommand{\twww}{{\bar{\www}}}
\newcommand{\trrr}{\bar{\rrr}}
\newcommand{\tE}{{\bar{\bE}}}
\newcommand{\tGam}{{\bar{\Gamma}}}
\newcommand{\tB}{{\bar{B}}}
\newcommand{\tV}{{\bar{V}}}
\newcommand{\tU}{{\bar{U}}}
\newcommand{\tF}{{\bar{F}}}
\newcommand{\tSigma}{\bar{\Sigma}}

\newcommand{\Tbbb}{{\breve{\bbb}}}
\newcommand{\Tuuu}{{\breve{\uuu}}}
\newcommand{\Tvvv}{{\breve{\vvv}}}
\newcommand{\TGam}{{\breve{\Gamma}}}

\newcommand{\brep}{\breve{\epsilon}}

\title[Algorithmic threshold for high-dimensional projection pursuit I]{Algorithmic threshold for high-dimensional projection pursuit I: general theory}
\date{\today}

\author[B.~Huang, M.~Sellke, and N.~Sun]{Brice Huang$^*$ \and Mark Sellke$^\bullet$ \and Nike Sun$^\circ$}

\thanks{$^*$Department of Statistics, Stanford University. \newline\indent $^\bullet$Department of Statistics, Harvard University. \newline\indent 
$^\circ$Department of Mathematics, Massachusetts Institute of Technology.}

\renewcommand{\emptyset}{\varnothing}

\newcommand{\bbb}{\mathfrak{b}}
\newcommand{\vvv}{\mathfrak{v}}
\newcommand{\uuu}{\mathfrak{u}}
\newcommand{\www}{\mathfrak{w}}
\newcommand{\rrr}{\mathfrak{r}}

\newcommand{\iB}{\hat{B}}
\newcommand{\iF}{\hat{F}}
\newcommand{\iV}{\hat{V}}
\newcommand{\iU}{\hat{U}}
\newcommand{\iw}{\hat{\www}}
\newcommand{\ir}{\hat{\rrr}}
\newcommand{\iX}{\hat{X}}
\newcommand{\iY}{\hat{Y}}
\newcommand{\iZ}{\hat{Z}}

\newcommand{\MAX}{\mathfrak{M}}
\newcommand{\SMAX}{\mathfrak{S}}

\newcommand{\TMAX}{\mathfrak{T}}

\def\lt{\left}
\def\rt{\right}
\def\de{{\rm d}}
\def\cN{{\mathcal N}}
\newcommand{\Law}{\mathscr{L}}

\newcommand{\vx}{\underline{\smash{x}}}
\newcommand{\vy}{\underline{\smash{y}}}
\newcommand{\vI}{\underline{\smash{1}}}

\renewcommand{\vec}{\underline}

\newcommand{\EmpDist}{\textup{\textsf{Law}}}
\newcommand{\ProgMsrbl}{\textup{\textsf{Prog}}}

\DeclareMathOperator{\Tr}{tr}

\def\eps{{\varepsilon}}

\newcommand{\ab}{\acute{b}}
\newcommand{\av}{\acute{v}}
\newcommand{\ar}{\acute{r}}
\newcommand{\aGam}{\acute{\Gamma}}
\newcommand{\au}{\acute{u}}
\newcommand{\aw}{\acute{w}}
\newcommand{\asig}{\acute{\sigma}}

\def\Adm{{\mathsf{Adm}}}
\def\unif{{\mathsf{unif}}}
\def\aux{{\mathsf{aux}}}
\newcommand{\cE}{{\mathcal{E}}}
\def\sP{{\mathscr{P}}}
\def\sS{{\mathscr{S}}}
\def\sW{{\mathscr{W}}}
\def\ALG{{\mathsf{ALG}}}

\def\ioeps{{\iota}}
\def\iotasol{{\iota_{\mathsf{sol}}}}
\def\iotamsr{{\iota_{\mathsf{msr}}}}
\def\iotaval{{\iota_{\mathsf{val}}}}

\newif\iffull\fulltrue

\newcommand{\concave}{\textit{c}}
\newcommand{\WERR}{\textup{\textsf{ERR}}}

\newcommand{\coarse}{\textup{c}}
\newcommand{\eaux}{\epsilon_\textup{aux}}

\begin{document}

\begin{abstract}
We study a null model of high-dimensional projection pursuit: we are given $M$ points sampled i.i.d. from a standard gaussian in $N$ dimensions, where $M,N\to\infty$ with $M/N\to\alpha\in(0,\infty)$. Our goal is to characterize the possible empirical distributions of these points' projections along a data-dependent direction $\bx$, which ranges over either the sphere $S_N=\sqrt{N}\mathbb{S}^{N-1}$ or cube $\Sigma_N=\{-1,+1\}^N$. We consider this problem in an algorithmic setting, where $\bx$ must be the output of an algorithm with dimension-free Lipschitz dependence on the input; this class of algorithms includes general gradient-based methods such as Langevin dynamics and approximate message passing (AMP). Our main result exactly characterizes the set of empirical distributions attainable by this class in terms of a one-dimensional stochastic control problem. As a consequence of our main result, we obtain exact algorithmic thresholds for optimizing the Hamiltonian of a spherical or Ising perceptron model with general bounded continuous activation. For the spherical problem, independent work of Montanari and Zhou \cite{montanari2024exceptional} characterized the empirical distributions attainable by a related two-stage AMP algorithm, also in terms of stochastic control.

Our proof of hardness builds on the branching overlap gap property introduced in earlier work by the first two authors. Our main innovation is to develop stochastic control theory within the branching OGP framework, significantly expanding the settings in which it locates an exact algorithmic threshold. Notably, our methods apply even though the non-algorithmic problem of characterizing all feasible projections remains a major outstanding challenge. For the matching algorithmic result, we construct a new incremental AMP algorithm that acts on a Brownian-bridge revelation of the gaussian disorder and simulates the same family of controlled SDEs.

All ideas in this paper are human-generated, and all the writing was done by the human authors. AI was used in the writing of this paper only for light proofreading and copy-editing.
\end{abstract}

\maketitle

\vspace{-0.3cm}
\setcounter{tocdepth}{1}
\tableofcontents

\vfill

\pagebreak\section{Introduction}\label{s:intro}

\def\fq{{\mathfrak{q}}}

\iffull
% !TEX root = main.tex

In this paper we study \textbf{projection pursuit}, a model of exploratory data analysis which can be described informally as follows.
If we imagine that we have $M$ data points in $N$ dimensions, then projection pursuit refers to the task of finding a low-dimensional projection of this data that looks ``interesting'' --- in particular, non-gaussian. The formal notation for the problem is introduced in \eqref{eq:proj-pursuit-def} below.
The idea of projection pursuit was put forth by Kruskal \cite{kruskal1969toward,kruskal1972linear}, and we review some of the surrounding literature in \S\ref{subsec:literature} below.

It was shown in \cite{diaconis1984asymptotics} that under fairly mild conditions, if $M/N\to\infty$ (the \textbf{low-dimensional regime}), then \emph{all} $k$-dimensional projections ($k\le N$ any constant) are gaussian, in the sense that the empirical measure of the projected data approximates a gaussian measure in $k$ dimensions
 (see also \cite{friedman1974projection,friedman1987exploratory}). In the \textbf{high-dimensional regime} where $M,N\to\infty$ with $M\asymp N$, this is no longer the case, and the set of possible projections is much less well understood. Note that a natural null model is that the data consists of $M$ points sampled independently from the standard gaussian distribution in $N$ dimensions. For this gaussian null model in the high-dimensional regime, we consider the set of achievable  \textbf{one-dimensional} projections. This question has been studied by prior works (\cite{bickel2018projection,montanari2022overparametrized}), which showed that non-gaussian projections do exist, but that all projections must be bounded in $2$-Wasserstein distance to a standard gaussian distribution.

Beyond the questions of which projections are possible, projection pursuit also poses computational challenges. In particular, given a target distribution $\mu$, the problem of finding a projection close to $\mu$ cannot be formulated via a convex objective. Hence, there may be measures $\mu$ such that there typically \emph{exists} a direction along which the empirical distribution of projected data is very close to $\mu$, but \emph{finding} any such direction is computationally hard. \textbf{In this paper we study the question of which one-dimensional projections are algorithmically feasible. Our main result Theorem~\ref{thm:main} (see also the informal version Theorem~\ref{t:main.INFORMAL}) gives a precise characterization of the set of measures achievable by Lipschitz algorithms.}

We now introduce some notation to formalize our discussion. As noted above, we assume the data are points $\bg^1,\ldots,\bg^M$ sampled independently from the standard gaussian distribution in $N$ dimensions. We can collect this data in the $M\times N$ matrix $\bG$, where each row corresponds to one data point. Let $S_N$ denote the sphere of radius $N^{1/2}$ inside $\R^N$. For any $\bx\in S_N$, define the empirical measure
\begin{equation}
    \label{eq:proj-pursuit-def}
    \mu_{\bG}(\bx)
    \equiv \EmpDist\bigg(\frac{\bG\bx}{N^{1/2}}\bigg)
    \equiv
    \frac{1}{M}\sum_{a=1}^M \delta\bigg\{
    \frac{(\bg^a,\bx)}{N^{1/2}}
    \bigg\}\,,
\end{equation}
where $\delta\{c\}$ denotes a Dirac mass at point $c$. Thus $\mu_{\bG}(\bx)$ is the empirical distribution of the data $\bG$ projected in the $\bx$ direction. The set of all possible empirical measures is given by
	\beq
\label{eq:proj-pursuit-set}
\cM_{\bG}=\Big\{\mu_{\bG}(\bx):\bx\in S_N\Big\}.
\eeq
In this paper, we study the question of
which measures $\mu$ can be approximated using directions $\bx$ which can be found \emph{by a large class of efficient algorithms}, given the data $\bG$ as input. We also study the case where $\bx$ is restricted to $\Sigma_N\equiv\{-1,+1\}^N$.

Let us also mention a closely related problem, the \textbf{perceptron model}. With $\bg^a$ as before, this can be defined in terms of the random Hamiltonian
\beq
    \label{e:hamiltonian}
    H_N(\bx) \equiv
    \sum_{a=1}^M
    \phi\bigg(
        \frac{(\bg^a,\bx)}{N^{1/2}}
    \bigg)
    \equiv M \phi(\mu_{\bG}(x))
    \,,
\eeq
where $\bx\in\R^N$, the activation function $\phi\in C_b(\bbR;\bbR)$ is fixed (independently of $N$),
and $\phi(\mu)$ denotes the integral of $\phi$ against $\mu$. One may consider the problem of (efficiently) optimizing $H_N(\bx)$
 over either $\bx\in S_N$ (the \textbf{spherical perceptron} problem), or $\bx\in\Sigma_N$ (the \textbf{Ising perceptron} problem). Recalling \eqref{eq:proj-pursuit-set}, it is immediate that the optimal value of the spherical perceptron Hamiltonian is given by
\[
\frac1{M} \sup \Big\{
H_N(\bx) : \bx\in S_N \Big\}
=
\sup \Big\{ 
\phi (\mu)
:\mu\in\cM_{\bG}\Big\}\,,
\]
and similarly in the Ising case (replacing $S_N$ with $\Sigma_N$). Thus, from our perspective, (algorithmic) projection pursuit contains the (algorithmic) perceptron problem as a special case, so we focus our attention on the former. The remainder of this section is organized as follows: 
\begin{itemize}
\item In \S\ref{ss:intro.results.informal} we introduce our main results, deferring some more technical definitions to later in the section. Our main results for the algorithmic projection pursuit problem are stated informally in Theorem~\ref{t:symm.INFORMAL} for the symmetrized setting,
and in Theorem~\ref{t:main.INFORMAL} for the general setting.

\item In \S\ref{subsec:literature} we review surrounding literature. A large amount of work in statistical physics, probability, computer science, and statistics has studied the above models and closely related variants.

\item In \S\ref{subsec:ideas} we mention some of the key ideas of the proof. 

\item In \S\ref{ss:intro.lip.algs} we present the formal definitions that were omitted from \S\ref{ss:intro.results.informal}, and give the formal statements of our main results: Theorem~\ref{thm:symmetric} for the symmetrized setting (formalizing Theorem~\ref{t:symm.INFORMAL}), and Theorem~\ref{thm:main} for the general setting (formalizing Theorem~\ref{t:main.INFORMAL}).

\item In \S\ref{ss:intro.organization} we describe how the proof is organized across the remainder of this paper.
\end{itemize}

\subsection{Main results}
\label{ss:intro.results.informal}

We now introduce our main results, deferring some of the more technical definitions to \S\ref{ss:intro.lip.algs} below. As mentioned above, our goal is to understand the power of efficient algorithms for projection pursuit and the perceptron problem. We view algorithms as maps
\beq\label{eq:algorithms-as-maps}
  \begin{aligned}  \cA :
     \R^{M\times N} \times \R^{N_{\aux}}
    &\to \bbR^N\,, \\
    (\bG, \bg^{\aux}) 
    &\mapsto \bx\,.
    \end{aligned}\eeq
In the above, $\bG$ is the gaussian disorder matrix appearing in the definition of the problem \eqref{eq:proj-pursuit-def}. Meanwhile, $\bg^{\aux}$ is an auxiliary input, given by a standard gaussian random vector in $N_{\aux}$ dimensions that is independent of $\bG$.

Our main results concern the class of  algorithms with \textbf{dimension-free Lipschitz} dependence on the input $(\bG,\bg^\aux)$; the precise specification is given by Definition~\ref{d:Lip} below. This class includes many of the gradient-based algorithms which are widely used on high-dimensional optimization and search problems (see e.g. \cite[Section 8]{HuangSellke2021}). In our setting, the goal of the algorithm will be to output $\bx$ such that $\mu_\bG(\bx)$ approximates a given target measure $\mu$. We will precisely characterize the set of measures $\mu$ which can be approximated by using Lipschitz algorithms for projection pursuit.

We will use a stochastic control problem to describe the set of achievable measures. We first describe this in the simpler setting of the \textbf{symmetrized} projection pursuit measures. For any measure $\mu$ on the real line, define the symmetrization
\beq
    \label{e:symmetrize.measure}
    \sym(\mu) = \cL(|x| : x \sim \mu)\,.
\eeq
Then, rather than the measures 
$\mu_{\bG}(\bx)$ from \eqref{eq:proj-pursuit-def}, suppose we are interested in the symmetrizations
	\beq
    \label{eq:symmetrized-empirical-measure-intro}
    \mu_{\bG,\sym}(\bx)
    \equiv
    \sym(\mu_{\bG}(\bx))
    \equiv\frac{1}{M}\sum_{a=1}^M 
    \delta\bigg\{
    \frac{|(\bg^a,\bx)|}{N^{1/2}}
    \bigg\}\,.
	\eeq
The \textbf{symmetrized projection pursuit} problem asks, for typical $\bG$, to understand the set of achievable measures $\mu_{\bG,\sym}(\bx)$.
This is closely related to the perceptron optimization problem \eqref{e:hamiltonian} with \emph{even} activation functions $\phi$, and to the symmetric perceptron constraint satisfaction problems studied in a line of recent work (reviewed in \S\ref{subsec:literature}). As mentioned above, we consider both the \textbf{spherical} setting where $\bx$ goes over all of $S_N \equiv N^{1/2}\mathbb{S}^{N-1}$, and the \textbf{Ising} setting where $\bx$ is restricted to $\Sigma_N \equiv \{-1,+1\}^N$. (See \S\ref{ss:intro.lip.algs} where we introduce relaxed versions of $S_N,\Sigma_N$ that are appropriate for the analysis of Lipschitz algorithms.)

\begin{dfn}[stochastic control and achievable measures for symmetrized projection pursuit]
\label{d:achievable-msrs-sym}
Let $B, W$ be one-dimensional standard Brownian motions, and let $U,U' \sim \unif([0,1])$ (with $B,W,U,U'$ mutually independent).
Define the filtrations 
\[\begin{aligned}
    \cF_X(t) &= \sigma(U, (B(s) : 0\le s\le t))\,, \\
    \cF_{\Ising}(t) &= \sigma(U', (W(s) : 0\le s\le t))\,.
\end{aligned}\]
Let $\ProgMsrbl(\cF_X)$ denote the space of stochastic processes that are progressively measurable with respect to $\cF_X$, and similarly $\ProgMsrbl(\cF_{\Ising})$ for $\cF_{\Ising}$. For $\sigma\in\ProgMsrbl(\cF_X)$ and $w\in\ProgMsrbl(\cF_{\Ising})$, let  
	\begin{align}
	\label{eq:main-process-intro-p=1}
	X(t)
	&\equiv\int_0^s \sigma_t\,dB(s)\,,\\
	X^\Ising(t)
	&\equiv\int_0^t w_s\,dW(s)
	\label{eq:X-ising-intro}
	\end{align}
Let $\cM^{\Ising,\sym}(\alpha)$ be the set of measures $\mu$ which can be expressed as $\mu=\Law(|X(1)|)$
(i.e., $\mu$ is the law of the real-valued random variable $|X(1)|$)
 for $X$ defined by \eqref{eq:main-process-intro-p=1}, and for controls $(\sigma,w)$ satisfying the following conditions:
\begin{itemize}
\item the nonnegativity condition
$\sigma_t\ge0$ for all $t\in[0,1]$,
\item the time parametrization condition
$\E[(w_t)^2]=1$ for all $t\in[0,1]$,
\item the \textbf{budget constraint}
	\beq\label{e:intro.sym.budget}
	\E\Big[(\sigma_t-1)^2\Big]
	\le \frac{(\E w_t)^2}{\alpha}
	\eeq
for all $t\in[0,1]$, and 
\item the \textbf{endpoint condition}
$|X^{\Ising}(1)|=1$.
\end{itemize}
 Let $\cM^{\sph,\sym}(\alpha)$ be the set of measures $\mu$ which can be expressed as $\mu=\Law(|X(1)|)$ for $X$ defined by \eqref{eq:main-process-intro-p=1}, and for controls $(\sigma,w\equiv1)$ satisfying the same conditions as above, \textbf{except for the endpoint constraint}. We use $\bar{\cM}^{\Ising,\sym}(\alpha)$,
$\bar{\cM}^{\sph,\sym}(\alpha)$ to denote the closures with respect to the $2$-Wasserstein metric $\bbW_2$ (defined by \eqref{def:W2} below).
\end{dfn}

\begin{thm}[\textbf{informal statement of main result for symmetrized projection pursuit}]\label{t:symm.INFORMAL}
 A measure $\mu$ on the real line can be approximated by a symmetrized projection pursuit measure $\mu_{\bG,\sym}(\bx)$ (as defined by \eqref{eq:symmetrized-empirical-measure-intro}), with $\bx$ the output of a Lipschitz algorithm \eqref{eq:algorithms-as-maps},
if and only if:
\begin{enumerate}[(a)]
\item (for the spherical setting, where $\bx$ must lie near $S_N$) 
$\mu$ lies in 
$\bar{\cM}^{\sph,\sym}(\alpha)$;
\item (for the Ising setting, where $\bx$ must lie near $\Sigma_N$)  $\mu$ lies in 
$\bar{\cM}^{\Ising,\sym}(\alpha)$.
\end{enumerate}
The sets
$\bar{\cM}^{\sph,\sym}(\alpha)$ and  $\bar{\cM}^{\Ising,\sym}(\alpha)$
are specified by Definition~\ref{d:achievable-msrs-sym}.
\end{thm}

In \S\ref{ss:intro.lip.algs} we present precise definitions for the various approximations glossed over in the statement above. \textbf{Theorem~\ref{t:symm.INFORMAL} is then formalized by Theorem~\ref{thm:symmetric} below.} We next present our main result in the general (non-symmetrized) setting, where the stochastic control problem becomes more complicated. To this end let $\sP$ denote the set of increasing, absolutely continuous functions $p:[0,1] \to [0,1]$ with $p(1)=1$, and let $\sP_c \subseteq \sP$ denote the set of such $\sP$ that are furthermore concave with $p(0) = 0$.
(See Remark~\ref{r:intro.interpretation.p} below for an initial comment on the interpretation of the function $p$.) We also introduce the abbreviation 
    \[(tp)'(t)\equiv 
    \frac{d}{dt} [t p(t)]
    \equiv p(t) + t p'(t)\,,
    \] 
which is defined almost-everywhere since the function $tp(t)$ is absolutely continuous.

\begin{dfn}[admissible controls]
\label{d:admissible.controls}
Let $B,W,U,U',\cF_X,\cF_{\Ising}$ be as in
Definition~\ref{d:achievable-msrs-sym}. Let
\begin{align*}
    \sS &= \lt\{
        \begin{array}{ll}
        (b,\sigma): & \textup{$b,\sigma\in\ProgMsrbl(\cF_X)$, with
        $\sigma_t \ge 0$ almost surely} \\
        &\textup{and
        $\max\{\bbE[(b_t)^2], \bbE[(\sigma_t)^2]\} < \infty$ for all $t\in[0,1]$}
        \end{array}
    \rt\}\,, \\
    \sW &= \lt\{
        \begin{array}{ll}
        w : & \textup{$w\in\ProgMsrbl(\cF_{\Ising})$, with 
        $\bbE [(w_t)^2] = 1$ and} \\
        &\textup{$X^{\Ising}(1) \in \{\pm 1\}$ almost surely} 
        \end{array}
    \rt\}\,,
\end{align*}
where $X^{\Ising}$ is as in \eqref{eq:X-ising-intro}.
    For $p\in \sP$, let $\Adm^{\Ising}(\alpha,p)$
 denote the set of controls $(b,\sigma,w)\in \sS \times \sW$ satisfying, for all $0\le t\le 1$, the \textbf{Ising budget constraint}
    \beq
        \label{eq:is-budget-constraint}
        \bbE \bigg[(b_t)^2 + \frac{(tp)'(t)}{p(t)}(\sigma_t-1)^2\bigg]
        \le \frac{(\bbE w_t)^2 }{ \alpha}\,.
    \eeq
Similarly, let $\Adm^{\sph}(\alpha,p)$ denote the set of controls $(b,\sigma)\in\sS$ satisfying, for all $0\le t\le 1$, the \textbf{spherical budget constraint}
    \beq
        \label{eq:sp-budget-constraint}
        \bbE \bigg[(b_t)^2 + \frac{(tp)'(t)}{p(t)}(\sigma_t-1)^2\bigg]
        \le \frac1\alpha\,.
    \eeq
Then $\Adm^{\Ising}(\alpha,p)$ and $\Adm^{\sph}(\alpha,p)$ are the \textbf{admissible controls} for the stochastic control problems corresponding to the Ising and spherical problems.
\end{dfn}

Note that by taking $p\equiv1$ and $b\equiv0$ in Definition~\ref{d:admissible.controls}, we can retrieve the control problems given by Definition~\ref{d:achievable-msrs-sym} for the symmetrized setting. This is a substantial simplification and will be exploited to precisely analyze the symmetric perceptron in our companion work \cite{bogpinprogress}, where we compute the asymptotic behavior of associated satisfiability thresholds. In the spherical symmetric case, this subclass of algorithms was also considered by \cite{montanari2024exceptional}: their treatment of the asymmetric case is \textit{a~priori} different, though it might be possible to connect with ours.

\begin{dfn}[achievable measures]
    \label{d:achievable-msrs}
    For $(b,\sigma,p) \in \sS \times \sP$, let $\mu(b,\sigma,p)$ denote the law of $X(1)$ given by the stochastic integral
    \beq
        \label{eq:main-process-intro-general}
        X(t)
        = \int_0^t p'(s)^{1/2}
        b_s\,ds
        +
        \int_0^t \Big[ (sp)'(s) \Big]^{1/2} \sigma_s\,dB(s)\,.
    \eeq
The \textbf{sets of achievable measures} 
for the Ising and spherical problems are defined to be 
    \begin{align*}
        \cM^{\Ising}(\alpha) 
        &= \Big\{
            \mu(b,\sigma,p) : p \in \sP, (b,\sigma,w) \in \Adm^{\Ising}(\alpha,p)
        \Big\}, \\
        \cM^{\sph}(\alpha) &= \Big\{
            \mu(b,\sigma,p) : p \in \sP, (b,\sigma) \in \Adm^{\sph}(\alpha,p)
        \Big\}.
    \end{align*}
Define $\cM^{\Ising, \concave}(\alpha)$ and $\cM^{\sph, \concave}(\alpha)$ analogously, but with the additional restriction that $p\in\sP_c$. Let $\ocM^{\Ising}(\alpha)$,
$\ocM^{\sph}(\alpha)$,
$\ocM^{\Ising, \concave}(\alpha)$, $\ocM^{\sph, \concave}(\alpha)$
 be the closures of these sets in the 
    $\bbW_2$ metric.
\end{dfn}

\begin{thm}[\textbf{informal statement of main result}]
\label{t:main.INFORMAL} A measure $\mu$ can be approximated by $\mu_{\bG}(\bx)$, with $\bx$ the output of a Lipschitz algorithm \eqref{eq:algorithms-as-maps},
if and only if:
\begin{enumerate}[(a)]
\item \label{t:main.INFORMAL.sph} (for the spherical setting, where $\bx$ must lie near $S_N$)  $\mu$ lies in 
$\bar{\cM}^{\sph}(\alpha)$;
\item \label{t:main.INFORMAL.Ising} (for the Ising setting, where $\bx$ must lie near $\Sigma_N$) $\mu$ lies in 
$\bar{\cM}^{\Ising}(\alpha)$.
\end{enumerate}
The sets
$\bar{\cM}^{\sph}(\alpha)$ and 
$\bar{\cM}^{\Ising}(\alpha)$
are specified by Definition~\ref{d:achievable-msrs}. Moreover they coincide with the sets
$\ocM^{\sph, \concave}(\alpha)$ and
$\ocM^{\Ising, \concave}(\alpha)$ respectively.
\end{thm}

Again, precise definitions for the approximations mentioned in the above statement are presented in \S\ref{ss:intro.lip.algs}. \textbf{Theorem~\ref{t:main.INFORMAL} is then formalized by Theorem~\ref{thm:main} below.} From Theorem~\ref{t:main.INFORMAL} it is also possible to read off the algorithmic threshold for the perceptron optimization problem \eqref{e:hamiltonian}, assuming $\phi\in C_b(\bbR;\bbR)$. Informally speaking, for Lipschitz algorithms, the achievable threshold in the Ising perceptron (where the algorithm must output $\bx$ near $\Sigma_N$) is the value
\beq
    \label{eq:def-ALG-Ising}
    \ALG^{\Ising}(\phi,\alpha) 
    \equiv \sup\bigg\{
    \phi(\mu) 
    : \mu \in \ocM^{\Ising}(\alpha)
    \bigg\}\,.
\eeq
Likewise, the achievable threshold in the spherical perceptron (where the algorithm must output $\bx$ near $S_N$) is the value $\ALG^{\sph}(\phi,\alpha)$, defined similarly as \eqref{eq:def-ALG-Ising} but with $\ocM^{\sph}(\alpha)$ in place of
$\ocM^{\Ising}(\alpha)$. \textbf{This assertion is formalized by Corollary~\ref{cor:alg-for-optimization} below.}

\begin{rmk}[achievable measures for spherical versus Ising problems]\label{r:ising.subset.sphere}
Recalling Definition~\ref{d:achievable-msrs}, the set $\ocM^{\sph}(\alpha)$ can be equivalently defined by dropping the endpoint condition $X^{\Ising}(1) \in \{\pm 1\}$ from the definition of $\ocM^{\Ising}(\alpha)$.
    Indeed, by the Cauchy--Schwarz inequality, 
    \beq
        \label{eq:cauchy-schwarz}
        (\bbE w_t)^2 \le \bbE [w_t^2] = 1.
    \eeq
So, without the endpoint condition, the optimal $w$ is $w_t\equiv 1$, and \eqref{eq:is-budget-constraint} reduces to \eqref{eq:sp-budget-constraint}.
    However, as a result of the endpoint condition, the projection of $\Adm^{\Ising}(\alpha,p)$ onto $\sS$ is a strict subset of $\Adm^{\sph}(\alpha,p)$: indeed, we refer to \cite{bogpinprogress} for explicit analytic conditions that hold for all $\mu\in\Adm^{\Ising}(\alpha,p)$ but not for all $\mu\in\Adm^{\sph}(\alpha,p)$.
\end{rmk} 

\begin{rmk}[correlation function $\chi$ and interpretation of $p$]
\label{r:intro.interpretation.p}
For an algorithm $\cA$, we define the \emph{correlation function} $\chi$ to be the expected overlap between the outputs of the algorithm given $t$-correlated inputs; see \eqref{e:p.corr.overlap}. In the hardness analysis, $p$ essentially enters as the functional inverse of $\chi$. The choice $p\equiv 1$ in the symmetric case corresponds to $\chi\equiv 0$, i.e., the outputs $\bbE[(\cA(\bG),\cA(\bG'))/N]\approx 0$ are approximately orthogonal for any correlation level $t<1$.\footnote{This approximate orthogonality cannot hold for $L$ fixed, but is possible if $L\to\infty$ slowly with $N$. Roughly speaking, this corresponds to algorithms taking $L$ sufficiently large depending on the error parameter $\iota$ of Definition~\ref{d:intro.relaxed.domains}.}
\end{rmk}

\subsection{Further background and related works} \label{subsec:literature}
In this subsection we review some of the prior literature surrounding our results. 

\subsubsection{Statistical motivation for projection pursuit}

Following work of Kruskal \cite{kruskal1972linear}, the projection pursuit problem was introduced by Friedman and Tukey in \cite{friedman1974projection} as a method of exploratory data analysis; see also the survey \cite{huber1985projection}.
The idea is that hidden signals of interest should correspond to directions with non-gaussian projections.
Several more recent works including \cite{blanchard2006search,virta2016projection,koldovsky2018gradient,goyal2019non} have developed methods to identify such special directions.
The study \cite{bickel2018projection} and subsequent work \cite{montanari2022overparametrized} showed that gaussian point clouds admit non-gaussian projections with high probability in the high-dimensional regime $M\asymp N$. This raises the question of which projections \textbf{are} suggestive of an underlying structure in the data. From the perspective of this paper, such projections are exactly those not contained in the sets $\cM(\alpha)$  described by Definition~\ref{d:achievable-msrs}.

\subsubsection{Satisfiability threshold and free energy}

The \textbf{(half-space) spherical perceptron} constraint satisfaction problem asks, for a fixed \emph{margin} parameter $\kappa\in\bbR$, to understand the random intersection
\beq
\label{eq:perceptron-intersection}
S^{\sph}(\bG;\kappa)
=
\bigg\{
\bx\in S_N: \frac{(\bg^a,\bx)}
	{N^{1/2}}\in [\kappa,\infty)~\forall 1\leq a\leq M
\bigg\}\,.
\eeq
The \textbf{(half-space) Ising perceptron}
 $S^{\Ising}(\bG;\kappa)$ is defined analogously, replacing $S_N$ with $\Sigma_N$. Note that \eqref{eq:perceptron-intersection} can be viewed as a zero-temperature variant of \eqref{e:hamiltonian}.\footnote{While the difference between the zero and positive temperatures models \eqref{eq:perceptron-intersection}, \eqref{e:hamiltonian} often poses technical challenges, it is benign for satisfiability thresholds thanks to a rounding trick introduced by \cite{ding2019capacity}. This is explained in our companion paper \cite{bogpinprogress}, which studies consequences for algorithmic satisfiability thresholds in zero-temperature perceptron models.} 
The spherical perceptron was introduced in \cite{wendel1962problem,cover1965geometrical}, where it was shown that $M/N\to 2$ is the critical threshold at which the random set $S^{\sph}(\bG;\kappa=0)$ transitions from being non-empty to empty with high probability (i.e., the satisfiability threshold). Subsequent work in statistical physics \cite{gardner1987maximum,gardner1988optimal,gardner1988space,krauth1989storage} used the non-rigorous replica method to formulate a variety of predictions on the satisfiability threshold and other aspects of the model, for both the spherical and Ising settings. More recently, \cite{franz2016simplest} has studied the full replica-symmetry breaking phase of the problem in a similar vein, with connections to sphere packing.

In terms of rigorous progress,
for the \emph{spherical} perceptron, an important work of
\cite{shcherbina2003rigorous} computed the asymptotic volume of $S^{\sph}(\bG;\kappa)$ for all $\kappa\geq 0$, and showed that it matches the replica symmetric physics heuristics. As detailed in \cite{talagrand2010mean,MR2731561}, a series of works by Talagrand established other fundamental results for this model, including replica symmetry under weak disorder and concentration estimates.
Later, \cite{stojnic2013another} gave a simple proof of the satisfiability threshold for the spherical perceptron for $\kappa\geq 0$, introducing what has become known as the convex Gordon minimax theorem \cite{thrampoulidis2018precise}; he also showed that the replica symmetric satisfiability threshold is incorrect for the spherical perceptron with negative $\kappa$ \cite{stojnic2013negative}.
More recently, the satisfiability threshold for the \emph{Ising} perceptron, as conjectured by \cite{krauth1989storage} via replica symmetric heuristics, was verified at $\kappa=0$ by \cite{ding2019capacity,huang2024capacity}, modulo an explicit numerical condition. Results on Ising perceptron models in the high-temperature regime were obtained by 
\cite{talagrand1999intersecting,MR3024566,bolthausen2022gardner}. A variety of sharp threshold results have also been recently established for Ising perceptron models, including \cite{xu2021sharp,nakajima2023sharp,minzer2023perfectly,altschuler2023zero}.

We also mention a separate line of work on high-dimensional logistic regression, which studies models closely related to the random perceptron but with a hidden signal term. \cite{candes2020phase} established a sharp phase transition for existence of the maximum likelihood estimator in the proportional regime, while \cite{sur2019modern} characterized its asymptotic bias and variance and the limiting law of the likelihood ratio in the regime where it exists.
\cite{deng2022model} derived exact asymptotic classification errors for gradient descent on logistic loss across the separability transition, relating the two regimes to the maximum likelihood and max-margin classifiers and demonstrating double descent.
For a teacher--student perceptron model, \cite{aubin2020generalization} obtained exact high-dimensional generalization errors for $\ell_2$-regularized convex losses and showed that logistic and hinge losses can approach the Bayes-optimal error.
Finally, \cite{chardon2024finite} proved sharp non-asymptotic guarantees for existence of the maximum likelihood estimator and its excess logistic risk. 

\subsubsection{Previous algorithms}

The spherical half-space perceptron with $\kappa\geq 0$
can be reformulated as a convex problem (see e.g. \cite{stojnic2013another}), and so solutions can be found efficiently whenever they exist. 

For the Ising half-space perceptron, Kim and Roche \cite{kim1998covering} gave the first non-trivial algorithm. Their procedure efficiently finds a solution in the regime $M/N \le \alpha_0$ (a small constant). Significantly, this algorithm falsified statistical physics heuristics from \cite{zdeborova2008constraint,huang2014origin} that solutions to the Ising perceptron should not be efficiently computable at any positive clause density. Indeed, the binary perceptron was believed to exhibit the \emph{frozen $1$-RSB} property: at any clause density $\alpha$, all but an exponentially small fraction of solutions have linear Hamming distance from all other solutions. A reconciliation of these phenomena was proposed in \cite{baldassi2015subdominant,baldassi2016unreasonable}, which argued that efficient algorithms find rare dense clusters that do not appear when analyzing the statics of the model. 

For the more tractable \emph{symmetric} binary perceptron where $[\kappa,\infty)$ is replaced by $[-\kappa,\kappa]$ in \eqref{eq:perceptron-intersection}, the frozen 1RSB physics picture is verified by 
 \cite{perkins2021frozen,abbe2022proof} following the initial work \cite{aubin2019storage} (see also the sharper results \cite{altschuler2023critical,sah2023distribution}). 
This problem is of interest in its own right, being equivalent to an average-case version of the discrepancy minimization problem.
At the same time, efficient algorithms \cite{bansal2020line,abbe2022binary} and hardness results \cite{gamarnik2022algorithms,gamarnik2023geometric} for the symmetric binary perceptron were developed, which revealed the presence of a large parameter regime in which solutions exist but cannot be found efficiently, at least using families of stable algorithms similar to those we consider.
We also mention \cite{li2024discrepancy}, which obtained improved algorithms and hardness results for a variety of parameter regimes; and \cite{barbier2024atypical}, which gave a replica-based study of the local entropy around typical solutions of positive margin.
Recent work \cite{fiedler2026mean} studies an online stochastic version of this vector balancing problem and characterizes its asymptotic limit through a mean-field stochastic control problem.

Our results lend further credence to the physics picture of algorithms finding rare dense clusters in general perceptron models, with a precise definition of ``rare dense cluster'' as dense ultrametric trees. Finally, we mention the independent work \cite{montanari2024exceptional} of Montanari and Zhou, which focuses on the spherical perceptron and studies a similar class of approximate message passing (AMP) algorithms.
In our language, their iterative AMP iterations only permit $p\equiv 1$, but may be augmented by a different type of first-stage AMP.
It is natural to conjecture, based on \cite{sellke2021optimizing,huang2024optimization}, that these classes of algorithms have the same power.
The analytical route taken by \cite{montanari2024exceptional} is different from ours, and focuses on analyzing a Hamilton-Jacobi PDE corresponding to the $p\equiv 1$ part of these algorithms.

\subsubsection{Hardness of random optimization problems}

Our proof of hardness for Lipschitz algorithms is based on the \textbf{overlap gap property (OGP)} framework introduced by \cite{gamarnik2017limits} to show suboptimality of local algorithms on sparse random graphs. 
This technique was subsequently developed in a series of work including \cite{gamarnik2017performance,rahman2017local,gamarnik2018finding,chen2019suboptimality,wein2022optimal,gamarnik2023algorithmic,bresler2022algorithmic,gamarnik2024hardness,alaoui2024near}, see also the survey \cite{gamarnik2021ogp}.
In \cite{HuangSellke2021}, the first two authors developed an extension called \textbf{branching OGP}, which gives a route to proving exact algorithmic thresholds.
This first work relied on the Guerra--Talagrand interpolation method which is possible only in even mixed $p$-spin models.
Beyond this setting, the ground state analysis required of the multi-replica spin glasses appearing in branching OGP arguments generally fell outside the scope of known techniques.\footnote{We believe the recent breakthrough \cite{chen2026free} provides another route to establish the branching OGP in the multi-species spherical spin glass setting considered in \cite{huang2023algorithmic}. In particular, \cite{chen2026free} gives the ground state energy of a general vector spin glass, which in principle upper bounds the average energy of the ultrametric constellations of solutions appearing in the branching OGP. 
(While \cite{chen2026free} studies Ising spins, spherical spins will be addressed in a forthcoming followup work by the same authors.) 
However these advances do not appear to address the projection pursuit problem studied here.}
However it was observed in \cite{huang2023algorithmic}, building on the ``uniform concentration'' idea of Subag \cite{subag2018free}, that the \emph{asymptotic family} of spin glass ensembles in branching OGP arguments could nevertheless be analyzed in the ``small step size limit'' by a different, more robust approach.
%%
% but \cite{huang2023algorithmic} later gave a more robust method based on iterative application of the ``uniform concentration'' idea of Subag \cite{subag2018free}.
See \cite{jones2022random,du2023algorithmic,bhamidi2025finding} for further uses of the branching OGP including to classes of online algorithms.

Our present work develops the technique of \cite{huang2023algorithmic} to another setting where the interpolation method is inapplicable. 
The new challenges are quite significant, and stem from the lack of spherical symmetry in the problem.
As explained in \S\ref{subsec:ideas}, the starting point is to track a Doob martingale for $\cA(\bG)$ as information about $\bG$ is gradually revealed through a gaussian channel.
Roughly speaking, in the setting of \cite{huang2023algorithmic}, the only relevant information to keep track of along this path is a finite-dimensional vector of norms in different subspaces.
In our setting, one must instead track an evolving probability measure on the real line (in fact, two of them in the Ising case).
We are able to analyze this by building the required stochastic control theory ``within'' the framework of the branching overlap gap property.
This is a substantial generalization of the iterative uniform concentration technique of \cite{huang2023algorithmic}, which allows us to derive sharp algorithmic thresholds in settings with Ising spins and more general Hamiltonians.

On the algorithmic side, the framework for optimal algorithms in glassy optimization problems was developed in ground-breaking work \cite{subag2021following,mon18}, and later extended by \cite{ams20,alaoui2022perceptron,sellke2021optimizing,alaoui2023local,chen2023local,jekel2025potential} (see also \cite{lopatto2026full}).
The idea behind these algorithms is that when the model is full RSB at zero temperature, one should find ground states by following the ultrametric tree of states, and in general one should proceed analogously to the extent possible.
Our algorithms follow this intuition, and as in \cite{huang2023algorithmic} we are able to read off such algorithms from the proof of hardness, so that the two directions match by construction.
However, as discussed in \S\ref{subsec:ideas} below, our implementation of these algorithms is new and involves a family of auxiliary disorder matrices, in order to closely mimic the proof of hardness.

\subsection{Some proof ideas}
\label{subsec:ideas}

Here we describe the main ideas in our proof, focusing primarily on the easier setting of spherical projection pursuit. Throughout our analysis, we let $(\bG(t) : 0\le t\le 1)$ be an $M\times N$ matrix-valued Brownian motion conditioned to satisfy $\bG(1)=\bG$.
Treating $\bG$ as fixed, this means
\beq\label{e:matrix.br.bridge}
\bG(t) = t\bG + \bG_{\circ}(t)
\eeq
where $\bG_{\circ}$ is an $M\times N$ matrix-valued Brownian bridge. (That is, the entries of $\bG_{\circ}$ are independent, and each entry of $\bG_{\circ}$ is a Brownian motion started from zero and conditioned to return to zero at time $t=1$.) 

Recall from \eqref{eq:proj-pursuit-def} that $\mu_{\bG}(\bx)$ denotes the empirical measure of the projection of data $\bG$ in the $\bx$ direction. Our hardness proof goes by considering the evolving family of probability measures
\beq\label{e:intro.evolving.measures}
\mu_{\bG(t)}(
\bbE[\cA_N(\bG)\,|\,\bG(t)]
)
= \frac1M\sum_{a=1}^M \delta\bigg\{
\frac{(\bg^a(t), \bbE[\cA_N(\bG)\,|\,\bG(t)])}{N^{1/2}}
\bigg\}
\bigg) \in\cP_2(\bbR)
\,,\eeq
where $\bg^a(t)$ denotes the $a$-th row of $\bG(t)$, and $\cA_N$ is an arbitrary $L$-Lipschitz algorithm (up to a small transformation implemented in Proposition~\ref{p:wlogable}).
Note that the disorder matrix $\bG(t)$ and the Doob martingale $\bbE[\cA_N(\bG)\,|\,\bG(t)]$ for the eventual solution evolve simultaneously. 
We aim to classify all possible subsequential limits of the measures \eqref{e:intro.evolving.measures} as $N\to\infty$. 

Note that \eqref{e:intro.evolving.measures} corresponds naturally to the real-valued process
	\beq\label{e:intro.X.process}
	X(t)
	= 
	\frac{(\bg^a(t), \bbE[\cA_N(\bG)\,|\,\bG(t)])}{N^{1/2}}\,,\eeq
where $a\in[M]$ is sampled uniformly at random. For a suitably preprocessed and rerandomized variant of $X$ (Section~\ref{s:rerand}), we show tightness within the space of semi-martingales (Appendix~\ref{s:kolmogorov}). We then construct discretized proxies for the drift and quadratic variation of the processes. 
These essentially capture the average drift and quadratic variation experienced by atoms at any time-space location $(t,x)\in [0,1]\times \bbR$. 
This naturally gives rise to limiting measures on time-space encoding the distribution of drift and quadratic variation.
Crucially, these measures reside in a compact space, so we may consider their $N\to\infty$ subsequential limits. 
(This argument is related to the theory of relaxed controls, which has a long history \cite{fleming1966existence,mcshane1967relaxed,gyongy1986mimicking,nicole1987compactification,haussmann1990existence,kushner1990controlled,kurtz1998existence,brunick2013mimicking}.)
Further, known results on Fokker--Planck equations show that under suitable regularity conditions, the limiting measures uniquely determine the limiting marginal distributions at each time. 
Thus, if we can characterize the constraints on the drift and quadratic variation profiles imposed by the condition that $\cA_N$ is Lipschitz, then we can expect to arrive at a stochastic control approximation for a general measure $\mu_{\bG}(\cA_N(\bG))$ output by a Lipschitz algorithm.
This is the route we follow: Appendix~\ref{a:freeprob} derives discretized analogs of these ``budget'' constraints using random matrix theory, while Section~\ref{s:sde} implements several technical smoothing procedures in order to pass the budget constraints to the $N\to\infty$ limiting processes.
In fact our analysis focuses on the harder Ising case, where a second control problem governs the Doob martingale's coordinate process 
	\beq\label{e:intro.X.process.ISING}
	X^\Ising(t)
	=
	(\bbE[\cA_N(\bG)|\bG(t)],\be_i)
	\,,\eeq
where $\be_i$ is $i$-th coordinate vector in $N$ dimensions, and $i\in[N]$ is sampled uniformly at random. The interaction between these two control problems complicates the analysis significantly, as does the consideration of a general correlation function $\chi$ (cf.\ Remark~\ref{r:intro.interpretation.p}). 
These ideas culminate in Theorem~\ref{thm:BOGP-hardness-main}, which is a close approximation to Theorem~\ref{thm:main}\ref{it:thm-main-BOGP}.

We remark that while our hardness proof builds heavily on the branching OGP of \cite{HuangSellke2021,huang2023algorithmic}, the technical connection is primarily contained within Appendix~\ref{a:freeprob}. In particular, the ultrametric trees of correlated problem instances and algorithm outputs from \cite{HuangSellke2021,huang2023algorithmic} do not explicitly appear in this paper. Instead, the construction of $\bG(t)$ and $\bbE[\cA(\bG) | \bG(t)]$ above can be understood as tracing a single root-to-leaf path of these trees. Appendix~\ref{a:freeprob} follows the approach of \cite{huang2023algorithmic}, controlling how the processes \eqref{e:intro.X.process}, \eqref{e:intro.X.process.ISING} can evolve in one (discretized) time step using the ``uniform concentration'' idea of \cite{subag2018free}.

The remainder of the paper addresses the new challenges specific to this model.
For our purposes, the upshot of the method from \cite{huang2023algorithmic} (again based on the uniform concentration idea from \cite{subag2018free}) is as follows.
To analyze the increments of $\bbE[\cA_N(\bG)|\bG(t)]$, one can naively pretend that $\bbE[\cA_N(\bG)|\bG(t)]$ is independent of $\bG(t)$ in some sense, perform explicit random matrix theory calculations to constrain its movement, and then argue that these calculations remain valid despite the lack of independence \emph{because $\cA_N$ is a Lipschitz algorithm}.
In our setting, these explicit calculations exactly yield the discretized budget constraints, which must then be pieced together and passed through various limits.

Following the hardness analysis, Section~\ref{sec:IAMP} implements matching algorithms based on incremental approximate message passing (IAMP) as introduced in \cite{mon18,ams20}. 
While the basics are recalled there, we mention for now that these algorithms consist of a combination of $O(1)$ multiplications by $N\times M$ and $M\times N$ gaussian matrices, and entry-wise application of Lipschitz functions between real spaces of constant dimension.
We essentially show that any Markovian diffusion with Lipschitz coefficients satisfying the budget constraints corresponds to an explicit IAMP algorithm.
Unlike previous AMP algorithms, our implementation involves the full family of disorder matrices $\bG(t)$ rather than just the input matrix $\bG$ (such a family with the correct law can be generated artifically given just the endpoint $\bG=\bG(1)$).
This is the first AMP algorithm with this property, and enables the algorithms and hardness results to match in a canonical way.

The hardness analysis in Sections~\ref{s:rerand} and \ref{s:sde} shows (Theorem~\ref{thm:BOGP-hardness-main}) that all empirical measures attainable by a Lipschitz algorithm are the endpoint measure of a controlled SDE, while the IAMP analysis in Section~\ref{sec:IAMP} shows (Theorem~\ref{thm:IAMP-main}) that any endpoint measure of another controlled SDE is attainable by a Lipschitz algorithm.
The controlled SDEs in these two theorems involve slightly different assumptions on the control processes and various error terms, and can both be thought of as approximate versions of the ``ideal'' SDE \eqref{eq:main-process-intro-general}, \eqref{eq:X-ising-intro}.
In Section~\ref{sec:alternate-diffusions}, we show that the sets of endpoint measures of these SDEs coincide as the various error parameters are sent to $0$, and that $p\equiv 1$ suffices to attain all feasible \emph{symmetrized} endpoint measures.
This proves most of our main results Theorems~\ref{thm:main} and \ref{thm:symmetric}, but with a somewhat weaker hardness claim: the hardness analysis in Theorem~\ref{thm:BOGP-hardness-main} controls the set of all subsequential $\bbW_2$-limits of feasible empirical measures, which is weaker than the notion of ``confinement'' in Theorems~\ref{thm:main} and \ref{thm:symmetric} (see Definitions~\ref{d:attain} and \ref{d:sym-attain}).
Section~\ref{s:confinement} upgrades these hardness claims to complete the proofs of Theorems~\ref{thm:main} and \ref{thm:symmetric}.
We refer the reader to Section~\ref{ss:intro.organization} for a more detailed summary of Sections~\ref{sec:alternate-diffusions}--\ref{s:confinement}. 

\subsection{Formal definitions and theorem statements}
\label{ss:intro.lip.algs}

In this subsection we present precise definitions for Lipschitz algorithms and approximations of target measures, as well as formal statements of the main results that were informally presented in \S\ref{ss:intro.results.informal}.

\begin{dfn}[Lipschitz algorithms]
\label{d:Lip}
For $\cA$ as in \eqref{eq:algorithms-as-maps},
we say that $\cA$ is an \textbf{$L$-Lipschitz algorithm} if it satisfies the following two conditions:
\begin{enumerate}[(a)]\item \label{it:L-lip-condition}
$\cA$ defines an $L$-Lipschitz mapping, that is,
    \[
        \Big\|\cA(\bG,\bg^{\aux})-\cA(\bG',(\bg')^{\aux})\Big\| 
        \le L \Big\|(\bG,\bg^{\aux}) - (\bG',(\bg')^{\aux})\Big\|\,,
    \]
where $\|\cdot\|$ always refers to Euclidean norm unless specifically indicated otherwise.
    \item \label{it:expectation-condition}
    We have 
    $\bbE[\|\cA(\bG,\bg^{\aux})\|^2] = N$.
\end{enumerate} We say $\cA$ is a \textbf{deterministic $L$-Lipschitz algorithm} if it satisfies the above conditions with $\cA=\cA(\bG)$, i.e., without dependence on any auxiliary randomness $\bg^\aux$. 
\end{dfn}
\begin{rmk}
    As discussed in Remark~\ref{r:deterministic-lipschitz} below, the power of this class of algorithms remains unchanged if we do not have the random seed $\bg^{\aux}$.
    We include this seed for additional flexibility: Proposition~\ref{p:wlogable} uses it to ensure additional properties of $\cA$ without loss of generality, while in Section~\ref{sec:IAMP} it will be technically useful to consider certain external random variables to be ``auxiliary''.
\end{rmk}

\begin{rmk}
\label{rmk:omega-general-seed}
    We may even more generally consider algorithms $\cA(\bG,\bg^{\aux}; \omega)$ for $\omega$ an arbitrary random variable independent of $(\bG,\bg^{\aux})$, where we require conditions
    \eqref{it:L-lip-condition} and \eqref{it:expectation-condition}
    from Definition~\ref{d:Lip} to hold conditional on almost all $\omega$ (with uniform choice of $L$).
    However, this generalization does not yield an increase in algorithmic power because we can fix the best seed $\omega$.
\end{rmk}

We next define relaxed domains for the spherical and Ising problems.

\begin{dfn}[relaxed domains]\label{d:intro.relaxed.domains}
Let $\dist(\cdot,\cdot)$ denote the Euclidean distance in $\R^N$, and define
\beq\label{eq:relaxed-domain}
\begin{aligned}
    S_N(\iota)
    &\equiv\bigg\{\bx\in\R^N :
    \frac{\dist(\bx,S_N)}{N^{1/2}} \le \iota
    \bigg\}\,, &
    \Sigma_N(\iota)
    &\equiv\bigg\{\bx\in\R^N :
    \frac{\dist(\bx,\Sigma_N)}{N^{1/2}} \le \iota
    \bigg\}\,.
\end{aligned}\eeq
We call $S_N(\iota)$ the \textbf{$\iota$-relaxed spherical domain},
and $\Sigma_N(\iota)$  the \textbf{$\iota$-relaxed Ising domain}.\footnote{While these relaxed domains are needed in the Ising setting for the notion of a Lipschitz algorithm to make sense, any point in $S_N(\iota)$ or $\Sigma_N(\iota)$ can be rounded to a point in $S_N$ or $\Sigma_N$ with nearly the same performance; see Remark~\ref{r:relaxed-wlog}.}
\end{dfn}

Throughout the following, we let $\cP_2(\bbR)$ denote the set of square-integrable Borel probability measures on $\bbR$, metrized by the $2$-Wasserstein distance
\beq
\label{def:W2}
\bbW_2(\mu,\mu')
=\inf\bigg\{
\bbE^{\pi}[(x-x')^2]^{1/2}
:\pi\in \Pi(\mu,\mu')
\bigg\}\,,
\eeq
where $\Pi(\mu,\mu')$ denotes the set of couplings of the measures $\mu$ and $\mu'$. (Recall that a coupling $\pi$ of measures $\mu,\mu'$ is any law of a pair of random variables $(X,X')$ such that marginally $X\sim\mu$ and $X'\sim\mu'$.) For $\mu \in \cP_2(\bbR)$ and $\cM \subseteq \cP_2(\bbR)$, we denote the point-to-set distance
\beq
\label{def:W2-to-set}
\bbW_2(\mu,\cM)
=\inf\bigg\{
\bbW_2(\mu,\mu')
:\mu'\in \cM\bigg\}\,.
\eeq
We now formalize what it means for an algorithm \eqref{eq:algorithms-as-maps} to succeed at the approximation task $\mu_{\bG}(\bx) \approx \mu$: 

\begin{dfn}\label{d:attain}
    In the spherical (resp.\ Ising) projection pursuit problem, we say that an algorithm $\cA$ (as in \eqref{eq:algorithms-as-maps}) \textbf{$(\ioeps,\gamma)$-attains} a target measure $\mu\in \cP_2(\bbR)$ if the output $\bx = \cA(\bG,\bg^{\aux})$ satisfies, with probability at least $\gamma$, \emph{both} of the following properties:
    \begin{enumerate}[(i)]
         \item \label{it:coord-profile-succeed} 
         $\bx \in S_N(\iota)$ (resp.\ $\Sigma_N(\iota)$);
         \item \label{it:proj-pursuit-succeed}
         $\bbW_2(\mu,\mu_{\bG}(\bx)) \le \ioeps$.
    \end{enumerate}
For $\fq \in [1,2)$, we say that $\cA$ is \textbf{$(\fq,\iotamsr,\iotasol,\gamma)$-confined} to a subset $\cM\subseteq \cP_2(\bbR)$ if the output $\bx$ satisfies, with probability at least $\gamma$, \emph{either} of the following properties:
    \begin{enumerate}[(i)]
        \item \label{it:coord-profile-confined} $\bx \notin S_N(\iotasol)$ (resp.\ $\Sigma_N(\iotasol)$);
        \item \label{it:proj-pursuit-confined} $\bbW_\fq(\mu_{\bG}(\bx),\cM) \le \iotamsr$.
    \end{enumerate}
\end{dfn}

We will be interested in the setting where the Lipschitz parameter $L$ of Definition~\ref{d:Lip} is a constant, which may be arbitrarily large but must be independent of the dimension $N$. Likewise, the parameters $\iota,\iotamsr,\iotasol$ of Definition~\ref{d:attain}
will be small constants, also independent of $N$. The success probability $\gamma$ will be exponentially close to $1$. In this regime, our main result for Ising projection pursuit (Theorem~\ref{t:main.INFORMAL}\eqref{t:main.INFORMAL.Ising})
identifies the set $\ocM^{\Ising}(\alpha)$ of achievable measures, in the sense that all measures in $\ocM^{\Ising}(\alpha)$ can be \emph{attained} by an $L$-Lipschitz algorithm, while all $L$-Lipschitz algorithms are confined to $\ocM^{\Ising}(\alpha)$:

\begin{thm}[\textbf{main result}]
\label{thm:main}
For the spherical projection pursuit problem \eqref{eq:proj-pursuit-def}, we have the following:
    \begin{enumerate}[(a)]
        \item 
        \label{it:thm-main-IAMP}
        If $\mu \in \ocM^{\sph}(\alpha)$, for any $\iota > 0$ there exist $L$ large enough and $c$ small enough such that for sufficiently large $N$, there exists an $L$-Lipschitz algorithm that $(\iota,1-e^{-cN})$-attains $\mu$.
        \item 
        \label{it:thm-main-BOGP}
        For any $\fq \in [1,2)$ and $\iotamsr > 0$, there exists small $\iotasol > 0$ such that the following holds.
        For any $L$, there exists small $c>0$ such that for sufficiently large $N$, all $L$-Lipschitz algorithms are $(\fq,\iotamsr,\iotasol,1-e^{-cN})$-confined to the set $\ocM^{\sph}(\alpha)$.
        \item 
        \label{it:thm-main-concave}
We have $\ocM^{\sph}(\alpha) = \ocM^{\sph, \concave}(\alpha)$, i.e., concave $p$ suffice to achieve all endpoint measures of \eqref{eq:main-process-intro-general}.
    \end{enumerate}
The analogous results hold in the Ising setting, replacing $\ocM^{\sph}(\alpha), \ocM^{\sph, \concave}(\alpha)$
with $\ocM^{\Ising}(\alpha), \ocM^{\Ising, \concave}(\alpha)$.
\end{thm}

We note the implications of Theorem~\ref{thm:main} for the perceptron model:
\begin{cor}[\textbf{algorithmic threshold for perceptron model}]
\label{cor:alg-for-optimization}
For $\phi\in C_b(\bbR;\bbR)$, recall the perceptron Hamiltonian $H_N$ from in \eqref{e:hamiltonian}. We have the following:
    \begin{enumerate}[(a)]
        \item \label{it:cor-optimization-IAMP} 
        For any $\iota> 0$, there exists $L$ large enough and $c$ small enough such that the following holds.
        For sufficiently large $N$, there exists a $L$-Lipschitz algorithm $\cA_N$ such that for $\bx = \cA_N(\bG,\bg^{\aux})$,
        \[
            \bbP\bigg(\bx \in \Sigma_N(\iota), H_N(\bx) \ge M(\ALG-\iota)\bigg) \ge 1-e^{-cN}\,.
        \]
        \item \label{it:cor-optimization-BOGP}
        For any $\iotaval > 0$, there exists sufficiently small $\iotasol > 0$ such that the following holds.
        For any $L$, there exists $c$ small enough that for any $L$-Lipschitz algorithm $\cA_N$, and $\bx = \cA_N(\bG,\bg^{\aux})$,
        \[
            \bbP\bigg(\bx \in \Sigma_N(\iotasol), H_N(\bx) \ge M(\ALG+\iotaval)\bigg) \le e^{-cN}\,.
        \]
    \end{enumerate}
In the above, $\ALG=\ALG^{\Ising}(\phi,\alpha)$ is defined by \eqref{eq:def-ALG-Ising} in the Ising setting, and analogously in the spherical setting
by replacing
$\ocM^{\Ising}(\alpha)$ with $\ocM^{\sph}(\alpha)$.
\end{cor}

We next formally present the specialization of our results in the symmetrized setting. We now let $\cP_2(\bbR_{\ge 0}) \subseteq \cP_2(\bbR)$ be the set of square-integrable Borel probability measures on $\bbR_{\ge 0}$, again metrized by the $2$-Wasserstein metric $\bbW_2$.

\begin{dfn}
    \label{d:sym-attain} In the symmetric spherical (resp.\ Ising) projection pursuit problem, we say that $\cA$ \textbf{$(\ioeps,\gamma)$-symmetrically attains} a target measure $\mu\in\cP_2(\bbR_{\ge 0})$ the output $\bx = \cA(\bG,\bg^{\aux})$ satisfies both of the following properties
    with probability at least $\gamma$:
    \begin{enumerate}[(i)]
         \item \label{it:coord-profile-succeed-sym} 
         $\bx \in S_N(\iota)$ (resp.\ $\Sigma_N(\iota)$);
         \item \label{it:proj-pursuit-succeed-sym}
         $\bbW_2(\mu,\mu_{\bG,\sym}(\bx)) \le \ioeps$.
    \end{enumerate}
We say that $\cA$ is \textbf{$(\fq,\iotamsr,\iotasol,\gamma)$-symmetrically confined} to a subset $\cM\subseteq \cP_2(\bbR_{\ge 0})$ if the output $\bx$ satisfies both of the following properties with probability at most $1-\gamma$:
    \begin{enumerate}[(i)]
        \item \label{it:coord-profile-confined-sym} $\bx \in S_N(\iotasol)$ (resp.\ $\Sigma_N(\iotasol)$);
        \item \label{it:proj-pursuit-confined-sym} $\bbW_\fq(\mu_{\bG,\sym}(\bx),\cM) \ge \iotamsr$.
    \end{enumerate}
\end{dfn}   

\begin{thm}[\textbf{main result for symmetrized projection pursuit}]
    \label{thm:symmetric}
The following holds in the Ising (resp. spherical) cases.
    \begin{enumerate}[(a)]
        \item 
        \label{it:thm-symmetric-IAMP}
        If $\mu \in \ocM^{\Ising,\sym}(\alpha)$ (resp. $\ocM^{\sph,\sym}(\alpha)$), for any $\iota > 0$ there exist $L$ large enough and $c$ small enough such that for sufficiently large $N$, there exists an $L$-Lipschitz algorithm that $(\iota,1-e^{-cN})$-symmetrically attains $\mu$.
        \item 
        \label{it:thm-symmetric-BOGP}
        For any $\fq \in [1,2)$ and $\iotamsr > 0$, there exists small $\iotasol > 0$ such that the following holds.
        For any $L$, there exists small $c>0$ such that for sufficiently large $N$, all $L$-Lipschitz algorithms are $(\fq, \iotamsr,\iotasol,1-e^{-cN})$-symmetrcally confined to $\ocM^{\Ising,\sym}(\alpha)$.
    \end{enumerate} 
\end{thm}

On the algorithmic achievability side, it suffices to use the class of incremental approximate message passing algorithms introduced in \cite{mon18} and used in \cite{ams20,sellke2021optimizing,alaoui2022perceptron,montanari2024exceptional} (see \S\ref{subsec:ideas}).
These algorithms are suitably Lipschitz and require only $O_{\iota,\epsilon}(N^2)$ arithmetic operations, assuming query access to suitable Lipschitz, Markovian control functions $(p(t),\sigma(t,x),b(t,x),w(t,x))$ which are independent of $N$.
They are introduced and analyzed in Section~\ref{sec:IAMP}.

\subsection{Organization of this paper} 
\label{ss:intro.organization}
The rest of the paper is structured as follows.

\begin{itemize}
\item In Section~\ref{s:rerand}, we consider the processes $X^N$ and $X^{N,\Ising}$ of \eqref{e:intro.X.process} and \eqref{e:intro.X.process.ISING},
     which can be associated to any Lipschitz algorithm $\cA_N$.
    The stochastic control descriptions $X$ and $X^\Ising$ in \eqref{eq:main-process-intro-general} and \eqref{eq:X-ising-intro} will later emerge as a suitable continuum limit of $X^N$ and $X^{N,\Ising}$.
    We then perform several preprocessing operations,
    which regularize $X^N, X^{N,\Ising}$ to a form $\vX^N$ suitable for taking a continuum limit.
    We state the tightness result Theorem~\ref{t:tightness}, which implies that $\vX^N$ has subsequential limits.
    We then introduce discrete local averages of the drift and quadratic variation of $\vX^N$, from which the coefficients of the limiting SDE will be extracted.
    We finally state Theorem~\ref{t:free.prob}, which shows that these discrete local averages satisfy a discrete version of the budget constraint \eqref{eq:is-budget-constraint} or \eqref{eq:sp-budget-constraint}, and provide a heuristic free-probability derivation of a simplified version of this constraint.
    The tightness result Theorem~\ref{t:tightness} is proved in Appendix~\ref{s:kolmogorov} via Kolmogorov estimates on $\vX^N$, while the discrete budget contraint Theorem~\ref{t:free.prob} is proved in Appendix~\ref{a:freeprob} via gaussian comparison inequalities.
    \item In Section~\ref{s:prelim} we collect preliminary estimates concerning gaussian matrices, subgaussian vectors, and Lipschitz functions on gaussian space.
    We then define and study the correlation function $\chi$ of a Lipschitz algorithm (cf.\ Remark~\ref{r:intro.interpretation.p}).  In Proposition~\ref{p:correlation-fn-ub}, we prove an upper bound on the derivative of $\chi$.
    In Proposition~\ref{p:wlogable}, we show that any Lipschitz algorithm can be perturbed, with little change in its output or performance, so that the derivative of $\chi$ is also bounded away from zero.
    This regularization ensures that the inverse correlation function $p \equiv \chi^{-1}$ used to parametrize the processes of Sections \ref{s:rerand} and \ref{s:sde} is well behaved.
    \item In Section~\ref{s:sde}, we prove Theorem~\ref{thm:BOGP-hardness-main}, which shows that for any sequence of $L$-Lipschitz algorithms $(\cA_N)^\circ$, if $(\cA_N)$ are their perturbations from Proposition~\ref{p:wlogable}, then any subsequential $\bbW_2$-limit of the annealed empirical measure $\mu(\cA_N) = \bbE \mu_\bG(\cA_N(\bG,\bg^\aux))$ is approximated by the endpoint law of a controlled SDE \eqref{e:sde.x}--\eqref{e:sde.x.ising}.
    This SDE resembles the SDE \eqref{eq:main-process-intro-general}, \eqref{eq:X-ising-intro} above; see the discussion of Section~\ref{sec:alternate-diffusions} below for the differences. The proof in Section~\ref{s:sde} is based on showing that any subsequential limit of the process $\vX^N$ from Section~\ref{s:rerand} (recall Theorem~\ref{t:tightness}) approximates this controlled SDE; note that $\vX^N$ is constructed from the perturbed algorithm $\cA_N$.
    We will argue that if the occupation measures for the position, drift, and quadratic variation of $\vX^N$ are regular enough, it is possible to extract coefficients of a Lipschitz SDE that has the same time marginals as $\vX^N$.
    We will spatially and temporally smooth the discrete-time process $\vX^N$ in order to ensure these occupation measures have the required regularity, and to pass the discrete budget constraint derived in Theorem~\ref{t:tightness} to the continuum limit.
    Using similar ideas, we prove Theorem~\ref{thm:SDE-smoothing-general}, which provides a related regularization result: a diffusion with general progressively measurable controls satisfying the budget constraint \eqref{eq:is-budget-constraint} can be approximated by a Markovian SDE with bounded Lipschitz coefficients that approximately satisfies the same budget constraint.
    \item In Section~\ref{sec:IAMP}, we construct matching Lipschitz algorithms using incremental approximate message passing (IAMP).
    We first prove Theorem~\ref{thm:IAMP-main}\ref{i:IAMP-main-main}, which provides an IAMP algorithm to achieve any endpoint measure of a controlled SDE \eqref{e:IAMP.X.SDE}--\eqref{e:IAMP.Y.SDE} that also resembles \eqref{eq:main-process-intro-general}, \eqref{eq:X-ising-intro}.
    We then prove Theorem~\ref{thm:IAMP-main}\ref{i:IAMP-main-centered}, which shows that the IAMP algorithm can be arranged to have mean zero.
    Finally we prove Theorem~\ref{thm:IAMP-main}\ref{i:IAMP-main-chaotic}, which shows that, if we are only interested in the \textbf{symmetrized} empirical measure $\mu_{\bG,\sym}(\cA(\bG,\bg^\aux))$ (as defined by \eqref{eq:symmetrized-empirical-measure-intro}), then we can arrange for the IAMP algorithm to have correlation function $\chi(1-\iota) \le \iota$, i.e. the algorithm maps $(1-\iota)$-correlated inputs to nearly orthogonal outputs.
    These last two facts will greatly simplify the final SDE characterization we obtain in Section~\ref{sec:alternate-diffusions}.
    \item In Section~\ref{sec:alternate-diffusions}, we prove Theorem~\ref{thm:control.problems.main}, which shows that the set of possible subsequential $\bbW_2$-limits of $\mu(\cA_N)$, for $(\cA_N)$ 
    a sequence of Lipschitz algorithms,  is precisely the set of endpoint measures of \eqref{eq:main-process-intro-general}, \eqref{eq:X-ising-intro} under the budget constraint \eqref{eq:is-budget-constraint}.
    The SDE \eqref{eq:main-process-intro-general}, \eqref{eq:X-ising-intro} can be thought of as the ``ideal'' version of the SDEs \eqref{e:sde.x}--\eqref{e:sde.x.ising} and \eqref{e:IAMP.X.SDE}--\eqref{e:IAMP.Y.SDE} considered in Sections \ref{s:sde}--\ref{sec:IAMP}: these SDEs approximate the ideal SDE in different ways and feature a starting time $q_0 \in [0,1)$ that is not necessarily $0$.
    In the ideal SDE, all the approximation errors and $q_0$ are set to $0$, and the controls are general progressively measurable processes.
    To prove Theorem~\ref{thm:control.problems.main}, we define a collection of control and measure classes interpolating between these three SDEs and prove a cycle of inclusions among these classes.
    Theorems~\ref{thm:BOGP-hardness-main} and~\ref{thm:IAMP-main} will justify key steps of this inclusion loop.
    Note that in Theorem~\ref{thm:control.problems.main}, $\cA_N$ is the original algorithm, not its perturbed algorithm; while Theorem~\ref{thm:BOGP-hardness-main} pertains to its perturbed variant, Proposition~\ref{p:wlogable} implies that $\cA_N$ and its perturbation achieve similar measures.
    Theorem~\ref{thm:IAMP-main}\ref{i:IAMP-main-centered} will crucially imply that $q_0=0$ suffices to attain all endpoint measures.
    We also prove Theorem~\ref{thm:control.problems.sym}, which states that in symmetrized projection pursuit, the set of subsequential $\bbW_2$-limits by Lipschitz algorithms is the set of endpoint measures of \eqref{eq:main-process-intro-p=1}, \eqref{eq:X-ising-intro}, where we have set $p\equiv 1$.
    This is proved by an analogous inclusion loop, where we use Theorem~\ref{thm:IAMP-main}\ref{i:IAMP-main-chaotic} to show that $p\equiv 1$ suffices to attain all symmetrized endpoint measures.
    \item In Section~\ref{s:confinement}, we deduce the main Theorems~\ref{thm:main} and~\ref{thm:symmetric} from Theorems~\ref{thm:control.problems.main} and~\ref{thm:control.problems.sym}.
    While the achievability parts of the main theorems essentially follow from the results of Section~\ref{sec:alternate-diffusions}, the confinement parts do not, as the sequence of measures $\mu(\cA_N)$ is not necessarily uniformly integrable in $\bbW_2$.
    To address this, we define a truncated version $\cA_{N[K]}$ of $\cA_N$ (Definition~\ref{d:truncated.lip.alg}) that truncates consistently large entries of $\mu_\bG(\cA_N(\bG,\bg^\aux))$.
    The truncation ensures that the measures $\mu(\cA_{N[K]})$ have subsequential $\bbW_2$-limits, which are characterized by the results of Section~\ref{sec:alternate-diffusions}.
    We then show (Lemma~\ref{l:lip.trunc.effect.on.Wq}) that the truncation has a small effect on $\mu(\cA_N)$ in $\bbW_{\fq}$ for any $\fq \in [1,2)$.
    This proves the confinement assertions in Theorems~\ref{thm:main} and~\ref{thm:symmetric}.
    Finally, we prove Corollary~\ref{cor:alg-for-optimization} by optimizing linear functionals over the characterized sets of measures.
\end{itemize}

\subsection*{Acknowledgements} 
We would like to thank Andrea Montanari and Kangjie Zhou for discussions about their work \cite{montanari2024exceptional}, and Christian Fiedler, Joe Jackson, Daniel Lacker, and Jonathan Niles-Weed, for discussions about their work \cite{fiedler2026mean}.
We are also grateful to Ahmed El Alaoui, Kavita Ramanan, Tselil Schramm, Youngtak Sohn, and Lenka Zdeborov\'a for motivating discussions.
We gratefully acknowledge the support of NSF CAREER grant DMS-1940092, 
NSF-Simons collaboration grant DMS-2031883,
the Solomon Buchsbaum Research Fund at MIT (BH and NS); 
the Stanford Science Fellowship, NSF Mathematical Sciences
Postdoctoral Fellowship and Google PhD Fellowship (BH); NSF grant DMS-2347177 (NS);
and 
NSF CAREER grant DMS-2540987 and a Sloan Research Fellowship (MS). 
All ideas in this paper are human-generated, and all the writing was done by the human authors. AI was used in the writing of this paper only for light proofreading and copy-editing.

\fi

\pagebreak\section{Overview of limiting procedure}
\label{s:rerand}

\iffull
% !TEX root = main.tex

In this section we introduce the discrete-time process from which we will extract the SDE limit described in our main result Theorem~\ref{thm:main}.
We consider a Brownian bridge on the time interval $0\le t\le 1$, whose terminal value is the gaussian disorder $(\bG,\bg^\aux)$ input to the algorithm $\cA$. Our process then tracks the expected output of the algorithm, conditional on the Brownian bridge up to time $t$. This section is organized as follows:
\begin{itemize}
\item In \S\ref{ss:X.decomp} we formally introduce the discrete-time processes, and decompose them into a small number of key components for further analysis.
\item In \S\ref{ss:preprocessing} we perform some technical preprocessing on the processes, then perform a crucial spatial rerandomization which imposes a natural filtration on the discrete-time processes.
\item In \S\ref{ss:tightness.statements} we state a tightness result, Theorem~\ref{t:tightness}, for the spatially rerandomized processes. We also introduce the (discrete-time) quantities $(\bbb,\vvv,\uuu,\rrr,\www)^{\frozen,N}$ from which we will extract the coefficients of the limiting SDE.

\item In \S\ref{ss:apriori} we collect some \textit{a~priori} estimates on the discrete-time processes, to be used in Section~\ref{s:sde}.

\item In \S\ref{ss:budget} we present Theorem~\ref{t:free.prob}, which gives budget constraints in a simplified setting without preprocessing considerations. We then bound the effect of the preprocessing operations, and thereby deduce budget constraints for the quantities of interest $(\bbb,\vvv,\uuu,\rrr,\www)^{\frozen,N}$.

\item Theorem~\ref{t:free.prob} is a central ingredient to our main result. However, its proof is quite long, and is deferred to Appendix~\ref{a:freeprob}. Instead, in \S\ref{ss:heuristic.budget} we present a heuristic derivation of a simplified version of Theorem~\ref{t:free.prob}, using free probability calculations.
\end{itemize}
Throughout what follows, we will abbreviate $\|\bG\|\equiv\|\bG\|_\op$ for the operator norm of a matrix $\bG$, and abbreviate $\|\bx\|\equiv\|\bx\|_2$ for the euclidean norm of a vector $\bx$.

\subsection{Processes tracking the algorithm}
\label{ss:X.decomp}

Let $\cA=\cA(\bG,\bg^{\aux})$ be an $L$-Lipschitz algorithm, as in Definition~\ref{d:Lip}. Let $\chi_{\cA}(p) \equiv \E R(p)$, where $R(p)$ denotes the overlap between the outputs of $\cA$ resulting from a pair of $p$-correlated inputs:
\beq
    \label{e:p.corr.overlap}
    \chi_{\cA}(p)
    \equiv \E R(p)
    \equiv \E\bigg[
        \frac{(\cA(\bG,\bg^{\aux}),\cA(p\bG+(1-p^2)^{1/2}\bG',p\bg^{\aux}+(1-p^2)^{1/2}(\bg')^{\aux}))}{N}
    \bigg]\,,
\eeq
where $\bG',(\bg')^{\aux}$ are independent copies of $\bG,\bg^{\aux}$.
We will write $\chi = \chi_{\cA}$ when $\cA$ is clear from context.  \textbf{Throughout this section, we will assume that $\chi=\chi_{\cA}$ satisfies}
    \beq\label{e:chi.deriv.bounds}
    \frac{1}{L^2}
    \le \chi'(p) \le L^2\eeq
\textbf{for all $p\in[0,1]$.} This assumption will later be justified by arguments in \S\ref{subsec:cor-func-bounds}. In particular, Proposition~\ref{p:correlation-fn-ub} tells us that the upper bound of \eqref{e:chi.deriv.bounds} is always satisfied. Proposition~\ref{p:wlogable} tells us that while a general Lipschitz algorithm $\cA^\circ$ may not satisfy the lower bound of \eqref{e:chi.deriv.bounds}, a small perturbation $\cA$ of it has similar performance and does satisfy \eqref{e:chi.deriv.bounds}. \textbf{The discussion of this section should be understood as applying to the perturbed algorithm $\cA$.}

Assuming \eqref{e:chi.deriv.bounds}, fix parameters 
    \beq\label{e:p.q}
    \begin{array}{l}
    0 \equiv p_{-1}  = p_0 < p_1 < \ldots < p_{\dmax} =1\,,\\
    0 \equiv q_{-1} \le q_0 < q_1 < \ldots < q_{\dmax}=1\,,
    \end{array}
 \eeq
such that $q_d=\chi(p_d)$ for all $d$. Let $\delta_d\equiv q_{d+1}-q_d$; we write
    \beq\label{e:delta.cutoffs}
    \delta\equiv\max\bigg\{
    \delta_d : 0\le d\le \dmax\bigg\}\,.
    \eeq
We let $p:[q_0,1]\to [0,1]$ denote the continuous function which satisfies $p(q_d)=p_d$ for $0\leq d\leq d_{\max}$ and is linear on each interval $[q_d,q_{d+1}]$; we also abbreviate $p'(q)$ for the right derivative of $p$ at $q$. 
Since $p$ is a discretized inverse of $\chi$, the mean value theorem together with \eqref{e:chi.deriv.bounds} implies
    \beq\label{e:deriv.p.bound}
    p'(q_d) \equiv 
    \frac{p(q_{d+1})-p(q_d)}{\delta_d} =
    \frac{p_{d+1}-p_d}{q_{d+1}-q_d}
    \le L^{O(1)}
    \eeq
for all $0\leq d\leq d_{\max}$. 

For each $\ell$ let $\bXi^\ell\equiv\bXi(p_\ell)$ be an $M\times N$ matrix with i.i.d.\ gaussian entries, independent of all else. We use this to define, for $0 \le d\le\dmax$,
    \beq\label{e:gaus.decomp}
    \bG(q_d)
    \equiv \sum_{\ell=0}^d (p_\ell-p_{\ell-1})^{1/2} \bXi^\ell\,.
    \eeq
Let $\bG \equiv \bG(1)$; note that this again has standard gaussian entries.
Equivalently, $\bG(t)$ for $0\le t\le 1$ is a Brownian bridge in the space of $M\times N$ matrices, started from the zero matrix and ending at $\bG(1)=\bG$, and we observe it at the discrete set of times $q_d$. We analogously define the Brownian bridge $\bg^\aux(t)$ for $0\le t\le 1$, started from the zero vector and ending that $\bg^\aux(1)=\bg^\aux$. Let
    \beq\label{e:x.q.mean}
    \bx(q_d)
    \equiv \E\Big[ \cA(\bG ,\bg^\aux)
        \,\Big|\, \cG(q_d)\Big]
    = \E\Big[ \cA(\bG ,\bg^\aux) \,\Big|\, 
        \bG(q_d)  ,\bg^\aux(q_d)  \Big]\,,
    \eeq
where $\cG(q_d)$ is the $\sigma$-algebra generated by the gaussian disorder up to level $d$:
    \beq\label{e:gaus.filt}
    \cG(q_d)\equiv \sigma\bigg(
    \bG(q_\ell) ,\bg^\aux(q_\ell) 
     : \ell \le d \bigg)\,.
    \eeq
Let $a\in[M]$ be a uniformly random index, and consider the process
    \beq\label{e:process.X}
    X\equiv X^N
    \equiv (X^N(q_d,a))_{0\le d\le \dmax}
    \equiv \bigg( \frac{(\bg^a(q_d),\bx(q_d))}
    {N^{1/2}}\bigg)_{0\le d\le \dmax}\,,
    \eeq
where $\bg^a(q_d)$ is the $a$-th row of $\bG(q_d)$.  In the Ising perceptron case, take also a uniformly random index $i\in[N]$, and consider the process
    \beq\label{e:process.X.Ising}
    X^{\Ising}
    \equiv X^{N,\Ising}
    \equiv
    \Big(X^{N,\Ising}(q_d,i)\Big)_{0\le d\le \dmax}
    \equiv 
    \Big( (\be_i,\bx(q_d))\Big)_{0\le d\le \dmax}\,,
    \eeq 
where $\be_i$ is the $i$-th standard basis vector in $\R^N$.     (We can view $X$ and $X^{\Ising}$ as being defined for all times $t\in  [q_0,1]$ by making them linear on each interval $[q_d,q_{d+1}]$.) \textbf{In Remark~\ref{rmk:tree.decomp} we explain that the above may be regarded as conceptually equivalent to the branching overlap gap property formalism introduced by \cite{HuangSellke2021}.}

\begin{rmk}[non-centered algorithms]
Recall $p_0=0$, which means that the process $X^N$ of \eqref{e:process.X} starts at $X^N(q_0,a)\equiv0$. Meanwhile, recall $\chi(p)=\E R(p)$ with $R(p)$ as in \eqref{e:p.corr.overlap}, and note that
    \[
    q_0 = \chi(p_0) = \chi(0) 
    = \frac{\|\bx(q_0)\|^2}{N}
    = \frac{\|\E \cA(\bG,\bg^\aux)\|^2}{N}\,.
    \]
Thus the process $X^{N,\Ising}$ of \eqref{e:process.X.Ising} need not start at zero, although it will be constrained to satisfy the second moment condition $\E[X^{N,\Ising}(q_0)^2]=q_0$. In Sections~\ref{sec:IAMP} and \ref{sec:alternate-diffusions}, we will argue that Lipschitz algorithms with zero mean have the same power as general Lipschitz algorithms. However, prior to Section~\ref{sec:IAMP} we 
do not make this assumption, meaning that $X^{N,\Ising}(q_0)$ can in general be non-zero.
\end{rmk}

\begin{rmk}[annealed versus quenched]\label{r:annealed.quenched}
Throughout, we write $\P$ and $\E$ for overall probability and expectation, including over the gaussian disorder.
We write $\bP\equiv \bP_\bG$ for probability conditional on the gaussian disorder matrices 
    \[\Big(\bG(q_d) , \bg^\aux(q_d)
    :0\le d\le\dmax
        \Big)\,,\]
and write $\bE\equiv\bE_\bG$ for expectation under $\bP$. We refer to estimates under $\P$ and $\E$ as \textbf{annealed}, and 
estimates under $\bP$ and $\bE$ as \textbf{quenched}.
 For example, for the Ising perceptron,
 we write $(X,X^{\Ising})\sim\bP$ to mean that we first sample $(a,i)$ uniformly from $[M]\times[N]$, and let $X=X(a,\cdot)$ and $X^{\Ising}=X^{\Ising}(i,\cdot)$ as in \eqref{e:process.X} and \eqref{e:process.X.Ising}, conditional on the gaussian disorder.
\end{rmk}

\begin{lem}\label{l:nearly.ultrametric}
Let $\cA$ be an $L$-Lipschitz algorithm, and let $\bx(q_d)$ be as defined by
\eqref{e:x.q.mean}. Then
    \[
    \P\bigg(\bigg|
    \frac{(\bx(q_d),\bx(q_\ell))}{N}-q_{\min\{d,\ell\}}
    \bigg|\ge
    \frac{1}{N^{1/3}}\bigg) 
    \le \exp(-N^{0.1})\,.
    \]

\begin{proof}
This follows easily from gaussian Lipschitz concentration (see e.g.\ Lemma~\ref{l:lip.subgaus}).
Indeed note that $\|\cA(\bG,  \bg^\aux)\|_2\leq N^{1/2}$ almost surely, while by definition $\bx(q_d)$
is an $L$-Lipschitz function of the gaussian vector $(\bG(q_d),\bg^\aux(q_d))$. It follows that for $\ell \le d$, the normalized inner product
    \[\frac{(\bx(q_d),\bx(q_\ell))}{N}\]
is an 
 $O(LN^{-1/2})$-Lipschitz function of 
    \[\Big(\bG(q_\ell),
    \bg^\aux(q_\ell),
    \bG(q_d) - \bG(q_\ell),\bg^\aux(q_d) - \bg^\aux(q_\ell)\Big)\,.\]
 Therefore it concentrates around its mean by Lemma~\ref{l:lip.subgaus}, which gives the conclusion.
\end{proof}
\end{lem}

\begin{rmk}[connection to branching overlap gap property in previous works] \label{rmk:tree.decomp} In this remark we explain that the analysis of this paper conceptually follows the branching overlap gap property introduced by \cite{HuangSellke2021}, even though some of the formalism may appear slightly different. For simplicity, we drop $\bg^\aux$ from this discussion.  Following \cite{HuangSellke2021}, 
for an algorithm $\cA$ we can define the overlap function $\chi(p)\equiv \E R(p)$ as in \eqref{e:p.corr.overlap}. As in \cite[Defn.~2.1]{HuangSellke2021}, we say that $\cA$ is \textbf{overlap concentrated} if $R(p)$ concentrates around its mean $\chi(p)$. Definition~\ref{d:Lip} is slightly \textbf{more} restrictive than overlap concentration: if $\cA$ is Lipschitz and takes values in $\{\|\bx\| \le 2N^{1/2}\}$, then the standard concentration bound for Lipschitz functions of gaussians (see e.g.\ \eqref{e:lip.subgaus}) implies that $\cA$ is overlap concentrated. Now fix parameters $p,q$ as in \eqref{e:p.q}.
Let $\mathbb{T}$ be a $k$-ary tree of depth $\dmax$, with $d$-th layer $V_d\cong [k]^d$. For each $a\in[M]$ and each vertex $u$ in $\mathbb{T}$, let $\bmeta^{a,u}$ be an independent standard gaussian vector in $\R^N$, and let
    \[\bg^{a,u}
    \equiv \sum_{\ell=1}^{|u|}
    \Big( p_\ell - p_{\ell-1} \Big)^{1/2}
    \bmeta^{a,u_{\le\ell}}\,,
    \]
where $u_{\le\ell}$ denotes the prefix of $u$ of length $\ell$. Let $\bG^u \equiv (\bg^{a,u})_{a\le M}$. 
Let $\cG^u$ denote the filtration generated by $(\bG^w)_{w\le u}$, where $w\le u$ means that $w$ is an ancestor of $u$ in the tree $\mathbb{T}$. If $v\in V_{\dmax}$ is any descendant of $u\in V_d$, then we define
    \beq\label{e:cond.exp.x}
    \bx^u
    \equiv \E\Big( 
    \cA( \bG^v ) 
    \,\Big|\, \cG^u\Big)
    = \E\Big( \cA( \bG^v ) 
    \,\Big|\, \bG^u\Big)\,.
    \eeq
Thus, for any $v\in V_{\dmax}$, the sequence of $\bx^u$ on the path from the root to $v$ is an $\R^N$-valued martingale. For any vertices $u,v$ in the tree, it follows by iterated expectations that
    \beq \label{e.def.Ruv}
    \E R^{u,v}
    \equiv \E \frac{(\bx^u,\bx^v)}{N}
    = \chi(p_{|u\wedge v|})
    = q_{|u\wedge v|}
    \,.
    \eeq
If $\cA$ is a Lipschitz algorithm in the sense of Definition~\ref{d:Lip}, then $R^{u,v}$ concentrates around its mean.  It follows that 
the configurations $(\bx^u)_u$ defined by the conditional expectations
\eqref{e:cond.exp.x} form a nearly \textbf{ultrametric tree}.
The process \eqref{e:process.X}
is equidistributed as the process
    \[
     \bigg( \frac{(\bg^{a,u_d},\bx^{u_d})}
    {N^{1/2}}\bigg)_{0\le d\le \dmax}
    \]
where $(u_0,\ldots,u_{\dmax})$ is a path in $\mathbb{T}$ from root to boundary. Previous works \cite{HuangSellke2021,huang2023algorithmic} on the branching OGP used overlap concentration of the algorithm $\cA$ to build a tree of algorithmic outputs $(\bx^u)_{u \in \mathbb{T}}$ with overlaps given by \eqref{e.def.Ruv}, and then characterized the maximum energy attainable by a typical leaf of this tree. In this paper, we instead use the Lipschitzness of $\cA$ directly to extract the SDE limit in Theorem~\ref{thm:main}, so this tree construction will not be necessary.
\end{rmk}

We first decompose $X$ into a few separate parts, as follows. 
Take indices $1\le d\le \dmax$ and $a\in[M]$, and abbreviate
    \begin{align}
    \label{e:normalized.vectors.g}
    \bar{\bg}^a \equiv \bar{\bg}^a(q_d)
    &\equiv \frac{\bg^{a}(q_d)}{(p_d)^{1/2}}\,,\\
    \label{e:normalized.vectors.x}
    \bar{\bx}
    \equiv \bar{\bx}(q_d)
    &\equiv \frac{\bx(q_d)}{(q_d)^{1/2}}\,,\\
    \label{e:normalized.vectors.y}
    \by
    \equiv \by(q_{d+1})
    &\equiv \frac{\bx(q_{d+1})
        -\bx(q_d)}{(\delta_d)^{1/2}}\,.
    \end{align}
In the case $p_d=0$, we define
$\bar{\bg}^a$ (for $a\le M$) to be a collection of i.i.d.\ standard gaussian vectors in $\R^N$ (independent of everything else).
In the case $q_d=0$ we define
$\bar{\bx} = \ind\in\R^N$. (For $p_d=0$ and $q_d=0$ the definitions of $\bar{\bg}^a$ and $\bar{\bx}$ are somewhat arbitrary, since they will always eventually be multiplied by zero, but we have made choices that are convenient for the analysis that follows.) We abbreviate $\bar{\bG}\equiv\bar{\bG}(q_d)$ for the matrix with rows $\bar{\bg}^a$. Let us also abbreviate $\bmeta^a= \bmeta^a(q_{d+1})$ for the $a$-th row of the matrix $\bXi \equiv \bXi^{d+1}$ from \eqref{e:gaus.decomp}. Note that $\bar{\bg}^a$ and $\bmeta^a$ are  standard gaussians in $\R^N$; and the norms $\|\bar{\bx}\|$ and $\|\by\|$ will be concentrated around $N^{1/2}$.

For later use, we now introduce a quantity which will appear repeatedly in some of our quenched estimates: we define the random variable
    \beq\label{e:MAX.N} 
    \MAX_N
    \equiv 1 + \max_{0\le d\le\dmax-1}
        \bigg\{
    \frac{\|\bar{\bG}(q_d)\|^2}{M}
    +\frac{\|\bXi(q_{d+1})\|^2}{M}
    +\frac{\|\bar{\bx}(q_d)\|^2}{N}
    +\frac{\|\by(q_{d+1})\|^2}{N}
    \bigg\}\,.
    \eeq
Note that this is a function of the gaussian disorder, so we can also denote it as $\MAX(\bXi)$. We will verify in Proposition~\ref{p:Y.kolmogorov} that $\MAX_N$ satisfies the bound
    \beq\label{e.MAX.bound}
    \P\Big(\MAX_N
    \ge L^{O(1)}\Big)
    \le \frac{1}{\exp(\Theta(N))}\,.
    \eeq
In some later parts of the proof we will restrict to the (high-probability) event $U_\MAX$ that $\MAX_N\le L^{O(1)}$.

Now recall \eqref{e:process.X}, and consider the increment
    \[\Delta X(q_d,a)
    \equiv 
    X(q_{d+1},a)-X(q_d,a)
    =
    \frac{(\bg^a(q_{d+1}),\bx(q_{d+1}))}{N^{1/2}}
        -\frac{(\bg^a(q_d),\bx(q_d))}{N^{1/2}}\,.\]
Recalling the notation \eqref{e:deriv.p.bound}, we shall decompose this as the sum of three terms:
    \begin{align}\nonumber
    \Delta X^\RomI(q_d,a)
    &\equiv\delta_d p'(q_d)^{1/2}
    \frac{(\bmeta^a(q_{d+1}),
        \by(q_{d+1}))}{N^{1/2}}\,,\\ 
    \nonumber
    \Delta X^\RomII(q_d,a)
    &\equiv
    [\delta_d p(q_d)]^{1/2}
    \frac{
    (\bar{\bg}^a(q_d),\by(q_{d+1}))}
        {N^{1/2}}\,,\\ 
    \Delta X^\RomIII(q_d,a)
    &\equiv
    [\delta_d q_d p'(q_d)]^{1/2}
    \frac{
    (\bmeta^a(q_{d+1}),\bar{\bx}(q_d))}{N^{1/2}}\,.
    \label{e:X.decomp}
    \end{align}
For $\sigma\in\{\RomI,\RomII,\RomIII\}$ we let $X^\sigma$ be the piecewise linear process with increments $\Delta X^\sigma$, with initial values $X^{\RomI}(q_0) =X^{\RomII}(q_0)=X^{\RomIII}(q_0)=0$.
This gives $X = X^{\RomI}+X^{\RomII}+X^{\RomIII}$. For the Ising perceptron, we additionally have, for each $i\in[N]$,
    \beq\label{e:Delta.X.Ising}
    \Delta X^{\Ising}(q_d,i)
    \equiv X^{\Ising}(q_{d+1},i)-X^{\Ising}(q_d,i)
    = \Big(\be_i, \bx(q_{d+1})-\bx(q_d)\Big)
    = (\delta_d)^{1/2} \by_i\,.
    \eeq
We hereafter write, for $ (a,i)\in[M]\times[N]$,
    \begin{align}\nonumber
    \vX(q_d, (a,i))
    &\equiv
    \Big(\vX^{\RomI:\RomIII}(q_d, a),
    X^\Ising(q_d,i)\Big)\\
    &\equiv
    \Big(X^{\RomI}(q_d,a),
    X^{\RomII}(q_d,a),
    X^{\RomIII}(q_d,a),
    X^{\Ising}(q_d, i)
    \Big)
    \label{e:vX}
    \end{align}
(for the spherical perceptron, we do not keep track of $X^\Ising$, so we do not need the $i\in[N]$ index).

\textbf{From the decomposition \eqref{e:X.decomp} and \eqref{e:Delta.X.Ising}, the intuition is that, in some appropriate limit, $\Delta X^\RomI$ is a drift term while 
$\Delta X^\sigma$ for $\sigma\in\{\RomII,\RomIII,\Ising\}
$ are martingale terms.} For example, $\Delta X^\RomII$ gives a martingale term because $\bar{\bg}^a$ is measurable with respect to $\cG(q_d)$ while $\by$ has mean zero conditional on $\cG(q_d)$, so overall the inner product $(\bar{\bg}^a,\by)$ has mean zero conditional on $\cG(q_d)$. By contrast, in the $\Delta X^\RomI$ term, both $\bmeta^a$ and $\by$ are not measurable with respect to $\cG(q_d)$, so their inner product can have non-zero mean conditional on $\cG(q_d)$, thus contributes a drift.  \textbf{In the appropriate limit, we expect the decomposition of $X$ to satisfy an SDE which is roughly of the form}
    \beq\label{e:intro.SDE}
    \begin{pmatrix}
    dX^{\RomI}(t)\\
    dX^{\RomII}(t)\\
    dX^{\RomIII}(t)\\
    dX^{\Ising}(t)
    \end{pmatrix}
    = \begin{pmatrix}
    p'(t)^{1/2} b_t\\
    0\\
    0\\
    0
    \end{pmatrix}\,dt
    +\begin{pmatrix}
    0&0&0&0\\
    0& p(t)^{1/2}[v_t-(u_t)^2]^{1/2}
        & p(t)^{1/2} u_t & 0 \\
    0& 0& [tp'(t)]^{1/2} & 0 \\
    0& 0 & 0 & w^{1/2} \\
    \end{pmatrix}
    \begin{pmatrix}
    dB^1(t)\\
    dB^2(t)\\
    dB^3(t)\\
    dB^4(t)
    \end{pmatrix}\,.
    \eeq
In the above, for instance, $b_t$ corresponds to a \textbf{``local average''} of the quantities
    \[
    \frac{(\bmeta^a,\by)}{N^{1/2}}\,,
    \]
corresponding to the first term in the decomposition \eqref{e:X.decomp}. This intuition will be formalized in the next few sections. In particular, to make sense of the idea of ``local average,'' in Definitions \ref{d:rerand} and \ref{d:rerand.Ising} we introduce a ``spatial rerandomization'' of the $\vX$ process. This allows us to define certain ``local averages'' 
$(\bbb,\vvv,\uuu,\rrr,\www)$
(see \eqref{e:def.b}, \eqref{e:def.v}, \eqref{e:def.u}, \eqref{e:def.r},
and \eqref{e:def.w})
which we will use to approximate the coefficients in
\eqref{e:intro.SDE}. The limiting martingale equations are given by Theorem~\ref{t:tightness} (see additionally Proposition~\ref{p:limiting.drift.qv}). 
\textbf{In Section~\ref{s:sde} we give a formal version of the limiting SDE \eqref{e:intro.SDE} --- see \eqref{e:actual.sde} and Proposition~\ref{p:sde.approx}.}
The spatial rerandomization goes through a few steps of processing, which we turn to next.

\subsection{Preprocessing and spatial rerandomization} 
\label{ss:preprocessing}

In this subsection, we first preprocess $\vX$ so that it will satisfy some technical requirements for the analysis of later sections. We then define the spatial rerandomization of the (preprocessed) $\vX$, which we will analyze in the following sections.

\begin{dfn}[truncation] \label{d:trunc}
For the spherical perceptron, let $\vX\equiv\vX^{\RomI:\RomIII}$ be the three-dimensional process defined by \eqref{e:X.decomp} and \eqref{e:vX} above. Given a (large constant) parameter $\trK$, we define a \textbf{truncated} process $\vX^\trstar$ by initializing
$\vX^\trstar(q_0,a)=\vX(q_0,a) = (0,0,0)\in\R^3$, and setting
    \[
    \vX^\trstar(q_{d+1},a)
    =\begin{cases}
    \vX(q_{d+1},a) &\textup{if $\vX^\trstar(q_d,a)\in[-\trK,\trK)^3$,}\\
    \vX^\trstar(q_d,a) &\textup{otherwise,}
    \end{cases}
    \]
so that the process $\vX^\trstar$ stays frozen upon exiting $[-\trK,\trK)^3$. For the Ising perceptron, we apply the above procedure to obtain $\vX^{\trstar,\RomI:\RomIII}$. We separately define $X^{\trstar,\Ising}$ to be the process $X^\Ising$ from \eqref{e:Delta.X.Ising}, frozen upon exiting $[-\trK,\trK)$. We then denote $\vX^\trstar\equiv(\vX^{\trstar,\RomI:\RomIII},X^{\trstar,\Ising})$.
\end{dfn}

\begin{dfn}[addition of frozen particles]\label{d:froze} We now define a further modification $\vX^\frstar$, which may be thought of as $\vX^\trstar$ together with an additional \textbf{positive density of frozen particles}. For the spherical perceptron, the process $\vX^\frstar\equiv\vX^{\frstar,\RomI:\RomIII}$ is defined on the augmented space
    \[[M^\frozen] 
    \equiv
    \bigg[ M\bigg(1 + \frac1{\trK^6}\bigg)\bigg]\,.
    \]
For $a\in[M]$ we set
$\vX^\frstar(q_d,a)\equiv\vX^\trstar(q_d,a)$.
The remaining particles
$\vX^\frstar(q_d,a)$ for $a\in[M^\frozen]\setminus[M]$ are spaced evenly over $[-2\trK,2\trK]^3$ and frozen for all time. For the Ising perceptron, we again apply the above procedure to obtain $\vX^{\frstar,\RomI:\RomIII}$. We separately define $X^{\frstar,\Ising}$ on the augmented space
    \[[N^\frozen] 
    \equiv
    \bigg[ N\bigg(1 + \frac1{\trK^3}\bigg)\bigg]\,,
    \]
such that the added particles are spaced evenly over $[-2\trK,2\trK]$ and frozen for all time. We then denote $\vX^\frstar\equiv(\vX^{\frstar,\RomI:\RomIII},X^{\frstar,\Ising})$.
\end{dfn}

\begin{dfn}[spatial blocks] \label{d:blocks}
Let $\vl\equiv (\ell_1,\ell_2,\ell_3)$, and 
partition $[-\trK,\trK)^3$ into \textbf{blocks}
    \beq\label{e:J.blocks}
    J_{\vl}
    \equiv 
    \prod_{i=1}^3
    \Big[ \ell_i\eta,
    (\ell_i+1)\eta\Big)\,.
    \eeq
(We assume without loss that $\trK/\eta$ is integer-valued.) For any $0\le d\le\dmax$, abbreviate
    \[J(q_{0:d})
    \equiv
    \Big(J(q_0),\ldots,J(q_d)\Big)
    \]
where each $J(q_\ell)$ is one of the blocks defined by \eqref{e:J.blocks}.
For the Ising perceptron, we additionally let
    \[K(q_{0:d})
    \equiv
    \Big(K(q_0),\ldots,K(q_d)\Big)
    \]
where the \textbf{blocks} $K_\ell \equiv [\ell\eta,(\ell+1)\eta)$ partition the set $[-\trK,\trK)$.
\end{dfn}

\begin{dfn}[buckets for spherical perceptron]
\label{d:buckets} We now define $\vX^\trunc$ and $\vX^\frozen$, which may be regarded as the processes $\vX^\trstar$ and $\vX^\frstar$ with \textbf{small spatial buckets frozen}. For the spherical perceptron, initialize $\vX^\trunc(q_0)\equiv\vX^\trstar(q_0)$
and $\vX^\frozen(q_0)\equiv\vX^\frstar(q_0)$, and suppose inductively that $\vX^\trunc(q_\ell)$
and $\vX^\frozen(q_\ell)$
have been defined for $0\le\ell\le d$. Then let
    \begin{align}\label{e:bucket}
    B(q_d, J(q_{0:d})) 
    &\equiv \Big\{ a\in[M]:
    \vX^\trunc(q_\ell,a)
    \in J(q_\ell)
    \textup{ for each }0\le\ell\le d
    \Big\}\,,\\
    \label{e:bucket.fr}B^\frozen(q_d, J(q_{0:d})) 
    &\equiv \Big\{ a\in[M^\frozen]:
    \vX^\frozen(q_\ell,a)
    \in J(q_\ell)
    \textup{ for each }0\le\ell\le d
    \Big\}\,.\end{align}
We emphasize that in both \eqref{e:bucket} and \eqref{e:bucket.fr}, the $J$ blocks are restricted to within $[-\trK,\trK)^3$. We then let
    \begin{align*}
    B(q_d,a)
    &= \begin{cases}
    B(q_d, J(q_{0:d}))
    &\textup{if $a\in B(q_d, J(q_{0:d}))$,}\\
    \{a\}
    &\textup{otherwise, meaning 
    $\vX^\trunc(q_d, a)
    \notin[-\trK,\trK)^3$;}
    \end{cases} \\
    B^\frozen(q_d,a)
    &= \begin{cases}
    B^\frozen(q_d, J(q_{0:d}))
    &\textup{if }a\in B(q_d, J(q_{0:d}))\,,\\
    \{a\}
    &\textup{otherwise, meaning }\vX^\frozen(q_d, a)\notin[-\trK,\trK)^3\,.
    \end{cases}
    \end{align*}
The sets $\mathcal{B}(q_d)\equiv\{B(q_d,a) : a\in[M]\}$ partition $[M]$, and the partition is progressive in the sense that each $\mathcal{B}(q_{d+1})$ refines the previous partition $\mathcal{B}(q_d)$.  Similarly, the sets
$\mathcal{B}^\frozen(q_d)\equiv\{B^\frozen(q_d,a) : a\in[M^\frozen]\}$ progressively partition $[M^\frozen]$. In either case we call these the \textbf{buckets} at time $q_d$.  We define the quantity
    \beq\label{e:buckets.scB}
    \BUCKETS
    \equiv \bigg(\frac{8\trK}{\eta}\bigg)^{6\dmax}\,,\eeq
so that $\BUCKETS$ upper bounds the total number of buckets within $[-\trK,\trK)^3$ over all times. We then let
    \begin{align}
    \label{e:small-bucket-freezing-trunc}
    \vX^\trunc(q_{d+1},a)
    &=\begin{cases}
    \vX^\trstar(q_{d+1},a)
    &\textup{if } |B(q_d,a)| 
        \ge M/\BUCKETS^2\,,\\
    \vX^\trunc(q_d,a)
    &\textup{otherwise}\,.
    \end{cases}\\
    \label{e:small-bucket-freezing-frozen}
    \vX^\frozen(q_{d+1},a)
    &=\begin{cases}
    \vX^\frstar(q_{d+1},a)
    &\textup{if } |B^\frozen(q_d,a)  \cap[M]  | 
        \ge M/\BUCKETS^2\,,\\
    \vX^\frozen(q_d,a)
    &\textup{otherwise}\,.
    \end{cases}
    \end{align}
Note that because the partitions $\mathcal{B}(q_d)$ and $\mathcal{B}^\frozen(q_d)$ are progressive, particles that are frozen in step $d+1$ (i.e. such that the second case of \eqref{e:small-bucket-freezing-trunc} or \eqref{e:small-bucket-freezing-frozen} holds) stay frozen in all future steps. Note also that when determining whether buckets for the $\vX^\frozen$ process are too small we only look at particle indices in $[M]$, that is, we ignore the added particles from Definition~\ref{d:froze} when determining bucket size. With this definition, one can verify by induction that the set of buckets $B\in \mathcal{B}(q_d)$ can be obtained by taking $B=B^\frozen\cap[M]$ for each $B^\frozen\in\mathcal{B}^\frozen(q_d)$.
\end{dfn} 

\begin{dfn}[buckets for Ising perceptron]\label{d:buckets.Ising}
For the Ising perceptron, first apply Definition~\ref{d:buckets}  to obtain $\vX^{\trunc,\RomI:\RomIII}$ and $\vX^{\frstar,\RomI:\RomIII}$. We separately apply an analogous procedure in the Ising coordinate: let 
$X^{\trunc,\Ising}(q_0) \equiv X^{\trstar,\Ising}(q_0)$
and
$X^{\frozen,\Ising}(q_0) \equiv X^{\frstar,\Ising}(q_0)$,
and suppose inductively that 
$X^{\trunc,\Ising}(q_\ell)$
and $X^{\frozen,\Ising}(q_0)$
have been defined for all $0\le\ell\le d$. Then let
    \begin{align*}
    B^\Ising(q_d, K(q_{0:d})) 
    &\equiv \Big\{ i\in[N]:
    X^{\trunc,\Ising}(q_\ell,i)
    \in K(q_{\ell})
    \textup{ for each }0\le\ell\le d
    \Big\}\,,\\
    B^{\frozen,\Ising}(q_d, K(q_{0:d})) 
    &\equiv \Big\{ i\in[N^\frozen]:
    X^{\frozen,\Ising}(q_\ell,i)
    \in K(q_\ell)
    \textup{ for each }0\le\ell\le d
    \Big\}\,,\end{align*}
analogously to \eqref{e:bucket} and \eqref{e:bucket.fr}. Then define the \textbf{buckets}
    \begin{align*}
    B^\Ising(q_d,i)
    &= \begin{cases}
    B^\Ising(q_d, K(q_{0:d}))
    &\textup{if $i\in B^\Ising(q_d, K(q_{0:d}))$,}\\
    \{i\}
    &\textup{otherwise, meaning $X^{\trunc,\Ising}(i)\notin[-\trK,\trK)$;}
    \end{cases} \\
    B^{\frozen,\Ising}(q_d,i)
    &= \begin{cases}
    B^\frozen(q_d, K(q_{0:d}))
    &\textup{if }i\in B(q_d, K(q_{0:d}))\,,\\
    \{i\}
    &\textup{otherwise, meaning }X^{\frozen,\Ising}(i)\notin[-\trK,\trK)\,.
    \end{cases}
    \end{align*}
Write $\mathcal{B}^\Ising(q_d)\equiv\{B^\Ising(q_d,i):i\in[N]\}$
and $\mathcal{B}^{\frozen,\Ising}(q_d)
\equiv\{B^{\frozen,\Ising}(q_d,i):i\in[N^\frozen]\}$. Lastly, let
    \begin{align*}
    X^{\trunc,\Ising}(q_{d+1},i)
    &=\begin{cases}
    X^{\trstar,\Ising}(q_{d+1},i)
    &\textup{if $|B^\Ising(q_d,i)| 
        \ge N/\BUCKETS^2$,}\\
    X^{\trunc,\Ising}(q_d,i)
    &\textup{otherwise.}
    \end{cases}\\
    X^{\frozen,\Ising}(q_{d+1},i)
    &=\begin{cases}
    \vX^\frstar(q_{d+1},i)
    &\textup{if $|B^{\frozen,\Ising}(q_d,i)  \cap[N]  | 
        \ge N/\BUCKETS^2$,}\\
    X^{\frozen,\Ising}(q_d,i)
    &\textup{otherwise.}
    \end{cases}
    \end{align*}
Note that when determining whether buckets are too small we only look at particle indices in $[N]$, that is, we ignore the added particles from Definition~\ref{d:froze} when determining bucket size. Altogether we denote
$\vX^\trunc\equiv(\vX^{\trunc,\RomI:\RomIII},X^{\trunc,\Ising})$ and $\vX^\frozen\equiv(\vX^{\frozen,\RomI:\RomIII},X^{\frozen,\Ising})$. One can verify by induction that the set of buckets $B^\Ising\in \mathcal{B}^\Ising(q_d)$ can be obtained by taking $B^\Ising = B^{\frozen,\Ising}\cap[N]$ for each $B^{\frozen,\Ising}\in\mathcal{B}^{\frozen,\Ising}(q_d)$.
\end{dfn}

We now finally define the spatial rerandomizations. In subsequent sections we will take limits of these processes.

\begin{dfn}[spatial rerandomization for spherical perceptron]\label{d:rerand}
We use the buckets from Definition~\ref{d:buckets} to
 define a new process $\vY\equiv\vY^{\RomI:\RomIII}$, as follows. First let $a(q_0)$ be sampled uniformly at random from $[M]$. Then, conditional on $a(q_0), \ldots, a(q_d)$, we let $a(q_{d+1})$ be sampled from the measure $\bP$ conditional on the bucket at time $q_d$, that is,
    \[a(q_{d+1})
    \sim \bP
    \bigg( a = \cdot \,\bigg|\, 
    a\in B(q_d,a(q_d)) \bigg)\,.
    \]
We then let $\vY(q_d)$ be the average of $\vX^\trunc$ over all the particles in the bucket containing $a(q_d)$, 
    \beq\label{e:def.Y}
    \vY(q_d)
    =
    \sum_{a\in B(q_d,a(q_d))}
    \frac{\vX^\trunc(q_d,a)}{|B(q_d,a(q_d))|}
    \,.
    \eeq
Informally, $\vY$ is the process $\vX^\trunc$ with ``spatial rerandomization at scale $\eta$.'' It is easy to check (by induction) that each $a(q_d)$ is marginally uniform on $[M]$, which implies that $\vX^\trunc(q_d)$ and $\vY(q_d)$ are within $W_{\infty}$ distance $O(\eta/\delta)$ for each fixed $d$.  Likewise, we define $\vY^\frozen\equiv\vY^{\frozen,\RomI:\RomIII}$ to be the spatial rerandomization of $\vX^\frozen$ at scale $\eta$, with random indices $a^\frozen(q_d) \in[M^\frozen]$
and buckets $\mathcal{B}^\frozen(q_d)$.
\end{dfn}

\begin{dfn}[spatial randomization for Ising perceptron]
\label{d:rerand.Ising} 
For the Ising perceptron, we first apply the above Definition~\ref{d:rerand} to obtain $\vY^{\RomI:\RomIII}$ and $Y^{\frozen,\RomI:\RomIII}$.  We separately use the buckets $B^\Ising(q_d,i)$ introduced in Definition~\ref{d:buckets.Ising} to construct $Y^\Ising$ as the spatial rerandomization of $X^{\trunc,\Ising}$, with random indices $i(q_d)\in[N]$. Likewise we use the buckets $B^{\frozen,\Ising}(q_d,i)$ from 
Definition~\ref{d:buckets.Ising} to construct $Y^{\frozen,\Ising}$
as the spatial rerandomization of $X^{\frozen,\Ising}$, with random indices $i^\frozen(q_d)\in[N^\frozen]$.  We denote
$\vY\equiv(\vY^{\RomI:\RomIII},Y^\Ising)$ where the $a(q_d)$ indices evolve independently from the $i(q_d)$ indices. Likewise we denote $\vY^\frozen\equiv(\vY^{\frozen,\RomI:\RomIII},Y^{\frozen,\Ising})$ where the $a^\frozen(q_d)$ indices evolve independently from the $i^\frozen(q_d)$ indices. 
\end{dfn}

Recall the $\sigma$-algebra $\cG(q_d)$ defined by \eqref{e:gaus.filt}. Note that
for all $0\le\ell\le d$, the partitions $\mathcal{B}(q_\ell)$,
$\mathcal{B}^\Ising(q_\ell)$,
$\mathcal{B}^\frozen(q_\ell)$, and
$\mathcal{B}^{\frozen,\Ising}(q_\ell)$
are all measurable with respect to $\cG(q_d)$. 
Recalling Remark~\ref{r:annealed.quenched}, from now on we also write $\bP\equiv\bP_{\bG}$ for the law of $(\vY,Y^\Ising)$ conditional on the gaussian disorder. We then see from \eqref{e:def.Y} that the processes $\vY$ and $\vY^\frozen$ 
are adapted with respect to the filtration $\cG=(\cG(q_d):0\le d\le \dmax)$.  In Appendix~\ref{s:kolmogorov}, we formally define and analyze (discrete-time) ``drift'' and ``quadratic variation'' processes based on this filtration; these will be used to help characterize the limit of the processes \eqref{e:process.X}. 

\begin{ass}[order of parameters]\label{a:params} 
Throughout this section and the analysis of Section~\ref{s:sde}, we will assume that parameters are sent to limits in an appropriate order. In particular, we will require
$\eta\ll\delta$ and $1/N^{0.01}\ll\delta\ll1$,  and 
$\BUCKETS\ll N^{0.01}\delta$, where
$\BUCKETS$ is defined by \eqref{e:buckets.scB}.  Eventually we will assume \beq
\label{eq:order-of-limits}
\frac1N \ll 
\eta\ll 
\delta\ll \tep \ll \lambda \ll
\frac{1}{\trK}
\ll \sep \ll 
\epsilon^\circ
\ll \frac{1}{L}
\ll
1,
\eeq
where each quantity is sufficiently small depending on all larger ones (the $\epsilon$ parameters will be introduced later). More precisely we will first take limits in the parameters $N,\eta,\delta$, obtain an SDE description, and then argue that in the limit of the remaining parameters this SDE gives a better and better approximation for $\mu_{\bG}(\cA)$.
 \end{ass}

\subsection{Tightness in continuum limit}
\label{ss:tightness.statements}

In this subsection we present
Theorem~\ref{t:tightness}, which states that the spatially rerandomized processes converge along subsequential limits, to continuous semimartingales. \textbf{Throughout the remainder of this section, ``$N\to\infty$'' means $N\to\infty$, $\eta\to0$, and $\delta\to0$ as specified in Assumption~\ref{a:params}.}

\begin{thm}\label{t:tightness}
Let $\cA_N=\cA_N(\bG,\bg^\aux)$ be a sequence of $L$-Lipschitz algorithms, as in Definition~\ref{d:Lip}. Let $\chi^N\equiv \chi_{\cA_N}$, and assume further that
    \[
    \frac{1}{L^2} \le (\chi^N)'(p) \le L^2
    \]
for all $0\le p\le 1$. Consider a Brownian bridge terminating at $(\bG,\bg^\aux)$, and let $\vX^N\equiv (\vX^{\RomI:\RomIII},X^\Ising)^N$ be the resulting $\R^4$-valued process defined by \eqref{e:X.decomp}--\eqref{e:vX}. Let $\vX^{\trunc,N}$ be the truncation given by Definitions~\ref{d:trunc}--\ref{d:buckets.Ising}, and let $\vY^N$ be the spatial rerandomization of $\vX^{\trunc,N}$ from Definitions~\ref{d:rerand} and \ref{d:rerand.Ising}. Let $\vD^N$ and $\vQ^N$ denote respectively its discrete-time drift and covariation with respect to the filtration given by the buckets --- see Definitions~\ref{d:D} and \ref{d:Q} below ($\vD^N$ is $\R^4$-valued, while $\vQ^N$ is $\R^{10}$-valued).

We make
$(\vY,\vD,\vQ)^N$ into continuous-time $\R^{18}$-valued processes by taking the linear interpolation on each time interval $[q_d,q_{d+1}]$. We also denote the piecewise constant variants $(\vY,\vD,\vQ)^{\bullet,N}$, so for example $\vY^{\bullet,N}(t)=\vY^N(q_d)$ for all $[q_d,q_{d+1})$. Let $\mu_{\bG^N}$ be the law of $(\vY,\vD,\vQ)^N$ conditional on the disorder $\bG^N$, and let $\Q_N$ be the law of $\mu_{\bG^N}$. Let $\mu_{\bullet,\bG^N}$ be the law of $(\vY,\vD,\vQ)^{\bullet,N}$ conditional on the disorder $\bG^N$, and let $\Q_{\bullet,N}$ be the law of $\mu_{\bullet,\bG^N}$. Note that 
$\mu_{\bG^N},\mu_{\bullet,\bG^N}$ are probability measures on the space of $\R^{18}$-valued paths equipped with the uniform topology, and we equip the space of such measures with the topology of weak convergence.

If $N\to\infty$ with parameters satisfying Assumption~\ref{a:params}, then the measures $\Q_N$ are tight, with respect to the weak convergence topology described above. If $\Q_N\to\Q$ along an integer subsequence $N\to\infty$, then we have $\Q_{\bullet,N}\Rightarrow\Q_N$ along the same subsequence. For any subsequential limit $\Q$, if $\mu\sim\Q$ and $(\vY,\vD,\vQ)\sim\mu$, then $\vY$ is a semimartingale with drift $\vD$ and covariation $\vQ$. Moreover, $Y^\RomI$ is a finite-variation process while $Y^{\RomII}$, $Y^{\RomIII}$, and $Y^{\Ising}$ are martingales: this means $D^\sigma\equiv0$ for $\sigma\in\{\RomII,\RomIII,\Ising\}$, and $Q^{\RomI,\sigma}\equiv0$ for all $\sigma$. We also have $Q^{\sigma,\Ising}\equiv0$ for $\sigma\in\{\RomI,\RomII,\RomIII\}$.

All the analogous statements hold with $\vY^{\frozen,N}$, $\vY^\frozen$, $\vD^\frozen$, $\vQ^\frozen$ in place of $\vY^N$, $\vY$, $\vD$, $\vQ$. 
\end{thm}

The \hyperlink{proof:t.tightness}{proof of Theorem~\ref{t:tightness}} is based on relatively straightforward Kolmogorov-type moment estimates, and we defer this to Appendix~\ref{s:kolmogorov}.

In our main result Theorem~\ref{thm:main}, the coefficients of the stochastic integral \eqref{eq:main-process-intro-general} will be obtained by taking suitable limits of discrete analogues, which are defined as local averages of the drift and quadratic variation processes. To write this in an explicit way, let $q_T$ be the first time $q$ that $\vY^{\RomI:\RomIII}(q)$ has either exited $[-\trK,\trK)^3$ (Definition~\ref{d:trunc}) or has landed in a small bucket (Definition~\ref{d:buckets}). 
Similarly, let $q_{\Ising,T}$ denote the first time $q$ that $Y^\Ising(q)$ has either exited $[-\trK,\trK)$ (Definition~\ref{d:trunc}) or has landed in a small bucket (Definition~\ref{d:buckets.Ising}). Recall also from
Definitions \ref{d:buckets} and \ref{d:buckets.Ising}
that particles that have exited $[-\trK,\trK)$ in any coordinate are frozen in singleton buckets. It follows that we can express
	\beq\label{e:indicator.qd.less.than.qT}
	\ind\{q_d<q_T\}
	= \ind\bigg\{
	\Big|B(q_d,a(q_d))\Big| \ge \frac{M}{\BUCKETS^2}
	\bigg\}\,,
	\eeq
that is to say, the indicator $\ind\{q_d<q_T\}$ can be determined as a function of the bucket $B(q_d,a(q_d))$ containing the random index $a(q_d)$.  Similarly, $\ind\{q_d<q_{\Ising,T}\}$ can be determined as a function of $B^\Ising\in\mathcal{B}^\Ising(q_d)$.  Recall the definitions of $\bar{\bg}^a,\bar{\bx},\by$ from 
\eqref{e:normalized.vectors.g}--\eqref{e:normalized.vectors.y}.
Recall also that
$\bmeta^a$ denotes the $a$-th row of $\bXi^{d+1}$ from \eqref{e:gaus.decomp}. We then define 
\begin{align}\label{e:def.b}
    \bbb^{\trunc,N}(q_d)
    &\equiv
    \frac{\Delta\bar{D}^{\RomI}(q_d)}{\delta_d p'(q_d)^{1/2}}
    = \ind\{ q_d < q_T\} 
    \frac{1}{|B|}
    \sum_{a\in B} 
    \frac{(\bmeta^a,\by)}{N^{1/2}}\,,\\
    \label{e:def.v}
    \vvv^{\trunc,N}(q_d)
    &\equiv\frac{\Delta\tilde{Q}^{\RomII}(q_d)}{\delta_d p(q_d)}
    = \ind\{q_d<q_T\} 
    \frac{1}{|B|}\sum_{a\in B}
    \frac{(\bar{\bg}^a,\by)^2}{N}\,,\\
    \uuu^{\trunc,N}(q_d)
    &\equiv
    \frac{\Delta\tilde{Q}^{\RomII,\RomIII}(q_d)}{\delta_d
    [p(q_d) q_d p'(q_d)]^{1/2}}
    = \ind\{q_d<q_T\} 
    \frac{1}{|B|}\sum_{a\in B}
    \frac{(\bmeta^a,\bar{\bx})
        (\bar{\bg}^a,\by)}{N} \,.
    \label{e:def.u}\\
    \label{e:def.r}
    \rrr^{\trunc,N}(q_d)
    &\equiv
    \frac{\Delta C^{\RomIII}(q_d)}{\delta_d}
    \equiv  \ind\{q_d<q_T\}  
    q_d p'(q_d) \,,\\
    \label{e:def.w}
    \www^{\trunc,N}(q_d)
    &\equiv \frac{\Delta\tilde{Q}^{\Ising}(q_d)}
        {\delta_d}
    = \ind\{q_d< q_{\Ising,T}  \} 
    \frac{1}{|B^{\Ising}|}
    \sum_{i\in B^{\Ising}}
    (\be_i,\by)^2\,.
    \end{align}
In \eqref{e:def.b}--\eqref{e:def.w}, for now the reader should take the right-most quantities as the definitions of the quantities $(\bbb,\vvv,\uuu,\rrr,\www)^{\trunc,N}$.
Each of these quantities is a function of the bucket at time $q_d$, and we sometimes emphasize this by denoting the quantities as
$\bbb^\trunc(q_d, B)$,
$\www^\trunc(q_d, B^\Ising)$, and so on. The intermediate quantities refer to processes $\bar{D}$, $\tilde{Q}$, $C^{\RomIII}$ that we have not yet discussed: for now, we will simply say that
\begin{itemize}
\item $\bar{D}^\RomI$ approximates the drift $D^\RomI$ of $Y^\RomI$,
\item $\tilde{Q}^\RomII$ approximates the quadratic variation $Q^\RomII$ of $Y^\RomII$, 
\item $\tilde{Q}^{\RomII,\RomIII}$ 
approximates the covariation $Q^{\RomII,\RomIII}$ between $Y^\RomII$ and $Y^\RomIII$,
\item $C^{\RomIII}$ approximates the quadratic variation $Q^\RomIII$ of $Y^\RomIII$, and
\item $\tilde{Q}^{\Ising}$ approximates the quadratic variation $Q^\Ising$ of $Y^\Ising$.
\end{itemize}
See Definitions~\ref{d:barD} and \ref{d:tQ}, and display \eqref{e:APPX.def.r}, for the precise definitions. We defer further discussion of these  processes to Appendix~\ref{s:kolmogorov}, since they will not be immediately relevant in what follows. 
We also define analogous quantities
$(\bbb,\vvv,\uuu,\rrr,\www)^{\frozen,N}$ corresponding to the $\vY^{\frozen,N}$ process: 
	\begin{align}
	\label{e:b.plus.hist}
	\bbb^{\frozen,N}(q_d,B^\frozen)
	&\equiv
	\frac{\Delta \bar{D}^{\frozen,\RomI,N}(q_d,B^\frozen)}
	{\delta_d (p^N)'(q_d)^{1/2}}
	= 
	\frac{|B|}{|B^\frozen|}
	\bbb^{\trunc,N}(q_d,B)
	\,, \\ \label{e:v.plus.hist}
	\vvv^{\frozen,N}(q_d,B^\frozen)
	&\equiv
	\frac{\Delta\tilde{Q}^{\frozen,\RomII,N}(q_d,B^\frozen)}{\delta_d p^N(q_d)}
	=
	\frac{|B|}{|B^\frozen|}
	\vvv^{\trunc,N}(q_d,B)
	\,, \\ \label{e:u.plus.hist}
	\uuu^{\frozen,N}(q_d,B^\frozen)
	&\equiv 
	\frac{\Delta\tilde{Q}^{\frozen,\RomII,\RomIII,N}(q_d,B^\frozen)}
		{ \delta_d
		[ p^N(q_d) \cdot q_d (p^N)'(q_d)]^{1/2}}
	=\frac{|B|}{|B^\frozen|}
	\uuu^{\trunc,N}(q_d,B)
		\,,\\
	\label{e:P.plus.hist}
	\rrr^{\frozen,N}(q_d,B^\frozen)
	&\equiv
	\frac{\Delta C^{\frozen,\RomIII}
		(q_d,B^\frozen)}{\delta_d}
	=\frac{|B|}{|B^\frozen|}
	\rrr^{\trunc,N}(q_d,B)\,,\\
	\label{e:w.plus.hist}
	\www^{\frozen,N}(q_d)
    &\equiv \frac{\Delta\tilde{Q}^{\frozen,\Ising}(q_d)}
        {\delta_d}
    =
    \frac{|B^{\Ising}|}{|B^{\frozen,\Ising}|}
    \www^{\trunc,N}(q_d)\,.
	\end{align}
We next define a ``coarsened'' version of the above coefficients, as follows:

\begin{dfn}[coarsened buckets]
\label{d:coarse.buckets}
Recall from Definition~\ref{d:buckets}
the progressive partition
$\mathcal{B}^\frozen(q_d)$. 
For each time $q_d$, and for each of the spatial blocks $J_{\vl}$ defined by \eqref{e:J.blocks}, we now define a single bucket
	\beq\label{e:coarsened.bucket}
	B^{\frozen,\coarse}
	\equiv B^{\frozen,\coarse}(q_d,J_{\vl})
	\equiv \bigg\{
	a\in [M^+]
	: \vX^\frozen(q_d,a) \in J_{\vl}
	\bigg\}\,.\eeq
Each such $B^{\frozen,\coarse}$ is a union of buckets $B^\frozen\in\mathcal{B}^\frozen(q_d)$.
For any $B^\frozen\in\mathcal{B}^\frozen(q_d)$ not covered by one of the above (i.e., singleton buckets corresponding to particles that have left $[-\trK,\trK)^3$), we define simply
$B^{\frozen,\coarse} = B^\frozen$.
We define $\mathcal{C}^\frozen(q_d)$
to be the coarsening of $\mathcal{B}^\frozen(q_d)$ given by all these buckets $B^{\frozen,\coarse}$. We analogously define
the coarsening $\mathcal{C}^{\frozen,\Ising}(q_d)$ of $\mathcal{B}^{\frozen,\Ising}(q_d)$.
\end{dfn}

We then use the coarsened buckets to define averaged coefficients
	\begin{align}
	\label{e:b.plus.coarse}
	\bbb^{\frozen,\coarse,N}(q_d,B^{\frozen,\coarse})
	&\equiv
	\sum_{B^\frozen \subseteq B^{\frozen,\coarse}}
	\frac{|B^\frozen|}
		{|B^{\frozen,\coarse}|}
	\bbb^{\frozen,N}(q_d,B^\frozen)
	\,, \\ \label{e:v.plus.coarse}
	\vvv^{\frozen,\coarse,N}(q_d,B^{\frozen,\coarse})
	&\equiv
\sum_{B^\frozen \subseteq B^{\frozen,\coarse}}
	\frac{|B^\frozen|}
		{|B^{\frozen,\coarse}|}
	\vvv^{\frozen,N}(q_d,B^\frozen)
	\,, \\ \label{e:u.plus.coarse}
	\uuu^{\frozen,\coarse,N}(q_d,B^{\frozen,\coarse})
	&\equiv
	\sum_{B^\frozen \subseteq B^{\frozen,\coarse}}
	\frac{|B^\frozen|}
		{|B^{\frozen,\coarse}|}
	\uuu^{\frozen,N}(q_d,B^\frozen)
	\,,\\
	\label{e:P.plus.coarse}
	\rrr^{\frozen,N}(q_d,B^{\frozen,\coarse})
	&\equiv\sum_{B^\frozen \subseteq B^{\frozen,\coarse}}
	\frac{|B^\frozen|}
		{|B^{\frozen,\coarse}|}
	\rrr^{\frozen,N}(q_d,B^\frozen)\,,\\
	\label{e:w.plus.coarse}
	\www^{\frozen,N}(q_d,B^{\frozen,\Ising})
	&\equiv\sum_{B^{\frozen,\Ising}
		\subseteq
		B^{\frozen,\coarse,\Ising}}
	\frac{|B^{\frozen,\Ising}|}
		{|B^{\frozen,\coarse,\Ising}|}
	\www^{\frozen,N}
	(q_d,B^{\frozen,\Ising})
	\end{align}
\textbf{In the stochastic integral \eqref{eq:main-process-intro-general},
the drift coefficient $b_s$ is obtained via a certain limit of the averaged coefficients $\bbb^{\frozen,\coarse,N}$.
The diffusivity coefficient $\sigma_s$ is obtained via a certain limit of the  averaged coefficients $\vvv^{\frozen,\coarse,N}$ and $\uuu^{\frozen,\coarse,N}$, together with a simple transform which will be presented in \S\ref{ss:sde.reparam.sigma}. Similarly in \eqref{eq:X-ising-intro}, the diffusivity coefficient $w_s$ is a certain limit of the averaged coefficients $\www^{\frozen,N}$.} The limiting procedure occupies most of Section~\ref{s:sde}.

\subsection{A priori estimates}\label{ss:apriori}

In this subsection we collect some \textit{a~priori} estimates that will be used in Section~\ref{s:sde}, although the proofs will be deferred to Appendix~\ref{s:kolmogorov}.

\begin{ppn}[Kolmogorov estimates for $\vY$ and $\vY^\frozen$]
\label{p:Y.kolmogorov}
As in Definitions~\ref{d:rerand} and \ref{d:rerand.Ising}, let $\vY$ and $\vY^\frozen$ be the spatial rerandomizations of $\vX^\trunc$ and $\vX^\frozen$. We then have the quenched moment bound
    \beq\label{e:Y.I.Kolmogorov}
    \frac{\bE[
    (Y^{\RomI}(t)
    -Y^{\RomI}(s))^2]}{(t-s)^2}
    \le L^{O(1)} (\MAX_N)^2\,,\eeq
where $\MAX_N$ is defined by \eqref{e:MAX.N}, and satisfies
$\E[(\MAX_N)^2] \le L^{O(1)}$ as well as the bound
\eqref{e.MAX.bound}. For $N$ large enough, we also have the annealed moment bounds
    \beq\label{e:Y.II.III.Is.Kolmogorov}
    \sum_{\sigma\in\{\RomII,\RomIII,\Ising\}}
    \frac{\E[(Y^{\sigma}(t)-Y^{\sigma}(s))^8]}{(t-s)^4}
    \le
    L^{O(1)}\,,
    \eeq
for all $q_0\le s < t\le 1$. The same bounds hold with $\vY^\frozen$ in place of $\vY$.
\end{ppn}

The \hyperlink{proof:p.Y.kolmogorov}{proof of Proposition~\ref{p:Y.kolmogorov}} appears at the end of \S\ref{ss:kolmogorov.Y}.

\begin{lem}[quenched second moment bounds for 
$\bar{D}^{\RomI}$ and $\bar{D}^{\frozen,\RomI}$]
\label{l:barD.I.second.quenched}
Let $\bar{D}^{\RomI}$ and $\bar{D}^{\frozen,\RomI}$ be defined by \eqref{e:barD} and \eqref{e:barD.plus}; and let $\bbb^\trunc$ and $\bbb^\frozen$ be defined by \eqref{e:def.b} and \eqref{e:def.b.plus}. We then have
    \[
    \bE[\bbb^\trunc(q_d)^2]
    =\frac{\bE[\Delta \bar{D}^{\RomI}(q_d)^2]}
        {(\delta_d)^2  p'(q_d)}
    \le
    (\MAX_N)^2\,,\]
with $\MAX_N$ as in \eqref{e:MAX.N}. The same bound holds with $(\bar{D}^{\frozen,\RomI},\bbb^\frozen)$ in place of
$(\bar{D}^{\RomI},\bbb^\trunc)$. 
\end{lem}

The \hyperlink{proof:l.barD.I.second.quenched}{proof of Lemma~\ref{l:barD.I.second.quenched}} appears in \S\ref{ss:kolmogorov.drift}.

\begin{lem}[quenched first moment bounds on $\tilde{Q}$ and $\tilde{Q}^\frozen$]\label{l:tQ.first.quenched}
With $\tilde{Q}$ the approximation of $Q$ given by Definition~\ref{d:tQ}, we have the quenched bounds
    \begin{align*}
    &\max\bigg\{
    \frac{\bE\Delta\tilde{Q}^{\RomI}(q_d)}{(\delta_d)^2 p'(q_d)},
    \frac{\bE\Delta\tilde{Q}^{\RomII}(q_d)}{\delta_d p(q_d)}
    =\bE \vvv^\trunc(q_d), \\
    &\qquad\qquad\qquad
    \frac{\bE|\Delta\tilde{Q}^{\RomII,\RomIII}|}
    {\delta_d [p(q_d) q_d p'(q_d)]^{1/2}}
    =\bE |\uuu^\trunc(q_d)|,
    \frac{\bE\Delta\tilde{Q}^{\Ising}(q_d)}{\delta_d}
    =\bE \www^\trunc(q_d)
    \bigg\}
    \le (\MAX_N)^2\,,
    \end{align*}
with $\MAX_N$ as in \eqref{e:MAX.N}. The same bounds hold for $\tilde{Q}^\frozen$, $\vvv^\frozen$, $\uuu^\frozen$, $\www^\frozen$.
\end{lem}

The \hyperlink{proof:l.tQ.first.quenched}{proof of Lemma~\ref{l:tQ.first.quenched}} appears in \S\ref{ss:kolmogorov.qv}. 
We next define some additional quantities which will also appear repeatedly in our quenched bounds. In Proposition~\ref{p:drift.quenched} and \eqref{e:SMAX.N}, we will define a random variable $\SMAX_N$
satisfying the bound
        \beq\label{e.SMAX.bound}
    \P\Big(\SMAX_N\ge \log N\Big)
    \le\exp\bigg(-\frac{(\log N)^2}{L^{O(1)}}\bigg)\,.
    \eeq
We then define
    \beq\label{e:def.TMAX}
    \TMAX_N
    \equiv
    \bigg\{
    (\MAX_N)^2+ \frac{(\SMAX_N)^2\BUCKETS}{N\delta}
    \bigg\}^{1/2}.
\eeq
Combining \eqref{e.MAX.bound} and \eqref{e.SMAX.bound}, it follows that 
\beq\label{e.JMAX.bound}
    \P[\TMAX_N\geq L^{O(1)}]\leq \exp\bigg(-\frac{(\log N)^2}{L^{O(1)}}\bigg).
\eeq
Note that in the limit $N\to\infty$, $\TMAX_N$ will behave like
$\MAX_N$. 

\begin{ppn}[quenched second moment bounds for $\vY$ and $\vY^\frozen$]\label{p:Y.notI.second.quenched}
For all $0\le s<t\le1$ we have
    \[
    \sum_{\sigma\in\{\RomII,\RomIII,\Ising\}}
    \frac{\bE[(Y^\sigma(t)-Y^\sigma(s))^2]}{(t-s)}
    \le L^{O(1)} \bigg\{
    (\MAX_N)^2+ \frac{(\SMAX_N)^2\BUCKETS}{N\delta}
    \bigg\}
    = L^{O(1)} (\TMAX_N)^2
    \,,
    \]
with $\TMAX_N$ defined by \eqref{e:def.TMAX}. The same bound holds with $\vY^\frozen$ in place of $\vY$.
\end{ppn}
The \hyperlink{proof:p.Y.notI.second.quenched}{proof of Proposition~\ref{p:Y.notI.second.quenched}} appears in \S\ref{ss:kolmogorov.qv}.

\subsection{Budget constraints}
\label{ss:budget}

In this subsection we present the  \textbf{budget constraints} for our processes, corresponding to
\eqref{eq:is-budget-constraint} and \eqref{eq:sp-budget-constraint} from Definition~\ref{d:admissible.controls}. Recall the quantities $(\bbb,\vvv,\uuu,\rrr,\www)^{\trunc,N}$
were previously defined in 
\eqref{e:def.b}--\eqref{e:def.w}, relating to the drift and covariation of the $\vY$ processes. 
The vectors
$\bar{\bg}^a$,
$\bar{\bx}$, and $\by$ are defined by 
\eqref{e:normalized.vectors.g}--\eqref{e:normalized.vectors.y}; and $\bmeta^a$ refers to  the $a$-th row of $\bXi^{d+1}$ from \eqref{e:gaus.decomp}.
Next recall from
\eqref{e:b.plus.hist}--\eqref{e:w.plus.hist} that $(\bbb,\vvv,\uuu,\rrr,\www)^{\frozen,N}$ are defined analogously, but for $\vY^\frozen$ in place of $\vY$. Finally recall the coarsened coefficients from \eqref{e:b.plus.coarse}--\eqref{e:w.plus.coarse}.

Clearly, all these quantities depend on $N$ and on other parameters from Assumption~\ref{a:params}, although we often suppress this from the notation. Below, we will say an event has \textbf{high probability} if it has probability tending to $1$ in the limiting regime described by  Assumption~\ref{a:params}. In this section we further introduce a new positive parameter $\epsilon^\circ\ll 1$, which is much larger than all the other parameters $\delta$, $\eta$, $1/\trK$ that have appeared so far. Recall that $\bE$ denotes expectation conditional on the gaussian disorder.

For $\epsilon\ge0$ let $\Cost_\epsilon : \R \times [0,\infty) \to\R$ be defined by
    \beq\label{e:cost}
    \Cost_\epsilon(b,v,u)
    \equiv b^2 + u^2 + 
        \bigg( \Big[(v-u^2 + \epsilon)_+\Big]^{1/2}-1\bigg)^2\,.
    \eeq
One can easily check that $\Cost_\epsilon$ is convex when restricted to $\{(b,v,u) : v\ge u^2-\epsilon\}$. We abbreviate $\Cost\equiv\Cost_0$. The next theorem presents the budget constraint in a slightly simplified setting. We then show that this implies constraints for the coefficients 
$(\bbb,\vvv,\uuu,\www)^{\trunc,N}$,
 $(\bbb,\vvv,\uuu,\www)^{\frozen,N}$, and
 $(\bbb,\vvv,\uuu,\rrr,\www)^{\frozen,N}$ in Propositions \ref{p:budget-constraints-without-frozen}--\ref{p:budget-constraints-coarsened}  below. 

\begin{thm}[budget constraint]
\label{t:free.prob}
Fix parameters $\alpha,L,\iota,\epsilon$. Let $M,N\to\infty$ with $M/N\to\alpha$. Let $\omega$ denote a random variable. Let $\bG,\bXi$ be $M\times N$ matrices with i.i.d.\ standard gaussian entries, with rows $\bg^1,\ldots,\bg^M$ and $\bmeta^1,\ldots,\bmeta^M$ respectively,  such that $\bXi$ is independent of $(\bG,\omega)$. Let $\bar{\bx}\equiv \bar{\bx}(\bG)$ and $\hat{\by} \equiv \hat{\by}(\bG,\bXi)$ be $L$-Lipschitz, $\R^N$-valued functions such that 
	\begin{align*}
	&\E(\|\bar{\bx}\|^2)=\E(\|\hat{\by}\|^2)=N,\\
	&\E(\hat{\by}\,|\,\bG)=0 \textup{ almost surely.}
	\end{align*}
Suppose for some $r,s \ge 1$ we have a partition, which is permitted to depend on $(\bG,\omega)$ but not on $\bXi$:
	 \begin{align*}
    [M] = \bigsqcup_{j=1}^r B_j \textup{ with } 
    \frac{|B_j|}{M}
    = \lambda_j(\bG,\omega)\,,\\
    [N] = \bigsqcup_{j=1}^s 
    K_j \textup{ with } 
    \frac{|K_j|}{N}
    = \mu_j(\bG,\omega)\,,
    \end{align*}
such that $\min_j\lambda_j\ge\iota$ and $\min_j\mu_j\ge\iota$ hold $(\bG,\omega)$-almost surely. Then define the random variables
     \begin{align}
    \label{e:freeprob.b}
    \bbb_j \equiv \bbb(B_j)
        &\equiv 
        \frac{1}{|B_j|}
        \sum_{a\in B_j}
        \frac{(\bmeta^a,\by)}{N^{1/2}}
        \,,\\
    \label{e:freeprob.v}
    \vvv_j \equiv \vvv(B_j)
        &\equiv 
        \frac{1}{|B_j|}
        \sum_{a\in B_j}
        \frac{(\bar{\bg}^a,\by)^2}{N}
        \,,\\
    \label{e:freeprob.u}
    \uuu_j \equiv \uuu(B_j)
        &\equiv 
        \frac{1}{|B_j|}
        \sum_{a\in B_j}
        \frac{(\bmeta^a,\bar{\bx})
        (\bar{\bg}^a,\by)}{N}
        \,,\\
    \label{e:freeprob.w}
    \www_j
    \equiv \www(K_j) 
    &\equiv \frac{1}{|K_j|}
        \sum_{i\in K_j}
        (\by_i)^2\,.\end{align}
Recall the definition of the function $\Cost\equiv\Cost_0$ from \eqref{e:cost}. Then, for all sufficiently large $N$, we have
    \[
    \sum_{j=1}^r
    \frac{|B_j|}{M}
    \Cost(\bbb_j,\vvv_j,\uuu_j)
    \le
    \frac{1}{\alpha}
    \bigg( \sum_{j=1}^s
        \frac{|K_j|}{N}
        (\www_j)^{1/2}\bigg)^2
     +  \epsilon\,,
    \]
with probability at least $1-e^{-cN}$, where $c$ is a positive constant depending on the parameters $\alpha,L,\iota,\epsilon$.\end{thm}

The proof of Theorem~\ref{t:free.prob} is deferred to Appendix~\ref{a:freeprob} (where it is restated as Proposition~\ref{p:fp}), although a heuristic derivation of a simplified statement is presented in \S\ref{ss:heuristic.budget} below. Although the constraint set established in Theorem~\ref{t:free.prob} is tight (see Remark~\ref{r:free.prob.tight}), we do not actually prove or require this tightness. In Section~\ref{sec:IAMP}, we present incremental AMP algorithms which explicitly construct any SDE satisfying the limiting constraints derived from Theorem~\ref{t:free.prob}.

\begin{lem}\label{l:w.average.1}
Recall $\by(q_{d+1})$ defined in \eqref{e:normalized.vectors.y}. With probability at least $1-\exp(-N^{0.1})$, we have 
    \[
    \bigg|\frac{\|\by(q_{d+1})\|^2}{N}-1\bigg|
    \le \frac{1}{N^{1/4}}\]
for all $0\leq d\leq \dmax-1$.

\begin{proof}
Recalling \eqref{e:normalized.vectors.y}, we have
    \[
    \frac{\|\by(q_{d+1})\|^2}{N}
    =
    \frac{1}{\delta_d}
    \bigg(
    \frac{(\bx(q_{d+1}),\bx(q_{d+1}))}{N}
    -
    \frac{2(\bx(q_{d+1}),\bx(q_{d}))}{N}
    +
    \frac{(\bx(q_{d}),\bx(q_{d}))}{N}
    \bigg)\,.\]
Applying Lemma~\ref{l:nearly.ultrametric} to each term on the right-hand side gives 
    \[
    \frac{\|\by(q_{d+1})\|^2}{N}
    =\frac{1}{\delta_d}
    \bigg\{q_{d+1}-q_d +\frac{O(1)}{N^{1/3}} \bigg\}
    = 1 + O\bigg(\frac{1}{\delta_d N^{1/3}}\bigg)
    \]
for all $0\le d\le \dmax$, with probability at least $1-\exp(-N^{0.1})$.
The claim follows by sending $N\to\infty$ fast enough depending on $\delta$, as discussed in Assumption~\ref{a:params}.
\end{proof}
\end{lem}

The next proposition applies Theorem~\ref{t:free.prob} to deduce ``preliminary'' budget constraints for the coefficients $(\bbb,\vvv,\uuu,\rrr,\www)^{\trunc,N}$. In the subsequent results, Proposition~\ref{p:budget-constraints-HIST} and Corollary~\ref{p:budget-constraints-coarsened}, we deduce the budget constraints for the quantities $(\bbb,\vvv,\uuu,\rrr,\www)^{\frozen,N}$ and  $(\bbb,\vvv,\uuu,\rrr,\www)^{\frozen,\coarse,N}$, from which we will ultimately extract an SDE limit in Section~\ref{s:sde}. 

\begin{ppn}[constraints on $(\bbb,\vvv,\uuu,\www)^{\trunc,N}$]
\label{p:budget-constraints-without-frozen} For the coefficients defined by \eqref{e:def.b}--\eqref{e:def.w}, with high probability, the following hold simultaneously for all $0\leq d\leq d_{\max}-1$.
\begin{enumerate}[(a)]
    \item 
    \label{it:domain-constraint}
    The preliminary ``domain constraint'': for all $B\in\mathcal{B}(q_d)$,
    \beq\label{e:domain.constraint}
    \vvv^{\trunc,N}(q_d,B) 
    \ge \uuu^{\trunc,N}(q_d,B)^2 
    - \frac{\eta^2}{2} 
    \,.\eeq
 \item
 \label{it:spherical-budget-constraint}
The preliminary ``budget constraint'' for the spherical perceptron:
    \beq
    \label{e:budget.constraint}
    \bE
    \bigg[ \Cost\Big(
        \bbb^{\trunc,N}(q_d),\vvv^{\trunc,N}(q_d),
            \uuu^{\trunc,N}(q_d)
        \Big) \bigg]
    \le \frac{1}{\alpha} 
    + \frac{\epsilon^\circ}{2} \,.
    \eeq
 \item
 \label{it:Ising-budget-constraint}
The preliminary ``budget constraint'' for the Ising perceptron:
    \beq
    \label{e:budget.constraint.Ising}
    \bE
    \bigg[ \Cost\Big(
        \bbb^{\trunc,N}(q_d),\vvv^{\trunc,N}(q_d),
        \uuu^{\trunc,N}(q_d)
        \Big) \bigg]
    \le \frac{\bE[ \www^{\trunc,N}(q_d)^{1/2}]^2}{\alpha} 
    + \frac{\epsilon^\circ}{2}\,.
    \eeq
\item \label{it:Ising-w-constraint} The preliminary ``diffusivity'' for the Ising perceptron:
    \[\Big|\bE(\www^{\trunc,N}(q_d)) -1\Big|
    \le \frac{\epsilon^\circ}{2}\,.\]
\end{enumerate}
(In the above, $\bE$ denotes probability conditional on the gaussian disorder.)
\end{ppn}

\begin{proof}
To simplify the notation, in the proof we suppress the dependence on $N$ whenever possible.

\eqref{it:domain-constraint} Recall that we write $\uuu^\trunc(q_d,B)$ and $\vvv^\trunc(q_d,B)$ to emphasize that the quantities depend on the bucket $B$ at time $q_d$. Recall from \eqref{e:indicator.qd.less.than.qT} that $\ind\{q_d<q_T\}$ can be determined as a function of $B$. 
If $q_d\ge q_T$, then $\uuu^\trunc(q_d)=0=\vvv^\trunc(q_d)$ and the constraint holds trivially. Assume therefore that $q_d<q_T$: by Definition~\ref{d:buckets}, we must then have $|B|\ge M/\BUCKETS^2$. The Cauchy--Schwarz inequality gives
    \begin{align}\nonumber
    \uuu^\trunc(q_d,B)^2
    &\stackrel{\eqref{e:def.u}}{=}
    \bigg\{ \frac{1}{|B|}\sum_{a\in B}
    \frac{(\bmeta^a,\bar{\bx})
        (\bar{\bg}^a,\by)}{N}\bigg\}^2
    \le 
    \bigg\{\frac{1}{|B|}\sum_{a\in B}
    \frac{(\bmeta^a,\bar{\bx})^2}{N}\bigg\}
    \bigg\{\frac{1}{|B|}\sum_{a\in B}
    \frac{(\bar{\bg}^a,\by)^2}{N}\bigg\} \\
    &\stackrel{\eqref{e:def.v}}{=}
    \bigg\{ \frac{1}{|B|}\sum_{a\in B}
    \frac{(\bmeta^a,\bar{\bx})^2}{N}\bigg\}
    \vvv^\trunc(q_d,B)
    \,.\label{e:u.v.cs}
    \end{align}
Next, $\vvv^\trunc(q_d,B)$ can be crudely bounded as
    \[
    \vvv^\trunc(q_d,B)
    \le \frac{\|\bar{\bG}(q_d)\|^2\|\by\|^2}{N|B|}
    \le \frac{\|\by\|^2}{N} 
    \BUCKETS^2 \MAX_N\,,
    \]
having used the definition of $\MAX_N$ from \eqref{e:MAX.N} and the assumption  $|B|\ge M/\BUCKETS^2$. We emphasize that the above bound holds over all $B\in\mathcal{B}(q_d)$.
Recall that $\MAX_N$ is controlled by the bound \eqref{e.MAX.bound} from Proposition~\ref{p:Y.kolmogorov}. 
Next, 
recalling the definition of $\by\equiv\by(q_{d+1})$ from \eqref{e:normalized.vectors.y}, it follows using Lemma~\ref{l:w.average.1} that $\|\by\|^2 \le 2N$ with probability at least $1-\exp(-N^{0.1})$. Substituting into the above bound on $\vvv^\trunc$, and taking a union bound over $1\le d\le\dmax$, we conclude that
    \[
    \max_{1\le d\le \dmax}
    \max\Big\{
    \vvv^\trunc(q_d, B)
    : B\in\mathcal{B}(q_d)
    \Big\}
    \le
    O(1)\BUCKETS^2
    \stackrel{\eqref{e:buckets.scB}}{=} 
    O(1) \bigg(\frac{8\trK}{\eta}\bigg)^{12\dmax}\,,
    \]
with very high probability. 

Now, returning to \eqref{e:u.v.cs}, the first factor
in the last expression is equidistributed as
    \[
    \frac{\|\bar{\bx}\|^2}{N} \frac{\Gamma_B}{|B|}
    \]
where $\Gamma_B$ denotes a chi-square random variable with $|B|$ degrees of freedom. Recalling the definition of $\bar{\bx}$ from \eqref{e:normalized.vectors.x}, it follows using
Lemma~\ref{l:nearly.ultrametric} that
    \[
    \bigg|\frac{\|\bar{\bx}\|^2}{N}-1\bigg|
    \stackrel{\eqref{e:normalized.vectors.x}}{=} 
    \bigg|\frac{\|\bx(q_d)\|^2}{Nq_d}-1\bigg|
    \le \frac{O(1)}{N^{1/3}\delta}\,,
    \]
except with probability at most $\exp(-N^{0.01})$. Recalling again that we assumed $|B| \ge N/\BUCKETS^2$, it follows by a Chernoff bound on the chi-squared random variable that
    \[
    \P\bigg(\frac{\Gamma_B}{|B|}
    \ge 1 + \frac{\BUCKETS}{N^{1/3}}\bigg)\le \frac{1}{\exp(N^{1/3})}\,.
    \]
Substituting into \eqref{e:u.v.cs} and taking a union bound gives, with high probability,
    \begin{align*}
    &\max_{1\le d\le\dmax}
    \max\Big\{
    \uuu^\trunc(q_d,B)^2 -\vvv^\trunc(q_d,B)
    : B\in\mathcal{B}(q_d)\Big\}
    \\
    &\qquad
    \le 
    O(1) \bigg\{ \frac{1}{N^{1/3}\delta} + \frac{\BUCKETS}{N^{1/3}}
    \bigg\} \bigg(\frac{8\trK}{\eta}\bigg)^{12\dmax}
    \stackrel{\eqref{e:buckets.scB}}{\le} \frac{O(1)}{N^{1/3}}
    \bigg(\frac{8\trK}{\eta}\bigg)^{18\dmax}
    \ll \eta^2\,.
    \end{align*}
The last bound above holds  by sending $N\to\infty$ fast enough depending on the other parameters, as explained in  Assumption~\ref{a:params}.

\eqref{it:spherical-budget-constraint}
We omit the proof of \eqref{it:spherical-budget-constraint}, since it is very similar to but simpler than the proof of \eqref{it:Ising-budget-constraint} below. In particular, we proceed identically to the proof of \eqref{it:Ising-budget-constraint}, but ignore the blocks $\mathcal{B}^\Ising(q_d)$ and take the singleton partition $K_1 = [N']$.

\eqref{it:Ising-budget-constraint} Fix $\eaux=1/\trK>0$ and denote $M'\equiv M(1+\eaux)$ and $N'\equiv N(1+\eaux)$. We will apply Theorem~\ref{t:free.prob} to an $M'\times N'$ system, as follows.  Recalling \eqref{e:gaus.decomp}, we let $\omega'\equiv(\bXi^1,\ldots,\bXi^d)$. We let $\bG'$ be an $M'\times N'$ matrix where the top left $M\times N$ block is given by $\bar{\bG}\equiv \bar{\bG}(q_d)$, and the remaining entries are i.i.d.\ standard gaussians (independent of all else). Similarly, we let $\bXi'$ be an $M'\times N'$ matrix where the top left $M\times N$ block is given by $\bXi\equiv\bXi^{d+1}$, and the remaining entries are i.i.d.\ standard gaussians (independent of all else). Thus $\bXi'$ is independent of $(\bG',\omega')$.
 With $\bar{\bx}$ and $\by$ as in \eqref{e:normalized.vectors.x} and \eqref{e:normalized.vectors.y}, define 
    \begin{align*}
    \bar{\bx}'
    &\equiv \bar{\bx}'(\bG')
    \equiv
    \bigg(\frac{N'}{N}\bigg)^{1/2}
     (\bar{\bx},\bzero)\in\R^{N'}
     \,,\\
    \by'
    &\equiv
    \by'(\bG',\bXi')
    \equiv
    \bigg(\frac{N'}{N}\bigg)^{1/2}
    (\by,\bzero) \in\R^{N'}\,.
    \end{align*}
where $\bzero$ stands for the zero vector in $N\eaux$ dimensions. Recalling 
\eqref{e:p.corr.overlap},
\eqref{e:p.q} and the relation $q_d=\chi(p_d)$, we can verify that 
$\E(\by'\,|\, \bG')=0$, and
	\begin{align*}
	\frac{\E(\|\bar{\bx}'\|^2)}{N'}
	&= \frac{\E(\|\bar{\bx}\|^2)}{N}
	= \frac{\chi(p_d)}{q_d} = 1\,,\\
	\frac{\E(\|\by'\|^2)}{N'}
	&=\frac{\E(\|\by\|^2)}{N}
	= \frac{\chi(p_{d+1}) + \chi(p_d)
		-2\chi(p_d)}{\delta_d}
	= \frac{q_{d+1}-q_d}{\delta_d}
    = 1\,.
	\end{align*}
We will also verify in Lemma~\ref{l:bar.x.subgaus} that $\bar{\bx}$ is an $L^{O(1)}$-Lipschitz function of $\bG$, and we verify in Lemma~\ref{l:y.subgaus} that $\by$ is an $L^{O(1)}$-Lipschitz function of $\bG,\bXi$.

Recall the blocks $\mathcal{B}(q_d)$
and $\mathcal{B}^\Ising(q_d)$  from Definitions~\ref{d:buckets} and \ref{d:buckets.Ising}.
Then, in keeping with the notation of Theorem~\ref{t:free.prob}, we take the partitions
	\[[M'] = \bigsqcup_{j=1}^r
		B_j\,,\quad
    [N'] = \bigsqcup_{j=1}^s 
    	K_j\,,
    \]
where $B_2,\ldots,B_r$ are the blocks of $\mathcal{B}(q_d)$ for which $q_d<q_T$,
and $B_1 = B_\trunc \sqcup B_\textup{aux}$ where $B_\trunc$ is the union of blocks with $q_d\ge q_T$, and $B_\aux \equiv [M']\setminus [M]$. 
Similarly, $K_2,\ldots,K_s$ are the blocks of $\mathcal{B}^\Ising(q_d)$ for which $q_d < q_{\Ising,T}$, and $K_1 = K_\trunc \sqcup K_\textup{aux}$ where $K_\trunc$ is the union of the blocks with $q_d\ge q_{\Ising,T}$, and $K_\textup{aux}\equiv[N']\setminus[N]$.  Since $\eaux>0$, it follows that every bucket $B_j$ occupies a positive fraction of $[M']$, and every bucket $K_j$ occupies a positive fraction of $[N']$. It follows from Definition~\ref{d:buckets}  that these random partitions can be viewed as functions of $(\bG',\omega')$.
 
Let $\bbb',\vvv',\uuu',\www'$ 
denote the quantities
defined by
\eqref{e:freeprob.b}--\eqref{e:freeprob.w}, but for the expanded $M\times N'$ system. Writing $\bmeta^a$ for the ($N'$-dimensional) rows of $\bXi'$,
for $j\ge2$ we have
    \[
    \bbb'(B_j)
    = \frac{1}{|B_j|}
    \sum_{a\in B_j}
    \frac{(\bmeta^a,\by')}{(N')^{1/2}}
    = \bbb(B_j)\,.
    \]
Similarly we have $\vvv'(B_j)=\vvv(B_j)$ and $\uuu'(B_j)=\uuu(B_j)$ for all $j\ge2$. We also have $\www'(K_j) = (N'/N)\www(K_j)$ for $2\le j\le r$, and
    \[
    \www'(K_1)
    = \frac{1}{|K_1|}
    \sum_{i\in K_1} 
    \frac{N'}{N}(\by_i)^2
    =
    \frac{|K_\trunc|
    	N'/N}
    	{|K_\trunc|
	+ |K_\textup{aux}|}
    \frac{1}{|K_\trunc|}
    \sum_{i\in K_\trunc} 
    (\by_i)^2
    = \frac{p_\trunc(1+\eaux)}
        {p_\trunc+\eaux}
    \www(K_\trunc)\,,
    \]
where $p_\trunc\equiv |K_\trunc|/N$.

The conditions of
Theorem~\ref{t:free.prob} are satisfied for the $M'\times N'$ system with randomness $(\omega',\bG',\bXi')$. Applying the theorem gives, with very high probability,
    \begin{align}\nonumber
    &\sum_{j=1}^r
    \frac{|B_j|}{M}
    \Cost\Big(
        \bbb(B_j),
        \vvv(B_j),
        \uuu(B_j)\Big)
    \le
    \frac{1+\eaux}{\alpha}
    \bigg( 
    \sum_{j=1}^s
        \frac{|K_j|}{N'}
        \www'(K_j)^{1/2}
        \bigg)^2
     + \frac{\epsilon^\circ}{4}
      \\
    &\qquad
    = \frac1{\alpha}
    \bigg( 
    \Big[
    (p_\trunc+\eaux)
    p_\trunc 
    \www(K_\trunc)
    \Big]^{1/2}
    +
    \sum_{j=2}^s
        \frac{|K_j|}{N}
        \www(K_j)^{1/2}
        \bigg)^2
     + \frac{\epsilon^\circ}{4}
    \label{e:eaux.application.freeprob}
    \end{align}
for all times $0\le d\le\dmax$.
The quantity $p_\trunc=|K_\trunc|/N$ has two sources: (i) particles that exited the interval $[-\trK,\trK)$ in some coordinate (Definition~\ref{d:trunc}), and (ii) particles that never exited
the interval $[-\trK,\trK)$, but landed in a small bucket (Definition~\ref{d:buckets}). For (ii), recall from the discussion around \eqref{e:buckets.scB} that
$\BUCKETS$ upper bounds the total number of buckets within $[-\trK,\trK)^3$ over all times. It follows that the total number of particles frozen in small buckets is upper bounded by
    \beq\label{e:total.number.of.particles.in.small.buckets}
    \BUCKETS
    \cdot \frac{M}{\BUCKETS^2}
    = \frac{M}{\BUCKETS}
    \le M o_\eta(1)\,,
    \eeq
recalling Assumption~\ref{a:params}.
For (i), Corollary~\ref{c:trunc.small} gives that the total number of particles exiting $[-\trK,\trK)$ in any coordinate is at most $Mo_{\trK}(1)$, with very high probability. It follows that $p_\trunc \le o_{\trK}(1)$ with very high probability. Meanwhile,
Lemma~\ref{l:w.average.1} gives that with very high probability,
    \[
    2 \ge \frac{\|\by\|^2}{N}
    =\frac{|K_\trunc|}{N}
    \www(K_\trunc)
    +
    \sum_{j=2}^s \frac{|K_j|}{N}
    \www(K_j)
    \ge p_\trunc 
    \www(K_\trunc)\,.
    \]
Recall also that we set $\eaux=1/\trK$.
Altogether we conclude that with very high probability, we have $(p_\trunc+\eaux)
    p_\trunc \www(K_\trunc)\le o_{\trK}(1)$ for all times $0\le d\le \dmax$. Combining with \eqref{e:eaux.application.freeprob} gives
    \begin{align*}
    \bE\bigg[ \Cost\Big(\bbb^\trunc(q_d),
        \vvv^\trunc(q_d),\uuu^\trunc(q_d)\Big)\bigg]
    &\le
    \frac{1}{\alpha}
    \bigg(
    o_{\trK}(1)
     + \bE[ \www^\trunc(q_d)^{1/2}]
    \bigg)^2
    + \frac{\epsilon^\circ}{4} \\
    &\le
    \frac{\bE[ \www^\trunc(q_d)^{1/2}]^2}{\alpha} +\frac{\epsilon^\circ}{2} \,,
    \end{align*}
where the last bound again uses that
$\bE[ \www^\trunc(q_d)^{1/2}]^2 \le \|\by\|^2/N \le2$ with high probability, as noted above. (Note that the left-hand side above corresponds to the left-hand side of \eqref{e:eaux.application.freeprob} after dropping the $j=1$ term, which can only make the sum smaller.)

\eqref{it:Ising-w-constraint} 
Recalling the definition \eqref{e:def.w}, we can write
    \[
    \bE\www^{\trunc,N}(q_d)
    = \sum_{B^\Ising}
    \frac{|B^\Ising|}{N}
    \www^{\trunc,N}(B^\Ising,q_d)
    =\frac{\|\by\|^2}{N}
    - \sum_{B^\Ising}\frac1{N}
    \sum_{i\in B^\Ising} \ind\{q_d\ge q_{\Ising,T}(B^\Ising)\}
    (\by,\be_i)^2
    \,,
    \]
where $q_d\ge q_{\Ising,T}(B^\Ising)$ indicates that particles in $B^\Ising$ have exited  $[-\trK,\trK)$ (Definition~\ref{d:trunc}) or have entered a small bucket (Definition~\ref{d:buckets.Ising}). We again recall from Lemma~\ref{l:w.average.1} that $\|\by\|^2/N$, the first term on the right-hand side above, is concentrated around one. By the Cauchy--Schwarz inequality, the second term on the right-hand side is bounded in absolute value by
    \[
    \bigg\{\frac{1}{N} \sum_{i=1}^N(\by,\be_i)^4\bigg\}^{1/2}
    (p_{\trunc,\Ising})^{1/2}
    = \frac{\|\by\|_4^2}{N^{1/2}}
    \cdot (p_{\trunc,\Ising})^{1/2}
    \,,
    \]
where $p_{\trunc,\Ising}$ is the total fraction of particles that ever exit $[-\trK,\trK)$ or enter a small bucket. By a similar argument as for $p_\trunc$ above, we have $p_{\trunc,\Ising}= o_\eta(1)$. Since $\by$ is a subgaussian vector by Lemma~\ref{l:y.subgaus}, its $\ell_4$ norm satisfies $(\|\by\|_4)^4 \le O(N)$ with high probability, so the claim follows.
\end{proof}

We next show that the above constraints for $(\bbb,\vvv,\uuu,\www)^\trunc$ imply similar constraints for 
$(\bbb,\vvv,\uuu,\www)^\frozen$.

\begin{ppn}[constraints on 
$(\bbb,\vvv,\uuu,\www)^{\frozen,N}$]
\label{p:budget-constraints-HIST}
For the coefficients defined by \eqref{e:b.plus.hist}--\eqref{e:w.plus.hist}, 
with high probability,
the following hold simultaneously for all $0\leq d\leq d_{\max}-1$:
\begin{enumerate}[(a)]
\item 
    \label{it:domain-constraint-plus-HIST}
   The ``domain constraint'':
   for all $B^\frozen\in\mathcal{B}^\frozen(q_d)$,
    \beq\label{e:domain.plus.HIST}
    \vvv^{\frozen,N}(q_d,
    	B^\frozen) 
    \ge \frac{|B^\frozen|}{|B|}
     \uuu^{\frozen,N}(q_d,
     	B^\frozen)^2 
    - \frac{\eta^2}{2}\,.
    \eeq
 \item
 \label{it:spherical-budget-constraint-plus-HIST}
The ``budget constraint'' for the spherical perceptron:
    \beq
    \label{e:budget.plus.constraint.HIST}
    \bE_{\bG^N}
    \bigg[ \Cost\Big(
        \bbb^{\frozen,N}(q_d),
        \vvv^{\frozen,N}(q_d),
            \uuu^{\frozen,N}(q_d)
        \Big) \bigg]
    \le \frac{1}{\alpha} 
    + \epsilon^\circ
 .
    \eeq
\item \label{it:Ising-budget-constraint-plus-HIST}
The ``budget constraint'' for the Ising perceptron:
    \beq
    \label{e:budget.constraint.Ising.plus.HIST}
    \bE_{\bG^N}\bigg[ \Cost\Big(
        \bbb^{\frozen,N}(q_d),
        \vvv^{\frozen,N}(q_d),
        \uuu^{\frozen,N}(q_d)
        \Big) \bigg]
    \le \frac{\bE_{\bG^N}[
    \www^{\frozen,N}(q_d)^{1/2}]^2}{\alpha} 
    + \epsilon^\circ\,.
    \eeq
\item \label{it:Ising-w-constraint-plus-HIST} The ``diffusivity'' for the Ising perceptron:
    \[\Big|\bE(\www^{\frozen,N}(q_d)) -1\Big|
    \le\epsilon^\circ\,.\]
\end{enumerate}

\begin{proof} Recall the stopping times $q_T$, $q_{\Ising,T}$ that appear in \eqref{e:def.b}--\eqref{e:def.w}. Define analogously the stopping times $q_{\frozen,T}$ and $q_{\frozen,\Ising,T}$ with respect to the processes $\vY^{\frozen}$. In the same way that the indicator $\ind\{q_d<q_T\}$ can be determined as a function of $B\in\mathcal{B}(q_d)$, the indicator $\ind\{q_d<q_{\frozen,T}\}$ can be determined as a function of $B^\frozen\in\mathcal{B}^\frozen(q_d)$. Moreover, recalling that $B\in\mathcal{B}(q_d)$ corresponds to $B^\frozen\cap[M]$ for $B^\frozen\in\mathcal{B}^\frozen(q_d)$, we note it follows from the definitions that
    \beq\label{e:q.T.versus.q.frozen.T}
    \ind\{q_d<q_T\}(B)
    = \ind\{q_d<q_{\frozen,T}\}(B^\frozen)
    \eeq
This will be convenient in the proof below.

\eqref{it:domain-constraint-plus-HIST} For any $B^\frozen\in\mathcal{B}^\frozen(q_d)$, let $B=B^\frozen\cap[M]$. It follows from the definitions (and using the correspondence \eqref{e:q.T.versus.q.frozen.T} between $q_T$ and $q_{\frozen,T}$) that
    \[
    \vvv^\frozen(q_d,B^\frozen)
    = \frac{|B|}{|B^\frozen|}
    \vvv^\trunc(q_d,B)
    \stackrel{\eqref{e:domain.constraint}}
        {\ge}
    \frac{|B|}{|B^\frozen|}
    \bigg\{
    \uuu^\trunc(q_d,B)^2-\frac{\eta^2}{2}\bigg\}
    \ge \frac{|B^\frozen|}{|B|}
        \uuu^\frozen(q_d,B^\frozen)^2
        -\frac{\eta^2}{2}\,,
    \]
where the bound \eqref{e:domain.constraint} holds with high probability over all buckets and all $0\le d\le \dmax-1$ by the result of Proposition~\ref{p:budget-constraints-without-frozen}\ref{it:domain-constraint}. This proves the claim.

\eqref{it:spherical-budget-constraint-plus-HIST} We omit the proof of \eqref{it:spherical-budget-constraint-plus-HIST}, since it is very similar to but simpler than the proof of \eqref{it:Ising-budget-constraint-plus-HIST} below.

\eqref{it:Ising-budget-constraint-plus-HIST} 
As above, for $B^\frozen\in\mathcal{B}^\frozen(q_d)$, write $B=B^\frozen\cap[M]$. Recalling the correspondence \eqref{e:q.T.versus.q.frozen.T} between $q_T$ and $q_{\frozen,T}$, we can express
    \[
    \bbb^\frozen(q_d,B^\frozen)
    =\frac{|B|}{|B^\frozen|}
    \bbb^\trunc(q_d,B)
    +\frac{|B^\frozen\setminus B|}{|B^\frozen|}0\,,
    \]
and similarly for $\vvv^\frozen(q_d,B^\frozen)$, $\uuu^\frozen(q_d,B^\frozen)$.
Note that $\Cost_\epsilon(0,0,0)\le1$ for any small nonnegative $\epsilon$. Recall also the $\Cost_\epsilon$ function is convex subject to the constraint $\{v\ge u^2-\epsilon\}$. 
By Proposition~\ref{p:budget-constraints-without-frozen}\ref{it:domain-constraint} and \eqref{e:domain.plus.HIST}, both $(\vvv^\trunc,\uuu^\trunc)$
and $(\vvv^\frozen,\uuu^\frozen)$ satisfy this constraint, with high probability, provided we take $\epsilon=\eta^2$. Applying Jensen's inequality then gives
    \begin{align*}
    &\bE \bigg[\Cost_{\eta^2}\Big(\bbb^\frozen(q_d),
        \vvv^\frozen(q_d),\uuu^\frozen(q_d)\Big)
        \bigg] 
    =
    \sum_{B^\frozen\in\mathcal{B}^\frozen(q_d)}
    \frac{|B^\frozen|}{M^\frozen}
    \Cost_{\eta^2}\Big(\bbb^\frozen(q_d,B^\frozen),
        \vvv^\frozen(q_d,B^\frozen),\uuu^\frozen(q_d,B^\frozen)\Big)\\
    &\qquad\le
    \sum_{B^\frozen\in\mathcal{B}^\frozen(q_d)}
    \frac{|B^\frozen|}{M^\frozen}
    \bigg\{\frac{|B|}{|B^\frozen|}
    \Cost_{\eta^2}\Big(\bbb^\trunc(q_d,B),
        \vvv^\trunc(q_d,B),\uuu^\trunc(q_d,B)\Big)
    +\frac{|B^\frozen\setminus B|}
        {|B^\frozen|}
        \Cost_{\eta^2}(0,0,0)\bigg\}\\
    &\qquad\le
    \frac{M}{M^\frozen}
    \bE\bigg[\Cost_{\eta^2}\Big(
        \bbb^\trunc(q_d),
        \vvv^\trunc(q_d),
        \uuu^\trunc(q_d)
        \Big)\bigg]
        +\frac{M^\frozen-M}{M^\frozen}\\
    &\qquad
    \stackrel{\eqref{e:budget.constraint.Ising}}{\le}   \frac{\bE[ \www^\trunc(q_d)^{1/2}]^2}{\alpha}+\frac{\epsilon^\circ}{2}
    + O(\eta) + \frac{O(1)}{\trK^6} \,,
    \end{align*}
having used the definition of $M^\frozen$ from Definition~\ref{d:froze}, along with the fact that $\Cost$ and $\Cost_{\eta^2}$ differ by at most $O(\eta)$. We then expand
    \begin{align*}
    &\bE
    \Big[ \www^\trunc(q_d)^{1/2}\Big]
    = \sum_{B^\Ising\in\mathcal{B}^\Ising(q_d)}
    \frac{|B^\Ising|}{N}
    \www^\trunc(q_d,B^\Ising)^{1/2}
    =\sum_{B^{\frozen,\Ising}\in\mathcal{B}^\frozen(q_d)}
    \frac{|B^{\frozen,\Ising}|}{N}
    \bigg(\frac{|B^\Ising|}{|B^{\frozen,\Ising}|}\bigg)^{1/2}
    \www^\frozen(q_d,B^{\frozen,\Ising})^{1/2}\\
    &\qquad\le 
    \sum_{B^{\frozen,\Ising}\in\mathcal{B}^\frozen(q_d)}
    \frac{|B^{\frozen,\Ising}|}{N}
    \www^\frozen(q_d,B^{\frozen,\Ising})^{1/2}
    = \frac{N^\frozen}{N} \bE \Big[
    \www^\frozen(q_d,B^{\frozen,\Ising})^{1/2}
    \Big]\,.
    \end{align*}
Recalling again that $\bE[ \www^\trunc(q_d)^{1/2}]^2 \le \|\by\|^2/N \le2$ with high probability, we have
\[
    \bE \Big[ \www^\trunc(q_d)^{1/2}\Big]
    \le \bE \Big[
        \www^\frozen(q_d,B^{\frozen,\Ising})^{1/2}
    \Big] + \frac{N^\frozen - N}{N^\frozen} \bE \Big[ \www^\trunc(q_d)^{1/2}\Big]
    \le \bE \Big[
        \www^\frozen(q_d,B^{\frozen,\Ising})^{1/2}
    \Big] + \frac{O(1)}{\trK^2}\,.
\]
Combining the above bounds gives the claim.

 \eqref{it:Ising-w-constraint-plus-HIST} It follows from \eqref{e:w.plus.hist}
 that
    \[\bE\www^{\frozen,N}(q_d)
    =\sum_{B^{\frozen,\Ising}}
    \frac{|B^{\frozen,\Ising}|}{N^\frozen}
    \www^{\frozen,N}(q_d,B^{\frozen,\Ising})
    = \frac{N}{N^\frozen}
    \bE\www^{\trunc,N}(q_d)\,,
    \]
so the claim follows by Proposition~\ref{p:budget-constraints-without-frozen}\ref{it:Ising-w-constraint} combined with the choice of $N^\frozen$ from Definition~\ref{d:froze}.
\end{proof}
\end{ppn}

\begin{ppn}[constraints on $(\bbb,\vvv,\uuu,\rrr,\www)^{\frozen,\coarse,N}$]
\label{p:budget-constraints-coarsened} For the coefficients defined by \eqref{e:b.plus.coarse}--\eqref{e:w.plus.coarse},
with high probability,
the following hold simultaneously for all $0\leq d\leq d_{\max}-1$:
\begin{enumerate}[(a)]
\item 
    \label{it:domain-constraint-plus-coarse} 
    The ``domain constraint'': for all 
    $B^{\frozen,\coarse} \in \mathcal{C}^\frozen(q_d)$,
    \beq\label{e:domain.plus.coarse}
    \vvv^{\frozen,\coarse,N}(q_d,B^{\frozen,\coarse}) 
    \ge 
    q_d p'(q_d)
    \frac{\uuu^{\frozen,\coarse,N}(q_d,B^{\frozen,\coarse})^2 }
    {\rrr^{\frozen,\coarse,N}(q_d,B^{\frozen,\coarse})}
    - \frac{\eta^2}{2}\,.
    \eeq
(If $\uuu^{\frozen,\coarse,N}$ and $\rrr^{\frozen,\coarse,N}$ are both zero, we interpret the ratio on the right-hand side as zero.)
 \item
 \label{it:spherical-budget-constraint-plus-coarse}
The ``budget constraint'' for the spherical perceptron:
    \beq
    \label{e:budget.plus.constraint.coarse}
    \bE_{\bG^N}
    \bigg[ \Cost\Big(
        \bbb^{\frozen,\coarse,N}(q_d),
        \vvv^{\frozen,\coarse,N}(q_d),
            \uuu^{\frozen,\coarse,N}(q_d)
        \Big) \bigg]
    \le \frac{1}{\alpha} 
    + \epsilon^\circ
 .
    \eeq
\item \label{it:Ising-budget-constraint-plus-coarse}
The ``budget constraint'' for the Ising perceptron:
    \beq
    \label{e:budget.constraint.Ising.plus.coarse}
    \bE_{\bG^N}\bigg[ \Cost\Big(
        \bbb^{\frozen,\coarse,N}(q_d),
        \vvv^{\frozen,\coarse,N}(q_d),
        \uuu^{\frozen,\coarse,N}(q_d)
        \Big) \bigg]
    \le \frac{\bE_{\bG^N}[
    \www^{\frozen,\coarse,N}(q_d)^{1/2}]^2}{\alpha} 
    + \epsilon^\circ\,.
    \eeq
\item \label{it:Ising-w-constraint-plus-coarse} The ``diffusivity'' for the Ising perceptron:
    \[\Big|\bE(\www^{\frozen,\coarse,N}(q_d)) -1\Big|
    \le\epsilon^\circ\,.\]
\end{enumerate}

\begin{proof}
Assertions
\eqref{it:spherical-budget-constraint-plus-coarse}--\eqref{it:Ising-w-constraint-plus-coarse} follow straightforwardly from Jensen's inequality. For assertion \eqref{it:domain-constraint-plus-coarse}, note 
from \eqref{e:indicator.qd.less.than.qT},
 \eqref{e:def.r}, and \eqref{e:P.plus.hist}
that Proposition~\ref{p:budget-constraints-HIST}\ref{it:domain-constraint-plus-HIST} can be rewritten as
	\[
	\vvv^{\frozen,N}(q_d,B^\frozen) 
    \ge  q_d p'(q_d) 
    I(B^\frozen)
	\frac{\uuu^{\frozen,N}(q_d,B^\frozen)^2}
		{\rrr^{\frozen,N}(q_d,B^\frozen)}
    - \frac{\eta^2}{2}\,,
	\]
where we again interpret a $0/0$ ratio as zero, and
	\[I(B^\frozen)\equiv \ind\bigg\{|B^\frozen \cap[M]| \ge 
		\frac{M}{\BUCKETS^2}
		\bigg\}\,.\]
If $I(B^\frozen)=0$ for all $B^\frozen\in B^{\frozen,\coarse}$,
then the bound \eqref{e:domain.plus.coarse} follows.
Otherwise, consider a discrete probability measure $\pi$, and random variables $v_j,u_j,r_j$ which are zero on a set of indices $Z$, and otherwise satisfy
	\[
	v_j \ge A' \frac{(u_j)^2}{r_j} - a'\]
Assume $\pi_+ \equiv 1-\pi(Z)\in(0,1]$. Then Jensen's inequality gives
	\[
	\frac{\sum_j \pi_j v_j}{\pi_+}
	\ge \frac{1}{\pi_+}
	\sum_{j\notin Z} \pi_j \bigg\{
	A' \frac{(u_j)^2}{r_j} - a'\bigg\}
	\ge A'
	\bigg(\frac{1}{\pi_+}\sum_j \pi_j u_j\bigg)^2\bigg/
	\bigg(
	\frac{1}{\pi_+}\sum_j \pi_j r_j
	\bigg) - a'
	\]
and rearranging gives
	\[
	\sum_j \pi_j v_j
	\ge
	A' \bigg(\sum_j \pi_j u_j\bigg)^2\bigg/
	\bigg(\sum_j \pi_j r_j
	\bigg)
	-a'\pi_+\,.
	\]
If we apply the above inequality with the probability measure $\pi$
over buckets $B^\frozen\in B^{\frozen,\coarse}$ with probabilities 
$\pi(B^\frozen)=|B^\frozen|/|B^{\frozen,\coarse}|$, then we obtain the desired conclusion.
\end{proof}
\end{ppn}

\subsection{Heuristic derivation of budget constraint}\label{ss:heuristic.budget}

In this subsection we present a heuristic derivation of (a simplified version of) the budget constraint \eqref{eq:is-budget-constraint} and \eqref{eq:sp-budget-constraint}, using free probability calculations.
For an $N\times N$ hermitian matrix $A$, we let $\mu_A$ be its empirical eigenvalue distribution, and denote its Cauchy transform
    \[
    G_A(z)
    = \int \frac{\mu_A(dt)}{z-t}
    = \frac{\Tr[(z-A)^{-1}]}{N}\,.
    \]
We let $(G_A)^{-1}$ denote the functional inverse, and denote the $R$-transform
    \[
    R_A(s)
    = (G_A)^{-1}(s)-\frac1s\,.
    \]
Then free probability theory (see e.g.\ \cite{MR2760897} or \cite{MR2266879} for expository accounts) tells us that if $A$ and $B$ are ``free,'' then in the large-$N$ limit the eigenvalue distribution of $A+B$ can be characterized by the relation
    \beq\label{e:free.conv}
    R_{A+B}\approx R_A+R_B\,.
    \eeq
We will use this relation to heuristically derive the budget constraint.

\begin{xmp}[Wishart matrix]
To review a well-known example, 
consider the Wishart matrix
    \beq\label{e:wishart} W= \sum_{a=1}^M \frac{(\bg^a)^{\otimes2}}{N}
    = \sum_{a=1}^M W^a
    \eeq
with $M=N\alpha$. It is well known that the eigenvalue distribution of $W$ converges to the Marchenko--Pastur distribution, supported on $[\lambda_-,\lambda_+]$ with $\lambda_\pm=(1\pm\alpha^{1/2})^2$. We now review how to predict this from a free probability calculation, following \cite[Exercise 5.3.27]{MR2760897}. The rank-one matrix $W^a$ has one non-zero eigenvalue around $1$, so its Cauchy transform is approximately given by
    \[
    G_{W^a}(z) \approx \frac{1}{N}\bigg\{
    \frac{1}{z-1}+\frac{N-1}{z}
    \bigg\}\,.
    \]
Solving the above for $z$ gives
    \[
    (G_{W^a})^{-1}(s)
    \approx \frac{1+ s}{2s}
    \pm \frac{|1- s|}{2s}
    \bigg( 1+ \frac{4 s}
    {N(1- s)^2} \bigg)^{1/2}\,.
    \]
We can choose the positive root using the heuristic that $G_a(z)\sim 1/z$ for large $z$, leading to the $R$-transform
    \beq\label{e:Wa.Rtransform}
    R_{W^a}(s)
    = (G_{W^a})^{-1}(s)
    -\frac1s
    \approx 
    \frac{1}{N(1- s)}
    \eeq
for $s<1$. Applying the free convolution formula \eqref{e:free.conv} gives
    \beq\label{e:wishart.R}
    R_W(s)
    \stackrel{\eqref{e:free.conv}}{\approx} \sum_{a=1}^M 
    R_{W^a}(s)
    = \frac{\alpha}{1-s}
    = (G_W)^{-1}(s)- \frac1s\,.
    \eeq
Again using the heuristic that $G_W(z)\sim1/z$ for large $z$, one can check that the appropriate domain of the function $(G_W)^{-1}(s)$
is the interval with endpoints $G_W(\lambda_\pm)$, and the range is precisely $\R\setminus[\lambda_-,\lambda_+]$. This shows how to use \eqref{e:free.conv} to heuristically recover the support of the Marchenko--Pastur distribution.
\end{xmp}

\begin{xmp}[spiked Wishart matrix]
Let $P=\bu^{\otimes2}$ where $\bu$ is a unit vector. By a slight generalization of the calculation \eqref{e:Wa.Rtransform}, the $R$-transform of $\theta P$ is given by
    \[
    R_{\theta P}(s)
    =(G_{\theta P})^{-1}(s)-\frac1s
    \approx 
    \frac{\theta}{N(1-\theta s)}\,.
    \]
Let $W$ be the Wishart matrix \eqref{e:wishart}, whose $R$-transform is given by \eqref{e:wishart.R}.
For the spiked Wishart matrix $A=W+\theta P$, it is well known that for $0\le\theta\le 1/G_W(\lambda_+)$ then the spectrum of $A$ is similar to that of $W$; otherwise, if $\theta>1/G_W(\lambda_+)$, then there is one outlier eigenvalue around $(G_W)^{-1}(1/\theta)$ (see \cite{MR2165575}). To predict this from a free probability calculation, note that the free convolution formula \eqref{e:free.conv} gives
    \[
    (G_A)^{-1}(s)
    = \frac1s+R_A(s)
    \approx\frac1s+
    R_W(s)+R_{\theta P}(s)
    =\frac1s+\frac{\alpha}{1-s}
    + \frac{\theta}{N(1-\theta s)}\,.
    \]
If $\theta<1/ G_W(\lambda_+)$, then for large $N$ the function $(G_A)^{-1}$ looks very similar to the function $(G_W)^{-1}$ on the appropriate domain (the interval with endpoints $[G_W(\lambda_\pm)]$, as noted above). This is consistent with the eigenvalue distribution being asymptotically unchanged. If $\theta>1/G_W(\lambda_+)$, then $(G_A)^{-1}$ still looks similar to the function $(G_W)^{-1}$, but with a new singularity around $s=1/\theta$ corresponding to the new outlier eigenvalue around $(G_W)^{-1}(1/\theta)$. 
\end{xmp}

\begin{xmp}[spiked Wishart matrix continued]
For $W$ and $P$ as above, and any vector $\by\in\R^N$, let
    \begin{align*}
    v &= \frac{(\by,W\by)}{M}\,,\\
    x &= \frac{(\by,P\by)}{M}\,.
    \end{align*}
As $\by$ ranges over all vectors with $\|\by\|^2\le N$, we can ask what values of $(v,x)$ can be attained. If we consider $v$ or $x$ alone, then clearly we must have
    \begin{align*}
    v &\le \frac{\|\by\|^2}{M}
    \lambda_{\max}(W)
    = \frac{(1+\alpha^{1/2})^2}{\alpha}
    \,,\\
    x &\le \frac{\|\by\|^2}{M}
    \lambda_{\max}(P)
    = \frac{1}{\alpha}\,.
    \end{align*}
Moreover, it follows from the previous calculation that in general we must have
    \[
    v+\theta x
    = \frac{(\by,(W+\theta P)\by)}{M}
    \le 
    \begin{cases}
    \displaystyle
    \frac{\|\by\|^2}{M}(G_W)^{-1}\bigg(\frac{1}{\theta}\bigg)
    = \frac{\theta}{\alpha}
    + \frac{\theta}{\theta-1}
    \equiv h(\theta)
    &\textup{if $\theta>
    1/G_W(\lambda_+)$,}\\
    \displaystyle
    \frac{\lambda_+}{\alpha} &\textup{if $0\le\theta\le1/G_W(\lambda_+)$.}
    \end{cases}\,.
    \]
We claim that the full range of feasible $(v,x)$ corresponds to the set 
    \beq\label{e:v.x.constraint}
    \bigg\{(v,x):v\ge0, x\ge0,
    x + (v^{1/2}-1)^2
    \le \frac{1}{\alpha}\bigg\}\,.
    \eeq
The constraint in \eqref{e:v.x.constraint} can be predicted from free probability calculations as follows: let us assume a constraint of the form $x \le f(v)$, and solve for $f(v)$ such that
    \[
    \sup\bigg\{v+\theta x : x \le f(v)\bigg\} = h(\theta)\,.
    \]
(To recover the full domain \eqref{e:v.x.constraint}, one must also consider the problem of minimizing $v+\theta x$ for negative $\theta$, but for simplicity we will omit this calculation.) The Lagrangian is $L = v+\theta x+\zeta( f(v)-x)$, and stationarity equations give
    \[
    0 = \frac{\partial L}{\partial x}
    =\theta-\zeta\,,\quad
    0 = \frac{\partial L}{\partial v}
    = 1+\zeta f'(v)\,.
    \]
Thus 
at the stationary point we have $v=v(\theta)$ satisfying $f'(v)=-1/\theta$, 
and 
$x=x(\theta)=f(v(\theta))$. Thus
    \[
    v+\theta x
    = v(\theta)+\theta f(v(\theta))
    = h(\theta) = \frac{\theta}{\alpha}
    + \frac{\theta}{\theta-1}\,.
    \]
Differentiating in $\theta$ gives
$f(v(\theta))= h'(\theta)$, and differentiating again gives
    \[
    v'(\theta)
    = \frac{h''(\theta)}{f'(v(\theta))}
    = -\theta h''(\theta)
    =\frac{2\theta}{(1-\theta)^3}\,.
    \]
For $\theta$ large we expect the optimal $\by$ to be close to the principal eigenvector $\bu$ of $P$, and as a result $v(\theta)\approx1$. Integrating $v'(\theta)$ with this boundary condition gives
    \[v(\theta)
    = \frac{\theta^2}{(\theta-1)^2}\,,\quad
    \theta=\theta(v)
    = \frac{v^{1/2}}{v^{1/2}-1}\,.\]
Substituting this into the above gives
    \[
    f(v)
    = \frac{h(\theta)-v}{\theta}
    \bigg|_{\theta=\theta(v)}
    =\frac{1}{\alpha}-(v^{1/2}-1)^2\,,
    \]
thereby recovering the constraint from \eqref{e:v.x.constraint}.
\end{xmp}

\begin{xmp}[spherical perceptron with single bucket]
Let $\bar{\bg}^a$ and $\bmeta^a$ be i.i.d.\ standard gaussian vectors in $\R^N$. Consider the quantities
    \begin{align*}
    \bbb
    &=\frac1M
    \sum_{a=1}^M \frac{(\bmeta^a,\by)}{N^{1/2}}
    = \frac{(\bmeta,\by)}{(MN)^{1/2}}\,,\\
    \vvv
    &= \frac1M
    \sum_{a=1}^M \frac{(\bar{\bg}^a,\by)^2}{N}
    \,,\\
    \uuu
    &= \frac1M
    \sum_{a=1}^M 
    \frac{(\bmeta^a,\bar{\bx})}{N^{1/2}}
    \frac{(\bar{\bg}^a,\by)}{N^{1/2}}
    = \frac1M
    \sum_{a=1}^M z_a 
    \frac{(\bar{\bg}^a,\by)}{N^{1/2}}\,.
    \end{align*}
If we ignore $\uuu$, then the previous calculation shows that the range of feasible $(\bbb^2,\vvv)$ values is given by
    \[\bigg\{(\bbb^2,\vvv):\vvv\ge0,
    \bbb^2 +( \vvv^{1/2}-1)^2 \le \frac{1}{\alpha}
    \bigg\}\,.
    \]
Even if we take $\uuu$ into account, the $\bbb$ term enters in a simple way, so for expository purposes we focus on the interaction between $\uuu$ and $\vvv$. Let $\bz = (z_1,\ldots,z_M)$,  $\bv^1=\bz/\|\bz\|$, and complete $\bv^1$ to an orthonormal basis $(\bv^1,\ldots,\bv^M)$ of $\R^M$. Then we can rewrite
    \begin{align*}
    \uuu^2
    &= 
    \frac{1}{M}
    \frac{\|\bz\|^2}{M}
    \by^\st
    \bigg\{
    \frac{(\bar{\bG}^\st \bv^1)^{\otimes2}}{N}
    \bigg\} \by\,,\\
    \vvv
    &= \frac1M 
    \sum_{a=1}^M 
    \by^\st
    \frac{(\bar{\bG}\bv^a)^{\otimes2}}{N}
    \by
    = \uuu^2 + \tilde{\vvv}\,,
    \end{align*}
where $(\uuu^2,\tilde{\vvv})$ is roughly equidistributed as $(\bbb^2,\vvv)$. It follows that the feasible reason of $(\uuu^2,\vvv)$ is given by
    \[\bigg\{
    (\uuu^2,\vvv)
    : 
    \tilde{\vvv}=\vvv-\uuu^2\ge0,
    \uuu^2+\Big(
    (\vvv-\uuu^2)^{1/2}-1\Big)^2 \le \frac{1}{\alpha}
    \bigg\}\,.\]
A similar but slightly more complicated calculation gives the feasible region of $(\bbb^2,\vvv,\uuu^2)$ to be
    \[\bigg\{
    (\bbb^2,\vvv,\uuu^2):
    \vvv\ge\uuu^2,
    \Cost(\bbb,\vvv,\uuu)
    = \bbb^2+
    \uuu^2+\Big(
    (\vvv-\uuu^2)^{1/2}-1\Big)^2 
    \le
    \frac{1}{\alpha}
    \bigg\}\,.\]
This explains how free probability calculations can be used to derive our cost function. 
\end{xmp}

\begin{rmk}\label{r:free.prob.tight}
The limiting constraint in Theorem~\ref{t:free.prob} is sharp, though we do not need this fact.
  To see the mechanism behind the matching construction, consider the statistics $\bbb_j^2,\vvv_j,\uuu_j^2$ for $j\le r$ and $\www_j$ for $j\le s$.
  The budget constraint in Theorem~\ref{t:free.prob} (with $\epsilon = 0$) defines a convex body in these variables, and its support function is obtained by maximizing an arbitrary linear functional in these variables.
  The calculations in this subsection show that maximizing each such linear functional amounts to identifying the limiting maximum eigenvalue of a random matrix, which reduces to a free probability calculation.
  Conversely, each exposed boundary point of this convex body is the maximizer of some linear functional, and can be attained by taking $\by$ to be top eigenvector of the corresponding random matrix.
  The remaining (non-exposed) boundary points follow by taking convex combinations of the $\by$ corresponding to exposed boundary points, and interior points follow by taking a further convex combination with an independent standard gaussian vector in $\R^N$.
\end{rmk}

\fi

\pagebreak\section{Lipschitz algorithms and a priori estimates}\label{s:prelim}

\iffull
% !TEX root = main.tex

This section contains some preliminary observations on Lipschitz algorithms, and other basic estimates. The section is organized as follows:
\begin{itemize}
    \item In \S\ref{ss:basic.estimates} we recall some standard facts concerning gaussian matrices, subgaussian random vectors, and Lipschitz functionals on gaussian space.
    \item In \S\ref{subsec:cor-func-bounds} we show that any Lipschitz algorithm has a perturbation whose correlation function has derivatives bounded above and below. The bounded derivative assumption will be used repeatedly in subsequent sections.
    
\end{itemize}
As before, for any vector $\bx$, we write $\EmpDist(\bx)$
for the empirical distribution of coordinates of $\bx$.
Recall from \eqref{eq:proj-pursuit-def} that we let $\mu_{\bG}(\bx)$ denote the empirical distribution of the projection of data $\bG$ in direction $\bx$, 
    \[
    \mu_{\bG}(\bx)
    \equiv\EmpDist\bigg(\frac{\bG\bx}{N^{1/2}}\bigg)\,.
    \]
Given an algorithm $\cA=\cA(\bG,\bg^\aux)$ as in \eqref{eq:algorithms-as-maps}, we write 
    \beq\label{e:mu.of.Alg}
    \mu(\cA)
    =\bbE\Big[\mu_{\bG}(\cA(\bG,\bg^{\aux}))\Big]
    =\bbE \bigg[\EmpDist
    \bigg( \frac{\bG\cA(\bG,\bg^{\aux})}{N^{1/2}}\bigg)
    \bigg]\,,
    \eeq
where $\E$ denotes the expectation over $(\bG,\bg^\aux)$. Similarly, we denote
    \beq\label{e:mu.Ising.of.Alg}
    \mu^{\Ising}(\cA_N)
    =\E\Big[\EmpDist(\cA(\bG,\bg^{\aux}))\Big]
    =\E\bigg[
    \frac{1}{N}\sum_{i=1}^N \delta\{ \cA(\bG,\bg^{\aux})_i \}
    \bigg]\,.
    \eeq
Given $\mu\in\cP(\bbR)$, we write $\mu^{\sym}$ for the law of $sx$, where $x\sim\mu$ and $s\in\{\pm1\}$ is an independent symmetric random sign.
We write $\cP^{\sym}(\bbR)$ for the set of symmetric probability measures $\mu$ with $\mu=\mu^{\sym}$.

\subsection{Basic estimates}
\label{ss:basic.estimates}

In this subsection we recall some standard results (on subgaussian random variables, concentration of measure, and random matrix operator norms) and record some basic implications for in our problem setup.
Recall that a real-valued random variable $X$ is \emph{subgaussian with variance proxy $\sigma^2$} if it satisfies
    \[
    \E \exp (sX) \le \exp\bigg\{ 
        \frac{s^2 \sigma^2}{2}\bigg\}
    \]
for all $s\in\R$. It follows by a Chernoff bound that
    \[
    \P(|X|\ge t)
    \le 2\exp\bigg\{ -\frac{t^2}{2\sigma^2}\bigg\}\,.
    \]
As a result $X$ satisfies the moment bound
    \beq\label{e:subgaus.mmt}
    \E[|X|^p]
    =
    \int_0^\infty 
    pt^{p-1} \P(|X|\ge t)\,dt
    \le
    2\int_0^\infty 
    pt^{p-1}
    \exp\bigg\{-\frac{t^2}{2\sigma^2}\bigg\}\,dt
    \le C_p \sigma^p 
    \eeq
where $C_p$ is a constant depending only on $p$.
Next, a random vector $\bX \in \R^N$ is subgaussian with variance proxy $\sigma^2$ if $(\theta,\bX)$ is subgaussian with variance proxy $\sigma^2$ for all $\theta$ in the unit ball $B_N \subseteq \R^N$.

\begin{lem}\label{l:sg.norm.bound}
Suppose the random vector $\bX$ in $\R^N$ is \emph{subgaussian with variance proxy $\sigma^2$}, meaning that $(\theta,\bX)$ is subgaussian with variance proxy $\sigma^2$ for all $\theta$ in the unit ball $B_N\subset\R^N$. Then
    \[
    \P\bigg(\frac{\|\bX\|}{N^{1/2}} \ge t\bigg)
    \le \exp\bigg\{ -\frac{Nt^2}{16\sigma^2} \bigg\}
    \]
provided $t\ge 6\sigma$. It follows that
$\E(\|\bX\|^p) \le  C_p \sigma^p N^{p/2}$
where $C_p$ is a constant depending only on $p$.

\begin{proof}
To bound the norm of a random vector $\bX$ in $\R^N$, the standard argument is to let $N_{1/2}$ be a $1/2$-net of the unit ball $B\subset \R^N$, and note that
    \begin{align*}
    \|\bX\|
    &= \sup \Big\{
    (\theta,\bX) : \theta\in B\Big\}
    \le 
    \sup\bigg\{
    (\theta_1+\theta_2,\bX):
    \theta_1\in N_{1/2},
    \theta_2\in \frac12 B\bigg\} \\
    &\le
    \sup\Big\{(\theta,\bX) :\theta\in N_{1/2}\Big\}
    + \frac{\|\bX\|}{2}
    \equiv S + \frac{\|\bX\|}{2} \,.
    \end{align*}
Rearranging the above gives
    \[\|\bX\|
    \le 2\sup\Big\{(\theta,\bX) :\theta\in N_{1/2}\Big\}
    = 2S\,.
    \]
We can arrange for $|N_{1/2}| \le 6^N$. If $(\theta,\bX)$ is subgaussian with variance proxy $\sigma^2$ for all $\theta\in B_N$, then a union bound over the net $N_{1/2}$ gives
    \beq\label{e:subgaus.norm.conc}
    \P\bigg(\frac{\|\bX\|}{N^{1/2}} \ge 2t\bigg)
    \le \P(S \ge N^{1/2} t)
    \le 6^N 
    \exp\bigg\{ -\frac{Nt^2}{2\sigma^2}\bigg\}
    \le \exp\bigg\{ -\frac{Nt^2}{4\sigma^2}\bigg\}\,,
    \eeq
where the last bound holds provided
$t\ge 3\sigma$. Integrating (similarly to \eqref{e:subgaus.mmt}) gives the moment bound
    \beq\label{e:subgaus.vec.mmt.bd}
    \frac{\E[ \|\bX\|^p]}{N^{p/2}}
    \le
    \frac{2^p \E(S^p)}{N^{p/2}}
    \le (6\sigma)^p
    +2^p \int_{3\sigma}^\infty
    pt^{p-1}
    \exp\bigg\{ -\frac{Nt^2}{4\sigma^2}\bigg\}\, dt
    \le C_p \sigma^p
    \eeq
for a constant $C_p$ depending only on $p$.
\end{proof}
\end{lem}

\begin{lem}\label{l:lip.subgaus}
If $\bX:\R^{N'}\to\R^N$ is $L$-Lipschitz, and $\bg$ is a standard gaussian in $\R^{N'}$, then $\bX(\bg)-\E[\bX(\bg)]$ is subgaussian with variance proxy $2L^2$.
(As a result, Lemma~\ref{l:sg.norm.bound} applies with $\sigma^2=2L^2$.)

\begin{proof}
Recall from \cite[Thm.~1.3.4]{MR2731561} that if the function $F: \R^N\to\R$ is $L$-Lipschitz, and $\bg$ is a standard gaussian in $\R^N$, then
    \[
    \E \exp\bigg\{
     s\Big( F(\bg)-\E f(\bg)\Big)\bigg\}
    \le \exp(s^2 L^2)
    \]
for all $s\in\R$, i.e., the random variable $F(g)-\E F(g)$ is subgaussian
with variance proxy $2L^2$. As a consequence we have the subgaussian tail bound
    \beq\label{e:lip.subgaus}
    \P\bigg(
    \Big|  F(\bg)-\E f(\bg) \Big| \ge t
    \bigg)
    \le 2 \exp\bigg\{ -\frac{t^2}{4L^2}
    \bigg\}
    \eeq
for all $t\ge0$. Now suppose $\bX:\R^{N'}\to\R^N$ is $L$-Lipschitz. If $\theta$ is a unit vector in $\R^N$, then it is easily seen that $\bg\mapsto (\theta,\bX(\bg))$ is $L$-Lipschitz. It follows from \eqref{e:lip.subgaus} that the scalar random variable $(\theta,X(\bg)-\E[X(\bg)])$ is subgaussian with variance proxy $2L^2$, as claimed.
\end{proof}
\end{lem}

\begin{lem}\label{l:lip.norm.conc}
If $\bX:\R^{N'}\to\R^N$ is $L$-Lipschitz, $\bg$ is a standard gaussian in $\R^{N'}$, and $\E[\|\bX(\bg)\|^2] = N$, then for any $\delta > 0$, there exists $c(\delta) > 0$ such that 
    \[
    \P\bigg(
    \bigg|\frac{\|\bX(\bg)\|^2}{N} - 1\bigg|
    \ge \delta\bigg)
    \le e^{-c(\delta)N}\,.\]

\begin{proof}
    Since $\bg \mapsto \|\bX(\bg)\|$ is $L$-Lipschitz, it follows from Lemma~\ref{l:lip.subgaus} that $\|\bX(\bg)\|-\E\|\bX(\bg)\|$ is subgaussian with variance proxy $2L^2$, and satisfies the concentration bound \eqref{e:lip.subgaus}:
    \[
    \P\bigg( \Big|
    \|\bX(\bg)\| - \E\|\bX(\bg)\|\Big| \ge t
    \bigg) \le 2\exp\bigg\{-\frac{t^2}{2L^2}\bigg\}\,.
    \]
It follows that $\|\bX(\bg)\|$ has variance $O(L^2)$, which implies
    \[
     \frac{O(L^2)}{N}
     \ge
    \frac{\E[\|\bX(\bg)\|^2] - ( \E\|\bX(\bg)\|)^2}{N}
    =1- \frac{ ( \E\|\bX(\bg)\|)^2}{N}\,.
    \]
Combining the last two displays gives the claim.
\end{proof}
\end{lem}

\begin{lem}[{see \cite[Theorem 4.4.5]{wainwright2019high}}]
\label{l:wishart}
If $\bG$ is an $M\times N$ matrix with i.i.d.\ gaussian entries, then
    \[\P\bigg(
        \|\bG\|_\op
        \ge C_0 \Big(M^{1/2}+N^{1/2}\Big)
        \bigg)
    \le e^{-Nc_0} \]
for absolute constants $0<c_0,C_0<\infty$.
\end{lem}

\begin{lem}[$\bbW_2$-stability of $\mu_\bG(\bx)$]
    \label{l:w2-stability}
    For all $\bx,\bx' \in \bbR^N$, we have
    \[
        \bbW_2(\mu_\bG(\bx),\mu_\bG(\bx'))
        \le \frac{\|\bG\|_\op \|\bx-\bx'\|}{(MN)^{1/2}}\,.
    \]
\begin{proof}
The $\bbW_2$ distance can be bounded by
    \[
        \bbW_2(\mu_\bG(\bx),\mu_\bG(\bx'))^2
        \le \frac{1}{M} \sum_{a=1}^M \frac{(\bg^a,\bx-\bx')^2}{N}
        = \frac{\|\bG (\bx-\bx')\|^2}{MN}
        \le \frac{(\|\bG\|_\op \|\bx-\bx'\|)^2}{MN}\,.
    \]
This proves the claim.
\end{proof}
\end{lem}

\begin{rmk}\label{r:relaxed-wlog}
As a consequence of the basic estimates above, we note that any algorithmic output in the relaxed domains $\Sigma_N(\iota)$, $S_N(\iota)$ can be rounded to an element of $\Sigma_N$ and $S_N$ with nearly the same performance. 
That is, if $\bx = \cA(\bG,\bg^{\aux})$ satisfies conditions \eqref{it:coord-profile-succeed} and \eqref{it:proj-pursuit-succeed} in the first part of   Definition~\ref{d:attain}, then we can round $\bx$ to $\bx'$ lying in $\Sigma_N$ or $S_N$ such that, on the complement of the event in Lemma~\ref{l:wishart},
\[
    \bbW_2(\mu,\mu_\bG(\bx'))
    \le \bbW_2(\mu,\mu_\bG(\bx)) + \bbW_2(\mu_\bG(\bx),\mu_\bG(\bx'))
    \le \epsilon + C (\alpha^{1/2} + 1) \iota\,.
\]
\end{rmk}

\subsection{Correlation function bounds}
\label{subsec:cor-func-bounds}

In this subsection we establish bounds on the correlation function that will be used throughout the paper. As in Definition~\ref{d:Lip},
let $\cA=\cA(\bG, \bg^{\aux})$ be an $L$-Lipschitz algorithm. Recall from \eqref{e:p.corr.overlap} the definition of the correlation function $\chi$. In this subsection, we will show that we can essentially assume, without loss of generality,
\beq
    \label{e:correlation-fn-bds}
    \frac{1}{L^2}
    \le \chi'(t) \le L^2
\eeq
for all $0\le t\le 1$.
By Hermite expansion of $\cA$ (see \cite[Propn.~3.1]{HuangSellke2021}), $\chi$ is increasing and convex, and therefore it suffices to bound $\chi'(0)$ and $\chi'(1)$, which is done by the next two propositions:

\begin{ppn}\label{p:correlation-fn-ub}
We have $\chi'(1) \le L^2$.

\begin{proof}
Throughout this proof we abbreviate $\bg\equiv(\bG,\bg^\aux)$.  Consider small positive $\delta$, and define
    \[\begin{aligned}
    \bx&=\cA( (1-\delta)^{1/2}\bg+\delta^{1/2}\bg')\,,\\
    \by &= \E(\bx\,|\, \bg)\,,\\
    \bz &= \bx-\by\,.
    \end{aligned}\]
By iterated expectations, we have
    \[
    \frac{\E(\by,\by)}{N}
    =\frac{\E(\bx,\by)}{N}
    =\chi(1-\delta)\,, 
    \]
as well as $\E(\by,\bz)=0$. It follows that
    \beq\label{eq:bz}
    \chi(1)
    =\frac{\E(\bx,\bx)}{N}
    =\frac{\E(\by,\by)}{N}+\frac{\E(\bz,\bz)}{N}
    =\chi(1-\delta) + \frac{\E(\bz,\bz)}{N}\,.\eeq
Let us treat $\bz$ as a function of $\bg'$ for fixed $\bg$. Note that $\E(\bz\,|\,\bg)=0$. Since $\cA$ is $L$-Lipschitz, we see that $\bz$ is $(L\delta^{1/2})$-Lipschitz as a function of $\bg'$. It follows (using e.g.\ Lemma~\ref{l:lip.subgaus}) that
    \[
    \frac{\E(\bz,\bz)}{N} \le O(L^2\delta)\,.
    \]
Recalling \eqref{eq:bz}, we see that $\chi(1-\delta) \ge \chi(1)-O(L^2\delta)$, which implies $\chi'(1) \le L^2$ as claimed. 
\end{proof}
\end{ppn}

\begin{ppn}\label{p:wlogable}
    In both the spherical and Ising settings, there exists $c>0$ and $C = C(\alpha) > 0$ such that the following holds.
    Let $\cA^\circ$ be a $L$-Lipschitz algorithm with $L\ge 1$.
    Let $\chi^\circ \equiv \chi_{\cA^\circ}$ and $q_0 \equiv \chi^\circ(0)$, and suppose $q_0 < 1$.
    Then, there exists an $L$-Lipschitz algorithm $\cA$ such that the following hold:
\begin{enumerate}[(a)]
\item\label{it:wlogable.PERTURBED.CHI} $\chi \equiv \chi_{\cA}$ satisfies $\chi(0) = q_0$ and
    \[
        \chi'(0) \ge  \frac1{L^2}\,.
    \]
In fact, we can take
    \[
        \chi(p) = \bigg(1-\frac{1}{(1-q_0)L^2}\bigg) \chi^\circ(p)+ \frac{q_0 + (1-q_0)p}{(1-q_0)L^2}\,.
    \]

\item \label{it:wlogable.PERTURBED.APPX}
The perturbed algorithm $\cA$ approximates the original algorithm $\cA^\circ$ in the sense that
    \[
        \frac{1}{N} \E\Big[\|\cA^\circ(\bG,\bg^{\aux}) - \cA(\bG,\bg^{\aux})\|^2\Big]
        \le \frac{2}{(1-q_0)L^2}\,.
    \]

\end{enumerate}
We remark that the above error should be considered small, since we will prove later in Lemma~\ref{l:sM.Lip.1.trivial} that a sequence of algorithms with $q_0$ tending to one can only achieve the gaussian distribution.

\begin{proof} 
Let $\dot{\bg}$ be a standard gaussian in $\bbR^N$ independent of $\bg^{\aux}$, and
let $\bv$ be a deterministic element of $S_N$ with 
 $(\bv,\bbE \cA^\circ(\bG,\bg^{\aux})) = 0$. The perturbed algorithm $\cA$ is defined by
    \[\cA(\bG, (\bg^{\aux}, \dot{\bg} )) 
        = \bigg\{1-\frac{1}{(1-q_0)L^2}
            \bigg\}^{1/2}
        \cA^\circ(\bG,\bg^{\aux}) 
        + \bigg(\frac{q_0}{1-q_0}\bigg)^{1/2} \frac{\bv}{L}
        + \frac{\dot{\bg}}{L} \,.
    \]
Then the gradient of $\cA$ can be bounded as
    \begin{align*}
        \|\nabla\cA
            (\bG,(\bg^{\aux}, \dot{\bg})) \|^2
        &= \|\nabla_{\bG,\bg^{\aux}}
            \cA(\bG,(\bg^{\aux}, \dot{\bg} )) \|^2
        + \|\nabla_{\dot{\bg}} 
            \cA(\bG,(\bg^{\aux},\dot{\bg})) \|^2 \\
        &\le \bigg\{1-\frac{1}{(1-q_0) L^2}   \bigg\} L^2 
            + \frac1{L^2}
        \le L^2\,.
    \end{align*}
Since $\bv$ is deterministic and orthogonal to $\bbE \cA^\circ(\bG,\bg^{\aux})$,
 $\dot{\bg}$ is independent of $(\bG,\bg^{\aux})$, and $\bbE[\|\dot{\bg}\|^2] = N$, we have
    \[\frac{\bbE[
             \|\cA(\bG,(\bg^{\aux},\dot{\bg}))\|^2]}{N}
    =\bigg\{ 1-\frac1{(1-q_0)L^2}\bigg\}
    	+\frac{q_0}{(1-q_0)L^2}
	+\frac{1}{L^2}
    = 1\,,
    \]
where the last step follows from the trivial inequality $1-q_0\le1\le L^2$.  The above two displays show that $\cA$ is an $L$-Lipschitz algorithm in the sense of Definition~\ref{d:Lip}. Moreover, recalling the abbreviations $\chi\equiv\chi_\cA$ and $\chi^\circ\equiv\chi_{\cA^\circ}$, for any $0\le p\le1$ we have
    \[\chi(p) 
    = \bigg(1-\frac1{(1-q_0)L^2}\bigg) \chi^\circ(p) + \frac{q_0+(1-q_0)p}{(1-q_0)L^2}\,,
    \]
which implies $\chi(0)=q_0$.
   As noted above, $\chi^\circ$ is nondecreasing, which implies $\chi'(0)\ge 1/L^2$. This proves part~\eqref{it:wlogable.PERTURBED.CHI}. 

Next, since $|a^{1/2}-b^{1/2}| \le |a-b|^{1/2}$ for all $a,b\ge0$, we can bound
    \begin{align*}
        \frac{1}{N} \E\bigg[\Big\|\cA^\circ(\bG,\bg^{\aux}) - \cA(\bG,\bg^{\aux})\Big\|^2\bigg]
        &= \bigg(1 - \bigg(1-\frac{1}{(1-q_0)L^2}\bigg)^{1/2}\bigg)^2
        + \frac{q_0}{(1-q_0)L^2} + \frac{1}{L^2} \\
        &\le \frac{1}{(1-q_0)L^2} + \frac{q_0}{(1-q_0)L^2} + \frac{1}{L^2}
        = \frac{2}{(1-q_0)L^2}\,.
    \end{align*}
This proves part~\eqref{it:wlogable.PERTURBED.APPX}.
\end{proof}
\end{ppn}

\begin{rmk}\label{r:deterministic-lipschitz} 
Lastly, we show that the power of the class of algorithms in Definition~\ref{d:Lip} does not change if we remove the  auxiliary randomness $\bg^\aux$. Recall from Definition~\ref{d:Lip} that we call $\cA$ a \emph{deterministic} $L$-Lipschitz algorithm if it satisfies the requirements for an $L$-Lipschitz algorithm, and does not depend on any auxiliary randomness $\bg^\aux$. 
Suppose that for $\mu \in \cP_2(\R)$, there exists a $L$-Lipschitz Ising (resp. spherical) algorithm $\cA$ that $(\iota,\gamma)$-attains $\mu$: then we claim that there exists a deterministic $L$-Lipschitz algorithm $\cA'$ that also $(\iota',\gamma')$-attains $\mu$, with slightly worse parameters $\iota'$, $\gamma'$. Indeed, let $\cA$ be a general $L$-Lipschitz algorithm, in the sense of Definition~\ref{d:Lip}. Let $E=E(\iota)$ denote the event that it satisfies attainment properties \eqref{it:coord-profile-succeed} and \eqref{it:proj-pursuit-succeed}  from Definition~\ref{d:attain}. From that definition, we say that $\cA$ $(\iota,\gamma)$-attains $\mu$ if the event $E$ occurs with probability at least $\gamma$ over the randomness of $(\bG,\bg^\aux)$. If we consider $f(\bg^\aux)\equiv \P(E \,|\, \bg^\aux)$, and let $q=\P(f(\bg^\aux) \ge \gamma/2)$, then we can bound
	\[
	\gamma \le \E f(\bg^\aux)
	\le (1-q) \frac{\gamma}{2} + q
	\le \frac{\gamma}{2}+q\,,
	\]
so $q\ge\gamma/2$. On the other hand, if 
$h(\bg^\aux)
	\equiv \E[\|\cA(\bG,\bg^\aux)\|^2\,|\,\bg^\aux]$, then $h$ is well-concentrated around its mean value $N$, so we can arrange for
	\[
	\P\bigg(
	\bigg|\frac{h(\bg^\aux)}{N}-1\bigg|
	\ge\iota''\bigg)
	\le \frac{\gamma}{4}\,.
	\]
It follows that we can choose a realization $\bg^{\aux,\star}$ for which
	\[
	f(\bg^\aux) \ge \frac{\gamma}{2}
	\quad\textup{and}\quad
	\bigg|\frac{h(\bg^\aux)}{N}-1\bigg|
	\le\iota''\,.
	\]
The deterministic algorithm 
	\[
	\cA'(\bG)\equiv
	\frac{N^{1/2}\cA(\bG,\bg^{\aux,\star})}{h(\bg^{\aux,\star})^{1/2}}
	\]
satisfies the claim.
\end{rmk}

\fi

\pagebreak\section{SDE characterization of limits}
\label{s:sde}

% !TEX root = main.tex

\newcommand{\Num}{\tilde{N}}
\newcommand{\Den}{\tilde{D}}

Recall the processes constructed in Section~\ref{s:rerand}.
In particular, we introduced ``$\frozen$'' processes which were ``smoothed'' by the addition of a positive $\trK^{-O(1)}$ density of frozen particles on $[-2\trK,2\trK]^4$.
\textbf{Throughout this section we work with the ``$\frozen$'' processes only.}
An outline of this section is as follows:
\begin{itemize}
\item In \S\ref{ss:sde.abstract} we give an overview of this section in a simplified abstract setting, where $\vY$ is an $\R^d$-valued semimartingale with drift $\vD$ and covariation $\vQ$. We define $\rho_{\vY}$, $\rho_D$, and $\rho_Q$ to be the occupation measures of $\vY$, $\vD$, and $\vQ$ respectively. Then, under some regularity assumptions, we show how to construct an SDE solution $\hvX$ such that the marginal law of $\hvX(t)$ approximates the marginal law of $\vY(t)$ for each fixed $t$. The formal statement is given in
Lemma~\ref{l:compare.fokker.planck}, and uses a result of \cite{MR2375067} on uniqueness of solutions to Fokker--Planck equations under mild conditions. The coefficients of the SDE defining $\hvX$ are given by a certain smoothing of $\rho_D$ and $\rho_Q$, which we refer to as the ``ideal'' smoothing.

\item In \S\ref{ss:discrete.coeffs} we return to the discrete-time  process $\vY^{\frozen,N}$. We define the occupation measures $\rho_{\vY,N}$, $\rho_{D,N}$, and $\rho_{Q,N}$ associated to $\vY^{\frozen,N}$,
$\vD^{\frozen,N}$, and $\vtQ^{\frozen,N}$ respectively; these measures are directly related to the coefficients $(\bbb,\vvv,\uuu,\www,\rrr)^{\frozen,N}$
analyzed in previous sections. For technical reasons, it is not easy to show that the budget constraints from Section~\ref{ss:budget} pass through the ideal smoothing of these measures. We therefore show instead that the budget constraints pass through a modified smoothing procedure. 

\item In \S\ref{ss:limit.coefs} we consider the subsequential limiting semimartingale $\vY^\frozen$ obtained from 
Theorem~\ref{t:tightness}, with drift $\vD^\frozen$ and covariation $\vQ^\frozen$. (Passing to a further subsequence if needed, we assume that $p^N$ also converges to a limiting function $p$ in an appropriate sense.)  We define $\rho_{\vY}$, $\rho_D$, and $\rho_Q$ to be the occupation measures associated to $\vY^{\frozen}$, $\vD^{\frozen}$, and $\vQ^{\frozen}$, and we apply the modified smoothing procedure from \S\ref{ss:discrete.coeffs} to these limiting measures. We then show that the modified smoothing operation commutes with taking the limit $N\to\infty$, and consequently that the limiting smoothed coefficients also satisfy the budget constraints from Section~\ref{ss:budget}.

\item In \S\ref{ss:limit.sde} we apply the abstract result from \S\ref{ss:sde.abstract} to obtain an SDE solution $\hvX^{\frozen}$ such that the marginal law of $\hvX^{\frozen}(t)$ approximates the marginal law of $\vY^{\frozen}(t)$ for each fixed $t$. As noted above, the coefficients of the SDE defining $\hvX^{\frozen}$ are given by an ``ideal'' smoothing of the occupation measures $\rho_D$ and $\rho_Q$ associated to the drift and quadratic variation of $\vY^{\frozen}$. We finally show that this is well approximated by the SDE with coefficients given by the modified smoothing procedure, which were shown to satisfy the budget constraints from Section~\ref{ss:budget}.
\item Lastly, in \S\ref{ss:sde.reparam.sigma} we reparametrize and simplify the SDE resulting from \S\ref{ss:limit.sde}. \textbf{This yields Theorem~\ref{thm:BOGP-hardness-main} below, which is the main result of this section. Although further technical work is needed, Theorem~\ref{thm:BOGP-hardness-main} constitutes the most important step in our proof of hardness for Lipschitz algorithms.}
\end{itemize}

\begin{dfn}\label{d:incr.p} For $q_0 \in [0,1)$, let $\incr([q_0,1];[0,1])$ denote the set of increasing and absolutely continuous functions $p:[q_0,1] \to [0,1]$, with $p(1)=1$. Note that this includes the set of increasing and concave functions with $p(1)=1$. We do not require $p(q_0) = 0$ unless explicitly stated. Note that  $\incr([0,1];[0,1])$ coincides with the class $\sP$ from Section~\ref{s:intro}.
\end{dfn}

In what follows it will be useful to abbreviate
	\beq\label{e:s.sqrt.fn}
	s(t)\equiv
	\Big( (tp)'(t)\Big)^{1/2} 
\equiv \Big(
p(t)+t p'(t) \Big)^{1/2}\eeq
for any $p\in\incr([q_0,1];[0,1])$.  Recall the definitions of the averaged inner product and coordinate distributions $\mu(\cA_N),\mu^{\Ising}(\cA_N)$ of an algorithm $\cA_N$ from \eqref{e:mu.of.Alg} and \eqref{e:mu.Ising.of.Alg}. The first main result of this section is the following: 

\begin{thm}
\label{thm:BOGP-hardness-main}
Fix $L,\alpha>0$ and let $(\cA_N)^\circ$ be a sequence of $L$-Lipschitz algorithms in the sense of Definition~\ref{d:Lip}. Let $\cA_N$ be the perturbation of $(\cA_N)^\circ$ from Proposition~\ref{p:wlogable}. By passing to a subsequence $N_j\to\infty$, we assume that $\chi_{\cA_N}$ converges to a limiting function $\chi$, pointwise and in $L^1$. Denote $q_0\equiv \chi(0)\in[0,1]$. The function $p\equiv\chi^{-1}$ belongs to $\incr([q_0,1];[0,1])$ and is concave, with
\[\frac1{L^2}
\le p'(q) \le L^2\] for all $q\in [q_0,1]$. Suppose $(\mu,\mu^\Ising)$ satisfies 
	\beq\label{e:BOGP-hardness-main-hypothesis}
	\bbW_2\big(\mu,\mu(\cA_N)\big)
	+\bbW_2\big(\mu^\Ising,
		\mu^\Ising(\cA_N)\big)
		\le \epsilon_\textup{apx}
		\,.
	\eeq
Then there exist $C(L) > 0$ and $q_*,b,\sigma,w,\zeta,\zeta^\Ising$ such that the following hold. Let $B(t)$ and $W(t)$ be standard Brownian motions and
\begin{align} \label{e:sde.x}
        dX(t)
        &= p'(t)^{1/2} b_t\,dt
        + s(t) \sigma_t\,dB(t)\,,
        \\
        dX^{\Ising}(t)
        &= w_t\,dW(t) \label{e:sde.x.ising}
        \end{align}
on the time interval $q_* \le t \le 1$, with $X(q_*)\sim\zeta$ and $X^\Ising(q_*)\sim\zeta^\Ising$.  Then: 
\begin{enumerate}[(i)]
\item \label{it:BOGP-hardness-main-Lip}
$b_t=b(t,X(t))$,
$\sigma_t=\sigma(t,X(t))$,
and $w_t=w(t,Y(t))$ 
for functions $b:[q_*,1]\times \bbR\to[-C(L),C(L)]$ and $\sigma,w:[q_*,1]\times \bbR\to[0,C(L)]$, such that the functions $p'(t)^{1/2} b(t,x), s(t) \sigma(t,x), w(t,x)$, are $C(L)$-Lipschitz in $(t,x)$.
\item \label{it:BOGP-hardness-main-diffusivity}
We have the diffusivity constraint
        \beq\label{eq:wt-squared-at-most-1}
            \max\bigg\{
            \Big|\bbE[w(t,X^\Ising(t))^2]-1
            	\Big|
	:q_* \le t\le 1
\bigg\} \leq  o(1)\,.
\eeq

\item \label{it:BOGP-hardness-main-budget}
    We have the approximate budget constraints
       \beq \label{eq:budget-constraint-in-main-approximation-statement}
            \int_{q_*}^1
            \bigg\{
            \bbE\bigg[ b(t,X(t))^2+ \frac{(tp)'(t)}{p(t)}
            	(\sigma(t,X(t))-1)^2\bigg]
            -
            \frac{(\bbE w(t,X^\Ising(t)) )^2}{\alpha}
            \bigg\}_+
            \de t
            \leq o(1)\,.\eeq

 \item \label{it:BOGP-hardness-main-initial} The initial measures $\zeta,\zeta^\Ising$ satisfy $|\bE[X(q_*)^2]|\leq o(1)$ and  $|\bE[X^\Ising(q_*)^2] - q_*|\leq o(1)$.
\item \label{it:BOGP-hardness-main-qstar-condition} The initial time $q_*$ lies in $[q_0,q_0+ o(1)]$ and satisfies $p(q_*) \le o(1)$. 
\end{enumerate}
Lastly, we have
\beq\label{eq:BOGP-hardness-main-endpt}
\begin{aligned}
	\bbW_2(\Law(X(1)),\mu)
	&\leq 
	O(\epsilon_\textup{apx}) 
	+ o(1)\,,\\ 
	\bbW_2(\Law(X^{\Ising}(1)),\mu^{\Ising})
	&\leq  O(\epsilon_\textup{apx})
	+ o(1)\,.
\end{aligned}
\eeq
In all the above assertions, the $o(1)$ denotes an error which can tend to zero as $L^{-\Theta(1)}$. Thus, the distributions $\mu,\mu^\Ising$ are well approximated by the distributions of the endpoints of the solution to the constrained SDE described above.
\end{thm}

The 
\hyperlink{proof:t.BOGP-hardness-main}{proof of Theorem~\ref{thm:BOGP-hardness-main}} appears near the end of this section. In Section~\ref{sec:alternate-diffusions} we will relate several classes of  limiting SDEs, and the result of Theorem~\ref{thm:BOGP-hardness-main} will be applied in Lemma~\ref{l:Lip.in.BOGP}. The arguments of this section will also be needed at another point in Section~\ref{sec:alternate-diffusions}, in Lemma~\ref{l:prf3.in.IAMP}, to use
(Markovian) SDEs with Lipschitz coefficients
to approximate diffusions with general coefficients
(all subject to budget constraints).
This is encapsulated by the following more abstract formulation; and the required \hyperlink{proof:t.SDE-smoothing-general}{proof modifications} will be explained at the end of the section.

\begin{thm}
\label{thm:SDE-smoothing-general}
Given any $\acute{b},\acute{\sigma},\acute{w},p,\acute{\zeta},
\acute{\zeta}^\Ising$, let  
\begin{align*}
    d Z(t) &= p'(t)^{1/2} \ab_t\,dt + s(t) \asig_t\,dB(t)\,,
    \\
    dZ^{\Ising}(t)
    &= \aw_t\,dW(t)\,,
    \end{align*}
on the time interval $q_0 \le t \le 1$, with $Z(q_0)\sim\acute{\zeta}$ and $Z^\Ising(q_0)\sim\acute{\zeta}^\Ising$. Assume the following:
\begin{enumerate}[(a)]
\item \label{i:SDE-smoothing-general.INPUT.1}
The coefficients
$\ab,\asig,\aw$ are 
progressively measurable processes on the time interval $[q_0,1]$;

\item \label{i:SDE-smoothing-general.INPUT.2b}
$\bE[(\aw_t)^2]=1$ for all $t\in[q_0, 1]$;

\item \label{i:SDE-smoothing-general.INPUT.implies.3a} $p\in\incr([q_0,1];[0,1])$ with
$p(q_0)\ge 1/L$
and $\|p\|_{C^2([q_0,1])}\le L$;

\item \label{e:Z.budget.assumption}
The coefficients $\ab,\asig,\aw$ obey the budget constraints \eqref{eq:is-budget-constraint}, that is,
	\[
	\bE \bigg[(\ab_t)^2 + \frac{(tp)'(t)}{p(t)}(\asig_t-1)^2\bigg]
        \le \frac{(\bE \aw_t)^2 }{ \alpha}\,,\]
for all $t\in [q_0,1]$;

\item \label{i:SDE-smoothing-general.INPUT.5b}
$|Z^{\Ising}(1)|=1$ almost surely;

\item \label{i:SDE-smoothing-general.INPUT.6b}
$Z(q_0)=0$ almost surely, and $\bE[Z^{\Ising}(q_0)^2]=q_0$;

\item \label{i:SDE-smoothing-general.INPUT.7b}
$\zeta^\Ising\in\cP([-1,1])$,
the space of Borel probability measures on $[-1,1]$.\footnote{We note that since $Z^\Ising$ is a martingale with $|Z^{\Ising}(1)|=1$ almost surely, the initial distribution cannot have support outside of $[-1,1]$.}
\end{enumerate}
Given any
$q_*,b,\sigma,w,p,\zeta,\zeta^\Ising$
(same $p$ as above),
let $X,X^\Ising$ be given by \eqref{e:sde.x} and \eqref{e:sde.x.ising} on the time interval $[q_*,1]$, with $X(q_*)\sim\zeta$ and $X^\Ising(q_*)\sim\zeta^\Ising$.
Then for any $\epsilon>0$ there exist parameters $q_*,b,\sigma,w,\zeta,\zeta^\Ising$ such that the following hold: 
\begin{enumerate}[(i)]
\item \label{it:SDE-smoothing-general.coefs}
$b_t= b(t,X_t)$, $\sigma_t=\sigma(t,X_t)$, and $w_t= w(t,X_t)$,
for $C$-Lipschitz functions
$b:[q_*,1]\times\R\to[-C,C]$
and $\sigma,w:[q_*,1]\times\R\to[0,C]$, for a constant $C=C(L,\epsilon)$.
 
\item \label{it:SDE-smoothing-general.OUTPUT.diffus}
 $|\bE[(w_t)^2] -1| \le\epsilon$
 for all $t\in[q_*,1]$. 
 
\item \label{it:SDE-smoothing-general.budget}
 $X$ and $X^{\Ising}$ approximately obey the budget constraints at each time $t$, that is,
	\[\bE \bigg[(b_t)^2 + \frac{(tp)'(t)}{p(t)}(\sigma_t-1)^2\bigg]
        \le \frac{(\bE w_t)^2 }{ \alpha}
        +\epsilon\]
for all $t\in[q_*,1]$.

\item \label{it:SDE-smoothing-general.initial} The initial distributions of $(X,X^\Ising)$ and $(Z,Z^\Ising)$ are close:
\[\begin{aligned}
\bbW_2(\zeta,\acute{\zeta}) &= 
\bbW_2(\Law(X(q_*)),\Law(Z(q_0)))
\leq \epsilon\,,\\
\bbW_2(\zeta^\Ising,\acute{\zeta}^\Ising)
&= \bbW_2(\Law(X^{\Ising}(q_*)),
        	\Law(Z^{\Ising}(q_0)))
\leq \epsilon\,.
\end{aligned}\]

\item \label{it:SDE-smoothing-general.OUTPUT.qstar} 
The initial time $q_*$ lies in $[q_0,q_0+\epsilon]$.
\end{enumerate}
Lastly, we have
	\beq\label{eq:SDE-smoothing-general-endpt}
    \begin{aligned}
	\bbW_2(\Law(X(1)),\Law(Z(1)))
	&\leq \epsilon\,,\\
	\bbW_2(\Law(X^{\Ising}(1)),
        	\Law(Z^{\Ising}(1)))
	&\leq \epsilon
    \end{aligned}
	\eeq
Thus the more general diffusions $(Z,Z^\Ising)$ are well approximated by the Markovian SDE solutions
$(X,X^\Ising)$. Finally, if $\Law(Z^{\Ising}(t))$ is even for each $t$, then we can take $w_t$ to be even for each $t$, that is, $w(t,x)=w(t,-x)$ for all $(t,x)\in [0,1]\times \bbR$. 
\end{thm}

The \hyperlink{proof:t.SDE-smoothing-general}{proof of Theorem~\ref{thm:SDE-smoothing-general}}
is an adaptation of the 
\hyperlink{proof:t.BOGP-hardness-main}{proof of Theorem~\ref{thm:BOGP-hardness-main}},
and appears at the end of this section.  Both the proofs are based on a common template given by Theorem~\ref{thm:Lipschitz-SDE-approx}. 

\subsection{Ideal smoothing in abstract setting}
\label{ss:sde.abstract}

In this subsection we give an introductory overview of this section. This allows us to introduce some key objects, particularly the occupation measures (Definition~\ref{d:intro.cts.occ}) 
and the ``ideal'' SDE coefficients (Definition~\ref{d:intro.SDE.coefs}). At the end of this subsection we explain some of the technical issues that require modifications in our approach, which we implement in the remainder of this section.

Abstractly, let $\vY$ be an $\R^d$-valued semimartingale with drift $\vD$ and covariation $\vQ$. That is, $Y^i-D^i$ is a martingale for each $1\le i\le d$; and $Y^i Y^j-Q^{ij}$ is a martingale for each $1\le i,j\le d$. This is a generalization of the situation arising from Theorem~\ref{t:tightness}.  \textbf{In this subsection we generally use $x,y$ to denote vectors in $\R^d$, unless explicitly indicated otherwise.}

\begin{dfn}[occupation measures]
\label{d:intro.cts.occ}
In the setting described above,
let $\rho_{\vY}$ be the occupation measure of the $\vY$ process, so that
	\[\int_{t,y} 
	f(t,y)\,\rho_{\vY}(dt\,dy)
	= \bE
	\int_0^1 f(t,\vY(t))\,dt\,.
	\]
Let $\rho_{D^i}$ and $\rho_{Q^{ij}}$
be the occupation measures of the drift and quadratic variation processes, so that
	\begin{align*}
	\int_{t,y}
	f(t,y) \rho_{D^i}(dt\,dy) 
	&= \bE
	\int_0^1 f(t,\vY(t))\,dD^i(t)\,,\\
	\int_{t,y} f(t,y) 
	\rho_{Q^{ij}}(dt\,dy) 
	&= \bE
	\int_0^1 f(t,\vY(t))\,dQ^{ij}(t)\,.
	\end{align*}
for any continuous bounded function $f$.
\end{dfn}

We first consider our ``ideal smoothing'' which is given simply by the process $\hvY(t) = \vY(t-\sss)+\xxx$
where $(\sss,\xxx)$ is distributed according to a smooth density $\varphi$.
The next (easy) lemma records how the occupation density of the process transforms under this smoothing operation:

\begin{lem}[occupation density under ideal smoothing]\label{l:ideal.smoothed.occ.msr}
Let
$\hvY(t) = \vY(t-\sss)+\xxx$
where $(\sss,\xxx)\sim\varphi$ is sampled independently of $\vY$. If $\vY$ has occupation measure $\rho_{\vY}$ as given by Definition~\ref{d:intro.cts.occ}, then $\hvY$ has occupation measure
$\rho_{\vY}*\varphi$ (the convolution of $\rho_{\vY}$ with $\varphi$).\footnote{Note that at times close to $1$ or $q_0$, $\hvY$ may not in general be well defined, since $t-\sss$ can can lie outside the time interval $[q_0,1]$. In the full proof, we work on a subinterval $[q_*,1]\subsetneq [q_0,1]$ and choose $\varphi$ such that $0\leq \sss\leq q_0-q_*\ll 1$ almost surely, thus circumventing these issues. See \eqref{e:qbar.and.qstar} and Definition~\ref{d:t.smooth}.}

\begin{proof}
For the purposes of this proof let us abbreviate $\rho=\rho_{\vY}$. 
We then calculate
\begin{align*}
&\bE \int_t f(t,\hvY(t))\,dt
= \int_{s,x} \varphi(s,x) 
	\bE \int_t
	f(t+s,\vY(t)+x)\,dt \,ds\,dx\\
&\qquad=
	\int_{s,x} \varphi(s,x)
	\int_{t,y}
	\ind\{t+s\in[0,1]\}
	f(t+s,y+x)
	\,\rho(dt\,dy)\,ds\,dx\\
&\qquad= \int_{s',x'}
	f(s',x')
	\bigg\{\int_{t,y}
	\varphi(s'-t,x'-y)\,\rho(dt\,dy)
	\bigg\}
	\,ds'\,dx'\\
&\qquad= \int_{t,y}
	f(t,y)
	\rho*\varphi(t,y)
	\,dt\,dy\,.
\end{align*}
This proves the claim.
\end{proof}
\end{lem}

\begin{dfn}[SDE coefficients with ideal smoothing]\label{d:intro.SDE.coefs}
Abbreviating $\rho\equiv\rho_{\vY}$ for the occupation density of $\vY$, let us define
	\begin{align}
	\label{e:intro.occ.ratio.drift}
	B^i(t,y)
	&\equiv
	\frac{\rho_{D^i}*\varphi(t,y)}{\rho*\varphi(t,y)}\,,\\
	A^{ij}(t,y)
	&\equiv \frac{\rho_{Q^{ij}}*\varphi(t,y)}{\rho*\varphi(t,y)}\,.
	\label{e:intro.occ.ratio.qv}
	\end{align}
Let $\Sigma(t,y)$ be a matrix square root of $A(t,y)$, and let 
$\hvX$ be the solution to the SDE (in $d$ dimensions)
	\beq\label{e:abstract.sde}
	d\hvX(t) = B(t,\hvX(t))\,dt+
		\Sigma(t,\hvX(t))\,dW(t)\,,
	\eeq
where $W$ is a $d$-dimensional standard Brownian motion, with the initial condition that $\hvX(0)$ equidistributed as $\hvY(0)$. \textbf{In the abstract discussion of this section, we assume $B$ and $\Sigma$ are sufficiently regular that $\hvX$ is the unique solution of the SDE.}
\end{dfn}

\begin{lem}[comparison of Fokker--Planck equations]
\label{l:compare.fokker.planck}
Let
$\hvY(t) = \vY(t-\sss)+\xxx$
where $(\sss,\xxx)\sim\varphi$.
Let $\hvX$ be the SDE solution from Definition~\ref{d:intro.SDE.coefs}.
Assume the coefficients $B^i$ and $\Sigma^{ij}$ of \eqref{e:intro.occ.ratio.drift} and \eqref{e:intro.occ.ratio.qv} are bounded and Lipschitz. Then the processes $\hvX$ and $\hvY$ satisfy the same Fokker--Planck equations. Consequently, for each fixed $t$, the marginal law of $\hvX(t)$ agrees with the marginal law of\/ $\hvY(t)$.

\begin{proof}
Abbreviate $\rho\equiv\rho_{\vY}$ for the occupation density of $\vY$, and recall from Lemma~\ref{l:ideal.smoothed.occ.msr} that $\hvY$ has occupation density $\nu\equiv\rho_{\vY}*\varphi$.
Let $\mu$ denote the occupation density of $\hvX$.
 Let $f:\R^d\to\R$ be a smooth compactly supported test function. 
\smallskip

\noindent
\textbf{Fokker--Planck equations for $\hvY(t) = \vY(t-\sss)+\xxx$.} We will calculate
	\[
	\bE\Big[ f(\vY(T))-f(\vY(0))\Big]
	= \int_{s,x} \varphi(s,x)
	\bE\Big[ f(Y(T-s)+x)
	-f(Y(-s)+x)\Big]
	\,ds\,dx
	=\textup{(dr)}+\textup{(qv)}
	\]
where $\textup{(dr)}$
and $\textup{(qv)}$
denote the contributions from the drift and covariation respectively:
	\begin{align*}
	\textup{(dr)}
	&= \sum_{i=1}^d \textup{(dr)}^i
	= \sum_{i=1}^d
	\int_{s,x}\varphi(s,x)
	\int_{-s}^{T-s}
	\partial_i f(Y(t)+x)\,dD^i(t)
	\,ds\,dx \,,\\
	\textup{(qv)}
	&= \sum_{i,j=1}^d \textup{(qv)}^{ij}
	= \frac12
	\sum_{i,j=1}^d
	\int_{s,x}\varphi(s,x)
	\int_{-s}^{T-s}
	\partial_i\partial_j
	f(Y(t)+x)\,dQ^{ij}(t)
	\,ds\,dx\,.
	\end{align*}
Abbreviate $I_T(s)\equiv \ind\{s\in[0,T]\}$. The above can be further simplified as follows:
	\begin{align*}
	\textup{(dr)}^i
	&= \int_{s,x}\varphi(s,x)
	\int_{t,y} I_T(t+s)
	\partial_i f(y+x)
	\,\rho_{D^i}(dt\,dy)\,ds\,dx \\
	&=
	\int_{s',x'}
	I_T(s')\partial_i f(x')
	\bigg\{\int_{t,y} 
	\varphi(s'-t,x'-y)
	\,\rho_{D^i}(dt\,dy)
	\bigg\}\,ds'\,dx'\\
	&=
	\int_{t,y}
	I_T(t)\partial_i f(y)
	\frac{\rho_{D^i}*\varphi(t,y)}
	{\rho*\varphi(t,y)}
	\nu(dt\,dy)\,,
	\end{align*}
recalling at the last step that $\nu=\rho*\varphi$ is the occupation density of $\hvY$.\smallskip

\noindent
\textbf{Fokker--Planck equations for SDE solution $\hvX$.}
 If $f:\R^d\to\R$ is a smooth compactly supported test function, define
	\[
	L_t f(x)
	= \sum_{i=1}^d
	B^i(t,x)
	\partial_i f(x)
	+\frac12\sum_{i,j=1}^d
	A^{ij}(t,x)
	\partial_i\partial_j f(x)
	\]
Then, by It\=o's formula, the process
	\[
	N(t)
	= f(\hvX(t))-f(\hvX(0)) - \int_0^t L_s f(\hvX(s))\,ds
	\]
is a local martingale. Since we assumed the coefficients are bounded, we can conclude that $N$ is a true martingale, and the optional stopping theorem gives
	\[
	\bE\Big[ f(\hvX(T))-f(\hvX(0))\Big]
	= \bE \int_0^T L_t f(\hvX(s))\,ds
	= \textup{(Dr)}+\textup{(QV)}\,,
	\]
where $\textup{(Dr)}$ and $\textup{(QV)}$ denote the contributions from the drift and covariation respectively:
	\begin{align*}
	\textup{(Dr)}
	&= \sum_{i=1}^d \textup{(Dr)}^i
	=\sum_{i=1}^d
	\bE \int_0^T \partial_i f(\hvX(t))
	B^i(t,\hvX(t))\,dt\,,\\
	\textup{(QV)}
	&= \sum_{i,j=1}^d \textup{(QV)}^{ij}
	=\frac12\sum_{i=1}^d
	\bE \int_0^T 
	\partial_i\partial_j f(\hvX(t))
	A^{ij}(t,\hvX(t))\,dt\,.
	\end{align*}
Again abbreviate $I_T(s)\equiv \ind\{s\in[0,T]\}$. The above can be further simplified as
	\begin{align*}
	\textup{(Dr)}^i
	&= \int I_T(t)
	\partial_i f(y) 
	\frac{\rho_{D^i}*\varphi(t,x)}
		{\rho*\varphi(t,x)}
	\mu(dt\,dy)\,,\\
	\textup{(QV)}^{ij}
	&= \int I_T(t)
	\partial_i\partial_j f(y) 
	\frac{\rho_{Q^{ij}}*\varphi(t,x)}
		{\rho*\varphi(t,x)}
	\mu(dt\,dy)
	\end{align*}
recalling that $\mu$ is the occupation density of $\hvX$, and using the definitions \eqref{e:intro.occ.ratio.drift} and \eqref{e:intro.occ.ratio.qv} of $B^i$ and $A^{ij}$.

Comparing the expressions for $\hvX$ and $\hvY$ shows that $\mu$ and $\nu$ satisfy the same Fokker--Planck equations. We now apply \cite[Lem.~2.3]{MR2375067}: this is a general result applying to $\R^d$-valued stochastic differential equations
such as \eqref{e:abstract.sde}, with a given law for the initial position.
By the assumption of Lipschitz coefficients, all solutions to the martingale problem corresponding to $L_t$ have the same marginal distribution at each time $t\in [0,1]$  (see e.g.\ \cite{MR2190038}). Then \cite[Lem.~2.3]{MR2375067} says that the associated Fokker--Planck equation also has a unique solution (as a probability measure-valued function of time), from which the conclusion follows.
\end{proof}
\end{lem}

\textbf{Our approach in this section is largely based on Lemma~\ref{l:compare.fokker.planck}.} However, it will be difficult to ensure sufficient regularity for the most obvious matrix square root $\Sigma(t,y)$ for the process $(Y^{\RomI},Y^{\RomII},Y^{\RomIII})$. For this reason, we will instead work with 
the coordinate sum $Y^{\RomI}+Y^{\RomII}+Y^{\RomIII}$ (which is ultimately the object of study), where the required regularity is easier to show. Thus, for later use we record how the objects introduced above transform under this  projection:

\begin{dfn}[coordinate sum]
\label{d:projected.coeffs}
For any measure $\rho$ on $[0,1]\times\R^d$, let $\bar{\rho}$ be the measure on $[0,1]\times\R$ that satisfies
	\[
	\int_{t,y} f\bigg(t, \sum_{i=1}^d y^i\bigg)
	\,\rho(dt\,dy)
	= \int_{t,x} f(t,x)
	\,\bar{\rho}(dt\,dx)\,,
	\]
where $t\in[0,1]$, $y\in\R^d$, and $x\in\R$.
As above, let $\vY$ be an $\R^d$-valued semimartingale with drift $\vD$ and covariation $\vQ$. We now consider the one-dimensional projections
	\[
	Y = \sum_{i=1}^d Y^i\,,\quad
	\sum_{i=1}^d D^i\,,\quad
	\sum_{i,j=1}^d Q^{ij}\,.
	\]
Then $Y$ is a (one-dimensional) semimartingale with drift $D$ and covariation $Q$. We hereafter abbreviate $\rho\equiv\rho_{\vY}$ for the occupation measure of $\vY$ defined in Definition~\ref{d:intro.cts.occ}.
Then the occupation measure of $Y$ is $\bar{\rho}$. The occupation measures associated to the drift and covariation of $Y$ are
	\[
	\sum_{i=1}^d
	\bar{\rho}_{D^i}\,,\quad
	\sum_{i,j=1}^d
	\bar{\rho}_{Q^{ij}}\,.	
	\]
Let $\hvY$ be as defined by Lemma~\ref{l:ideal.smoothed.occ.msr},
and let $\iY$ be the one-dimensional projection, so $\iY(t)=Y(t-\sss)+\mathfrak{x}$ where $(\sss,\mathfrak{x})\sim\bar{\varphi}$ (the density $\varphi$ integrated over $x^1 + \ldots + x^d=\mathfrak{x}$). It follows that $\iY$ has occupation measure $\overline{\rho*\varphi}=\bar{\rho}*\bar{\varphi}$. Similarly, the occupation measures associated to the drift and covariation of $\iY$ have densities
	\[
	\sum_{i=1}^d
	\overline{\rho_{D^i}*\varphi}
	\,,\quad
	\sum_{i,j=1}^d
	\overline{\rho_{Q^{ij}}
		*\varphi}\,.	
	\]
Then, analogously to \eqref{e:intro.occ.ratio.drift} and \eqref{e:intro.occ.ratio.qv}, we define (with $x$ now a scalar)
	\begin{align}
	\label{e:proj.coef.drift}
	B(t,x) 
	&\equiv
	\sum_{i=1}^D \frac{
	\overline{\rho_{D^i}*\varphi}(t,x)
	}
		{\overline{\rho*\varphi}(t,x)}\,,\\
	A(t,x)
	&\equiv
	\sum_{i,j=1}^D \frac{
	\overline{\rho_{Q^{ij}}*\varphi}(t,x)
	}
	{\overline{\rho*\varphi}(t,x)}\,.
	\label{e:proj.coef.qv}
	\end{align}
Let $\iX$ be the solution to the (one-dimensional) SDE
	\[
	d\iX(t)
	= B(t,\iX)\,dt
	+A(t,\iX)^{1/2}\,dW(t)\,,
	\]
where $W$ is a one-dimensional standard Brownian motion.
\end{dfn}

\begin{lem}[SDE for coordinate sum]
\label{l:projected.fokker.planck}
Let
$\hvY(t) = \vY(t-\sss)+\xxx$
where $(\sss,\xxx)\sim\varphi$, and let 	\[
	\iY(t)
	\equiv \sum_{i=1}^d \iY^i(t)\,.
	\]
Let $\iX$ be the SDE solution from Definition~\ref{d:projected.coeffs}.
Assume the coefficients $B$ and $A^{1/2}$ of \eqref{e:proj.coef.drift} and \eqref{e:proj.coef.qv} are bounded and Lipschitz.
Then the processes $\iX$ and $\iY$ satisfy the same Fokker--Planck equations. 
Consequently, for each fixed $t$, the marginal law of $\iX(t)$ agrees with the marginal law of\/ $\iY(t)$.
\end{lem}

\begin{proof}
This is a direct consequence of the multidimensional result
Lemma~\ref{l:compare.fokker.planck}. 
\end{proof}

As mentioned above, in our final application we will not apply Lemma~\ref{l:compare.fokker.planck} directly, but will aim to apply the one-dimensional version Lemma~\ref{l:projected.fokker.planck} where regularity of coefficients is easier to guarantee. Moreover, we will ultimately work with a modification of the ``ideal'' coefficients \eqref{e:proj.coef.drift} and \eqref{e:proj.coef.qv}, which make it easier for us to pass the budget constraints of Section~\ref{ss:budget} to the limiting SDE. This plan is carried out in the remainder of this section.

\subsection{Discrete coefficients and smoothing}
\label{ss:discrete.coeffs}

In this subsection, we define and analyze discrete-time analogues of the occupation measures and SDE coefficients discussed in the previous section, now in reference to the processes $\vY^{\frozen,N}$ from Definition~\ref{d:rerand.Ising}. 
\textbf{We again emphasize that throughout this section we work with the ``$\frozen$'' processes only.}
The subsection is organized as follows:
\begin{itemize} 
\item
We first consider ``ideal smoothing'' as in Definition~\ref{d:intro.SDE.coefs}, but \textbf{in the spatial coordinate only}, then combined with a coordinate sum projection as in Definition~\ref{d:projected.coeffs} ---
see Definition~\ref{d:sp.smooth}. We then show that the constraints of  Proposition~\ref{p:budget-constraints-coarsened} pass to the spatially smoothed coefficients; see Lemmas \ref{l:sp.smooth.domain} and \ref{l:sp.smooth.budget}. We also show some regularity estimates, Lemmas \ref{l:sp.smooth.reg} and \ref{l:sp.smooth.reg.SDE}, to be applied later.
\item We then introduce a \textbf{temporal smoothing} rule, Definition~\ref{d:t.smooth}, which deviates from the ``ideal smoothing,'' but which allows us to easily pass the constraints from 
Lemmas \ref{l:sp.smooth.domain} and \ref{l:sp.smooth.budget} to the smoothed quantities. In Lemma~\ref{l:ideal.vs.actual.coeffs}
we will address the discrepancy between the ``ideal'' and the ``actual'' coefficients (Definition~\ref{d:intro.SDE.coefs} versus Definition~\ref{d:t.smooth}).
\end{itemize}

We take a moment to describe some preliminary reductions in the \hyperlink{proof:t.BOGP-hardness-main}{proof of Theorem~\ref{thm:BOGP-hardness-main}}; this will be occasionally be relevant in the analysis that follows. Given $(\cA_N)^\circ$ as in the statement of the theorem, we immediately pass to the perturbed algorithms $\cA_N$ from Proposition~\ref{p:wlogable}. Let $\chi^N\equiv\chi_{\cA^N}$ and $p^N\equiv(\chi^N)^{-1}$. By Propositions~\ref{p:correlation-fn-ub} and \ref{p:wlogable}, we have
\beq\label{e:sdelimits.pprime}
	\frac{1}{L^2} \leq (p^N)'(q)
	\leq L^2
	\eeq
for all $q\in [q_0,1]$.
Apply Theorem~\ref{t:tightness} to the sequence $\cA_N$, and let $\Q_N\to \Q$ be any subsequential limit. Passing to a further subsequence if needed, we assume that $\chi_{\cA_N}$ converges to a limiting function $\chi$, pointwise and in $L^1$. We also always assume $q_0=\chi(0)<1$ since Theorems~\ref{thm:BOGP-hardness-main} and \ref{thm:SDE-smoothing-general} become trivial when $q_0=1$. In the setting of Theorem~\ref{thm:BOGP-hardness-main}, by construction, the function $p\equiv\chi^{-1}$ is concave and nondecreasing, with $p(q_0)=p_0=0$.
By concavity, for almost all $q\in [q_0,1]$, the derivative $p'(q)$ exists and agrees with the limit $\lim_{N\to\infty} (p^N)'(q)$: for such $q$, the value of $p'(q)$ must also satisfy the bound \eqref{e:sdelimits.pprime}. Most of the analysis of this section applies to this setting: see in particular the statement of Theorem~\ref{thm:Lipschitz-SDE-approx}, which we will use to deduce Theorem~\ref{thm:BOGP-hardness-main}.

In the setting of Theorem~\ref{thm:BOGP-hardness-main}, to avoid issues of dividing by $p(q_0)=0$, we slightly shrink the interval $[q_0,1]$ as follows. Abbreviate $\brep\equiv 1/L^{\breve{c}}$, where $\breve{c}$ is a large absolute constant to be determined (see \eqref{e:use.of.breve.c}). Choose $\breve{q}$ such that $p(\breve{q})=\brep$. The bounds on $p'$ imply
	\beq\label{e:breve.q.bounds}
	\frac{1}{L^{\breve{c}+2}}
	\le \breve{q}-q_0
	\le \frac{L^2}{L^{\breve{c}}}\,.
	\eeq
We then define the parameters
\beq\label{e:qbar.and.qstar}
\begin{aligned}
\bar q &=q_0+2(\breve{q}-q_0),\\
q_* &=q_0+4(\breve{q}-q_0).
\end{aligned}\eeq
In some of the analysis below, we will restrict to the interval $[q_*,1]$ to avoid issues of dividing by zero.

Let $p^N:[q_0,1]\to[0,1]$ be the continuous function which satisfies
$p^N(q_d)=p_d$, and is linear on each interval $[q_d,q_{d+1}]$. We will discuss only the Ising case, since the spherical case is simpler.  \textbf{For clarity, in what follows we generally write $x,y$ for scalars and $\vx,\vy$ for vectors.}
	
\begin{dfn}[discrete-time coefficients]\label{d:disc.coefficients}
Let $\vY^N$ and $\vY^{\frozen,N}$ be the spatial rerandomization processes from Definition~\ref{d:rerand.Ising}. 
Recall the coefficients $(\bbb,\vvv,\uuu,\rrr,\www)^{\trunc,N}$ from \eqref{e:def.b}--\eqref{e:def.w}, which were defined with respect to the process $\vY^\trunc$. We also have the analogous coefficients
$(\bbb,\vvv,\uuu,\rrr,\www)^{\frozen,N}$ from \eqref{e:b.plus.hist}--\eqref{e:w.plus.hist}, which were defined with respect to the process $\vY^\frozen$. Finally, the coarsened coefficients
$(\bbb,\vvv,\uuu,\rrr,\www)^{\frozen,\coarse,N}$ are defined by \eqref{e:b.plus.coarse}--\eqref{e:w.plus.coarse}.
 We will view these coarsened coefficients as functions on $\R^5$, as follows. Partition time into intervals $I_d=[q_d,q_{d+1})$, and 
space into rectangles
	$R_{\vl,j}
	\equiv J_{\vl}
	\times K_j$, following the notation of Definition~\ref{d:blocks} (the partitions depend on $N$). With a mild abuse of notation, we let
$(\bbb,\vvv,\uuu,\rrr)^{\frozen,N}$
denote the piecewise constant functions on $\R^5$ such that, for all $(t,\vy) = (t,y^1,y^2,y^3,y^4)\in I_d\times R_{\vl,j}$, we have
	\begin{align}\label{e:review.b}
	\bbb^{\frozen,N}(t,\vy) &\equiv 
	\bbb^{\frozen,\coarse,N}(q_d,B^{\frozen,\coarse})
	= \sum_{B^\frozen
	\subseteq
	 B^{\frozen,\coarse}}
	 \frac{|B^\frozen|}
	 	{|B^{\frozen,\coarse}|}
	\frac{\Delta \bar{D}^{\frozen,\RomI,N}(q_d,B^\frozen)}
	{\delta_d (p^N)'(q_d)^{1/2}}\,,\\ \label{e:review.v}
	\vvv^{\frozen,N}(t,\vy) 
	&\equiv 
	\vvv^{\frozen,\coarse,N}(q_d,B^{\frozen,\coarse})
	= \sum_{B^\frozen
	\subseteq
	 B^{\frozen,\coarse}}
	 \frac{|B^\frozen|}
	 	{|B^{\frozen,\coarse}|} 
	\frac{\Delta\tilde{Q}^{\frozen,\RomII,N}(q_d,B^\frozen)}{\delta_d p^N(q_d)}
	\,, \\ \label{e:review.u}
	\uuu^{\frozen,N}(t,\vy) &\equiv
	\uuu^{\frozen,\coarse,N}
	(q_d,B^{\frozen,\coarse})
	\equiv 
	\sum_{B^\frozen
	\subseteq
	 B^{\frozen,\coarse}}
	 \frac{|B^\frozen|}
	 	{|B^{\frozen,\coarse}|} 
	\frac{\Delta\tilde{Q}^{\frozen,\RomII,\RomIII,N}(q_d,B^\frozen)}
		{ \delta_d
		[ p^N(q_d) \cdot q_d (p^N)'(q_d)]^{1/2}}
		\,,\\
	\label{e:review.P}
	\rrr^{\frozen,N}(t,\vy)
	&\equiv 
	\rrr^{\frozen,\coarse,N}
		(q_d,B^{\frozen,\coarse})
	\equiv
	\sum_{B^\frozen
	\subseteq
	 B^{\frozen,\coarse}}
	 \frac{|B^\frozen|}
	 	{|B^{\frozen,\coarse}|} 
	\frac{\Delta C^{\frozen,\RomIII}
		(q_d,B^\frozen)}{\delta_d}\,,
	\end{align}
where $B^{\frozen,\coarse}=B^{\frozen,\coarse}(q_d,J_{\vl})$
as defined by \eqref{e:coarsened.bucket}. The last equality in \eqref{e:review.b} follows from combining  \eqref{e:b.plus.coarse} and \eqref{e:b.plus.hist}; the last equality in each of 
\eqref{e:review.v}--\eqref{e:review.P} is similar. Note that $B^{\frozen,\coarse}$ does not depend on $y^4$, so the functions in \eqref{e:review.b}--\eqref{e:review.P} also do not depend on $y^4$. Likewise, we let $\www^{\frozen,N}$ denote the piecewise constant function on $\R^5$ such that, for all $(t,\vy)\in I_d\times R_{\vl,j}$, we have
	\beq\label{e:review.w}
	\www^{\frozen,N}(t,\vy) \equiv 
	\www^{\frozen,\coarse,N}(q_d,B^{\frozen,\coarse,\Ising})
	\equiv
	\sum_{B^{\frozen,\Ising}
	\subseteq
	 B^{\frozen,\coarse,\Ising}}
	 \frac{|B^{\frozen,\Ising}|}
	 {|B^{\frozen,\coarse,\Ising}|}
	\frac{\Delta \tilde{Q}^{\Isplus,N}(q_d,B^\Isplus)}{\delta_d}
	\,,
	\eeq
where $B^{\frozen,\coarse,\Ising}=B^{\frozen,\coarse,\Ising}(q_d,K_j)$
is defined analogously to \eqref{e:coarsened.bucket}. Note that $B^{\frozen,\coarse,\Ising}$ does not depend on $(y^1,y^2,y^3)$, so the function in \eqref{e:review.w} also does not depend on $(y^1,y^2,y^3)$. Finally, we take the functions $(\bbb,\vvv,\uuu,\rrr,\www)^\frozen$ to be identically zero on $\R^5$ outside of
$(t,\vy) \in [0,1]\times[-\trK,\trK]^4$.
We hereafter abbreviate
$\Gamma^{\frozen,N}\equiv(\bbb,\vvv,\uuu)^{\frozen,N}$.
\end{dfn}

We define the following discrete-time analogues of the occupation measures from Definition~\ref{d:intro.cts.occ}: 

\begin{dfn}[discrete-time occupation measures]
\label{d:discrete.occ}
Partition time into intervals $I_d \equiv [q_d,q_{d+1})$, and partition space into rectangles 
	$R_{\vl,j}
	\equiv J_{\vl}
	\times K_j$, as above.
Let $\rho_{\vY,N}$ be the occupation density for $\vY^{\frozen,\bullet,N}$: that is, we let $\rho_{\vY,N}$ be uniform on each $I_d\times R_{\vl,j}$, with
	\[
	\frac{\rho_{\vY,N}(I_d\times R_{\vl,j})}{|I_d|}
	= \mu_{\bullet,\bG^N}
	\Big( \vY^{\frozen,\bullet,N}(q_d) \in R_{\vl,j}\Big)
	\,,\]
where $\mu_{\bullet,\bG^N}$ is as defined in Theorem~\ref{t:tightness}.  Recall from Definition~\ref{d:rerand.Ising} that in the spatial rerandomization, the $a(q_d)\in[M]$ indices and the $i(q_d)\in[N]$ indices evolve completely independently, so the processes $\vY^{\frozen,\bullet,N,\RomI:\RomIII}$ and $Y^{\frozen,\bullet,N,\Ising}$ are independent under $\mu_{\bullet,\bG^N}$. This implies the above is equal to
	\beq\label{e:a.i.factorization}
	\mu_{\bullet,\bG^N}
	\bigg( 
	\vY^{\frozen,\bullet,
	N,\RomI:\RomIII}(q_d)
	 \in J_{\vl}
	 \bigg)
	\cdot
	\mu_{\bullet,\bG^N}
	\Big( Y^{\frozen,\bullet,N,\Ising}(q_d) 
	\in K_j \Big)\,.
	\eeq
We similarly let $\rho_{D,N}$ be the signed measure which is uniform on each $I_d\times R_{\vl,j}$, with
	\begin{align}\nonumber
	&\frac{\rho_{D,N}
	(I_d\times R_{\vl,j})}{|I_d|}
	= 
	\bE_{\bG^N} \bigg[
	\frac{\Delta \bar{D}^{\frozen,\RomI,N}(q_d)}
	{\delta_d}
	\ind\Big\{
	\vY^{\frozen,\bullet,N}(q_d) \in R_{\vl,j}
	\Big\}\bigg] \\
	&\qquad\stackrel{\eqref{e:review.b}}{=}
	(p^N)'(q_d)^{1/2}
	\bigg\{ \sum_{B^\frozen \in 
		B^{\frozen,\coarse}}
		\frac{|B^\frozen|}
			{|B^{\frozen,\coarse}|}
	\bbb^{\frozen,N}\Big(q_d, 
		B^\frozen \Big)
	\bigg\}
	\mu_{\bullet,\bG^N}
		\Big( 
		\vY^{\frozen,\bullet,N}(q_d) \in R_{\vl,j}\Big) \nonumber\\
	&\qquad\stackrel{\eqref{e:review.b}}{=}
	(p^N)'(q_d)^{1/2}
	\bbb^{\frozen,\coarse,N}(q_d, B^{\frozen,\coarse}) 
	\mu_{\bullet,\bG^N}
		\Big( 
		\vY^{\frozen,\bullet,N}(q_d) \in R_{\vl,j}\Big)
	\label{e:def.rho.b.N}
	\end{align}
with $B^{\frozen,\coarse}=B^{\frozen,\coarse}(q_d,J_{\vl})$ as defined by \eqref{e:coarsened.bucket}, and recalling that the last term above factorizes as in \eqref{e:a.i.factorization}. We similarly define measures $\rho_{V,N}$, $\rho_{U,N}$, 
$\rho_{\www,N}$,
and $\rho_{\rrr,N}$
 using increments of 
$\tilde{Q}^{\frozen,\RomII,N}$, $\tilde{Q}^{\frozen,\RomII,\RomIII,N}$, $\tilde{Q}^{\frozen,\Ising,N}$, and 
$C^{\frozen,\RomIII,N}$ respectively. 
\end{dfn}

Now recall the ``ideal'' coefficients defined by 
\eqref{e:intro.occ.ratio.drift} and \eqref{e:intro.occ.ratio.qv} in Definition~\ref{d:intro.SDE.coefs}. In the discrete setting, since the occupation densities are piecewise constant, we have enough regularity to define analogous coefficients without any smoothing, as follows: 

\begin{dfn}[discrete-time coefficients]
\label{d:disc.rescaled.coeffs}
We let $B^N$ be the piecewise constant function on $\R^5$ such that
for all $(t,\vy)\in I_d\times R_{\vl,j}$, we have
	\beq\label{e:discrete.ratio.occ.B}
	B^N (t,\vy) 
	\equiv \frac{\rho_{D,N}(I_d\times R_{\vl,j})}
		{\rho_{\vY,N}(I_d\times R_{\vl,j})}
	\stackrel{\eqref{e:review.b}}{=}
	(p^N)'(t)^{1/2} \cdot
	\bbb^{\frozen,N}(t,\vy)\,.\eeq
Likewise, let $V^N$ and $U^N$ be the piecewise constant functions on $\R^5$ such that
	\begin{align}
	\label{e:discrete.ratio.occ.V}
	V^N(t,\vy) 
	&\equiv 
	\frac{\rho_{V,N}(I_d\times R_{\vl,j})}
		{\rho_{\vY,N}(I_d\times R_{\vl,j})}
	\stackrel{\eqref{e:review.v}}{=}
	p^N(t) 
	\cdot \vvv^{\frozen,N}(t,\vy)\,,\\
	U^N(t,\vy) 
	&\equiv
	\frac{\rho_{U,N}(I_d\times R_{\vl,j})}
		{\rho_{\vY,N}(I_d\times R_{\vl,j})}
	\stackrel{\eqref{e:review.u}}{=}
	\Big[ p^N(t)\cdot
	t (p^N)'(t) 
	\Big]^{1/2}
	\cdot
	\uuu^{\frozen,N}(t,\vy)
	\label{e:discrete.ratio.occ.U}
	\end{align}
for all $(t,\vy)\in I_d\times R_{\vl,j}$. Similarly, for all $(t,\vy)\in I_d\times R_{\vl,j}$ we have the identities
	\beq\label{e:discrete.ratio.occ.r.w}
	\rrr^{\frozen,N}(t,\vy)
	= \frac{\rho_{\rrr,N}(I_d\times R_{\vl,j})}
		{\rho_{\vY,N}(I_d\times R_{\vl,j})}\,,\quad
	\www^{\frozen,N}(t,\vy)
	= \frac{\rho_{\www,N}(I_d\times R_{\vl,j})}
		{\rho_{\vY,N}(I_d\times R_{\vl,j})}\,.
	\eeq
It follows from Definitions~\ref{d:disc.coefficients} and \ref{d:discrete.occ} that $(B,V,U,\rrr)^N$ do not depend on $y^4$; we will often abuse notation and regard them as functions of $(t,y^1,y^2,y^3)\in\R^4$. Likewise, $\www^{\frozen,N}$ does not depend on $(y^1,y^2,y^3)$; we will often abuse notation and regard it as a function of $(t,y^4)\in\R^2$. We let  $\Psi^N\equiv(B,V,U)^N$.
\end{dfn}

\textbf{We next introduce and analyze smoothed variants of the above coefficients, in order to carry out (a modification of) the approach outlined in Lemmas \ref{l:compare.fokker.planck} and \ref{l:projected.fokker.planck}.} 
Recalling Definition~\ref{d:discrete.occ}, let $\rho_t\equiv \rho_{\vY,N,t}$ denote the measure $\rho_{\vY,N}$ conditioned on time $t$; define similarly the measures $\rho_{D,N,t}$, $\rho_{V,N,t}$ $\rho_{U,N,t}$, $\rho_{\rrr,N,t}$, and $\rho_{\www,N,t}$. We first consider the effect of \textbf{spatial smoothing alone}, along with the projection onto the coordinate sum. The below coefficients are analogous to \eqref{e:proj.coef.drift} and \eqref{e:proj.coef.qv} from Definition~\ref{d:projected.coeffs}, with the only difference that we convolve with a smooth kernel $\sphi$ in the spatial coordinate only, rather than in both space and time.

\begin{dfn}[SDE coefficients with spatial smoothing
	and coordinate sum projection] 
\label{d:sp.smooth}
Let $\sphicirc:\R\to[0,\infty)$ be a smooth nonnegative function that integrates to one and has support $[-1,1]$. We also assume that we have
$\sphicirc(\vx)\ge \Omega(1)$ for all $\vx\in[-1/2,1/2]$. With an abuse of notation, for $\vx\in\R^d$ for any $d\ge1$ we will denote
	\beq
	\label{e:spatial.smoothing.product}
	\sphicirc(\vx)
	= \prod_{i=1}^d
	\sphicirc(x^i)\,,\quad
	\sphi(\vx)\equiv 
	\frac{1}{\sep^d}
	\sphicirc\bigg(
	\frac{\vx}{\sep}\bigg)\,,
	\eeq
so that $\sphi$ also integrates to one and has support $[-\sep,\sep]^d$.
Recall from above that we use the abbreviation 
$\rho_t\equiv \rho_{\vY,N,t}$. It will also be convenient to factorize $\rho_t=\varrho_t\otimes\varsigma_t$ where $\varrho_t$ refers to the first three spatial coordinates while $\varsigma_t$ refers to the fourth spatial coordinate (cf.\ \eqref{e:a.i.factorization}). We then have $\rho_{D,t} = \varrho_{D,t}\otimes\varsigma_t$, etc., as well as $\rho_{\www,t}= \varrho_t\otimes\varsigma_{\www,t}$.
 Recall from Definition~\ref{d:projected.coeffs} the notation $\rho\mapsto\bar{\rho}$ for the coordinate sum projection.
Then, analogously to
\eqref{e:discrete.ratio.occ.B}--\eqref{e:discrete.ratio.occ.U}, we let
	\begin{align*}
	\sB^N(t,y)
	&\equiv
	\frac{\overline{\varrho_{D,N,t}*\sphi}(y)}
		{\overline{\varrho_{\vY,N,t}*\sphi}(y)}
	\equiv
	(p^N)'(t)^{1/2} \cdot \sbbb^N(t,y)\,,\\
	\sV^N(t,y)
	&\equiv
	\frac{\overline{\varrho_{V,N,t}*\sphi} (y)}
		{\overline{\varrho_{\vY,N,t}*\sphi} (y)}
	\equiv
	p^N(t) \cdot \svvv^N(t,y)\,,\\
	\sU^N(t,y)
	&\equiv
	\frac{\overline{\varrho_{U,N,t}*\sphi} (y)}
		{\overline{\varrho_{\vY,N,t}*\sphi} (y)}
	\equiv
	\Big[ p^N(t) \cdot 
	t (p^N)'(t) \Big]^{1/2}
		\cdot \suuu^N(t,y)\,.
	\end{align*}
Likewise, analogously to \eqref{e:discrete.ratio.occ.r.w}, we let
	\[
	\srrr^N(t,y)\equiv
	\frac{\overline{\varrho_{\rrr,N,t}*\sphi}(y)}
		{\overline{\varrho_{\vY,N,t}*\sphi}(y)}\,,\quad
	\swww^N(t,y^4)
	\equiv
	\frac{\varsigma_{\www,N,t}*\sphi(y^4)}
		{\varsigma_{\vY,N,t}*\sphi(y^4)}\,.
	\]
We emphasize again that we have combined spatial smoothing with coordinate sum projection in this definition, so the above coefficients are all functions on $\R^2$ rather than $\R^5$. Lastly, we shall abbreviate
$\sPsi^N\equiv(\sB,\sV,\sU)^N$ and $\sGam^N\equiv(\sbbb,\svvv,\suuu)^N$; these are variants of the quantities
$\Psi^N$ and $\Gamma^N$ from Definitions \ref{d:disc.coefficients} and \ref{d:disc.rescaled.coeffs}.
\end{dfn}

Recalling the notations introduced in Definitions \ref{d:disc.coefficients}--\ref{d:disc.rescaled.coeffs}, the domain constraint \eqref{e:domain.plus.coarse} from Proposition~\ref{p:budget-constraints-coarsened} can be rewritten as
	\begin{align}
	\nonumber
	&\frac{V^N(t,\vy)}{p^N(t)}
	\stackrel{\eqref{e:discrete.ratio.occ.V}}{=}
	\vvv^{\frozen,N}(t,\vy)
 \stackrel{\eqref{e:domain.plus.coarse}}{\ge}
	q_d p'(q_d)
    \frac{\uuu^{\frozen,N}(t,\vy)^2 }
    {\rrr^{\frozen,N}(t,\vy)}
    - \frac{\eta^2}{2}
	\\
	&\qquad
	\stackrel{\eqref{e:discrete.ratio.occ.U}}{=}
	\frac{U^N(t,\vy)^2}{p^N(t)
		\cdot \rrr^{\frozen,N}(t,\vy)}
		- \frac{\eta^2}{2}
	\ge \frac{U^N(t,\vy)^2}{p^N(t)
		\cdot t (p^N)'(t)}
		- \frac{\eta^2}{2}\,,
	\label{e:domain.plus.rewrittten}
	\end{align}
for all $q_0 \le t\le 1$. The next lemma says that this constraint passes to the spatially smoothed quantities. Hereafter we will often suppress $N$ from the notation; so for example we will further abbreviate $\rho\equiv \rho_{\vY,N}$ and $\rho_t\equiv \rho_{\vY,N,t}$.

\begin{lem}[spatial smoothing of domain constraint]\label{l:sp.smooth.domain}
Let $\sV^N$ and $\sU^N$ be the spatially smoothed quantities from Definition~\ref{d:sp.smooth}. Then, analogously to \eqref{e:domain.plus.coarse}, with high probability we have
	\[
	\svvv^N(t,y)
	=\frac{\sV^N(t,y)}{p^N(t)}
	\ge
	\frac{\sU^N(t,y)^2}{p^N(t)
		\cdot  \srrr^N(t,y)}
	-\frac{\eta^2}{2}
	\ge 
	\frac{\sU^N(t,y)^2}{p^N(t)
		\cdot  t (p^N)'(t)}
	-\frac{\eta^2}{2}
	=\suuu^N(t,y)^2
	-\frac{\eta^2}{2}
	\]
for all $(t,y) \in[q_0,1]\times\R$.

\begin{proof}
We use the abbreviations introduced above. In particular, recall from Definition~\ref{d:disc.rescaled.coeffs} that we can regard $V$ as a function of $(t,\vy)$ for $\vy\in\R^3$. 
Using 
\eqref{e:discrete.ratio.occ.V}, we can rewrite
	\[
	\rho_{V,t}*\sphi(\vy)
	= \int_{\R^3} \sphi(\vy-\vx) V(t,\vx)\,\varrho_t(d\vx)\,.
	\]
It follows that we can express 
$\sV(t,y)\equiv 
\sE_{t,y} V(t,\xxx)$ where $\sE_{t,y}$ 
refers to expectation with respect to the law of $(\vy,\xxx)\sim\tilde{\bP}_{t,y}$ defined by
	\beq\label{e:tilde.bP}
	\tilde{\bP}_{t,y}(d\vy\,d\vx)
	= \frac{ 
		\varrho_t*\sphi(\vy)}
		{\overline{\varrho_t*\sphi}
		(y)}
		\cdot 
		\frac{\sphi(\vy-\vx)\,
		\varrho_t(d\vx)}
		{\varrho_t*\sphi(\vy)} \, d\vy
	\eeq
for $\vx\in\R^3$ and $\vy\in H(y)$ where
	\beq\label{e:hspace.y}
	H(y)
	\equiv \bigg\{\vy\in\R^3 : (\vI,\vy)=y\bigg\}\,.\eeq
Here $H(y)$ is equipped with the Hausdorff measure inherited from Lebesgue measure on $\R^3$, and in \eqref{e:tilde.bP} $d\vy$ denotes integration with respect to this measure. Similarly, $\sU(t,y) = \sE_{t,y} U(t,\xxx)$ and $\srrr(t,y)=\sE_{t,y} \rrr(t,\xxx)$. One can check that the function
	\[
	(u,r) \mapsto \begin{cases}
	 \frac{u^2}{r} & r>0\,,\\
	0 & u=r=0
	 \end{cases}
	\]
is convex on the domain $\{u\in\R,r>0\}  \cup\{u=r=0\}$, so Jensen's inequality gives
	\[
	\frac{\sV(t,y)}{p(t)}
	= \frac{\sE_{t,y} V(t,\xxx)}{p(t)}
	\stackrel{\eqref{e:domain.plus.rewrittten}}{\ge}
	\sE_{t,y} \bigg[
		\frac{U(t,\xxx)^2}{p(t) 
		\cdot \rrr(t,\xxx)}
		- \frac{\eta^2}{2}\bigg]
	\ge 
	\frac{\sU(t,y)^2}{p(t) 
	\cdot  \srrr(t,y)}
		- \frac{\eta^2}{2}\,,
	\]
as claimed.
\end{proof}
\end{lem}

Next we define the following slight modifications of the original cost function \eqref{e:cost}: let $0<\bar{\lambda}\ll\epsilon^\circ$ (to be determined; see \eqref{e:def.bar.lambda} below), and define 
	\beq\label{e:cost.t.eps}
	\Cost_{t,2\bar{\lambda}}(B,V,U)
	\equiv \frac{B^2}{p'(t)}
	+ \frac{U^2}{p(t) \cdot  t p'(t)}
	+\bigg(
	\frac{1}{p(t)^{1/2}} 
	\bigg[\bigg(V-
	\frac{U^2}{ t p'(t) } 
	+ 2\bar{\lambda}
	\bigg)_+\bigg]^{1/2}-1\bigg)^2\,.
	\eeq
As in Definition~\ref{d:disc.rescaled.coeffs}, let
$\Psi^N\equiv(B,V,U)^N$.
Since $\bar{\lambda}\ll\epsilon^\circ$, 
\eqref{e:budget.constraint.Ising.plus.coarse} implies that, with high probability,
	\begin{align}\nonumber
	&C^N(t)
	\equiv
	\int_{\R^4} \Cost_{t,2\bar{\lambda}}
	(\Psi^N(t,\vy))
		\,\rho_{\vY,N,t}(d\vy)\\ 
	&\qquad\le
	\frac{1}{\alpha}
	\bigg(
	\int_{\R^4} \www^{\frozen}(t,\vy)^{1/2} \,\rho_{\vY,N,t}(d\vy)
	\bigg)^2
	+ 3\epsilon^\circ
	\equiv \frac{W^N(t)^2}{\alpha}+ 3\epsilon^\circ\,.
	\label{e:budget.Ising.plus.eps}
	\end{align}
The next lemma says that \eqref{e:budget.Ising.plus.eps} behaves well with respect to spatial smoothing:

\begin{lem}[spatial smoothing of budget constraint]
\label{l:sp.smooth.budget}
Analogously to \eqref{e:budget.Ising.plus.eps}, with high probability we have
	\begin{align*}\nonumber
	&\tilde{C}^N(t)
	\equiv
	\int_{\R}
	\Cost_{t,2\bar{\lambda}}(\sPsi^N(t,y)) \,
		\overline{\varrho_{\vY,N,t}*\sphi}(y)\,dy\\ 
	&\qquad\le
	\frac{1}{\alpha}
	\bigg(
	\int_{\R} \swww(t,y)^{1/2} 
		\,\varsigma_{\vY,N,t}*\sphi(y)\,dy
	\bigg)^2
	+ 3\epsilon^\circ
	\equiv \frac{\tilde{W}^N(t)^2}{\alpha}+ 3\epsilon^\circ
	\end{align*}
for all $q_0\le t\le 1$. 

\begin{proof}
We continue to use the abbreviations introduced above.\smallskip

\noindent\textbf{Relation between $\tilde{C}^N(t)$ and $C^N(t)$.}
 As in Definition~\ref{d:sp.smooth}, let $\sPsi\equiv(\sB,\sV,\sU)$. The function $\Cost_{t,2\bar{\lambda}}$ is convex on the  expanded domain 
	\[\bigg\{(B,V,U)
	: \frac{V}{p(t)}
	\ge \frac{U^2}{p(t) \cdot t p'(t)} 
	-  \frac{2\bar\lambda}{p(t)}\bigg\}\,.\]
It follows from Lemma~\ref{l:sp.smooth.domain}
that $\sPsi$ lies in this expanded domain with high probability. 
Recall from Definition~\ref{d:disc.rescaled.coeffs} that we can regard $\Psi=(B,V,U)$ as functions of $(t,\vy)$ with $\vy\in\R^3$.
It then follows by Jensen's inequality that 
	\begin{align*}
	&\Cost_{t,2\bar{\lambda}}\Big(\sPsi(t,y)\Big)
	= \Cost_{t,2\bar{\lambda}}
	\Big( \sE_{t,y} \Psi(t,\xxx)\Big)
	\le \sE_{t,y}
	\bigg[ \Cost_{t,2\bar{\lambda}}
		\Big(\Psi(t,\xxx)\Big)\bigg] \\
	&\qquad
	\stackrel{\eqref{e:tilde.bP}}{=}
	\int_{H(y)}\int_{\R^3}
	\Cost_{t,2\bar{\lambda}} (\Psi(t,\vx))
	\frac{\varrho_t*\sphi(\vy)}
		{\overline{\varrho_t*\sphi}(y)}
		\cdot \frac{\sphi(\vy-\vx)\,\varrho_t(d\vx)}
			{\varrho_t*\sphi(\vy)} \,d\vy\,.
	\end{align*}
Integrating the above bound with respect to $\overline{\varrho_t * \sphi}(y)\,dy$ gives
	\begin{align*}
	\tilde{C}^N(t) &=
	\int_{\R} \Cost_{t,2\bar{\lambda}}\Big(\sPsi(t,y)\Big)
		\overline{\varrho_t * \sphi}(y)\,dy
	\le
	\int_{\R^3}
	\Cost_{t,2\bar{\lambda}} (\Psi(t,\vx)) \bigg[
	\int_{\R^3}
	\sphi(\vy-\vx) \,d\vy\bigg] \,\varrho_t(d\vx) \\
	&= \int_{\R^3}
	\Cost_{t,2\bar{\lambda}} (\Psi(t,\vx))\, \varrho_t(d\vx) = C^N(t)\,,
	\end{align*}
having used that $\sphi$ integrates to one.\smallskip

\noindent\textbf{Relation between $\tilde{W}^N(t)$ and $W^N(t)$.}
Similarly, recall from Definition~\ref{d:disc.rescaled.coeffs} that we can regard $\www^\frozen$ as a function of $(t,y^4)\in\R^2$. We can then express
$\swww(t,y) = \sE_{t,y,\Ising} \www^\frozen(t,\xxx)$ where $\sE_{t,y,\Ising}$ denotes expectation with respect to the law of $(\vy,\xxx)\sim\tilde{\bP}_{t,y,\Ising}$ defined by
	\[
	\tilde{\bP}_{t,y,\Ising}(dy\,dx)
	= \frac{\sphi(y-x)\,\varsigma_t(dx)}
			{\varsigma_t*\sphi(y)} \, dy\,.
	\]
for $x,y\in\R$.
Then, since $w\mapsto w^{1/2}$ is concave in $w\ge0$, we have
	\begin{align*}
	&\swww(t,y)^{1/2}
	=\Big( \sE_{t,y,\Ising} 
	\www^\frozen(t,\xxx) \Big)^{1/2}
	\ge\sE_{t,y,\Ising} 
	\Big[\www^\frozen(t,\xxx)^{1/2}\Big] \\
	&\qquad\stackrel{\eqref{e:tilde.bP}}{=}\int_{\R}\int_{\R}
		\www^\frozen(t,x)^{1/2}
	\frac{\sphi(y-x)\,\varsigma_t(dx)}
		{\varsigma_t*\sphi(y)} \,dy
	\,.
	\end{align*}
Integrating the above bound with respect to $\varsigma_t * \sphi(y)\,dy$ gives
	\begin{align*}
	\tilde{W}^N(t)
	&=
	\int_{\R} \swww(t,y)^{1/2} 
		\varsigma_t * \sphi(y)\,dy
	\ge
	\int_{\R} \www^\frozen(t,x)^{1/2}
	\bigg[\int_{\R}\sphi(y-x)\,dy\bigg] \, \varsigma_t(dx)\\
	&=\int_{\R} 
	\www^\frozen(t,x)^{1/2}\, \varsigma_t(d\vx)
	= W^N(t)\,.
	\end{align*}
Combining with \eqref{e:budget.Ising.plus.eps} gives
	\[
	\tilde{C}^N(t)
	\le C^N(t)
	\stackrel{\eqref{e:budget.Ising.plus.eps}}{\le}
	\frac{W^N(t)^2}{\alpha}+ 3\epsilon^\circ
	\le
	\frac{\tilde{W}^N(t)^2}{\alpha}+ 3\epsilon^\circ\,,
	\]
as claimed.
\end{proof}
\end{lem}

For later use (for instance, in the proof of Lemma~\ref{l:t.smooth.budget}, and in ensuring uniqueness of SDE solutions), we next argue that the spatially smoothed coefficients are well-behaved. \textbf{At this point it is important to remember that throughout this section we work with the ``$\frozen$'' processes from Definition~\ref{d:trunc}.} In the next lemma the added density of frozen particles is essential to ensure bounded coefficients.

\begin{lem}[regularity for spatially smoothed coefficients]
\label{l:sp.smooth.reg}
The coefficients from Definition~\ref{d:sp.smooth} are continuously differentiable in the space coordinate, and satisfy the bounds
	\[
	\Big(
	\|\sB^N\|_\infty
    +
    \|\partial_y\sB^N\|_\infty\Big)^2
	+ \max\bigg\{
    \|\tilde{A}\|_\infty
    +\|\partial_y\tilde{A}\|_\infty
    : \tilde{A}\in\{ \sV,\sU,\srrr,\swww\}^N
    \bigg\}
	\le 
	\bigg(\frac{\trK L}{\sep}\bigg)^{O(1)} 
	(\MAX_N)^2\,,
	\]
with $\MAX_N$ as in \eqref{e:MAX.N} and Proposition~\ref{p:Y.kolmogorov}.

\begin{proof}
We continue to use the abbreviations introduced above. To simplify notation, we mostly suppress the dependence on $N$ inside the proof. Consider
	\[
	\sB(t,y)
	=\frac{\overline{\varrho_{D,t}*\sphi}(y)}
		{\overline{\varrho_t*\sphi}(y)}
	\equiv \frac{\Num(t,y)}{\Den(t,y)}\,.
	\]
We first aim to lower bound the denominator $\Den$.
Recall from Definition~\ref{d:sp.smooth} that $\sphi$ has support $[-\sep,\sep]^3$, and is $\ge \Omega(1)/\sep^3$ on $[-\sep/2,\sep/2]^3$. 
From the construction of $\vX^\frozen$ and $\vY^\frozen$
(Definitions \ref{d:trunc}--\ref{d:rerand.Ising}), we see that $\varrho_t$ gives measure at least $\Omega( \sep^3 /\trK^{O(1)} )$ to any box of side length $\sep$ contained in $[-2\trK,2\trK]^3$. It follows that, for all $\vy\in[-2\trK+1,2\trK-1]^3$,
	\[
	\varrho_t*\sphi(\vy)
	\ge \int_{[-\sep/2,\sep/2]^3}
		\sphi(\vy-\vx)\,\varrho_t(d\vx)
	\ge \frac{1}{\trK^{O(1)}}\,.
	\]
Integrating over $\vy\in H(y)$ (as defined by \eqref{e:hspace.y}) gives
	\beq\label{e:denominator.lbd}
	\Den(t,y)
	= \overline{\varrho_t*\sphi}(y)
	\ge \frac{1}{\trK^{O(1)}}
	\eeq
for all $y\in[-5\trK,5\trK]$. Next note that since particles outside $[-\trK,\trK]^3$ are frozen (see Definition~\ref{d:trunc}), the numerator $\Num(t,y)$ must be zero for $y\notin[-3\trK-1,3\trK+1]$. For $y\in[-3\trK-1,3\trK+1]$, we bound
	\begin{align*}
	|\Num(t,y)|
	&=\Big|\overline{\varrho_{D,t}*\sphi}(y)\Big|
	=\bigg| \int_{H(y)}
	\int_{\R^3} \sphi(\vy-\vx)B(t,\vx)\,\varrho_t(d\vx) \,d\vy\bigg|\\
	&=\bigg| \int_{\R^3}
	B(t,\vx) \overline{\sphi}(y-(\underline{1},\vx))\,\varrho_t(d\vx)\bigg|
	\le \|\overline{\sphi}\|_\infty \cdot \int_{\R^3}
	|B(t,\vx)|\,\varrho_t(d\vx)
	\le  \bigg(\frac{L}
		{\sep}\bigg)^{O(1)}
		\MAX_N\,,
	\end{align*}
where the last step uses the quenched moment bound on $\bbb$ from Lemma~\ref{l:barD.I.second.quenched}, together with the bound $(p^N)' \le  L^2$ from \eqref{e:sdelimits.pprime}. Combining with \eqref{e:denominator.lbd} gives
	\[\|\sB\|_\infty 
	\le \bigg(\frac{\trK L}
		{\sep}
		\bigg)^{O(1)}
		\MAX_N\,.
	\]
We can also apply \eqref{e:denominator.lbd} to obtain
	\beq\label{e:sB.quotient.rule}
	\Big|\partial_y \tilde{B}(t,y)\Big|
	\le \trK^{O(1)}
	\bigg|
	\frac{\partial}{\partial y} \Num(t,y)
	- \tilde{B}(t,y)
	\frac{\partial}{\partial y} \Den(t,y)
	\bigg|\,.\eeq
The first term on the right-hand side of \eqref{e:sB.quotient.rule} can be bounded by
	\begin{align*}
	\bigg|\frac{\partial}{\partial y} \Num(t,y)\bigg|
	&=\bigg|\frac{\partial}{\partial y}
	\overline{\varrho_{D,t}*\sphi}(y)\bigg|
	=\bigg| 
	\int_{\R^3} B(t,\vx)
	\overline{\sphi}' (y-(\vec{1},\vx))
	\,\varrho_t(d\vx)
	\bigg|\\
	&\le
	\|\overline{\sphi}' \|_\infty
	\int_{\R^3} | B(t,\vx) | \,\varrho_t(d\vx)
	\le \bigg(\frac{L}
		{\sep}\bigg)^{O(1)} \MAX_N\,.
	\end{align*}
The second term on the right-hand side of \eqref{e:sB.quotient.rule} can be controlled using the preceding bound on $\tilde{B}$, together with a similar  calculation as above:
	\[
	\bigg|\frac{\partial}{\partial y}
		\Den(t,y) \bigg|
	=\bigg|\frac{\partial}{\partial y}
	\overline{\varrho_t*\sphi}(y)\bigg|
	\le  \| \overline{\sphi}' \|_\infty \cdot
	\int_{\R^3} \varrho_t(d\vx)
	\le \frac{1}{\sep^{O(1)} }\,.
	\]
Combining the above gives the claimed bound on $\partial_y \sB$.
The bounds on $\sV$, $\sU$, $\srrr$, and $\swww$ follow by the same argument, using the quenched moment bounds from Lemma~\ref{l:tQ.first.quenched}.
\end{proof}
\end{lem}

Recalling
Definition~\ref{d:projected.coeffs}
and Lemma~\ref{l:projected.fokker.planck}, the quadratic variation of the combined process
	\beq\label{e:Y.plus}
	Y^\frozen \equiv 
	\sum_{\sigma\in\{\RomI,
		\RomII,\RomIII\}} 
	Y^{\frozen,\sigma}
	\eeq
will be approximately captured by
	\beq\label{e:sde.F.N}F^N(t,\vy)
	\equiv
	\Big\{ V+2U+\rrr^\frozen\Big\}
	(t,\vy)+2\eta^{1/2}\,.
	\eeq
It follows from \eqref{e:domain.plus.rewrittten} that for each $t\geq q_0$, with high probability, 
	\[
	U(t,\vy)^2
	\le p(t) \rrr^\frozen(t,\vy)
	\bigg( \frac{V(t,\vy)}{p(t)} +\frac{\eta^2}{2}\bigg)
	\le V(t,\vy)
	\rrr^\frozen(t,\vy) +
	O(L^{O(1)} \eta^2)\,,
	\]
where the last step uses 
\eqref{e:sdelimits.pprime}. We claim that this implies $F^N\ge0$: 
by rearranging, this can be reduced to the inequality
	\[
	\bigg(V+\rrr^\plus 
	+2\eta^{1/2}\bigg)^2
	\ge
	(V+\rrr^\plus )^2
	+ 4\eta
	\ge 
	4 \bigg(
	V\rrr^\plus + 
	O( L^{O(1)} \eta^2)
	\bigg)
	\ge 4U^2\,,
	\]
having used that $\eta\ll 1/L$
from Assumption~\ref{a:params}.
It follows that $F^N(t,\vy)\ge0$ for all $t\ge q_0$. 
Let $\rho_F$ be the discrete occupation density associated with $F$. Then, recalling Definition~\ref{d:sp.smooth} and \eqref{e:sde.F.N}, let us define
	\beq\label{e:sp.smoothed.F.N}
	\sF^N(t,y)
	\equiv \frac{\overline{\varrho_{F,N,t}*\sphi}(y)}
		{\overline{\varrho_{\vY,N,t}*\sphi}(y)}
	=
	\Big\{ \sV
		+2\sU+\srrr \Big\}(t,y)
		+2\eta^{1/2}\,.
	\eeq
Since $F^N\ge0$ and $\sF^N$ is obtained by smoothing $F^N$, we must also have $\sF^N\ge0$ on $t\geq q_0$. Lemma~\ref{l:sp.smooth.reg.SDE} below gives regularity in the spatial coordinate for $(\sF^N)^{1/2}$ and $(\swww^N)^{1/2}$, which will later be passed to the limit $N\to\infty$ in Corollary~\ref{c:reg.in.limit}. 
Towards this end we first record an elementary fact:

\begin{lem}\label{l:bump.sqrt.lipschitz}
Let $f:\R\to[0,\infty)$ be a smooth function which integrates to one and has compact support. Then the derivative of $f^{1/2}$ is uniformly bounded, that is,
	\[ \|(f^{1/2})'\|_\infty
	= \frac12 \bigg\|\frac{f'}{f^{1/2}}\bigg\|_\infty
	\le O(1)\,,\]
where $O(1)$ indicates a constant depending on $f$.

\begin{proof}
Since $f$ is smooth, we have $\|f''\|_\infty \le C$ where $C$ is a constant depending only on $f$. Assume without loss that $C$ is large. Now suppose $f(y)>0$ is small while $|f'|$ is large, say
$f'(y) \ge 2 C f(y)^{1/2}$. The bound on $f''$ implies $f'(x) \ge C f(y)^{1/2}$ for all $|x-y| \le f(y)^{1/2}$. It follows that
	\[
	f\Big(y-f(y)^{1/2}\Big)
	\le f(y) - Cf(y)^{1/2} \cdot f(y)^{1/2}
	< 0\,,
	\]
contradicting the assumption that $f\ge0$. The case $f'(y) \le - 2 C f(y)^{1/2}$ leads to a similar contradiction, so we conclude 
	\[\bigg\|\frac{f'}{f^{1/2}}\bigg\|_\infty
	\le 2C\,,
	\]
which proves the claim.
\end{proof}
\end{lem}

\begin{lem}[regularity for SDE coefficients]
\label{l:sp.smooth.reg.SDE}
For $t\geq q_0$ we have the bounds 
	\[
    \bigg\|\frac{\partial}{\partial y}
    	\Big[(\sF^N)^{1/2}\Big]\bigg\|_\infty
    +\bigg\|\frac{\partial}{\partial y}
    	\Big[(\swww^N)^{1/2}\Big]\bigg\|_\infty
	\le 
	\bigg(\frac{\trK L}{\sep}\bigg)^{O(1)} \MAX_N\,,
	\]
with $\MAX_N$ as in \eqref{e:MAX.N} and Proposition~\ref{p:Y.kolmogorov}. 

\begin{proof}
We again simplify notation by suppressing the dependence on $N$ inside the proof. We only consider $\sF$, since $\swww$ is handled identically.
Similarly as in the proof of Lemma~\ref{l:sp.smooth.reg}, we express
	\[
	\sF(t,y)
	=\frac{\overline{\varrho_{F,t}*\sphi}(y)}
		{\overline{\varrho_t*\sphi}(y)}
	\equiv \frac{\Num(t,y)}{\Den(t,y)}\,.
	\]
(Note the abuse of notation, in that we now use $\Num$ to refer to the numerator of $\sF$ rather than the numerator of $\sB$.)
Recall from
\eqref{e:denominator.lbd} and the surrounding discussion that $\Den(t,y)\ge 1/\trK^{O(1)}$ for all $|y|\le 5\trK$, while 
$\Num(t,y)=0$ for $|y|\ge3\trK+1$. Therefore, similarly to the derivation of \eqref{e:sB.quotient.rule}, we have
	\beq\label{e:sF.quotient.rule}
	\bigg|\frac{\partial}{\partial y} \sF(t,y)\bigg|
	\le \trK^{O(1)}
	\bigg| 
	\frac{\partial_y\Num(t,y)}{\Den(t,y)^{1/2}}
	- \sF(t,y) \partial_y \Den(t,y)
	\bigg|\,,\eeq
and we aim to show that $\partial_y\sF/\sF^{1/2}$ is bounded. To this end, we calculate
	\[\bigg|\frac{\partial}{\partial y}
	\Num(t,y)\bigg|
	=\bigg|\frac{\partial}{\partial y}
	\overline{\varrho_{F,t}*\sphi}(y)\bigg|
	= \bigg| \int_{\R^3} F(t,\vx)
		\overline{\sphi}'
		(y-(\vec{1},\vx))
		\,\varrho_t(d\vx)\bigg|\,.\]
Recalling the explicit definition \eqref{e:spatial.smoothing.product} of $\sphi$ from Definition~\ref{d:sp.smooth}, we have 
	\begin{align*}
	\bigg|
	\frac{d}{dy}\overline{\sphi}(y)
	\bigg|
	&=\bigg|
	\lim_{\delta\to0}
	\int_{H(y)}
	\frac{\sphi(\vy+(\delta,0,0))
		-\sphi(\vy)}{\delta}
	\,d\vy\bigg|\\
	&\le  \frac{1}{\sep^4}\int_{H(y)}
	\bigg| (\sphicirc)'\bigg(
		\frac{y_1}{\sep}
		\bigg) \bigg|
	\sphicirc\bigg(\frac{y_2}{\sep}
		\bigg)
	\sphicirc\bigg(\frac{y_3}{\sep}
		\bigg)
	\,d\vy \\
	&\le 
	\frac{O(1)}{\sep^4}
	\int_{H(y)} 
	\sphicirc\bigg(
		\frac{y_1}{\sep}
		\bigg)^{1/2}
	\sphicirc\bigg(\frac{y_2}{\sep}
		\bigg)
	\sphicirc\bigg(\frac{y_3}{\sep}
		\bigg) 
	\,d\vy\\
	&\le
	\frac{O(1)}{\sep^4}\bigg\{
	\int_{H(y)} \sphicirc\bigg(\frac{\vy}{\sep}\bigg)\,d\vy\bigg\}^{1/2}
	\bigg\{
		\int_{H(y)} \boldsymbol{1}\{\|\vy\|_\infty\le\sep\}\,d\vy
		\bigg\}^{1/2}
	\le \frac{O(1)}{\sep^{3/2}} \overline{\sphi}(y)^{1/2}
	\end{align*}
where the third line above follows from Lemma~\ref{l:bump.sqrt.lipschitz} applied to the function $\sphicirc:\R\to[0,\infty)$, and the first inequality in the last line above is by the Cauchy--Schwarz inequality. We then apply the Cauchy--Schwarz inequality to bound 
	\begin{align*}
	\bigg|\frac{\partial}{\partial y}
	\Num(t,y)\bigg|
	&=\bigg|\frac{\partial}{\partial y}
	\overline{\varrho_{F,t}*\sphi}(y)\bigg|
	= \bigg| \int_{\R^3} F(t,\vx)
		\overline{\sphi}' (y-(\vec{1},\vx))
		\,\varrho_t(d\vx)\bigg|\\
	&\le \frac{1}{\sep^{O(1)}}
		 \bigg| \int_{\R^3} F(t,\vx)
		 \overline{\sphi} (y-(\vec{1},\vx))^{1/2}
		\,\varrho_t(d\vx)\bigg|\\
	&\le\frac{1}{\sep^{O(1)}}
	\bigg\{
	\int_{\R^3} F(t,\vx)\,\varrho_t(d\vx)
	\bigg\}^{1/2}
	\bigg\{  \int_{\R^3} F(t,\vx)
		 \overline{\sphi} (y-(\vec{1},\vx))
		\,\varrho_t(d\vx)
	\bigg\}^{1/2} \\
	&=
	\frac{1}{\sep^{O(1)}}
	\bigg\{
	\int_{\R^3} F(t,\vx)\,\varrho_t(d\vx)
	\bigg\}^{1/2} \Num(t,y)^{1/2}\,,
	\end{align*}
where the last equality is by recalling the definition of $\Num$. Recalling the definition of $F$ from \eqref{e:sde.F.N}, we can argue similarly as 
in the proof of Lemma~\ref{l:sp.smooth.reg}
(using the quenched bounds of Lemma~\ref{l:tQ.first.quenched}) to bound 
	\[\bigg\{\int_{\R^3} F(t,\vx)\,\varrho_t(d\vx)\bigg\}^{1/2}
	\le  L^{O(1)}
	\MAX_N\,.
	\]
Altogether this bounds the first term in \eqref{e:sF.quotient.rule}, giving
	\[
	\bigg|\frac{\partial}{\partial y}\Num(t,y)\bigg|
	\le \frac{L^{O(1)}\MAX_N}
		{\sep^{O(1)}}
	\cdot \Num(t,y)^{1/2}
	= 
	\frac{L^{O(1)}\MAX_N}
		{\sep^{O(1)}}
		\Den(t,y)^{1/2} \sF(t,y)^{1/2}
	\,.
	\]
For the second term in \eqref{e:sF.quotient.rule}, we note that Lemma~\ref{l:sp.smooth.reg} implies
	\[
	\sF(t,y)
	\le \bigg(\frac{\trK L}{\sep}\bigg)^{O(1)}\MAX_N
	 \cdot \sF(t,y)^{1/2}\,.
	\]
We also note that
	\[
	\bigg|\frac{\partial}{\partial y}\Den(t,y)\bigg|
	=\bigg| \int_{\R^3} \overline{\sphi}'(y-(\vec{1},\vx))\varrho_t(d\vx)\bigg|
	\le \frac{1}{\sep^{O(1)}}\,.
	\]
Combining the above bounds and substituting into \eqref{e:sF.quotient.rule} gives
	\[
	\bigg\|\frac{\partial}{\partial y}(\sF^{1/2})\bigg\|_\infty
	= \bigg\|
	\frac{\partial_y\sF(t,y)}{\sF(t,y)^{1/2}}
	\bigg\|_\infty
	\le
	\bigg(\frac{\trK L}
		{\sep}\bigg)^{O(1)}\MAX_N\,,
	\]
as claimed.
\end{proof}
\end{lem}

\begin{lem}[continuity in time for spatially smoothed occupation densities]
\label{l:sp.smooth.cty}
Recall from Definition~\ref{d:discrete.occ} that $\rho_{\vY,N}$ denotes the occupation density of $\vY^{\frozen,N}$. Abbreviate $\rho_{N,s}$ for the density conditional on time $s$, and recall that we factorize
$\rho_{N,s}=\varrho_{N,s}
\otimes \varsigma_{N,s}$. With high probability, for $t\geq q_0$ the spatially smoothed densities satisfy the bounds
	\begin{align*}
	\Big\|\overline{\sphi*\varrho_{N,s}}
		-\overline{\sphi*\varrho_{N,t}}
		\Big\|_\infty 
	&\le
	\bigg(\frac{L}{\sep}\bigg)^{O(1)}
	\TMAX_N
	(t-s)^{1/2}\,,\\
	\Big\|
	\sphi*\varsigma_{N,s}
	-\sphi*\varsigma_{N,t}
	\Big\|_\infty
	&\le
	\bigg(\frac{L}{\sep}\bigg)^{O(1)}
	\TMAX_N
	(t-s)^{1/2}\,.
	\end{align*}
for all $0\le s\le t\le 1$, with $\TMAX_N$ as in \eqref{e:def.TMAX}. 

\begin{proof}
We again simplify notation by suppressing the dependence on $N$ inside the proof. From the definition of the spatially smoothed densities, we have
	\begin{align*}
	&\Big| \overline{\sphi*\varrho_s}(y)
		-\overline{\sphi*\varrho_t}(y)\Big|
	=\bigg|\int_{H(y)}
	\bE\Big[ \sphi(\vy-\vY^{\frozen,\RomI:\RomIII}(s)) -
		\sphi(\vy-\vY^{\frozen,\RomI:\RomIII}(t)) \Big] \,d\vy
	\bigg| \\
	&= \bigg|\bE\Big[ \overline{\sphi}
		(y-(\vec{1},\vY^{\frozen,\RomI:\RomIII}(s)))
	- \overline{\sphi}
		(y-(\vec{1},\vY^{\frozen,\RomI:\RomIII}(t)))\Big]\bigg|\\
	&\le
	O(1) \|\overline{\sphi}'\|_\infty
	\bE\bigg[ \Big\|\vY^{\frozen,\RomI:\RomIII}
		(s)-\vY^{\frozen,\RomI:\RomIII}(t)
		\Big\|^2\bigg]^{1/2}
	\le \bigg(\frac{L}{\sep}\bigg)^{O(1)} 
	\TMAX_N\,,
	\end{align*}
where the last step uses the bound
\eqref{e:Y.I.Kolmogorov} from Proposition~\ref{p:Y.kolmogorov} for increments of $Y^{\frozen,\RomI}$,
together with the bound from Proposition~\ref{p:Y.notI.second.quenched} for increments of $Y^{\frozen,\RomII}$ and $Y^{\frozen,\RomIII}$.
A similar argument applies with $\varsigma$ in place of $\varrho$, using the bound from
Proposition~\ref{p:Y.notI.second.quenched} for increments of $Y^{\frozen,\Ising}$.
\end{proof}
\end{lem}

Up to now, we have studied the effect of ``ideal smoothing'' --- in the sense of Definition~\ref{d:intro.SDE.coefs} --- in the spatial coordinate only. \textbf{We next turn to temporal smoothing. Here, rather than the ``ideal smoothing,'' we use a different smoothing rule, Definition~\ref{d:t.smooth}, which allows us to pass the constraints from Lemmas \ref{l:sp.smooth.domain} and \ref{l:sp.smooth.budget} more easily to the smoothed quantities.} We will later need to address the discrepancy between ``ideal smoothing'' and
Definition~\ref{d:t.smooth}; see Lemma~\ref{l:ideal.vs.actual.coeffs}. We recall from Assumption~\ref{a:params} that $\tep \ll \sep$ --- that is, the bandwidth for temporal smoothing is much smaller than the bandwidth for spatial smoothing.

\begin{dfn}[temporal smoothing]\label{d:t.smooth}
Recall the spatially smoothed quantities specified by Definition~\ref{d:sp.smooth}. Recall the function $\sphicirc$ from Definition~\ref{d:sp.smooth}, and let
	\[\tphi(t)
	\equiv
	\frac{2}{\tep}
	\sphicirc\lt(\frac{2
		t}{\tep}-1\rt)\,,\] so that $\tphi$ integrates to one and has support $[0,\tep]$. We define temporally smoothed coefficients by simply convolving with $\tphi$ in the time coordinate:
	\[\tB^N(t,y)
	\equiv\int_{\R} 
        \tphi(s) 
        \sB^N(t-s,y)\,ds
	= \tE\Big[ \sB^N(t-\sss,y)\Big]	\,,\]
where $\tE$ refers to expectation over the law of the random variable $\sss$ sampled from density $\tphi$. We define analogously $\tV^N$, $\tU^N$, $\trrr^N$,and $\twww^N$.  Recalling \eqref{e:qbar.and.qstar}, we will view all these quantities as being defined only for $t\in [q_*,1]$, so its definition only involves the values $\sbbb^N(t-s)$ for $t-s\in [\bar q,1]$, thus avoiding any temporal boundary issues. We then define 
	\begin{align*}
	\tbbb^N(t,y)
	&\equiv\frac{\tB^N(t,y)}{(p^N)'(t)^{1/2} }\,,\\
	\tvvv^N(t,y)
	&\equiv \frac{\tV^N(t,y)}{p^N(t) } 
	\,,\\
	\tuuu^N(t,y)
	&\equiv
	\frac{\tU^N(t,y)}{[ p^N(t)
		\cdot t (p^N)'(t)]^{1/2}} \,.
	\end{align*}
We abbreviate $\tGam^N\equiv(\tbbb,\tvvv,\tuuu)^N$. Recalling \eqref{e:sp.smoothed.F.N} we let
	\beq\label{e:temp.smoothed.F.N}
	\tF^N(t,y)
\equiv\Big\{ \tV^N+2\tU^N+\trrr^N
\Big\}(t,y)+2\eta^{1/2}\,,
	\eeq
and we note that $\tF^N$ can be obtained by convolving $\sF^N$ with $\tphi$ in the time coordinate. In particular, since $\sF^N\ge0$, it follows that $\tF^N\ge0$ also.
\end{dfn}

\begin{dfn}[variants of temporal smoothing]
\label{d:t.smooth.variant}
For comparison with Definition~\ref{d:t.smooth}, we define
$\TGam^N\equiv(\Tbbb,\Tvvv,\Tuuu)^N$ by directly convolving $\sGam^N$ with $\tphi$:
	\[\Tbbb^N(t,y)
	\equiv
	\int_{\R} 
        \tphi(s) 
        \sbbb^N(t-s,y)\,ds
	= \tE\Big[ \sbbb^N(t-\sss,y)\Big]
	= \tE\bigg[
	\frac{\sB^N(t-\sss,y)}{(p^N)'(t-\sss)}
	\bigg]\,,
	\]
and analogously $\Tvvv^N$ and $\Tuuu^N$. We abbreviate $\TGam^N\equiv(\Tbbb,\Tvvv,\Tuuu)^N$.
Recalling Definition~\ref{d:intro.SDE.coefs}, we also define the ``ideal'' smoothed coefficients 
	\[\iB^N(t,y)
	\equiv\frac{ \displaystyle
	\int_{\R} \tphi(t-s) \overline{\varrho_{D,N,s}*\sphi}(y)\,dy}
		{\displaystyle
		 \int_{\R} \tphi(t-s) \overline{\varrho_{N,s}*\sphi}(y)\,dy}\,,
	\]
and similarly $\iV$, $\iU$, $\ir$. We similarly define
	\[\iw^N(t,y)
	\equiv
	\frac{\displaystyle
	\int_{\R} \tphi(t-s) \varsigma_{\www,N,s}*\sphi(y)\,dy
	}{\displaystyle
	\int_{\R} \tphi(t-s) \varsigma_{N,s}*\sphi(y)\,dy
	}\,.
	\]
Recalling \eqref{e:sde.F.N}, \eqref{e:sp.smoothed.F.N}, and \eqref{e:temp.smoothed.F.N}, we let 
	\beq\label{e:ideal.smoothed.F.N}
	\iF^N(t,y)\equiv
	\Big\{ \iV^N+\iU^N+\ir^N\Big\}(t,y)
	+2\eta^{1/2}\,.\eeq
Note that $\iF^N$ is the ideal smoothing of $F^N$ from \eqref{e:sde.F.N}, and since $F^N\ge0$ we also have $\iF^N\ge0$.
\end{dfn}

\begin{cor}[regularity for temporally smoothed coefficients]
\label{c:t.smooth.regularity.N} 
Recall $\MAX_N$ from \eqref{e:MAX.N} and Proposition~\ref{p:Y.kolmogorov}.
The coefficients from Definitions~\ref{d:t.smooth} and \ref{d:t.smooth.variant} satisfy the following bounds:
\begin{enumerate}[(a)]
\item\label{c:t.smooth.regularity.N.a}
Similarly to Lemma~\ref{l:sp.smooth.reg}, the coefficients are continuously differentiable in space and time, and satisfy
	\begin{align}\nonumber
	&\Big(
	\|\tB^N\|_\infty
    +
    \|\partial_y\tB^N\|_\infty
 + \tep \|\partial_t \tB^N\|_\infty\Big)^2\\
	&\qquad+ \max\bigg\{
    \|\bar{A}\|_\infty
    +\|\partial_y\bar{A}\|_\infty
    + \tep\|\partial_t\bar{A}\|_\infty
    : \bar{A}\in\{ \tV,\tU,\trrr,
	\twww\}^N
    \bigg\} 
	\le 
	\bigg(\frac{\trK L}{\sep}\bigg)^{O(1)} 
	(\MAX_N)^2\,.
	\label{e:t.smooth.regularity.N}
	\end{align}
The analogous bound holds for the coefficients $\{\iB,\iV,\iU,\ir,\iw\}^N$.

\item\label{c:t.smooth.regularity.N.b}
Similarly to Lemma~\ref{l:sp.smooth.reg.SDE}, we have 
\beq\label{e:t.smooth.regularity.N.b}
\max\bigg\{
	\bigg\|\frac{\partial}{\partial y}
    \Big[\bar{A}^{1/2}\Big]\bigg\|_\infty
    +\tep
    \bigg\|\frac{\partial}{\partial t}
    \Big[\bar{A}^{1/2}\Big]\bigg\|_\infty
	:\bar{A}\in\{\tF,\twww\}^N\bigg\}
	\le 
	\bigg(\frac{\trK L}
		{\sep}
		\bigg)^{O(1)} \MAX_N\,.
	\eeq
The analogous bound holds for $\iF^N$ and $\iw^N$.
\end{enumerate} 

\begin{proof} 
\eqref{c:t.smooth.regularity.N.a}
Recall that the temporally smoothed coefficients from Definition~\ref{d:t.smooth} are obtained by taking the spatially smoothed coefficients from Definition~\ref{d:sp.smooth} and simply convolving with the kernel $\tphi$ in the time coordinate, so they satisfy the same spatial regularity bounds as in Lemma~\ref{l:sp.smooth.reg}. As for regularity in the time coordinate, we can calculate for example 
	\begin{align*}
    \bigg|
	\frac{\partial}{\partial t}
	\tB^N(t,y)
    \bigg|
	&= \bigg(\frac{2}{\tep}\bigg)^2
    \bigg|
	\int_{\R}
	(\sphicirc)'\bigg(\frac{2(t-s)}{\tep}-1\bigg)
	\sB^N(s,y)\,ds
    \bigg|\\
	&= \frac{2}{\tep}
    \bigg|
	\int_{\R}
	(\sphicirc)'(u-1)
	\sB^N\bigg(t-\frac{\tep u}{2},y\bigg)\,du
    \bigg|
	\le \frac{2\|(\sphicirc)'\|_1}{\tep} \|\sB^N\|_\infty
	\le \frac{O(1)}{\tep}\|\sB^N\|_\infty\,,
	\end{align*}
where the last term is bounded by Lemma~\ref{l:sp.smooth.reg}. This implies \eqref{e:t.smooth.regularity.N}. To obtain the analogous bound for the ``ideal'' smoothed coefficients, recall the notation $\Num(t,y)$ and $\Den(t,y)$ from the proof of Lemma~\ref{l:sp.smooth.reg}, and note that for example we can express
	\[
	\iB^N(t,y)
	=\frac{\displaystyle \int_{\R} \tphi(t-s)\Num(s,y)\,ds}
	{\displaystyle \int_{\R} \tphi(t-s)\Den(s,y)\,ds}
	\equiv \frac{\bar{N}(t,y)}{\bar{D}(t,y)}\,.
	\]
It is easy to conclude that $\iB^N$ satisfies similar spatial regularity bounds as $\sB^N$. Lastly, we have
	\[
	\frac{\partial}{\partial t} \iB^N(t,y)
	= \frac{\partial_t \bar{N}(t,y)}{\bar{D}(t,y)} 
	- \iB(t,y) \frac{\partial_t \bar{D}(t,y)}{\bar{D}(t,y)}\,,
	\]
and from this it is straightforward to conclude that $\iB^N$ satisfies similar temporal regularity bounds as $\tB^N$. This proves \eqref{c:t.smooth.regularity.N.a}.

\eqref{c:t.smooth.regularity.N.b} Recall the proof of Lemma~\ref{l:sp.smooth.reg.SDE}, where we decomposed $\sF=\Num/\Den$. We now again abuse notation by letting $\Num$ stand for the numerator of $\sF$ rather than of $\sB$. Then
	\[\iF(t,y)
	=\frac{ \displaystyle\int_{\R}\tphi(t-s)\Num(s,y)\,ds}
		{ \displaystyle\int_{\R}\tphi(t-s)\Den(s,y)\,ds}
	\equiv \frac{\bar{N}(t,y)}{\bar{D}(t,y)}\,.
	\]
The bound on the spatial derivative of $\iF^{1/2}$ follows by essentially the same argument as in 
Lemma~\ref{l:sp.smooth.reg.SDE}, so we consider only the time derivative.
Following the proof of Lemma~\ref{l:sp.smooth.reg.SDE}, 
first recall from \eqref{e:denominator.lbd} that $\Den(t,y)\ge 1/\trK^{O(1)}$ for all $|y|\le 5\trK$. Then, similarly to
\eqref{e:sF.quotient.rule}, we have
	\[\bigg|\frac{\partial}{\partial t} \iF(t,y)\bigg|
	\le \trK^{O(1)}
	\bigg|
	\frac{\partial}{\partial t}\bar{N}(t,y)
	- \iF(t,y) 
	\frac{\partial}{\partial t} \bar{D}(t,y)
	\bigg|\,,\]
and we wish to show that $\partial_t\iF/\iF^{1/2}$ is bounded. To this end, recalling the explicit definition of $\tphi$ from Definition~\ref{d:t.smooth},
we bound (using Lemma~\ref{l:bump.sqrt.lipschitz})
	\[ \bigg|\frac{d}{dt} \tphi(t)\bigg|
	= \bigg|\bigg(\frac{2}{\tep}\bigg)^2
	(\sphicirc)'\bigg(
	\frac{2t}{\tep}-1
	\bigg) \bigg|
	\le \frac{O(1)}{\tep^2}
	\sphicirc
    \bigg(
	\frac{2t}{\tep}-1
	\bigg)^{1/2}
	= \frac{O(1)}{\tep^{3/2}}
	\tphi(t)^{1/2}\,.
	\]
It follows from this that
	\begin{align*}
	\bigg|\frac{\partial}{\partial t}\bar{N}(t,y)\bigg|
	&= \int_{\R} (\tphi)'(t-s) \Num(s,y)\,ds
	\le \frac{O(1)}{\tep^{3/2}}
	\int_{\R} \tphi(t-s)^{1/2}\Num(s,y)\,ds
	\\
	&\le \frac{O(1)}{\tep^{3/2}}
	\bigg\{\int_{\R} \tphi(t-s)\Num(s,y)\,ds\bigg\}^{1/2}
	\bigg\{\int_{\R} 
	\ind\{|t-s|\le\tep\}
	\Num(s,y)\,ds\bigg\}^{1/2} \\
	&\le
	\frac{O(1)}{\tep}
	 \bigg(\frac{\trK L}{\sep}\bigg)^{O(1)}\MAX_N
	\cdot \bar{N}(t,y)^{1/2}\,,
	\end{align*}
where the last inequality uses that $\Num \le \trK^{O(1)} \sF$ is bounded by part \eqref{c:t.smooth.regularity.N.a}. Now the rest of the argument of Lemma~\ref{l:sp.smooth.reg.SDE} goes through essentially without change, yielding the claimed bound for the time derivative of $\iF^{1/2}$.
As for $\tF^{1/2}$, the bound on the spatial derivative is again straightforward, so we consider only the time derivative. Similarly to the above calculation we can bound
	\begin{align*}
	\bigg|
	\frac{\partial}{\partial t}
	\tF(t,y)\bigg|
	&\le
	\frac{O(1)}{\tep^{3/2}}
	\bigg\{\int_{\R} \tphi(t-s)
		\sF(s,y)\,ds\bigg\}^{1/2}
	\bigg\{\int_{\R} 
	\ind\{|t-s|\le\tep\}
	\sF(s,y)\,ds\bigg\}^{1/2}\\
	&\le
	\frac{O(1)}{\tep}
	\bigg(\frac{\trK L}{\sep}\bigg)^{O(1)}\MAX_N
	\cdot \tF(t,y)^{1/2}\,,
	\end{align*}
which implies \eqref{e:t.smooth.regularity.N.b}.
This concludes the proof of \eqref{c:t.smooth.regularity.N.b}.
\end{proof}
\end{cor}

For any concave function $p:[q_0,1]\to[0,1]$, let
\beq\label{e:bad.times.J}
J=J(p,\tep,\lambda)
\equiv 
\Big\{q_0\leq t\leq 1~:~ p'(t+\tep)-p'(t-\tep) < -\lambda \Big\}.
\eeq
The next lemmas show that the domain and budget constraints from Lemmas \ref{l:sp.smooth.domain} and \ref{l:sp.smooth.budget} pass through the temporal smoothing rule of Definition~\ref{d:t.smooth}, outside a set of ``bad'' times, as defined by \eqref{e:bad.times.J} with $p^N$ in place of $p$. We will show in Lemma~\ref{l:bad.times.JN.J} below that the set of bad times for $p^N$ can be controlled by an analogous set $J$ of bad times for $p$, and we will show in Lemma~\ref{l:bad.J} that $J$ is a small subset of $[0,1]$. We will show in Lemma~\ref{l:t.smooth.domain} that the different smoothing rules of Definitions~\ref{d:t.smooth} and \ref{d:t.smooth.variant} behave similarly outside the set of bad times, and we will pass this result to the limit in Corollary~\ref{c:limit.coefs.variants.of.temp.smoothing}. Likewise we will show in Lemma~\ref{l:t.smooth.budget} that the budget constraint
from Lemma~\ref{l:sp.smooth.budget} passes to the temporally smoothed quantities outside the set of bad times, and we will pass this result to the limit in Corollary~\ref{c:limit.budget}.

\begin{lem}[temporal smoothing of domain constraint]
\label{l:t.smooth.domain}
Let $\tvvv^N$ and $\tuuu^N$ be the temporally smoothed quantities from Definition~\ref{d:t.smooth}. Let $J^N\equiv J(p^N,\tep,\lambda)$ be as defined by \eqref{e:bad.times.J}. Then, analogously to 
\eqref{e:domain.plus.coarse} and the bound from Lemma~\ref{l:sp.smooth.domain}, with high probability we have
	\[
	\Tvvv^N(t,y)
	\ge
	\Tuuu^N(t,y)^2
	-\frac{\eta^2}{2}\]
for all $(t,y)\in[q_*,1]\times\R$. We have additionally
	\[\max\left\{
	\begin{array}{c}
	|\tbbb^N(t,y)^2
		-\Tbbb^N(t,y)^2|,\\
	|\tvvv^N(t,y)-\Tvvv^N(t,y)|,\\
	|\tuuu^N(t,y)^2-\Tuuu^N(t,y)^2|
	\end{array}
	\right\}
	\le \bar{\lambda}_N
	\equiv \bigg(\frac{ \trK L}{\sep}\bigg)^{O(1)} 
	(\MAX_N)^2 \lambda\,.
	\]
for all $(t,y)\in ([q_*,1]  \setminus J^N)\times \bbR$.

\begin{proof}
As discussed in Definition~\ref{d:t.smooth}, 
the coefficients $\tuuu^N$, $\tvvv^N$,  $\Tuuu$, $\Tvvv$, etc.\ are well-defined for all times $t\in[q_*,1]$. The first conclusion follows from Lemma~\ref{l:sp.smooth.domain}: for all $t\in[q_*,1]$,
    \begin{align*}
    \Tvvv^N(t,y)
    &=\tE\bigg[ \frac{\sV^N(t-\sss,y)}
        {p^N(t-\sss)}\bigg]
    \ge
    \tE\bigg[ 
    \frac{\sU^N(t-\sss,y)^2}
        {p^N(t-\sss)
        \cdot (t-\sss)
        (p^N)'(t-\sss)
        }
        \bigg]
         - \frac{\eta^2}{2}\\
    &= \tE[ \suuu^N(t-\sss,y)^2]
    - \frac{\eta^2}{2}
    \ge (\tE[\suuu^N(t-\sss,y)])^2
    - \frac{\eta^2}{2}
    = \Tuuu^N(t,y)^2 -\frac{\eta^2}{2}
    \,.
    \end{align*} 
From Definition~\ref{d:t.smooth}, we have 
	\[\tvvv^N(t,y)
	= \frac{\tV^N(t,y)}{p^N(t)}
	= \frac{\tE[\sV^N(t-\sss,y)]}{p^N(t)}\]
Since $(p^N)'\le L^2$ by \eqref{e:sdelimits.pprime},
the difference between
$p^N(t)$ and $p^N(t-\sss)$ is at most
$L^2\tep$. Moreover, since we restrict to $t\ge q_*$ as defined by \eqref{e:qbar.and.qstar}, it follows by recalling \eqref{e:breve.q.bounds} that 
\beq\label{e:sde.lbd.on.p.N.assump}
p^N(t)\ge \brep=\frac1{L^{\breve{c}}}\,.\eeq
Combining with the bound on $\sV^N$ from Lemma~\ref{l:sp.smooth.reg} gives, for $t\in [q_*,1]$, 
	\[
	\Big|\tvvv^N(t,y)
	-\Tvvv^N(t,y)\Big|
	=\bigg|
	\frac{\tE[\sV^N(t-\sss,y)]}
		{p^N(t)}
	-\tE\bigg[ \frac{\sV^N(t-\sss,y)}
		{p^N(t-\sss)}\bigg]\bigg|
	\le \bigg(\frac{ \trK L}{\sep}\bigg)^{O(1)} 
	(\MAX_N)^2 \tep\,.\]
Now applying Lemma~\ref{l:sp.smooth.domain} gives, for all $t\in[q_*,1]$,
	\[
	\frac{\sU^N(t,y)^2}
		{p^N(t)
		\cdot t
		(p^N)'(t)
		} - \frac{\eta^2}{2}
	\le \frac{\sV^N(t,y)}
		{p^N(t)}
	\le \bigg(\frac{ \trK L}{\sep}\bigg)^{O(1)} 
	(\MAX_N)^2\,,
	\]
which also implies a bound on $\sU^N$. 
Finally we note that if $t\in [q_*,1]  \setminus J^N$, then
	\begin{align*}
	\Big|\tuuu^N(t,y)
	-\Tuuu^N(t,y)\Big|
	&=\bigg|
	\frac{\tE[\sU^N(t-\sss,y)]}
		{p^N(t) \cdot t (p^N)'(t)}
	-\tE\bigg[\frac{\sU^N(t-\sss,y)}
		{p^N(t-\sss) \cdot
		(t-\sss) (p^N)'(t-\sss)}
		\bigg]
	\bigg| \\
	&\le
	\bigg(\frac{ \trK L}{\sep}\bigg)^{O(1)} 
	(\MAX_N) \lambda\,,
	\end{align*}
using the bound on $\sU^N$ noted above. It follows that
for $t\in [q_*,1] \setminus J^N$ we have
	\[
	\Big|\tuuu^N(t,y)^2
	-\Tuuu^N(t,y)^2\Big|
	\le
	\Big|\tuuu^N(t,y)
	-\Tuuu^N(t,y)\Big|
	\Big|\tuuu^N(t,y)^2
	+\Tuuu^N(t,y)^2\Big|
	\le \bigg(\frac{ \trK L}{\sep}\bigg)^{O(1)} 
	(\MAX_N)^2 \lambda\,.
	\]	
The discrepancy between $\tbbb^N$ and $\Tbbb^N$ can be bounded similarly.
Combining the above bounds proves the claim, recalling the order of parameters from Assumption~\ref{a:params}.
\end{proof}
\end{lem}

With $\bar{\lambda}_N$ as in the statement of Lemma~\ref{l:t.smooth.domain}, we let
	\beq\label{e:def.bar.lambda}
	\bar{\lambda}\equiv
	\limsup_{N\to\infty}\bar{\lambda}_N
	\le \bigg(\frac{ \trK L}{\sep}\bigg)^{O(1)} \lambda\,.
	\eeq
Then, similarly to \eqref{e:cost} and \eqref{e:cost.t.eps}, we now define
	\[
	\Cost_{2\bar{\lambda}}(b,u,v)
	\equiv b^2 + u^2 + 
		\bigg( 
        \Big[(v-u^2 
        + 2\bar{\lambda})_+\Big]
        ^{1/2}-1\bigg)^2\,.
	\]
We then define temporally smoothed ``cost'' and ``budget'' quantities
	\begin{align}
	\label{e:cost.temp.smoothed}
	\bar{C}^N(t)
	&\equiv\int_{\R} 
	\Cost_{2\bar{\lambda}}\Big(\tGam^N(t,y)\Big)
	\overline{\varrho_{\vY,N,t}*\sphi(y)}
	\,dy\,,\\
	\bar{W}^N(t)
	&\equiv
	\int_{\R} \twww^N(t,y)^{1/2} 
		\varsigma_{\vY,N,t}*\sphi(y)\,dy\,.
		\label{e:budget.temp.smoothed}
	\end{align}
We also let $\breve{C}^N$ be defined analogously to $\bar{C}^N$, but with $\TGam^N$ in place of $\tGam^N$.

\begin{lem}[temporal smoothing of budget constraint]
\label{l:t.smooth.budget}
Analogously to \eqref{e:budget.constraint.Ising.plus.coarse} and the bound from Lemma~\ref{l:sp.smooth.budget}, with high probability we have 
	\begin{align*}
	\breve{C}^N(t)
	- \frac{\bar{W}^N(t)^2}{\alpha} 
	\le
	\bigg(\frac{\trK L}{\sep}\bigg)^{O(1)} 
	(\MAX_N)^2 \TMAX_N \tep^{1/2}
	+ 3\epsilon^\circ
	\end{align*}
for all $t\in[q_*,1]$. Additionally, with $J^N\equiv J(p^N,\tep,\lambda)$ as  defined by \eqref{e:bad.times.J}, we have
	\[
	\Big|\breve{C}^N(t)
	-\bar{C}^N(t)
	\Big|
	\le \bigg(\frac{ \trK L}{\sep}\bigg)^{O(1)} 
	(\MAX_N)^2
	\lambda^{1/2}
	\]
for all $t\in[q_*,1] \setminus J^N$.

\begin{proof}
We continue to use the abbreviations introduced above. The function $\Cost_{2\bar{\lambda}}$ is convex on the expanded domain $\{(b,v,u) : v \ge u^2-2\bar{\lambda}\}$. Recall $\tGam^N\equiv(\tbbb,\tvvv,\tuuu)^N$
and $\TGam^N\equiv(\Tbbb,\Tvvv,\Tuuu)^N$: we have from Lemma~\ref{l:t.smooth.domain} that
(with high probability)
 $\TGam^N(t,y)$ lies in the expanded domain for all $(t,y)\in[q_*,1]\times \R$, while
$\tGam^N(t,y)$ lies in the expanded domain for all $(t,y)\in ([q_*,1]\setminus J^N)\times \R$. Lemma~\ref{l:t.smooth.domain} implies, for $(t,y)\in ([q_*,1]\setminus J^N)\times \R$,
	\[
	\bigg|
	\Cost_{2\bar{\lambda}}
	\Big( \tGam^N(t,y)\Big)
	-\Cost_{2\bar{\lambda}}
	\Big( \TGam^N(t,y)\Big)
	\bigg|
	\le
	\bigg(\frac{ \trK L}{\sep}\bigg)^{O(1)} 
	(\MAX_N)^2
	\lambda^{1/2}\,.
	\]
Integrating with respect to $\overline{\varrho_{\vY,N,t}*\sphi(y)}\,dy$ proves the second conclusion. For all $(t,y)\in[q_*,1]\times\R$, we can bound
	\begin{align*}
	& \Cost_{2\bar{\lambda}}
	\Big(\TGam^N(t,y)\Big)
	=
	\Cost_{2\bar{\lambda}}\bigg( \tE \Big[\sGam(t-\sss,y)\Big]\bigg) \\
	&\qquad\le \tE\bigg[
	\Cost_{2\bar{\lambda}}\Big( \sGam(t-\sss,y)\Big)\bigg]
	= \int_{\R} \tphi(t-s)
	\Cost_{2\bar{\lambda}}( \sGam(s,y) )\,ds\,.
	\end{align*}
Integrating with respect to $\overline{\varrho_t*\sphi}(y)\,dy$ gives
	\begin{align*}
	&\breve{C}^N(t)
	=
	\int_{\R} \Cost_{2\bar{\lambda}}\Big(\TGam^N(t,y)\Big)
	\overline{\varrho_t*\sphi}(y)\,dy
	\le
	\int_{\R}
	\bigg[
	 \int_{\R} \tphi(t-s)
	\Cost_{2\bar{\lambda}}( \sGam(s,y) )\,ds\bigg]
	\overline{\varrho_t*\sphi}(y)\,dy \\
	&\qquad=
	\int_{\R} \tphi(t-s)
	\bigg[\int_{\R}
	\Cost_{2\bar{\lambda}}( \sGam(s,y) )
	\overline{\varrho_s*\sphi(y)} \,dy\bigg]\,ds
	+O\bigg[ \bigg(\frac{\trK L}{\sep}\bigg)^{O(1)} 
	(\MAX_N)^2 \TMAX_N
	\tep^{1/2}\bigg]\\
	&\qquad=
	\int_{\R} \tphi(t-s)
	\tilde{C}^N(s)\,ds
	+O\bigg[ \bigg(\frac{\trK L}{\sep}\bigg)^{O(1)} 
	(\MAX_N)^2 \TMAX_N
	\tep^{1/2}\bigg]
	\,,
	\end{align*}
where the error estimate uses the bounds on $\sGam=(\sbbb,\svvv,\suuu)$ proved in Lemma~\ref{l:sp.smooth.reg}, along with the regularity of the map $t\mapsto \overline{\varrho_s*\sphi}$ proved in Lemma~\ref{l:sp.smooth.cty}. 
 Similarly, note that
	\[
	S_t(w)
	\equiv
	\bigg( \int_{\R} w(y)^{1/2} \varsigma_t*\sphi(y)\,dy\bigg)^2
	\]
is a concave function. It follows that
	\begin{align*}
	&\bar{W}^N(t)^2 
	=\bigg(\int_{\R} \twww(t,y)^{1/2} \varsigma_t*\sphi(y)\,dy\bigg)^2
	= S_t\bigg(\tE\Big[ \swww(t-\sss,\cdot)\Big]\bigg)
	\ge \tE\bigg[ S_t\Big( \swww(t-\sss,\cdot)\Big)\bigg]\\
	&\qquad= \int_{\R} \tphi(t-s)
	\bigg( \int_{\R} \swww(s,y)^{1/2} 
		\varsigma_t*\sphi(y)\,dy\bigg)^2\,ds\\
	&\qquad=
	\int_{\R} \tphi(t-s)
	\bigg( \int_{\R} \swww(s,y)^{1/2} 
		\varsigma_s*\sphi(y)\,dy\bigg)^2
	\,ds
	+
	O\bigg[ \bigg(\frac{\trK L}{\sep}\bigg)^{O(1)} 
	(\MAX_N)^2 \TMAX_N \tep^{1/2}\bigg]
	\\
	&\qquad=
	\int_{\R} \tphi(t-s)
	\tilde{W}^N(s)^2\,ds
	+
	O\bigg[ \bigg(\frac{\trK L}{\sep}\bigg)^{O(1)} 
	(\MAX_N)^2 \TMAX_N \tep^{1/2}\bigg]
	\end{align*}
where the error estimate again follows from Lemmas \ref{l:sp.smooth.reg} and \ref{l:sp.smooth.cty}. The claim then follows by combining with
Lemma~\ref{l:sp.smooth.budget}.
\end{proof}
\end{lem}

\subsection{Limiting coefficients and constraints}\label{ss:limit.coefs}

In the previous subsection, we showed that the domain and budget constraints
 from Proposition~\ref{p:budget-constraints-coarsened} behave well with respect to spatial and temporal smoothing. \textbf{In the remainder of this subsection, we write ``$N\to\infty$'' to refer to the limit $N\to\infty$, $\eta\to0$, and $\delta\to0$, as indicated by Assumption~\ref{a:params}. In the current subsection we will use the smoothing to show that the above constraints can be passed through this limit.}

We first review the continuous-time occupation measures: these were previously introduced in Definition~\ref{d:intro.cts.occ} for the abstract setting of \S\ref{ss:sde.abstract}, and they are the continuous-time analogues of the occupation measures of Definition~\ref{d:discrete.occ}:

\begin{dfn}[occupation measures, continuous-time]
\label{d:cts.occ} Let $\vY^\frozen$ be the subsequential limit from Theorem~\ref{t:tightness}. Similarly to Definitions~\ref{ss:sde.abstract}
and \ref{d:discrete.occ},
let $\rho_{\vY}$ be the occupation measure for $\vY^\frozen=(\vY^{\frozen,\RomI:\RomIII},Y^{\frozen,\Ising})$. We also let
$\rho_D$, $\rho_V$, $\rho_U$, $\rho_\rrr$, $\rho_\www$
be the occupation measures for
$D^{\frozen,\RomI}$, $Q^{\frozen,\RomII}$, $Q^{\frozen,\RomII,\RomIII}$, $Q^{\frozen,\RomIII}$, and $Q^{\frozen,\Ising}$ respectively: for instance, $\rho_D$ is the signed measure such that
	\[\int_{\R^5}
	f(t,\vy)\rho_D(dt\,d\vy)
	= \bE\int_0^1
	f(t,\vY^\frozen,(t))
	\,dD^\frozen(t)\]
for any bounded continuous function $f$ on $\R^5$. Note it 
follows from Theorem~\ref{t:tightness} that $D^\frozen$, $Q^{\frozen,\RomII}$, $Q^{\frozen,\RomII,\RomIII}$, $Q^{\frozen,\RomIII}$, and $Q^{\frozen,\Ising}$ are finite-variation processes, so the integrals with respect to these processes are well-defined. Denote the conditional measures $\rho_t\equiv \rho_{\vY,t}$, $\rho_{D,t}$, etc. Note that since $(Y^{\frozen,\RomI},Y^{\frozen,\RomII},Y^{\frozen,\RomIII})$ evolves independently of $Y^{\frozen,\Ising}$,
we can factorize
 $\rho_t=\varrho_t\otimes\varsigma_t$,
 $\rho_{D,t}=\varrho_{D,t}\otimes\varsigma_t$, and so on.
\end{dfn}

We next define continuous-time analogues of the smoothed coefficients from Definitions \ref{d:sp.smooth} and \ref{d:t.smooth}.

\begin{dfn}[continuous-time smoothed coefficients]
\label{d:cts.coeffs}
Analogously to Definition~\ref{d:sp.smooth}, we define the spatially smoothed continuous-time coefficients
	\[\srrr(t,y)\equiv
	\frac{\overline{\varrho_{\rrr,t}*\sphi}(y)}
		{\overline{\varrho_t*\sphi}(y)}\,,\quad
	\swww(t,y)\equiv
	\frac{\varsigma_{\www,t}*\sphi(y)}
		{\varsigma_t*\sphi(y)}\,,\]
as well as
	\begin{align*}
	\sB(t,y)
	&\equiv
	\frac{\overline{\varrho_{D,t}*\sphi}(y)}
		{\overline{\varrho_t*\sphi}(y)}
	\equiv
	p'(t)^{1/2} \cdot \sbbb(t,y)\,,\\
	\sV(t,y)
	&\equiv
	\frac{\overline{\varrho_{V,t}*\sphi}(y)}
		{\overline{\varrho_t*\sphi}(y)}
	\equiv
	p(t) \cdot \svvv(t,y)\,,\\
	\sU(t,y)
	&\equiv
	\frac{\overline{\varrho_{U,t}*\sphi}(y)}
		{\overline{\varrho_t*\sphi}(y)}
	\equiv
	\Big[ p(t) \cdot t p'(t) \Big]^{1/2}
		\cdot \suuu(t,y)\,.
	\end{align*}
Analogously to Definition~\ref{d:t.smooth}, we define the 
temporally smoothed continuous-time coefficients
	\beq
	\label{eq:tB-def}
	\tB(t,y)
	\equiv
	\int_{\R} \tphi(s) \sB(t-s,y)\,ds
	= \tE\Big[ \sB(t-\sss,y)\Big]
	\eeq
and similarly $\tV$, $\tU$, $\trrr$, $\twww$. We define also the scaled coefficients
	\begin{align*}
	\tbbb(t,y)
	&\equiv\frac{\tB(t,y)}
		{p'(t)^{1/2} }\,,\\
	\tvvv(t,y)
	&\equiv \frac{\tV(t,y)}{p(t) } 
	\,,\\
	\tuuu(t,y)
	&\equiv
	\frac{\tU(t,y)}{[ p(t)
		\cdot t p'(t)]^{1/2}} \,.
	\end{align*}
and abbreviate $\tGam\equiv(\tbbb,\tvvv,\tuuu)$. Likewise, analogously to Definition~\ref{d:t.smooth.variant}, we define $\TGam\equiv(\Tbbb,\Tvvv,\Tuuu)$ as well as $(\iB,\iV,\iU,\ir,\iw)$.
We then let
	\begin{align*}
	\sF(t,y)
	&\equiv\Big\{ \sV+2\sU+\srrr\Big\}(t,y)\,,\\
	\tF(t,y)
	&\equiv\Big\{ \tV+2\tU+\trrr\Big\}(t,y)\,,\\
	\iF(t,y)
	&\equiv\Big\{ \iV+2\iU+\ir\Big\}(t,y)\,,
	\end{align*}
analogously to \eqref{e:sp.smoothed.F.N}, \eqref{e:temp.smoothed.F.N}, and \eqref{e:ideal.smoothed.F.N}.
\end{dfn}

The next lemma confirms that the measures
from Definition~\ref{d:cts.occ}  indeed correspond to the $N\to\infty$
weak limits of the measures from Definition~\ref{d:discrete.occ}. Recall the objects $\Q_{\bullet,N}$, $\Q$, $\mu_{\bullet,\bG^N}$, $\mu$ from Theorem~\ref{t:tightness}.

\begin{lem}[weak convergence of occupation measures]
\label{l:occ.wk.conv}
Suppose $\Q_{\bullet,N}\Rightarrow\Q$. There exists a Skorohod coupling of the disorder matrices $\bG^N$ for which, with probability one over $(\bG^N)_{N\ge1}$, 
we have as $N\to\infty$ the following weak convergence, restricted to the time interval $[q_*,1]$:
\[
(\mu_{\bullet,\bG^N},\rho_{\vY,N},\rho_{D,N},\rho_{V,N},\rho_{U,N},\rho_{\rrr,N},\rho_{\www,N})
\Longrightarrow
(\mu,\rho_{\vY},\rho_D,\rho_V,\rho_U,\rho_\rrr,\rho_\www).
\]
Here the former measures are as in Definition~\ref{d:discrete.occ}, and the latter are as in Definition~\ref{d:cts.occ}.

\begin{proof} 
It follows from Definition~\ref{d:discrete.occ} that for all bounded continuous $f$, we have
	\[
	\int_{\R^4}
	f(t,y) \rho_{\vY,N,t}(d\vy)
	= \bE_{\bG^N} 
	\Big[f(t, \vY^{\frozen,\bullet,N}(t))
	\Big]
	\]
where $\bE_{\bG^N}$ is expectation with respect to $\mu_{\bullet,\bG^N}$, and $Y^{\frozen,\bullet,N}$ is the piecewise constant variant of $Y^{\frozen,N}$.  By Theorem~\ref{t:tightness}, we have a Skorohod coupling along which $\mu_{\bullet,\bG^N}\Rightarrow\mu$. This implies, for each fixed $t$,
	\begin{align*}
	&\lim_{N\to\infty}\int_{\R^4}
		f(t,y) \rho_{\vY,N,t}(d\vy)
	=\lim_{N\to\infty}
	\bE_{\bG^N} f(t,
		\vY^{\frozen,\bullet,N}(t))
	= \bE f(t,\vY^\frozen(t))
	=
	\int_{\R^4}f(t,y) \rho_{Y,t}(d\vy)
	\end{align*}
almost surely, which implies 
$\rho_{\vY,N,t}\Rightarrow\rho_{Y,t}$ for each $t$. This can be integrated over $0\le t\le 1$ using the dominated convergence theorem, so 
	\begin{align*}
	&\lim_{N\to\infty}
	\int_{[0,1]\times\R^4}f(t,\vy) \rho_{\vY,N}(dt\,d\vy)
	-\int_{[0,1]\times\R^4}f(t,\vy) \rho_{\vY}(dt\,d\vy) \\
	&\qquad=
	\lim_{N\to\infty} \int_0^1\bE_{\bG^N} f(t,
		\vY^{^\frozen,\bullet,N}(t))
	-\int_0^1 \bE f(t,\vY^\frozen(t))
	=0\,,
	\end{align*}
which implies $\rho_{\vY,N}\Rightarrow\rho_Y$. 

We next argue that
$\rho_{D,N}\Rightarrow\rho_D$.
Since $\mu_{\bullet,\bG^N}\Rightarrow\mu$, there is a coupling of random paths along which
$\vY^{\frozen,\bullet,N}$ converges uniformly to $\vY^\frozen$, and
 $\bar{D}^{\bullet,\frozen,\RomI,N}$ converges uniformly to $D^\frozen$. Since $D^\frozen$ is a finite-variation process, it defines a (random) signed measure $\nu_D$ on $[0,1]$, namely the measure satisfying
 \[
    D^\frozen(t) = \nu_D([0,t))
 \]
 for all $t\in [0,1]$.  Likewise, $\bar{D}^{\bullet,\frozen,\RomI,N}$ defines a signed measure $\nu_{D,N}$ on $[0,1]$, and the uniform convergence implies $\nu_{D,N}\Rightarrow\nu_D$. Note from \eqref{e:def.rho.b.N} that we can rewrite
 	\[\int_{[0,1]\times\R^4} f(t,\vy) \rho_{D,N}(dt\,d\vy)
	= \bE_{\bG^N} \int_0^1
	f(t,\vY^{\frozen,\bullet,N}(t))
	\,\nu_{D,N}(dt)\,.
	\]
It follows that, for any bounded continuous function $f$, we have
 	\begin{align*}
	&\bigg|
	\int f(t,\vy) \rho_{D,N}(dt\,d\vy)
	-\int f(t,\vy) \rho_D(dt\,d\vy)\bigg| \\
	&\quad=
	\bigg|\bE \int_0^1 
		f(t,\vY^{\frozen,\bullet,N}(t))
	\,\nu_{D,N}(dt)
	-\bE\int_0^1
	f(t,\vY^\frozen(t))
	\,\nu_D(dt)\bigg|
	\le \textup{(a)}+\textup{(b)}
	\end{align*}
where (a) accounts for the difference in the integrand, and (b) accounts for the difference in the measures:
 	\begin{align*}
	\textup{(a)}
	&=
	\bigg|\bE \int_0^1
	\Big[ 
	f(t,\vY^{\frozen,\bullet,N}(t))
	-f(t,\vY^\frozen(t))\Big]
	\,\nu_{D,N}(dt)\bigg|
	\,,\\
	\textup{(b)}
	&=\bigg|
	\bE\int_0^1
	f(t,\vY^\frozen(t))
	\,\nu_{D,N}(dt)
	-\bE\int_0^1
	f(t,\vY^\frozen(t))
	\,\nu_D(dt)\bigg|\,.
	\end{align*}
Recall that $\nu_D$ is a signed measure, so it has a Jordan decomposition
$\nu_D = (\nu_D)_+-(\nu_D)_-$. Denote $|\nu_D|\equiv (\nu_D)_++(\nu_D)_-$; this is the total variation measure associated to $D^\frozen$. Denote analogously $|\nu_{D,N}|$ for the 
total variation measure associated to $\bar{D}^{\frozen,\RomI,\bullet,N}$. We can then bound
	\begin{align*}
	\textup{(a)}
	&\le
	\bE\bigg[
	\sup_{0\le t\le 1}
	\Big| f(t, 
	\vY^{\frozen,\bullet,N}(t))
	-f(t,\vY^\frozen(t))
	\Big|
	\cdot \Big|\nu_{D,N}([0,1])\Big|
	\bigg] \\
	&\le
	\bE\bigg[
	\sup_{0\le t\le 1}
	\Big| f(t, 
	\vY^{\frozen,\bullet,N}(t))
	-f(t,\vY^\frozen(t))
	\Big|^2\bigg]^{1/2}
	\bE \bigg[
	\Big|\nu_{D,N}([0,1])\Big|^2\bigg]^{1/2}\,.
	\end{align*}
The last line above is a product of two terms: we claim that the first term tends to zero while the second stays bounded. Indeed, for the first term,
 the expression inside the expectation tends to zero almost surely by the Skorohod coupling, and is bounded because we assume $f$ is bounded, so the expectation tends to zero by the dominated convergence theorem. For the second term, we note that
 	\[
	\bE \bigg[
	\Big|\nu_{D,N}([0,1])\Big|^2\bigg]
	= \bE\bigg[\bigg(
	\sum_{0\le d\le \dmax-1}
	|\Delta\bar{D}^{\frozen,\RomI,N}(q_d)|\bigg)^2\bigg]
	\le L^{O(1)} (\MAX_N)^2
	\]
by Lemma~\ref{l:barD.I.second.quenched}. This proves that (a) tends to zero as $N\to\infty$. Meanwhile, (b) tends to zero as $N\to\infty$ by the weak convergence $\nu_{D,N}\to\nu_D$. This proves that $\rho_{D,N}\Rightarrow\rho_D$. 
The analogous claims for $V$, $U$, $\rrr$, and $\www$ follow similarly, using Lemma~\ref{l:tQ.first.quenched} in place of Lemma~\ref{l:barD.I.second.quenched}.
This concludes the proof.
\end{proof}
\end{lem}

\begin{lem}[convergence of smoothed coefficients]
\label{l:limit.coefs}
The coefficients 
$(\tB,\tV,\tU,\trrr,\twww)^N$,
$(\Tbbb,\Tvvv,\Tuuu)^N$, and 
$(\iB,\iV,\iU,\ir,\iw)^N$
from Definitions~\ref{d:t.smooth} 
and \ref{d:t.smooth.variant}
converge pointwise as $N\to\infty$ to the coefficients
$(\tB,\tV,\tU,\trrr,\twww)$,
$(\Tbbb,\Tvvv,\Tuuu)$, and 
$(\iB,\iV,\iU,\ir,\iw)$
from Definition~\ref{d:cts.coeffs}. Since $\twww^N,\iw^N,\tF^N,\iF^N$ are all nonnegative, it follows that 
$\twww,\iw,\tF,\iF$ are also nonnegative.  Additionally, analogously  to \eqref{e:domain.plus.coarse} and the bounds from Lemmas \ref{l:sp.smooth.domain} and \ref{l:t.smooth.domain}, we have almost surely
	\[
	\Tvvv(t,y)
	\ge
	\Tuuu(t,y)^2 
	\]
for all 
$(t,y)\in [q_*,1]\times \bbR$. 

\begin{proof} 
From Definitions \ref{d:sp.smooth} and \ref{d:t.smooth},
	\begin{align*}
	\Tuuu^N(t,y)
	&= \int_{\R}\tphi(t-s) \suuu^N(s,y)\,ds
	=
	\int_{\R}
	\frac{\tphi(t-s)}{p^N(s) 
	\cdot s (p^N)'(s)}
	\frac{\overline{\varrho_{U,N,s}*\sphi}(y)}
	{\overline{\varrho_{\vY,N,s}*\sphi}(y)}
	\,ds \\
	&=\int_{H(y)} \int_{\R^4}
	\frac{\tphi(t-s)}{p^N(s) 
	\cdot s (p^N)'(s)}
	\frac{\sphi(\vy-\vx)}
	{\overline{\varrho_{\vY,N,s}*\sphi}(y)}
	\,\rho_{U,N}(ds\,d\vx)\,d\vy\,.
	\end{align*}
We will argue that the above converges pointwise to
	\[\Tuuu(t,y)
	= \int_{H(y)} \int_{\R^4}
	\frac{\tphi(t-s)}{p(s) 
	\cdot s p'(s)}
	\frac{\sphi(\vy-\vx)}
	{\overline{\varrho_{\vY,s}*\sphi}(y)}
	\,\rho_U(ds\,d\vx)\,d\vy\,.
	\]
In the above display for $\Tuuu^N(t,y)$, consider the rightmost expression in the first line. We will first argue that $\Tuuu^N(t,y)$ is close to
	\begin{align*}
	\Tuuu^{N,(1)}(t,y)
	&\equiv
	\int_{\R}
	\frac{\tphi(t-s)}{p(s) \cdot  s p'(s)}
	\frac{\overline{\varrho_{U,N,s}*\sphi}(y)}
	{\overline{\varrho_{\vY,N,s}*\sphi}(y)}
	\,ds\\
	&= \int_{H(y)} \int_{\R^4}
	\frac{\tphi(t-s)}{p(s) \cdot  s p'(s)}
	\frac{\sphi(\vy-\vx)}
	{\overline{\varrho_{\vY,N,s}*\sphi}(y)}
	\,\rho_{U,N}(ds\,d\vx)\,d\vy
	\,.
	\end{align*}
To this end, note it follows from 
Definition~\ref{d:sp.smooth} and Lemma~\ref{l:sp.smooth.reg}  that
	\[
	\|\sU^N\|_\infty
	= \sup_{s\in[q_0,1]} 
	\bigg\|
	\frac{\overline{\varrho_{U,N,s}*\sphi}}
	{\overline{\varrho_{\vY,N,s}*\sphi}}\bigg\|_\infty
	\le
	\bigg(
	\frac{\trK L}{\sep}\bigg)^{O(1)}
	(\MAX_N)^2
	\le
	\bigg(
	\frac{\trK L}{\sep}\bigg)^{O(1)}
	\,,
	\]
where the last bound holds for $N$ large enough by 
the bound \eqref{e.MAX.bound} from Proposition~\ref{p:Y.kolmogorov}.
The functions $p^N$ are concave and converge almost everywhere to $p$, and it is well-known that consequently $(p^N)'$ converges almost everywhere to $p'$. By the definition \eqref{e:qbar.and.qstar} of $\bar q$, both $1/p^N$ and $1/(p^N)'$ are uniformly bounded by $L^2$ on $[\bar q,1]$. For each $t\ge q_*$ and $s$ such that $t-s \in \mathsf{supp}(\tphi)$, we have $s\ge t-\tep \ge \bar q$. So, it follows by the bounded convergence theorem that
	\[\Big|
	\Tuuu^N(t,y)
	-\Tuuu^{N,(1)}(t,y)\Big|
	\le
	\bigg(
	\frac{\trK L}{\sep}\bigg)^{O(1)}
	\int_{\R}
	\bigg|\frac{\tphi(t-s)}
		{p^N(s) \cdot s (p^N)'(s)}
	-\frac{\tphi(t-s)}
		{p(s) \cdot s p'(s)}\bigg| \,ds
	\le o_N(1)\,.
	\]
Thus
$\Tuuu^N(t,y)$ is close to $\Tuuu^{N,(1)}(t,y)$, as claimed. We will next argue that $\Tuuu^{N,(1)}(t,y)$ is close to
	\[
	\Tuuu^{N,(2)}(t,y)
	\equiv
	\int_{H(y)} \int_{\R^4}
	\frac{\tphi(t-s)}{p(s) \cdot  s p'(s)}
	\frac{\sphi(\vy-\vx)}
	{\overline{\varrho_{\vY,s}*\sphi}(y)}
	\,\rho_{U,N}(ds\,d\vx)\,d\vy\,.
	\]
For this it suffices to note that the functions
$\overline{\varrho_{\vY,N,s}*\sphi}$ and 
$\overline{\varrho_{\vY,s}*\sphi}$ are bounded away from zero by the argument from Lemma~\ref{l:sp.smooth.reg} (cf.\ \eqref{e:denominator.lbd}), and are
uniformly close: 
	\begin{align*}
	&\sup_{s,y}\bigg\{ \Big|\overline{\varrho_{\vY,N,s}*\sphi}(y)
	-\overline{\varrho_{\vY,s}*\sphi}(y)\Big|\bigg\}\\
	&\qquad= \sup_{s,y}\bigg|
		\int_{H(y)} \bE\bigg\{ \sphi\Big(\vy-\vY^{\frozen,N}(s)\Big)
		-\sphi\Big(\vy-\vY^\frozen(s)\Big)\bigg\}\,d\vy\bigg| 
		\le o_N(1)\,,
	\end{align*}
where the last bound follows by considering the Skorohod coupling in which $\vY^{\frozen,N}$ converges uniformly to $\vY^\frozen$. It follows that
$\Tuuu^{N,(1)}(t,y)$ and $\Tuuu^{N,(2)}(t,y)$ are close. Finally, it follows by the weak convergence result of Lemma~\ref{l:occ.wk.conv} that
$\Tuuu^{N,(2)}(t,y)$ is close to
$\Tuuu(t,y)$:  this uses the fact that $p'$ is monotone, so it has countably many discontinuities, and the set of discontinuities has measure zero under the limiting measure $\rho_U$. This proves that $\Tuuu^N(t,y)$ converges pointwise to $\Tuuu(t,y)$, and a very similar argument shows that $\tU^N(t,y)$ converges pointwise to $\tU(t,y)$. The other convergence claims follow by very similar arguments.
\end{proof}
\end{lem}

Recall that
$\MAX_N$ is defined by \eqref{e:MAX.N}, and satisfies the bound \eqref{e.MAX.bound}. 
Recall also the random variable 
$\SMAX_N$ from \eqref{e:SMAX.N};
it is formally defined by
Proposition~\ref{p:drift.quenched} and \eqref{e:SMAX.N}.
Lastly recall that $\TMAX_N$ is defined by \eqref{e:def.TMAX}, and satisfies the bound \eqref{e.JMAX.bound}. It follows that
	\[\limsup_{N\to\infty}\MAX_N
	=\limsup_{N\to\infty}\TMAX_N
	\le L^{O(1)}\]
almost surely. (The bound \eqref{e.SMAX.bound} together with Assumption~\ref{a:params} implies that $\limsup\MAX_N$ agrees with $\limsup\TMAX_N$. The bound  \eqref{e.MAX.bound} then implies that the limit must be at most $L^{O(1)}$, almost surely.) The next corollary gives the limiting version of Corollary~\ref{c:t.smooth.regularity.N} and Lemma~\ref{l:sp.smooth.cty}.

\begin{cor}[regularity for limiting coefficients]\label{c:reg.in.limit}
Consider the limiting coefficients from
Definition~\ref{d:cts.coeffs}. Recall also that $\rho_{\vY}$ denotes the occupation measure of the limiting process $\vY(t)$, and let $\rho_s$ be the measure conditional on time $s$.  The bounds from  Corollary~\ref{c:t.smooth.regularity.N} and 
Lemma~\ref{l:sp.smooth.cty} pass to these limiting objects:
\begin{enumerate}[(a)]
\item\label{c:reg.in.limit.a} Similarly to  Corollary~\ref{c:t.smooth.regularity.N}\ref{c:t.smooth.regularity.N.a}, 
we have
	\begin{align*}
	&\Big(
	\|\tB\|_\infty
    +
    \|\partial_y\tB\|_\infty
 + \tep \|\partial_t \tB\|_\infty\Big)^2\\
	&\qquad+ \max\bigg\{
    \|\bar{A}\|_\infty
    +\|\partial_y\bar{A}\|_\infty
    + \tep\|\partial_t\bar{A}\|_\infty
    : \bar{A}\in\{ \tV,\tU,\trrr,
	\twww\}
    \bigg\} 
	\le 
	\bigg(\frac{\trK L}{\sep}\bigg)^{O(1)} \,.
	\end{align*}
The analogous bound holds for $(\iB,\iV,\iU,\ir,\iw)$.
\item \label{c:reg.in.limit.b} Similarly to  Corollary~\ref{c:t.smooth.regularity.N}\ref{c:t.smooth.regularity.N.b}, 
we have 
	\[\max\bigg\{
	\bigg\|\frac{\partial}{\partial y}
    \Big[\bar{A}^{1/2}\Big]\bigg\|_\infty
    +\tep\bigg\|\frac{\partial}{\partial t}
    \Big[\bar{A}^{1/2}\Big]\bigg\|_\infty
	:\bar{A}\in\{\tF,\twww\}\bigg\}
	\le 
	\bigg(\frac{\trK L}
		{\sep}
		\bigg)^{O(1)}\,.
		\]
The analogous bound holds for $\iF$ and $\iw$.
\item \label{c:reg.in.limit.c} Similarly to Lemma~\ref{l:sp.smooth.cty},
for all $0\le s\le t\le1$ we have
	\begin{align*}
	\Big\|\overline{
	\sphi*\varrho_s}-
		\overline{
		\sphi*\varrho_t}\Big\|_\infty
	&\le
	\bigg(\frac{L}{\sep}\bigg)^{O(1)}
	(t-s)^{1/2}\,\\
	\Big\|\sphi*\varsigma_s-\sphi*\varsigma_t\Big\|_\infty
	&\le
	\bigg(\frac{L}{\sep}\bigg)^{O(1)}
	(t-s)^{1/2}\,.
	\end{align*}
\end{enumerate}

\begin{proof}
Corollary~\ref{c:t.smooth.regularity.N}, 
together with the 
convergence result of Lemma~\ref{l:limit.coefs}, implies claims \eqref{c:reg.in.limit.a} and \eqref{c:reg.in.limit.b}. Claim \eqref{c:reg.in.limit.c} follows using Fatou's lemma.
\end{proof}
\end{cor}

\begin{dfn}[L\'evy--Prohorov convergence of derivatives] \label{d:levy.prohorov.derivs}
The function $(p^N)'$ is nonincreasing, so there exists a nonnegative measure $\nu_N$ on $[q_0,1]$ with
	\[
	\nu_N( (s,t]) = -\Big[(p^N)'(t)-(p^N)'(s)\Big]
	\]
for all $q_0 \le s\le t\le1$, and with $\nu_N(\{0\})=0$. By passing to a subsequence, we may assume that $\nu_N$ converges in the L\'evy--Prohorov metric to a limiting measure $\nu$. Since $(p^N)'$ converges to $p'$ almost everywhere, we immediately conclude that 
	\[\nu((s,t])
	=-\Big[ p'(t)-p'(s)\Big]
	\]
for almost all $s,t$. Now recall that $p'$ denotes the derivative of $p$ from the right, so in particular $p'$ is right-continuous: this implies that the above equality holds for all $s,t\in[0,1]$. For this form of convergence, we hereafter say for short that \textbf{$(p^N)'$ converges to $p'$ in the L\'evy--Prohorov metric.}
\end{dfn}

Recall from \eqref{e:bad.times.J} the definition of the set of ``bad'' times $J(p,\sep,\lambda)$. We now show that under the assumption of L\'evy--Prohorov convergence, the set of bad times for $p^N$, with $N$ large, is contained in a set of bad times for $p$.

\begin{lem}\label{l:bad.times.JN.J}
Recall the definition \eqref{e:bad.times.J}. 
If $(p^N)'$ converges to $p'$ in the L\'evy--Prohorov metric in the sense of Definition~\ref{d:levy.prohorov.derivs}, then  we have
	\[
	J^N\equiv
	J(p^N,\tep,\lambda)
		\subseteq 
		J\bigg(p,2\tep,\frac{\lambda}{2}\bigg)
	\equiv J
	\] for all $N$ large enough. 

\begin{proof}
The function $(p^N)'$ is nonincreasing, so there exists a nonnegative measure $\nu_N$ on $[q_0,1]$ with
	\[
	\nu_N( (s,t]) = -\Big[(p^N)'(t)-(p^N)'(s)\Big]
	\]
for all $q_0 \le s\le t\le1$, and with $\nu_N(\{0\})=0$. Define analogously the measure $\nu$ corresponding to the limiting function $p$.
The L\'evy--Prohorov convergence assumption implies that for all $N\ge N(\tep)$ we have
	\[
	\nu_N(A) \le \nu(A^{\tep})+\tep
	\]
where $A\subseteq[q_0,1]$ is any Borel set, and $A^{\tep}$ is the $\tep$-fattening of $A$. (The same inequality also holds with $\nu$ and $\nu_N$ exchanged.) It follows that
	\begin{align*}
	&-\Big[(p^N)'(t+\tep)-(p^N)'(t-\tep)\Big]
	= \nu_N([t-\tep,t+\tep])\\
	&\qquad\le \nu([t-2\tep,t+2\tep]) + \tep
	=-\Big[p'(t+2\tep)-p'(t-2\tep)\Big]+ \tep
	\le \frac{\lambda}{2}+\tep\,,
	\end{align*}
where the last bound holds provided $t\notin J$. Since $\tep\ll\lambda$ by Assumption~\ref{a:params}, this shows that if $t\notin J$ then $t\notin J^N$, as claimed.
\end{proof}
\end{lem}

\begin{cor}\label{c:limit.coefs.variants.of.temp.smoothing}
Let $J\equiv J(p,2\tep,\lambda/2)$ as  defined by \eqref{e:bad.times.J}. Then the coefficients from Definition~\ref{d:cts.coeffs} satisfy the estimates 
	\[
	\max\left\{
	\begin{array}{c}
	|\tbbb(t,y)^2
		-\Tbbb(t,y)^2|,\\
	|\tvvv(t,y)-\Tvvv(t,y)|,\\
	|\tuuu(t,y)^2-\Tuuu(t,y)^2|
	\end{array}
	\right\}
	\le \bigg(\frac{\trK L}{\sep}\bigg)^{O(1)} \lambda
	\]
for all $(t,y)\in( [q_*,1]  \setminus J)\times\R$.

\begin{proof}
Lemma~\ref{l:limit.coefs} gives
	\[\begin{aligned}
	\tB(t,y)&=\tB^N(t,y)+o_N(1)\,,\\
	\Tbbb(t,y)&=\Tbbb^N(t,y)+o_N(1)
	\end{aligned}\]
for all $(t,y)\in[q_*,1]\times\R$. Recall that $(p^N)'(t)$ converges to $p'(t)$ for almost every $t$. Therefore, for \emph{every} $t$, we can find $s\in[t-\tep,t+\tep]$ such that $(p^N)'(s)$ converges to $p'(s)$. If $t\notin J$, then $t\notin J^N$ for sufficiently large $N$ by Lemma~\ref{l:bad.times.JN.J}, so $p'(s)$ differs from $p'(t)$ by at most $\lambda/2$, while $(p^N)'(s)$ differs from $(p^N)'(t)$ by at most $\lambda$. Therefore, for $t\in[q_*,1]\setminus J$, we have
	\[
	\Big|\tbbb(t,y)-\tbbb^N(t,y)\Big|
	=\bigg|\frac{\tB(t,y)}{p'(t)^{1/2}}
	-\frac{\tB^N(t,y)}{(p^N)'(t)^{1/2}}\bigg|
	\le o_N(1)
	+ \bigg( \frac{\trK L}{\sep}\bigg)^{O(1)}(\MAX_N)^2
	\lambda\,,
	\]
using the bounds on $\tB^N$ and $\tB$ implied by Lemma~\ref{l:sp.smooth.reg} and Corollary~\ref{c:reg.in.limit}.
Additionally, since $t\notin J^N$,
Lemma~\ref{l:t.smooth.domain} bounds the discrepancy between $\tbbb^N(t,y)$ and $\Tbbb^N(t,y)$.
Sending $N\to\infty$ gives the claimed bound on the discrepancy between $\tbbb(t,y)$ and $\Tbbb(t,y)$. The other bounds follow by similar arguments:
	\begin{align*}
	\Big|\tvvv(t,y)-\tvvv^N(t,y)\Big|
	&=\bigg|
	\frac{\tV(t,y)}{p(t)}
	-\frac{\tV^N(t,y)}{p^N(t)}\bigg|\,,\\
	\Big|\tuuu(t,y)-\tuuu^N(t,y)\Big|
	&=\bigg|\frac{\tU(t,y)}{[ p(t)
		\cdot t p'(t)]^{1/2}}
		-\frac{\tU^N(t,y)}{[ p^N(t)
		\cdot t (p^N)'(t)]^{1/2}}
		\bigg|\,.
	\end{align*}
Note in these cases that in the factor $L^{O(1)}$ appearing in the final bound, the $O(1)$ exponent can depend on the constant $\breve{c}$ from \eqref{e:breve.q.bounds}. We will remark on this point again in the proof of Theorem~\ref{thm:Lipschitz-SDE-approx}.
\end{proof}
\end{cor}

Analogously to \eqref{e:cost.temp.smoothed} and \eqref{e:budget.temp.smoothed}, let us define
\begin{align}
	\label{e:cost.temp.smoothed.LIMIT}
	\bar{C}(t)
	&\equiv\int_{\R} 
	\Cost_{2\bar{\lambda}}\Big(\tGam(t,y)\Big)
	\overline{\varrho_{\vY,t}*\sphi(y)}
	\,dy\,,\\
	\bar{W}(t)
	&\equiv
	\int_{\R} \twww^N(t,y)^{1/2} 
		{\varsigma_{\vY,t}*\sphi(y)}\,dy\,.
		\label{e:budget.temp.smoothed.LIMIT}
	\end{align}
We also let $\breve{C}(t)$ be defined analogously to $\bar{C}(t)$, but with $\TGam$ in place of $\tGam$.

\begin{cor}[limiting budget constraints]
\label{c:limit.budget}
Analogously  to \eqref{e:budget.constraint.Ising.plus.coarse} and the bounds from Lemmas \ref{l:sp.smooth.budget} and \ref{l:t.smooth.budget}, we have
almost surely 
	\begin{align*}
	&\breve{C}(t)
	- \frac{\bar{W}(t)^2}{\alpha} 
	\le
	\bigg(\frac{\trK L}{\sep}\bigg)^{O(1)}  
	\tep^{1/2}
	+ 3\epsilon^\circ
	\end{align*}
for all $t\in[q_*,1]$. Additionally, with $J\equiv J(p,2\tep,\lambda/2)$ as  defined by \eqref{e:bad.times.J}, we have
	\[
	\Big|\bar{C}(t)-\breve{C}(t)\Big|
	\le \bigg(\frac{\trK L}{\sep}\bigg)^{O(1)} \lambda^{1/2}
	\]
for all $t\in[q_*,1]\setminus J$. 

\begin{proof}
Recall $\TGam^N\equiv(\Tbbb,\Tvvv,\Tuuu)^N$
and $\tGam\equiv(\tbbb,\tvvv,\tuuu)$.
We showed in Lemma~\ref{l:limit.coefs} above that
the coefficients 
$(\Tbbb,\tB,\Tvvv,\tV,\Tuuu,\tU,\trrr,\twww)^N$
from Definition~\ref{d:t.smooth} 
converge pointwise to the coefficients
$(\Tbbb,\tB,\Tvvv,\tV,\Tuuu,\tU,\trrr,\twww)$
from Definition~\ref{d:cts.coeffs}.  Since the coefficients are defined on a compact domain and are uniformly Lipschitz by  Corollary~\ref{c:t.smooth.regularity.N}, the convergence in fact holds uniformly in $y$ for any fixed $t$. Let us abbreviate $\rho_{\vY,N,t} \equiv \rho_{N,t} \equiv \varrho_{N,t}\otimes \varsigma_{N,t}$ and $\rho_{Y,t} \equiv \rho_t \equiv \varrho_t\otimes \varsigma_t$. Then
	\begin{align*}
	\Big|\breve{C}^N(t)
	-\breve{C}(t)\Big|
	&\le \bigg|\int_{\R} \Cost_{2\bar{\lambda}}
	\Big(\TGam^N(t,y)\Big)
	\overline{\varrho_{N,t}*\sphi}(y)
	\,dy
	-\int_{\R} \Cost_{2\bar{\lambda}}
	\Big(\TGam(t,y)\Big)
	\overline{\varrho_{N,t}*\sphi}
	(y)\,dy\bigg| \\
	&\qquad+\bigg|\int_{\R} \Cost_{2\bar{\lambda}}
	\Big(\TGam(t,y)\Big)
	\overline{\varrho_{N,t}*\sphi}
	(y)\,dy
	-\int_{\R} \Cost_{2\bar{\lambda}}
	\Big(\TGam(t,y)\Big)
	\overline{\varrho_t*\sphi}(y)
	\,dy\bigg|\,.
	\end{align*}
On the right-hand side, the first term tends to zero by the convergence result from Lemma~\ref{l:limit.coefs}, making use of the above observation that the convergence holds uniformly over $y$.
For the second term, it follows from Corollary~\ref{c:reg.in.limit}\ref{c:reg.in.limit.a} that
$\TGam$ is a continuous function. It is compactly supported, from which we conclude that $\Cost_{2\bar{\lambda}}\circ\TGam$ is a bounded continuous function. The second term then tends to zero by the weak convergence result from Lemma~\ref{l:occ.wk.conv}. A similar analysis applies to the difference between $\bar{W}(t)^2$ and 
$\bar{W}^N(t)^2$. The first assertion then follows from the corresponding assertion from Lemma~\ref{l:t.smooth.budget}. The bound on the discrepancy between $\bar{C}(t)$ and $\breve{C}(t)$ follows from Corollary~\ref{c:limit.coefs.variants.of.temp.smoothing}.
\end{proof}
\end{cor}

\subsection{Limiting SDE}
\label{ss:limit.sde}

In this subsection, we first apply the abstract result from \S\ref{ss:sde.abstract} to show that the limiting process $\vY^\frozen$ is well approximated by an SDE solution $\hvX^\frozen$ with ``ideal'' coefficients.
We then show that $\hvX^\frozen$ in turn is well approximated by an SDE solution $\vX^\frozen$ with the smoothed coefficients defined in \S\ref{ss:limit.coefs}. This yields the main result of this subsection, Proposition~\ref{p:sde.approx}. To define $\hvX^\frozen$,
let $B$ and $B^\Ising$ be independent standard one-dimensional Brownian motions, and consider the SDE
    \beq\label{e:iX.sde}
	\begin{aligned}
	d\iX^\frozen(t) 
    &= 
    \iB(t,\iX^\frozen(t))\,dt
		+ \iF(t,\iX^\frozen(t))^{1/2}\,dB(t)\,,
    \\
    d\iX^{\frozen,\Ising}(t)
    &=
    \iw(t,\iX^{\frozen,\Ising}(t))^{1/2} dB^{\Ising}(t).
    \end{aligned}
	\eeq
We then have the following:

\begin{ppn}\label{p:fokker.planck}
Take independent random variables $\sss\sim\tphi$ and $\xxx\sim\sphi$.
Let $\vY^\frozen$ be a subsequential limit of the $\vY^{\frozen,N}$ processes, and let $\hvY^\frozen(t)\equiv \vY^\frozen(t-\sss)+\xxx$. Take the coordinate sum
	\beq\label{e:iY.plus}
	\iY^\frozen\equiv
	\sum_{\sigma\in\{\RomI,
		\RomII,\RomIII\}}
	\iY^{\frozen,\sigma}\,.\eeq
Let $\iX^\frozen$ and $\iX^{\frozen,\Ising}$ be the unique solutions of the SDEs \eqref{e:iX.sde} such that  $\iX^\frozen(q_*)$ is equidistributed as $\iY^\frozen(q_*)$, and 
$\iX^{\frozen,\Ising}(q_*)$ is equidistributed as $\iY^{\frozen,\Ising}(q_*)$. Then the (one-dimensional) processes 
$\iX^\frozen$ and $\iY^\frozen$
satisfy the same Fokker--Planck equations on the time interval $[q_*,1]$.
As a consequence, for each fixed $t$ the marginal law of $\iX^\frozen(t)$ agrees with the marginal law of $\iY^\frozen(t)$.
The same statements hold for $\iX^{\frozen,\Ising}$ and $\iY^{\frozen,\Ising}$.

\begin{proof}
This follows from Lemma~\ref{l:projected.fokker.planck},  using the regularity estimates from Corollary~\ref{c:reg.in.limit}. 
\end{proof}
\end{ppn}

We next state and prove an elementary lemma which helps us to control the large jumps of the monotone function $p'$:

\begin{lem}\label{l:bad.J}
Fix $\tep\ll\lambda\ll 1$, and let $J=J(p,2\tep,\lambda/2)$ as defined by \eqref{e:bad.times.J}. Then there exists a set $\oJ$ such that
$J\subseteq \oJ\subseteq [q_0,1]$,
and $\oJ$ can be expressed as the (not necessarily disjoint) union of at most $L^{O(1)}/\lambda$ open intervals each of width $10\tep$. 
In particular, $\oJ$ has one-dimensional Lebesgue measure
	\[\mu_{\textup{Leb}}(\oJ)
	\leq \frac{L^{O(1)}\tep}{\lambda}
	\ll 1\,.\]

\begin{proof}
Subdivide $[q_0,1]$ into intervals of width $2\tep$. (If $(1-q_0)/(2\tep)$ is not integer-valued, then simply let the last interval have width less than $2\tep$.) For the purposes of this proof let us say an interval is ``bad'' if $p'$ decreases  by more than $\lambda/8$ over the interval; otherwise we say the interval is ``good.'' Since $p'$ is monotone and assumed bounded by $L^{O(1)}$, the total number of bad intervals is at most $L^{O(1)}/\lambda$. Let $B\subseteq[0,1]$ be the union of all the bad intervals, and let $\oJ$ be the $(4\tep)$-neighborhood of $B$. We then claim that $J\subseteq\oJ$: to see this, note that if $x\in \oJ^c$, then the interval $[x-2\tep,x+2\tep]$ is disjoint from the bad intervals $B$, and is covered by at most three good intervals. It follows that $p'$ decreases by at most $3\lambda/8$ over the interval $[x-2\tep,x+2\tep]$, which means $x\in J^c$. This implies the claim.
\end{proof}
\end{lem}

\begin{lem}\label{l:ideal.vs.actual.coeffs}
The smoothed coefficients and ``ideal'' coefficients from
Definition~\ref{d:cts.coeffs} are uniformly close to one another: 
	\[\sup\bigg\{
	\Big|\bar{B}(t,y)-\hat{B}(t,y)\Big|
	: (t,y) \in [q_*,1] \times\R \bigg\}
	\le
	\bigg(\frac{\trK L}{\sep}\bigg)^{O(1)} \tep^{1/2}\,.
	\]
The same bound holds for 
$\bar{V}-\hat{V}$,
$\bar{U}-\hat{U}$, and 
$\bar{\rrr}-\hat{\rrr}$.
\end{lem}

\begin{proof}
The different coefficients are handled in the same way, so we give the proof only for the case of $\bar{B}$ and $\iB$. From Definition~\ref{d:cts.coeffs}, $\bar{B}$ is given by
	\[
	\bar{B}(t,y)
	=\int_{\R}\tphi(t-s)
	\frac{\overline{\varrho_{D,s}*\sphi}(y)}
	{\overline{\varrho_s*\sphi}(y)}
	\,ds\,.
	\]
Meanwhile, from Definition~\ref{d:cts.coeffs}, $\iB$ is given by
	\[
	\iB(t,y)
	= \frac{\displaystyle
	\int_{\R}\tphi(t-s)
	\overline{\varrho_{D,s}*\sphi}(y)\,ds}
	{\displaystyle
	\int_{\R}\tphi(t-s)
	\overline{\varrho_s*\sphi}(y)\,ds}\,.
	\]
 Recall from Definition~\ref{d:t.smooth} that the function
$\tphi$ integrates to one and has support $[0,\tep]$. Corollary~\ref{c:reg.in.limit}\ref{c:reg.in.limit.c}
and the bound \eqref{e:denominator.lbd} together imply 
	\[
	\sup\bigg\{
	\bigg|
	\frac{\overline{\varrho_s*\sphi}(y)}
	{\overline{\varrho_{s'}*\sphi}(y)}
	-1
	\bigg|
	: s,s' \in [t-\tep, t\bigg\}
	\leq 
	\bigg(\frac{\trK L}{\sep}\bigg)^{O(1)}
	 \tep^{1/2}
	\]
Comparing $\tB$ 
and $\iB$ with this bound gives
	\[
	\bigg|
	\frac{\tB 
	(t,y)}{\iB(t,y)}
	-1
	\bigg|
	\leq 
	\bigg(\frac{\trK L}{\sep}\bigg)^{O(1)}
	\tep^{1/2}\,.
	\]
Finally we use Corollary~\ref{c:reg.in.limit}\ref{c:reg.in.limit.a}  to bound the maximum value of $\bar{B}$, so the above multiplicative error bound implies the additive error bound claimed in the statement of this lemma.
\end{proof}

\begin{lem}\label{l:gronwall}
Let $X^1$ and $X^2$ be solutions to SDEs
    \[
    dX^i(t)
    =
    b^i(t,X^i(t))\,dt
    +
    \sigma^i(t,X^i(t))\,dB_t\,,
    \]
such that $X^1(q_*)$ and $X^2(q_*)$ are equidistributed.
Suppose for some finite constant $\oL$ that the SDE coefficients $b^i$ and $\sigma^i$ are uniformly bounded by $\oL$, and moreover are $\oL$-Lipschitz in the space coordinate. Lastly suppose
	\[
    \sup_{(t,x)\in [q_*
    ,1]\times \mathbb R}
    \bigg\{
    \Big|b^1(t,x) - b^2(t,x)\Big|
    +
 	\Big|\sigma^1(t,x)-\sigma^2(t,x)\Big|
    \bigg\}
    \le \epsilon\,.
    \]
Then the marginal laws of $X^1$ and $X^2$ are close in $L^2$-Wasserstein distance for all time: more precisely, there is a coupling of $X^1$ with $X^2$ under which
	\[
	\sup\bigg\{
	\E\bigg[\Big(X^1(t)-X^2(t)\Big)^2\bigg] 
	: t\in[q_*
        ,1]
	\bigg\} \le 
	C(\oL) \epsilon^2\,,
	\]
where $C(\oL)$ denotes a constant depending only on $\oL$.

\begin{proof}
Couple $X^1$ and $X^2$ with the same Brownian motion $B$, and same initial value (almost surely). For this coupling let $Z(t) = X^1(t)-X^2(t)$, and $E(t) = \E[Z(t)^2]$. By It\=o's formula, for any
$0\le t\le t'\le1$, we have
	\[
	E(t')-E(t)
	=\E \int_t^{t'} \bigg\{
		2 Z(s)
		\Big(b^1(s,X^1(s))-b^2(s,X^2(s))
		\Big)
	+\Big( \sigma^1(s,X^1(s))
		-\sigma^2(s,X^2(s))
		\Big)^2
		\bigg\}\,ds\,.
	\]
For $s\in[q_*,1]$, using that $b^1$ and $b^2$ are $\epsilon$-close to one another, and $\oL$-Lipschitz in the space coordinate, gives
	\[
	\Big|b^1(s,X^1(s))-b^2(s,X^2(s))
	\Big|
	\le \oL\Big| X^1(s)-X^2(s)\Big|
		+ \epsilon
	=\oL |Z(s)|
		+ \epsilon\,.
	\]
A similar bound holds with $\sigma^i$ in place of $b^i$. It follows that 
for $[t,t']\subseteq[q_*,1]$, we have
	\[
	E(t')-E(t)
	\le 3
	\E \int_t^{t'}
	\Big(\oL |Z(s)| + \epsilon\Big)^2\,ds
	\le 6 \bigg\{
	\oL^2 \int_t^{t'} E(s)\,ds
	+ \epsilon^2 (t'-t)
	\bigg\}\,.
	\]
The claimed bound follows from Gr\"onwall's inequality. More precisely, we can use the form stated in \cite[Lem.~6.1]{MR1071170} which says that if for all $t\in[q_*
    ,1]$ we have
	\[
	E(t) \le a(t) +\int_{q_0}^t b(s) E(s)\,ds
	\]
with $a,b$ continuous nonnegative functions, then it holds that
	\[
	E(1)
	\leq 
	E(q_*)
	+
	\int_{q_*
        }^1
	a(s)b(s)
	\exp
	\lt(
	\int_s^1 b(r)\,dr
	\rt)\,ds.
	\]
In our setting $E(q_*)=0$, $b(t)=6\oL^2$, and $a(t)=6\epsilon^2(t-q_*)$. Therefore
$a(1) \le 6\epsilon^2$, and the claimed bound follows.
\end{proof}
\end{lem}

Let $\vY^\frozen=(Y^{\frozen,\RomI},Y^{\frozen,\RomII},Y^{\frozen,\RomIII},Y^{\frozen,\Ising})$ be the limiting process given by Theorem~\ref{t:tightness}, and recall from \eqref{e:Y.plus} that
	\[
	Y^\frozen
	\equiv \sum_{\sigma\in\{
		\RomI,\RomII,\RomIII\}}
		Y^{\frozen,\sigma}
	\,.\]
In Definition~\ref{d:sp.smooth} we introduced a spatial smoothing kernel $\sphi$ with bandwidth $\sep$, and in Definition~\ref{d:t.smooth} we introduced a temporal smoothing kernel $\tphi$ with bandwidth $\tep$. Let $\hvY^\frozen(t)\equiv \vY^\frozen(t-\sss)+\xxx$ where $\xxx$ is sampled from $\sphi$ while $\sss$ is sampled from $\tphi$; this is what we called the ``ideal'' smoothing.
Recall that $\iX^\frozen$ and $\iX^{\frozen,\Ising}$
are the solutions of the SDE \eqref{e:iX.sde}, with coefficients $\iB$, $\iF^{1/2}$, and $\iw^{1/2}$ from Definition~\ref{d:cts.coeffs}.
 Similarly, let $X^\frozen$ and $X^{\frozen,\Ising}$ be the solution of the SDE 
	\beq\label{e:actual.sde}
	\begin{aligned}
	dX^\frozen(t)
	&= \tB(t,X^\frozen(t)) \,dt 
	+ \tF(t,X^\frozen(t))^{1/2} \,dB(t)\,,\\
	dX^{\frozen,\Ising}(t)
	&= \twww^{1/2}(t,
		X^{\frozen,\Ising}(t))\,dB^\Ising(t)\,,
	\end{aligned}
	\eeq
with coefficients $\tB$, $\tF^{1/2}$,
and $\twww^{1/2}$ from Definition~\ref{d:cts.coeffs}.
We then have the following:

\begin{ppn}\label{p:sde.approx}
Let $X^\frozen$ and $X^{\frozen,\Ising}$ be the unique solutions of the SDEs \eqref{e:actual.sde} 
such that   $X^\frozen(q_*)$ is equidistributed as $Y^\frozen(q_*)$, and 
$X^{\frozen,\Ising}(q_*)$ is equidistributed as $Y^{\frozen,\Ising}(q_*)$. The processes $X^\frozen$ and $Y^\frozen$ have marginals which are close in $L^2$-Wasserstein distance: for each fixed $t$, there exists a coupling of $X^\frozen(t)$ with $Y^\frozen(t)$ under which we have the bound
	\[
	\bE\bigg[\Big( X^\frozen(t)-Y^\frozen(t)\Big)^2\bigg]
	\le O(\sep)+
	C\bigg( \bigg(\frac{\trK L}{\sep}\bigg)^{O(1)} 
	\bigg) \tep^{1/2}\,.
	\]
Likewise with the processes $Y^{\frozen,\Ising}$ and $X^{\frozen,\Ising}$.
 
\begin{proof}
First, $\vY^\frozen$ is close to 
$\hvY^\frozen(t)= \vY^\frozen(t-\sss)+\xxx$ because $\sss$ and $\xxx$ are both small: $\xxx$ is sampled from the spatial kernel $\sphi$ (Definition~\ref{d:sp.smooth}), so $\|\xxx\|_\infty \le\sep$. Meanwhile $\sss$ is sampled from the temporal kernel $\tphi$ (Definition~\ref{d:t.smooth}), so $|\sss| \le\tep$. We can apply  the bound \eqref{e:Y.I.Kolmogorov} from Proposition~\ref{p:Y.kolmogorov} for increments of $Y^{\frozen,\RomI}$, and the bound from Proposition~\ref{p:Y.notI.second.quenched}
for increments of $Y^{\frozen,\RomII}$, $Y^{\frozen,\RomIII}$, and $Y^{\frozen,\Ising}$. By Fatou's lemma, in the $N\to\infty$ limit we have the quenched second moment bound
	\[
	\bE\bigg[\Big\| \vY^\frozen(t)-\hvY^\frozen(t)\Big\|^2\bigg]
	\le O(\sep)+ L^{O(1)}
	 \tep\,.
	\]
Next recall that $Y^\frozen$ is the coordinate sum \eqref{e:Y.plus}, and 
recall that $\iY^\frozen$ is the coordinate sum  \eqref{e:iY.plus}.
Meanwhile, $\iX^\frozen$ solves the SDE \eqref{e:iX.sde}. By Proposition~\ref{p:fokker.planck},
$\iY^\frozen(t)$ and $\iX^\frozen(t)$ have the same marginal law for all $t$. We finally compare $\iX^\frozen$ to the solution $X^\frozen$ of the SDE \eqref{e:actual.sde}. By Corollary~\ref{c:reg.in.limit}\ref{c:reg.in.limit.a},  the coefficients are uniformly bounded by
$(\frac{\trK L}{\sep})^{O(1)}$. 
By Lemma~\ref{l:ideal.vs.actual.coeffs},
the coefficients of the two SDEs
are close:
	\[
	\sup\bigg\{
	\Big|\tF(t,y)-\iF(t,y)
	\Big|
	: y\in\R, t\in[q_*,1]\setminus\oJ
	\bigg\}
	\le
	\bigg(\frac{\trK L}{\sep}\bigg)^{O(1)} \tep^{1/2}\,.
	\] 
Since for any $x,y\ge0$ we have $|x^{1/2}-y^{1/2}| \le |x-y|^{1/2}$, we immediately deduce a bound on $\tF^{1/2}-\iF^{1/2}$.
Lemma~\ref{l:gronwall} then implies the existence of a coupling under which
	\[
	\bE\bigg[ \Big( X^\frozen(t)-\iX^\frozen(t)\Big)^2\bigg]
	\le
	C\bigg( \bigg(\frac{\trK L}{\sep}\bigg)^{O(1)} 
	\bigg) \tep^{1/2}
	\]
for all $q_*\le t\le 1$. Combining the above bounds gives the conclusion.
\end{proof}
\end{ppn}

\subsection{Final budget constraints for limiting SDE}
\label{ss:sde.reparam.sigma}

We now complete our Lipschitz SDE approximation by simplifying the budget constraint and showing it remains approximately satisfied. The main source of error we address here comes from the discrepancy between $\trrr(t,y)$ and $tp'(t)$. However we first prove an elementary inequality which lets us simplify $(\tuuu,\tvvv)$ to a single parameter $\sigma$: here $u$ and $v$ correspond to $\tuuu$ and $\tvvv$, while $P_1$ and $P_2$ correspond to $p(t)$ and $tp'(t)$.

\begin{ppn}
\label{ppn:budget-constraint-reparametrize} 
Let $u,v$ be real numbers such that $v\geq u^2 - \bar{\lambda}$, and let $P_1,P_2>0$.
Set
\beq\label{e:sigma.of.u.v}
    \sigma 
    \equiv
    \bigg(
    \frac{P_1v + 2
    (P_1P_2)^{1/2}
    u+P_2}{P_1+P_2}
    \bigg)^{1/2}
\eeq
Then we have the bound
\beq\label{eq:u-v-to-sigma-bound}
  u^2 + \Big((v-u^2 + 2\bar{\lambda})^{1/2}-1\Big)^2 
  \ge \bigg(1+\frac{P_2}{P_1}\bigg)(\sigma-1)^2\,.
\eeq
Given any $\sigma\ge0$, if we let $(u,v)$ be defined by the equations
	\[
	u^2 = \frac{P_2}{P_1}(\sigma-1)^2 
	= v-\sigma^2 + 2\bar{\lambda}\,,
	\]
then $(u,v)$ achieves equality in both \eqref{e:sigma.of.u.v} and \eqref{eq:u-v-to-sigma-bound}. 

\begin{proof}
Denote $\tilde{v}\equiv v-u^2 + 2\bar{\lambda} \ge0$. The above definition of $\sigma$ can then be rewritten as
	\beq\label{eq:LM-constraint}    
	(P_1+P_2)\sigma^2
	= P_1 \tilde{v}  +
	\Big( (P_1)^{1/2} u + (P_2)^{1/2}
	\Big)^2\,.
	\eeq
Consider optimizing the left-hand side of \eqref{eq:u-v-to-sigma-bound} subject to the constraint 
\eqref{eq:LM-constraint}: the Lagrangian is
	\[
	\mathscr{L}
	= u^2 + \Big(\tilde{v}^{1/2}-1\Big)^2 
	+ \lambda
	\bigg\{
	(P_1+P_2)\sigma^2
	- P_1 \tilde{v}  -
	\Big( (P_1)^{1/2} u + (P_2)^{1/2}
	\Big)^2
	\bigg\}\,.
	\]
The Lagrangian stationarity conditions are
	\begin{align*}
	0&=\frac{\partial\mathscr{L}}
		{\partial u}
	= 2u - 2\lambda (P_1)^{1/2}
	\Big( (P_1)^{1/2} u + (P_2)^{1/2}
	\Big)\,\\
	0&=\frac{\partial\mathscr{L}}
		{\partial\tilde{v}}
	= 1-\frac{1}{\tilde{v}^{1/2}}
	- \lambda P_1\,.
	\end{align*}
Simplifying the above equations gives
	\[
	u
	= \bigg(\frac{P_2}{P_1}\bigg)^{1/2}
	 \frac{\lambda P_1}{1-\lambda P_1}
	=
	\bigg(\frac{P_2}{P_1}\bigg)^{1/2}
	(\tilde{v}^{1/2}-1)
	\,.
	\]
Combining this with \eqref{eq:LM-constraint} gives $\sigma^2=\tilde{v}$. Substituting this back into the left-hand side of \eqref{eq:u-v-to-sigma-bound} gives
	\[
	u^2 + (\tilde{v}^{1/2}-1)^2
	= \bigg(1+ \frac{P_2}{P_1}\bigg)
	(\tilde{v}^{1/2}-1)^2
	= \bigg(1+ \frac{P_2}{P_1}\bigg)
	(\sigma-1)^2\,,
	\]
so \eqref{eq:u-v-to-sigma-bound} holds with equality in the case of an interior stationary point. It remains to check the boundary case $\tilde{v}=0$. In this case we have
	\[
	\sigma
	=\bigg|\bigg( \frac{P_1}{P_1+P_2}\bigg)^{1/2} u
	+\bigg( \frac{P_2}{P_1+P_2}\bigg)^{1/2}\bigg|
	\equiv
	\Big| (A_1)^{1/2} u + (A_2)^{1/2}
	\Big|\,,
	\]
with $A_1+A_2=1$. The desired bound \eqref{eq:u-v-to-sigma-bound} then simplifies to
	\begin{align*}
	u^2 + 1
	&= \bigg(\frac{-(A_2)^{1/2}\pm\sigma}{(A_1)^{1/2}}\bigg)^2+1
	= \frac{A_2 + \sigma^2
		\pm 2 (A_2)^{1/2}\sigma}{A_1}
		+1
	= \frac{1+\sigma^2
		\pm 2 (A_2)^{1/2}\sigma}{A_1}
	\\
	&\ge
	\frac{1 + \sigma^2-2\sigma}{A_1}
	= \bigg(1+\frac{P_2}{P_1}\bigg)(\sigma-1)^2\,,
	\end{align*}
which gives the bound for the boundary case $\tilde{v}=0$. This concludes the proof.
\end{proof} 
\end{ppn}

\begin{lem}\label{lem:trrr-error-small}
    For any $t\notin \oJ$ (recall Lemma~\ref{l:bad.J}) and parameters as in \eqref{eq:order-of-limits}: 
    \[\bE\Big[|tp'(t) - \trrr_t|\Big]
    \leq  
    o_{\trK}(1)
    \]
for $ \trrr_t\equiv \trrr(t,X^\plus(t))$. 

\begin{proof}
Let $\bar{P}(t)$ be the temporal smoothing of $tp'(t)$. Then for $t\notin\bar{J}$ we have
	\[
    \bE\Big[|tp'(t) - \trrr_t|\Big]
    =
    \bE\Big[|tp'(t) - \bar P(t)|\Big]
    +\bE\Big[|\bar P(t) - \trrr_t|\Big]
    = O(\lambda)
    +\bE\Big[|\bar P(t) - \trrr_t|
    \Big]\,.
    \]
Let $\bar{P}^N$ be the temporal smoothing of $t(p^N)'(t)$. It follows from Lemma~\ref{l:limit.coefs} that $\bar{P}^N$ converges pointwise to $\bar{P}$, and $\trrr^N$ converges pointwise to $\trrr$. It follows by the bounded convergence theorem that we also have convergence in expectation, so
	\[
	\bE\Big[|\bar P(t) - \trrr_t|\Big]
	=\bE\Big[|\bar P^N(t) - (\trrr^N)_t|\Big] + o_N(1)\,.
	\] 
We then note that $\bar{P}^N-\trrr^N$ is the smoothing (via Definitions~\ref{d:sp.smooth} and \ref{d:t.smooth}) of 
	\begin{align} \nonumber
	A^N(t,\vy)
	&\equiv t(p^N)'(t)-\rrr^{\plus,N}(t,\vy)
	= t(p^N)'(t)-\sum_{B^\frozen \subseteq B^{\frozen,\coarse}}
    \frac{|B|}{|B^{\frozen,\coarse}|}
    \rrr^{\trunc,N}(t,B) 
    \\
	&\qquad= \bigg(
	1-\frac{|B^{\frozen,\coarse} \cap [M]|}{|B^{\frozen,\coarse}|} 
	\bigg)
	t(p^N)'(t)
	+ \sum_{B^\frozen \subseteq B^{\frozen,\coarse}} \frac{|B|}{|B^{\frozen,\coarse}|}
	\bigg(
	t(p^N)'(t) -
	\rrr^{\trunc,N}(t,B)\bigg)\,,
	\label{e:P.N.r.frozen.twoterms}
	 \end{align}
where the correspondence between $\vy$ and $B^{\frozen,\coarse}$ is as described around \eqref{e:review.P}. For the first term in
\eqref{e:P.N.r.frozen.twoterms}, at any fixed time $t$, we have
\[\bE\bigg[ 1-\frac{|B^{\frozen,\coarse} \cap [M]|}{|B^{\frozen,\coarse}|} \bigg]
    = 1 - \sum_{B^{\frozen,\coarse}}
    \frac{|B^{\frozen,\coarse}|}{M^\frozen}
    \frac{|B^{\frozen,\coarse} \cap [M]|}{|B^{\frozen,\coarse}|}
    = 1 - \frac{M}{M^\frozen}
    = o_{\trK}(1)\,,
\]
by recalling the choice of $M^\frozen$ from Definition~\ref{d:froze}. For the second term in \eqref{e:P.N.r.frozen.twoterms},  we argued in the proof of Proposition~\ref{p:budget-constraints-without-frozen}\ref{it:Ising-budget-constraint} that the total fraction of particles ever frozen is with high probability at most $o_{\trK}(1)$, so this term is also small on average at any fixed time for the process $Y$. Altogether this shows
	\beq\label{e:AN.small.finiteN}
	\bE\Big[ |A^N(t,\vY^{\frozen,N,\RomI:\RomIII}(t)) |\Big]
	= \int_{\R^3} |A^N(t,\vy)| \varrho_{N,t}(d\vy)
	\le o_{\trK}(1)\,.
	\eeq
Let $\tilde{A}^N$ be the spatial smoothing of $A^N$, via Definition~\ref{d:sp.smooth}. Let $\tvY^{\frozen,N}(t)\equiv \vY^{\frozen,N}(t)+\xxx$, and denote the projection onto the sum of the first three coordinates as
 $\tilde{Y}^{\frozen,N} = \tilde{Y}^{\frozen,N,\RomI}+\tilde{Y}^{\frozen,N,\RomII}+\tilde{Y}^{\frozen,N,\RomIII}$. Then
	\begin{align*} &\bE\Big[ |\tilde{A}^N(t, \tilde{Y}^{\frozen,N}(t)) |\Big]
	= \int_{\R} |\tilde{A}^N(t,y)| \overline{\varrho_{N,t}*\sphi}(y)\,dy \\
	&\qquad = \int_{\R} \bigg|
	\int_{H(y)} \int_{\R^3} \sphi(\vy-\vx) A^N(t,\vx)
		\,\rho_{N,t}(d\vx)\,d\vy \bigg| \,dy \\
	&\qquad\le  \int_{\R^3}
	|A^N(t,\vx)| \bigg( \int_{\R^3}
		 \sphi(\vy-\vx) \,d\vy\bigg)
		\,\rho_{N,t}(d\vx) =  \int_{\R^3}
	|A^N(t,\vx)|\rho_{N,t}(d\vx) \le o_{\trK}(1)\,,
	\end{align*}
using the previous estimate \eqref{e:AN.small.finiteN}. This shows that the bound \eqref{e:AN.small.finiteN} can be passed to the spatially smoothed quantity $A^N$. Similarly, let $\bar{A}^N$ be the temporal smoothing of $\tilde{A}^N$, via Definition~\ref{d:t.smooth}. Then
	\begin{align*}
	&\bE\Big[ |\bar{A}^N(t,\tilde{Y}^{\frozen,N}(t))| \Big]
	=\int_{\R} |\bar{A}^N(t,y)| \overline{\varrho_{N,t}*\sphi}(y)\,dy\\
	&\qquad\le \int_{\R}  
	\bigg( \int_{\R}\tphi(t-s) |  \tilde{A}^N(s) |\,ds  \bigg) \,
	 \overline{\varrho_{N,t}*\sphi}(y)\,dy \\
	&\qquad=
	\int_{\R} \tphi(t-s)
	\bigg[ \int_{\R} | \tilde{A}^N(s) | 
	 \overline{\varrho_{N,s}*\sphi}(y)
	 \,dy \bigg] \,ds + 
	 	\bigg(\frac{L}{\sep}\bigg)^{O(1)}
	 	\TMAX_N \tep^{1/2} \,,
	\end{align*}
where the last estimate comes from Lemma~\ref{l:sp.smooth.cty}. In the last line above, the integral inside the square brackets can be bounded by the previous estimate, giving
	\[
	\bE\Big[ |\bar{A}^N(t,\tilde{Y}^{\frozen,N}(t))| \Big]
	\le  o_{\trK}(1) + 
	 	\bigg(\frac{L}{\sep}\bigg)^{O(1)}
	 	\TMAX_N \tep^{1/2}  \le  o_{\trK}(1)\,,
	\]
where the last bound follows by  \eqref{e.JMAX.bound} and recalling the order of limits from Assumption~\ref{a:params}. It follows by the dominated convergence theorem that the same bound holds in the $N\to\infty$ limit:
	\[\bE\Big[ |\bar{A}(t,\tilde{Y}^\frozen(t))| \Big]
	\le  o_{\trK}(1)\,,\]
where $\tilde{Y}^\frozen$ is the (subsequential) limit of $\tilde{Y}^{\frozen,N}$. Finally, we argue that the same estimate holds with $X^\plus(t)$ in place of $\tilde{Y}^\frozen(t)$. To this end, note that the argument of Proposition~\ref{p:sde.approx} also gives
	\beq\label{e:tilde.Y.X.discrep}
	\bE\bigg[\Big( X^\frozen(t)-\tilde{Y}^\frozen(t)\Big)^2\bigg]
	\le	C\bigg( \bigg(\frac{\trK L}{\sep}\bigg)^{O(1)} 
	\bigg) \tep^{1/2}\,.
	\eeq 
Indeed, in the proof of Proposition~\ref{p:sde.approx}, the error term $O(\sep)$ comes from approximating $Y^\frozen(t)$ by its spatial rerandomization $\tilde{Y}^\frozen(t)$, and this term does not arise if we bound $X^\frozen(t)-\tilde{Y}^\frozen(t)$ rather than $X^\frozen(t)-Y^\frozen(t)$. Corollary~\ref{c:reg.in.limit}\ref{c:reg.in.limit.a}
bounds the Lipschitz norm of $\bar{A}$. It follows that
	\[\bE\bigg[ \Big|\bar{A}(t,\tilde{Y}(t))
		- \bar{A}(t, X^\plus(t)) \Big| \bigg]
	\le
	\bigg( \frac{\trK L}{\sep}\bigg)^{O(1)} \cdot
	\bE\bigg[\Big( X^\frozen(t)-\tilde{Y}^\frozen(t)\Big)^2\bigg]^{1/2}
	\le o_{\trK}(1)\,,
	\]
again recalling Assumption~\ref{a:params}. This proves the claim.
\end{proof}
\end{lem}

Using Proposition~\ref{ppn:budget-constraint-reparametrize}, we now reparametrize the SDEs $X^\frozen,X^{\frozen,\Ising}$ and show they approximately obey the corresponding budget constraints.
In what follows we use $o(1)$ to denote an error tending to zero in the asymptotic regime described by Assumption~\ref{a:params}.

\begin{thm}\label{thm:Lipschitz-SDE-approx}
Let $\cA_N$ be a sequence of $L$-Lipschitz algorithms such that $\chi^N\equiv\chi_{\cA_N}$ satisfies
\[\frac{1}{L^2} \le (\chi^N)'(p) \le L^2\]
for all $p\in[0,1]$, and converges to a limiting function $\chi$ both pointwise and in $L^1$. 
Consider a Brownian bridge terminating at $(\bG,\bg^\aux)$, and let $\vX^N\equiv (\vX^{\RomI:\RomIII},X^\Ising)^N$ be the resulting process defined by \eqref{e:X.decomp}--\eqref{e:vX}. Let $\vX^{\frozen,N}$ be the modification given by Definitions~\ref{d:trunc}--\ref{d:buckets.Ising}, and let $\vY^{\frozen,N}$ be the spatial rerandomization of $\vX^{\frozen,N}$ from Definitions~\ref{d:rerand} and \ref{d:rerand.Ising}. By passing to a subsequence, assume that we have a subsequential limit $\vY^\frozen$ as in Theorem~\ref{t:tightness}. Writing $p^N\equiv (\chi^N)^{-1}$ and $p\equiv \chi^{-1}$, and passing to a further subsequence as needed, assume that $(p^N)'$ converges in the L\'evy--Prohorov metric to $p'$, in the sense of Definition~\ref{d:levy.prohorov.derivs}. Denote $q_0\equiv \chi(0)$. Then, for some $q_*\in[q_0,q_0+L^{-\Theta(1)}]$, we have an SDE
	\beq\label{e:actual.sde.reparametrized}   \begin{aligned}
	dX^{\frozen}(t)
    &= \tB(t,X^{\frozen}(t))\,dt
    + \tSigma(t,X^{\frozen}(t))\,dB(t)
    ;
    \\ dX^{\frozen,\Ising}(t)
    &=
    w(t,X^{\frozen,\Ising}(t))\,d B^\Ising(t)
    \end{aligned}\eeq
on the time interval $[q_*,1]$, satisfying the following conditions:
\begin{enumerate}[(a)]
\item 
\label{it:Lipschitz-SDE-approx.Lip}
The functions $\tB,\tSigma, w:[q_*,1]\times\bbR\to\bbR$
 are $C(L)$-bounded, and $C(L)$-Lipschitz in space and time.
 
 \item \label{thm:Lipschitz-SDE-approx.b}
In addition, $w_t\equiv w(t,X^{\frozen,\Ising}(t))$ must satisfy, for all $t\in[q_*,1]$, the constraint
 \[
  \Big|\bE[(w_t)^2]-1\Big| \le o(1)\,.\]

\item\label{thm:Lipschitz-SDE-approx.a}
 Denote 
$b\equiv \tB/p'(t)^{1/2}$ and 
$\sigma\equiv \tSigma/[p(t) + t p'(t)]^{1/2}$. Abbreviate $b_t\equiv b(t,X^\frozen(t))$, etc. 
Then
    \[
        \int_{q_*}^1
        \bigg(\bE\bigg[
        (b_t)^2+\lt(\frac{p(t)+tp'(t)}{p(t)}\rt)
        (\sigma_t-1)^2
        \bigg]
        -
        \frac{(\bE w_t)^2}{\alpha}
        \bigg)_+ \,dt
        \leq 
        o(1)\,.\]

\item \label{thm:Lipschitz-SDE-approx.c.INIT}
The SDE initializations satisfy 
    \begin{align}
    \label{eq:SDE-init-1}
    \bE[X^{\frozen}(q_*)^2]&\leq o(1)\,,
    \\
    \label{eq:SDE-init-2}
|\bE[X^{\frozen,\Ising}(q_*)^2-q_*]| &\leq o(1)\,,
    \end{align}

\end{enumerate}
Finally, we have the $\bbW_2$ approximation estimates
\begin{align*}
\bbW_2(\Law(X^{\frozen}(1)),\Law(Y^\plus(1))) &\le o(1)\,,\\
\bbW_2(\Law(X^{\frozen,\Ising}(1)),\Law(Y^{\plus,\Ising}(1))) &\le o(1)\,.
\end{align*}
In all the above assertions, the $o(1)$ denotes an error which can tend to zero as $L^{-\Theta(1)}$. 

\begin{proof}
Apply the procedure of this section to obtain the define the processes $X^\frozen$ and $X^{\frozen,\Ising}$ from the SDE \eqref{e:actual.sde}, where we recall that the coefficients are given by Definition~\ref{d:cts.coeffs}. 
Recall the order of parameters from Assumption~\ref{a:params}, and note that we are considering the $N\to\infty$ limit while all the other parameters are kept finite. 
As indicated above,
we define
	\[b(t,y)\equiv\tbbb(t,y)\equiv
	\frac{\tB(t,y)}{p'(t)^{1/2}} \,,\]
so that the drift terms for $X^\plus(t)$ in 
\eqref{e:actual.sde} and
\eqref{e:actual.sde.reparametrized} coincide exactly. 
Recall from Lemma~\ref{l:limit.coefs}
that the quantities $\twww,\iw,\tF,\iF$ are all nonnegative.
We also define
 $w\equiv\twww^{1/2}$, so that the diffusive terms for $X^{\plus,\Ising}(t)$ in \eqref{e:actual.sde} and
\eqref{e:actual.sde.reparametrized} also coincide exactly.
Then, as suggested by Proposition~\ref{ppn:budget-constraint-reparametrize}, we define
        \begin{align*}\sigma
        &\equiv
        \bigg(
        \frac{p(t)\tvvv+2[p(t) \cdot tp'(t)]^{1/2}\tuuu+\trrr}{p(t)+tp'(t)}
        \bigg)^{1/2}
        =\frac{\tF^{1/2}}{[p(t)+tp'(t)]^{1/2}}
        \,,
        \\
        \ddot{\sigma}
        &=\bigg(
        \frac{p(t)\tvvv+2[p(t) \cdot tp'(t)]^{1/2}
        \tuuu+tp'(t)}{p(t)+tp'(t)}
        \bigg)^{1/2}
        \equiv\frac{\ddot{F}^{1/2}}{[p(t)+tp'(t)]^{1/2}}
       \,.
    \end{align*}
We also write $\tSigma\equiv\tF^{1/2}$ and
 $\ddot{\Sigma}\equiv \ddot{F}^{1/2}$. 
For this definition of $\sigma$ and $\tSigma$, the diffusive terms for $X^\plus(t)$ in
\eqref{e:actual.sde} and
\eqref{e:actual.sde.reparametrized} coincide exactly.
\begin{itemize}
\item Proof of \eqref{it:Lipschitz-SDE-approx.Lip}:
For the drift coefficient $\tB$,
the supremum norm and Lipschitz constant are  bounded by Corollary~\ref{c:reg.in.limit}\ref{c:reg.in.limit.a}. For the diffusivity coefficients
$w=\twww^{1/2},
\tSigma=\tF^{1/2},
\ddot{\Sigma}=\ddot{F}^{1/2}$,
the supremum norm and Lipschitz constant are bounded by Corollary~\ref{c:reg.in.limit}\ref{c:reg.in.limit.b}.
 
\item 
Proof of \eqref{thm:Lipschitz-SDE-approx.b}: 
Recall from Proposition~\ref{p:budget-constraints-coarsened}\ref{it:Ising-w-constraint-plus-coarse} that with high probability, $\bE\www^{\frozen,N}(t)$ is close to $1$ for all times $t=q_d$. Now recall from Definition~\ref{d:sp.smooth} that 
	\[\swww^N(t,y)
	\equiv
	\frac{\varsigma_{\www,N,t}*\sphi(y)}
		{\varsigma_{\vY,N,t}*\sphi(y)}
	=
	\frac{1}
		{\varsigma_{\vY,N,t}*\sphi(y)}
	\int 
	\sphi(y-x)
	\www^{\frozen,N}(t,x) \varsigma_{\vY,N,t}(dx)
	\]
with $y\in\R$. It follows that 
	\begin{align*}
	&\bE\swww^N(t,\tilde{Y}^{\frozen,N,\Ising}(t))
	= \int_{\R}
	\swww^N(t,y)
	\varsigma_{\vY,N,t}*\sphi(y)\,dy \\
	&\qquad
	= \int_{\R}
	\bigg( \int_{\R}
	\sphi(y-x) \,dy\bigg)
	\www^{\frozen,N}(t,x) \varsigma_{\vY,N,t}(dx) 
	= \bE\www^{\frozen,N}(t)\,.
	\end{align*}
It therefore follows directly from Proposition~\ref{p:budget-constraints-coarsened}\ref{it:Ising-w-constraint-plus-coarse} that with high probability, 
	\beq\label{eq:swww-1st-moment}
	\Big|\bE\swww^N(t,\tilde{Y}^{\frozen,N,\Ising}(t))-1\Big| 
	\le \epsilon^\circ 
	\eeq
for all times $t=q_d$. Next, recalling Definition~\ref{d:t.smooth}, we expand
	\[
	\bE\twww^N(t,\tilde{Y}^{\frozen,N,\Ising}(t))
	=\int_{\R}\bigg(
	\int_{\R}
	\tphi(t-s)
	\swww^N(s,y)\,ds\bigg)
	\varsigma_{\vY,N,t}*\sphi(y)\,dy
	\,.
	\]
Recall the kernel $\tphi(t-s)$ is only non-zero for $|t-s|\le\tep$, and Lemma~\ref{l:sp.smooth.cty} controls the regularity of the mapping
$s\mapsto \varsigma_{\vY,N,s}*\sphi(y)$. Lemma~\ref{l:sp.smooth.reg} bounds the size of $\swww^N$. Altogether it gives
	\begin{align*}
	&\Big|\bE\twww^N(t,\tilde{Y}^{\frozen,N,\Ising}(t))-1\Big|\\
	&\qquad\le \bigg| \int_{\R}
	\tphi(t-s)
	\bigg(\int_{\R}
	\swww^N(s,y)
	\varsigma_{\vY,N,s}*\sphi(y)\,dy
	\bigg)
	\,ds -1\bigg| + \bigg(\frac{\trK L}{\sep}\bigg)^{O(1)} 
	(\MAX_N)^2 \tep^{1/2}\\
	&\qquad\le
	\int_{\R}
	\tphi(t-s)
	\Big| \bE\swww^N(t,\tilde{Y}^{\frozen,N,\Ising}(t))
	-1\Big|
	\,ds
	+ \bigg(\frac{\trK  L}{\sep}\bigg)^{O(1)} 
	(\MAX_N)^2 \tep^{1/2}
	\,.
	\end{align*}
Applying \eqref{eq:swww-1st-moment} and recalling Assumption~\ref{a:params} gives
	\[\Big|\bE\twww^N(t,\tilde{Y}^{\frozen,N,\Ising}(t))-1\Big|
	\le  2\epsilon^\circ\,.\]
We now send $N\to\infty$. Passing to a subsequence as needed, we can assume we have a coupling for which $\tilde{Y}^{\frozen,N,\Ising}$ converges uniformly to the limiting process 
$\tilde{Y}^{\frozen,\Ising}$. We also have $\twww^N$ converging pointwise to $\twww$ by Lemma~\ref{l:limit.coefs}. It follows by the bounded convergence theorem that
	\[
	\Big|\bE\twww(t,
	\tilde{Y}^{\frozen,\Ising}(t))-1\Big| 
	= \lim_{N\to\infty}
	\Big|\bE\twww^N(t,
		\tilde{Y}^{\frozen,N,\Ising}(t))-1\Big| 
	\le 2\epsilon^\circ\,.
	\]
Now recall that Corollary~\ref{c:reg.in.limit}
gives regularity for $\twww$, while
 \eqref{e:tilde.Y.X.discrep} bounds the discrepancy between
$X^\plus(t)$ and $\tilde{Y}^\frozen(t)$. It follows that a similar bound holds with
$X^\plus(t)$ in place of $\tilde{Y}^\frozen(t)$, concluding the proof of \eqref{thm:Lipschitz-SDE-approx.b}.

\item
Proof of \eqref{thm:Lipschitz-SDE-approx.a}: Let $J\equiv J(p,2\tep,\lambda/2)$ as defined by \eqref{e:bad.times.J}. It follows by combining Lemma~\ref{l:limit.coefs} and    Corollary~\ref{c:limit.coefs.variants.of.temp.smoothing} that
    	\[
	\tvvv(t,y)-\tuuu^2(t,y) \ge -\bar{\lambda} \equiv - \bigg(\frac{\trK L}{\sep}\bigg)^{O(1)}\lambda
	\]
provided $(t,y)\in([q_*,1]\setminus J)\times\R$. In this case we can apply Proposition~\ref{ppn:budget-constraint-reparametrize} to obtain
	\begin{align*}
	&\bigg\{
	b^2 + \bigg(1+\frac{tp'(t)}{p(t)}\bigg)(\ddot{\sigma}-1)^2
	\bigg\}(t,y) \\
	&\qquad\le\bigg\{ b^2 + \tuuu^2 + 
	\Big((\tvvv-\tuuu^2 +2 \bar{\lambda} )^{1/2}-1\Big)^2 
	\bigg\}(t,y)
	=\Cost_{ 2\bar{\lambda}}\Big(\tGam(t,y)\Big)
	\end{align*}
for all $(t,y)\in([q_*,1]\setminus J)\times\R$. Combining with
Corollary~\ref{c:limit.budget} and recalling Assumption~\ref{a:params} gives
	\[\bE \Big[
		\Cost_{2\bar{\lambda}}\Big(
			\tGam(t,
			\tilde{Y}^\frozen(t))
			\Big)\Big]
	- \frac{(\bE w(t,\tilde{Y}^\frozen(t)) )^2}{\alpha} 
	\le \bigg(\frac{\trK L}{\sep}\bigg)^{O(1)}  
	(\tep^{1/2} + \lambda^{1/2})
	+ 3\epsilon^\circ
	\le 4\epsilon^\circ
	\]
for all $t\in[q_*,1] \setminus J$, where $\tilde{Y}^\frozen$ is the same as in the proof of Lemma~\ref{lem:trrr-error-small}. Combining the last two displays gives,
again for all $t\in[q_*,1] \setminus J$, the bound
	\[\bE\bigg[ 
	b(t,\tilde{Y}^\frozen(t))^2
	+ \bigg(1+\frac{tp'(t)}{p(t)}\bigg)(
		\ddot{\sigma}(t,\tilde{Y}^\frozen(t))
		-1)^2\bigg]
	\le \frac{(\bE w(t,\tilde{Y}^\frozen(t)))^2}{\alpha} 
	+ 4\epsilon^\circ\,.
	\]
Now \eqref{e:tilde.Y.X.discrep} bounds the discrepancy between
$X^\plus(t)$ and $\tilde{Y}^\frozen(t)$,
while Corollary~\ref{c:reg.in.limit} gives regularity for $b,\ddot{\sigma},w$ (here we only use regularity in the space coordinate). It follows that a similar bound holds with $X^\plus(t)$ in place of $\tilde{Y}^\frozen(t)$,  which we write in shorthand notation as
	\[\bE\bigg[ 
	(b_t)^2 + \bigg(1+\frac{tp'(t)}{p(t)}\bigg)(
		\ddot{\sigma}_t-1)^2\bigg]
	\le \frac{(\bE w_t)^2}{\alpha}     
	+ 5\epsilon^\circ
	\]
for all $t\in[q_*,1] \setminus J$. We now need to replace $\ddot{\sigma}$ with $\sigma$: it follows from the definitions that $\sigma$ and $\ddot{\sigma}$ satisfy
	\[
	\bigg(1+ \frac{tp'(t)}{p(t)}\bigg)\bE\Big[
		\Big| (\sigma_t)^2-(\ddot{\sigma}_t)^2 \Big| \Big]
	\le  L^{O(1)} \bE\Big[| t p'(t)- \trrr_t|\Big]\,.
	\]
We note additionally that we can bound
	\[
	\bigg(\bE\Big[\Big|\sigma_t-\ddot{\sigma_t}\Big|\Big]\bigg)^2
	\le \bE\Big[ (\sigma_t-\ddot{\sigma_t})^2\Big]
	\le \bE\Big[ \Big| (\sigma_t)^2 -(\ddot{\sigma}_t)^2 \Big|\Big]\,.
	\]
Combining the preceding estimates gives
	\begin{align*}
	&\bigg(1+ \frac{tp'(t)}{p(t)}\bigg)\bigg|
		\bE\Big[(\sigma_t-1)^2\Big]
		-\bE\Big[(\ddot{\sigma}_t-1)^2\Big]\bigg|\\
	&\le L^{O(1)}
	\max\bigg\{
	\bE\Big[| t p'(t)- \trrr_t|\Big],
	\bE\Big[| t p'(t)- \trrr_t|\Big]^{1/2}\bigg\}\,.
	\end{align*}
Lemma~\ref{lem:trrr-error-small} shows that the right-hand side above is $o_{\trK}(1)$ for $t\notin \oJ$. Combining the above estimates gives
	\[
	\int_{[q_*,1]\setminus\oJ}
	\bigg(
	\bE\bigg[
        (b_t)^2+\lt(\frac{p(t)+tp'(t)}{p(t)}\rt)
        (\sigma_t-1)^2
        \bigg]
        -
        \frac{(\bE w_t)^2}{\alpha}
        \bigg)_+\,dt
   \le 
   5\epsilon^\circ
   +o_{\trK}(1)\,.
   \]
On the other hand, for all $t\in[q_*,1]$ (including $t\in\oJ$), we can crudely bound
	\[\bE\bigg[
        (b_t)^2+\lt(\frac{p(t)+tp'(t)}{p(t)}\rt)
        (\sigma_t-1)^2
        \bigg]
        -
        \frac{(\bE w_t)^2}{\alpha}
    \le \bigg(\frac{\trK L}{\sep}\bigg)^{O(1)}\,,
    \]
by dropping the term involving $w$ and using the bounds on $b,\sigma$ implied by Corollary~\ref{c:reg.in.limit}, together with the bound $1/p(q_*)\le L^{O(1)}$ from the discussion around \eqref{e:breve.q.bounds}. Applying Lemma~\ref{l:bad.J} gives
	\[
	\int_{\oJ}
	\bigg(
	\bE\bigg[
        (b_t)^2+\lt(\frac{p(t)+tp'(t)}{p(t)}\rt)
        (\sigma_t-1)^2
        \bigg]
        -
        \frac{(\bE w_t)^2}{\alpha}\bigg)_+
      \le\bigg(\frac{\trK L}{\sep}\bigg)^{O(1)}
      \frac{\tep}{\lambda}\,.
	\]
Combining the integrals over $[q_*,1]\setminus\oJ$ and $\oJ$ concludes the proof of \eqref{thm:Lipschitz-SDE-approx.a}.

\item Proof of \eqref{thm:Lipschitz-SDE-approx.c.INIT}: 
Recall from
\eqref{e:process.X} and \eqref{e:process.X.Ising} that 
	\begin{align*}
	X^N(q_d,a)
	&= \frac{(\bg^a(q_d),\bx(q_d))}
		{N^{1/2}}\,,\\
	X^{N,\Ising}(q_d,i)
	&= \Big(\be_i,\bx(q_d)\Big)\,.
	\end{align*}
Since $p_0=0$, the vector $\bg^a(q_0)$ is simply zero, so $X^N(q_0,a)=0$ for all $a\in[M]$. On the other hand, by Lemma~\ref{l:nearly.ultrametric}, it holds with very high probability that
	\[\Big|
	\bE[X^{N,\Ising}(q_0)^2]-q_0\Big|
	=\bigg| \frac{\|\bx(q_0)\|^2}{N}-q_0
		\bigg|
	\le \frac1{N^{1/4}}\,.\]
Now recall we truncate to form $\vX^{\trstar}$ (Definition~\ref{d:trunc}),
 add frozen particles
 to form $\vX^{\frstar,N}$ (Definition~\ref{d:froze}), and freeze small spatial buckets to form $\vX^{\frozen,N}$ (Definitions~\ref{d:buckets} and \ref{d:buckets.Ising}). In particular, since most of these steps do not affect what happens at time $q_0$,
 we observe that
 $\vX^{\frozen,N}(q_0)$ differs from
 $\vX^N(q_0)$ only by
the addition of frozen particles (Definition~\ref{d:froze}). Thus we conclude
	\[\max
	\bigg\{\bE[X^{\frozen,N}(q_0)^2],
	\Big| \bE[X^{\frozen,N,\Ising}(q_0)^2] -q_0\Big| 
	\bigg\}
	 \le o_{\trK}(1)
	\]
Recalling the spatial rerandomization of 
Definitions~\ref{d:rerand} and \ref{d:rerand.Ising}, the same estimates also hold for $\vY^{\frozen,N}$. By Fatou's lemma for $Y^{\frozen,N}(q_0)$, and the bounded convergence theorem for $Y^{\frozen,N,\Ising}(t)$, we can send $N\to\infty$ to conclude
	\[\max
	\bigg\{\bE[Y^\frozen(q_0)^2],
	\Big|\bE[Y^{\frozen,\Ising}(q_0)^2] -q_0\Big| 
	\bigg\}
	 \le o_{\trK}(1)\,.
	\]
Recall the definition of $q_*$ from \eqref{e:qbar.and.qstar}. It follows using Propositions~\ref{p:Y.kolmogorov} and \ref{p:Y.notI.second.quenched} that
	\beq
	\max
	\bigg\{\bE[Y^\frozen(q_*)^2],
	\Big| \bE[Y^{\frozen,\Ising}(q_*)^2] -q_*\Big| 
	\bigg\}\\
	\le o_{\trK}(1)
	+ L^{O(1)}\brep^{1/2}
	\le o(1)\,,\label{e:use.of.breve.c}
	\eeq
where for the last bound we take $\breve{c}$ to be a sufficiently large absolute constant in the definition 
$\brep=1/L^{\breve{c}}$. In this last step, it is essential to note the following: in earlier occurrences of 
$L^{O(1)}$ inside this proof, the $O(1)$ exponent can depend on the parameter $\breve{c}$ from \eqref{e:breve.q.bounds}. \textbf{However, in the $L^{O(1)}\brep^{1/2}$ term in  
\eqref{e:use.of.breve.c}, the $O(1)$ exponent comes from Propositions~\ref{p:Y.kolmogorov} and \ref{p:Y.notI.second.quenched}, and does not depend on $\breve{c}$ --- therefore it is indeed possible to take $\breve{c}$ large enough to make $L^{O(1)}\brep^{1/2} \le o(1)$.} Finally, it follows from Proposition~\ref{p:sde.approx} that a similar estimate holds with 
$\vX^\plus$
in place of $\vY^\plus$, and this concludes the proof of \eqref{thm:Lipschitz-SDE-approx.c.INIT}.  
\end{itemize}
The final $\bbW_2$ error bounds are  a direct consequence of Proposition~\ref{p:sde.approx}. We also remark that in all the $o(1)$ error bounds obtained in this theorem, 
the $o(1)$ can be taken to tend to zero at least polynomially in $L$, with the largest error coming from \eqref{e:use.of.breve.c}.
\end{proof}\end{thm}

We finally conclude this section with the proofs of Theorems~\ref{thm:BOGP-hardness-main}
and \ref{thm:SDE-smoothing-general}.

\begin{proof}[\hypertarget{proof:t.BOGP-hardness-main}{Proof of Theorem~\ref{thm:BOGP-hardness-main}}]
Given the general $L$-Lipschitz algorithms $(\cA_N)^\circ$, we first define the perturbations $\cA_N$ as in Proposition~\ref{p:wlogable}. Write $\chi^{N,\circ}\equiv \chi_{(\cA_N)^\circ}$ and $\chi^N\equiv\chi_{\cA_N}$.
Recall from Proposition~\ref{p:wlogable} that
	\[
	\chi^N(p) 
	=  
    \bigg(1-\frac{1}{(1-q_0)L^2}\bigg) \chi^\circ(p)+ \frac{q_0 + (1-q_0)p}{(1-q_0)L^2}\,.
	\]
Together with Proposition~\ref{p:correlation-fn-ub},
we see that we  have
	\[
	\frac{1}{L^2}
	\le (\chi^N)'(p) \le L^2
	\]
for all $0\le p \le 1$. Since the $\chi^N$ are monotone functions, it follows by the Helly selection principle that $\chi^N\to\chi$ pointwise along a subsequence $N_j\to\infty$. Since the $\chi^N$ are bounded functions on the bounded interval $[0,1]$, this also implies convergence in $L^1$. Denote $p^N\equiv (\chi^N)^{-1}$. By passing to a further subsequence, we can assume that $(p^N)'$ converges to $p'$ in the L\'evy--Prohorov sense. Apply $\cA_N$ to define processes $\vX^N$ as in \eqref{e:X.decomp}--\eqref{e:vX}, and then consider the resulting $\vX^{\frozen,N}$ (Definitions~\ref{d:trunc}--\ref{d:buckets.Ising}) and $\vY^{\frozen,N}$
(Definitions~\ref{d:rerand} and \ref{d:rerand.Ising}). By passing to a further subsequence, assume that we have a subsequential limit $\vY^\frozen$ as in Theorem~\ref{t:tightness}. Applying Theorem~\ref{thm:Lipschitz-SDE-approx} gives an SDE \eqref{e:actual.sde.reparametrized}, which we now relabel as $(X,X^\Ising)$ (dropping the ``$\plus$'' notation), such that 
	\beq\label{e:sde.w2.apx.guarantee.repeated}
	\bbW_2\Big(
	\Law(X(1)),\Law(Y^\frozen(1))\Big)
	+\bbW_2\Big(
	\Law(X^\Ising(1)),
		\Law(Y^{\frozen,\Ising}(1))
		\Big)
	\le o(1)\,.
	\eeq
Then:
\begin{itemize}
\item conclusion~\eqref{it:BOGP-hardness-main-Lip} follows from Theorem~\ref{thm:Lipschitz-SDE-approx}\ref{it:Lipschitz-SDE-approx.Lip}
(recalling 
$b\equiv \tB/p'(t)^{1/2}$, 
$\sigma\equiv \tSigma/s(t)$);
\item conclusion~\eqref{it:BOGP-hardness-main-diffusivity}
follows from Theorem~\ref{thm:Lipschitz-SDE-approx}\ref{thm:Lipschitz-SDE-approx.b};
\item conclusion~\eqref{it:BOGP-hardness-main-budget}
follows from Theorem~\ref{thm:Lipschitz-SDE-approx}\ref{thm:Lipschitz-SDE-approx.a};
\item conclusion~\eqref{it:BOGP-hardness-main-initial} follows from Theorem~\ref{thm:Lipschitz-SDE-approx}\ref{thm:Lipschitz-SDE-approx.c.INIT};
\item conclusion~\eqref{it:BOGP-hardness-main-qstar-condition}
follows from the choice of $q_*$ in Theorem~\ref{thm:Lipschitz-SDE-approx}; see the discussion around \eqref{e:breve.q.bounds}.
\end{itemize}
Lastly, we bound
	\[
	\bbW_2(\mu,\Law(X(1)))
	\le
	\bbW_2(\mu,\mu(\cA_N))
	+\bbW_2\Big(\mu(\cA_N),
		\Law(Y^\frozen(1))\Big)
	+\bbW_2\Big(
		\Law(Y^\frozen(1)),
		\Law(X(1))
		\Big)
	\]
On the  above right-hand side, the first term is at most $\eps_\textup{apx}$ by assumption; the last term is at most $o(1)$ by Theorem~\ref{thm:Lipschitz-SDE-approx}; and the second term is at most $O(\eps_\textup{apx}) + o_{\trK}(1)$ by Proposition~\ref{p:summarize.W2.errors.preprocessing} (recalling from the spatial rerandomization that $\vY^\frozen(1)$ has the same law as $\vX^\frozen(1)$). The conclusion follows.
\end{proof}

\begin{proof}[\hypertarget{proof:t.SDE-smoothing-general}{Proof of Theorem~\ref{thm:SDE-smoothing-general}}]
The proof is an adaptation of the \hyperlink{proof:t.BOGP-hardness-main}{proof of Theorem~\ref{thm:BOGP-hardness-main}}, and primarily amounts to repeating the arguments of this section for the modified setting.
One convenient simplification is that we begin with continuous-time processes $Z$ and $Z^\Ising$, instead of sequences of processes indexed by $N$. Since we now assume $p(q_0)\ge 1/L$, we no longer have issues with division by zero, so we no longer need to choose $q_*$ as in \eqref{e:breve.q.bounds}--\eqref{e:qbar.and.qstar}. Instead, we can take any $q_*\in[q_0,q_0+\epsilon]$, provided $q_*-q_0 \gg\tep$ so that the temporal smoothing step of Definition~\ref{d:t.smooth} still goes through. Additionally, thanks to the  regularity  assumption on $p$, we will be able to take $J=\emptyset$ in \eqref{e:bad.times.J}. We use the same parameters as in Assumption~\ref{a:params}, \textbf{which now depend on the given process $(p,\au,\ab,\asig,\Law(Z),\Law(Z^{\Ising}))$,  as well as on the error tolerance $\epsilon$ from the statement of the theorem}. See e.g.\ step \eqref{i:W2.nontrunc.to.trunc} below where $\trK$ needs to depend on $\Law(Z)$. In this setting \eqref{eq:order-of-limits} is adapted to 
	\[\ldots
	\ll \epsilon^\circ
	\ll \min\bigg\{\frac{1}{L},\epsilon
	\bigg\}
	\ll
	1\,.\] 
We enumerate the adaptations to the proof below:
\begin{enumerate}[(A)]
\item Similarly as in Proposition~\ref{ppn:budget-constraint-reparametrize}, let $\av$, $\au$ be defined by
	\[
	\au^2 = \frac{tp'(t)}{p(t)} (\asig-1)^2
	=\av-\asig^2\,.
	\]
Let $B^1,B^2$ be independent one-dimensional Brownian motions, and let
	\begin{align*}
	dZ^\RomI(t)
	&= p'(t)^{1/2}\ab_t\,dt\,,\\
	dZ^\RomII(t)
	&= p(t)^{1/2} [\av_t-(\au_t)^2]^{1/2}\,dB^1(t)
		+ p(t)^{1/2} \au_t \, dB^2(t)\,, \\
	dZ^\RomIII(t)
	&=[ tp'(t)]^{1/2} \, dB^2(t)\,.
	\end{align*}
Then one easily verifies that $Z$ is equidistributed as $Z^\RomI+Z^\RomII+Z^\RomIII$. 
Let $\vD$ denote its drift, and $\vQ$ its quadratic variation. Note that we have
	\beq\label{e:Z.drift.qv}
	\begin{aligned}
	dD^\RomI(t) &= p'(t)^{1/2}\ab_t\,dt\,,\\
	dQ^{\RomII}(t) 
		&= p(t) \av_t\,dt\,,\\
	dQ^{\RomII,\RomIII}(t) 
		&= [p(t) \cdot t p'(t)]^{1/2} \au_t\,dt\,,\\
	dQ^{\RomIII}(t) &=t p'(t)\,dt\,,\\
	d Q^{\Ising}(t) &= (\aw_t)^2\,dt
	\end{aligned}\eeq
We hereafter abbreviate $\vZ\equiv (\vZ^{\RomI:\RomIII},Z^\Ising)$.

\item 
Analogously to Definition~\ref{d:trunc}, let $\vZ^{\RomI:\RomIII,\trunc}$ be the process $\vZ^{\RomI:\RomIII}$ stopped at the first time $q_T$ that it exits $[-\trK,\trK]^3$. (Meanwhile, $Z^\Ising$ is a martingale that reaches $\{-1,+1\}$ at time $t=1$, so we must have $|Z^\Ising(t)|\le1$ for all $0\le t\le1$.) Denote $\vZ^\trunc \equiv (\vZ^{\RomI:\RomIII,\trunc},Z^{\Ising,\trunc})$. Denote its drift and quadratic variation $\vD^\trunc$ and $\vQ^\trunc$. Note then that
$\vD^\trunc$ and $\vQ^\trunc$ have similar expressions as \eqref{e:Z.drift.qv} above:
	\beq\label{e:Z.tr.drift.qv}
	\begin{aligned}
	dD^{\RomI,\trunc}(t)
	&= p'(t)^{1/2}\ab_{\trunc,t}\,dt\,,
		& \ab_{\trunc,t} &\equiv \ind\{t<q_T\} \ab_t\\
	dQ^{\RomII,\trunc}(t)
	&= p(t) \av_{\trunc,t}\,dt
		& \av_{\trunc,t} &\equiv \ind\{t<q_T\} \av_t\\
	dQ^{\RomII,\RomIII,\trunc}(t)
	&= [p(t) \cdot t p'(t)]^{1/2} \au_{\trunc,t}\,dt\,,
		& \au_{\trunc,t} &\equiv \ind\{t<q_T\} \au_t\\
	dQ^{\RomIII,\trunc}(t)
	&= \ar_{\trunc,t}\,dt,
		& \ar_{\trunc,t} &\equiv \ind\{t<q_T\} tp'(t)\,,\\
	dQ^{\Ising}(t)
	&= (\aw_{\trunc,t})^2\,dt,
		& \aw_{\trunc,t} &\equiv \ind\{t<q_{\Ising,T}\} \aw_t\\
	\end{aligned}\eeq
Note that this is analogous to \eqref{e:def.b}--\eqref{e:def.w}.
Similarly as in Definition~\ref{d:cts.occ}, let $\rho_\trunc$ be the occupation measure for the process $\vZ^\trunc$. Then let $\rho_{\trunc,D}$, $\rho_{\trunc,V}$, $\rho_{\trunc,U}$, $\rho_{\trunc,\rrr}$, $\rho_{\trunc,\www}$ be the occupation measures for 
$D^{\trunc,\RomI}$, $Q^{\trunc,\RomII}$, $Q^{\trunc,\RomII,\RomIII}$, $Q^{\trunc,\RomIII}$, and $Q^{\trunc,\Ising}$ respectively. As in Definition~\ref{d:cts.occ}, we can factorize the time-$t$ distributions as $\rho_{\trunc,t}=\varrho_{\trunc,t}\otimes\varsigma_{\trunc,t}$, $\rho_{\trunc,D,t}=\varrho_{\trunc,D,t}\otimes\varsigma_{\trunc,t}$, and so on.

\item \label{i:W2.nontrunc.to.trunc} 
Given the process $(Z,Z^{\Ising})$ and the error tolerance $\epsilon$, we claim that we can choose $\trK$ large enough that
\[
\bbW_2\Big(\Law(Z(1)),
	\Law(Z^{\trunc}(1)) \Big)
	+\bbW_2\Big(
	\Law(Z^{\trunc}(1),
	\Law(Z^{\frozen}(1)
	\Big) 
	\leq \frac{\epsilon}{2}\,.\]
Indeed, it follows from the budget constraints
\eqref{e:Z.budget.assumption} that $Z(1)$ and $Z^{\trunc}(1)$ have finite second moments. Moreover, we have $\lim_{\trK\to\infty} \bbP(Z(1)=Z^{\trunc}(1))=1$ because $\sup\{|Z_t|:0\le t\le1\}$ is almost surely finite. It follows by the monotone convergence theorem that as $\trK\to\infty$,
	\[
	\bE\bigg[ Z(1)^2 ; \sup\{|Z_t|:0\le t\le1\}\geq \trK\bigg]
	= o_{\trK}(1)\,,
	\]
which bounds the $\bbW_2$ error between $Z(1)$ and 
$Z^\trunc(1)$. The error between $Z^\trunc(1)$ and $Z^\frozen(1)$ can be handled very similarly. The analogous bound holds for the error between $Z^\Ising(1)=Z^{\trunc,\Ising}(1)$ and $Z^{\frozen,\Ising}(1)$. 

\item  Analogously to Definition~\ref{d:froze}, let $\vZ^\frozen$ be the process which equals $\vZ^\trunc$ with probability $1-\trK^{-10}$, and otherwise is frozen to a value chosen uniformly at random from the interval $[-2\trK,2\trK]^4$. Denote its drift and quadratic variation $\vD^\frozen$ and $\vQ^\frozen$. Define the corresponding occupation measures
$\rho_{\frozen,D}$, $\rho_{\frozen,V}$, $\rho_{\frozen,U}$, $\rho_{\frozen,\rrr}$, $\rho_{\frozen,\www}$. Again factorize $\rho_{\frozen,t}=\varrho_{\frozen,t}\otimes\varsigma_{\frozen,t}$, $\rho_{\frozen,D,t}=\varrho_{\frozen,D,t}\otimes\varsigma_{\frozen,t}$, and so on.
Note then that $\vD^\frozen$ and $\vQ^\frozen$ have similar expressions as \eqref{e:Z.tr.drift.qv} above, provided we take
	\beq\label{e:Z.frozen.versus.trunc}
	\ab_{\frozen,t}
	\equiv 
	\bigg(1-\frac{1}{\trK^{10}}\bigg)
	\frac{d\varrho_{\trunc,t}}{d\varrho_{\frozen,t}}
	\cdot
	\ab_{\trunc,t}
	\equiv 
	(1-p_{\frozen,t})\cdot
	\ab_{\trunc,t}
	\,,
	\eeq
and similarly for the other coefficients $(\av,\au,\ar,\aw)_{\frozen,t}$. Note that this is analogous to
\eqref{e:def.b.plus}.

\item We specify smoothed coefficients
as in Definitions \ref{d:t.smooth.variant} and \ref{d:cts.coeffs}. For example, we define the spatially smoothed coefficient
	\[
	\sB(t,y)
	= \frac{\overline{\varrho_{\frozen,D,t}*\sphi}(y)}
		{ \overline{\varrho_t*\sphi}(y) }
	\equiv p'(t)^{1/2} \cdot \sbbb(t,y)\,,
	\]
and from this we define the temporally smoothed variants 
$\tbbb$, $\tB$, $\Tbbb$, $\iB$. We denote as before 
$\tGam\equiv(\tbbb,\tvvv,\tuuu)$ and $\TGam\equiv(\Tbbb,\Tvvv,\Tuuu)$.

\item Recall the assumption that the original coefficients $(\ab,\av,\au,\aw)$ satisfy the budget constraint \eqref{e:Z.budget.assumption}. They also satisfy the domain constraint $\av-\au^2=\asig^2\ge0$. We argue that these constraints passes to the spatially smoothed quantities, in the following steps: 
\begin{itemize}
\item Comparison of
$(\ab,\av,\au,\aw)$ with $(\ab,\av,\au,\aw)_\trunc$: this step is analogous to the arguments for
 Proposition~\ref{p:budget-constraints-without-frozen}\ref{it:domain-constraint}-\ref{it:Ising-budget-constraint}. First, from \eqref{e:Z.tr.drift.qv}, it is immediate that $\av_\trunc$ and $\au_\trunc$ continue to satisfy the domain constraint:
	\[
	\av_{\trunc,t}-(\au_{\trunc,t})^2
	= \ind\{t<q_T\}\Big[ \av_t-(\au_t)^2\Big] \ge0\,.
	\]
We next argue that $(\ab,\av,\au,\aw)_\trunc$ approximately satisfies the budget constraint \eqref{e:Z.budget.assumption}: writing $\aGam\equiv(\ab,\av,\au)$ and 
$\aGam_\trunc\equiv(\ab,\av,\au)_\trunc$, we have
	\begin{align*}
	&\bE\Big[\Cost(\aGam_{\trunc,t})\Big]
	=\bE\bigg[ \Cost(\aGam_t)
		+\ind\{t\ge q_T\} \Big(
			\Cost(0,0,0)-\Cost(\aGam_t)\Big)\bigg]\\
	&\qquad\le\bE\Cost(\aGam_t)+\bP(t\ge q_T)
	\le \bE\Cost(\aGam_t)+o_{\trK}(1)\,.
	\end{align*}
Meanwhile, by the Cauchy--Schwarz inequality
and the assumption $\bE [(\aw_t)^2]=1$, we have
	\begin{align*}
	\bE \aw_t 
	&=\bE \aw_{\trunc,t}+ \bE( \ind\{t\ge q_T\} \aw_t )\\
	&\le\bE \aw_{\trunc,t}
	+ \bP(t\ge q_T) \bE [(\aw_t)^2]
	\le \bE \aw_{\trunc,t} + o_{\trK}(1)\,.
	\end{align*}
This proves that the coefficients $(\ab,\av,\au,\aw)_\trunc$ satisfy the budget constraint \eqref{e:Z.budget.assumption} up to additive error $o_{\trK}(1)$. 

\item Comparison of $(\ab,\av,\au,\aw)_\trunc$
with $(\ab,\av,\au,\aw)_\frozen$:
this step is analogous to the arguments for Proposition~\ref{p:budget-constraints-HIST}\ref{it:domain-constraint-plus-HIST}-\ref{it:Ising-budget-constraint-plus-HIST}. From \eqref{e:Z.frozen.versus.trunc} we obtain a modified domain constraint,
	\[
	\av_{\frozen,t}
	\ge \frac{(\au_{\frozen,t})^2}{1-p_{\frozen,t}}
	\ge (\au_{\frozen,t})^2\,.
	\]
For the budget constraint, denoting $\aGam_{\frozen}\equiv(\ab,\av,\au)_{\frozen}$, we use Jensen's inequality to bound
	\begin{align*}
	&\bE\Big[\Cost(\aGam_{\frozen,t})\Big]
	=\int \Cost\Big( 
	(1-p_{\frozen}(t,\vy))
	\aGam_\trunc(t,\vy)
	+p_{\frozen}(t,\vy)
	(0,0,0)
	\Big) 
	\,\varrho_{\frozen,t}(d\vy)
	\\
	&\qquad\le \int
	\bigg[ (1-p_{\frozen}(t,\vy))
		\Cost(\aGam_{\trunc}(t,\vy))
	+p_{\frozen}(t,\vy)\Cost (0,0,0) \bigg]
	\,\varrho_{\frozen,t}(d\vy)\\
	&\qquad=
	\bigg(1-\frac{1}{\trK^{10}}\bigg)
	\int \Cost(\aGam_{\trunc}(t,\vy))
	\varrho_{\trunc,t}(d\vy)
	+ \frac{1}{\trK^{10}}
	\le \bE\Big[\Cost(\aGam_{\trunc,t})\Big] + o_{\trK}(1)\,.
	\end{align*}
Similarly, we have
	\begin{align*}
	\bE \aw_{\frozen,t}
	&= \int (1-p_{\frozen}(t,\vy))
	 \aw_{\trunc}(t,\vy)\,\varrho_{\frozen,t}(d\vy)\\
	&= \bigg(1-\frac{1}{\trK^{10}}\bigg)
	\int\aw_{\trunc}(t,\vy)\,\varrho_{\trunc,t}(d\vy)
	= \bigg(1-\frac{1}{\trK^{10}}\bigg)
	\bE \aw_{\frozen,t}
	\end{align*}
This proves that the coefficients $(\ab,\av,\au,\aw)_\frozen$ also satisfy the budget constraint \eqref{e:Z.budget.assumption} up to additive error $o_{\trK}(1)$:
	\beq\label{e:Z.budget.frozen.approx}
	\bE\Big[\Cost(\aGam_{\frozen,t})\Big]
	\le\frac{( \bE \aw_{\frozen,t})^2}{\alpha} + o_{\trK}(1)\,.
	\eeq 

\item Comparison of
$(\ab,\av,\au,\aw)_\frozen$ with $(\sbbb,\svvv,\suuu,\swww)$:
this can be done using the arguments of Lemmas \ref{l:sp.smooth.domain} and \ref{l:sp.smooth.budget}. From this we can conclude
	\[
	\tilde{C}(t) - \frac{\tilde{W}(t)^2}{\alpha}
	\le 3\epsilon^\circ
	\]
for all $q_0 \le t \le 1$.
\end{itemize}

\item Some of the following steps will require some analogues of the \textit{a~priori} estimates from \S\ref{ss:apriori}. In the current setting, the required analogues will all be derived from the budget constraint \eqref{e:Z.budget.assumption}, together with the assumption 
$\bE[(\aw_t)^2]=1$. We first note that, by calculus, for $v\ge u^2$ we have
	\[
	\Cost(b,v,u)
	=b^2 +v -2(v-u^2)^{1/2}+1
	\ge b^2 + v - 2v^{1/2}+1
	\ge b^2 + \frac{v}{2}-1\,.
	\]
Combined with the approximate budget
\eqref{e:Z.budget.frozen.approx} derived above, it implies
	\beq\label{e:Z.apriori}
	\begin{aligned}
	&\bE\bigg[ (\ab_{\frozen,t})^2 +
		\frac{(\au_{\frozen,t})^2}{2} \bigg]
	\le 
	\bE\bigg[ (\ab_{\frozen,t})^2 + \frac{\av_{\frozen,t}}{2}\bigg]\\
	&\qquad\le \bE\Big[ \Cost(\ab_{\frozen,t},\av_{\frozen,t},\au_{\frozen,t})
		\Big]+1
	\le \frac{(\bE\aw_{\frozen,t})^2}{\alpha} + 1 + o_{\trK}(1)
		\le \frac{1}{\alpha} + 1 + o_{\trK}(1)\,.
	\end{aligned}
	\eeq
In the current setting, the estimate \eqref{e:Z.apriori} can be used in place of Lemmas \ref{l:barD.I.second.quenched} and \ref{l:tQ.first.quenched}. Meanwhile, analogues of the Kolmogorov estimates of Propositions \ref{p:Y.kolmogorov} and \ref{p:Y.notI.second.quenched} can be derived 
from the stochastic integral representation of $\vZ^\frozen$, together with \eqref{e:Z.apriori}. For example, the stochastic integral representation of $Z^{\frozen,\RomII}$, combined with \eqref{e:Z.apriori}, gives
	\beq\label{e:Z.holder}
	\bE\bigg[\Big( Z^{\frozen,\RomII}(t)-Z^{\frozen,\RomII}(s)\Big)^2\bigg]
	= \bE \int_s^t p(r) \av_{\frozen,r} \,dr
	\le O(1) (t-s)\,.
	\eeq
Similar estimates apply for the other coordinates.

\item \label{i:Z.regularity} The smoothed coefficients also satisfy regularity estimates:
\begin{itemize}
\item The spatially smoothed coefficients satisfy regularity estimates as in Lemmas \ref{l:sp.smooth.reg}--\ref{l:sp.smooth.reg.SDE}, using
\eqref{e:Z.apriori} in place of Lemmas \ref{l:barD.I.second.quenched} and \ref{l:tQ.first.quenched}.

\item The spatially smoothed occupation densities are also regular in time in the manner of Lemma~\ref{l:sp.smooth.cty}, using
\eqref{e:Z.holder} (and its analogues for the other coordinates) in place of Propositions \ref{p:Y.kolmogorov} and \ref{p:Y.notI.second.quenched}.

\item The temporally smoothed coefficients satisfy regularity estimates by the argument of Corollary~\ref{c:t.smooth.regularity.N}.

\end{itemize}

\item \label{i:Z.tsmooth.domain.budget}
The above domain and budget constraints can then be passed to the temporally smoothed coefficients
$(\tbbb,\tvvv,\tuuu,\twww)$, adapting the arguments of Lemmas \ref{l:t.smooth.domain} and \ref{l:t.smooth.budget}. In the current setting, an important simplification is the regularity assumption on $p$, which allows us to take $J=\emptyset$ in \eqref{e:bad.times.J}. It follows
from the argument of
Lemma~\ref{l:t.smooth.domain} that
$\tGam\equiv(\tbbb,\tvvv,\tuuu)$ and $\TGam\equiv(\Tbbb,\Tvvv,\Tuuu)$ are uniformly close over all $(t,y)\in[q_*,1]\times\R$. One notable difference is that in the original application of Lemma~ \ref{l:t.smooth.domain}, towards proving  Theorem~\ref{thm:BOGP-hardness-main}, we had $p(q_0)=0$, and chose $q_*$ according to \eqref{e:breve.q.bounds}--\eqref{e:qbar.and.qstar} to ensure the bound \eqref{e:sde.lbd.on.p.N.assump}, which was used in the analysis of Lemma~ \ref{l:t.smooth.domain}. In the setting of the current theorem, we already assume $p(q_0)\ge 1/L$, so the bound  \eqref{e:sde.lbd.on.p.N.assump} is already satisfied without any further constraints on $q_*$. Next, analogously to 
\eqref{e:cost.temp.smoothed.LIMIT}
and \eqref{e:budget.temp.smoothed.LIMIT}, denote
	\begin{align*}
	\bar{C}(t)
	&\equiv\int_{\R} 
	\Cost\Big(\tGam(t,y)\Big)
	\overline{\varrho_t*\sphi(y)}
	\,dy\,,\\
	\bar{W}(t)
	&\equiv
	\int_{\R} \twww(t,y)^{1/2} 
		\varsigma_t*\sphi(y)\,dy\,.
	\end{align*}
The argument of Lemma~\ref{l:t.smooth.budget} then gives
	\[
	\bar{C}(t) - \frac{\bar{W}(t)^2}{\alpha}
	\le 4\epsilon^\circ
	\]
for all times $t\in[q_*,1]$.
At this point we have obtained analogues of all the results of
\S\ref{ss:discrete.coeffs} for the setting of the current theorem.

\item Since in the current setting there is no dependence on $N$, we do not require any analogues of  the results from
\S\ref{ss:limit.coefs}:
\begin{itemize}
\item Lemmas \ref{l:occ.wk.conv} and \ref{l:limit.coefs} hold trivially.
\item
Corollary~\ref{c:reg.in.limit} is redundant with step~\eqref{i:Z.regularity} above.
\item Lemma \ref{l:bad.times.JN.J} holds trivially.
\item Corollaries \ref{c:limit.coefs.variants.of.temp.smoothing} and \ref{c:limit.budget} are redundant with step \eqref{i:Z.tsmooth.domain.budget} above.
\end{itemize}

\item We next note that we have straightforward analogues of the results of \S\ref{ss:limit.sde}:
\begin{itemize}
\item Let $\hvZ^\frozen(t)\equiv \vZ^\frozen(t-\sss)+\xxx$, and define the coordinate sum $\hat{Z}^\frozen$ analogously to
\eqref{e:iY.plus}. Define $\iX^\frozen$ and $\iX^{\frozen,\Ising}$ to be the solutions of the SDE with ``ideal'' coefficients. Then $\iZ^\frozen$ and $\iX^\frozen$ satisfy the same Fokker--Planck equations on $[q_*,1]$, and likewise for
$\iZ^{\frozen,\Ising}$ and $\iX^{\frozen,\Ising}$. 
This is the analogue of Proposition~\ref{p:fokker.planck},
using step~\eqref{i:Z.regularity} in place of 
Corollary~\ref{c:reg.in.limit}.

\item Lemma~\ref{l:bad.J} is not needed since $J=\emptyset$ as noted above.

\item The discrepancy between the temporally smoothed coefficients and the ``ideal'' coefficients can be bounded as in Lemma~\ref{l:ideal.vs.actual.coeffs}, again using step~\eqref{i:Z.regularity} in place of Corollary~\ref{c:reg.in.limit}.

\item Lemma~\ref{l:gronwall} is simply a general statement of Gr\"onwall's inequality.

\item An analogue of Proposition~\ref{p:sde.approx} follows by collecting the above: if $X^\frozen$ is the solution of the actual SDE, then it is close to $\iX^\frozen$ by the analogue of Lemma~\ref{l:ideal.vs.actual.coeffs},
along with the Gr\"onwall inequality from Lemma~\ref{l:gronwall}.
Then $\iX^\frozen$ and $\iZ^\frozen$ have matching marginals at any fixed time by the analogue of Proposition~\ref{p:fokker.planck}.
Lastly $\iZ^\frozen$ is close to $Z^\frozen$. Here we again use step~\eqref{i:Z.regularity} in place of 
Corollary~\ref{c:reg.in.limit}, and use \eqref{e:Z.holder} (and its analogues for the other coordinates) in place of Propositions \ref{p:Y.kolmogorov} and \ref{p:Y.notI.second.quenched}.

\end{itemize}

\item We next collect the analogues of the results of \S\ref{ss:sde.reparam.sigma} leading to the  \hyperlink{proof:t.BOGP-hardness-main}{proof of Theorem~\ref{thm:BOGP-hardness-main}}:
\begin{itemize}
\item
Proposition~\ref{ppn:budget-constraint-reparametrize} is a general calculus result.

\item The error $tp'(t)-\trrr_t$ can be bounded by a similar argument as for
Lemma~\ref{lem:trrr-error-small}: let $\bar{P}(t)$ be the temporal smoothing of $tp'(t)$, and note that $tp'(t)$ is uniformly close to $\bar{P}(t)$, since $J=\emptyset$. We therefore turn to bounding $\bar{A}_t=\bar{P}(t)-\trrr_t$, which we regard as the spatio-temporal smoothing of 
	\[A_t=tp'(t) -\ar_{\frozen,t}
	=p_{\frozen,t} tp'(t)+ (1-p_{\frozen,t}) (tp'(t)-\ar_{\trunc,t} )\,.
	\]
For the first term it suffices to recall that
	\[\bE p_{\frozen,t} = \frac1{\trK^{10}} = o_{\trK}(1)\,.\] 
For the second term it suffices to bound
	\[\bE\Big[|tp'(t)-\ar_{\trunc,t}|\Big]
	= tp'(t) \bP(t\le q_T) \le o_{\trK}(1)\,.\]
The rest of the argument of Lemma~\ref{lem:trrr-error-small} goes through similarly. Since $J=\emptyset$, the final bound holds uniformly over all $t\in[q_*,1]$.

\item 
An analogue of Theorem~\ref{thm:Lipschitz-SDE-approx} follows by collecting the above. An important difference again is that, since $J=\emptyset$ and Lemma~\ref{lem:trrr-error-small} holds uniformly in time, we can improve the budget constraint \eqref{it:BOGP-hardness-main-budget} to hold at each time $t$ rather than only after integrating over time.  This gives SDE approximations to $(Z^\frozen,Z^{\frozen,\Ising})$. We denote these SDE approximations as $(X,X^\Ising)$.
\end{itemize}
\end{enumerate}
The result of Theorem~\ref{thm:SDE-smoothing-general} now follows:
\begin{itemize}
\item conclusion~\eqref{it:SDE-smoothing-general.coefs}
follows by the analogue of Theorem~\ref{thm:Lipschitz-SDE-approx}\ref{it:Lipschitz-SDE-approx.Lip};
\item conclusion~\eqref{it:SDE-smoothing-general.OUTPUT.diffus}
follows by the analogue of Theorem~\ref{thm:Lipschitz-SDE-approx}\ref{thm:Lipschitz-SDE-approx.b};
\item conclusion~\eqref{it:SDE-smoothing-general.budget}
follows by the analogue of Theorem~\ref{thm:Lipschitz-SDE-approx}\ref{thm:Lipschitz-SDE-approx.a} --- as discussed above, the bound now holds uniformly over $t\in[q_*,1]$ because we can take $\oJ=\emptyset$;
\item conclusion~\eqref{it:SDE-smoothing-general.initial}
follows by the analogue of Theorem~\ref{thm:Lipschitz-SDE-approx}\ref{thm:Lipschitz-SDE-approx.c.INIT},
using step \eqref{i:W2.nontrunc.to.trunc}  to handle the error between $(Z(q_0),Z^\Ising(q_0))$ and $(Z^\frozen(q_0),Z^{\frozen,\Ising}(q_0))$;
\item conclusion~\eqref{it:SDE-smoothing-general.OUTPUT.qstar}
follows from the choice of $q_*$ discussed above.
\end{itemize}
The $\bbW_2$ approximation estimates 
follow from the corresponding estimates given by 
Theorem~\ref{thm:Lipschitz-SDE-approx},  using step \eqref{i:W2.nontrunc.to.trunc}  to handle the error  between $(Z(1),Z^\Ising(1))$ and $(Z^\frozen(1),Z^{\frozen,\Ising}(1))$.  The final assertion --- that $w_t$ can be taken even if $\Law(Z^{\Ising}(t))$ is even for each $t$ --- holds by symmetry, as we may assume the mollifier $\sphi$ used for spatial smoothing is also an even function.
\end{proof}

\newcommand{\chaos}{{\textup{ch}}}
\newcommand{\linit}{{\ell_{\init}}}
\newcommand{\init}{\textup{init}}
%PROXY
\newcommand{\prxby}{\textup{\texttt{y}}}
\newcommand{\prxbx}{\textup{\texttt{x}}}
\newcommand{\prxX}{\textup{\texttt{X}}}
\newcommand{\prxY}{\textup{\texttt{Y}}}
\newcommand{\prxgamma}{\Gamma}
\newcommand{\prxpsi}{\uppsi}
\newcommand{\upp}{\textup{\texttt{p}}}

\def\I{{\rm I}}
\def\II{{\rm II}}
\def\III{{\rm III}}

\pagebreak\section{Matching IAMP algorithms}
\label{sec:IAMP}

\iffull
% !TEX root = main.tex

In this section we present Lipschitz \textbf{incremental approximate message passing (IAMP)} algorithms which simulate any SDE of the form we have identified. The section is organized as follows: 
\begin{itemize}
\item In \S\ref{subsec:state-evolution-multiple} we review known results (\cite{gerbelot2021graph}) on the state evolution limit of a sufficiently general form of AMP algorithms.
\item In \S\ref{subsec:IAMP-setup} we introduce an explicit IAMP iteration
(Definition~\ref{d:iamp.v1}) which is designed to mimic the SDE coming from the BOGP analysis. This forms the basis of the main algorithmic result of this section, Theorem~\ref{thm:IAMP-main}\ref{i:IAMP-main-main}, which gives an IAMP-based Lipschitz algorithm that mimics the SDE appearing in our hardness result Theorem~\ref{thm:BOGP-hardness-main}. 
\item In \S\ref{ss:proof.iamp.diffus.limit} we present the \hyperlink{proof:t.IAMP.main.main}{proof of Theorem~\ref{thm:IAMP-main}\ref{i:IAMP-main-main}}, based on an analysis of the IAMP iteration of Definition~\ref{d:iamp.v1}.
\item In \S\ref{subsec:center-IAMP} we present the \hyperlink{proof:t.IAMP.main.centered}{proof of Theorem~\ref{thm:IAMP-main}\ref{i:IAMP-main-centered}}. This theorem states that, under additional symmetry assumptions on the SDE's control processes, the desired behavior of the IAMP algorithm can be achieved while enforcing the additional constraint that the algorithm's output has mean $\bzero\in\R^N$.

\item In \S\ref{subsec:p=1-IAMP} we present the \hyperlink{proof:t.IAMP.main.chaotic}{proof of Theorem~\ref{thm:IAMP-main}\ref{i:IAMP-main-chaotic}}. This theorem states that any feasible \textbf{symmetrized} endpoint measure can be achieved by a ``chaotic IAMP'' algorithm $\cA_N$ whose correlation function (recall \eqref{e:p.corr.overlap}) satisfies $\chi_{\cA_N}(1-\iota) \le \iota$. 

\item Lastly, in \S\ref{ss:even.solution.trees} we show that in the symmetrized setting, chaotic IAMP algorithms are able to construct large, approximately ultrametric solution trees, matching the BOGP intuition. 
\end{itemize}
The following is the main result of this section.
\begin{thm}\label{thm:IAMP-main}
  For any $\iota > 0$, there exists $\epsilon_0 = \epsilon_0(\alpha,\iota) > 0$ such that for all $\epsilon \in (0,\epsilon_0)$ and $L > 0$, there exists $C(L,\epsilon)$ such that the following holds.
  Consider any $q_0 \in [0,1)$. Given controls $(b,\sigma,w,p,\zeta,\zeta^\Ising)$,
recall from \eqref{e:s.sqrt.fn} that we denote $s(t) = ((tp)'(t))^{1/2}$. Let $B$ and $W$ be standard Brownian motions on the time interval $q_0 \le t\le 1$, and consider the SDEs (cf.\ Theorem~\ref{thm:BOGP-hardness-main})   \begin{align}
    \label{e:IAMP.X.SDE}
    \de X(t) &= p'(t)^{1/2} b(t,X(t)) \,\de t + s(t) \sigma(t,X(t)) \,\de B(t)\,, \\
    \label{e:IAMP.Y.SDE}
    \de Y(t) &= w(t,Y(t)) \,\de W(t)\,,
  \end{align}
with initial values $X(q_0) \sim \zeta$ and $Y(q_0) \sim \zeta^\Ising$ independent of $B,W$. Suppose the controls satisfy the following conditions:
  \begin{enumerate}[(i)]
    \item \label{i:IAMP-main-coefs} $b(t,x) : [q_0,1] \times \bbR \to [-L,L]$ and $\sigma(t,x), w(t,x) : [q_0,1] \times \bbR \to [0,L]$ are $L$-Lipschitz in $(t,x)$.
    \item \label{i:IAMP-main-diffusivity} $|\bbE[(w_t)^2] -1| \le \epsilon$ for all $t\in [q_0,1]$.
    \item \label{i:IAMP-main-p} $p\in \incr([q_0,1];[0,1])$ (see Definition~\ref{d:incr.p}) is a twice differentiable function, with $p(q_0) \in [1/L,\epsilon]$ and $\|p\|_{C^2([q_0,1])} \le L$.
    \item \label{i:IAMP-main-budget} For all $t\in [q_0,1]$, the controls satisfy the approximate budget constraint:
    \begin{align}\nonumber
      \Cost(t)
      &\equiv
      \bbE\bigg[
      b_t(X(t))^2 + \bigg(
        \frac{s(t)^2}{p(t)}
        \bigg) \big(\sigma_t(X(t))-1\big)^2
      \bigg]\\
      &\leq 
      \budget(t)+\epsilon
      \equiv
      \frac{(\bbE w_t(Y(t)))^2}
        {\alpha}
      +\epsilon\,.\label{eq:budget-constraint-ising}
    \end{align}
    \item \label{i:IAMP-main-endpt} The endpoint of the $Y$ process satisfies
    \beq\label{eq:endpoint-constraint-IAMP}
    \bbE\Big[\big(|Y(1)|-1\big)^2\Big]\leq \epsilon\,.
    \eeq
    \item \label{i:IAMP-main-init} The initial measures $\zeta, \zeta^\Ising \in\cP_2(\R)$ are such that $|\bbE[X(q_0)^2]|\le \epsilon$ and $|\bbE[Y(q_0)^2] - q_0|\le \epsilon$.
  \end{enumerate}
We then have the following:
  \begin{enumerate}[(a)]
    \item \label{i:IAMP-main-main} There exists a $C(L,\epsilon)$-Lipschitz algorithm $\cA_N$ (in the sense of Definition~\ref{d:Lip}) which agrees with the output of an efficient IAMP iteration with probability $1-e^{-cN}$, such that
    \begin{align}
      \label{e:IAMP-main-main-coords}
      \bbP\lt(\cA_N(\bG,\bg^\aux) \in \Sigma(\iota) \rt) \ge 1-e^{-cN}\,, \\
      \label{e:IAMP-main-main-inner-prods}
      \bbP\lt(
        \bbW_2\Big(\mu_\bG(\cA_N(\bG,\bg^\aux)), \Law(X(1)) \Big) \le \iota
      \rt) \ge 1-e^{-cN}\,.
    \end{align}
    Moreover, for $\mu(\cA_N)$ and $\mu^{\Ising}(\cA_N)$ defined in \eqref{e:mu.of.Alg} and \eqref{e:mu.Ising.of.Alg}, we have
    \begin{align}
      \label{e:IAMP-main-main-inner-prods-averaged}
      \bbW_2\Big(\mu(\cA_N), \Law(X(1)) \Big) \le \iota\,, \\
      \label{e:IAMP-main-main-coords-averaged}
      \bbW_2\Big(\mu^{\Ising}(\cA_N), \cP(\{\pm 1\}) \Big) \le \iota\,,
    \end{align}
    where in \eqref{e:IAMP-main-main-coords-averaged}, $\bbW_2$ denotes the point-to-set $\bbW_2$ distance.
    \item \label{i:IAMP-main-centered} Suppose that, in addition to the assumptions (\ref{i:IAMP-main-coefs})--(\ref{i:IAMP-main-init}) above, we have $w(t,x) = w(t,-x)$ for all $t,x$ and $\zeta^\Ising$ is symmetric, i.e. $\Law(-x : x\sim \zeta^\Ising) = \zeta^\Ising$. Then, there exists a $C(L,\epsilon)$-Lipschitz algorithm $\cA_N$ satisfying the conclusions of part (\ref{i:IAMP-main-main}) and furthermore $\bbE \cA_N(\bG,\bg^\aux) = \bzero$.
    \item \label{i:IAMP-main-chaotic} 
    Under only the assumptions (\ref{i:IAMP-main-coefs})--(\ref{i:IAMP-main-init}) above, 
    there exists a $C(L,\epsilon)$-Lipschitz algorithm $\cA_N$, which agrees with the output of an efficient IAMP iteration with probability $1-e^{-cN}$, such that \eqref{e:IAMP-main-main-coords} holds,
    \beq\label{e:IAMP-main-main-inner-prods-sym}
      \bbP\Big(
        \bbW_2\lt(\mu_{\bG,\sym}(\cA_N(\bG,\bg^\aux)), \Law(|X(1)|) \rt) \le \iota
      \Big) \ge 1-e^{-cN}\,,
    \eeq
    and $\chi_{\cA_N}$ defined in \eqref{e:p.corr.overlap} satisfies $\chi_{\cA_N}(0) = 0$ and $\chi_{\cA_N}(1-\iota) \le \iota$.
    Moreover, \eqref{e:IAMP-main-main-coords-averaged} holds, and
    \beq\label{e:IAMP-main-main-inner-prods-sym-averaged}
      \bbW_2\Big(\sym(\mu(\cA_N)), \Law(|X(1)|) \Big) \le \iota\,,
    \eeq
    where $\sym(\cdot)$ is defined in \eqref{e:symmetrize.measure}.
  \end{enumerate}
\end{thm}
 
As noted at the beginning of this section, the \hyperlink{proof:t.IAMP.main.main}{proof of Theorem~\ref{thm:IAMP-main}\ref{i:IAMP-main-main}} is given at the end of \S\ref{ss:proof.iamp.diffus.limit}; the \hyperlink{proof:t.IAMP.main.centered}{proof of Theorem~\ref{thm:IAMP-main}\ref{i:IAMP-main-centered}}
is given at the end of \S\ref{subsec:center-IAMP}; and the
\hyperlink{proof:t.IAMP.main.chaotic}{proof of Theorem~\ref{thm:IAMP-main}\ref{i:IAMP-main-chaotic}}
is given at the end of \S\ref{subsec:p=1-IAMP}.

\begin{rmk}[order of parameters]
\label{r:Lip.epsilon.discussion}
  All the arguments in this section will observe the order of parameters
  \beq
  \label{eq:IAMP-order-of-limits}
  0\ll \delta\ll \frac1L
    \ll \epsilon\ll \iota \ll 1\,,
  \eeq
  where $\delta$ is the time discretization of the IAMP iteration described in \S\ref{subsec:IAMP-setup}. This should not be confused with the parameters of Assumption~\ref{a:params}, which only applies to Sections \ref{s:rerand} and \ref{s:sde}.  Note that this is consistent with the order of quantifiers in Theorem~\ref{thm:IAMP-main}.
  In Theorem~\ref{thm:IAMP-main}, we only assume the controls $(b,\sigma,w,p,\zeta,\zeta^\Ising)$ satisfy (\ref{i:IAMP-main-coefs})--(\ref{i:IAMP-main-init}) with an error tolerance $\epsilon > 0$, so that there exists a Lipschitz $w$ satisfying assumption (\ref{i:IAMP-main-endpt}). The Lipschitz property of $b,\sigma,w$ is in turn necessary to make the SDE solutions of \eqref{e:IAMP.X.SDE}, \eqref{e:IAMP.Y.SDE} well-defined, and to analyze our IAMP iteration using state evolution.
  In Section~\ref{sec:alternate-diffusions}, we will show that setting $\epsilon=0$ does not make a difference in the limiting stochastic control problems when progressively measurable $b,\sigma,w$ are permitted.
\end{rmk}

\begin{rmk}[spherical vs.\ Ising]\label{r:IAMP.spherical}
  Theorem~\ref{thm:IAMP-main} is written for the Ising perceptron.
  The spherical case is subsumed in our analysis by fixing $w_t\equiv 1$, so that (\ref{i:IAMP-main-diffusivity}) trivially holds, and ignoring (\ref{i:IAMP-main-endpt}).
  As the analysis is the same, we focus on the Ising case below. 
\end{rmk}

\begin{rmk}[external randomness] Our IAMP algorithms will depend on many external sources of randomness, including $\bG_{\circ}$.
We sometimes view these random variables as a random seed $\omega$ in the sense of Remark~\ref{rmk:omega-general-seed}, and sometimes as part of the gaussian vector $\bg^{\aux}$ as in \eqref{eq:algorithms-as-maps}.
We will always make explicit which of these views we take; see e.g.\ \eqref{e:iamp.random.seed} or \eqref{e:iamp.chaotic.omega}.
\end{rmk}

Previous IAMP algorithms \cite{mon18,ams20,alaoui2022perceptron,sellke2021optimizing,montanari2022overparametrized,huang2024optimization,montanari2024exceptional} took the form of first-order iterations where vector-valued iterates were multiplied by the input matrix $\bG$.
These algorithms were restricted to simulating martingale SDEs.\footnote{Several of these algorithms include an initial ``root-finding'' phase using a simpler non-incremental form of AMP. This combination is also expected to constitute an optimal Lipschitz algorithm; \cite{sellke2021optimizing,huang2024optimization} essentially prove this in their respective settings, but it is non-obvious in general. The intuitive picture is that similarly to the statics, the ``algorithmically accessible maxima'' are arranged into an approximately ultrametric tree, which may have a non-trivial root at positive distance from the origin. This picture suggests that once the root is located, a martingale-type AMP suffices to descend the tree.}
In order to simulate the more general semimartingale diffusions \eqref{e:IAMP.X.SDE}, \eqref{e:IAMP.Y.SDE} within the IAMP framework, our IAMP algorithm multiplies its iterates by a correlated family of matrices.
More precisely, we let $(\bG(t) : 0\le t\le 1)$ be an $M\times N$ matrix-valued Brownian motion conditioned to satisfy $\bG(1)=\bG$. Treating $\bG$ as fixed, this means
\beq\label{e:matrix.br.bridge.IAMP}
\bG(t) = t\bG + \bG_{\circ}(t)
\eeq
where $\bG_{\circ}$ is an $M\times N$ matrix-valued Brownian bridge. (That is, the entries of $\bG_{\circ}$ are independent, and each entry of $\bG_{\circ}$ is a Brownian motion started from zero and conditioned to return to zero at time $t=1$.) In our IAMP algorithm, the $\ell$-th iteration multiplies by $\bG(t_\ell)$, for a suitably chosen increasing sequence $t_\ell \in [0,1]$.
The analysis of this algorithm requires a sufficiently general version of state AMP state evolution, which we review in \S\ref{subsec:state-evolution-multiple}.

While our main algorithmic result is Theorem~\ref{thm:IAMP-main}\ref{i:IAMP-main-main}, parts \ref{i:IAMP-main-centered}--\ref{i:IAMP-main-chaotic} of this theorem will substantially simplify the final stochastic control description we obtain in Section~\ref{sec:alternate-diffusions}.
In particular, part~\ref{i:IAMP-main-centered} is used in the proof of Lemma~\ref{l:IAMP.in.Lip}, and will imply that setting $q_0 = 0$ suffices to attain all feasible endpoint measures.
Part~\ref{i:IAMP-main-chaotic} is used in the proof of Lemma~\ref{l:IAMP.in.Lip.sym}, and will imply that setting $p\equiv1$ suffices to achieve any \textbf{symmetrized} endpoint measure. The latter implication is analogous to prior results from \cite{auffinger2015properties,ams20,huang2023algorithmic} on the behavior of the Parisi formula and algorithmic thresholds in mean-field spin glasses, shown using variational calculus.
Here we take a different approach: Theorem~\ref{thm:IAMP-main}\ref{i:IAMP-main-chaotic} will imply that any Lipschitz algorithm can be simulated by a ``chaotic IAMP'' algorithm whose intrinsic correlation function satisfies $\chi(1-\iota)\le \iota$.
This will mean that $p(\iota)\geq 1-\iota$, which for small $\iota$ approximates the constant function $p\equiv1$. 

\subsection{AMP state evolution for multiple correlated matrices}
\label{subsec:state-evolution-multiple}

We use the results of \cite{gerbelot2021graph}, which considered AMP algorithms on constant-sized directed graphs. Their general setup is that each vertex $v$ in the graph is associated with a vector space $U_v$, with dimension $N_v$ such that $\lim_{N\to\infty} N_v/N=\alpha_v\in (0,\infty)$.
Each directed edge $v\to w$ is associated with a random linear map $U_v\to U_w$, see 
\cite[\S2]{gerbelot2021graph}.
The basic case of a single symmetric matrix is a vertex with self-loop, while an asymmetric $M\times N$ matrix $\bG$ is a two-vertex graph with a single bi-directional edge, corresponding to $\bG$ in one orientation and $\bG^{\top}$ in the other.

Our IAMP algorithm will produce iterates $\bu^{\ell,k}\in\bbR^N$ and $\bv^{\ell,k}\in\bbR^M$ for $0\leq k\leq\ell\leq\ubl$. It involves $\ubl$ disorder matrices $\bJ_k$, where $1\le k\le \ubl$, which are i.i.d. samples from $\cN(0, I_{M\times N}/N)$ (which represent rescaled increments of the matrix Brownian motion $\bG(t)$, see \eqref{eq:Jk}).
In the setting of \cite{gerbelot2021graph}, this corresponds to a multi-graph on two vertices, labeled by $\bbR^M$ and $\bbR^N$, with $\ubl$ bi-directional edges between them each associated to a different $\bJ_k$.
We also let $\bar{\bg}^1,\dots\bar{\bg}^{\ubl}\in\R^N$ be independent random vectors (usually gaussian). Lastly, we let $\bx[0]\in\R^{M\times d_x}$ and $\by[0]\in\R^{N\times d_y}$ be auxiliary random initial data, for constant $d_x,d_y$, such that
\begin{equation}\begin{aligned}\label{e:iamp-side-information-w2-convergence}
  \frac{1}{M} \sum_{m=1}^M \delta(x[0]_m) &\stackrel{\bbW_2}{\longrightarrow} \cL(X[0]) \in \cP(\bbR^{d_x})\,, \\
  \frac{1}{N} \sum_{n=1}^N \delta(y[0]_n) &\stackrel{\bbW_2}{\longrightarrow} \cL(Y[0]) \in \cP(\bbR^{d_y})
\end{aligned}\end{equation}
in probability. 

Let $\bu[[\ell]]$ denote the collection of vectors $\bu^{i,j}$ indexed by $ 0\le j\le i\le\ell$. We also let $\bar{\bg}[\ell]$ denote the collection of vectors $\bar{\bg}^i$ indexed by $0\le i\le\ell$. Then, for our purposes, a sufficiently general AMP algorithm is given by the recursive definition
\beq
\label{eq:IAMP-fully-general-fixed}
\begin{aligned}
  \bu^{\ell+1,k}
  &= (\bJ_k)^{\top} 
  f_{\ell,k}\big(\bv[[\ell]]
  ,\bx[0]
  \big)
  +
  \sum_{i=0}^\ell
  a_{i,\ell,k}\bar{\bg}^i
  -\ons_{u,\ell,k}\in\bbR^N\,,
  \quad
  0 \le k\le \ell+1\,,
  \\
  \bv^{\ell,k}
  &=
  \bJ_{k}
  h_{\ell,k}
  \Big( \bu[[\ell]], \by[0]
  \Big)
  -\ons_{v,\ell-1,k}\in\bbR^M\,,
  \quad 
  0\le k\le \ell\,.
  \end{aligned}
\eeq
Here the functions $f_{\ell,k}$ and $h_{\ell,k}$ are applied coordinatewise,
i.e., the $n$-th coordinate of $f_{\ell,k}(\bv[[\ell]],\bx[0])$ is $f_{\ell,k}$ applied to the $n$-th coordinates of each of the vectors
in $\bv[[\ell]]$, and to each of the columns of $\bx[0]$. The Onsager correction terms $\ons_{u,\ell,k}$ and $\ons_{v,\ell-1,k}$ will be defined in \eqref{eq:IAMP-fully-general-fixed.Ons} below. 

\begin{rmk}\label{rem:external-noise-is-ok}
Technically, \cite{gerbelot2021graph} considers only AMP iterations without long-term memory, meaning that $f_{\ell}$ depends only on $\bv^{\ell}$. However, the extension to AMP iterations with long-term memory
is encompassed by matrix-valued iterates as used in \cite[\S3.3]{gerbelot2021graph}, and thus follows by the same proof. The exogenous randomness from the $\bar{\bg}^\ell,\bx[0],\by[0]$ is also not explicitly handled in \cite{gerbelot2021graph}, but also has exactly the same proof (see e.g.\ \cite[Appendix A]{huang2024optimization} for a similar extension to exogenous randomness for tensor AMP).
\end{rmk}

The corresponding \textbf{state evolution} description is in terms of centered jointly gaussian  random vectors $\vec {\bar{G}} = \bar{G}[\ubl]$, $\vec U = U[[\ubl]]$, $\vec V = V[[\ubl]]$, as well as initial data $X[0]\in\R^{d_x},Y[0]\in\R^{d_y}$. (Similarly to the notation used above, $\bar{G}[\ell]$ denotes the random variables $\bar{G}^i$ indexed by $0\le i\le\ell$, and $U[[\ell]]$ denotes the random variables $U^{i,j}$ indexed by $0\le j\le i \le\ell)$.) First of all, we can view $X[0]$, $Y[0]$, $(\vec{U},\vec{\bar{G}})$, and $\vec{V}$ as mutually independent. The random variables $(\vec{U},\vec{\bar{G}})$ are jointly gaussian, as are the random variables $\vec{V}$. The $\bar{G}^i$ are i.i.d.\ one-dimensional standard gaussian random variables. Beyond this, the covariances are defined recursively by the equations
\begin{align}\nonumber
  \bbE[U^{\ell_1+1,k_1}U^{\ell_2+1,k_2}]
  &= \ind\{k_1=k_2\}\bigg\{
  \alpha\,
  \bbE 
  \Big[
    f_{\ell_1,k_1}
    \big( V[[\ell_1]], X[0]
    \big)
    \,
    f_{\ell_2,k_2}
    \big( V[[\ell_2]], X[0]\big)
  \Big]\\ \nonumber
  &\qquad\qquad\qquad\qquad+
  \sum_{i=0}^{\min\{\ell_1,\ell_2\}}
  a_{i,\ell_1,k_1}
  a_{i,\ell_2,k_2}
  \bigg\}\,,\\
\label{eq:state-evolution-setup}
  \bbE[V^{\ell_1,k_1}V^{\ell_2, k_2}]
  &=
  \ind\{k_1=k_2\}
  \bbE 
  \bigg[
    h_{\ell_1,k_1}\Big(
    U[[\ell_1]] , Y[0] 
    \Big)
    \,
    h_{\ell_2,k_2}
    \Big(
    U[[\ell_2]] , Y[0] 
    \Big)
  \bigg]\,,
\end{align}
and $\bbE[U^{\ell_1+1,k}\bar{G}^i] = a_{i,\ell,k}$.\footnote{We remark that on the right-hand side of the first equality in \eqref{eq:state-evolution-setup}, the expectation is over both $\vec{V}$ and $X[0]$: thus the distribution of $\vec{U}$ depends on the \emph{distribution} of $X[0]$, not on the realization of $X[0]$ itself.
Likewise, on the right-hand side of the second equality in \eqref{eq:state-evolution-setup}, the expectation is over both $\vec{U}$ and $Y[0]$.}
The Onsager correction terms are then defined by
  \beq\label{eq:IAMP-fully-general-fixed.Ons}
  \begin{aligned}
  \ons_{u,\ell,k}
  &=
  \alpha
  \sum_{j=0}^\ell
  \bbE\bigg[
  \frac{\partial f_{\ell,k}}
  {\partial v^{j,k}}
  (V[[\ell]])
  \bigg]
  \cdot
  h_{j,k} (\bu[[j]], \by[0] 
	) \in\R^N\,,\\
  \ons_{v,\ell-1,k}
  &=
  \sum_{j=0}^\ell
  \bbE\bigg[
  \frac{\partial h_{\ell,k}}
  {\partial u^{j,k}}
    (U[[\ell]])
  \bigg]
  \cdot
  f_{j-1,k} (\bv[[j-1]]
  , \bx[0]
  )
  \in\R^M\,.
  \end{aligned}\eeq
We say $\psi:\bbR^d\to\bbR$ is \textbf{$1$-pseudo-Lipschitz} if 
  \[
  |\psi(x)-\psi(y)|\leq (|x|+|y|)\cdot (|x-y|)\,.\]
The resulting AMP characterization is as follows, where we note that the proper normalization is given as \cite[Assump.~(A1)]{gerbelot2021graph}.

\begin{ppn}[{\cite{gerbelot2021graph}}]
\label{prop:graph-AMP-cor}
In the AMP recursion \eqref{eq:IAMP-fully-general-fixed}, for any $\ell\geq 1$ and $1$-pseudo-Lipschitz $\psi$, we have the convergence in probability to the limiting random variables
\eqref{eq:state-evolution-setup} in the sense
  \begin{equation}
    \label{eq:state-evolution-limit-graph}
    \begin{aligned}
    \plim_{N\to\infty}
    \frac{1}{N}
    \sum_{n=1}^N
    \psi\Big(
    u_n[[\ell]],
    \bar{g}_n[\ell]
    , y_n[0] 
    \Big)
    &=
    \bbE\psi\Big( 
    \vec U, \vec{\bar{G}}
    , Y[0] 
     \Big)\,,
    \\
    \plim_{N\to\infty}
    \frac{1}{M}
    \sum_{m=1}^M
    \psi\Big(
    v_m[[\ell]]
    , x_m[0] 
    \Big)
    &=
    \bbE\psi\Big( \vec V,
    	X[0] \Big)\,.
    \end{aligned}
  \end{equation}
We also have convergence in Wasserstein $\bbW_2$ distance,
  \beq
  \label{eq:state-evolution-limit-W2}
    \frac{1}{N}
    \sum_{n=1}^N
    \delta_{u_n[[\ell]],
    \bar{g}_n[\ell],y_n[0]
    }
    \to \Law( \vec U,
      \vec {\bar{G}}
      , Y[0] 
      ), \qquad 
    \frac{1}{M}
    \sum_{m=1}^M
    \delta_{v_m[[\ell]], x_m[0]} 
    \to \Law(\vec V, X[0] )\,,
  \eeq
where again the convergence holds in probability as $N\to\infty$.
\end{ppn}

Throughout this section, we will use the term \textbf{state evolution limit} to indicate convergence of the form \eqref{eq:state-evolution-limit-graph}, using capital letters to indicate limiting distributions and lower case bold to denote iterates in $N$ or $M$ dimensions (as in \eqref{eq:w-sigma-b-setup}). We similarly use the term \textbf{$\bbW_2$ state evolution limit} to indicate convergence of the form 
\eqref{eq:state-evolution-limit-W2}. It follows from Proposition~\ref{prop:graph-AMP-cor} that each observable $\psi$ can similarly be said to have a $\bbW_2$ state evolution limit given by its evaluation on the limit variables $V$ or $Y$. We state this formally as follows, where $\psi_\#\mu$ denotes the pushforward of the measure $\mu$ by the function $\psi$.

\begin{cor}
\label{cor:psi-state-evolution}
In the setting of Proposition~\ref{prop:graph-AMP-cor}, suppose $\psi$ preserves $\bbW_2$ convergence in the sense that if $\mu_k\to \mu$ in $\bbW_2$ as $k\to\infty$, then $\psi_\#\mu_k\to\psi_\#\mu$ in $\bbW_2$ as $k\to\infty$. Let
  \[
  \boldsymbol{\psi}
  \equiv
  \psi(\bu[[\ell]],\bar{\bg}[\ell], \by[0])
  \equiv
  \Big( \psi
  ( u_n[[\ell]],
  \bar{g}_n[\ell]
  , y_n[0] 
   )\Big)_{n\le N} \in \R^N\,.\]
Then the vector $\boldsymbol{\psi}$ also obeys $\bbW_2$ state evolution, meaning that in $\bbW_2$ distance we have
  \[\frac{1}{N}
    \sum\limits_{n=1}^N
    \delta_{
      \psi_n
    }
    \to \Law(\psi( \vec U,
    \vec {\bar{G}}
    , Y[0] 
    ))\]
in probability as $N\to\infty$.
A similar statement holds for
$(\bv[[\ell]],\bx[0])$ and $(\vec V,X[0])$.
\end{cor}

\begin{proof}
  Immediate from Proposition~\ref{prop:graph-AMP-cor}.
\end{proof}

\begin{ppn}
\label{prop:W2-convergence}
If $\psi:\bbR^{d_1}\to\bbR^{d_2}$ is Lipschitz, then $\psi$ preserves $\bbW_2$ convergence.
If $d_1=2,d_2=1$ and $\psi(X,Y)=XY$, then $\psi$ preserves $\bbW_2$ convergence for sequences on which $X$ is uniformly bounded.

\begin{proof}
The first statement follows because $L$-Lipschitz maps increase $\bbW_2$ distances by at most a factor of $L$.
For the second statement, suppose 
$X_N\to X$ and $Y_N\to Y$ in $L^2$, with $|X_N|\leq C$ almost surely for all $N$. 
We will argue that this implies $X_NY_N\to XY$ in $L^2$, and the assertion then follows. First,
the uniform bound on $X_N$ gives
  \[
  \E\Big[(X_NY_N-X_NY)^2\Big]
  \le C^2 \E\Big[(Y_N-Y)^2\Big]
  = o_N(1)\,,
  \]
so $X_NY_N$ and $X_NY$ are close in $L^2$. For $T$ growing sufficiently slowly with $N$, we then bound
  \[
  \E\Big[ (X_NY-XY)^2\Big]
  \le T^2\E\Big[ (X_N-X)^2\Big]
  + 4C^2\E\Big[Y^2; |Y| \ge T\Big]
  \le T^2 o_N(1) + o_T(1)\,,
  \]
where the last term is because $Y^2\in L^1$. Therefore $X_NY$ converges to $XY$ in $L^2$, and the claim follows.
\end{proof}
\end{ppn}

The above constructions allow us to define two probability spaces. 
The first corresponds to the $N$-dimensional state evolution limit of $(\vec U,\vec{\bar G}, Y[0])$, and the second to the $M$-dimensional limit of $(\vec V, X[0])$. 
It may be useful to implicitly think of state evolution computation as being taken in the product of these two probability spaces, in which state evolution limits of $N$-dimensional vectors are independent of limits of $M$-dimensional vectors. However it is not formally necessary to couple together the spaces at all.

\subsection{General form IAMP iterations}
\label{subsec:IAMP-setup}

 For the rest of this section, we fix $q_0 \in [0,1)$ and $(b,\sigma,w,p,\zeta,\zeta^\Ising)$ satisfying assumptions (\ref{i:IAMP-main-coefs})--(\ref{i:IAMP-main-init}) of Theorem~\ref{thm:IAMP-main}. We will construct an IAMP algorithm whose coordinate marginals follow a discrete-time approximation of the SDE solutions $X,Y$ to \eqref{e:IAMP.X.SDE} and \eqref{e:IAMP.Y.SDE}. Throughout this section we focus on the Ising perceptron; as explained in Remark~\ref{r:IAMP.spherical}, the spherical perceptron can be treated as a special case. Recall from \eqref{e:matrix.br.bridge.IAMP} that we set
\[
  \bG(t) = t\bG + \bG_{\circ}(t)
\]
where $\bG_{\circ}$ is an $M\times N$ matrix-valued Brownian bridge independent of $\bG$. Let $p \in \incr([q_0,1];[0,1])$ be as in Theorem~\ref{thm:IAMP-main}.
Fix a small $\delta > 0$ such that $\ubl \equiv (1 - q_0)/\delta$ is an integer, and abbreviate throughout
\[
  p^k \equiv p(q_0 + k\delta)\,.
\]
Define the matrices
\beq\label{eq:A-ell}
\bA_{\ell}
\equiv \frac{\bG(p(q_0+\ell\delta))}{N^{1/2}}
= \frac{\bG(p^\ell)}{N^{1/2}},
\quad 0\leq \ell\leq \ubl\,,
\eeq
and note in particular that $\bA_{\ubl} = \bG/N^{1/2} \equiv \bA$. 
Define also the normalized increments
\beq\label{eq:Jk}
  \bJ_k
  =
  \frac{
  \bA_k-\bA_{k-1}
  }
  {(p^k-p^{k-1})^{1/2}
  }\,,\quad
  0\leq k\leq \ubl,
  \eeq
with the convention $p^k=0$ for $k \leq 0$. The matrices $\bJ_k$ are i.i.d. samples from $\cN(0,I_{M\times N}/N)$, as in the setup of \S\ref{subsec:state-evolution-multiple}. 
They will thus form our basis for state evolution computations.

The AMP iteration encoding the above diffusions involves four primary sequences of high-dimensional vectors, denoted by $\bx^{\ell},\prxbx^{\ell},\by^{\ell},\prxby^{\ell}$, for $0\leq \ell\leq \ubl$.
The iterates $\prxbx^{\ell}\in\R^M$ and $\prxby^{\ell}\in\R^N$ will converge in the state evolution sense to $X_t$ and $Y_t$ respectively as $\delta\to 0$, but may depend badly on $L$. 
Meanwhile the errors $\|\bx^{\ell}-\prxbx^{\ell}\|/N^{1/2}$ and $\|\by^{\ell}-\prxby^{\ell}\|/N^{1/2}$ will be $O(\epsilon^{1/2})$, coming from the various $\epsilon$-errors in assumptions (\ref{i:IAMP-main-diffusivity})--(\ref{i:IAMP-main-init}) of Theorem~\ref{thm:IAMP-main}. This two-stage approximation ensures the overall error vanishes in the limiting regime \eqref{r:Lip.epsilon.discussion}, even though $L$ can be large depending on $\epsilon$. On a first reading it may help to set $\epsilon=0$, so that we have $\bx^{\ell}=\prxbx^{\ell}$ and $\by^{\ell}=\prxby^{\ell}$, though as mentioned in Remark~\ref{r:Lip.epsilon.discussion} the formal result requires $\epsilon > 0$.

\begin{dfn}[IAMP mimicking BOGP SDE]
\label{d:iamp.v1}
  Recall from Theorem~\ref{thm:IAMP-main} 
  the initial conditions $X(q_0)\sim\zeta$ and $Y(q_0)\sim\zeta^\Ising$. Now let
  \[
  \breve{X}(q_0)\equiv0\,,\quad
  \breve{Y}(q_0)\equiv
    \frac{Y(q_0) \sqrt{q_0}}{(\bbE[Y(q_0)^2])^{1/2}}\,.\]
If $Y(q_0)\equiv0$ then we (somewhat arbitrarily) set $\breve{Y}(q_0)\sim\cN(0,q_0)$.
Then $\E[\breve{Y}(q_0)^2]=q_0$ in either case,
and it follows from
assumption~(\ref{i:IAMP-main-init}) of Theorem~\ref{thm:IAMP-main} 
that
  \beq\label{eq:tilde-Y-approx-Y}
    \bbE\bigg[
      \Big(
  Y(q_0)-\breve{Y}(q_0)
  \Big)^2\bigg]
  =\Big(
   \E[Y(q_0)^2]^{1/2}
   -  (q_0)^{1/2}\Big)^2
    \leq \Big|
    \bbE[Y(q_0)^2]-q_0\Big|
    \leq \epsilon\,.
    \eeq
We initialize vectors $\bx^0,\prxbx^0\in\R^M$ and 
$\by^0,\prxby^0\in\R^N$ as follows:
\begin{itemize}
\item The coordinates of $\prxby^0$ are sampled i.i.d.\ from $\Law(Y(q_0))=\zeta^\Ising$. If $\bbE[Y(q_0)^2]>0$, then we set 
	\[\by^0 = \bigg(
	\frac{q_0}{\bbE[Y(q_0)^2]}
	\bigg)^{1/2} \prxby^0\,.\]
If $\bbE[Y(q_0)^2]=0$, then we let $\by^0\sim\mathcal{N}(0,q_0 I_N)$. 

\item the coordinates of $\prxbx^0$ are sampled i.i.d.\ from $\Law(X(q_0))=\zeta$, and we set $\bx^0 = \bzero \in \R^M$.
\end{itemize}
For fixed $\delta$ we abbreviate $ b^\ell(x)\equiv b_{q_0+\ell\delta}(x) \equiv b(q_0+\ell\delta,x)$, etc. We then introduce the vectors
\beq\label{eq:w-sigma-b-setup}
\begin{aligned}
\bb^{\ell}
&\equiv
b^\ell(\prxbx^{\ell})
\equiv
b_{q_0+\ell\delta}(\prxbx^{\ell}) \in \bbR^M\,,\\
\bsig^{\ell}
&=
\sigma^\ell(\prxbx^{\ell})
\equiv
\sigma_{q_0+\ell\delta}(\prxbx^{\ell}) \in \bbR^M,
\\
\bw^{\ell}
&\equiv
w^\ell(\prxby^{\ell})
\equiv w_{q_0+\ell\delta}(\prxby^{\ell})
 \in \bbR^N, 
\end{aligned}\eeq
where the functions $b^\ell,\sigma^\ell,w^\ell$ are applied coordinate-wise. Let 
  \beq\label{e:c.ell.budget.cost}c_\ell(\epsilon)
  \equiv
  \bigg[\alpha\delta
  \Big( \budget(q_0+\ell\delta)+\epsilon-\Cost(q_0+\ell\delta) \Big)
    \bigg]^{1/2}\,.
  \eeq
As in \eqref{eq:IAMP-fully-general-fixed}, each $\bar{\bg}^\ell$ is an independent standard gaussian vector in $\R^N$. We let $\gamma(\ell,\delta)$ and $\prxgamma(\ell,\delta)$ be certain positive scalars, defined in \eqref{eq:proxy-C-ell-delta-def} below.  The IAMP iterates are as follows, with $\bu^0=\bzero \in \R^N$ and $\bm^{-1} = \bm^0=\bzero\in\R^M$: 
\begin{align}
\label{eq:proxy-y-amp-def}
  \prxby^{\ell+1}
  -
  \prxby^{\ell}
  &=
  \prxgamma(\ell,\delta)^{-1}
  \bw^{\ell}\odot(\bu^{\ell+1}-\bu^{\ell})
  \in\bbR^N,
  \\
\label{eq:y-amp-def}
  \by^{\ell+1}
  -
  \by^{\ell}
  &=
  \gamma(\ell,\delta)^{-1}
  \bw^{\ell}\odot(\bu^{\ell+1}-\bu^{\ell})\in\bbR^N\,,
  \\
\label{eq:u-amp-def}
  \bu^{\ell+1}
  -
  \bu^{\ell}
  &\equiv
  (\bA_\ell)^{\top}
  \big(
    \bm^{\ell}
    -
    \bm^{\ell-1}
  \big)
  + \delta^{1/2}
  (\bJ_{\ell+1})^\top\bb^\ell
  + c_\ell(\epsilon)
  \bar{\bg}^{\ell}
  -
  \ons_{\bu,\ell}
  \in\R^N\,,\\
\label{eq:v-amp-def}
  \bv^{\ell}&=\bA_{\ell}\by^{\ell}-\ons_{\bv,\ell-1}\in\bbR^M\,,
  \\
\label{eq:m-amp-def}
  \bm^{\ell+1}-\bm^{\ell}
  &=p(q_0+\ell\delta)^{-1}
  \big(\bsig^{\ell}
  -
  1
  \big)
  \odot
  \big(\bv^{\ell+1}-\bv^{\ell}\big)
  \in \bbR^M\,,
  \\
\label{eq:x-amp-def}
  \prxbx^{\ell+1}
  -
  \prxbx^{\ell}
  =
  \bx^{\ell+1}
  -
  \bx^{\ell}
  &=
  \bsig^{\ell}\odot
  \big(\bv^{\ell+1}-\bv^{\ell}\big)
  +
  \delta
  p'(q_0+\ell\delta)^{1/2}
  \,
  \bb^{\ell}
  \in 
  \bbR^M
  \,.
\end{align}
The Onsager terms $\ons_{\bu,\ell}$ and $\ons_{\bv,\ell-1}$ will be defined in \eqref{e:u.combined.ons} and \eqref{e:v.combined.ons} below. 
\textbf{We denote the above iteration by $\IAMP_N(\sigma,b,w,p,\zeta,\zeta^\Ising,\delta)$.} 
In the notation of Remark~\ref{rmk:omega-general-seed}, 
this algorithm will have random seed given by
	\beq\label{e:iamp.random.seed}
	\omega= \Big(
	\bG_\circ,
	(\bar\bg^\ell)_{1\le\ell\le\ubl},
	\bx[0], \by[0]\Big)\,.
	\eeq
For a fixed realization of $\omega$, this IAMP iteration will be a $C(L,\epsilon)$-Lipschitz function of the random matrix $\bG$ on an event with probability $1-e^{-cN}$, and the algorithm witnessing Theorem~\ref{thm:IAMP-main}\ref{i:IAMP-main-main} is its globally Lipschitz extension. 
\end{dfn}

\textbf{We refer to Remarks~\ref{r:IAMP.intuition} and \ref{r:proxy} below for some heuristic discussion of the above recursion.} The \hyperlink{proof:t.IAMP.main.main}{proof of Theorem~\ref{thm:IAMP-main}\ref{i:IAMP-main-main}} is based on analyzing this AMP iteration, and appears at the end of \S\ref{ss:proof.iamp.diffus.limit}. Towards this end, we first rewrite it
in the standard form \eqref{eq:IAMP-fully-general-fixed} we introduced for the relevant AMP. Recall from \eqref{eq:Jk} that we can decompose
\beq
\label{eq:A-ell-decompose}
  \bA_{\ell}
  = \sum_{k=0}^{\ell}
  (p^k-p^{k-1})^{1/2}
  \bJ_k
  \equiv \sum_{k=0}^{\ell}
  a_k
  \bJ_k\,.
  \eeq
Applying \eqref{eq:u-amp-def}, \eqref{eq:v-amp-def}, \eqref{eq:A-ell-decompose} yields 
\begin{align*}\nonumber
  \bu^{\ell+1}
  &= \sum_{j=0}^\ell (\bu^{j+1} - \bu^j)
  = \sum_{j=0}^\ell \bigg\{
    \bigg(\sum_{k=0}^j a_k \bJ_k\bigg)^\top (\bm^j - \bm^{j-1})
    + \delta^{1/2} (\bJ_{j+1})^\top \bb^j 
    + c_j(\epsilon) \bar{\bg}^j - \ons_{u,j}
  \bigg\} \\
  \nonumber
  &= \sum_{k=0}^{\ell+1}
  \bigg\{
  (\bJ_k)^\top\Big[
  \ind\{k\le\ell\}
  a_k(\bm^\ell-\bm^{k-1})
  +\ind\{k\ge1\}
  \delta^{1/2}\bb^{k-1}
  \Big]
  +\ind\{k\le\ell\} \Big[
    c_k(\epsilon)\bar{\bg}^k - \ons_{\bu,k}
  \Big]
  \bigg\}\,, \\
  \bv^\ell 
  &= \sum_{k=0}^\ell \bJ_k a_k \by^\ell - \ons_{\bv,\ell-1}\,.
\end{align*}
We express this in the format of \eqref{eq:IAMP-fully-general-fixed} by writing
\begin{equation}\label{e:u.ellp1.decomp}
\begin{aligned}
  \bu^{\ell+1}
  &=\sum_{k=0}^{\ell+1}
  \bu^{\ell+1,k}\,, \\
	\bv^{\ell}
  &=\sum_{k=0}^{\ell}
  \bv^{\ell,k}\,,
\end{aligned}\end{equation}
for the following choices of parameters in \eqref{eq:IAMP-fully-general-fixed}. We take $\bx[0]$ to be the $M\times 2$ matrix with columns $\prxbx^0$, $\bx^0$.
We take $\by[0]$ to be the $N\times 2$ matrix with columns $\prxby^0$, $\by^0$. Note that these satisfy the assumption \eqref{e:iamp-side-information-w2-convergence} because, by the construction in Definition~\ref{d:iamp.v1}, the rows of $\bx[0]$ and $\by[0]$ are i.i.d. draws from a $L^2$-bounded distribution in $\cP(\R^2)$. We take $a_{i,\ell,k}=\ind\{i=k\le\ell\}c_k(\epsilon)$, and 
\begin{align}\nonumber
  f_{\ell,k}(\bv[[\ell]], \bx[0] )
  &= \ind\{k\le\ell\}
  a_k
  \sum_{s=k-1}^{\ell-1}
  (\bm^{s+1}-\bm^s)
  +\ind\{k\ge1\}
  \delta^{1/2}\bb^{k-1}
  \\
  \label{e:f.l.k.explicit} 
  &= \ind\{k\le\ell\}
  a_k
  \sum_{s=k-1}^{\ell-1}
  \frac{1}{p^s}
  (\bsig^s-1)
  \odot(\bv^{s+1}-\bv^s)
  +\ind\{k\ge1\}
  \delta^{1/2} \bb^{k-1}\,, \\
  \label{e:h.l.k.explicit}
  h_{\ell,k}(\bu[[\ell]], \by[0] )
  &= a_k\by^\ell
  = a_k\sum_{s=0}^{\ell-1}
  (\by^{s+1}-\by^s) 
  + a_k\by^0 \\
  &=a_k\sum_{s=0}^{\ell-1}
  \frac{1}{\gamma(s,\delta)}
  \bw^s
  \odot(\bu^{s+1}-\bu^s)
  + a_k\by^0 \,.
\end{align}
Note the above does not immediately express $f_{\ell,k}(\cdot,\bx[0])$ and $h_{\ell,k}(\cdot,\by[0])$ as functions of the $\bv[[\ell]]$ and $\bu[[\ell]]$ alone. 
However, because $\bsig^s$ and $\bb^s$ are functions of $\prxbx^s$, by recursive expansion we can write $f_{\ell,k}(\cdot,\bx[0])$ as a function of $\bv[[\ell]]$ only.
Similarly, we can write $h_{\ell,k}(\cdot,\by[0])$ as a function of the $\bu[[\ell]]$. We remark also that in the above iteration, $f_{\ell,k}$ does not depend on $\bx[0]$, and $h_{\ell,k}$ depends only on $\by^0$ (but not on $\prxby^0$). 

We can then define the Onsager corrections $\ons_{\bu,\ell,k}$ and $\ons_{\bv,\ell-1,k}$ for the standard form AMP iterates $\bu^{\ell+1,k}$, $\bv^{\ell,k}$ as in \eqref{eq:IAMP-fully-general-fixed.Ons}. We emphasize that in taking partial derivatives in \eqref{eq:IAMP-fully-general-fixed.Ons}, we use the versions of $f_{\ell,k}$ and $h_{\ell,k}$ that are written as functions of the variables $\bv[[\ell]]$ and $\bu[[\ell]]$ alone.
Finally, we define the Onsager corrections in Definition~\ref{d:iamp.v1} by
\begin{align}\label{e:u.combined.ons}
  \ons_{\bu,\ell}
  &= \sum_{k=0}^{\ell+1}
  \ons_{\bu,\ell,k}
  -\sum_{k=0}^{\ell}
  \ons_{\bu,\ell-1,k}\,,\\
  \ons_{\bv,\ell-1}
  &=\sum_{k=0}^\ell
  \ons_{\bv,\ell-1,k}\,.
  \label{e:v.combined.ons}
\end{align}
The state evolution limits $(\vec{U}, \bar{G}^\ell,Y[0])$ and $(V^{\ell,k}, X[0])$ are defined as in Corollary~\ref{cor:psi-state-evolution} for suitable $\psi$. In particular, we have
\begin{equation}\label{e:state.evol.initialization}
\begin{aligned}
X[0] &\equiv (\prxX^0,X^0) 
	\equiv (X(q_0),\breve{X}(q_0))
	= (X(q_0), 0)\,,\\
Y[0] &\equiv (\prxY^0,Y^0) 
\equiv (Y(q_0),\breve{Y}(q_0))\,.
\end{aligned}
\end{equation}
(We point out again that the $\bar{G}^\ell$ here are one-dimensional i.i.d.\ standard gaussian random variables, in contrast with $\bG$ which usually denotes an $M\times N$ gaussian matrix.) We emphasize that in defining the state evolution limit, we have sent $N\to\infty$ \textbf{but not} $\delta\to 0$.
We let
  \beq\label{e:state.evol.U}
  U^\ell
  =\sum_{k=0}^\ell U^{\ell,k}\,,
  \eeq
and similarly $V^\ell$. In particular, $U^0\equiv0$ while $V^0\sim\cN(0,p(q_0)q_0)$. We then define recursively
  \beq
  \label{e:state.evol.x.y.increments}
  \begin{aligned}
  Y^{\ell+1}-Y^\ell
  &= \gamma(\ell,\delta)^{-1}
  w^\ell
  (\prxY^\ell)
  (U^{\ell+1}-U^\ell)\,,\\
  \prxY^{\ell+1}-\prxY^\ell
  &= \prxgamma(\ell,\delta)^{-1}
  w^\ell
  (\prxY^\ell)
  (U^{\ell+1}-U^\ell)\,,\\
  X^{\ell+1}-X^\ell
  =\prxX^{\ell+1}-\prxX^\ell
  &= \sigma^\ell(\prxX^\ell)
  (V^{\ell+1}-V^{\ell}\big)
  +
  \delta
  p'(q_0+\ell\delta)^{1/2}
  b^\ell(\prxX^\ell)\,,\\
  M^{\ell+1}-M^\ell
  &= (p^\ell)^{-1}
  (\sigma^\ell(\prxX^\ell)-1)
  (V^{\ell+1}-V^\ell)\,,
  \end{aligned}
  \eeq
where we recall $p^\ell\equiv p(q_0+\ell\delta)$, $b^\ell\equiv b_{q_0+\ell\delta}$, etc.
The scalars $\gamma$ and $\prxgamma$ are defined by
  \beq\label{eq:proxy-C-ell-delta-def}
  \begin{aligned}
  \prxgamma^\ell
  \equiv
  \prxgamma(\ell,\delta)
  &\equiv 
  \delta^{-1/2}
  \bbE\big[(U^{\ell+1}-U^\ell)^2\big]^{1/2}\,,\\
  \gamma^\ell
  \equiv \gamma(\ell,\delta)
  &\equiv 
  \delta^{-1/2}
  \bbE\big[(U^{\ell+1}-U^\ell)^2\big]^{1/2}
    \bbE\big[
  w^\ell
  (\prxY^\ell)^2\big]^{1/2}
  \,.
  \end{aligned}
  \eeq
The validity of these state evolution limits (i.e., the fact that the IAMP iterates defined by \eqref{eq:proxy-y-amp-def}--\eqref{eq:x-amp-def} converge to these limits in the manner of Proposition~\ref{prop:graph-AMP-cor}) will be shown in Proposition~\ref{ppn:state-evolution-valid-for-IAMP} below.
It will follow by repeatedly applying Proposition~\ref{prop:W2-convergence}, after verifying that $\Gamma^{\ell}$ and $\gamma^{\ell}$ are strictly positive.

\begin{rmk}[IAMP intuition]\label{r:IAMP.intuition}
We present here some overview and heuristic intuition on the IAMP defined by \eqref{eq:proxy-y-amp-def}--\eqref{eq:x-amp-def} and \eqref{eq:proxy-C-ell-delta-def}. For simplicity, in this discussion we consider only the spherical perceptron with $\epsilon=0$, for which we can use the simplified IAMP
  \begin{align} \label{eq:u-amp-def.spherical}
  \bu^{\ell+1}
  -
  \bu^{\ell}
  &\equiv
  (\bA_\ell)^{\top}
  \big(
    \bm^{\ell}
    -
    \bm^{\ell-1}
  \big)
  + \delta^{1/2}
  (\bJ_{\ell+1})^\top\bb^\ell
  + c_\ell(\epsilon)
  \bar{\bg}^{\ell}
  -
  \ons_{\bu,\ell}
  \in\R^N\,, \\ \nonumber
  \bv^{\ell}&=\bA_{\ell} \bu^{\ell}-\ons_{\bv,\ell-1}\in\bbR^M\,,\\
  \nonumber
  \bm^{\ell+1}-\bm^{\ell}
  &=p(q_0+\ell\delta)^{-1}
  \big(\bsig^{\ell}
  -
  1
  \big)
  \odot
  \big(\bv^{\ell+1}-\bv^{\ell}\big)
  \in \bbR^M\,,\\ 
  \bx^{\ell+1}
  -
  \bx^{\ell}
  &=
  \bsig^{\ell}\odot
  \big(\bv^{\ell+1}-\bv^{\ell}\big)
  +
  \delta
  p'(q_0+\ell\delta)^{1/2}
  \,
  \bb^{\ell}
  \in 
  \bbR^M
  \,, \nonumber
  \end{align}
where $\bb^\ell,\bsig^\ell$ correspond to $b^\ell,\sigma^\ell$ evaluated on $\bx^\ell$. Recall that the $\bA_\ell$ are (normalized) observations of the matrix Brownian bridge
\eqref{e:matrix.br.bridge.IAMP}: the basic idea is that the iterates $\bu^\ell\in\R^N$ should mimic the behavior of the Doob martingale
  \[
  \bx^\textup{alg}(t)
  \equiv \E[\cA(\bG)\,|\,\bG(t)]
  \in\R^N\,,
  \]
and the coordinate profile of $\bu^\ell$ should mimic the evolution of the SDE
solution $Y(t)$ from \eqref{e:IAMP.Y.SDE}. Similarly, the iterates $\bx^\ell\in\R^M$ should mimic the behavior of the $\R^M$-valued process
  \[ 
  \frac{\bG(t) \bx^\textup{alg}(t)}{N^{1/2}}
  \in\R^M\,,
  \]
and the coordinate profile of $\bx^\ell$ should mimic the evolution of the SDE solution $X(t)$ from \eqref{e:IAMP.X.SDE}. Most importantly, the two sequences will be linked by the relation $\bA_\ell\bu^\ell \approx \bx^\ell$. The algorithm finally outputs (a Lipschitz approximation of) $\bu^{\ubl}$, so that $\bG\bu^{\ubl}\approx\bx^{\ubl}$ will approximate $\mathscr{L}(X_1)$, as desired. 
The key calculations involved are as follows:
\begin{itemize}
\item Using the relation
\eqref{eq:v-amp-def}, it is straightforward to calculate
  \[
  \frac{\E[(V^{\ell+1}-V^\ell)^2]^{1/2}}{\delta^{1/2}}
  =  s_\delta(t)
  \approx s(t)\,,
  \]
for $s_\delta(t)$ defined by \eqref{e:s.sqrt.fn.delta} below.
Therefore we can let $B$ be the standard Brownian motion with increments
  \[
  B(q_0+(\ell+1)\delta)
  -B(q_0+\ell\delta)
  = \frac{V^{\ell+1}-V^\ell}{s_\delta(t)}\,.
  \]
Then \eqref{eq:x-amp-def} suggests
that $\bx^\ell$ mimics the behavior of the SDE \eqref{e:IAMP.X.SDE} for $X(t)$, with driving Brownian motion $B$.
This is proved in Proposition~\ref{prop:Ising-IAMP-analysis}\ref{it:m.x.barx.stateevol}.

\item \textbf{It remains to explain the role of the most complicated step \eqref{eq:u-amp-def.spherical}. This is specifically designed to achieve the goal $\bA_\ell\bu^\ell\approx \bx^\ell$.} More precisely, it follows from \eqref{eq:u-amp-def.spherical} that
  \[
  \bA_\ell\bu^\ell-\bx^\ell
  =\ons_{\bv,\ell-1} -(\bx^\ell-\bv^\ell)\,,
  \]
where $\bx^\ell$ can be expanded using \eqref{eq:x-amp-def}, and $\ons_{\bv,\ell-1}$ is given explicitly by 
\eqref{eq:IAMP-fully-general-fixed.Ons} and \eqref{e:v.combined.ons}.
A careful expansion of the above quantities suggests the approximate form of $f_{\ell,k}$. The details of this calculation are carried out in the \hyperlink{proof:p.IAMP-to-SDE-limit.Ax.y}{proof of Proposition~\ref{prop:IAMP-to-SDE-limit.Ax.y}}. 
On a more intuitive but heuristic level, the recursion \eqref{eq:u-amp-def.spherical} may be interpreted as follows: recall from \eqref{e:X.decomp} that we decompose
  \[
  \bA_{\ell+1}\bu^{\ell+1}-\bA_\ell\bu^\ell
  = a_{\ell+1}\bJ_{\ell+1}(\bu^{\ell+1}-\bu^\ell)
  + \bA_\ell(\bu^{\ell+1}-\bu^\ell) +a_{\ell+1}\bJ_{\ell+1} \bu^\ell
  \]
where the three terms on the right-hand side correspond to the labels $\RomI,\RomII,\RomIII$ from \S\ref{ss:X.decomp}. Imagine that we have already reached $\bu^\ell$, and want to choose $\bu^{\ell+1}$.
In \eqref{eq:u-amp-def.spherical}, the middle term
$(\bJ_{\ell+1})^\top\bb^\ell$  is the natural choice for term $(\RomI)$ to create the desired drift. The first and third terms combined can be viewed as the diffusive term, where the diffusivity depends on the spatial location through the definition of $\bm^\ell-\bm^{\ell-1}$. 
\end{itemize}
Combining the above steps leads to the \hyperlink{proof:t.IAMP.main.main}{proof of Theorem~\ref{thm:IAMP-main}\ref{i:IAMP-main-main}}.
\end{rmk}

\begin{rmk}[role of proxy processes]
\label{r:proxy}
In the previous remark we took for simplicity $\epsilon=0$, so that there was no distinction between $\bx^\ell,\by^\ell$ and the ``proxy'' processes $\prxbx^\ell,\prxby^\ell$. We now make some comments on the distinction between the two processes. The proxy processes
$\prxbx^\ell,\prxby^\ell$ are designed to be an extremely close approximation to the solutions of the SDEs \eqref{e:IAMP.X.SDE} and \eqref{e:IAMP.Y.SDE}: see in particular Proposition~\ref{prop:IAMP-to-SDE-limit}, where we emphasize that the error bounds in \eqref{eq:X-diffusive-limit} and \eqref{eq:Y-diffusive-limit} go to $0$ with $\delta$ (with possibly bad dependence on $L$) and do not depend on $\epsilon$. The process $\by^\ell$ arises from slightly rescaling the increments of $\prxby^\ell$ to achieve a precise radius schedule, which is convenient for some of the analysis; however, we expect $\by^\ell$ to be a worse approximation of the SDE, with error bounds involving $\epsilon$ rather than $\delta$ alone (see \eqref{eq:by-prxby-close}). Thus, in \eqref{eq:w-sigma-b-setup} it is essential that the coefficient functions $b,\sigma,w$ are evaluated on $\prxbx^\ell,\prxby^\ell$
rather than on $\bx^\ell,\by^\ell$, so that the overall error tends to $0$ in the regime \eqref{r:Lip.epsilon.discussion}. 
Ultimately, in the below 
\hyperlink{proof:t.IAMP.main.main}{proof of Theorem~\ref{thm:IAMP-main}\ref{i:IAMP-main-main}},
we can take either
$\cA(\bG)=\bar{c}_N\by^{\ubl}$
or
$\cA(\bG)=\bar{c}_N\prxby^{\ubl}$, for deterministic $\bar{c}_N$ chosen so that $\bbE[\|\cA(\bG)\|^2] = N$ as required by Definition~\ref{d:Lip}.
\end{rmk}

In order to compare the state evolution limits with the solutions of the SDEs \eqref{e:IAMP.X.SDE} and \eqref{e:IAMP.Y.SDE}, we introduce some further notation. Denote
    \[t_{\delta}\equiv q_0
    +\delta\bigg\lfloor \frac{t-q_0}{\delta} \bigg\rfloor\,,\quad
    t_{\delta+}\equiv t_{\delta}+\delta\,,\quad k(t)
        \equiv \bigg\lfloor\frac{t-q_0}{\delta}\bigg\rfloor\,.\]
As functions of $t$ these are all piecewise constant and right-continuous.  Recalling \eqref{e:s.sqrt.fn}, abbreviate
  \beq\label{e:s.sqrt.fn.delta}
  s_\delta(t_\delta)
  \equiv 
   \bigg(   \frac{
  t_{\delta+} p(t_{\delta+})
  -
  t_{\delta} p(t_{\delta})
  }{\delta}
  \bigg)^{1/2}
  \eeq
We will verify in Lemma~\ref{lem:(tp)'-stable} below that $s_\delta$ is a good approximation of
$s(t)$.
Define now the piecewise constant continuous-time processes
  \beq\label{e:state.evol.piecewise}
  \begin{aligned}
  \prxY(t)
  &\equiv \prxY^{k(t)}\,,\quad
  \breve{Y}(t)
  \equiv Y^{k(t)}\,,\\
  \prxX(t)
  &\equiv \prxX^{k(t)}\,,\quad
  \breve{X}(t)
  \equiv X^{k(t)}\,,
  \end{aligned}
  \eeq
where we recall that 
$\prxY^k$, $Y^k$, 
$\prxX^k$, $X^k$ are the state evolution limits defined by \eqref{e:state.evol.x.y.increments} above. The next few results will be proved in \S\ref{ss:proof.iamp.diffus.limit}.

\begin{ppn}[\hyperlink{proof:p.Ising-IAMP-analysis}{proved in \S\ref{ss:proof.iamp.diffus.limit}}]
\label{prop:Ising-IAMP-analysis} In the setting of Theorem~\ref{thm:IAMP-main}\ref{i:IAMP-main-main}, and recalling the notation \eqref{e:state.evol.piecewise}, there exist standard Brownian motions $B$ and $W$ such that the behavior of the state evolution limit \eqref{e:state.evol.x.y.increments} can be described as follows:
\begin{enumerate}[(a)]
\item 
  \label{it:Ising-BM-limit}
The random variables $V^k,Y^k,\prxY^k$ can be written as
  \begin{align*}
  V^k
  &=
        V^0
        +
  \int_{q_0}^{q_0+k\delta}
  s_\delta(t_\delta)
  \,\de B(t),
  \\
  Y^k
  = \breve{Y}(q_0+k\delta)
  &=
        \breve{Y}(q_0)
        +
  \int_{q_0}^{q_0+k\delta}
  \frac{w_{t_\delta}(\prxY(t))}
  {
  \bbE[w_{t_\delta}(\prxY(t))^2]^{1/2}
  }
  \,\de W(t),
  \\
  \prxY^k
  = \prxY(q_0+k\delta)
  &=\prxY^0+
  \int_{q_0}^{q_0+k\delta}
  w_{t_\delta}(\prxY(t))
  \,\de W(t)
  .
  \end{align*}

\item
\label{it:m.x.barx.stateevol}
The random variables $M^k$, $X^k$, $\prxX^k$ can be written as
  \begin{align*}
  M^k
  &=
  \int_{q_0}^{q_0+k\delta} \frac{s_\delta(t_\delta)}{p(t_{\delta})}
  \Big(\sigma_{t_{\delta}}(\prxX(t))-1\Big)
  \,\de B(t);
  \\
  X^k = \breve{X}(q_0+k\delta)
  &=
  \int_{q_0}^{q_0+k\delta}
  s_\delta(t_\delta)
  \sigma_{t_{\delta}}(
  \prxX(t))
  \,\de B(t)
  +
  \int_{q_0}^{q_0+k\delta}
  p'(t_{\delta})^{1/2}
  b_{t_{\delta}}(\prxX(t))
  \,\de t
  \\
  \prxX^k
  =\prxX(q_0+k\delta)
  &=
  \prxX^0
  +
  \int_{q_0}^{q_0+k\delta}
  s_\delta(t_\delta)
  \sigma_{t_{\delta}}
  (\prxX(t))
  \,\de B(t)
  +
  \int_{q_0}^{q_0+k\delta}
  p'(t_{\delta})^{1/2}
  b_{t_{\delta}}
  (\prxX(t))\,\de t
  .\end{align*}
(Recall from Definition~\ref{d:iamp.v1} that $\breve{X}(q_0)=0$.)
\end{enumerate}
Recall $\prxX^0$ and $\prxY^0$
are equidistributed respectively as $X(q_0)$ and $Y(q_0)$ from the SDEs
\eqref{e:IAMP.X.SDE} and \eqref{e:IAMP.Y.SDE}. Moreover $\prxY^0$ is independent of $W$, and $\prxX^0$ is independent of $B$.\footnote{We recall from the discussion after Corollary~\ref{cor:psi-state-evolution} that $N$-dimensional and $M$-dimensional state evolution variables such as $W$ and $B$ are defined on different probability spaces, and can thus be  thought of as independent.}
\end{ppn}

\begin{ppn}[\hyperlink{proof:p.IAMP-to-SDE-limit}{proved in \S\ref{ss:proof.iamp.diffus.limit}}]
\label{prop:IAMP-to-SDE-limit}
Let $X(t)$ and $Y(t)$ be the solutions of the SDEs \eqref{e:IAMP.X.SDE} and \eqref{e:IAMP.Y.SDE}. Let $\prxX(t)$ and $\prxY(t)$ be the (piecewise constant) processes \eqref{e:state.evol.piecewise}, as characterized by Proposition~\ref{prop:Ising-IAMP-analysis}. If we couple them with the same driving Brownian motions $B$ and $W$, then we have
  \begin{align}
  \label{eq:X-diffusive-limit}
  \sup\bigg\{
  \bbE\Big[
  \Big(\prxX(q_0+\ell\delta)
    -X(q_0+\ell\delta)\Big)^2\Big]
  : 0\leq \ell \leq \ubl
  \bigg\}
  &\leq 
  C\delta,
  \\
  \label{eq:Y-diffusive-limit}
  \sup\bigg\{
  \bbE\Big[
  \Big(
  \prxY(q_0+\ell\delta)
  -Y(q_0+\ell\delta)
  \Big)^2
  \Big]
  : 0\leq \ell \leq \ubl
  \bigg\}
  &\leq 
  C\delta
  .
  \end{align}
where $C=C(L,\alpha)$.
\end{ppn}

\begin{ppn}[\hyperlink{proof:p.IAMP-to-SDE-limit.Ax.y}{proved in \S\ref{ss:proof.iamp.diffus.limit}}]
\label{prop:IAMP-to-SDE-limit.Ax.y}
 In the setting of Theorem~\ref{thm:IAMP-main}\ref{i:IAMP-main-main},  the IAMP defined by \eqref{eq:proxy-y-amp-def}--\eqref{eq:x-amp-def} satisfies the bound
  \[
  \sup\bigg\{
  \plim_{N\to\infty}
  \frac{\|\bA_{\ell}\by^\ell-\bx^\ell\|}{N^{1/2}}
  : 0\le \ell\le\ubl
  \bigg\}
  \leq 
  C\delta^{1/4}
  +C_0\epsilon^{1/2}\,.\]
where $C=C(L,\alpha)$ and $C_0=C_0(\alpha)$.
\end{ppn}

\subsection{Proof of IAMP diffusion limit}
\label{ss:proof.iamp.diffus.limit}

The main goal of this subsection is to complete the \hyperlink{proof:t.IAMP.main.main}{proof of Theorem~\ref{thm:IAMP-main}\ref{i:IAMP-main-main}}. Recalling Remark~\ref{rem:external-noise-is-ok}, we will treat the vectors $\bar{\bg}^k$ as external gaussian noise in the state evolution limits.
First we prove a lemma on orthogonal increments, which will also be relevant for later generalizations. Define the $\sigma$-fields
	\begin{align*}
	 \cF^{X,\prxX}(\ell)
	& \equiv \sigma(
	X^0,\ldots,X^\ell,
	\prxX^0,\dots,\prxX^{\ell}
	)\,,\\
	\cF^{Y,\prxY}(\ell)
	&\equiv \sigma(
	Y^0,\ldots,Y^\ell,
	\prxY^0,\dots,\prxY^{\ell}
	)\,,\\
	\cF^U(\ell)
	&\equiv \sigma( Y[0], U^0,\dots,U^{\ell})\,,\\
	\cF^V(\ell)
	&\equiv \sigma( X[0], V^0,\ldots,V^\ell)\,,\\
	\cF^M(\ell)
	&\equiv \sigma( X[0], M^0,\ldots,M^\ell)\,.
	\end{align*}
We emphasize that these $\sigma$-fields are defined only with respect to the state evolution random variables defined by \eqref{e:state.evol.x.y.increments}.

\begin{lem}\label{lem:orthogonal-increments}
  For all $0\leq \ell\leq \ubl$ we have:
  \begin{align}
  \label{eq:u-increments-orthogonal}
  \bbE[
  U^\ell-U^{\ell-1}
  \,|\,\cF^U(\ell-1)]
  &=0,
  \\
  \label{eq:y-increments-orthogonal}
  \bbE[
  Y^\ell-Y^{\ell-1}
  \,|\,
  \cF^{Y,\prxY}(\ell-1)]
  =
   \bbE[
  Y^\ell-Y^{\ell-1}
  \,|\,
  \cF^U(\ell-1)]
  &=
  0,
  \\
  \label{eq:v-increments-orthogonal}
  \bbE[
  V^\ell-V^{\ell-1}
  \,|\,
  \cF^V(\ell-1)]
  &=0,
  \\
  \label{eq:m-increments-orthogonal}
  \bbE[
  M^\ell-M^{\ell-1}
  \,|\,
  \cF^M(\ell-1)]
  =
  \bbE[M^\ell-M^{\ell-1}
  \,|\,
  \cF^V(\ell-1)]
   &=
  0
  .
  \end{align}

\begin{proof}
We will show inductively that
\eqref{eq:u-increments-orthogonal}
implies
\eqref{eq:y-increments-orthogonal},
which implies
\eqref{eq:v-increments-orthogonal},
which implies 
\eqref{eq:m-increments-orthogonal},
which finally implies
\eqref{eq:u-increments-orthogonal}
with $\ell+1$ in place of $\ell$.
In the base case, \eqref{eq:u-increments-orthogonal} with $\ell=1$ is trivial, since $U^0=0$, $\mathcal{F}^U(0)$ is the $\sigma$-field generated by $Y[0]$, and $U^1$ is a centered gaussian random variable independent of $Y[0]$. We therefore proceed with the inductive argument:

\begin{itemize}
\item
\eqref{eq:u-increments-orthogonal} $\Rightarrow$ \eqref{eq:y-increments-orthogonal}: from the definition
\eqref{e:state.evol.x.y.increments},
  \[
  \E\Big(Y^\ell-Y^{\ell-1} 
  \,\Big|\,
  \cF^{Y,\prxY}(\ell-1)\Big)
  = \frac{w^{\ell-1}
  (\prxY^{\ell-1})}{\gamma(\ell-1,\delta)}
  \E \Big(U^\ell-U^{\ell-1} 
  \,\Big|\,
  \cF^{Y,\prxY}(\ell-1)\Big)\,.
  \]
We also claim
that $\cF^{Y,\prxY}(\ell)$
is contained in 
$\cF^U(\ell)$ for all $\ell$: indeed, in the base case we have $\cF^{Y,\prxY}(0) =\cF^U(0)$, so if we suppose inductively that 
$\cF^{Y,\prxY}(\ell)$ is contained in $\cF^U(\ell)$, then
\eqref{e:state.evol.x.y.increments}  implies 
\[\cF^{Y,\prxY}(\ell+1)
\subseteq \sigma\Big(
\cF^{Y,\prxY}(\ell),
\cF^U(\ell+1)
\Big)
\subseteq \cF^U(\ell+1)
\,,\]
which verifies the inductive hypothesis. Thus \eqref{eq:u-increments-orthogonal} implies that the above equals zero, which in turn verifies \eqref{eq:y-increments-orthogonal}.

\item 
\eqref{eq:y-increments-orthogonal} $\Rightarrow$ \eqref{eq:v-increments-orthogonal}: we apply the state evolution recursion \eqref{eq:state-evolution-setup}, where
we recall from \eqref{e:h.l.k.explicit}
that $h_{\ell,k}=a_k Y^\ell$.  It follows that for $k\le s\le \ell -1$,
  \begin{align*}
  &\E\Big[(V^{\ell,k}-V^{\ell-1,k}) V^{s,k}\Big]\\
  &\qquad\stackrel{\eqref{eq:state-evolution-setup}}{=}
  \E\bigg[\Big( 
  h_{\ell,k}(U[[\ell]],
  	 Y[0] 
	)
  -h_{\ell-1,k}(U[[\ell-1]],
  	 Y[0] 
	\Big)
  h_{s,k}( U[[s]],
  	 Y[0] 
	)\bigg]
  \\
  &\qquad
  \stackrel{\eqref{e:h.l.k.explicit}}{=}
  (a_k)^2
  \E\Big[(Y^\ell-Y^{\ell-1})
  Y^s
  \Big]
  =0\,,
  \end{align*}
where the last step is by the inductive hypothesis \eqref{eq:y-increments-orthogonal}. Recall also from \eqref{eq:state-evolution-setup} that $V^{\ell_1,k_1}$ and $V^{\ell_2,k_2}$ are uncorrelated for $k_1\ne k_2$, so the above implies
  \[
  \E\Big[(V^\ell-V^{\ell-1}) 
  V^s\Big]=0
  \]
for $s\le \ell-1$. Since the $V^{\ell,k}$ are jointly gaussian, it follows that $V^\ell-V^{\ell-1}$ is independent of $V[[\ell-1]]$, thus yielding \eqref{eq:v-increments-orthogonal}.

\item \eqref{eq:v-increments-orthogonal} $\Rightarrow$ \eqref{eq:m-increments-orthogonal}: this follows by similar reasoning as for \eqref{eq:u-increments-orthogonal} $\Rightarrow$ \eqref{eq:y-increments-orthogonal} above. First, we claim that $\cF^{X,\prxX}(\ell)$ is contained in $\cF^V(\ell)$ for all $\ell\ge0$: indeed, in the base case we have 
$\cF^{X,\prxX}(0)=\cF^V(0)$, so if we suppose inductively that $\cF^{X,\prxX}(\ell)$ is contained in $\cF^V(\ell)$, then \eqref{e:state.evol.x.y.increments} implies 
	\[
	\cF^{X,\prxX}(\ell+1)
	\subseteq
	\sigma\Big(
	\cF^{X,\prxX}(\ell),
	\cF^V(\ell+1)
	\Big)
	\subseteq 
	\cF^V(\ell+1)\,,
	\]
which verifies the inductive hypothesis.
We also note that 
$\cF^M(0) = \cF^V(0)$, and another application of \eqref{e:state.evol.x.y.increments} gives
	\[
	\cF^M(\ell+1)
	\subseteq
	\sigma\Big(
	\cF^{X,\prxX}(\ell),
	\cF^V(\ell+1)
	\Big)
	\subseteq \cF^V(\ell+1)
	\]
for all $\ell\ge0$. From \eqref{e:state.evol.x.y.increments}  we have
	\[
	\bbE\Big( M^\ell-M^{\ell-1}
  \,\Big|\,
  \cF^V(\ell-1)\Big)
  	= 
	\frac{\sigma^{\ell-1}(\prxX^{\ell-1})-1}{p^{\ell-1}}
	\bbE\Big( V^\ell-V^{\ell-1}
  \,\Big|\,
  \cF^V(\ell-1)\Big)
  =0\,,
	\]
where the last equality holds by the hypothesis \eqref{eq:v-increments-orthogonal}. The same therefore holds if we replace $\cF^V(\ell-1)$ by the its sub-$\sigma$-field $\cF^M(\ell-1)$, which verifies \eqref{eq:m-increments-orthogonal}.

\item \eqref{eq:m-increments-orthogonal}
$\Rightarrow$
\eqref{eq:u-increments-orthogonal}:
we again apply the state evolution recursion \eqref{eq:state-evolution-setup}, where $f_{\ell,k}$ is given explicitly by \eqref{e:f.l.k.explicit}.  It follows that for $k\le  s+1\le \ell$,
  \begin{align*}
  &\E\Big[
  (U^{\ell+1,k} -U^{\ell,k})
  U^{s+1,k}
  \Big]
  = \alpha \, \E\bigg[a_k
  (M^\ell-M^{\ell-1})
   f_{s,k}(V[[s]] , X[0] )
   \bigg]
  \\
  &\qquad
  =
  \alpha a_k\E\bigg[
  (M^\ell-M^{\ell-1})
  \bigg\{
  \ind\{k\le s\}
  a_k(M^s-M^{k-1})
  +\ind\{k\ge1\}
  \delta^{1/2}
  b^{k-1}(\prxX^{k-1})
  \bigg\}
  \bigg]
  =0\,,
  \end{align*}
where the last step is by the inductive hypothesis \eqref{eq:m-increments-orthogonal}. Recall also from \eqref{eq:state-evolution-setup}
that $V^{\ell_1,k_1}$ and $V^{\ell_2,k_2}$ are uncorrelated for $k_1\ne k_2$, so the above implies
  \[
  \E\Big[(U^{\ell+1}-U^\ell) 
  U^{s+1}\Big]=0
  \]
for $s\le \ell-1$. Since the $U^{\ell,k}$ are jointly gaussian, it follows that $U^{\ell+1}-U^\ell$ is independent of $U[[\ell]]$, thus yielding \eqref{eq:u-increments-orthogonal} with $\ell+1$ in place of $\ell$.
\end{itemize}
The above steps give one complete round of the induction, and the claim follows.
\end{proof}
\end{lem}

\begin{proof}[\hypertarget{proof:p.Ising-IAMP-analysis}{Proof of Proposition~\ref{prop:Ising-IAMP-analysis}}] We begin with the proof of
part~\eqref{it:Ising-BM-limit}:
From Lemma~\ref{lem:orthogonal-increments},
$U^{\ell+1}-U^{\ell}$
is independent of 
$w^\ell(\prxY^{\ell})$.
Combining this fact with \eqref{e:state.evol.x.y.increments}
and \eqref{eq:proxy-C-ell-delta-def}  allows us to write
  \begin{align*}
  Y^{\ell+1}-Y^\ell
  &= 
  \frac{w(q_0+\ell\delta,\prxY^\ell)}
  {\E[w(q_0+\ell\delta,
    \prxY^\ell)^2]^{1/2}}\cdot
  \frac{\delta^{1/2}(U^{\ell+1}-U^\ell)}
    {\E[(U^{\ell+1}-U^\ell)^2]^{1/2}}\,,\\
  \prxY^{\ell+1}-\prxY^\ell
  &=w(q_0+\ell\delta,\prxY^\ell)
  \cdot 
  \frac{\delta^{1/2} (U^{\ell+1}-U^\ell)}{\E[(U^{\ell+1}-U^\ell)^2]^{1/2}}\,.
  \end{align*}
We therefore define $W$ to be a Brownian motion with increments
  \beq\label{e:disc.B.y}
  W(q_0+(\ell+1)\delta)
  -W(q_0+\ell\delta)
  = \frac{\delta^{1/2} (U^{\ell+1}-U^\ell)}{\E[(U^{\ell+1}-U^\ell)^2]^{1/2}}\,,
  \eeq
which yields the claimed representations for $Y$ and $\prxY$.
Next recalling the definition of $a_k$ in \eqref{eq:A-ell-decompose}, we use the state evolution recursion \eqref{eq:state-evolution-setup} to calculate, for $\ell_1\le\ell_2$,
  \begin{align*}
  \E[ V^{\ell_1,k}
  V^{\ell_2,k}]
  &=\E\Big[
    h_{\ell_1,k}(U[[\ell_1]], Y[0] ) 
    h_{\ell_2,k}(U[[\ell_2]], Y[0] )
    \Big]\\
  &= (a_k)^2
  \E\bigg[\bigg\{Y^0 + 
  \sum_{s_1=0}^{\ell_1-1} 
    (Y^{s_1+1}-Y^{s_1}) \bigg\}
  \bigg\{Y^0 + 
  \sum_{s_2=0}^{\ell_2-1}
    (Y^{s_2+1}-Y^{s_2})
    \bigg\}\bigg]
  \\
  &= (a_k)^2\bigg\{
  \E[(Y^0)^2] + 
  \sum_{s=0}^{\ell_1-1}\delta
  \bigg\}
  =( q_0 + \ell_1\delta) (p^k-p^{k-1})
  \,,
  \end{align*}
recalling again that $Y^0 \equiv \breve{Y}(q_0)$ has second moment $q_0$. Recall also from \eqref{eq:state-evolution-setup}
that $V^{\ell_1,k_1}$  and $V^{\ell_2,k_2}$ are uncorrelated for $k_1\ne k_2$. It follows that, for $\ell_1\le \ell_2$,
  \[
  \E[ V^{\ell_1}V^{\ell_2}]
  = \sum_{k=0}^{\ell_1}
  \E[ V^{\ell_1,k}V^{\ell_2,k}]
  = (q_0 + \ell_1\delta)
  \sum_{k=0}^{\ell_1}
  (p^k-p^{k-1})
  = (q_0+ \ell_1\delta)
   p(q_0+\ell_1\delta) \,.
  \]
Recalling Lemma~\ref{lem:orthogonal-increments} again, we conclude
  \beq\label{e:V-increment-var}
  \frac{\E[ (V^{\ell+1}-V^\ell)^2]}{\delta}
  = \frac{\E[(V^{\ell+1})^2]-\E[(V^\ell)^2]}{\delta}
  = s_\delta(q_0+\ell\delta)^2\,,
  \eeq
with $s_\delta$ as defined by \eqref{e:s.sqrt.fn.delta}. Then, similarly to 
\eqref{e:disc.B.y}, we define
$B$ to be the standard Brownian motion with increments
  \[
  B(q_0+(\ell+1)\delta)
  -B(q_0+\ell\delta)
  = \frac{\delta^{1/2} 
  (V^{\ell+1}-V^\ell)}
  {\E[(V^{\ell+1}-V^\ell)^2]^{1/2}}
  =\frac{V^{\ell+1}-V^\ell}{s_\delta(q_0+\ell\delta)}\,.
  \]
This implies the result of
part~\eqref{it:Ising-BM-limit}.
Part~\eqref{it:m.x.barx.stateevol}
follows directly from part \eqref{it:Ising-BM-limit}
combined with
\eqref{e:state.evol.x.y.increments}.
\end{proof}

The following simple lemma will be useful to analyze the $X_t$ diffusive limit.

\begin{lem}
\label{lem:(tp)'-stable}
Recall $s$ and $s_\delta$ from 
\eqref{e:s.sqrt.fn} and \eqref{e:s.sqrt.fn.delta}.
  There exists $C=C(L)$ such that for all $t\in [q_0,1]$,
  \[
     \Big|s(t)
      -s_{\delta}(t)\Big|
  \le C\delta.
  \]

\begin{proof}
Since $p(q_0)\geq 1/L$ and $p$ is increasing, we have
\[
  s(t) + s_\delta(t) \ge s(t) \ge L^{-1/2}\,.
\]
By adjusting $C$ appropriately it suffices to show
  \[
     \Big|s(t)^2
      -s_{\delta}(t)^2\Big|
  \le C\delta.
  \]
Note that
$s_\delta(t)^2$ is the average value of $s(t)^2$ on the interval $[t_\delta,t_{\delta+}]$:
  \[
  s_\delta(t)^2
  = \frac{t_{\delta+} p(t_{\delta+})
    -
    t_{\delta} p(t_{\delta})
    }{\delta}
    = \frac{1}{\delta}\int_{t_\delta}^{t_{\delta_+}}
    (tp)'(t)\,dt
    = \frac{1}{\delta}\int_{t_\delta}^{t_{\delta_+}} s(t)^2\,dt\,.
  \]
By the mean-value theorem, it suffices to show that $s(t)^2=(tp)'(t)$ is Lipschitz, which is clear since assumption (\ref{i:IAMP-main-p}) of Theorem~\ref{thm:IAMP-main} gives that  $\|p\|_{C^2([0,1])}\leq L$.
    \end{proof}
\end{lem}

\begin{proof}[\hypertarget{proof:p.IAMP-to-SDE-limit}{Proof of Proposition~\ref{prop:IAMP-to-SDE-limit}}] 
Throughout, let $C$ be a constant depending on $L$. It may vary from one occurrence to the next, but ultimately can be bounded in terms of $L$ alone.
\medskip

\noindent\textbf{Diffusive limit for $\prxY$.} We will compare the state evolution limit $\prxY$, as described by Proposition~\ref{prop:Ising-IAMP-analysis}\ref{it:Ising-BM-limit}, with the SDE solution $Y$ from \eqref{e:IAMP.Y.SDE}, assuming they are coupled via the same driving Brownian motion $W$.
Define the error function
  \[
  D^Y(t)
  \equiv \E\bigg[
  \Big( Y(t)-\prxY(t) \Big)^2
  \bigg] + \delta\,,
  \]
where we note that $D^Y(q_0)=\delta$ since we initialized $\prxY(q_0)=Y(q_0)$.
Recall also that $w$ is assumed to be $L$-Lipschitz in both space and time.
By It\=o's formula,
  \begin{align*}
  D^Y(t)
  &= \int_{q_0}^{t_\delta}
  \E\bigg[\Big(
  w(u,Y(u))-w(u_\delta,\prxY(u))
  \Big)^2\bigg]
  \,du
  + \int_{t_\delta}^t
  \E\Big[w(u,Y(u))^2\Big]
  \,du \\
  &\le C 
  \int_{q_0}^{t_\delta}
  \E \bigg[
  \Big(w(u,Y(u))
  -w(u_\delta,Y(u))\Big)^2
  +\Big( w(u_\delta,Y(u))
  -w(u_\delta,\prxY(u))
  \Big)^2
  \bigg]
  \,du + C\delta \\
  &\le 
  C\int_{q_0}^{t_\delta}\bigg\{
  (u-u_\delta)^2 
  +\E\Big[ (Y(u)-\prxY(u))^2\Big]
  \bigg\} \,du
  + C\delta 
  \le C \int_{q_0}^t D^Y(u)\,du + C\delta\,.
  \end{align*}
It follows by Gr\"onwall's inequality that $D^Y(t) \le C\delta \exp(Ct)$ for all $t$.
This implies \eqref{eq:Y-diffusive-limit}.\medskip

\noindent\textbf{Diffusive limit for $\prxX$.} We next compare the state evolution limit $\prxX$, described by Proposition~\ref{prop:Ising-IAMP-analysis}\ref{it:m.x.barx.stateevol},
with the SDE solution $X$ from \eqref{e:IAMP.X.SDE}, assuming they are coupled via the same driving Brownian motion $B$. Similarly as above, let
  \[
  D^X(t)
  \equiv \E\bigg[
  \Big( X(t)-\prxX(t) \Big)^2
  \bigg] + \delta\,,
  \]
where we note that $D^X(q_0)=\delta$ since we initialized $\prxX(q_0)=X(q_0)$.  By It\=o's formula,
\[
D^X(t) \le C D^{X,1}(t)+C D^{X,2}(t)+C\delta,
\]
where (using the notation from Lemma~\ref{lem:(tp)'-stable}):
  \begin{align*}
  D^{X,1}(t)
  &= \int_{q_0}^{t_\delta}
  \E\bigg[
    \Big( 
    s(u) \sigma(u,X(u))
    -s_\delta(u_\delta)
    \sigma(u_\delta,\prxX(u))\Big)^2
    \bigg]
    \,du
    +\int_{t_\delta}^t
      \E\Big[
      s(u)^2 \sigma(u,X(u))^2
      \Big]\,du\,,\\
  D^{X,2}(t)&=
    \E\bigg[\bigg(
    \int_{q_0}^t 
    \Big(
    p'(u)^{1/2} b(u,X(u))
    -p'(u_\delta)^{1/2}
    b(u_\delta,\prxX(u_\delta))
    \Big)\,du
    +\int_{t_\delta}^t p'(u) b(u,X(u))\,du
    \bigg)^2\bigg]\,.
  \end{align*}
Recall that the functions $\sigma$ and $b$ are assumed to be $L$-Lipschitz in both space in time. We then bound
  \begin{align*}
  D^{X,1}(t)
  &\le
  C\int_{q_0}^{t_\delta}
  \E\bigg[
  \Big(
  s(u)
  -s_\delta(u_\delta)
  \Big)^2\sigma(u,X(u))^2
  +s_\delta(u_\delta)^2
    \Big(\sigma(u,X(u))
  -\sigma(u_\delta,X(u))
    \Big)^2\\
  &\qquad\qquad\qquad
  +s_\delta(u_\delta)^2
    \Big(
    \sigma(u_\delta,X(u))
    -\sigma(u_\delta,\prxX(u))
    \Big)^2\bigg]\,du
  +C\delta\\
  &\le
  C\int_{q_0}^t D^X(u)\,du
  +C\delta\,,
  \end{align*}
having used Lemma~\ref{lem:(tp)'-stable} to compare $s(u)$ with $s_\delta(u_\delta)$. Similarly, we have
  \begin{align*}
  D^{X,2}(t)
  &\le
  C \bigg\{\E
  \int_{q_0}^t
  \Big( p'(u)^{1/2}
  -p'(u_\delta)^{1/2}\Big) b(u,X(u))
  +p'(u_\delta)^{1/2} \Big(
    b(u,X(u))- b(u_\delta,X(u))
  \Big)\\
  &\qquad
  +p'(u_\delta)^{1/2}
    \Big(
    b(u_\delta,X(u))
    -b(u_\delta,\prxX(u))
    \Big)
  \,du\bigg\}^2
  + C\delta
  \le 
  C \int_{q_0}^t D^X(u)\,du
  +C\delta\,.
  \end{align*}
Applying Gr\"onwall to the function $D^X$ 
gives \eqref{eq:X-diffusive-limit}.
\end{proof}

We next compute the variance of the $U$-increments, to be used in the proof of Proposition~\ref{prop:IAMP-to-SDE-limit.Ax.y}:

\begin{lem}\label{l:gamma.estimates}
Let $U^k$ be the state evolution limits from \eqref{e:state.evol.U}, and let $Y$ be the solution of the SDE \eqref{e:IAMP.Y.SDE}.  Then, for all $0\le \ell\le\ubl-1$, we have
  \beq\label{e:U.stdev.Ew}
  \bigg|\frac{\E[(U^{\ell+1}-U^\ell)^2
  ]^{1/2}}{\delta^{1/2}}
  - \E w^\ell
    (Y(q_0+\ell\delta)) \bigg|
  \le C\delta^{1/4}
  +C_0\epsilon^{1/2}
  \,,\eeq
where $C=C(L,\alpha)$ while $C_0=C_0(\alpha)$.
As a consequence, for $\gamma$ and $\prxgamma$ as defined by \eqref{eq:proxy-C-ell-delta-def}, we have
  \[\max\bigg\{
  \Big|\gamma(\ell,\delta)-
  \E w^\ell
    (Y(q_0+\ell\delta))\Big|
  ,\Big|\prxgamma(\ell,\delta)-
  \E w^\ell
    (Y(q_0+\ell\delta))\Big|
    \bigg\}
  \le C\delta^{1/4}
  +C_0\epsilon^{1/2}
  \]
for all $0\le\ell\le\ubl-1$.

\begin{proof}
We apply the state evolution recursion \eqref{eq:state-evolution-setup}, where $f_{\ell,k}$ and $a_{i,\ell,k}$ are defined by \eqref{e:u.ellp1.decomp} and \eqref{e:f.l.k.explicit}:
  \begin{align*}
  f_{\ell,k}(V[[\ell]], X[0])
  &= \ind\{k\le\ell\}
  a_k (M^\ell-M^{k-1})
  +\ind\{k\ge1\}
  \delta^{1/2} 
  b^{k-1}(\prxX^{k-1})\,,\\
  a_{i,\ell,k}
  &= \ind\{i=k\le\ell\} c_k(\epsilon)\,,
  \end{align*}
for $0\le k\le\ell+1$. It follows 
using \eqref{eq:state-evolution-setup}
and Lemma~\ref{lem:orthogonal-increments} that
  \[
  \E\Big[(U^{\ell+1,k})^2\Big]
  =\ind\{k\le \ell\}\bigg\{
  \alpha (a_k)^2 \E 
  \Big[ (M^\ell-M^{k-1})^2\Big]
  + c_k(\epsilon)^2\bigg\}
  +\ind\{k\ge1\}
  \alpha \delta \E \Big[
    b^{k-1}
      (\prxX^{k-1})^2\Big]\,.
  \]
Recall also from \eqref{eq:state-evolution-setup} that $U^{\ell_1,k_1}$ and $U^{\ell_2,k_2}$ are uncorrelated for $k_1\ne k_2$.
Summing the above calculation over $k$, and again recalling 
Lemma~\ref{lem:orthogonal-increments}, gives
  \begin{align}
  \nonumber
  &\E\Big[(U^{\ell+1}-U^\ell)^2\Big]
  = \E\Big[(U^{\ell+1})^2\Big]-\E\Big[(U^\ell)^2\Big]
  = \sum_{k=0}^{\ell+1}
  \E[(U^{\ell+1,k})^2]
  -\sum_{k=0}^\ell
  \E[(U^{\ell,k})^2]\\ 
  \nonumber
  &\qquad
  = \sum_{k=0}^\ell
  \bigg\{
  \alpha (a_k)^2 \E 
  \Big[ (M^\ell-M^{k-1})^2\Big]
  + c_k(\epsilon)^2\bigg\}
  + \sum_{k=1}^{\ell+1} 
  \alpha \delta\, \E\Big[ 
  b^{k-1}
  (\prxX^{k-1})^2\Big]\\ 
  \nonumber
  &\qquad\qquad- 
  \sum_{k=0}^{\ell-1}
  \bigg\{
  \alpha (a_k)^2 \E 
  \Big[ (M^{\ell-1}-M^{k-1})^2\Big]
  + c_k(\epsilon)^2\bigg\}
  - \sum_{k=1}^{\ell} 
  \alpha \delta\, \E\Big[ b^{k-1}(\prxX^{k-1})^2\Big] 
  \\ 
  \nonumber
  &\qquad=\alpha\,
  \E \Big[(M^\ell-M^{\ell-1})^2\Big]
  \sum_{k=0}^\ell (a_k)^2
  +\alpha\delta\,\E\Big[ b^\ell(\prxX^\ell)^2\Big]
  + c_\ell(\epsilon)^2 
  \\
  &\qquad=
  \alpha p^\ell
  \E \Big[(M^\ell-M^{\ell-1})^2\Big]
  + \alpha\delta\,\E\Big[ b^\ell(\prxX^\ell)^2\Big]
  + c_\ell(\epsilon)^2\,.
  \label{e:state.evol.var.increment.u}
  \end{align}
\textbf{Note that \eqref{e:state.evol.var.increment.u} is precisely what one would naturally expect from the form of the iteration \eqref{eq:u-amp-def}, and we highlight this calculation for later use.}
Recalling \eqref{e:state.evol.x.y.increments}, we can further simplify
  \beq\label{e:variance.Delta.M.calculation}
  \E\Big[(M^{\ell+1}-M^\ell)^2\Big]
  = \frac{\E[(\sigma^\ell(\prxX^\ell)-1)^2]}{(p^\ell)^2}
  \E\Big[(V^{\ell+1}-V^\ell)^2\Big]
  = 
  \frac{\E[(\sigma^\ell(\prxX^\ell)-1)^2]}{(p^\ell)^2}
  \delta 
  s_\delta(q_0+\ell\delta)^2\,,
  \eeq
where the variance of the $V$-increments was computed in the proof of Proposition~\ref{prop:Ising-IAMP-analysis}. Now recall the definition of $c_\ell(\epsilon)$ from \eqref{e:c.ell.budget.cost}; recall also the definition of $\budget(q_0+\ell\delta)$ and $\Cost(q_0+\ell\delta)$ from
\eqref{eq:budget-constraint-ising}.
We then have 
  \begin{align*}
  &\E\Big[(U^{\ell+1}-U^\ell)^2\Big]
  =
  \frac{\alpha \delta 
  s_\delta(q_0+(\ell-1)\delta)^2}
    {(p^{\ell-1})^2 / p^\ell}
  \E\Big[(\sigma^{\ell-1}
  (\prxX^{\ell-1})-1)^2\Big]
  \\
     &\qquad\qquad 
    + \alpha \delta 
  \E\Big[b^\ell(\prxX^\ell)^2\Big]
  + 
  \Big(
  \E w^\ell(Y(q_0+\ell\delta)) \Big)^2
  + \alpha\delta \epsilon
    \\
    & \qquad\qquad
  - \alpha\delta\,\E\bigg[
  b^\ell(X(q_0+\ell\delta))^2 + 
  \bigg(\frac{
  s(q_0+\ell\delta)^2}{p^\ell}\bigg)
  (\sigma^\ell(X(q_0+\ell\delta))-1)^2\bigg]\\
  &\quad
 = \delta
  \Big(
  \E w^\ell(Y(q_0+\ell\delta)) \Big)^2
  +\alpha\delta\epsilon
  +\alpha\delta\, \E \bigg[
    b^\ell(\prxX^\ell)^2 
    -
    b^\ell(X(q_0+\ell\delta))^2
    \bigg]\\
  &\qquad
  + \alpha\delta\bigg\{
  \frac{s_\delta(q_0+(\ell-1)\delta)^2}
    {(p^{\ell-1})^2 / p^\ell}
  \E\Big[(\sigma^{\ell-1}(\prxX^{\ell-1})-1)^2\Big]
  - \frac{s(q_0+\ell\delta)^2}
    {p^\ell}
  \E\Big[
  (\sigma^\ell(X(q_0+\ell\delta))-1)^2
  \Big]
  \bigg\}\,.
  \end{align*}
Now recall the assumption that both $p'$ and $p''$ are bounded by $L$, while the functions $b$ and $\sigma$ are $L$-Lipschitz in both time and space. We can use these assumptions together with the result of Proposition~\ref{prop:IAMP-to-SDE-limit} to bound the above. The main contributions to the error are
  \begin{align*}
  &\alpha\delta
  \bigg| \E \bigg[
    b^\ell(\prxX^\ell)^2 
    -
    b^\ell(X(q_0+\ell\delta))^2
    \bigg]\bigg| \\
  &\qquad=
  \alpha\delta\,
  \bigg|\E\bigg[
  \Big(b^\ell(\prxX^\ell)
    +
    b^\ell(X(q_0+\ell\delta))\Big)\Big(b^\ell(\prxX^\ell)
    -
    b^\ell(X(q_0+\ell\delta))\Big)
    \bigg]\bigg|
  \le C\delta^{3/2}\,,
  \end{align*}
and similarly with the terms involving $\sigma$. 
Combining the above estimates gives
        \beq
        \label{eq:U-squared-bound}
        -C\delta^{1/2}
        \leq 
        \frac{\E[(U^{\ell+1}-U^\ell)^2
  ]}{\delta}
  - \Big(\E w^\ell
    (Y(q_0+\ell\delta))\Big)
    ^2
  \le 
  C\delta^{1/2} + C_0\epsilon\,,
        \eeq
and this implies \eqref{e:U.stdev.Ew}.
Now we return to the definitions of $\gamma$ and $\prxgamma$ from
\eqref{eq:proxy-C-ell-delta-def}.
Recall from assumptions (\ref{i:IAMP-main-coefs})--(\ref{i:IAMP-main-diffusivity}) of Theorem~\ref{thm:IAMP-main} that $|\E[ w^\ell(Y(q_0+\ell\delta))^2] -1| \le \epsilon$, and that $w$ is bounded by $L$ and spatially $L$-Lipschitz. It follows that
        \beq  
        \label{eq:w-close-for-prxY}
        \begin{aligned}
  &\bigg|\E\Big[ w^\ell(\prxY^\ell)^2
  - w^\ell(Y(q_0+\ell\delta))^2\Big]\bigg|
  \\
        &= \bigg|
  \E \Big[
  \Big(
  w^\ell(\prxY^\ell)
  + w^\ell(Y(q_0+\ell\delta))
  \Big)
  \Big(w^\ell(\prxY^\ell)
  - w^\ell(Y(q_0+\ell\delta))
  \Big)
  \Big]
  \bigg|
  \\
  &\lesssim C
  \E \Big|
  w^\ell(\prxY^\ell)
  - w^\ell(Y(q_0+\ell\delta))
  \Big|
  \le C \bigg\{
  \E\Big[
  (\prxY^\ell-Y(q_0+\ell\delta))^2\Big]
  \bigg\}^{1/2}
  \le C\delta^{1/2}\,,
  \end{aligned}
        \eeq
where the last step is by Proposition~\ref{prop:IAMP-to-SDE-limit}. Substituting these bounds into the definition \eqref{eq:proxy-C-ell-delta-def} gives the claimed estimates for $\gamma$ and $\prxgamma$.
\end{proof}
\end{lem}

\begin{ppn}
\label{ppn:state-evolution-valid-for-IAMP} In the setting of Theorem~\ref{thm:IAMP-main}\ref{i:IAMP-main-main}, the $\IAMP_N(\sigma,b,w,p,\zeta,\zeta^\Ising,\delta)$ iterates 
	\[(\bx^k, \prxbx^k, \by^k,\prxby^k,\bu^k, \bv^k,\bm^k)\] described by \eqref{eq:w-sigma-b-setup}--\eqref{eq:x-amp-def} 
converge in the manner of Proposition~\ref{prop:graph-AMP-cor} to the state evolution
 limits \[(X^k, \prxX^k,Y^k,\prxY^k,U^k,V^k,M^k)\] described by \eqref{e:state.evol.initialization}--\eqref{e:state.evol.x.y.increments}.
\end{ppn}

\begin{proof}
We first argue that the scaling constants $\Gamma^\ell$ and $\gamma^\ell$, as defined in \eqref{eq:proxy-C-ell-delta-def}, are bounded away from zero. Recall from assumptions (\ref{i:IAMP-main-coefs})--(\ref{i:IAMP-main-diffusivity}) of Theorem~\ref{thm:IAMP-main} that $|\bbE[w_t(Y(t))^2] - 1| \le \epsilon$  and $w_t\in [0,L]$. Therefore
    \begin{equation}
    \label{eq:E-wt-nonzero}
    \bbE[ w_t(Y(t))]
    \ge \frac{\bbE[w_t(Y(t))^2]}{L}
    \ge \frac{1-\epsilon}{L}
    \ge \frac{1}{2L}. 
    \end{equation} 
It follows from the left-hand estimate in \eqref{eq:U-squared-bound} above that
    \beq\label{e:U-increment-var-lb}
    \prxgamma^\ell
    =\frac{\bbE[(U^{\ell+1}-U^{\ell})^2]}{\delta}
    \geq 
    \frac{1}{4L^2}
    -
    C(L,\alpha)\delta^{1/2}
    \geq 
    \frac{1}{5L^2}
    \eeq
since $\delta$ is small depending on $L$. To show a similar bound for
$\gamma^{\ell}$, it suffices to note that $\bbE[w^\ell(\prxY^{\ell})^2]$ is bounded away from zero, which follows by combining
\eqref{eq:w-close-for-prxY} with the assumption (\ref{i:IAMP-main-diffusivity}) from Theorem~\ref{thm:IAMP-main}. 

Having shown $\prxgamma^{\ell}$ and $\gamma^{\ell}$ are bounded away from zero, we now see from \eqref{e:h.l.k.explicit} and \eqref{e:f.l.k.explicit} (and the fact that $p$ is bounded away from $0$) that the functions $h_{\ell,k}$ and $f_{\ell,k}$ are compositions of Lipschitz operations and multiplication by $b,\sigma,w$.
The coefficients
$b,\sigma,w$ were assumed in Theorem~\ref{thm:IAMP-main} to be uniformly bounded and spatially Lipschitz. Thus Corollary~\ref{cor:psi-state-evolution} justifies the state evolution convergence in $\bbW_2$ of $\by^{\ell}$ to $Y^{\ell}$, and similarly for all the other quantities.
\end{proof}

\begin{proof}[\hypertarget{proof:p.IAMP-to-SDE-limit.Ax.y}{Proof of Proposition~\ref{prop:IAMP-to-SDE-limit.Ax.y}}]
Recalling \eqref{eq:v-amp-def}, we have
  \beq\label{eq:ons-hopefully}
  \bA_\ell\by^\ell-\bx^\ell
  =\ons_{\bv,\ell-1}-(\bx^\ell-\bv^\ell)\,,
  \eeq
and we will show that this quantity is small. Recall that $\ons_{\bv,\ell-1}$ is defined by \eqref{e:v.combined.ons} and \eqref{eq:IAMP-fully-general-fixed.Ons}, where the functions $h$ and $f$ are given explicitly by \eqref{e:h.l.k.explicit} and \eqref{e:f.l.k.explicit}. Recall that we abbreviate $w^\ell \equiv w_{q_0+\ell\delta}$, etc. It follows from \eqref{e:h.l.k.explicit} that
  \[
  \frac{\partial h_{\ell,k}}{\partial u^{j,k}}
  = a_k
  \sum_{s=0}^{\ell-1}
  \frac{u^{s+1}-u^s}{\gamma^s}
  \frac{\partial w^s}{\partial u^{j,k}}
  + 
  a_k
  \sum_{s=0}^{\ell-1}
  \bigg(\frac{\partial u^{s+1}}{\partial u^{j,k}}
  -\frac{\partial u^s}{\partial u^{j,k}}\bigg)
  \frac{w^s}{\gamma^s}
  \,,
  \]
where all the functions are evaluated on the $U$ variables. The first term has mean zero by Lemma~\ref{lem:orthogonal-increments}. For the second term, recall from \eqref{e:u.ellp1.decomp} that $\bu^{s+1}$ is the sum of $\bu^{s+1,k}$ for $0\le k\le s+1$, so
  \[ \frac{\partial u^{s+1}}{\partial u^{j,k}}
  = \ind\{s+1=j\}
  \ind\{k\le s+1\}
  = \ind\{s+1=j\}
  \ind\{k\le j\}
  \,.
  \]
Therefore the expected partial derivative of $h$ is given by
  \begin{align*}
  \E \bigg[\frac{\partial h_{\ell,k}}{\partial u^{j,k}}(U[[\ell]])\bigg]
  &=\ind\{k\le j\}
  a_k
  \sum_{s=0}^{\ell-1}
  \frac{\bbE [w^s]}{\gamma^s}
  \bigg\{
  \ind\{s=j-1\}
  -\ind\{s=j\}
  \bigg\} \\
  &= \ind\{k\le j\}
  a_k
  \bigg\{
  \frac{\bbE [w^{j-1}]}{\gamma^{j-1}}
  \ind\{1\le j\le\ell\}
  -\frac{\bbE [w^j]}{\gamma^j}
  \ind\{0\le j\le\ell-1\}
  \bigg\}\,.
  \end{align*}
Substituting into
\eqref{eq:IAMP-fully-general-fixed.Ons}
(and recalling from \eqref{eq:IAMP-fully-general-fixed} that $f_{j,k}$ is only non-zero for $0\le k\le j+1$) gives
  \begin{align*}
  \ons_{\bv,\ell-1,k}
  &=a_k \sum_{j=0}^\ell
  \ind\{k\le j\}
  \bigg\{
  \frac{\bbE [w^{j-1}]}{\gamma^{j-1}}
  \ind\{1\le j\le\ell\}
  -\frac{\bbE [w^j]}{\gamma^j}
  \ind\{0\le j\le\ell-1\}
  \bigg\}
  f_{j-1,k}\\
  &= a_k
  \sum_{j=0}^{\ell-1}
  \frac{\E w^j}{\gamma^j}
  (f_{j,k}-f_{j-1,k})
  = a_k
  \sum_{j=k}^{\ell-1}
  \frac{\E w^j}{\gamma^j}
  (f_{j,k}-f_{j-1,k})
  + a_k 
  \frac{\E w^{k-1}}{\gamma^{k-1}}
  f_{k-1,k}\,.
  \end{align*}
Up to now we have not used the explicit definition~\eqref{e:f.l.k.explicit} of $f_{\ell,k}$. Recall from \eqref{eq:ons-hopefully} that our goal is for $\ons_{\bv,\ell-1}$ to be close to
$\bx^\ell-\bv^\ell$, which by \eqref{eq:x-amp-def} can be rewritten as
  \[\bx^\ell-\bv^\ell
        =
        \bx^0-\bv^0 + 
  \sum_{j=1}^\ell\bigg\{
  (\bsig^{j-1}-1)
  \odot(\bv^j-\bv^{j-1})
  +\delta p'(q_0+(j-1)\delta)^{1/2} \bb^{j-1}
  \bigg\}\,.
    \]
\textbf{This is the main reason for the precise definition~\eqref{e:f.l.k.explicit} of $f_{\ell,k}$:} it gives
  \[
  f_{j,k}-f_{j-1,k}
  = \frac{a_k(\bsig^{j-1}-1)
    \odot(\bv^j-\bv^{j-1})}
    {p^{j-1}}
  \]
for all $0\le k\le j$, as well as
$f_{k-1,k}
= \ind\{k\ge1\}
\delta^{1/2} \bb^{k-1}$.
Combing \eqref{e:v.combined.ons} and the above calculations gives
  \begin{align*}
  \ons_{\bv,\ell-1}
  &=\sum_{j=1}^{\ell-1}
  \frac{\E w^j}{\gamma^j p^{j-1}}
  \bigg\{\sum_{k=0}^j
  (a_k)^2 
  \bigg\}
  (\bsig^{j-1}-1)
    \odot(\bv^j-\bv^{j-1})
  + \sum_{k=1}^\ell
  \frac{\E w^{k-1}}{\gamma^{k-1}}
  (\delta^{1/2} a_k)
  \bb^{k-1} \\
  &= \sum_{j=1}^{\ell-1}
  \frac{(\E w^j) p^j}{\gamma^j p^{j-1}}
  (\bsig^{j-1}-1)
    \odot(\bv^j-\bv^{j-1})
  + \sum_{k=1}^\ell
  \frac{\E w^{k-1}}{\gamma^{k-1}}
  (\delta^{1/2} a_k)
  \bb^{k-1}
  \end{align*}
It follows that 
$\ons_{\bv,\ell-1}-(\bx^\ell-\bv^\ell)
=\textbf{err}^\textup{A}+\textbf{err}^\textup{B} - (\bx^0-\bv^0)$ 
where
  \begin{align*}
  \textbf{err}^\textup{A}
  &\equiv - (\bsig^{\ell-1}-1)
  \odot(\bv^\ell-\bv^{\ell-1})
  + \sum_{j=1}^{\ell-1}
    \bigg(
  \frac{(\E w^j)p^j}
    {\gamma^j p^{j-1}}
    -1\bigg)
    (\bsig^{j-1}-1)
    \odot
    (\bv^j-\bv^{j-1})\,,\\
  \textbf{err}^\textup{B}
  &\equiv \sum_{j=1}^{\ell-1}
  \delta 
  \bigg\{
  \frac{a_j}{\delta^{1/2}}
  \frac{\E w^{j-1}}
    {\gamma^{j-1}}
  - p'(q_0+(j-1)\delta)^{1/2}
  \bigg\}\bb^{j-1}\,.
  \end{align*}
For the $\bx^0$ and $\bv^0$ terms, we recall that $\bx^0=\bzero$, while \eqref{eq:v-amp-def} and \eqref{e:v.combined.ons} imply $\bv^0=\bA_0 \by^0 = a_0 \bJ_0 \by^0$, from which it follows that
	\[
	\frac{\bbE[\|\bv^0\|^2]}{M}
	= (a_0)^2\frac{\bbE[\|\by^0\|^2]}{N} 
	= (a_0)^2 q_0 = p(q_0) q_0
	\le q_0\epsilon \le\epsilon\,,
	\]
using assumption (\ref{i:IAMP-main-p}) of Theorem~\ref{thm:IAMP-main}. Since $p'$ and $p''$ are both bounded by $L$, it follows that
  \[
  \frac{(a_j)^2}{\delta}
  = \frac{p(q_0+j\delta)-p(q_0+(j-1)\delta)}{\delta}
  = p'(q_0+(j-1)\delta) + O(\delta)\,.
  \]
Again recall from assumption~(\ref{i:IAMP-main-coefs}) of Theorem~\ref{thm:IAMP-main} that $b,\sigma,w$ are bounded.
Combining with the result of Lemma~\ref{l:gamma.estimates} gives
  \[
  \frac{\|\textbf{err}^\textup{B}\|}{N^{1/2}}
  \le C\delta 
  (\delta^{1/2}+\epsilon)\sum_{k=0}^{\ell-1}
  \frac{\|\bb^k\|}{N^{1/2}}
  \le  C\delta^{1/4}
  +C_0\epsilon^{1/2}\,.
  \]
Combining 
Proposition~\ref{prop:graph-AMP-cor},
Lemma~\ref{lem:orthogonal-increments},
and Lemma~\ref{l:gamma.estimates} gives
  \[
  \plim_{N\to\infty} \frac{\|\textbf{err}^\textup{A}\|}{N^{1/2}}
  \le C\delta^{1/4}
  +C_0\epsilon^{1/2}\,.
  \]
Taking all the above estimates together proves the claim.
\end{proof}

\begin{proof}[\hypertarget{proof:t.IAMP.main.main}{Proof of Theorem~\ref{thm:IAMP-main}\ref{i:IAMP-main-main}}]
For sufficiently small $\delta=\delta(L,\epsilon)$, consider the map $\tilde\cA(\bG)=\by^{\ubl}$, where $\by^{\ubl}$ is given  by the IAMP iteration \eqref{eq:w-sigma-b-setup}--\eqref{eq:x-amp-def}. By Lemma~\ref{l:wishart}, with $C_0$ given therein, it holds with probability $1-\exp(-cN)$ that 
the Brownian bridge
 $\bG_{\circ}$ satisfies 
   \beq\begin{aligned}
   \label{eq:Brownian-bridge-condition}
   \Big\|\bG_{\circ}(p(q_0+\ell\delta))-\bG_{\circ}(p(q_0+(\ell-1)\delta))\Big\|_{\op}\leq 
    C_0 \big(N^{1/2}+M^{1/2}\big)
    \end{aligned}
    \eeq
for all $1\le \ell\le \ubl$. Then the Lipschitz constant of $\tilde\cA$, for a fixed Brownian bridge, can be shown inductively to be at most $C(L,\epsilon,\delta)$ on the subset 
	\[
	\mathcal{U}
	\equiv \Big\{
	\bG : 
	\|\bG\|_\textup{op}\le
	C_0 (M^{1/2}+N^{1/2})
	\Big\}\,.
	\]
(See for example the discussion of ``standard optimization algorithms'' in \cite{HuangSellke2021} for the inductive argument.) The Kirszbraun extension theorem then gives a globally $C(L,\epsilon,\delta)$-Lipschitz 
$\bar{\cA}$ agreeing with $\tilde\cA$ on the set $\mathcal{U}$. We also know from Lemma~\ref{l:wishart} that $\bG\in\mathcal{U}$ with probability $1-\exp(-cN)$. Lastly we set
  \[
  \cA(\bG)
  = \frac{N^{1/2}}{\E_{\bG'}[\|\bar\cA(\bG')\|^2]^{1/2}}\bar{\cA}(\bG)
  \equiv \bar{c}_N \bar{\cA}(\bG)\,,
  \]
so that $\cA$ satisfies condition~\eqref{it:expectation-condition} of Definition~\ref{d:Lip}. We will argue below that $\bar{c}_N$ is close to $1$, so that $\cA$ will also satisfy condition~\eqref{it:L-lip-condition} of Definition~\ref{d:Lip}. 

To this end, let $\iota'$ be such that $\epsilon \ll \iota' \ll \iota$;
combining with \eqref{eq:IAMP-order-of-limits} gives the overall order of parameters
\beq\label{eq:IAMP-order-of-limits-tmp}
  \delta \ll \frac{1}{L} \ll \epsilon \ll \iota' \ll \iota\,.
\eeq 
First we argue that the map $\cA$ takes values in $\Sigma_N(\iota')$ 
(or $S_N(\iota')$ in the case of the spherical perceptron) 
with probability $1-\exp(-cN)$.
The relaxed domains $\Sigma_N(\iota')$ and $S_N(\iota')$ were defined in \eqref{eq:relaxed-domain}. Recall that $q_0+\ubl\delta=1$. By Proposition~\ref{prop:Ising-IAMP-analysis}\ref{it:Ising-BM-limit}
and
\eqref{eq:tilde-Y-approx-Y},
  \begin{align*}
    &\bbE[(Y^{\ubl}-\prxY^{\ubl})^2]
    =
    \int_{q_0}^1
    \bigg(
    \frac{1}{\bbE[w_{t_{\delta}}(\prxY(t))^2]^{1/2}}-1
    \bigg)^2 \bbE[w_{t_{\delta}}(\prxY(t))^2]\,dt
    +\E[(Y^0-\prxY^0)^2]
    \\
    &\qquad\le 
    \int_{q_0}^1
 \Big|   \bbE[w_{t_{\delta}}(\prxY(t))^2]-1\Big|\,dt
 + \epsilon\,,
    \end{align*}
having used the inequality $(a-b)^2 \le |a^2-b^2|$. Now recall from the definition \eqref{e:state.evol.piecewise} that 
$\prxY(t)$ is piecewise constant, with 
$\prxY(t)=\prxY(t_\delta)$. It follows that 
    \beq
    \label{eq:by-prxby-close}
    \begin{aligned}
   \bbE[(Y^{\ubl}-\prxY^{\ubl})^2]
   &\le
  \int_{q_0}^1
    \bigg|\bbE\Big[\Big(
    w_{t_{\delta}}(\prxY(t_{\delta}))
    -
    w_{t_{\delta}}(Y(t_{\delta}))\Big)
    \Big(
    w_{t_{\delta}}(\prxY(t_{\delta}))
    +
    w_{t_{\delta}}(Y(t_{\delta}))\Big)
    \Big]
    \bigg|\,dt\\
    &\qquad+
    \int_{q_0}^1
    \Big|
    \bbE[w_{t_{\delta}}(Y(t_{\delta}))^2]-1
    \Big|\,dt + \epsilon
    \le C(L,\alpha)\delta^{1/2}
      + C_0\epsilon\,.
   \end{aligned}
    \eeq
In the above, the first integral is at most $C(L,\alpha)\delta^{1/2}$ because assumption~(\ref{i:IAMP-main-coefs}) of Theorem~\ref{thm:IAMP-main} gives that  $w$  is bounded by $L$ and is
$L$-Lipschitz, while the discrepancy between $\prxY(t)$ and $Y(t)$ is bounded by \eqref{eq:Y-diffusive-limit} from Proposition~\ref{prop:IAMP-to-SDE-limit}. The second integral is at most $C_0\epsilon$ using assumption~(\ref{i:IAMP-main-diffusivity}) from Theorem~\ref{thm:IAMP-main}. Combining the above with the state evolution result of 
Proposition~\ref{ppn:state-evolution-valid-for-IAMP} gives
 \[
    \plim_{N\to\infty}
    \frac{\|\by^{\ubl}-\prxby^{\ubl}\|}{N^{1/2}} 
    \leq C(L,\alpha)\delta^{1/4}+C_0\epsilon^{1/2}\,,
    \]
Moreover, Proposition~\ref{ppn:state-evolution-valid-for-IAMP} gives the state evolution convergence of $\prxby^{\ubl}$ to $\prxY^{\ubl}$, and by Proposition~\ref{prop:IAMP-to-SDE-limit},
  \[
    \bbE\lt[\Big(\prxY^{\ubl} - Y(1)\Big)^2\rt] \le C(L,\alpha) \delta\,.
  \]
  By \eqref{eq:IAMP-order-of-limits-tmp}, all these errors are $\ll \iota'$.
For the \textbf{spherical perceptron}, the SDE \eqref{e:IAMP.Y.SDE} with $w_t\equiv1$ implies $\E[Y(1)^2]=1$, and it follows that
$\tilde{\cA}$ takes values in the relaxed spherical domain $S_N(\iota')$ with high probability. For the \textbf{Ising perceptron}, assumption~(\ref{i:IAMP-main-endpt}) on the law of $Y(1)$ implies that 
$\tilde{\cA}$ takes values in $\Sigma_N(\iota')$ with high probability.
Since $\cA$ and $\tilde{\cA}$ agree with probability $1-e^{-cN}$, it follows that $\cA$ also takes values in $\Sigma_N(\iota')$ (or $S_N(\iota')$) with high probability.
Since $\cA$ is Lipschitz, a routine concentration argument then shows that $\cA(\bG) \in \Sigma_N(\iota)$ with probability $1-e^{-cN}$.

Next we argue that $\bar{c}_N=1+O(\iota')$. First note that 
  \[
  \E\Big[\|\bar{\cA}(\bG)\|^2\Big] -o(1)
  =\E\Big[\|\bar{\cA}(\bG)\|^2;\mathcal{U}\Big] 
  =\E\Big[\|\tilde{\cA}(\bG)\|^2;\mathcal{U}\Big]\,,
  \]
where the first equality holds because $\bar{\cA}$ is Lipschitz, and the second holds because $\bar{\cA}$ and $\tilde{\cA}$ agree on $\mathcal{U}$. It follows from the above that $|\|\tilde{\cA}(\bG)\|-N^{1/2}| \le N^{1/2} \iota'$ with very high probability, so the same holds for $\bar{\cA}(\bG)$. It follows from Lipschitz concentration that 
  \[\bar{c}_N
  = \frac{N^{1/2}}{\E[\|\bar{\cA}(\bG)\|^2]^{1/2}}
  = 1+O(\iota')\,.\]
Since $\bar{\cA}$ takes values in $\Sigma_N(\iota')$ (or $S_N(\iota')$) with very high probability, it follows that $\cA=\bar{c}_N\bar{\cA}$ takes values in $\Sigma_N(\iota')$ (or $S_N(\iota')$) with very high probability. This proves \eqref{e:IAMP-main-main-coords}. 

Finally we show $\bbW_2(\mu_\bG(\cA_N(\bG,\bg^\aux)),\Law(X(1)))\leq \iota$ with probability $1-e^{-cN}$. Again recall from 
\eqref{eq:A-ell}  that
$q_0+\ubl\delta = 1$, so $\bA_{\ubl}=\bA =\bG/N^{1/2}$.
    Thus Proposition~\ref{prop:IAMP-to-SDE-limit.Ax.y} directly yields
    \[
    \plim_{N\to\infty}
    \frac{\|\bA\by^{\ubl}-\bx^{\ubl}\|}{N^{1/2}}
    \leq 
    C\delta^{1/4}+C_0\epsilon^{1/2}.
    \] 
Arguing analogously to \eqref{eq:by-prxby-close} shows, similarly to the above estimate for $\by^{\ubl}-\prxby^{\ubl}$, that
\[
    \plim_{N\to\infty}
    \frac{\|\bx^{\ubl}-\prxbx^{\ubl}\|}{N^{1/2}} 
    \leq C(L,\alpha)\delta^{1/4}+C_0\epsilon^{1/2}\,,
    \]
Proposition~\ref{ppn:state-evolution-valid-for-IAMP} gives the convergence of $\prxbx^{\ubl}$ to $\prxX^{\ubl}$ in the manner of Proposition~\ref{prop:graph-AMP-cor}, where $\prxX^{\ubl}$ is a good approximation to $X(q_0+\ubl\delta)=X(1)$ by 
Proposition~\ref{prop:IAMP-to-SDE-limit}. Combining these estimates shows
\[
  \bbW_2(\mu_\bG(\cA_N(\bG,\bg^\aux)),
   \Law(X(1))) \le O(\iota')
  \le \frac{\iota}{2}
\]
with high probability.
Since $\cA$ is Lipschitz, a routine concentration argument implies 
\[\bbW_2(\mu_\bG(\cA_N(\bG,\bg^\aux)), \Law(X(1))) \le \iota\] with probability $1-e^{-cN}$. This proves \eqref{e:IAMP-main-main-inner-prods}. 

The remaining conclusions also follow from the Lipschitzness of $\cA$.
Indeed, since $\bbW_2(\cdot,\cdot)^2$ is convex,
\begin{align}
  \label{e:IAMP-main-main-inner-prods-averaged-jensen}
  \bbW_2\Big(\mu(\cA_N), \Law(X(1)) \Big)^2
  &\le \bbE \Big[\bbW_2\Big(\mu_\bG(\cA_N(\bG,\bg^\aux)), \Law(X(1)) \Big)^2\Big]\,, \\
  \label{e:IAMP-main-main-coords-averaged-jensen}
  \bbW_2\Big(\mu^{\Ising}(\cA_N), \cP(\{\pm 1\}) \Big)
  &\le \bbE \Big[\bbW_2\Big(\EmpDist(\cA(\bG,\bg^{\aux})), \cP(\{\pm 1\})\Big)^2\Big]\,.
\end{align}
The conclusions \eqref{e:IAMP-main-main-coords} and \eqref{e:IAMP-main-main-inner-prods} imply that the quantities inside the right-hand side expectations are bounded by $\iota^2$ with probability $1-e^{-cN}$.
A direct calculation shows that the functions
\begin{align*}
  (\bG,\bg^\aux) &\mapsto \bigg(1 + \bbW_2\Big(\mu_\bG(\cA_N(\bG,\bg^\aux)), \Law(X(1)) \Big)^2\bigg)^{1/4}\,, \\
  (\bG,\bg^\aux) &\mapsto \bigg(1 + \bbW_2\Big(\EmpDist(\cA(\bG,\bg^{\aux})), \cP(\{\pm 1\})\Big)^2\bigg)^{1/2}\,,
\end{align*}
are $O(1/N^{1/2})$-Lipschitz, and thus these quantities are subgaussian with variance proxy $O(1/N)$.
This provides enough control to integrate the tails of the right-hand side expectations in \eqref{e:IAMP-main-main-inner-prods-averaged-jensen}, \eqref{e:IAMP-main-main-coords-averaged-jensen}, and implies that they are $\le 2\iota^2$.
Adjusting $\iota$ proves \eqref{e:IAMP-main-main-coords-averaged} and \eqref{e:IAMP-main-main-inner-prods-averaged}.
\end{proof}

\subsection{Centering of general IAMP algorithms}
\label{subsec:center-IAMP}

In this subsection we give the \hyperlink{proof:t.IAMP.main.centered}{proof of Theorem~\ref{thm:IAMP-main}\ref{i:IAMP-main-centered}}.
Here, we take as given $(b,\sigma,w,p,\zeta,\zeta^\Ising)$ satisfying assumptions \eqref{i:IAMP-main-coefs}--\eqref{i:IAMP-main-init} from Theorem~\ref{thm:IAMP-main}, and additionally assume $\zeta^\Ising$ is symmetric and the functions $w_t = w(t,\cdot)$ are even.
We will exhibit a modified IAMP algorithm $\cA^\dagger$ achieving the guarantees of Theorem~\ref{thm:IAMP-main}\ref{i:IAMP-main-main} and $\bbE \cA(\bG,\bg^\aux) = \bzero$.

We introduce a new auxiliary variable $\bg^{\init} \sim \cN(0, I_N)$.
\textbf{For this algorithm, we will treat the auxiliary gaussian data $\bg^\aux$ from \eqref{eq:algorithms-as-maps} as consisting of $\bg^{\init}$, the vectors $\bar{\bg}^\ell$ from before, and the Brownian bridge $\bG_\circ$ (evaluated at the finite set of times $\{q_0+\ell\delta:0\le\ell\le\ubl\}$)}. 
Thus we will provide an IAMP algorithm $\cA^{\dagger}$ satisfying
\[
\bbE\bigg[\cA^{\dagger}\Big(\bG,
  \bg^\aux
  \equiv\big(\bG_{\circ},
    \bg^{\init},\bar\bg^0,\dots,\bar\bg^{\ubl}\big)
  \Big)
  \bigg]= \bzero\,,
\]
where the map is $C(L,\epsilon,\delta)$-Lipschitz with respect to $(\bG,\bg^\aux)$, and the expectation is over $(\bG,\bg^\aux)$. Let
  \[
  N^\aux
  \equiv MN(\ubl+1)
    + N(\ubl+2)
  \]
denote the total dimension of $\bg^\aux$. First, we argue that there exists a  Lipschitz odd function $\psi:\bbR\to\bbR$ such that if $Z \sim \cN(0,1)$, then the law of $\psi(Z)$ approximates the measure $\zeta^\Ising$:

\begin{lem}
\label{lem:even-manual-init}
Let $Z$ be a standard gaussian random variable, and suppose $Y(q_0)\sim\zeta^\Ising$ satisfies assumption~\eqref{i:IAMP-main-init} of Theorem~\ref{thm:IAMP-main}. Then for any $\eta>0$ there exists an odd function $\prxpsi:\bbR\to\bbR$ with bounded Lipschitz norm $\|\prxpsi\|_{\textup{Lip}}\leq \overline{L}(\eta,\zeta^{\Ising})$ such that $\bbW_2\big(\prxpsi(Z),\zeta^{\Ising}\big)\leq \eta$. Additionally, we have $|\bbE[\prxpsi(Z)^2] - q_0| \le 4\epsilon^{1/2}$.

\begin{proof}
    Using a standard quantile coupling, there exists an odd measurable function $\hat\psi\in L^2(\bbR)$ such that the law of $\hat\psi(Z)$ is precisely $\zeta^{\Ising}$. Since bounded Lipschitz functions are dense in $L^2$, there exists a bounded Lipschitz function $\tilde\psi:\bbR\to\bbR$ with $\bbW_2\big(\tilde\psi(Z),\zeta^{\Ising}\big)\leq \min\{\eta, \epsilon^{1/2}\}$.
The first assertion of the lemma immediately follows by taking $\prxpsi=\tilde\psi$.  For the second assertion, recall that under assumption~\eqref{i:IAMP-main-init} of Theorem~\ref{thm:IAMP-main}, the random variable $Y(q_0)\sim\zeta^\Ising$ satisfies
	\[\Big|\E[Y(q_0)^2-q_0]
	\Big|\le\epsilon\,.\]
Under the $\bbW_2$-optimal coupling of $\prxpsi(Z)$ with $Y(q_0)$, we have
	\begin{align*}
	&\bigg|\E[\prxpsi(Z)^2]
	-\E[Y(q_0)^2] \bigg|
	=\bigg|\E\Big[
	(\prxpsi(Z)-Y(q_0))
	(\prxpsi(Z)-Y(q_0)
		+2Y(q_0))
	\Big]\bigg|\\
	&\qquad=\bigg| \E\Big[
	(\prxpsi(Z)-Y(q_0))^2
	\Big]
	+ 2\E\Big[
	(\prxpsi(Z)-Y(q_0)) Y(q_0)\Big]
	\bigg|
	\le \epsilon + 2 \epsilon^{1/2}
		\E[Y(q_0)^2]^{1/2} \le 3\epsilon^{1/2}\,.
	\end{align*}
Combining these gives
	\[
	\Big|\E[\prxpsi(Z)^2-q_0]
	\Big|
	\le 4\epsilon^{1/2}\,,
	\]
as claimed.\end{proof}
\end{lem}

\begin{dfn}[centered IAMP] \label{d:centered.iamp}
  Let $\prxpsi$ be as in Lemma~\ref{lem:even-manual-init}, where $\eta$ is sufficiently small depending on $(\delta,L,\epsilon)$, and define
  \[
    \psi(z)=\frac{(q_0)^{1/2} \prxpsi(z)}{\bbE[\prxpsi(Z)^2]^{1/2}}\,,
  \]
  with $\psi(z) = (q_0)^{1/2} z$ if $\prxpsi(Z) = 0$ almost surely.
  In either case, $\bbE[\psi(Z)^2] = q_0$. 
  We now initialize the iteration with the vectors
  \beq
  \label{eq:mean-0-initialization}
  \begin{aligned}
  \by^{\dagger,0}
  &\equiv \psi(\bg^{\init})\in\R^N\,,\\
  \prxby^{\dagger,0}
  &\equiv\prxpsi(\bg^{\init})\in\R^N\,,
  \end{aligned}
  \eeq
  and $\bx^{\dagger,0}=\bx^0\in\R^M$, $\prxbx^{\dagger,0}=\prxbx^0\in\R^M$ exactly as before. We then follow 
  \textbf{exactly the same} IAMP iteration, as specified in \eqref{eq:proxy-y-amp-def}--\eqref{eq:x-amp-def}.
  These modified iterates are denoted $\by^{\dagger,k}$, $\prxby^{\dagger,k}$, etc., with state evolution limits $Y^\dagger(t)$, $\prxY^\dagger(t)$, etc.
  Finally, let $\tilde{\cA}^\dagger(\bG,\bg^\aux)\equiv \by^{\dagger,\ubl}$.
\end{dfn}
The key point of this construction is the following symmetry property which ensures that the resulting algorithm is centered. We emphasize that we have this symmetry exactly in $\bbR^N$ and $\bbR^M$, not just in the state evolution limit.

\begin{lem}
\label{lem:symmetrized-amp-is-symmetric}
Fix $(b,\sigma,w,p,\zeta,\zeta^\Ising)$ satisfying the conditions of Theorem~\ref{thm:IAMP-main}\ref{i:IAMP-main-centered}. 
Fix any vectors $\bx^{\dagger,0},\prxbx^{\dagger,0}\in\bbR^M$, and let $\tilde{\cA}^\dagger$ be the centered IAMP specified by Definition~\ref{d:centered.iamp}. 
Then, for $C\equiv C(L,\epsilon,\delta)$, 
there exists a $C$-Lipschitz map $\cA^\dagger$ which agrees
with $\tilde{\cA}^\dagger$ with probability $1-\exp(-cN)$, and satisfies
	    \[
    \cA(-\bG,-\bg^\aux)
    =
    -\cA(\bG,\bg^\aux)\,.
    \]
In particular, it follows that $\bbE[\cA(\bG,\bg^{\aux})~|~\bx^{\dagger,0},\prxbx^{\dagger,0}]=\bzero\in\bbR^N$.
\begin{proof}
    Let $(\hat\by^{\dagger,k},\hat\prxby^{\dagger,k},\hat\bx^{\dagger,k},\hat\prxbx^{\dagger,k})$ be the IAMP iterates corresponding to $(-\bG,- \bg^{\aux})$.
    We claim that for all $0\leq k\leq \ubl$ we have
    \begin{align*}
    (\hat\prxby,\hat\by,\hat\bu,\hat\bw)^{\dagger,k}
    &=
    (-\prxby,-\by,-\bu,\bw)^{\dagger,k}
      \,\\
    (\hat\prxbx,\hat\bx,\hat\bv,
    \hat\bm,\hat\bsig,\hat\bb)^{\dagger,k}
    &=
    (\prxbx,\bx,\bv,
    \bm,\bsig,\bb)^{\dagger,k}\,.
    \end{align*}
Indeed, the base case for the first statement above holds easily since
$w$ is even by the assumptions of Theorem~\ref{thm:IAMP-main} \eqref{i:IAMP-main-centered}, the functions
 $\psi$ and $\prxpsi$ are odd (Lemma~\ref{lem:even-manual-init}), 
the vectors $(\by,\prxby,\hat{\by},\hat{\prxby})^{\dagger,0}$ are defined by \eqref{eq:mean-0-initialization},  and  we defined $\bu^{\dagger,0}=\hat{\bu}^{\dagger,0}=\bzero\in\R^M$.
The remaining statements follow by 
induction, using the formulas \eqref{eq:proxy-y-amp-def}--\eqref{eq:x-amp-def}. 

The above shows the desired symmetry for the mapping $\tilde{\cA}^\dagger$, and we must now transfer it to $\cA^\dagger$. This proceeds similarly as in the \hyperlink{proof:t.IAMP.main.main}{proof of Theorem~\ref{thm:IAMP-main}\ref{i:IAMP-main-main}}, using the Kirzsbraun extension theorem, except that we now need to ensure that the symmetry requirement is satisfied. 
First, we have by easy induction that $\tilde{\cA}^\dagger$ is $C$-Lipschitz on a subset $\mathcal{U}^{\dagger}\subseteq \bbR^{MN+N_{\aux}}$, similarly to the \hyperlink{proof:t.IAMP.main.main}{proof of Theorem~\ref{thm:IAMP-main}\ref{i:IAMP-main-main}}. To be precise, 
    \[\mathcal{U}^{\dagger}
    =\bigg\{
    \begin{array}{c}
    (\bG,\bg^\aux) : 
    \textup{$\bG_\circ$
    satisfies \eqref{eq:Brownian-bridge-condition},
    $\|\bG\|_\textup{op}\le C_0
    (M^{1/2}+N^{1/2})$,}\\
    \textup{$\|\bar\bg^{\ell}\|\leq 2\sqrt{N}$
    for all $0\le\ell\le\ubl$, and
    $\|\bar\bg^{\init}\|\leq 2\sqrt{N}$}
    \end{array}
    \bigg\}\,.
    \]
In particular, this choice of $\mathcal{U}^{\dagger}$ is origin-symmetric. 
(Again, see e.g.\ the discussion on ``standard optimization algorithms'' in \cite{HuangSellke2021} for the Lipschitz bound on such a domain.)

Thus, let $S=\{\bg_1,\bg_2,\dots\}$ be a countable dense set in $\bbR^{MN+N_{\aux}}$.
    For each $\bg_i$, we assume $\tilde\cA^{\dagger}$ has been extended to a $C$-Lipschitz function on $\mathcal{U}^{\dagger,i}\equiv \mathcal{U}^{\dagger}\cup \{\bg_1,-\bg_1,\dots,\bg_i,-\bg_i\}$.
    Then by the Kirzsbraun extension theorem, there exists a $C$-Lipschitz extension to $\mathcal{U}^{\dagger,i}\cup \{\bg_{i+1}\}$.
    Since 
  $\mathcal{U}^{\dagger,i}$ is symmetric,  
the unique odd extension to $\mathcal{U}^{\dagger,i+1}$ is still $C$-Lipschitz. Continuing in this way yields a $C$-Lipschitz extension to $ \mathcal{U}^\dagger\cup S\cup -S$. This set is dense by definition of $S$, so the result follows by continuously extending to all of $\R^{MN+N_{\aux}}$.
\end{proof}
\end{lem}

Finally we show centered IAMP behaves the same as our original IAMP algorithm.

\begin{proof}[\hypertarget{proof:t.IAMP.main.centered}{Proof of Theorem~\ref{thm:IAMP-main}\ref{i:IAMP-main-centered}}]
Let $\hat\epsilon=4\epsilon^{1/2}\geq \epsilon$.
We claim that $(b,\sigma,w,p,\zeta,\Law(\prxpsi(Z))$ satisfies assumptions (\ref{i:IAMP-main-coefs})--(\ref{i:IAMP-main-init}) of Theorem~\ref{thm:IAMP-main} with parameters $(\alpha,q_0,L,\hat{\epsilon})$ in place of $(\alpha,q_0,L,\epsilon)$.

First, Lemma~\ref{lem:even-manual-init} shows $|\bbE[\prxpsi(Z)^2] - q_0| \le \hat\epsilon$.
Moreover, since $\eta$ is sufficiently small, a simple Gr{\"o}nwall estimate shows that the SDE solution $Y$ does not change much under small perturbations of $Y_0$, so 
the other assumptions of Theorem~\ref{thm:IAMP-main} are satisfied. Clearly, the order of limits in \eqref{eq:IAMP-order-of-limits} is still valid with $\epsilon$ replaced by $\hat\epsilon$. 

The result now follows directly from  Lemma~\ref{lem:symmetrized-amp-is-symmetric} together with Theorem~\ref{thm:IAMP-main}\ref{i:IAMP-main-main}. Indeed one can just treat the initializations \eqref{eq:mean-0-initialization} as given (i.e. forget $\bg^{\init}$); $\by^{\dagger,0},\prxby^{\dagger,0}$ then have i.i.d.\ coordinates exactly as required in Theorem~\ref{thm:IAMP-main}\ref{i:IAMP-main-main}, so the same conclusion applies.
\end{proof}

\subsection{Chaotic IAMP}
\label{subsec:p=1-IAMP}

In this subsection we give the \hyperlink{proof:t.IAMP.main.chaotic}{proof of Theorem~\ref{thm:IAMP-main}\ref{i:IAMP-main-chaotic}}.
We take as given $q_0 \in [0,1)$ and $(b,\sigma,w,p,\zeta,\zeta^\Ising)$ satisfying assumptions (\ref{i:IAMP-main-coefs})--(\ref{i:IAMP-main-init}) from Theorem~\ref{thm:IAMP-main}, and will exhibit a ``chaotic'' IAMP algorithm achieving the guarantees of Theorem~\ref{thm:IAMP-main}\ref{i:IAMP-main-chaotic}.

To implement the chaotic initialization, we again augment the input of the algorithm (recall \eqref{eq:algorithms-as-maps}) with auxiliary gaussian variables. This time the auxiliary input will consist of the Brownian bridge $\bG_\circ$, along with another $M\times N$ matrix $\bG^{\aux}$ with i.i.d.\ standard gaussian entries. Meanwhile, the ``external noise'' gaussians $\bar\bg^k$, which appeared in \eqref{eq:u-amp-def} in the $\IAMP_N$ iteration, will be artificially replaced below, and thus do not appear in this subsection. Thus the chaotic IAMP algorithm will be a Lipschitz function of inputs 
  \beq\label{e:chaotic.IAMP.inputs}
  \bz\equiv \Big(\bG,\bg^{\aux} \equiv (\bG_\circ, \bG^\aux)
    \Big)\,,\eeq
where the dimension of the auxiliary input can be understood as
$N^{\aux}= MN(\ell_{\max}+1) + MN$.

Given inputs $\bz\equiv(\bG,\bg^\aux)$ and $\tilde{\bz}\equiv(\tilde{\bG},\tilde{\bg}^\aux)$ which are $\upp$-correlated, where $\upp \le1-\epsilon$, we want to arrange that the outputs $\cA(\bz)$ and $\cA(\tilde{\bz})$ are nearly uncorrelated. The algorithm has two phases:
\begin{itemize}
\item \textbf{Phase I.} First use the auxiliary matrices $\bG^\aux$ and $\tilde{\bG}^\aux$ to generate a large amount of decorrelated noise --- the vectors $(\hat{\bg},\hat{\bs})$ and $(\tilde{\bg},\tilde{\bs})$ below, which are only $o_\eta(1)$-correlated with one another.
\item \textbf{Phase II.} Then run a modified IAMP algorithm: in the case of $\cA(\bz)$,  the matrices appearing in the iteration are based on $(\bG,\bG_\circ)$, but the noise vectors $(\hat{\bg},\hat{\bs})$ also appear in the iteration. In the case of $\cA(\tilde{\bz})$,  the matrices appearing in the iteration are based on $(\tilde{\bG},\tilde{\bG}_\circ)$, but the noise vectors $(\tilde{\bg},\tilde{\bs})$ also appear in the iteration.
\end{itemize}
The idea is that the decorrelated noise will be enough to decorrelate the outputs $\cA(\bz)$ and $\cA(\tilde{\bz})$, even though $\bz$ and $\tilde{\bz}$ are themselves quite correlated.

\begin{dfn}[chaotic IAMP, phase I]
\label{d:chaotic.IAMP.phase.I}
Similarly to \eqref{eq:A-ell}, define $\bA^{\aux}=\bG^{\aux}/N^{1/2}$.
Take a large integer $\linit\geq 10(\ubl)$ with 
\beq
\label{eq:eta-delta-eps}
0  \ll \frac{1}{\linit} \ll 
1-\upp\ll \eta\ll\delta\ll 
  \frac{1}{L}\ll \epsilon \ll \iota' \ll\iota \ll 1\,.
\eeq
(Again, this should not be confused with Assumption~\ref{a:params}, which only applies to Sections~\ref{s:rerand} and \ref{s:sde}.) Consider the IAMP iterates $(\hat{\by}^{\ell},\hat{\bv}^{\ell})_{-\linit\leq \ell\leq -1}$ defined by
\beq\label{eq:chaotic-init-setup}
\begin{aligned}
  \hat{\by}^{\ell+1}&=
  \alpha^{-1/2}
  (\bA^{\aux})^{\top}
    \hat{\bv}^{\ell}
  -\alpha^{1/2}
    \hat{\by}^{\ell}\in\bbR^N,
  \\
  \hat{\bv}^{\ell}
  &=\bA^{\aux}\hat{\by}^{\ell}
  -\alpha^{-1/2}\hat{\bv}^{\ell-1}
  \in\bbR^M\,,
\end{aligned}
\eeq
started from
$\hat{\by}^{-\linit}$ a standard gaussian vector in $\R^N$,
and $\hat{\bv}^{-\linit}=\bzero\in\R^M$. As mentioned above, we also consider another input $\tilde{\bz}$
which is $\upp$-correlated with $\bz$ for some constant $\upp\in (0,1)$. Generate iterates $(\tilde{\by}^{\ell},\tilde{\bv}^{\ell})$ for
$-\linit\leq \ell\leq -1$ using the matrix $\tilde{\bG}^{\aux}$,
started from
  \[\begin{aligned}
  \tilde{\by}^{-\linit} &=\hat{\by}^{-\linit} \in \R^N\,, \\ 
  \tilde{\bv}^{-\linit} &=\bzero\in\R^M\,.\end{aligned}
  \]
\textbf{In particular, the initialization $\tilde{\by}^{-\linit}=\hat{\by}^{-\linit}$ is treated as fixed in the algorithm, meaning it will be part of the random seed $\omega$ (see \eqref{e:iamp.chaotic.omega} below).} 
\end{dfn}

Before proceeding further, we note that the iteration
\eqref{eq:chaotic-init-setup}, considered jointly for 
$\bz$ and $\tilde{\bz}$,
 falls under the general AMP framework reviewed in \S\ref{subsec:state-evolution-multiple}. Indeed, we can decompose
\beq
\label{eq:AMP-for-correlated-disorder}
(\bG,\tilde{\bG})^{\aux} = 
\Big(\upp^{1/2}\bZ
+(1-\upp)^{1/2}\dot{\bZ},
\upp^{1/2}\bZ
+(1-\upp)^{1/2}\ddot{\bZ}
\Big)
\eeq
where $\bZ,\dot{\bZ},\ddot{\bZ}$ are i.i.d.\ $M\times N$ gaussian matrices.
Therefore we can model \eqref{eq:chaotic-init-setup} using a multi-graph with two vertices $v,y$ and three bidirectional edges. The first bidirectional edge consists of a directed edge from $v$ to $y$ labelled $\bZ^\top$, as well as a directed edge from $y$ to $v$ labelled $\bZ$. The second bidirectional edge is labelled similarly with $\dot{\bZ}$ and its transpose, and the third bidirectional edge is labelled similarly with $\ddot{\bZ}$ and its transpose. By using \eqref{eq:AMP-for-correlated-disorder} to reparametrize  in terms of $(\bZ,\dot{\bZ},\ddot{\bZ})$, it is straightforward to verify for instance that the last term on the right-hand side of each line of \eqref{eq:chaotic-init-setup} is the correct Onsager term. Similarly, after the reparametrization \eqref{eq:chaotic-init-setup}, we can apply the state evolution formulas from \S\ref{subsec:state-evolution-multiple}. This implies for instance
\beq
\label{eq:example-AMP-for-correlated-disorder}
\plim_{N\to\infty}
\frac{(\hat{\bv}^{\ell+1},\tilde{\bv}^{j+1})}{M}
=
\upp\cdot \plim_{N\to\infty}
\frac{(\hat{\by}^{\ell},\tilde{\by}^{j})}{N}\,.
\eeq
A similar calculation appears in Lemma~\ref{lem:decorrelate} below.

\begin{lem}
\label{lem:independent-init}
For any fixed $\linit$, 
as $N\to\infty$ we have the state evolution limit 
\[
(\hat{\by}^{-\linit},\ldots,
  \hat{\by}^{-1})
\to 
(Y^{-\linit},\dots,Y^{-1})
\sim \cN(0,I_{\linit})
\]
to an i.i.d. standard gaussian sequence in the manner of Proposition~\ref{prop:graph-AMP-cor}.
\begin{proof}
This is a straightforward consequence 
of Proposition~\ref{prop:graph-AMP-cor}; note that \eqref{eq:chaotic-init-setup} includes the correct Onsager terms for the AMP iteration.  
\end{proof}
\end{lem}

\begin{lem}\label{lem:decorrelate}
Given any $\upp,\eta\in (0,1)$, if $\linit$ is sufficiently large
depending on $\upp,\eta$, then
  \[
        \sup_{1\leq \ell\leq 3\linit/4}
        \lt\{
        \plim_{N\to\infty} 
  \bigg\{
  \frac{|(
  \hat{\by}^{-\ell},
  \tilde \by^{-\ell})|}{N}
  +\frac{| (
  \hat{\bv}^{-\ell},
  \tilde\bv^{-\ell})|}{M}\bigg\}
        \rt\}
  \leq \eta\,.
  \]

\begin{proof} 
From \eqref{eq:example-AMP-for-correlated-disorder} and the surrounding discussion, we have
        \begin{align}
        \label{eq:correlation-decay-for-initialization} \nonumber
    & \frac{(\hat{\by}^{-1},\tilde{\by}^{-1})}{N}
    \simeq
    \upp
    \cdot
    \frac{(\hat{\bv}^{-2},\tilde{\bv}^{-2})}{M}
    \simeq
    \upp^2
    \cdot  \frac{(
    \hat{\by}^{-2},\tilde\by^{-2})}{N}  \\
    &\qquad\simeq 
    \dots
    \simeq
    \upp^{2(\linit-1)}
    \frac{(\hat{\by}^{-\linit},
    \tilde\by^{-\linit})}{N}
  \simeq
  \upp^{2(\linit-1)} \,.
  \end{align}
The claim follows.
\end{proof}
\end{lem}

For what follows, define
 $\sign(x)\equiv\sign_0(x)\equiv \ind\{x\ge0\}-\ind\{x<0\}$. For $\eta>0$ define
  \[\sign_\eta(x)
  \equiv
  \begin{cases}
  -1 & \textup{for $x \le -\eta$,}\\
  x/\eta & \textup{for 
    $-\eta\le x\le \eta$,}\\
  1 &\textup{for $x \ge \eta$,}
  \end{cases}
  \]
and note that this function is $(1/\eta)$-Lipschitz.

\begin{dfn}[chaotic IAMP, phase II]\label{d:chaotic.IAMP}
Recall that phase I of $\cA(\bz)$ (Definition~\ref{d:chaotic.IAMP.phase.I})
has generated 
$(\hat{\by},\hat{\bv})$
for $-\linit\le\ell\le-1$. Apply $\sign_\eta$ entrywise 
to define the ``random sign vectors''
	\begin{align*}
  \hat{\bs}^\eta
  &\equiv\sign_\eta
  (\hat{\bv}^{-\linit/2}) \in\R^M\,,\\
  \hat{\br}^\eta
  &\equiv\sign_\eta
  (\hat{\by}^{-\linit/2-1}) \in\R^N\,,
  \end{align*}
which are well defined for all $\eta\ge 0$. Likewise, phase I of of $\cA(\tilde{\bz})$ 
has generated
$(\tilde{\by},\tilde{\bv})^\ell$
for $-\linit\le\ell\le-1$, and we use these to define 
  \begin{align*}
  \tilde{\bs}^\eta
  &\equiv\sign_\eta
  (\tilde{\bv}^{-\linit/2}) \in\R^M\,,\\
  \tilde{\br}^\eta
  &\equiv\sign_\eta
  (\tilde{\by}^{-\linit/2-1}) \in\R^N\,.
  \end{align*}
Phase II of $\cA(\bz)$ will  define iterates
$(\hat{\bx},\hat{\prxbx},\hat{\by},\hat{\prxby})^\ell$, while phase II of $\cA(\tilde{\bz})$ will define iterates
$(\tilde{\bx},\tilde{\prxbx},
\tilde{\by},\tilde{\prxby})^\ell$.
Recall that we abbreviate $b^\ell\equiv b_{q_0+\ell\delta}$, etc. Then, in place of 
\eqref{eq:w-sigma-b-setup}, we let
  \beq\label{eq:sign-nonlinearities}
  \begin{aligned}
  \hat{\bw}^{\chaos,\ell}
  \equiv
  \hat{\bw}^{\chaos,\ell,\eta}
  &\equiv w^\ell
    (\hat{\br}^\eta \odot \hat{\prxby}^\ell) \in\R^N\,, \\
  \hat{\bsig}^{\chaos,\ell}
  \equiv
  \hat{\bsig}^{\chaos,\ell,\eta}
  &\equiv \sigma^\ell
    (
    \hat{\bs}^\eta
    \odot\hat{\prxbx}^\ell
    ) \in\R^M\,,\\
  \hat{\bb}^{\chaos,\ell}
  \equiv \hat{\bb}^{\chaos,\ell,\eta}
  &\equiv
  \hat{\bs}^\eta\odot b^\ell
    (\hat{\bs}^\eta
    \odot
    \hat{\prxbx}^\ell) \in\R^M\,,
  \end{aligned}
  \eeq
and similarly 
$(\tilde{\bw},\tilde{\bsig},
\tilde{\bb})^{\chaos,\ell}$. 
\textbf{We also replace $\bar{\bg}^{\ell}$ from before with ``artificial'' external noise}
  \beq\label{eq:fake-external-noise}
  \begin{aligned}
  \hat{\bg}^{\ell}
  &=
  \hat{\by}^{-\linit/2 + 
        \ell
        } \in \R^N
  \,,\\
  \tilde{\bg}^\ell
  &=\tilde{\by}^{-\linit/2 + 
        \ell
        } \in\R^N\,.
  \end{aligned}
  \eeq
Recall that we chose $\linit \ge 10(\ubl)$, so the indices
$-\linit/2 + \ell$ will stay negative for all $0\le\ell\le\ubl$.
Then, recalling 
\eqref{eq:chaotic-init-setup}, the vectors $\hat{\bg}^{\ell}$ depend on the (fixed) initialization $\hat{\by}^{-\linit}$ and the random matrix $\bA^\aux$. Similarly, the vectors
$\tilde{\bg}^{\ell}$ depend on the 
same initialization
$\tilde{\by}^{-\linit}=\hat{\by}^{-\linit}$, and the random matrix $\tilde{\bA}^\aux$. Further, from Lemma~\ref{lem:independent-init} we have the state evolution limit
\begin{align}
  \nonumber
  (\hat{\br}^\eta, \hat{\bg}^0, \ldots, \hat{\bg}^{\ubl})
  &= (\sign_\eta(\hat{\by}^{-\linit/2-1}), \hat{\by}^{-\linit/2}, \ldots, \hat{\by}^{-\linit/2 + \ubl}) \\
  \label{e:aux-init-independent}
  &\rightarrow
  (\sign_\eta(Y^{-\linit/2-1}), Y^{-\linit/2}, \ldots, Y^{-\linit/2 + \ubl})
\end{align}
where the entries of this vector are mutually independent. Now recall the initialization of the original IAMP, below \eqref{eq:tilde-Y-approx-Y}.  We initialize phase II of the the chaotic IAMP $\cA(\bz)$ as follows: 
\begin{itemize}
\item As before, we 
let $\prxbx^0\in\R^M$ have coordinates i.i.d.\ from $\zeta$,
and set $\bx^0=\bzero\in\R^M$. We then set
\begin{equation}\label{eq:chaotic-init}
\begin{aligned}
	\hat{\prxbx}^0 
	&=\hat{\bs}^{\eta}\odot \prxbx^0
		\in\R^M\,,\\
	\hat{\bx}^0 &\equiv
	\hat{\bs}^{\eta}\odot \bx^0
	= \bzero\in\R^M\,.\\
\end{aligned}
\end{equation}

\item As before, we let $\prxby^0\in\R^N$ have coordinates i.i.d.\ from $\zeta^\Ising$.   If $\bbE[Y(q_0)^2]>0$, then we set
	\[\by^0 = \bigg(
	\frac{q_0}{\bbE[Y(q_0)^2]}
	\bigg)^{1/2} \prxby^0
	\in\R^N\,.\]
If $\bbE[Y(q_0)^2]=0$, then we let $\by^0\sim\mathcal{N}(0,q_0 I_N)$. In either case, we set
	\beq\label{eq:chaotic-init-2}
  \begin{aligned}
    \hat{\prxby}^0 &= \hat{\br}^\eta \odot \prxby^0
    \in\R^N\,, \\
    \hat{\by}^0 &= \hat{\br}^\eta \odot \by^0 \in\R^N\,, \,.
  \end{aligned}
\eeq
\end{itemize}
In the notation of \eqref{eq:IAMP-fully-general-fixed}, we take $\bx[0]$ to be the matrix with columns $(\prxbx^0,\bx^0)$, and we take $\by[0]$  to be the matrix with columns $(\prxby^0,\by^0)$. In the notation of 
Remark~\ref{rmk:omega-general-seed}, we view the algorithm $\cA(\bz)$ (phases I and II combined) as having random seed 
\beq\label{e:iamp.chaotic.omega}
\omega=\Big(
\tilde{\by}^{-\linit}
=\hat{\by}^{-\linit},
\bx[0]=(\prxbx^0,\bx^0),
\by[0]=(\prxby^0,\by^0)\Big)\,.\eeq
We note that in contrast with the random seed \eqref{e:iamp.random.seed} of the original IAMP, the chaotic IAMP treats the Brownian bridge $\bG_\circ$ as part of the input \eqref{e:chaotic.IAMP.inputs}, so it is no longer part of the random seed $\omega$.

For phase II of $\cA(\tilde{\bz})$, we use the same $\bx[0]$ and $\by[0]$ as for $\cA(\bz)$, since these are part of the random seed $\omega$. Then:
\begin{itemize}
\item Analogously to \eqref{eq:chaotic-init}, we define 
	\[
\begin{aligned}
\tilde{\prxbx}^0 
	&=\tilde{\bs}^{\eta}
		\odot \prxbx^0,\\
\tilde{\bx}^0 & \equiv
	\tilde{\bs}^{\eta}\odot \bx^0
	= \bzero.\\
\end{aligned}\,.
	\]
\item 
Likewise, analogously to \eqref{eq:chaotic-init-2} we define
	\[
  \begin{aligned}
    \tilde{\prxby}^0 &= \tilde{\br}^\eta \odot \prxby^0
    \in\R^N\,, \\
    \tilde{\by}^0 &= \tilde{\br}^\eta \odot \by^0 \in\R^N\,.
  \end{aligned}\]
\end{itemize}
This finishes our discussion of the initialization for phase II of $\cA(\bz)$ and $\cA(\tilde{\bz})$.

We then run the Ising IAMP iteration \eqref{eq:proxy-y-amp-def}--\eqref{eq:x-amp-def} with the modified nonlinearities from \eqref{eq:sign-nonlinearities}, and with $\hat{\bu}^0\equiv\bzero\in\R^N$, $\hat{\bm}^{-1} = \hat{\bm}^0 = \bzero \in \R^M$:
\begin{align}
\notag
  \hat{\prxby}^{\ell+1}
  -\hat{\prxby}^{\ell}
  &=
  \hat{\prxgamma}(\ell,\delta)^{-1}
  \hat{\bw}^{\chaos,\ell}
  \odot
  (\hat{\bu}^{\ell+1}
  -\hat{\bu}^\ell)
  \in\bbR^N\,,\\
\notag
  \hat{\by}^{\ell+1}
  -\hat{\by}^\ell
  &=\hat{\gamma}(\ell,\delta)^{-1}
  \hat{\bw}^{\chaos,\ell}
  \odot
  (\hat{\bu}^{\ell+1}
  -\hat{\bu}^{\ell})
  \in\bbR^N\,,\\
\label{eq:hat-u-IAMP}
  \hat{\bu}^{\ell+1}
  -\hat{\bu}^{\ell}
  &\equiv (\bA_\ell)^{\top}
  \big(
    \hat{\bm}^\ell
    -\hat{\bm}^{\ell-1}
  \big)
  + \delta^{1/2}(\bJ_{\ell+1})^\top
  \hat{\bb}^{\chaos,\ell}
  + c_\ell(\epsilon)\hat{\bg}^\ell
  -\ons_{\hat\bu,\ell}
  \in\bbR^N\,,\\
\label{eq:hat-v-def-IAMP}
  \hat{\bv}^\ell
  &=
  \bA_{\ell}\hat{\by}^\ell
  -\ons_{\hat\bv,\ell-1}
  \in\bbR^M
  \,,\\
\notag
  \hat{\bm}^{\ell+1}
  -\hat{\bm}^\ell
  &=p(q_0 + \ell\delta)^{-1}
  (\hat{\bsig}^{\chaos,\ell}
  -1)
  \odot
  (\hat{\bv}^{\ell+1}
  -\hat{\bv}^\ell)
  \in \bbR^M\,,
  \\
\label{eq:hat-bx-IAMP}
  \hat{\prxbx}^{\ell+1}
  -\hat{\prxbx}^{\ell}
  &=\hat{\bx}^{\ell+1}
  -\hat{\bx}^{\ell}
  =\hat{\bsig}^{\chaos,\ell}\odot
  \big(\hat{\bv}^{\ell+1}
    -\hat{\bv}^\ell\big)
  +
  \delta
  p'(q_0 + \ell\delta)^{1/2}
  \hat{\bb}^{\chaos,\ell}
  \in 
  \bbR^M\,,
\end{align}
where the constants $\gamma$ and $\prxgamma$ will be defined in \eqref{e:gamma.chaotic} below.
We analogously define $(\tilde{\bx},\tilde{\prxbx},\tilde{\by},\tilde{\prxby},\tilde{\bu},\tilde{\bv},\tilde{\bm})^\ell$, using the $\upp$-correlated matrices $\tilde{\bA}_\ell$, $\tilde{\bJ}_\ell$ derived from $\tilde{\bz}$. We call the above \textbf{chaotic $\eta$-symmetrized IAMP}, denoted $\IAMP_{\chaos,N}(p,\sigma,b,w,\delta,\eta,\zeta,\zeta^{\Ising})$.   
 (Note the algorithm does not depend on the scalar $\upp\in (0,1)$, which enters only to analyze $\chi$.) We emphasize that the above iterates depend on $\eta$ through \eqref{eq:sign-nonlinearities}, although we often suppress this from the notation. In what follows, we sometimes make the $\eta$-dependence explicit by writing $\hat{\bx}^\ell=\hat{\bx}^{\ell,\eta}$, etc.
\end{dfn}

Generalizing the discussion around \eqref{eq:AMP-for-correlated-disorder},
if we consider the combined iteration of phases I and II (as specified by Definitions~\ref{d:chaotic.IAMP.phase.I} and \ref{d:chaotic.IAMP}),
jointly for both $\bz$ and $\tilde{\bz}$, this is again encapsulated by the general framework from \S\ref{subsec:state-evolution-multiple}, which also specifies the Onsager terms. Indeed, similarly to \eqref{eq:AMP-for-correlated-disorder}, we can decompose
  \[
  (\bJ_\ell,\tilde{\bJ}_\ell) = 
\Big(\upp^{1/2}\bZ_\ell
+(1-\upp)^{1/2}\dot{\bZ}_\ell,
\upp^{1/2}\bZ_\ell
+(1-\upp)^{1/2}\ddot{\bZ}_\ell
\Big)\,,\]
and introduce additional bidirectional edges labelled with $\bZ_\ell,\dot{\bZ}_\ell,\ddot{\bZ}_\ell$ and their transposes. This allows us to handle the correlation between $\bz$ and $\tilde{\bz}$. In order to rewrite \eqref{eq:hat-u-IAMP} into the general form \eqref{eq:IAMP-fully-general-fixed}, the first two terms on the right-hand side of \eqref{eq:hat-u-IAMP} can be handled as before, but the third term $c_\ell(\epsilon)\hat{\bg}^\ell$ needs to be handled slightly differently: for this, we need to recall from \eqref{eq:fake-external-noise} and \eqref{eq:chaotic-init-setup} that
  \[
  \hat{\bg}^\ell
  = \hat{\by}^{-\linit/2+\ell}
  = \alpha^{-1/2}(\bA^\aux)^\top \hat{\bv}^{-\linit/2+\ell-1}
    -\alpha^{1/2} \hat{\by}^{-\linit/2+\ell-1}\,,
  \]
where the last term is the Onsager term. Thus, to formally fit \eqref{eq:hat-u-IAMP} into the form \eqref{eq:IAMP-fully-general-fixed}, we can first re-index \eqref{eq:IAMP-fully-general-fixed} to include an additional index $k=-1$ corresponding to $\bJ_{-1}\equiv\bA^\aux$, and then define $\bu^{\ell+1,-1}$ to capture the contribution to $\bu^{\ell+1}$ from the terms $c_j(\epsilon)\hat{\bg}^j$, for $0\le j\le\ell$: 
  \beq\label{e:chaotic.iamp.hatg.standard.form} 
  \bu^{\ell+1,-1}
  \equiv \bJ_{-1}^\top \bigg\{ 
  \frac{1}{\alpha^{1/2}}
  \sum_{j=0}^\ell c_j(\epsilon)
  \hat{\bv}^{-\linit/2+j-1} \bigg\}
  - \ons_{\bu,\ell,-1}
  \eeq
Note moreover that since $\hat{\bg}^\ell$ already incorporates its own Onsager term, this Onsager term does not appear in the $\ons_{\hat\bv,\ell-1}$ term in \eqref{eq:hat-u-IAMP}. As a result, the definition of $\ons_{\hat\bu,\ell}$ looks essentially similar to the definition of $\ons_{\bu,\ell}$ from \eqref{eq:IAMP-fully-general-fixed.Ons} and
\eqref{e:u.combined.ons} --- one should not include any extra $k=-1$ term corresponding to $\bA^\aux$. Similarly, the remaining parts of the chaotic IAMP are straightforward to rewrite into the general form \eqref{eq:IAMP-fully-general-fixed}, and all the corresponding Onsager terms are thereby determined. For the $k\ge0$ terms of \eqref{eq:IAMP-fully-general-fixed},
the functions $h$ and $f$ are given similarly to \eqref{e:h.l.k.explicit} and \eqref{e:f.l.k.explicit}:
  \beq\label{e:f.h.l.k.explicit.CHAOTIC}
  \begin{aligned}
  h_{\chaos,\ell,k}(\hat{U}[[\ell],Y[0] )
  &= a_k \hat{Y}^\ell\,,\\
  \hat{S}^\eta
  f_{\chaos,\ell,k}(
    \hat{V}[[\ell],
    X[0] )
  &=\ind\{k\le\ell\} a_k
    \hat{S}^\eta (\hat{M}^\ell
    -\hat{M}^{k-1})
  +
  \ind\{k\ge1\}
  \delta^{1/2} 
  (\hat{S}^\eta)^2
  b^{k-1}(\hat{S}^\eta
    \hat{\prxX}^{k-1})\,.
  \end{aligned}
  \eeq
As before, $\hat{U}^\ell$ is the sum of $\hat{U}^{\ell,k}$ over $0\le k\le\ell$, and $\hat{V}^\ell$ is the sum of $\hat{V}^{\ell,k}$ over $0\le k\le\ell$.

We denote the corresponding state evolution limits $\hat{X}^\ell\equiv \hat{X}^{\ell,\eta}$, and similarly $(\hat{\prxX},\hat{Y},\hat{\prxY},\hat{U},\hat{G},\hat{V},\hat{M})^\ell$, along with the random variables 
	\begin{equation}
	\label{e:hat.S.hat.R.eta}
	\begin{aligned}
	\hat{S}^\eta
	&= \sign_\eta(\hat{V}^{-\linit/2})\,, \\
	\hat{R}^\eta 
	&= \sign_\eta(\hat{Y}^{-\linit/2-1})
	\end{aligned}
	\end{equation}
representing the state evolution limits of $\hat{\bs}^\eta$ and $\hat{\br}^\eta$. 
(We again point out that here $\hat G^{\ell}$ denotes the state evolution limit of $\hat\bg^{\ell}$ and is a one-dimensional standard gaussian random variable.) Recalling \eqref{eq:sign-nonlinearities}, we abbreviate
\begin{equation}
\label{e:chaotic.w.sigma.b.state.evol}
\begin{aligned}
\hat{w}^{\chaos,\ell,\eta}
	&\equiv
	w^\ell(\hat{R}^\eta \prxY^\ell)\,,\\
\hat{\sigma}^{\chaos,\ell,\eta}
	&\equiv
	\sigma^\ell(\hat{S}^\eta
		\prxX^\ell)\,,\\
\hat{b}^{\chaos,\ell,\eta}
	&\equiv \hat{S}^\eta b^\ell(\hat{S}^\eta \prxX^\ell)\,.
\end{aligned}\end{equation}
Likewise we have the state evolution limits $(\tilde{X},\tilde{\prxX},\tilde{Y},\tilde{\prxY},\tilde{U},\tilde{G},\tilde{V},\tilde{M})^{\ell,\eta}$, 
$(\tilde{S},\tilde{R})^\eta$, and
$(\tilde{w},\tilde{\sigma}, \tilde{b})^{\chaos,\ell,\eta}$.  \textbf{(See below, around Lemma~\ref{lem:chaotic-symmetrization}, for comments about the case $\eta=0$.)} We note that \eqref{e:aux-init-independent}  
implies the state evolution limits $(\hat{R}^\eta,\hat{G}^1,\dots,\hat{G}^{\ell})$ are jointly independent, as are $(\tilde{R}^\eta,\tilde {G}^1,\dots,\tilde {G}^{\ell})$, even for positive $\eta$. Similarly $\hat{G}^i$ and $\tilde{G}^j$ are independent, unless $i=j$ in which case the correlation is $o_\eta(1)$ by Lemma~\ref{lem:decorrelate}. Thus, the $(\hat{U}^{\ell,k},\tilde{U}^{\ell,k})$ are independent of the $(\hat{V}^{\ell,k},\tilde{V}^{\ell,k})$, and the covariances among the $U$'s and $V$'s are described by a generalization of \eqref{eq:state-evolution-setup} --- for instance, the correct generalization of \eqref{e:state.evol.var.increment.u} is given by (recall $p^\ell = p(q_0 + \delta \ell)$; this should not be confused with
the correlation $\upp$ of $\bz$ and $\tilde{\bz}$)
\beq\label{eq:chaotic-U-covariance-recursion}
\begin{aligned}
\bbE\Big[(\hat U^{\ell+1}-\hat U^{\ell})(\tilde U^{\ell+1}-\tilde U^{\ell})\Big]
&=
\upp  \cdot \alpha p^\ell 
\bbE\Big[(\hat M^{\ell}-\hat M^{\ell-1})(\tilde M^{\ell}-\tilde M^{\ell-1})
  \Big] \\
&\qquad +
\upp \cdot \alpha \delta  \bbE[\hat{b}^{\chaos,\ell} \tilde{b}^{\chaos,\ell}]
+
\upp^{\linit+2\ell} c_{\ell}(\epsilon)^2\,.
\end{aligned}
\eeq
Here the last coefficient $\upp^{\linit+2\ell}$ follows from the calculation \eqref{eq:correlation-decay-for-initialization} from Lemma~\ref{lem:decorrelate}, and is $o_\eta(1)$.

Recalling the initialization  from Definition~\ref{d:chaotic.IAMP}, for $\cA(\bz)$ we set
\begin{itemize}
\item $\prxX^0\sim\zeta$, $X^0=0$, and 
	\[(\hat{\prxX}^0,\hat{X}^0)
	=\hat{S}^\eta(\prxX^0,X^0)
	=\hat{S}^\eta(\prxX^0,0)\,.
	\]
\item $\prxY^0\sim\zeta^\Ising$,
$Y^0 = (q_0)^{1/2} \prxY^0 / \E[Y(q_0)^2]^{1/2}$ if $\E[Y(q_0)^2]>0$,
$Y^0 \sim\cN(0,q_0)$
 if $\E[Y(q_0)^2]=0$, and
 	\[(\hat{\prxY}^0,\hat{Y}^0)
	=\hat{R}^\eta(\prxY^0,Y^0)\,.
	\]
\end{itemize}
We denote $X[0]=(\prxX^0,X^0)$ and $Y[0]=(\prxY^0,Y^0)$, and we emphasize again that this is used for both $\cA(\bz)$ and $\cA(\tilde{\bz})$, since it is part of the random seed \eqref{e:iamp.chaotic.omega}. For $\cA(\tilde{\bz})$, we define analogously
\begin{itemize}
\item $(\tilde{\prxX}^0,\tilde{X}^0)
	=\tilde{S}^\eta(\prxX^0,X^0)
	=\tilde{S}^\eta(\prxX^0,0)$, and
	
\item $(\tilde{\prxY}^0,\tilde{Y}^0)
	=\tilde{R}^\eta(\prxY^0,Y^0)$.
\end{itemize} 
Then, similarly to \eqref{e:state.evol.x.y.increments}, we define recursively
\beq\label{e:state.evol.x.y.increments.CHAOTIC}
\begin{aligned}
\hat{\prxY}^{\ell+1}-\hat{\prxY}^{\ell}
&= \hat{\prxgamma}(\ell,\delta,\eta)^{-1}
w^\ell(\hat{R}^\eta\hat{\prxY}^\ell)
(\hat{U}^{\ell+1}-\hat{U}^\ell)\,,\\
\hat{Y}^{\ell+1}-\hat{Y}^\ell
&= \hat{\gamma}(\ell,\delta,\eta)^{-1}
w^\ell(\hat{R}^\eta\hat{\prxY}^\ell)
(\hat{U}^{\ell+1}-\hat{U}^\ell)\,,
\\
\hat{M}^{\ell+1}-\hat{M}^\ell
&=p(q_0+\ell\delta)^{-1}
(\sigma^\ell(\hat{S}^\eta \hat{\prxX}^\ell)-1)
\cdot 
(\hat{V}^{\ell+1}-\hat{V}^{\ell})\,,
\\
\hat{\prxX}^{\ell+1}-\hat{\prxX}^\ell
=
\hat{X}^{\ell+1}-\hat{X}^\ell
&=\sigma^\ell(\hat{S}^\eta
  \hat{\prxX}^\ell)
\cdot 
(\hat{V}^{\ell+1}-\hat{V}^\ell)
+\delta p'(q_0+\ell\delta)^{1/2} 
\hat{S}^\eta
b^\ell(\hat{S}^\eta
  \hat{\prxX}^\ell)\,.
\end{aligned}
\eeq
Similarly to \eqref{eq:proxy-C-ell-delta-def}, we define
  \beq\label{e:gamma.chaotic}
  \begin{aligned}
  \hat{\prxgamma}^\ell
  \equiv \hat{\prxgamma}
    (\ell,\delta, \eta)
  &\equiv
  \delta^{-1/2} 
  \bbE\big[(\hat{U}^{\ell+1}
    - \hat{U}^\ell)^2\big]^{1/2}\,,
  \\
  \hat{\gamma}^\ell\equiv
  \hat{\gamma}(\ell,\delta,\eta)
  &\equiv 
  \delta^{-1/2} 
  \bbE\big[(\hat{U}^{\ell+1}
    - \hat{U}^\ell)^2\big]^{1/2}
    \bbE[w^\ell\big(\hat{R}^\eta\hat{\prxY}^\ell\big)^2]^{1/2}\,.
  \end{aligned}
  \eeq
When $\eta=0$, similarly to \eqref{e:hat.S.hat.R.eta}, we take
	\[
	\begin{aligned}
	\hat{S}
	&= \sign(\hat{V}^{-\linit/2})\,, \\
	\hat{R}^\eta 
	&= \sign(\hat{Y}^{-\linit/2-1})\,,
	\end{aligned}
	\]
and similarly $\tilde{S},\tilde{R}$. We then note that $(\hat S,\tilde S, \hat R, \tilde R)$ 
are uniformly random signs, independent of one another and of all the
 $\hat G^{\ell},\tilde{G}^\ell$.
Then the recursion \eqref{e:state.evol.x.y.increments.CHAOTIC} is well-defined for any $\eta\ge0$. 
However, for $\eta=0$, it is not immediately clear whether the random variables produced by \eqref{e:state.evol.x.y.increments.CHAOTIC}
 actually represent the state evolution limit, e.g., whether $\hat{X}^{\ell,\eta=0}$ approximates $\bx^{\ell,\eta=0}$ in the sense of \eqref{eq:state-evolution-limit-graph}. Indeed, the result of Proposition~\ref{prop:graph-AMP-cor} requires all non-linearities to be Lipschitz --- clearly, this holds for $\sign_{\eta}$ with $\eta>0$, but not for $\sign=\sign_0$ itself. For the same reason, taking $\eta=0$  in the chaotic IAMP does not give a Lipschitz algorithm. On the other hand, formally setting $\eta=0$ in the state evolution limits gives the following symmetrization property, which is one of the main purposes of the above construction:

\begin{lem}
\label{lem:chaotic-symmetrization}
For the collection of random variables described by \eqref{e:state.evol.x.y.increments.CHAOTIC} with $\eta=0$, let us abbreviate $\hat{S}\equiv \hat{S}^{\eta=0}$, $\hat{R}\equiv \hat{R}^{\eta=0}$, 
$\hat{X}^\ell\equiv \hat{X}^{\ell,\eta=0}$, and so on. On the other hand, let $(X^{\ell},M^{\ell},V^{\ell},\dots)_{\ell\leq \ubl}$ be the state evolution limit of our earlier iteration $\IAMP(p,\sigma,b,w,\delta)$, as described around \eqref{e:state.evol.x.y.increments}. Then
  \[\begin{aligned}
  (\hat{S},\hat{X}^\ell,
    \hat{\prxX}^\ell,
    \hat{M}^\ell,\hat{V}^\ell)_{0\leq \ell\leq \ubl}
  &\stackrel{\textit{d}}{=}
  (S', S' X^{\ell}, S' \prxX^\ell,
    S' M^{\ell}, S' V^{\ell})_{0\leq \ell\leq \ubl}\,,\\
  (\hat{R},\hat Y^\ell,\hat{\prxY}^\ell,
  \hat{U}^\ell,\hat{G}^{\ell})_{0\leq\ell\leq\ubl}
  &\stackrel{\textit{d}}{=}(R',R' Y^{\ell},R' \prxY^{\ell}, R' U^{\ell}, R' \bar{G}^{\ell})_{0\leq \ell\leq \ubl}. 
  \end{aligned}\]
where $S'$ and $R'$ are 
symmetric random signs independent of everything else. 
In particular, at $\eta=0$ the constants $\hat{\prxgamma}$, $\hat{\gamma}$ from \eqref{e:gamma.chaotic} are the same as the constants $\prxgamma$, $\gamma$ from \eqref{eq:proxy-C-ell-delta-def}. Additionally $\hat G^{\ell}$ is independent of $(\hat Y^k,\hat\prxY^k,\hat U^k)_{0\leq k\leq \ell-1}$. 
\end{lem}

\begin{proof}
Lemma~\ref{lem:independent-init} shows the equivalence in distribution of $\hat G^{\ell}$ and $\bar{G}^{\ell}$, the state evolution limits of $\hat\bg^{\ell}$ and $\bar\bg^{\ell}$, where $\bar{\bg}^\ell$ is the external gaussian noise from the original IAMP iteration \eqref{eq:u-amp-def}.
The dynamics
\eqref{e:state.evol.x.y.increments.CHAOTIC}
are exactly the same as 
\eqref{e:state.evol.x.y.increments}, but for the symmetrized random variables 
\[\hat S\hat X^{\ell},
\hat S\hat{\prxX}^\ell, \hat S\hat{M}^\ell, \hat S\hat{V}^\ell, \hat{R}\hat{Y}^\ell, \hat{R}\hat{\prxY}^\ell, \hat{R}\hat{U}^\ell, \hat{R}\hat{G}^\ell
\,.\]
(This is the reason for the initializations \eqref{eq:chaotic-init}--\eqref{eq:chaotic-init-2}.) Moreover, since $(\hat{V}^\ell)_\ell$ is a centered gaussian process independent of $\hat S$, it is equidistributed as $(\hat{S}\hat{V}^\ell)_\ell$. Similarly $(\hat{G}^\ell)_\ell$ is equidistributed as $(\hat{R}\hat{G}^\ell)_\ell$. The claim then follows by induction on $\ell$, using
the modified non-linearities \eqref{eq:sign-nonlinearities}.
\end{proof}

We next show stochastic continuity in $\eta$ of state evolution limits down to $\eta=0$.
This will let us apply Lemma~\ref{lem:chaotic-symmetrization} to analyze the behavior at small positive $\eta$.
(A similar approximation argument shows that state evolution is valid even for $\eta=0$, but will not be necessary.)

\begin{lem}
\label{lem:eta-to-0-chaotic-iamp}
For parameters as in \eqref{eq:eta-delta-eps} and $0\leq \ell\leq \ubl$,
  \beq
  \label{eq:eta-to-0-chaotic-iamp}
  \lim_{\eta\downarrow 0}
  \bbE\Big[
  \big(\hat{X}^{\ell,\eta}
  -\hat{X}^{\ell,0}\big)^2
  \Big]
  =0.\eeq
Similarly for $\hat{M}^{\ell},\hat{V}^{\ell},\hat{Y}^{\ell},\hat{U}^{\ell}$.
\end{lem}

\begin{proof}
Since the state evolution limits $Y^{-\linit/2- 1}$ and $V^{-\linit/2}$ are standard gaussian, we clearly have 
  \[
        \lim_{\eta\downarrow 0}
        \bbE\Big[\big(
        \hat{S}^\eta
        -\hat{S}^0\big)^2\Big]
        = \lim_{\eta\downarrow 0}
        \bbE\Big[\big(
        \hat{R}^\eta
        -\hat{r}^0\big)^2\Big]
        =0\,.\]
Then \eqref{eq:eta-to-0-chaotic-iamp} follows easily by inductively estimating the errors incurred by the state evolution recursion \eqref{e:state.evol.x.y.increments.CHAOTIC} across iterations $0\leq\ell\leq\ubl$.
\end{proof}

\begin{lem}
\label{lem:chaotic-is-ok}
  For parameters as in \eqref{eq:eta-delta-eps}, the state evolution limits of $\IAMP_{\chaos,N}(p,\sigma,b,w,\delta,\eta
  )$ satisfy:
  \begin{equation}
  \label{eq:chaotic-iamp-still-works}
  \begin{aligned}
  \bbW_2\Big( |X^{\ubl}|, |\hat X^{\ubl,\eta}|\Big)
  &\le
  o_{\eta\to 0}(1), \\
  \quad 
  \bbW_2\Big(|Y^{\ubl}|,
    |\hat Y^{\ubl,\eta}|\Big)
  &\le
        o_{\eta\to 0}(1)\,,
  \end{aligned}
  \end{equation}
where $X^{\ubl}$ and $Y^{\ubl}$ are as in \eqref{e:state.evol.x.y.increments}
from 
\S\ref{subsec:IAMP-setup}.
\end{lem}

\begin{proof}
Lemma~\ref{lem:chaotic-symmetrization} shows that $|X^{\ubl}|$ and $|\hat{X}^{\ubl,\eta=0}|$ have exactly the same distribution, so the first statement follows by taking $\eta$ small enough and applying Lemma~\ref{lem:eta-to-0-chaotic-iamp}. The second is shown similarly since Lemma~\ref{lem:eta-to-0-chaotic-iamp} also holds for the $Y$ processes.
\end{proof}

We now define a preliminary $C(L,\epsilon)$-Lipschitz algorithm $\breve\cA_N$ and argue it satisfies all the requirements of Theorem~\ref{thm:IAMP-main}\ref{i:IAMP-main-chaotic} except $\chi_{\cA_N}(0) = 0$ and $\chi_{\cA_N}(1-\iota) \le \iota$.
Let $\tilde\cA_N(\bz=(\bG,\bg^\aux))=\hat\by^{\ubl}$, as defined by the chaotic IAMP 
of Definitions \ref{d:chaotic.IAMP.phase.I} and \ref{d:chaotic.IAMP}. Let $\bar\cA_N$ be the Kirszbraun extension of $\tilde\cA_N$ constructed in the \hyperlink{proof:t.IAMP.main.main}{proof of Theorem~\ref{thm:IAMP-main}\ref{i:IAMP-main-main}}.
Then $\bar\cA_N$ is globally $C(L,\epsilon)$-Lipschitz and $\tilde\cA_N(\bz) = \bar\cA_N(\bz)$ with probability $1-e^{-cN}$.
Finally we set
\beq\label{e:breve-cA-N}
  \breve{\cA}_N(\bz) 
  = \frac{N^{1/2}}{\bbE_{\bz'}[\|\bar\cA_N(\bz')\|^2]^{1/2}} \bar\cA_N(\bz) 
  = \bar{c}_N \bar\cA_N(\bz)\,,
\eeq
so that $\breve{\cA}_N$ satisfies condition~\eqref{it:expectation-condition} of Definition~\ref{d:Lip}. 

\begin{ppn}
\label{ppn:chaotic-amp-still-works}
In the setting of Theorem~\ref{thm:IAMP-main}\ref{i:IAMP-main-chaotic}, the map $\breve{\cA}_N$ defined by \eqref{e:breve-cA-N}
  is a $C(L,\epsilon)$-Lipschitz algorithm in the sense of  Definition~\ref{d:Lip}, and satisfies the estimates 
  \begin{align*}
    \bbP\lt(\breve{\cA}_N(\bz) \in \Sigma(\iota) \rt) &\ge 1-e^{-cN}\,, \\
    \bbP\Big(
      \bbW_2\lt(\mu_{\bG,\sym}(\breve{\cA}_N(\bz)), \Law(|X(1)|) \rt) \le \iota
    \Big) &\ge 1-e^{-cN}\,,  \\
    \bbW_2\Big(\mu^{\Ising}(\breve{\cA}_N), \cP(\{\pm 1\}) \Big) &\le \iota\,, \\
    \bbW_2\Big(\sym(\mu(\breve{\cA}_N)), \Law(|X(1)|) \Big) &\le \iota\,.
  \end{align*}

\begin{proof}
Recall the order of parameters \eqref{eq:eta-delta-eps}. Thanks to Lemma~\ref{lem:chaotic-symmetrization}, the preceding proofs of Propositions~\ref{prop:Ising-IAMP-analysis} and \ref{prop:IAMP-to-SDE-limit} and Lemma~\ref{l:gamma.estimates} apply unchanged on the level of state evolution limits for $\eta=0$, adding the appropriate factors of $\hat S$ and $\hat R$ as indicated by Lemma~\ref{lem:chaotic-symmetrization}. By Lemmas~\ref{lem:eta-to-0-chaotic-iamp}--\ref{lem:chaotic-is-ok}, at $\eta > 0$ the state evolution limits only differ from the $\eta = 0$ case by $o_\eta(1)$ in $\bbW_2$ distance.

We will now imitate the proof of Theorem~\ref{thm:IAMP-main}\ref{i:IAMP-main-main}; note that since we are interested in the symmetrized measures $\mu_{\bG,\sym}(\breve{\cA}_N(\bG,\bg^\aux))$ and $\sym(\mu(\breve{\cA}_N))$, the argument will be insensitive to the factors of $\hat S$.
Similarly, both the event $\breve{\cA}_N(\bz) \in \Sigma(\iota)$ and the distance $\bbW_2(\mu^{\Ising}(\breve{\cA}_N), \cP(\{\pm 1\}))$ do not depend on the signs of the entries of $\breve{\cA}_N(\bz)$, so the argument will be insensitive to the factors of $\hat R$.

By a computation identical to the one in the \hyperlink{proof:p.IAMP-to-SDE-limit.Ax.y}{proof of Proposition~\ref{prop:IAMP-to-SDE-limit.Ax.y}}, we can show
\[
  \sup \bigg\{
        \plim_{N\to\infty}
  \frac{\|\bA\hat\by^{k,\eta=0}-\hat\bx^{k,\eta=0}\|}
    {N^{1/2}}
    : 0\leq k\leq \ubl
  \bigg\}
  \leq 
  C \delta^{1/4}
  +C_0 \epsilon^{1/2}\,.
\]
In particular, as mentioned around \eqref{e:chaotic.iamp.hatg.standard.form}, the formulas for the Onsager terms $\ons_{\hat{\bu},\ell}$ and $\ons_{\hat{\bv},\ell}$ do not account for the replacement of $\bar{\bg}^\ell$ with $\hat{\bg}^\ell$, so the computations of the \hyperlink{proof:p.IAMP-to-SDE-limit.Ax.y}{proof of Proposition~\ref{prop:IAMP-to-SDE-limit.Ax.y}} go through essentially unchanged. Then, approximating the $\eta > 0$ state evolution with $\eta = 0$ yields
\[
  \sup \bigg\{
        \plim_{N\to\infty}
  \frac{\|\bA\hat\by^{k,\eta}-\hat\bx^{k,\eta}\|}
    {N^{1/2}}
    : 0\leq k\leq \ubl
  \bigg\}
  \leq 
  C \delta^{1/4}
  +C_0 \epsilon^{1/2}
  + o_\eta(1)\,.
\]
The remainder is identical to the \hyperlink{proof:t.IAMP.main.main}{proof of Theorem~\ref{thm:IAMP-main}\ref{i:IAMP-main-main}}. The final algorithm is $C(L,\epsilon)$-Lipschitz because the IAMP parameters $\delta, \eta$ are set sufficiently small as functions of $L,\epsilon$.
\end{proof}
\end{ppn}

The symmetry of the distributions of $\hat X^\ell$ and $\hat{\prxX}^\ell$ shown in Lemma~\ref{lem:chaotic-symmetrization} for $\eta=0$, together with the continuity of the state evolution as $\eta \downarrow 0$ shown in Lemma~\ref{lem:eta-to-0-chaotic-iamp}, has another useful consequence concerning the average drift:

\begin{cor}
\label{cor:average-drift-zero}
  For all $\ell\leq \ubl$, we have 
  \begin{equation}
  \label{eq:symmetrized-b-mean-zero}
  \lim_{\eta\downarrow 0}
  \bbE\big[
  \hat{b}^{\chaos,\ell}(\hat X^{\ell,\eta})
  \big]
  =
  \lim_{\eta\downarrow 0}
  \bbE\big[
  \hat{b}^{\chaos,\ell}(\hat{\prxX}^{\ell,\eta}) \big]
  =
  0\,.
  \end{equation}
\end{cor}

\begin{proof}
For $\eta=0$ it follows by the sign symmetry shown in Lemma~\ref{lem:chaotic-symmetrization} that
  \[\bbE\big[
  \hat{b}^{\chaos,\ell}(\hat X^{\ell,\eta=0})
  \big]
  =\bbE\big[
  \hat{b}^{\chaos,\ell}(\hat{\prxX}^{\ell,\eta=0}) \big]
  =
  0\,.\]
The result for $\eta\to0$ then follows by applying Lemma~\ref{lem:eta-to-0-chaotic-iamp}.
\end{proof}

We have shown in
Proposition~\ref{ppn:chaotic-amp-still-works} that $\breve{\cA}_N$ satisfies all requirements of Theorem~\ref{thm:IAMP-main}\ref{i:IAMP-main-chaotic}, except those involving the correlation function $\chi$. Writing 
	\[\breve{\chi}_N(\upp)
	\equiv\chi_{\breve{\cA}_N}(\upp)
	\equiv 
		\frac{\E(\breve{\cA}_N(z),\breve{\cA}_N(z))}{N}
	\]
for the correlation function of $\breve{\cA}_N$, we next prove (in Lemma~\ref{l:chi-approx-by-amp-correlation} and Theorem~\ref{thm:chaotic-iamp-is-chaotic} below) that the nondecreasing function $\breve{\chi}_N$ approximates zero in that $\breve{\chi}_N(1-\iota)\leq \breve{\chi}_(\upp)\leq \iota$, where we recall the order of limits \eqref{eq:eta-delta-eps}. 
We again recall that we view $\breve{\cA}_N$ as a function of inputs $\bz$ as specified by \eqref{e:chaotic.IAMP.inputs}, with random seed $\omega$ as specified by \eqref{e:iamp.chaotic.omega}.

\begin{lem}\label{l:chi-approx-by-amp-correlation}
Under the assumptions of Theorem~\ref{thm:IAMP-main}\ref{i:IAMP-main-chaotic}, in the limiting regime \eqref{eq:eta-delta-eps} we have
  \[
    \sup_{\upp \in [0,1]}
    \lt|
      \lim_{N\to\infty} \chi_{\breve{\cA}_N}(\upp) - \plim_{N\to\infty} \frac{(\hat{\by}^{\ubl},\tilde{\by}^{\ubl})}{N}
    \rt| \le \iota'\,,
  \]
  where $\hat{\by}^{\ubl},\tilde{\by}^{\ubl}$ are the outputs of $\IAMP_{\chaos,N}(p,\sigma,b,w,\delta,\eta)$ on $\upp$-correlated $\bz, \tilde{\bz}$.
\begin{proof}
  For each $\upp \in [0,1]$, we have
  \[
    \chi_{\breve{\cA}_N}(\upp)
    = \frac{\bbE (\breve{\cA}_N(\bz), \breve{\cA}_N(\tilde \bz) )}{N}
    = \frac{N}{\bbE_{\bz'}[\|\bar\cA_N(\bz')\|^2]} \cdot \frac{\bbE (\bar{\cA}_N(\bz), \bar{\cA}_N(\tilde \bz) )}{N}\,.
  \]
  As argued in the \hyperlink{proof:t.IAMP.main.main}{proof of Theorem~\ref{thm:IAMP-main}\ref{i:IAMP-main-main}},
  \[
    \bigg|\frac{N}{\bbE_{\bz'}[\|\bar\cA_N(\bz')\|^2]} - 1\bigg| \le \iota'\,.
  \]
  We can write
  \[
    \frac{(\bar{\cA}_N(\bz), \bar{\cA}_N(\tilde \bz))}{N}
    = \frac{\|\bar{\cA}_N(\bz) + \bar{\cA}_N(\tilde \bz)\|^2}{2N} 
    - \frac{\|\bar{\cA}_N(\bz)\|^2}{2N} 
    - \frac{\|\bar{\cA}_N(\tilde\bz)\|^2}{2N}
  \]
  Since $\bar{\cA}$ is a $C(L,\epsilon)$-Lipschitz function of standard gaussians, each of the terms on the right-hand side has variance $O(1/N)$.
  Thus $(\bar{\cA}_N(\bz), \bar{\cA}_N(\tilde \bz)) / N$ has variance $O(1/N)$, which implies
  \[
    \lim_{N\to\infty} \frac{\bbE (\bar{\cA}_N(\bz), \bar{\cA}_N(\tilde \bz) )}{N}
    = \plim_{N\to\infty} \frac{(\bar{\cA}_N(\bz), \bar{\cA}_N(\tilde \bz) )}{N}\,.
  \]
  Finally, because $\bar{\cA}_N(\bz)$ agrees with $\tilde{\cA}_N(\bz)$ with probability $1-e^{-cN}$,
  \[
    \plim_{N\to\infty} \frac{(\bar{\cA}_N(\bz), \bar{\cA}_N(\tilde \bz) )}{N}
    = \plim_{N\to\infty} \frac{(\tilde{\cA}_N(\bz), \tilde{\cA}_N(\tilde \bz) )}{N}
    = \plim_{N\to\infty} \frac{(\hat{\by}^{\ubl},\tilde{\by}^{\ubl})}{N}\,. 
  \]
  This concludes the proof.
\end{proof}
\end{lem}

\begin{thm}
\label{thm:chaotic-iamp-is-chaotic}
Under the assumptions of Theorem~\ref{thm:IAMP-main}\ref{i:IAMP-main-chaotic}, for fixed $L,\epsilon,\delta$, we have 
  \beq\label{eq:chaotic-iamp-is-chaotic}
     \lim_{\eta\downarrow 0}
     \lim_{\upp\uparrow 1}
     \lim_{\linit\to\infty}
    \plim_{N\to\infty} \frac{(\hat{\by}^{\ubl},\tilde{\by}^{\ubl})}{N}
     =
     0\,.\eeq
(We note that the order of parameters is consistent with \eqref{eq:eta-delta-eps}.)
\end{thm}

Towards the \hyperlink{proof:thm.chaotic-iamp-is-chaotic}{proof of Theorem~\ref{thm:chaotic-iamp-is-chaotic}}, we recall from Lemma~\ref{lem:orthogonal-increments} that
  \begin{equation}
  \label{eq:V-self-correlation}
  \begin{aligned}
  \bbE[\hat V^0 (\hat V^{i+1}-\hat V^{i})] &= 0\,, & 
  \bbE[
  (\hat V^{i+1}-\hat V^{i})
  (\hat V^{j+1}-\hat V^j)]
  &=0\,,\\
  \bbE[\tilde V^0 (\tilde V^{i+1}-\tilde V^{i})] &= 0\,, & 
  \bbE[
  (\tilde V^{i+1}-\tilde V^{i})
  (\tilde V^{j+1}-\tilde V^j)]
  &=0\,,
  \end{aligned}
  \end{equation}
for all $i\ne j$.
Lemma~\ref{lem:pair-orthogonality} 
below is a pair generalization which implies, for all $i\neq j$,
\begin{align}
\label{eq:V-correlation}
  \bbE[\hat V^0 (\tilde V^{j+1}-\tilde V^j)] &= 0\,, & 
  \bbE[(\hat V^{i+1}-\hat V^{i})
  (\tilde V^{j+1}-\tilde V^j)]
  &=0,
  \\
\label{eq:Y-correlation}
  \bbE[\hat Y^0 (\tilde Y^{j+1}-\tilde Y^j)] &= 0\,, & 
  \bbE[(\hat Y^{i+1}-\hat Y^{i})
  (\tilde Y^{j+1}-\tilde Y^j)]
  &=0.
\end{align}
We define the $\sigma$-fields \begin{align*}
  \cF^{\hat U,\tilde U}(\ell)
  &= \sigma( Y[0], \hat{R}^\eta, \tilde{R}^\eta, 
    \hat{U}^0, \ldots, \hat{U}^\ell, 
    \tilde{U}^0, \ldots, \tilde{U}^\ell)\,, \\
  \cF^{\hat V,\tilde V}(\ell)
  &= \sigma( X[0], \hat{S}^\eta, \tilde{S}^\eta, 
    \hat{V}^0, \ldots, \hat{V}^\ell, 
    \tilde{V}^0, \ldots, \tilde{V}^\ell)\,,
\end{align*}
and analogously $\cF^{\hat Y,\tilde Y,\hat{\prxY},\tilde{\prxY}}(\ell)$ and $\cF^{\hat M,\tilde M}(\ell)$, which now include $\hat{R}^\eta, \tilde{R}^\eta$ and $\hat{S}^\eta, \tilde{S}^\eta$, respectively.
\begin{lem}
\label{lem:pair-orthogonality}
For all $1 \leq \ell\leq \ubl$ we have:
  \begin{align}
  \label{eq:pair-u-increments-orthogonal}
  \bbE[\hat U^{\ell}-\hat U^{\ell-1}
  \,|\,
  \cF^{\hat U,\tilde U}(\ell-1)]
  &=0\,,
  \\
  \label{eq:pair-y-increments-orthogonal}
  \bbE[\hat Y^{\ell}-\hat Y^{\ell-1}
  ~|~
  \cF^{\hat Y,\tilde Y,\hat{\prxY},\tilde{\prxY}}(\ell-1)]
  =
  \bbE[\hat Y^{\ell}-\hat Y^{\ell-1}
  ~|~
  \cF^{\hat U,\tilde U}(\ell-1)]
  &=
  0\,,
  \\
  \label{eq:pair-v-increments-orthogonal}
  \bbE[\hat V^{\ell}-\hat V^{\ell-1}
  ~|~
  \cF^{\hat V,\tilde V}(\ell-1)]
  &=0\,,
  \\
  \label{eq:pair-m-increments-orthogonal}
  \bbE[\hat M^{\ell}-\hat M^{\ell-1}
  ~|~
  \cF^{\hat M,\tilde M}(\ell-1)]
  =
  \bbE[\hat M^{\ell}-\hat M^{\ell-1}
  ~|~
  \cF^{\hat V,\tilde V}(\ell-1)]
  &=
  0
  \,.
  \end{align}

\begin{proof}
  We will show inductively that
  \eqref{eq:pair-u-increments-orthogonal} implies
  \eqref{eq:pair-y-increments-orthogonal}, which implies
  \eqref{eq:pair-v-increments-orthogonal}, which implies
  \eqref{eq:pair-m-increments-orthogonal}, which implies
  \eqref{eq:pair-u-increments-orthogonal} with $\ell+1$ in place of $\ell$.
We therefore proceed with the inductive argument: 

\begin{itemize}
\item \eqref{eq:pair-u-increments-orthogonal} $\Rightarrow$
  \eqref{eq:pair-y-increments-orthogonal}: 
for this we use
the recursion \eqref{e:state.evol.x.y.increments.CHAOTIC}. Note that $w^{\ell-1}(\hat{\prxY}^{\ell-1})$ is measurable with respect to the $\sigma$-field $\cF^{\hat Y,\tilde Y,\hat{\prxY},\tilde{\prxY}}(\ell-1)$, which in turn is contained in the $\sigma$-field 
$\cF^{\hat U,\tilde U}(\ell-1)$. Therefore
  \[
  \bbE[\hat Y^{\ell}-\hat Y^{\ell-1}
    ~|~
    \cF^{\hat U,\tilde U}(\ell-1)]
  = \frac{w^{\ell-1}(\hat{R}^\eta\hat{\prxY}^{\ell-1}) }
    {\hat{\gamma}(\ell,\delta,\eta)}
  \E[\hat{U}^\ell-\hat{U}^{\ell-1}
  ~|~ \cF^{\hat U,\tilde U}(\ell-1)
  ]=0\,,
  \]
where the last step is by the inductive hypothesis \eqref{eq:pair-u-increments-orthogonal}. This implies
\eqref{eq:pair-y-increments-orthogonal}.

\item \eqref{eq:pair-y-increments-orthogonal} $\Rightarrow$ \eqref{eq:pair-v-increments-orthogonal}:
it follows from the state evolution definition that the random variables $\hat{V}^{\ell,k}$ and $\tilde{V}^{\ell,k}$ are jointly gaussian. Using the covariance formula \eqref{eq:state-evolution-setup}, for $0\le k\le s \le \ell-1$ we have
  \begin{align*}
  \E[(\hat{V}^{\ell,k}-\hat{V}^{\ell-1,k})
  \tilde{V}^{s,k}]
  &\stackrel{\eqref{eq:state-evolution-setup}}{=} \E\bigg[
   \Big(h_{\chaos,\ell,k}(\hat{U}[[\ell]])-
  h_{\chaos,\ell-1,k}
  (\hat{U}[[\ell-1]])
  \Big)
    h_{\chaos,s,k}
    (\tilde{U}[[s]])
    \bigg]\\
  &\stackrel{
  \eqref{e:f.h.l.k.explicit.CHAOTIC}
  }{=} (a_k)^2
  \E[(\hat{Y}^\ell-\hat{Y}^{\ell-1})
    \tilde{Y}^s]
  = 0\,,
  \end{align*}
where the last step is by the inductive hypothesis \eqref{eq:pair-y-increments-orthogonal}. A similar calculation gives the same identity but with $\hat{V}^{s,k}$ in place of $\tilde{V}^{s,k}$. Since the $V$ are jointly gaussian, this implies \eqref{eq:pair-v-increments-orthogonal}.

\item \eqref{eq:pair-v-increments-orthogonal} $\Rightarrow$ \eqref{eq:pair-m-increments-orthogonal}: this follows from the recursion \eqref{e:state.evol.x.y.increments.CHAOTIC}, by the same reasoning as for \eqref{eq:pair-u-increments-orthogonal} $\Rightarrow$ \eqref{eq:pair-y-increments-orthogonal} above.

\item
\eqref{eq:pair-m-increments-orthogonal}
$\Rightarrow$
\eqref{eq:pair-u-increments-orthogonal}:
this follows by similar reasoning as for
\eqref{eq:pair-y-increments-orthogonal}
$\Rightarrow$
\eqref{eq:pair-v-increments-orthogonal} above.
\end{itemize}
The above steps give one complete round of the induction, and the claim follows.
\end{proof}
\end{lem}

We will exploit these orthogonality relations using the following lemma.

\begin{lem}\label{lem:conditional-gaussian-decorrelation}
Let $I$ be a finite, nonempty index set. Let 
  \[
    \lt(\Theta, 
    \vec{\hat{Z}}\equiv
    (\hat Z_i)_{i\in I},
    \vec{\tilde{Z}}
    \equiv(\tilde Z_i)_{i\in I}\rt)
  \]
  be a random variable where
  $(\vec{\hat{Z}},\vec{\tilde{Z}})$
  are centered, jointly gaussian, and independent of $\Theta$. 
  Assume that $\vec{\hat{Z}}$ and $\vec{\tilde{Z}}$ are equidistributed,
 and that, for distinct $i,k\in I$,
  \[
    \bbE[\hat Z_i\hat Z_k]
    =\bbE[\tilde Z_i\tilde Z_k]
    =\bbE[\hat Z_i\tilde Z_k]=0.
  \]
Let $F$ and $G$ be measurable functions
such that $F(\Theta,\hat Z)$ and 
$G(\Theta,\hat Z)$ have finite variance, and such that
	\beq\label{e:conditional-centering}
    \bbE\Big[F(\Theta,\hat Z) 
    	\,\Big|\, \Theta\Big] = 0 = 
    \bbE\Big[G(\Theta,\tilde Z) 
    	\,\Big|\, \Theta\Big]\,,
    \eeq
  $\Theta$-almost surely.
Then the covariance between $F(\Theta,\hat{Z})$ and $G(\Theta,\tilde{Z})$ is bounded by
  \beq
    \lt| \bbE[F(\Theta,\hat Z) G(\Theta,\tilde Z)] \rt|
    \le
    \rho_{\max} \cdot
    \bbE[F(\Theta,\hat Z)^2]^{1/2}
    \bbE[G(\Theta,\tilde Z)^2]^{1/2}\,, 
    \label{eq:conditional-gaussian-decorrelation-max}
  \eeq
where $\rho_i\in[-1,+1]$ denotes the correlation between $\hat{Z}_i$ and $\tilde{Z}_i$, and 
$\rho_{\max} 
   \equiv\max\{|\rho_i|: i\in i\}$.
 
\begin{proof}
Let $\sigma_i$ denote the standard deviation of $\hat{Z}_i$, and denote
  \[
    \hat W_i \equiv 
    \frac{\hat Z_i}{\sigma_i}\,, \quad
    \tilde W_i 
    \equiv \frac{\tilde Z_i}{\sigma_i}\,,
  \]
which are marginally standard gaussians. It follows from the assumptions that
$(\hat W_i,\tilde W_i)$ are independent over $i\in I$, and $\bbE[\hat W_i\tilde W_i]=\rho_i$.
Let $(h_m)_{m\geq 0}$ denote the orthonormal Hermite polynomials for the
one-dimensional standard gaussian measure. We note that if $W$ and $\tilde{W}$ are marginally standard gaussians, and $(W,\tilde{W})$ is jointly gaussian and $\rho$-correlated, then it is straightforward to verify that
	\beq\label{e:hermite.corr.ip}
	\E[ h_n(W) h_m(\tilde{W})]
	= \ind\{n=m\} \rho^n\,.\eeq
 For a multi-index $\vec n=(n_i)_{i\in I}\in\bbZ_{\ge 0}^I$, define the multivariate Hermite polynomial
  \[
    H_{\vec n}(w)
    \equiv\prod_{i\in I}h_{n_i}(w_i).
  \]
  We expand $F,G$ in the Hermite bases of $\hat W$ and $\tilde W$ respectively, conditional on $\Theta$, to obtain Hermite coefficients $a_{\vec n}(\Theta)$ and $b_{\vec n}(\Theta)$ such that
  \begin{align*}
    F(\Theta,\hat Z)
    &= \sum_{\vec n}a_{\vec n}(\Theta)
        H_{\vec n}(\hat W)\,, \\
    G(\Theta,\tilde Z)
    &= \sum_{\vec n}b_{\vec n}(\Theta)
        H_{\vec n}(\tilde W)\,.
  \end{align*}
  The assumption \eqref{e:conditional-centering} ensures $a_{\vec 0}(\Theta)=b_{\vec 0}(\Theta)=0$, $\Theta$-almost surely.
  Consequently, expanding in the Hermite bases and applying \eqref{e:hermite.corr.ip} gives
  \begin{align*}
    &\bigg|
      \bbE\Big[
        F(\Theta,\hat Z)G(\Theta,\tilde Z) \,\Big|\, \Theta
      \Big]
    \bigg|
    \le
      \sum_{\vec n\neq\vec 0}
         |a_{\vec n}(\Theta)b_{\vec n}(\Theta)|
         \prod_{i\in I}|\rho_i|^{n_i}\\
    &\qquad\le
      \rho_{\max}
      \bigg(\sum_{\vec n\neq\vec 0}
        a_{\vec n}(\Theta)^2
        \bigg)^{1/2}
      \bigg(\sum_{\vec n\neq\vec 0}
        b_{\vec n}(\Theta)^2
        \bigg)^{1/2}
    =\rho_{\max}
      \bbE\Big[F(\Theta,\hat Z)^2 
      	\,\Big|\, \Theta\Big]^{1/2}
      \bbE\Big[G(\Theta,\tilde Z)^2 
      	\,\Big|\, \Theta\Big]^{1/2}\,.
  \end{align*}
  Taking expectation over $\Theta$ and again applying the Cauchy--Schwarz inequality proves the conclusion.
\end{proof}
\end{lem}

\begin{proof}[\hypertarget{proof:thm.chaotic-iamp-is-chaotic}{Proof of Theorem~\ref{thm:chaotic-iamp-is-chaotic}}] Fix $\eta>0$, and recall
the order of parameters \eqref{eq:eta-delta-eps}. Throughout this proof, $o_\eta(1)$ denotes a quantity small in the in the limit $\linit \uparrow \infty$, then $\upp \uparrow 1$, then $\eta \downarrow 0$, all after the state evolution limit $N \uparrow \infty$, with $L,\epsilon,\delta$ fixed.
  By Lemma~\ref{lem:pair-orthogonality}, for any $j\le \ubl$, the hypothesis of Lemma~\ref{lem:conditional-gaussian-decorrelation} holds for
  \begin{align*}
    &\Big(
      \Theta^U, 
      (\hat{Z}_{U,R}, \hat{Z}_{U,1}, \ldots, \hat{Z}_{U,j}),
      (\tilde{Z}_{U,R}, \tilde{Z}_{U,1}, \ldots, \tilde{Z}_{U,j})
    \Big) \\
    &\equiv \Big( Y[0],
      (\hat Y^{-\linit/2-1}, \hat U^1, \hat U^2 - \hat U^1, \ldots, \hat U^j - \hat U^{j-1}),
      (\tilde Y^{-\linit/2-1}, \tilde U^1, \tilde U^2 - \tilde U^1, \ldots, \tilde U^j - \tilde U^{j-1})
    \Big)\,,
  \end{align*}
as well as for
  \begin{align*}
    &\Big(
      \Theta^V, 
      (\hat{Z}_{V,S}, \hat{Z}_{V,0}, \ldots, \hat{Z}_{V,j}),
      (\tilde{Z}_{V,S}, \tilde{Z}_{V,0}, \ldots, \tilde{Z}_{V,j})
    \Big) \\
    &\equiv \Big( X[0],
      (\hat V^{-\linit/2}, \hat V^0, \hat V^1 - \hat V^0, \ldots, \hat V^j - \hat V^{j-1}),
      (\tilde V^{-\linit/2}, \tilde V^0, \tilde V^1 - \tilde V^0, \ldots, \tilde V^j - \tilde V^{j-1})
    \Big)\,.
  \end{align*}
For $j\ge 0$, let $\Delta \hat{M}^j = \hat{M}^j - \hat{M}^{j-1}$, and analogously define $\Delta \tilde{M}^j$; recall from above \eqref{eq:hat-u-IAMP} that we initialize $\hat{M}^{-1} = \hat{M}^0 = \tilde{M}^{-1} = \tilde{M}^0 = 0$. We claim that for fixed $L,\epsilon,\delta$, there is a constant $C=C(L,\epsilon,\delta)<\infty$ such that, for all sufficiently small $\eta$ and all $0\le j \le \ubl$,
  \beq\label{e:2mt-bound-Y-M}
    \bbE[(\hat{Y}^j)^2]
    +\bbE[(\tilde{Y}^j)^2]
    +\bbE[(\Delta \hat M^j)^2]
    +\bbE[(\Delta \tilde M^j)^2]
    \le C
  \eeq
--- this follows from Lemmas~\ref{lem:chaotic-symmetrization}--\ref{lem:eta-to-0-chaotic-iamp}, recalling that the second moment of $Y^j$ is bounded as a consequence of Proposition~\ref{prop:Ising-IAMP-analysis}, while the second moment of $\Delta M^j$ is bounded as a consequence of the calculation \eqref{e:variance.Delta.M.calculation} from the proof of Lemma~\ref{l:gamma.estimates}, using assumptions \eqref{i:IAMP-main-coefs} and \eqref{i:IAMP-main-p} of Theorem~\ref{thm:IAMP-main}.
Recalling the notation \eqref{e:chaotic.w.sigma.b.state.evol},  since $b(t,x)$ takes values in $[-L, L]$ by assumption~\eqref{i:IAMP-main-coefs} of Theorem~\ref{thm:IAMP-main}, we have
  \beq\label{e:2mt-bound-b}
    \E[(\hat{b}^{\chaos,j})^2]
    +\E[(\tilde{b}^{\chaos,j})^2]
    \le 2L^2\,.
  \eeq 
  We will prove that
  \begin{align}\label{e:cross-correlations-small}
    \rho^V_j &\equiv \frac{\bbE[\hat{Z}_{V,j} \tilde{Z}_{V,j}]}{\bbE[(\hat{Z}_{V,j})^2]} = o_\eta(1)\,, &
    \rho^U_j &\equiv \frac{\bbE[\hat{Z}_{U,j} \tilde{Z}_{U,j}]}{\bbE[(\hat{Z}_{U,j})^2]} = o_\eta(1)\,,
  \end{align}
  for respectively $0\le j\le \ubl$ and $1\le j\le \ubl$.
  We first explain how \eqref{e:cross-correlations-small} implies the conclusion.
  We take in Lemma~\ref{lem:conditional-gaussian-decorrelation}
  \begin{align*}
    F(\Theta^U, \hat{Z}_{U,R}, \hat{Z}_{U,1}, \ldots, \hat Z_{U,\ubl})
    &= \hat{Y}^{\ubl}\,, &
    G(\Theta^U, \tilde{Z}_{U,R}, \tilde{Z}_{U,1}, \ldots, \tilde Z_{U,\ubl})
    &= \tilde{Y}^{\ubl}\,.
  \end{align*}
  Note that $\hat{Y}^{\ubl}$ has symmetric distribution, even conditional on $\Theta^U$, as its value is negated if we negate $(\hat{Z}_{U,R}, \hat{Z}_{U,1}, \ldots, \hat Z^U_{\ubl})$.
  Similarly $\tilde{Y}^{\ubl}$ has symmetric distribution conditional on $\Theta^U$, so \eqref{e:conditional-centering} is satisfied.
  Then, by state evolution  combined with Lemma~\ref{lem:conditional-gaussian-decorrelation} and \eqref{e:2mt-bound-Y-M}, we have
  \[
    \plim_{N\to\infty} \frac{(\hat{\by}^{\ubl},\tilde{\by}^{\ubl})}{N}
    = \bbE[\hat{Y}^{\ubl} \tilde{Y}^{\ubl}]
    = o_\eta(1)\,,
  \]
giving the conclusion of the theorem.
  It remains to prove \eqref{e:cross-correlations-small}.
  By Lemma~\ref{lem:decorrelate},
  \begin{align*}
    \rho^V_S &\equiv 
    \frac{\bbE[\hat{Z}_{V,S} \tilde{Z}_{V,S}]}
    	\bbE[(\hat{Z}_{V,S})^2]
	 = o_\eta(1)\,, & 
    \rho^U_R &\equiv 
    \frac{\bbE[\hat{Z}_{U,R} \tilde{Z}_{U,R}]}
	\bbE[(\hat{Z}_{U,R})^2]
     = o_\eta(1)\,.
  \end{align*}
  We will prove by induction on $j$ that
  \begin{align*}
    D_V(j) &\equiv \max\bigg\{
      |\rho^V_S|, \max_{0\le i\le j} |\rho^V_i|
    \bigg\} = o_\eta(1)
    \textup{ for all $0\le j\le \ubl$,}
    \\
    D_U(j) &\equiv \max\bigg\{
      |\rho^U_R|, \max_{1\le i\le j} |\rho^U_i|
    \bigg\} = o_\eta(1)
    \textup{ for all $ 1\le j\le \ubl$.}
  \end{align*}
For the base case of $D_V(0)$, note that $\rho^U_R = o_\eta(1)$ implies
  \begin{align*}
    |\rho^V_0| 
    &= \frac{\bbE[\hat V^0 \tilde V^0]}{\bbE[(\hat V^0)^2]} 
    = \upp \frac{\bbE[\hat Y^0 \tilde Y^0]}{\bbE[(\hat Y^0)^2]} 
    = \upp \frac{\bbE[\hat R^\eta \tilde R^\eta (Y^0)^2]}{\bbE[(\hat R^\eta)^2 (Y^0)^2]} 
    = \upp \frac{\bbE[\hat R^\eta \tilde R^\eta]}{\bbE[(\hat R^\eta)^2]}\\
    &= \upp \frac{\bbE[
    \sign_\eta(\hat{Y}^{-\linit/2-1})
    \sign_\eta(\tilde{Y}^{-\linit/2-1})
    ]}{\bbE[
    (\sign_\eta(\hat{Y}^{-\linit/2-1}))^2
    ]}
    = o_\eta(1)\,.
  \end{align*}
For the inductive step, suppose that $D_V(j)=o_\eta(1)$ for some $0\le j < \ubl$. Note that both $\Delta\hat M^j$ and $\hat b^{\chaos,j}$
  are functions of $(\Theta^V, \hat{Z}_{V,S}, \hat{Z}_{V,0}, \ldots, \hat{Z}_{V,j})$, and have symmetric distribution conditional on $\Theta^V = X_0$, as their values are negated if we negate $(\hat{Z}_{V,S}, \hat{Z}_{V,0}, \ldots, \hat{Z}_{V,j})$.
  So, Lemma~\ref{lem:conditional-gaussian-decorrelation} and \eqref{e:2mt-bound-Y-M}--\eqref{e:2mt-bound-b} give
  \[
    \Big|\bbE[\Delta\hat M^j\Delta\tilde M^j]\Big|
    +
    \Big|\bbE[\hat b^{\rm ch,j}\tilde b^{\rm ch,j}]\Big|
    =o_\eta(1)\,.
  \]
Applying \eqref{eq:chaotic-U-covariance-recursion} gives,
  \[
    \bbE[\hat Z_{U,j+1} \tilde Z_{U,j+1}] = \upp \alpha p^j \cdot \bbE[\Delta\hat M^j\Delta\tilde M^j]
    +\upp \alpha\delta \bbE[\hat b^{\chaos,j}\tilde b^{\chaos,j}]+\upp^{\linit+2j}c_j(\epsilon)^2
    = o_\eta(1)\,.
  \]
By Lemmas~\ref{lem:chaotic-symmetrization}--\ref{lem:eta-to-0-chaotic-iamp} and \eqref{e:U-increment-var-lb},
  \[
    \bbE[(\hat Z_{U,j+1})^2]
    = \bbE[(U^{j+1}-U^j)^2] + o_\eta(1)
    \ge \frac{\delta}{5L^2} - o_\eta(1)
  \]
  is bounded away from $0$ independently of $\eta$ (recalling the limiting regime \eqref{eq:eta-delta-eps}). Combining the last two displays shows $\rho^U_{j+1} = o_\eta(1)$, which implies $D_U(j+1) = o_\eta(1)$.

  Next, note that for $0\le i\le j+1$, $\hat Y^i$ is a function of $(\Theta^U, \hat{Z}_{U,R}, \hat{Z}_{U,1}, \ldots, \hat Z_{U,i})$, and has symmetric distribution conditional on $\Theta^U = Y[0]$ by the same argument as above.
  Thus Lemma~\ref{lem:conditional-gaussian-decorrelation} combined with \eqref{e:2mt-bound-Y-M} gives
  \[
    \Big|\bbE[\hat Y^i \tilde Y^i]\Big| = o_\eta(1)
  \]
  for all such $i$.
  Summing over \eqref{eq:state-evolution-setup} and using \eqref{e:f.h.l.k.explicit.CHAOTIC}, we find
  \[
    \bbE[\hat V^r\tilde V^s]
    =
    \upp p^{\min(r,s)}
    \bbE[\hat Y^r\tilde Y^s]\,.
  \]
  Together with Lemma~\ref{lem:pair-orthogonality}, this gives
  \begin{align*}
    \bbE[\hat Z_{V,j+1} \tilde Z_{V,j+1}] 
    &= \bbE[(\hat V^{j+1} - \hat V^j)(\tilde V^{j+1} - \tilde V^j)]
    = \bbE[\hat V^{j+1} \tilde V^{j+1}] - \bbE[\hat V^j \tilde V^j] \\
    &= \upp p^{j+1} \bbE[\hat Y^{j+1} \tilde Y^{j+1}]
    - \upp p^j \bbE[\hat Y^j \tilde Y^j]
    = o_\eta(1)\,.
  \end{align*}
By Lemmas~\ref{lem:chaotic-symmetrization}--\ref{lem:eta-to-0-chaotic-iamp} and \eqref{e:V-increment-var},
  \[
    \bbE[(\hat Z_{V,j+1})^2]
    = \bbE[(V^{j+1}-V^j)^2] + o_\eta(1)
    \ge \delta s_\delta(q_0+\ell\delta)^2 - o_\eta(1)
    \ge \frac{\delta}{L} - o_\eta(1)\,,
  \]
  where the last inequality holds because
  \[
    s_\delta(t)^2
    = p(t+\delta)
    +\frac{t}{\delta}(p(t+\delta)-p(t))
    \ge p(t+\delta)
    \ge p(0) \ge \frac1L\,.
  \]
  Thus $\bbE[(\hat Z_{V,j+1})^2]$ is bounded away from $0$ independently of $\eta$. Hence $\rho^V_{j+1} = o_\eta(1)$, which implies $D_V(j+1) = o_\eta(1)$.
  This completes the induction and implies the claim.
\end{proof}

\begin{proof}[\hypertarget{proof:t.IAMP.main.chaotic}{Proof of Theorem~\ref{thm:IAMP-main}\ref{i:IAMP-main-chaotic}}] Recall $\bz = (\bG, \bg^\aux)$.
  As in the discussion just above Proposition~\ref{ppn:chaotic-amp-still-works}, let $\tilde{\cA}_N$ denote the IAMP iteration $\IAMP_{\chaos,N}(p,\sigma,b,w,\delta,\eta)$, and let $\breve{\cA}_N$ be the the $C(L,\epsilon)$-Lipschitz algorithm constructed from $\tilde{\cA}_N$ as in \eqref{e:breve-cA-N}.
  We can ensure the conclusion of Proposition~\ref{ppn:chaotic-amp-still-works} holds with $\iota/2$ in place of $\iota$ by taking $\delta, \eta$ sufficiently small in the limiting regime \eqref{eq:eta-delta-eps}.
  We then take
  \begin{align} \nonumber
    \check{\cA}_N(\bz) &= \breve{\cA}_N(\bz) - \bbE[\breve{\cA}_N(\bz)]\,, \\ 
    \cA_N(\bz) &= \frac{N^{1/2}}{\bbE[\|\check{\cA}_N(\bz)\|^2]^{1/2}} \cdot \check{\cA}_N(\bz)
    \equiv c_N \check{\cA}_N(\bz)
    \,.
    \label{e:chaotic-iamp-final-construction}
  \end{align}
  We abbreviate $\chi \equiv \chi_{\cA_N}$ and $\breve\chi \equiv \chi_{\breve{\cA}_N}$, and note that these are related via
  \[
    \chi(\upp) = \frac{\breve\chi(\upp) - \breve\chi(0)}{1 - \breve\chi(0)}\,.
  \]
  The construction \eqref{e:chaotic-iamp-final-construction} ensures $\bbE [\|\cA_N(\bz)\|^2] = N$ and $\chi(0) = 0$.

As previously, let $\hat{\by}^{\ubl}$ and $\tilde{\by}^{\ubl}$ be the outputs of $\IAMP_{\chaos,N}(p,\sigma,b,w,\delta,\eta)$ on $\upp$-correlated inputs $\bz$ and $\tilde{\bz}$, as specified by Definitions \ref{d:chaotic.IAMP.phase.I} and \ref{d:chaotic.IAMP}.
  Since in the limiting regime \eqref{eq:eta-delta-eps} we have $1-\upp \ll \eta \ll \iota'$, Lemma~\ref{l:chi-approx-by-amp-correlation} and Theorem~\ref{thm:chaotic-iamp-is-chaotic} imply that
  \[
    \breve\chi(1-\iota') \le \breve\chi(\upp) \le \plim_{N\rightarrow\infty} \frac{(\hat{\by}^{\ubl},\tilde{\by}^{\ubl})}{N} + \iota' + o_N(1) \le 2\iota'\,.
  \]
  Since $\breve\chi$ is nondecreasing we also have $\breve\chi(0) \le 2\iota'$.
  Thus (since $\iota' \ll \iota$)
  \[
    \chi(1-\iota)
    \le \chi(1-\iota') 
    \le \frac{\breve\chi(1-\iota')}{1 - \breve\chi(0)} 
    \le 2\breve\chi(1-\iota') 
    \le 4\iota'
    \le \iota\,.
  \]
  Similarly, we note that
  \[
  \frac{1}{(c_N)^2}
   =\frac{\bbE[\|\check{\cA}_N(\bz)\|^2]}{N}
    = 1 - \breve\chi(0) \ge 1-2\iota'.
  \]
  Since $\check{\cA}_N$ is $C(L,\epsilon)$-Lipschitz, the construction \eqref{e:chaotic-iamp-final-construction} implies that $\cA_N$ is $2C(L,\epsilon)$-Lipschitz.
  Finally, we will check that
  \begin{align}
    \bbP\Big(\cA_N(\bz) \in \Sigma(\iota) \Big)
    &\ge 1-e^{-cN}\,, \\
    \label{e:chaotic-AN-2}
    \bbP\Big(
      \bbW_2\lt(\mu_{\bG,\sym}(\cA_N(\bz)), \Law(|X(1)|) \rt) \le \iota
    \Big) &\ge 1-e^{-cN}\,,  \\
    \bbW_2\Big(\mu^{\Ising}(\cA_N), \cP(\{\pm 1\}) \Big) &\le \iota\,, \\
    \label{e:chaotic-AN-4}
    \bbW_2\Big(\sym(\mu(\cA_N)), \Law(|X(1)|) \Big) &\le \iota\,.
  \end{align}
  These will each follow from the fact that $\breve{\cA}_N$ satisfies the conclusion of Proposition~\ref{ppn:chaotic-amp-still-works} with $\iota/2$ in place of $\iota$.
  We will only prove \eqref{e:chaotic-AN-2} and \eqref{e:chaotic-AN-4}, as the other two estimates are only simpler.
To this end, write
$\bm \equiv \bbE[\breve{\cA}_N(\bz)]$.
  It was proved above that $1 \le (c_N)^2 \le (1-2\iota')^{-1}$.
  Note that
  \[
    \mu_{\bG}(\cA_N(\bz))
    = \frac{1}{M} \sum_{a=1}^M \delta \lt(\frac{(\cA_N(\bz), \bg^a)}{\sqrt{N}}\rt)
    = \frac{1}{M} \sum_{a=1}^M \delta \lt(\frac{c_N(
    \breve{\cA}_N
    (\bz) - \bm, \bg^a)}{\sqrt{N}}\rt)
  \]
It follows that
  \begin{align*}
    \bbW_2(\mu_{\bG}(\cA_N(\bz)), \mu_{\bG}(\breve\cA_N(\bz)))
    &\le \frac{1}{M} 
    \sum_{a=1}^M \bigg(
    \frac{(c_N-1) |(\breve \cA_N(\bz),\bg^a)| + |(\bm, \bg^a)|}
    {N^{1/2}}\bigg)^2 \\
    &\le \frac{2}{MN} \bigg(
      (c_N - 1)^2
      \|\bG \breve \cA_N(\bz)\|^2
      + \|\bG \bm\|^2
    \bigg) \\
    &\le \frac{2}{MN} 
    (\|\bG\|_{\op})^2 \bigg(
      (c_N - 1)^2
      \|\breve \cA_N(\bz)\|^2
      + \|\bm\|^2
    \bigg)\,.
  \end{align*}
  Note that $\|\bm\|^2 / N = \breve\chi(0) \le \iota'$.
  By Lemma~\ref{l:wishart}, $\bbP(\|\bG\|_{\op} \le C_0 (1+\alpha^{1/2}) \sqrt{N}) \ge 1-e^{-cN}$ for $C_0$ given therein.
  Also, since $\breve \cA_N$ is $C(L,\epsilon)$-Lipschitz, $\bbP(\|\breve \cA_N(\bz)\| \le 2\sqrt{N}) \ge 1-e^{-cN}$.
  On the intersection of these events we have
  \[
    \bbW_2\Big(\mu_{\bG,\sym}(\cA_N(\bz)), \mu_{\bG,\sym}(\breve\cA_N(\bz))\Big)
    \le \bbW_2\Big(\mu_{\bG}(\cA_N(\bz)), \mu_{\bG}(\breve\cA_N(\bz))\Big)
    \le O(\iota')\,.
  \]
  Since $\iota' \ll \iota$ and we constructed $\breve\cA_N$ in Proposition~\ref{ppn:chaotic-amp-still-works} so that 
  \[
    \bbP\bigg(
      \bbW_2\Big(\mu_{\bG,\sym}(\breve\cA_N(\bz)), \Law(|X(1)|) \Big) \le \frac{\iota}{2}
    \bigg) \ge 1-e^{-cN}\,,
  \]
  the conclusion \eqref{e:chaotic-AN-2} follows.
  Further, by convexity of $(\bbW_2)^2$,
  \begin{align*}
   &\bbW_2\Big(\mu(\cA_N), \mu(\breve \cA_N)\Big)^2
    = \bbW_2\Big(\bbE \mu_{\bG}(\cA_N(\bz)), \bbE \mu_{\bG}(\breve \cA_N(\bz))\Big)^2\\
    &\qquad\le \bbE \bigg[
    \bbW_2\Big(\mu_{\bG}(\cA_N(\bz)), \mu_{\bG}(\breve \cA_N(\bz))\Big)
    ^2\bigg]\,.
  \end{align*}
  We showed above that $\bbW_2(\mu_{\bG}(\cA_N(\bz)), \mu_{\bG}(\breve \cA_N(\bz))) \le O(\iota')$ with probability $1-e^{-cN}$.
  A standard calculation also shows that
  \[
    \bz \mapsto \lt(1 + \bbW_2^2(\mu_{\bG}(\cA_N(\bz)), \mu_{\bG}(\breve \cA_N(\bz)))\rt)^{1/4}
  \]
  is $O(1/N)$-Lipschitz, and therefore the right-hand side has subgaussian fluctuations of order $N^{-1/2}$.
  From this it follows that
  \begin{align*}
    \bbW_2^2(\sym(\mu(\cA_N)), \sym(\mu(\breve \cA_N)))
    &\le \bbW_2^2(\mu(\cA_N), \mu(\breve \cA_N)) \\
    &\le \bbE \bbW_2^2(\mu_{\bG}(\cA_N(\bz)), \mu_{\bG}(\breve \cA_N(\bz)))
    \le O((\iota')^2)\,.
  \end{align*}
  Since $\iota' \ll \iota$ and we constructed $\breve\cA_N$ so that 
  \[
    \bbW_2\lt(\sym(\mu(\breve \cA_N)), \Law(|X(1)|) \rt) \le 
    \frac{\iota}{2}\,,
  \]
  the conclusion \eqref{e:chaotic-AN-4} follows.
\end{proof}

\subsection{Constructing solution trees for even activations}
\label{ss:even.solution.trees}

In this subsection we show that the chaotic IAMP algorithms of \S\ref{subsec:p=1-IAMP} are also able to construct large approximately ultrametric trees of outputs in the symmetric setting. The presence of such trees, possibly within an ensemble of correlated disorders, is characteristic of the algorithmically feasible side of the branching OGP threshold. However, we note that this subsection is not used elsewhere in the paper.  

\begin{cor}
\label{cor:b-iamp-tree}
Fix $0<q_1<q_2<\dots<q_m<1$. For small enough $\delta,\epsilon>0$ there exists $c>0$ such that for $N$ sufficiently large and $K=e^{cN}$, with probability $1-e^{-cN}$ the following holds for any $(L,\epsilon)$-good $(p,\sigma,b,w)$: 
 there exist points $(\by_w :w\in [K]^m)$, taking values in $\Sigma_N(\iota)$ (or $S_N(\iota)$) such that
  \beq\label{eq:all-outputs-good}
  \bbW_2\Big(
  \mu^{\sym}_{\bG}(\by_w),
  \mu^{\sym}(p,\sigma,b,w)
  \Big) \le\epsilon\eeq
for all $w\in [K]^m$, and
  \beq\label{eq:overlaps-ultrametric}
  \lt|\frac{( \by_w,\by_{w'})}{N}-q_{w\wedge w'}\rt|
  \leq \epsilon\eeq
for all $w,w'\in [K]^m$. In fact, we can take $\by_w=\cA_{\chaos}(\bz_w)$ where $(\bz_w:w\in [K]^m)$ is an ensemble of gaussians that are suitably correlated with $\bz\equiv(\bG,\bg^\aux)$.

\begin{proof}
Choose $\epsilon\ll q_1$ small. We will take $c$ small in order to union bound over $e^{O(cN)}$ exponentially unlikely events in the proof. Let $\cA_{\chaos}$ be as defined in 
Proposition~\ref{ppn:chaotic-amp-still-works}, taking inputs $\bz=(\bG,\bg^\aux)$. Let $\chi$ be the correlation function associated to $\cA_{\chaos}$. By Theorem~\ref{thm:chaotic-iamp-is-chaotic}, we have
  \[p_i\equiv \chi^{-1}(q_i)
  \ge 1-\epsilon\]
for all $1\leq i\leq m$. Let $T$ denote the $K$-ary tree of depth $m$, and identify the leaves with $[K]^m$. Generate a jointly gaussian ensemble
$(\bz_w: w\in [K]^m)$ such that $\bz_{1^m}=\bz=(\bG,\bg^\aux)$, and such that the correlation between $\bz_u$ and $\bz_v$ is $p_{|u\wedge v|+1}$ for any $u,v\in [K]^m$. Let $\by_w\equiv \cA_{\chaos}(\bz_w)$ for $w\in[K]^m$. 
By taking $\eta$ small enough in the 
definition of $\cA_{\chaos}$, and 
applying
Proposition~\ref{ppn:chaotic-amp-still-works} together with a union bound over $w\in[K]^m$, we obtain
that with probability at least $1- e^{-10cN}$, we have 
  \begin{equation}
  \label{eq:outputs-good-for-perturbed-disorder}
  \bbW_2\Big(
  \mu^{\sym}_{\bG_w}(\by_w),
  \mu^{\sym}(p,\sigma,b,w)
  \Big)\leq \epsilon
  \end{equation}
for all $w\in[K]^m$. Another union bound gives that with probability at least 
 $1-e^{-10cN}$, we have
  \[
  \begin{aligned}
  \max\bigg\{
  \frac{\|\by_w\|}{N^{1/2}}
  : w\in[K]^m\bigg\}
  &\le 2\,,\\
  \max\bigg\{
  \frac{\|\bG_w-\bG_{w'}\|_{\op}}{M^{1/2} + N^{1/2}}
  : w,w'\in[K]^m
  \bigg\}
  &\leq 10\epsilon^{1/2}\,.
  \end{aligned}
  \]
On the event that the above conditions hold, we have deterministically
  \[
  \bbW_2\Big(\mu^{\sym}_{\bG_w}(\by_w),\mu^{\sym}_{\bG}(\by_w)\Big)
  \leq 
  \bbW_2\Big(\mu_{\bG_w}(\by_w),\mu^{\sym}_{\bG}(\by_w)\Big)
  \leq 
  \frac{\|(\bG_w-\bG)\by_w\|}{
    (MN)^{1/2}}
  \leq 
  O(\epsilon^{1/2}).
  \]
Combining with \eqref{eq:outputs-good-for-perturbed-disorder} and adjusting $\epsilon$ gives \eqref{eq:all-outputs-good}. Next we note that
  \[
  \bbE\bigg[\frac{( \by_w,\by_{w'})}{N}\bigg]
  =
  \chi(p_{w\wedge w'})=q_{w\wedge w'}\,,
  \]
so \eqref{eq:overlaps-ultrametric} follows from concentration of measure (again taking a union bound with $c$ small enough).  
\end{proof}
\end{cor}

\begin{rmk}
A related approach to constructing multiple solutions is to introduce explicit branching into the IAMP algorithm. This type of ``branching IAMP'' was proposed in \cite{alaoui2020algorithmic} for mean-field spin glasses, and analyzed in \cite[Sec.\ 4]{sellke2021optimizing} and \cite[Sec.\ 3.3]{huang2024optimization}, which we refer to for further details.
Branching IAMP has the advantage of applying in situations where the optimal function $p$ equals one on some nontrivial interval $[q_1,1]$; informally this means the root of the ``algorithmic ultrametric solution tree'' has self-overlap strictly between zero and one --- we expect this to hold when $\psi$ is mildly asymmetric in some sense. In our setting, it works as follows: suppose $p(\tilde q_1)=1$ for some $\tilde q_1\in (0,1)$, and fix 
  \[
  0<\tilde q_1<\tilde q_2<\dots<\tilde q_k<1.
  \]
  and choose $\delta$ small compared to these constants. Then set $q_i=\delta\lfloor \tilde{q}_i/\delta\rfloor$ for each $i$, and let $\ell^\delta(q_i)
\equiv \ell(q_i) \equiv q_i/\delta
= \lfloor \tilde{q}_i/\delta\rfloor$. 
For each $0\leq i\leq k$, let $\hat\bg^i\sim \cN(0,I_N)$ be a standard gaussian ``perturbation'' vector.
At most times $\ell$, we perform the basic updates
\eqref{eq:proxy-y-amp-def}--\eqref{eq:x-amp-def} from Definition~\ref{d:iamp.v1}. The only exception is the following: at times $\ell=\ell(q_i)$, we perform the same updates \eqref{eq:u-amp-def}--\eqref{eq:x-amp-def}, but the updates \eqref{eq:proxy-y-amp-def}--\eqref{eq:y-amp-def} are replaced by
        \beq
  \label{eq:branching-IAMP-perturbed-update}
        \begin{aligned}
  \by^{\ell(q_i)+1}
  &=
  \by^{\ell(q_i)}
  +
  \frac{\bw^{\ell(q_i)}\odot \hat{\bg}^i}
  {
  \gamma(\ell(q_i),\delta)
  }
  \bbE[(U^{\ell(q_i)+1}-U^{\ell(q_i)})^2]^{1/2}
    ;
        \\
        \quad\quad
        \prxby^{\ell(q_i)+1}
  &=
  \prxby^{\ell(q_i)}
  +
  \frac{\bw^{\ell(q_i)}\odot \hat{\bg}^i}
  {
  \Gamma(\ell(q_i),\delta)
  }
  \bbE[(U^{\ell(q_i)+1}-U^{\ell(q_i)})^2]^{1/2}
  .
        \end{aligned}
        \eeq
The \hyperlink{proof:t.IAMP.main.main}{proof of Theorem~\ref{thm:IAMP-main}\ref{i:IAMP-main-main}} still holds here with minor modifications.
Further, for $0\leq j\leq k$ one can take 
  \[
  \begin{aligned}
  \hat\bg^{i,1}=\hat\bg^{i,2}\sim\cN(0,I_N) 
    &\textup{ for $i\le j$,}\\
  \hat\bg^{i,1},\hat\bg^{i,2}
    \stackrel{\textup{i.i.d.}}{\sim}\cN(0,I_N)
    &\textup{ for $i>j$.}
  \end{aligned}
  \]
Additionally, recalling \eqref{eq:u-amp-def},
we can take
  \[\begin{aligned}
  \bar{\bg}^{\ell,1}=\bar{\bg}^{\ell,2}
  \sim\cN(0,I_N) 
  &\textup{ for $\ell \le \ell(q_j)$,}\\
  \bar{\bg}^{\ell,1},\bar{\bg}^{\ell,2}
    \stackrel{\textup{i.i.d.}}{\sim}\cN(0,I_N)
  &\textup{ for $\ell > \ell(q_j)$.}
  \end{aligned}
  \]
Then the resulting algorithmic outputs $\by^{\ubl,1},\by^{\ubl,2}$ can be shown to have asymptotic overlap $q_j$, so long as $p(\tilde q_j)=1$ holds.
This is because for any pair of of trajectories which first differ at step $q_j$, one may show by induction that the future increments of the trajectories remain approximately orthogonal, i.e.,
    \[
    \plim_{N\to\infty}
    \frac{(\by_1^{\ell_1+1}-\by_1^{\ell_1}, \by_2^{\ell_2+1}-\by_2^{\ell_2} )}
      {N}=0
  \quad\textup{for all } \ell_1,\ell_2\geq \ell(q_j).
    \]
Using this and Lipschitz concentration, one may construct a branching tree of outputs, with combinatorial depth $k$ and width $e^{cN}$, such that pairs with least common ancestor at depth $j$ have overlap in the interval $[\tilde q_j-\epsilon,\tilde q_j+\epsilon]$, for any desired $\epsilon>0$ (when e.g. $\delta\ll\epsilon$).\end{rmk}

\fi

\pagebreak\section{Equivalence and continuity of stochastic control problems}
\label{sec:alternate-diffusions}

\def\cB{{\mathcal{B}}}
\def\Law{{\mathsf{Law}}}
\def\BUDGETS{{\mathsf{BUDGETS}}}
\def\Lip{{\mathsf{Lip}}}

\def\sM{{\mathscr{M}}}

\def\si{{\mathsf{i}}}
\def\sii{{\mathsf{ii}}}
\def\siii{{\mathsf{iii}}}
\def\siv{{\mathsf{iv}}}
\def\sv{{\mathsf{v}}}
\def\svi{{\mathsf{vi}}}
\def\svii{{\mathsf{vii}}}
\def\sviii{{\mathsf{viii}}}
\def\six{{\mathsf{ix}}}
\def\sx{{\mathsf{x}}}
\def\sxi{{\mathsf{xi}}}
\def\sxii{{\mathsf{xii}}}
\def\sideal{{\mathsf{ideal}}}
\def\sidealast{{\mathsf{ideal}^*}}

\iffull
% !TEX root = main.tex

For an algorithm $\cA_N$, recall again the definitions \eqref{e:mu.of.Alg} and \eqref{e:mu.Ising.of.Alg} of the averaged inner product and coordinate distributions, $\mu(\cA_N)$ and $\mu^{\Ising}(\cA_N)$.
Theorem~\ref{thm:BOGP-hardness-main} implies that for sequences of $O(1)$-Lipschitz $\cA_N$, all subsequential limits of $\mu(\cA_N)$ in $\bbW_2$ can be approximated by an endpoint measure $\Law(X(1))$, where $X(t)$ solves an SDE from a certain class. Conversely, Theorem~\ref{thm:IAMP-main} shows that the endpoint measures of a related class of SDEs can all be attained in the same sense by a sequence of $O(1)$-Lipschitz algorithms.
In this section, we show that these two sets of endpoint measures in fact coincide, in an appropriate limiting regime of parameters. Moreover, we show that they coincide with the set of endpoint measures $\ocM^{\Ising}(\alpha)$ (or $\ocM^{\sph}(\alpha)$ in the spherical setting) introduced in Definition~\ref{d:achievable-msrs}.
This is the set of endpoint measures of a much cleaner control problem, where all $o_L(1)$ terms in the aforementioned SDEs are set to $0$, the coefficients $(b,\sigma,w)$
can be general progressively measurable processes subject to the budget constraints, and the starting time $q_0$ is $0$. These endpoint measures are thus the achievable $\bbW_2$-subsequential limits of $\mu(\cA_N)$ for $O(1)$-Lipschitz $\cA_N$. We further show that $\ocM^{\Ising}(\alpha)$ coincides with $\ocM^{\Ising, \concave}(\alpha)$ from Definition~\ref{d:achievable-msrs}: this is the set of endpoint measures of the same control problem where the function $p$ must be concave with $p(0)=0$; thus such $p$ suffice to attain all endpoint measures.

We next state the results of this section formally. Recall $\cP_2(\bbR)$ denotes the set of square-integrable Borel probability measures on $\bbR$.
We let $\cP(\{\pm 1\})$ denote the set of Borel probability measures supported on $\{\pm 1\}$, and recall the point-to-set $\bbW_2$ distance defined in \eqref{def:W2-to-set}. For a Lipschitz algorithm $\cA$ and a probability measure $\mu\in\cP_2(\R)$, it will be useful to abbreviate
	\beq\label{e:Wasserstein.alg.err.msr}
	\WERR(\cA,\mu)
	\equiv\max\bigg\{
	\bbW_2\Big(\mu(\cA_{N_j}), \mu\Big),
        \bbW_2\Big(\mu^{\Ising}(
        	\cA_{N_j}), \cP(\{\pm 1\})\Big)
		\bigg\}\,.
	\eeq 
For $L,\epsilon > 0$, define
\begin{align}
    \label{e:cM.Lip}
    \cM^{\Lip}(\alpha;L,\epsilon) 
    &= \bigg\{
        \begin{array}{ll}
        \mu \in \cP_2(\bbR): & \text{$\exists$ $L$-Lipschitz algorithms $(\cA_{N_j})_{j\ge 1}$ with} \\
        & \WERR(\cA_{N_j},\mu) \le \epsilon
        \end{array}
    \bigg\}\,,\\
    \label{e:sM.Lip}
    \sM^{\Lip}(\alpha) 
    &\equiv\adjustlimits
    \bigcap_{\epsilon > 0} \bigcup_{L > 0} \cM^{\Lip}(\alpha;L,\epsilon)\,.
	\end{align}
The main result of this section is the following: 

\begin{thm}[$\bbW_2$-limits of Lipschitz algorithms for the Ising perceptron]
\label{thm:control.problems.main}
Let $\ocM^{\Ising}(\alpha)$ and $\ocM^{\Ising, \concave}(\alpha)$ be as specified by Definition~\ref{d:achievable-msrs}.
We then have
\[
    \sM^{\Lip}(\alpha) 
    = \ocM^{\Ising}(\alpha) = \ocM^{\Ising, \concave}(\alpha)
\]
for $\sM^{\Lip}(\alpha)$ as defined by \eqref{e:sM.Lip}. 
\end{thm}

\begin{rmk}\label{rmk:control.problems.spherical}
    For the spherical perceptron, we have the analogues of Theorem~\ref{thm:control.problems.main} and Theorem~\ref{thm:control.problems.sym} (see below), where we omit the constraint $\bbW_2(\mu^{\Ising}(\cA_{N_j}), \cP(\{\pm 1\})) \le \epsilon$ in \eqref{e:cM.Lip} and \eqref{e:cM.Lip.sym} and replace $\ocM^{\Ising}(\alpha)$ and $\ocM^{\Ising, \concave}(\alpha)$ with $\ocM^{\sph}(\alpha)$ and $\ocM^{\sph, \concave}(\alpha)$. The proof is strictly simpler: we set $w_t \equiv 1$ in Definition~\ref{d:SDE.conds} and drop the process $Y$ of  \eqref{e:formal.sde.ising} from the analysis. Thus we focus on the Ising setting in this section.
\end{rmk}

\begin{rmk}    \label{r:control.problems.interpretation}    In the \hyperlink{p:thm.control.problems.main}{proof of Theorem~\ref{thm:control.problems.main}} below, we show $\sM^{\Lip}(\alpha)$ coincides with several other classes of measures 
    $\sM^\star(\alpha,0)$, for $\star$ in \eqref{e:star.options}, specified by Definitions~\ref{d:SDE.conds}--\ref{d:measure.classes} below. Out of these, the $\star=\BOGP$ and $\star=\IAMP$ classes are precisely those described by Theorems~\ref{thm:BOGP-hardness-main} and \ref{thm:IAMP-main}, respectively (see Lemmas~\ref{l:IAMP.in.Lip} and \ref{l:Lip.in.BOGP}), while the $\star = \sideal$ and $\star = \sidealast$ classes coincide with $\ocM^{\Ising}(\alpha)$ and $\ocM^{\Ising, \concave}(\alpha)$ (see Remark~\ref{r:measure.classes.simplification}).
    The latter two formulations are by far the most convenient for direct analysis of the control problem, and are exclusively used in \cite{bogpinprogress}.
\end{rmk}

Most of this section is devoted to the proof of Theorem~\ref{thm:control.problems.main}.
The next theorem is proved in \S\ref{ss:sdes.symmetric.endpoints} by adapting this argument, and shows that if we are only interested in the \textbf{symmetrized} endpoint measure $\Law(|X(1)|)$, then the set of controls we need to consider simplifies even further. Analogously to \eqref{e:Wasserstein.alg.err.msr}, for a Lipschitz algorithm $\cA$ and a measure $\mu\in\cP_2(\R_{\ge0})$, denote
	\beq\label{e:WERR.sym.notation}
	\WERR^\sym(\cA,\mu)
	\equiv
	 \max\bigg\{
        \bbW_2\Big(\sym(\mu(\cA_{N_j})), \mu\Big),
        \bbW_2\Big(\mu^{\Ising}(\cA_{N_j}),
        \cP(\{\pm 1\})\Big)
        \bigg\}\,.
	\eeq
Then, analogously to \eqref{e:cM.Lip} and \eqref{e:sM.Lip}, we define
\begin{align}
    \label{e:cM.Lip.sym}
    \cM^{\Lip,\sym}(\alpha;L,\epsilon) 
    &\equiv \bigg\{
        \begin{array}{ll}
        \mu \in \cP_2(\bbR_{\ge 0}): & \text{$\exists$ $L$-Lipschitz algorithms $(\cA_{N_j})_{j\ge 1}$ with} \\
        & \WERR^\sym(\cA_{N_j},\mu)\le \epsilon
        \end{array}
    \bigg\} \,,\\
    \label{e:sM.Lip.sym}
    \sM^{\Lip,\sym}(\alpha) 
    &\equiv \adjustlimits
    \bigcap_{\epsilon > 0}
    \bigcup_{L > 0} \cM^{\Lip}(\alpha;L,\epsilon)\,.
\end{align}
We then have the following:

\begin{thm}
\label{thm:control.problems.sym}
Let $\ocM^{\Ising,\sym}(\alpha)$ be as specified in Definition~\ref{d:achievable-msrs-sym}.
We then have
\[
    \sM^{\Lip,\sym}(\alpha) 
    = \ocM^{\Ising,\sym}(\alpha)
\]
for $\sM^{\Lip,\sym}(\alpha)$ as defined by \eqref{e:sM.Lip.sym}.
\end{thm}
\begin{rmk}
    \label{r:control.problems.sym.interpretation}
    The class 
    $\ocM^{\Ising,\sym}(\alpha)$
    corresponds to ``ideal'' controls similarly as above, but with the further restriction $p\equiv 1$; see also Remark~\ref{r:measure.classes.simplification}. As a result, the drift control $b$ no longer participates in the problem, leaving a much simpler control problem over only $(\sigma,w)$.
    This is convenient for the analysis of symmetric perceptron models in \cite{bogpinprogress}.
    Similarly to Remark~\ref{r:control.problems.interpretation}, the \hyperlink{p:thm.control.problems.sym}{proof of Theorem~\ref{thm:control.problems.sym}} below also shows $\sM^{\Lip,\sym}(\alpha)$ coincides with several other classes 
    $\sM^{\star,\sym}(\alpha,0)$ defined in Definitions~\ref{d:SDE.conds} and \ref{d:SDE.classification.sym}--\ref{d:measure.classes.sym}.
\end{rmk}

Finally, we derive the following corollary, which establishes continuity of the sets  $\ocM^{\Ising}(\alpha)$ and $\ocM^{\Ising,\sym}(\alpha)$ in $\alpha$.
Recall that we equip $\cP_2(\bbR)$ and $\cP_2(\bbR_{\ge 0})$ with the $\bbW_2$ metric.
Let $P(\cP_2(\bbR))$ and $P(\cP_2(\bbR_{\ge 0}))$ denote the spaces of subsets of $\cP_2(\bbR)$ and $\cP_2(\bbR_{\ge 0})$, metrized by the Hausdorff distance
\beq
    \label{e:hausdorff}
    d_{\cH}(\sM,\sM')=\max\bigg\{
    \adjustlimits
    \sup_{\mu\in \sM} \inf_{\mu'\in \sM'} \bbW_2(\mu,\mu'),\,
    \adjustlimits
    \sup_{\mu'\in \sM'}
    \inf_{\mu\in \sM} \bbW_2(\mu,\mu')
    \bigg\}\,.
\eeq

\begin{cor}    \label{cor:feasible.distributions.continuity}
    The map $\alpha \mapsto \ocM^{\Ising}(\alpha)$
    from $(0,+\infty)$ to $P(\cP_2(\bbR))$ is continuous,
    as is the map $\alpha \mapsto  \ocM^{\Ising,\sym}(\alpha)$
    from $(0,+\infty)$ to $P(\cP_2(\bbR_{\ge 0}))$.
\end{cor}

To prove Theorem \ref{thm:control.problems.main}, we will introduce several classes of measures denoted $\sM^\star(\alpha, q_0)$, where $\star$ indicates the class type and $q_0 \in [0,1)$ is an a priori arbitrary starting time (see Definition \ref{d:measure.classes}).
These are the sets of attainable endpoint measures of a sequence of progressively cleaner control problems, which interpolate between the cases $\star = \IAMP,\BOGP$ and $\star = \sideal,\sidealast$.

We then prove a cycle of inclusions showing that the sets $\sM^\star(\alpha,0)$ coincide for all choices of $\star$.
This cycle also implies that each of these sets contains every $\sM^\star(\alpha, q_0)$ with $q_0 \in [0,1)$.

Most steps in this inclusion loop are established through approximation arguments: given a feasible control for one control problem in the loop, we show how to approximate it with a feasible control for the next problem so that the resulting endpoint measures are close in $\bbW_2$.

However, \textbf{several key inclusions rely crucially on the algorithmic interpretation of these stochastic control problems developed in the preceding sections}. 
For example, to show that the sets $\sM^\star(\alpha,0)$ contain $\sM^\star(\alpha,q_0)$ for all $q_0 \in [0,1)$, we use Theorem~\ref{thm:IAMP-main}\ref{i:IAMP-main-centered}, which implies that a general control process, with arbitrary $q_0$, can be simulated by a sequence of IAMP algorithms $\cA_N$ whose correlation functions satisfy $\chi_{\cA_N}(0)=0$.
Applying the results of \S\ref{s:sde}, we then extract control processes associated with these algorithms, which automatically have $q_0=0$. 
Thus all feasible endpoint measures can already be attained with $q_0=0$. We do not know a proof of this fact that works purely at the level of stochastic control.

We prove Theorem~\ref{thm:control.problems.sym} using a similar, though simpler, inclusion loop. 
This argument again depends on our algorithmic interpretation of the stochastic control problems. 
To show that $p \equiv 1$ suffices to attain all symmetrized endpoint measures, we use Theorem~\ref{thm:IAMP-main}\ref{i:IAMP-main-chaotic}, which implies that in the symmetric setting, any general control process can be simulated by a sequence of IAMP algorithms $\cA_N$ whose correlation functions satisfy $\chi_{\cA_N}(1-\epsilon) \le \epsilon$.
Extracting control processes associated with these algorithms using the results of \S\ref{s:sde} then yields functions $p \equiv \chi^{-1}$ that are close to $p \equiv 1$. 
Consequently, $p \equiv 1$ is sufficient to attain all symmetrized endpoint measures.

This section is structured as follows.

\begin{itemize}
    \item In \S\ref{ss:control.problems.measure.classes}, we introduce the control process classes and the corresponding measure classes $\sM^\star(\alpha, q_0)$ described above. The \hyperlink{p:thm.control.problems.main}{proof of Theorem~\ref{thm:control.problems.main}} appears at the end of this subsection and assumes most steps of the inclusion cycle, whose proofs are deferred to the remainder of the section.
    \item In \S\ref{ss:control.problems.Lip.approx}--\ref{ss:exact.endpoint.ising}, we prove all inclusions used in the proof of Theorem~\ref{thm:control.problems.main}. 
    These subsections are organized by the techniques used to show the inclusions: \S\ref{ss:control.problems.Lip.approx}--\ref{ss:control.problems.earlier.results} are based on algorithmic considerations, while \S\ref{ss:simple.adjustments}--\S\ref{ss:exact.endpoint.ising} use approximation arguments for the control processes. Specifically:
    \item In \S\ref{ss:control.problems.Lip.approx}, we prove the inclusions that follow from explicit perturbations to a Lipschitz algorithm $\cA_N$.
    \item In \S\ref{ss:control.problems.earlier.results}, we prove the inclusions that follow from our algorithmic results in \S\ref{s:sde}--\ref{sec:IAMP}.
    \item In \S\ref{ss:simple.adjustments} we prove the inclusions that follow from simple adjustments to the controls.
    \item In \S\ref{ss:continuity.in.budget} and \S\ref{ss:continuity.in.p}, we show that the feasible endpoint measures are suitably continuous in the budget and the function $p$, and prove the inclusions that follow from adjusting these quantities.
    \item In \S\ref{ss:exact.endpoint.ising}, we show how to approximate controls where the Ising process $Y(t)$ (see \eqref{e:formal.sde.ising}) satisfies $|Y(1)| \approx 1$ with controls where it satisfies $|Y(1)| = 1$ almost surely.
    \item Finally, in \S\ref{ss:sdes.symmetric.endpoints}, we present the \hyperlink{p:thm.control.problems.sym}{proof of Theorem~\ref{thm:control.problems.sym}}, which argues through a similar inclusion loop.
    Most steps of this loop are analogous to a step in the proof of Theorem~\ref{thm:control.problems.main}, and we provide proofs for those that are not.
\end{itemize}

\subsection{Definitions of measure classes}
\label{ss:control.problems.measure.classes}

We next define several collections of measures $\sM^\star(\alpha,q_0)$, for $q_0 \in [0,1)$.
These involve parameters $L,\epsilon,L_0,\epsilon_0 > 0$, which will eventually be sent to limits $L\to\infty$, $\epsilon \to 0$, $L_0 \to \infty$, $\epsilon_0 \to 0$ in that order.
In other words,
\beq
    \label{e:control.problems.limiting.regime}
    \frac{1}{L} \ll \epsilon \ll \frac{1}{L_0} \ll \epsilon_0 \ll 1\,.
\eeq

\begin{dfn}[controlled processes and endpoint measures]
    For $q_0 \in [0,1)$ and $q_* \in [q_0,q_0+\epsilon]$, define the following mutually independent random variables.
    Let $B,W$ be Brownian motions on the time interval $t\in [q_*,1]$ with $B(q_*) = W(q_*) = 0$. Let $U,U' \sim \unif([0,1])$.
    For measures $\zeta, \zeta^\Ising \in \cP_2(\R)$, let $X(q_*) \sim \zeta$ and $Y(q_*) \sim \zeta^\Ising$.
    Define the filtrations $\cF_X$ and $\cF_Y$ by
    \beq\label{e:filt.X.Y}
        \begin{aligned}
            \cF_X(t) &= \sigma\Big(
                (B(s): q_* \le s\le t),
                U,X(q_*)
            \Big)\,,\\
            \cF_Y(t) &= \sigma\Big(
                (W(s) : q_* \le s\le t),
                U',Y(q_*)
            \Big)\,. 
        \end{aligned}
    \eeq
Recall from \eqref{e:s.sqrt.fn} that $s(t) \equiv [(tp)'(t)]^{1/2}$.  Recall the notations of Definition~\ref{d:achievable-msrs-sym}.
    For $p \in \incr([q_0,1];[0,1])$
    (see Definition~\ref{d:incr.p}) and control processes $b,\sigma \in \ProgMsrbl(\cF_X)$, $w\in \ProgMsrbl(\cF_Y)$, define
    \begin{align}
        \label{e:formal.sde.main}
        X(t)
        &= X(q_*) + \int_{q_*}^t p'(u)^{1/2}\,b_u\,du
        +\int_{q_*}^t s(u)\sigma_u \,dB(u)\,,
        \\
        \label{e:formal.sde.ising}
        Y(t)
        &=Y(q_*) + \int_{q_*}^t w_u \,dW(u)\,,
    \end{align}
    with initial conditions $X(q_*) \sim \zeta$ and  $Y(q_*) \sim \zeta^\Ising$ as described above. Define
\beq
\label{e:control.problems.endpoint.measures} 
    \begin{aligned}
\mu(q_*,b,\sigma,p,\zeta) 
&= \Law(X(1))\,, \\
            \mu^{\Ising}(q_*,w,\zeta^{\Ising}) &= \Law(Y(1))\,.
    \end{aligned}\eeq
Note the similarity with
Definition~\ref{d:achievable-msrs}.
\end{dfn}

We stop to make a few comments:
\begin{itemize}
\item In the case that $b,\sigma$ take the form $b_t = b(t,X(t))$ and $\sigma_t = \sigma(t,X(t))$, such that the functions $p(t)^{1/2} b(t,x)$ and $s(t) \sigma(t,x)$ are Lipschitz in $(t,x)$, it is well known that the SDE \eqref{e:formal.sde.main} has a unique strong solution.
The analogous statement holds for the SDE \eqref{e:formal.sde.ising} when $w_t=w(t,Y(t))$ where $w$ is Lipschitz in time and space.  This case corresponds to condition~\ref{it:1.1bogp} in Definition~\ref{d:SDE.conds} below, and condition~\ref{it:1.1} if $p$ is also continuously differentiable.

\item  In the setting where $b,\sigma$ are progressively measurable with
respect to $\cF_X$, it is important to note that the definition \eqref{e:filt.X.Y} of $\cF_X(t)$ contains the information of the Brownian motion $B$ up to time $t$, but \textbf{not} the process $X$ up to time $t$. Thus, in this setting, \eqref{e:formal.sde.main} defines $X$ by \textbf{direct stochastic integration}, and is \textbf{not an SDE}. The analogous statement holds in the case that $w$ is progressively measurable with respect to $\cF_Y$. This case corresponds to condition~\ref{it:1.2} in Definition~\ref{d:SDE.conds} below.
\end{itemize}
Next define the cost and budget quantities (cf.\ \eqref{eq:is-budget-constraint}): 
\begin{align*}
	\Cost(t)
    &\equiv
    \bbE\bigg[(b_t)^2 + \frac{(tp)'(t)}{p(t)} \big(\sigma_t-1\big)^2
    	\bigg]\,,
    \\
    \budget(t;\alpha)
    &=
    \frac{(\bbE w_t)^2}
   {\alpha}\,.
\end{align*}
In this section we study different classes of solutions to \eqref{e:formal.sde.main}--\eqref{e:formal.sde.ising}, subject to
different
 variations of the $\Cost(t) \le\budget(t;\alpha)$ constraint
 and different technical conditions, as follows:

\begin{dfn}[conditions on control parameters]\label{d:SDE.conds}
For the control parameters $q_*,b,\sigma,w,p,\zeta,\zeta^\Ising$ introduced above, we define the following conditions:
\begin{enumerate}[(1)]
    \item Conditions on $b,\sigma,w$
    and symmetry in law of $Y$: 
    \label{it:1}
    \begin{enumerate}[(a)]
\item
\label{it:1.1}
$b_t=b(t,X(t))$,
$\sigma_t=\sigma(t,X(t))$, and
$w_t=w(t,Y(t))$ for $L$-Lipschitz functions $b:[q_*,1]\times \bbR\to[-L,L]$ and $\sigma,w:[q_*,1]\times \bbR\to[0,L]$. Moreover, $w(t,x)=w(t,-x)$ for all $t,x$ and $\zeta^\Ising$ is symmetric.
\item
        \label{it:1.1bogp}
$b_t=b(t,X(t))$,
$\sigma_t=\sigma(t,X(t))$,
and $w_t=w(t,Y(t))$ 
for functions $b:[q_*,1]\times \bbR\to[-L,L]$ and $\sigma,w:[q_*,1]\times \bbR\to[0,L]$, such that the functions $p'(t)^{1/2} b(t,x), s(t) \sigma(t,x), w(t,x)$, are $L$-Lipschitz in $(t,x)$.
        \item
        \label{it:1.2}
        $\sigma_t,b_t$ are progressively measurable processes with respect to $\cF_X$, and $w_t$ is progressively measurable with respect to $\cF_Y$, for filtrations $\cF_X,\cF_Y$ as defined by \eqref{e:filt.X.Y}.
        \item 
        \label{it:1.2even}
        Condition~\ref{it:1.2} holds, and moreover $\Law(Y(t))$ is even for all $t\in [q_*,1]$.
    \end{enumerate}
    \item Conditions on $\bbE[(w_t)^2]$:
    \label{it:2}
    \begin{enumerate}[(a)]
        \item
        \label{it:2.1}
        $\bbE[(w_t)^2]\in [1-\epsilon,1+\epsilon]$ for all $t\in [q_*,1]$.
        \item
        \label{it:2.2}
        $\bbE[(w_t)^2]=1$ for all $t\in [q_*,1]$.
    \end{enumerate}
    \item
    Conditions on $p \in \incr([q_0,1];[0,1])$ (see Definition~\ref{d:incr.p}): 
    \label{it:3}
    \begin{enumerate}[(a)]
        \item
        \label{it:3.1}
        $p$ is twice differentiable, with $1/L \le p(q_0)
        \le p(q_*) \le \epsilon$ and $\|p\|_{C^2([q_0,1])}\le L$.
        \item 
        \label{it:3.1concave}
        $p$ is concave and twice differentiable, with $1/L \le p(q_0) \le \epsilon$ and $\|p\|_{C^2([q_0,1])}\le L$. 
        \item
        \label{it:3.2}
        $p$ is concave with $p(q_0) = 0$, $p(q_*) \le \epsilon$, and $\|p\|_{C^1([q_0,1])}\le L$. (Note the condition $p(q_*) \le \epsilon$ is vacuous under condition~\ref{it:8.2} below.)
        \item
        \label{it:3.4}
        $p$ is concave and twice differentiable, with $1/L_0 \le p(q_0) \le \epsilon_0$ and $\|p\|_{C^2([q_0,1])}\le L_0$.
        \item
        \label{it:3.3} $p$ is concave and $p(q_0) = 0$. 
        \item
        \label{it:3.5withzero}
        $p(q_0) = 0$. 
        \item
        \label{it:3.5}
        no conditions (beyond the ones already imposed by the restriction $p\in\incr([q_0,1];[0,1])$).
        \item
        \label{it:3.6} $p(q_0) = 0$, $p(q_*) \le \epsilon$, $p(q_0 + \epsilon) \ge 1-\epsilon$, and $\|p\|_{C^1([q_0,1])}\le L$.
        \item
        \label{it:3.7} $p\equiv 1$.
    \end{enumerate}
    \item Budget constraint:
    \label{it:4}
    \begin{enumerate}[(a)]
        \item
        \label{it:4.1}
        $\int_{q_*}^1 (\Cost(t)-\budget(t;\alpha)))_+\de t\leq \epsilon$.
        \item
        \label{it:4.2}
        $\Cost(t)\leq \budget(t;\alpha)+\epsilon$ for all $t\in [q_*,1]$.
        \item
        \label{it:4.3}
        $\Cost(t)\leq \budget(t;\alpha)$ for all $t\in [q_*,1]$.
    \end{enumerate}
    \item
    \label{it:5}
    Conditions on endpoint $\Law(Y(1)) = \mu^{\Ising}(q_*,w,\zeta^{\Ising})$:
    \begin{enumerate}[(a)]
        \item
        \label{it:5.1}
        $\bbE[(|Y(1)|-1)^2]\leq \epsilon$.
        \item
        \label{it:5.2}
        $|Y(1)|=1$ almost surely.
    \end{enumerate}
    \item Conditions on $L^2$ norms of initializations $X(q_*)\sim\zeta$ and $Y(q_*)\sim\zeta^\Ising$ (with $\zeta,\zeta^\Ising\in\cP_2(\R)$): 
        \label{it:6}
    \begin{enumerate}[(a)]
        \item
        \label{it:6.1}
        $\bbE[X(q_*)^2] \le \epsilon$ and $|\bbE [Y(q_*)^2] - q_*| \le \epsilon$.
        \item
        \label{it:6.2}
        $X(q_*) = 0$ and $\bbE [Y(q_*)^2] = q_*$.
    \end{enumerate}
    \item Conditions on support of Ising initialization $Y(q_*)\sim\zeta^\Ising$:
    \label{it:7}
    \begin{enumerate}[(a)]
        \item \label{it:7.1}
        $\zeta^\Ising \in \cP_2(\R)$.
        \item \label{it:7.2}
        $\zeta^\Ising \in \cP([-1,1])$, the space of Borel probability measures on $[-1,1]$.
    \end{enumerate}
    \item Conditions on $q_*$:
    \label{it:8}
    \begin{enumerate}[(a)]
        \item \label{it:8.1}
        $q_* \in [q_0,q_0+\epsilon]$.
        \item \label{it:8.2}
        $q_* = q_0$.
    \end{enumerate}
\end{enumerate}
(We remark that conditions \ref{it:3.6} and \ref{it:3.7} do not appear in the \hyperlink{p:thm.control.problems.main}{proof of Theorem~\ref{thm:control.problems.main}}, but will be used in the \hyperlink{p:thm.control.problems.sym}{proof of Theorem~\ref{thm:control.problems.sym}} in \S\ref{ss:sdes.symmetric.endpoints}.)
\end{dfn}

Using the above conditions, we now define several classes of control parameters. We will show in the \hyperlink{p:thm.control.problems.main}{proof of Theorem~\ref{thm:control.problems.main}}
that the resulting sets of achievable laws of $X(1)$ are all equivalent in the limit \eqref{e:control.problems.limiting.regime}.

\colorlet{colA}{orange!80!yellow}
\colorlet{colZ}{orange!70!yellow}
\colorlet{colBb}{orange!60!yellow}
\colorlet{colB}{orange!50!yellow}
\colorlet{colBa}{orange!40!yellow}
\colorlet{colD}{orange!20!yellow}
\colorlet{colE}{blue!40!white}
\colorlet{colF}{blue!20!white}

\newcommand{\MSRheader}{ & coefs & $\E[(w_t)^2]$ & $p$ & budget & endpt & init & support & {{$q_*$}}}

\DeclareDocumentCommand{\IAMPdefn}{ g }{$(\alpha,
	\IfValueF{#1}{q_0}%
	\IfValueT{#1}{0},
	L,\epsilon)$-$\IAMP$
    & \cellcolor{colD}\ref{it:1.1} 
    & \cellcolor{colD}\ref{it:2.1}  
    & \cellcolor{colD}\ref{it:3.1} 
    & \cellcolor{colB}\ref{it:4.2} 
    & \cellcolor{colD}\ref{it:5.1} 
    & \cellcolor{colD}\ref{it:6.1}
    & \cellcolor{colD}\ref{it:7.1}
    & \cellcolor{colD}\ref{it:8.1}}
    % previously was 8b!

\DeclareDocumentCommand{\BOGPdefn}{ g }{$(\alpha,
	\IfValueF{#1}{q_0}%
	\IfValueT{#1}{0},
	L,\epsilon)$-$\BOGP$
    & \cellcolor{colBa}\ref{it:1.1bogp}
    & \cellcolor{colD}\ref{it:2.1} 
    & 
    \IfValueF{#1}{\cellcolor{colB}\ref{it:3.2}}%
    \IfValueT{#1}{\cellcolor{colE}\ref{it:3.6}}
    & \cellcolor{colD}\ref{it:4.1} 
    & \cellcolor{colD}\ref{it:5.1} 
    & \cellcolor{colD}\ref{it:6.1} 
    & \cellcolor{colD}\ref{it:7.1} 
    & \cellcolor{colD}\ref{it:8.1}}

\DeclareDocumentCommand{\FAIRdefn}{ g }{$(\alpha,
	\IfValueF{#1}{q_0}%
	\IfValueT{#1}{0},
	L,\epsilon)$-$\si$
    & \cellcolor{colA}\ref{it:1.2} 
    & \cellcolor{colD}\ref{it:2.1} 
    & \IfValueF{#1}{\cellcolor{colB}\ref{it:3.2}}%
    \IfValueT{#1}{\cellcolor{colE}\ref{it:3.6}}
    & \cellcolor{colD}\ref{it:4.1} 
    & \cellcolor{colD}\ref{it:5.1} 
    & \cellcolor{colD}\ref{it:6.1} 
    & \cellcolor{colD}\ref{it:7.1} 
    & \cellcolor{colD}\ref{it:8.1}}
    
\DeclareDocumentCommand{\TMPONEdefn}{ g }{$(\alpha,
	\IfValueF{#1}{q_0}%
	\IfValueT{#1}{0},
	L,\epsilon)$-$\sii$
    & \cellcolor{colA}\ref{it:1.2} 
    & \cellcolor{colD}\ref{it:2.1} 
    & \IfValueF{#1}{\cellcolor{colB}\ref{it:3.2}}%
    \IfValueT{#1}{\cellcolor{colE}\ref{it:3.6}}
    & \cellcolor{colD}\ref{it:4.1} 
    & \cellcolor{colD}\ref{it:5.1} 
    & \cellcolor{colD}\ref{it:6.1} 
    & \cellcolor{colA}\ref{it:7.2} 
    & \cellcolor{colD}\ref{it:8.1}}
    
\DeclareDocumentCommand{\TMPTWOdefn}{ g }{$(\alpha,
	\IfValueF{#1}{q_0}%
	\IfValueT{#1}{0},
	L,\epsilon)$-$\siii$
    & \cellcolor{colA}\ref{it:1.2} 
    & \cellcolor{colD}\ref{it:2.1} 
    & \IfValueF{#1}{\cellcolor{colB}\ref{it:3.2}}%
    \IfValueT{#1}{\cellcolor{colE}\ref{it:3.6}}
    & \cellcolor{colD}\ref{it:4.1} 
    & \cellcolor{colD}\ref{it:5.1} 
    & \cellcolor{colD}\ref{it:6.1} 
    & \cellcolor{colA}\ref{it:7.2} 
    & \cellcolor{colA}\ref{it:8.2}}

\DeclareDocumentCommand{\GOODdefn}{ g }{$(\alpha,
	\IfValueF{#1}{q_0}%
	\IfValueT{#1}{0},
	L,\epsilon)$-$\siv$
    & \cellcolor{colA}\ref{it:1.2} 
    & \cellcolor{colD}\ref{it:2.1} 
    & \IfValueF{#1}{\cellcolor{colB}\ref{it:3.2}}%
    \IfValueT{#1}{\cellcolor{colE}\ref{it:3.6}}
    & \cellcolor{colD}\ref{it:4.1} 
    & \cellcolor{colD}\ref{it:5.1} 
    & \cellcolor{colA}\ref{it:6.2} 
    & \cellcolor{colA}\ref{it:7.2} 
    & \cellcolor{colA}\ref{it:8.2}}

\DeclareDocumentCommand{\TMPTHRdefn}{ g }{$(\alpha,
	\IfValueF{#1}{q_0}%
	\IfValueT{#1}{0},
	L,\epsilon)$-$\sv$
    & \cellcolor{colA}\ref{it:1.2} 
    & \cellcolor{colA}\ref{it:2.2} 
    & \IfValueF{#1}{\cellcolor{colB}\ref{it:3.2}}%
    \IfValueT{#1}{\cellcolor{colE}\ref{it:3.6}}
    & \cellcolor{colD}\ref{it:4.1} 
    & \cellcolor{colD}\ref{it:5.1} 
    & \cellcolor{colA}\ref{it:6.2} 
    & \cellcolor{colA}\ref{it:7.2} 
    & \cellcolor{colA}\ref{it:8.2}}

\DeclareDocumentCommand{\GREATdefn}{ g }{$(\alpha,
	\IfValueF{#1}{q_0}%
	\IfValueT{#1}{0},
	L,\epsilon)$-$\svi$ % (``great'')
    & \cellcolor{colA}\ref{it:1.2} 
    & \cellcolor{colA}\ref{it:2.2} 
    & \IfValueF{#1}{\cellcolor{colB}\ref{it:3.2}}%
    \IfValueT{#1}{\cellcolor{colE}\ref{it:3.6}}
    & \cellcolor{colA}\ref{it:4.3} 
    & \cellcolor{colD}\ref{it:5.1} 
    & \cellcolor{colA}\ref{it:6.2} 
    & \cellcolor{colA}\ref{it:7.2} 
    & \cellcolor{colA}\ref{it:8.2}}

\DeclareDocumentCommand{\SPRdefn}{ g }{$(\alpha,
	\IfValueF{#1}{q_0}%
	\IfValueT{#1}{0},
	\epsilon
	\IfValueF{#1}{,L_0,\epsilon_0}
	)$-$\svii$ % (``superior'')
    & \cellcolor{colA}\ref{it:1.2} 
    & \cellcolor{colA}\ref{it:2.2} 
    & \IfValueF{#1}{\cellcolor{colBb}\ref{it:3.4}}%
    \IfValueT{#1}{\cellcolor{colF}\ref{it:3.7}}
    & \cellcolor{colA}\ref{it:4.3} 
    & \cellcolor{colD}\ref{it:5.1} 
    & \cellcolor{colA}\ref{it:6.2} 
    & \cellcolor{colA}\ref{it:7.2} 
    & \cellcolor{colA}\ref{it:8.2}}

\newcommand{\PRFONEdefn}{$(\alpha,q_0,L_0,\epsilon_0)$-$\sviii$
    & \cellcolor{colA}\ref{it:1.2} 
    & \cellcolor{colA}\ref{it:2.2} 
    & \cellcolor{colBb}\ref{it:3.4} 
    & \cellcolor{colA}\ref{it:4.3} 
    & \cellcolor{colA}\ref{it:5.2} 
    & \cellcolor{colA}\ref{it:6.2} 
    & \cellcolor{colA}\ref{it:7.2} 
    & \cellcolor{colA}\ref{it:8.2}}

\newcommand{\PRFONEHALFdefn}{$(\alpha,q_0,L,\epsilon)$-$\six$
    & \cellcolor{colA}\ref{it:1.2} 
    & \cellcolor{colA}\ref{it:2.2} 
    & \cellcolor{colBa}\ref{it:3.1concave} 
    & \cellcolor{colA}\ref{it:4.3} 
    & \cellcolor{colA}\ref{it:5.2} 
    & \cellcolor{colA}\ref{it:6.2} 
    & \cellcolor{colA}\ref{it:7.2} 
    & \cellcolor{colA}\ref{it:8.2}}
    
\DeclareDocumentCommand{\PRFTWOdefn}{ g }{$(\alpha,
	\IfValueF{#1}{q_0}%
	\IfValueT{#1}{0},
	L,\epsilon)$-$\sx$
    & \cellcolor{colA}\ref{it:1.2} 
    & \cellcolor{colA}\ref{it:2.2} 
    & \cellcolor{colD}\ref{it:3.1} 
    & \cellcolor{colA}\ref{it:4.3} 
    & \cellcolor{colA}\ref{it:5.2} 
    & \cellcolor{colA}\ref{it:6.2} 
    & \cellcolor{colA}\ref{it:7.2} 
    & \cellcolor{colA}\ref{it:8.2}}

\DeclareDocumentCommand{\PRFTHRdefn}{ g }{$(\alpha,
	\IfValueF{#1}{q_0}%
	\IfValueT{#1}{0},
	L,\epsilon)$-$\sxi$
    & \cellcolor{colBb}\ref{it:1.2even} 
    & \cellcolor{colA}\ref{it:2.2} 
    & \cellcolor{colD}\ref{it:3.1} 
    & \cellcolor{colA}\ref{it:4.3} 
    & \cellcolor{colA}\ref{it:5.2} 
    & \cellcolor{colA}\ref{it:6.2} 
    & \cellcolor{colA}\ref{it:7.2} 
    & \cellcolor{colA}\ref{it:8.2}}

\newcommand{\PRFFOURdefn}{$(\alpha,q_0)$-$\sxii$
    & \cellcolor{colA}\ref{it:1.2} 
    & \cellcolor{colA}\ref{it:2.2} 
    & \cellcolor{colZ}\ref{it:3.5withzero} 
    & \cellcolor{colA}\ref{it:4.3} 
    & \cellcolor{colA}\ref{it:5.2} 
    & \cellcolor{colA}\ref{it:6.2} 
    & \cellcolor{colA}\ref{it:7.2} 
    & \cellcolor{colA}\ref{it:8.2}}

\DeclareDocumentCommand{\SUPERPERFECTdefn}{ g }{$(\alpha,
	\IfValueF{#1}{q_0}%
	\IfValueT{#1}{0})$-$\sideal$
    & \cellcolor{colA}\ref{it:1.2} 
    & \cellcolor{colA}\ref{it:2.2} 
    & \IfValueF{#1}{\cellcolor{colA}\ref{it:3.5}}%
     \IfValueT{#1}{\cellcolor{colF}\ref{it:3.7}}
    & \cellcolor{colA}\ref{it:4.3} 
    & \cellcolor{colA}\ref{it:5.2} 
    & \cellcolor{colA}\ref{it:6.2} 
    & \cellcolor{colA}\ref{it:7.2} 
    & \cellcolor{colA}\ref{it:8.2}}

\newcommand{\SUPERPERFECTASTdefn}{$(\alpha,q_0)$-$\sidealast$
    & \cellcolor{colA}\ref{it:1.2} 
    & \cellcolor{colA}\ref{it:2.2} 
    & \cellcolor{colA}\ref{it:3.3} 
    & \cellcolor{colA}\ref{it:4.3} 
    & \cellcolor{colA}\ref{it:5.2} 
    & \cellcolor{colA}\ref{it:6.2} 
    & \cellcolor{colA}\ref{it:7.2}
    & \cellcolor{colA}\ref{it:8.2}}

\begin{dfn}[classification]
\label{d:SDE.classification}
We classify the parameters $(q_*,b,\sigma,w,p,\zeta,\zeta^\Ising)$ as follows:
\begin{center}
\begin{tabular}[h]{r|cccccccc}
\MSRheader\\
\hline
\IAMPdefn \\
\BOGPdefn \\
\FAIRdefn \\
\TMPONEdefn \\
\TMPTWOdefn \\
\GOODdefn \\
\TMPTHRdefn \\
\GREATdefn \\
\SPRdefn \\
\PRFONEdefn \\
\PRFONEHALFdefn \\
\PRFTWOdefn \\
\PRFTHRdefn \\
\PRFFOURdefn \\
\SUPERPERFECTdefn \\
\SUPERPERFECTASTdefn
\end{tabular}\end{center}
For example, we say the set of parameters $(q_*,b,\sigma,w,p,\zeta,\zeta^\Ising)$ is \emph{$(\alpha,q_0,L,\epsilon)$-$\si$} if it satisfies conditions \ref{it:1.2}, \ref{it:2.1}, \ref{it:3.2}, \ref{it:4.1}, \ref{it:5.1}, \ref{it:6.1}, \ref{it:7.1}, \ref{it:8.1}. For
\beq
    \label{e:star.options}
    \star \in \bigg\{
    \begin{array}{c}
    \IAMP,\BOGP,\si,
    \sii,\siii,\siv,
    \sv,\svi,\svii,\\
    \sviii,
    \six,
    \sx,
    \sxi,
    \sxii,
    \sideal,
    \sidealast
    \end{array}\bigg\}\,,
\eeq
we shall denote
	\[\Adm^\star(\alpha,q_0;L,\epsilon,L_0,\epsilon_0)\]
for the set of all $(q_*,b,\sigma,w,p,\zeta,\zeta^\Ising)$ satisfying conditions $\star$. Note that some of the parameters $L,\epsilon,L_0, \epsilon_0$ do not appear in the definitions of some of the above classes: in such cases, we take $\Adm^\star(\alpha,q_0;L,\epsilon,L_0,\epsilon_0)$ to be the same for all values of those parameters, and sometimes omit those parameters from the notation. 
For example,
	\[\Adm^{\BOGP}(\alpha,q_0;L,\epsilon,L_0,\epsilon_0) = \Adm^{\BOGP}(\alpha,q_0;L,\epsilon)\]
is the same for all $L_0,\epsilon_0$. 
\end{dfn}
Recall the point-to-set $\bbW_2$ distance defined in \eqref{def:W2-to-set}.
For $\cM \subseteq \cP_2(\bbR)$ and $\iota > 0$, define the $\bbW_2$-ball
\beq\label{e:W2.ball.notation}
    \cB_\epsilon(\cM) = \Big\{
        \mu \in \cP_2(\bbR) : 
        \bbW_2(\mu,\cM) \le \epsilon
    \Big\}.
\eeq 
We now use the classification of Definition~\ref{d:SDE.classification} to define corresponding sets of achievable endpoint measures, as follows:

\begin{dfn}[measure classes]
    \label{d:measure.classes}
Recall the definition of the endpoint measure $\mu(q_*,b,\sigma,p,\zeta)$ from \eqref{e:control.problems.endpoint.measures}.
    For $\star$ as in \eqref{e:star.options}, we define
    \[
        \cM^\star(\alpha,q_0;L,\epsilon,L_0,\epsilon_0) 
        \equiv \bigg\{
            \mu(q_*,b,\sigma,p,\zeta) : (q_*,b,\sigma,w,p,\zeta,\zeta^\Ising) \in \Adm^\star(\alpha,q_0;L,\epsilon,L_0,\epsilon_0)
        \bigg\}\,.
    \]
When $\Adm^\star(\alpha,q_0;L,\epsilon,L_0,\epsilon_0)$ does not depend on some of these parameters we will sometimes omit them from the notation, writing e.g. $\cM^\BOGP(\alpha,q_0;L,\epsilon) = \cM^\BOGP(\alpha,q_0;L,\epsilon,L_0,\epsilon_0)$ for all $L_0,\epsilon_0$. Define
    \[
        \sM^\star(\alpha,q_0) 
        = \adjustlimits 
        \bigcap_{\epsilon_0 > 0} 
        \bigcup_{L_0>0}
        \bigcap_{\epsilon>0}
        \bigcup_{L>0}
        \cB_{\epsilon_0}(\cM^\star(\alpha,q_0;L,\epsilon,L_0,\epsilon_0))\,.
    \]
Note that this order of quantifiers corresponds to the limit \eqref{e:control.problems.limiting.regime}. 
\end{dfn}

\begin{rmk}
    \label{r:measure.classes.simplification}
For all values of $\star$ except $\star \in \{\svii,\sviii\}$, the class $\Adm^\star(\alpha,q_0;L,\epsilon,L_0,\epsilon_0)$ does not depend on $L_0,\epsilon_0$.
    For such $\star$, we have the simplification
    \beq
        \label{e:measure.classes.simplification.without.L0}
        \sM^\star(\alpha,q_0)
        = \adjustlimits
        \bigcap_{\epsilon_0 > 0}
        \bigcap_{\epsilon>0}
        \bigcup_{L>0}
        \cB_{\epsilon_0}(\cM^\star(\alpha,q_0;L,\epsilon))
        = \adjustlimits
        \bigcap_{\epsilon>0}
        \bigcup_{L>0}
        \cB_{\epsilon}(\cM^\star(\alpha,q_0;L,\epsilon))\,,
    \eeq
    since $\cB_{\epsilon_0}(\cM^\star(\alpha,q_0;L,\epsilon))$ is decreasing in both $\epsilon$ and $\epsilon_0$.
    In the cases $\star \in \{\sxii,\sideal,\sidealast\}$, the parameters $L,\epsilon$ are also not relevant, so
    \[
        \sM^\star(\alpha,q_0)
        = \bigcap_{\epsilon>0}
        \cB_{\epsilon}(\cM^\star(\alpha,q_0))
        = \overline{\cM^\star(\alpha,q_0)}\,,
    \]
is the closure of $\cM^\star(\alpha,q_0)$ with respect to the $\bbW_2$ metric. In particular, we have
    \begin{align*}
    \sM^{\sideal}(\alpha,0) 
    &= \ocM^{\Ising}(\alpha)\,,\\
    \sM^{\sidealast}(\alpha,0) 
    &= \ocM^{\Ising, \concave}(\alpha)\,,
    \end{align*}
for $\ocM^{\Ising}(\alpha)$ and $\ocM^{\Ising, \concave}(\alpha)$ specified in Definition~\ref{d:achievable-msrs}.
\end{rmk}

\begin{rmk}\label{r:true.mg}
Under conditions~\ref{it:1.1} or \ref{it:1.1bogp},
since $\sigma$ and $w$ are bounded, it is immediate that the process $Y$ defined by \eqref{e:formal.sde.ising} is a true martingale (not only a local martingale), and that the local martingale term in \eqref{e:formal.sde.main} is likewise a true martingale. We note that the same holds for any of the classes in Definition~\ref{d:SDE.classification}, as a consequence of basic principles of stochastic calculus: as explained for example in \cite[Ch.~5]{MR3497465}, if $H$ is a progressive process and $M$ is a continuous local martingale with quadratic variation $\langle M\rangle$, such that we have the condition
    \[
    \E\int_0^\infty (H_t)^2\,d\langle M\rangle_t<\infty\,,
    \]
then the stochastic integral of $H$ with respect to $M$ is in fact an $L^2$-bounded true martingale. Conditions~\ref{it:2.1} or \ref{it:2.2} clearly guarantee
    \[\E\int_{q_*}^1 (w_t)^2\,d\langle W\rangle_t
    = \int_{q_*}^1 \E[(w_t)^2] \,dt \le 1+\epsilon <\infty\,,\]
so that the process $Y$ from \eqref{e:formal.sde.ising} is an $L^2$-bounded true martingale. Combining \ref{it:2.1} or \ref{it:2.2} with any of the conditions \ref{it:4.1}--\ref{it:4.3} gives
    \begin{align}\nonumber
    &\E\int_{q_*}^1 s(t)^2 (\sigma_t)^2\,d\langle B\rangle_t
    \le2 \E \int_{q_*}^1 s(t)^2 \Big( (\sigma_t-1)^2 + 1\Big)\,dt
    \le 2\bigg\{\E \int_{q_*}^1 \frac{(tp)'(t)}{p(t)}(\sigma_t-1)^2\,dt
     +1\bigg\} \\
    &\le 2\bigg\{\int_{q_*}^1 \Cost(t)\,dt + 1\bigg\}
    \le 2\bigg\{\int_{q_*}^1 \budget(t;\alpha)\,dt +\epsilon+ 1\bigg\}
    \le 2\bigg\{\frac{1+\epsilon}{\alpha} +\epsilon+ 1\bigg\}<\infty\,,
    \label{e:loc.mg.L2.bounded}
    \end{align}
so the local martingale term in \eqref{e:formal.sde.main} is a true martingale, as claimed.
\end{rmk}
We note the following simple estimate, which implies the sets $\sM^\star(\alpha,q_0)$ are $\bbW_2$-bounded.
\begin{lem}
    \label{l:L2-bound-on-stochastic-control}
    For any $\star$ as in \eqref{e:star.options} and $\mu \in \sM^\star(\alpha,q_0)$, we have
        \[\|\mu\|_{L^2}^2 \le
    \frac{10}{\alpha}
    +6 + O(\epsilon)\,.\]
\begin{proof}
Let $X$ be given by \eqref{e:formal.sde.main}. Let us write
$X(1)= X(q_*)+\textup{(a)}+\textup{(b)}$ where (a) denotes the drift term and (b) denotes the martingale term.
We have $\bbE[X(q_*)^2] \le \epsilon$ under either of the initialization conditions \ref{it:6.1} or \ref{it:6.2}. Next, using the Cauchy--Schwarz inequality, we have
    \begin{align*}
    \E[\textup{(a)}^2]
    &= \E\bigg[
    \bigg(\int_{q_*}^1
    p'(t)^{1/2} b_t\,dt\bigg)^2
    \bigg]
    \le \E\bigg[
    \bigg(\int_{q_*}^1 p'(t)\,dt\bigg)
    \bigg(\int_{q_*}^1 (b_t)^2 \,dt
        \bigg)\bigg]
    \le
    \int_{q_*}^1 \E[(b_t)^2] \,dt\\
    &\le \int_{q_*}^1 \Cost(t)\,dt
    \le \int_{q_*}^1 \budget(t;\alpha)\,dt + \epsilon
    \le \frac{1+\epsilon}{\alpha} + \epsilon\,,
    \end{align*}
where the last two inequalities hold under any of the budget conditions
\ref{it:4.1}, \ref{it:4.2}, or \ref{it:4.3}
combined with either condition
\ref{it:2.1} or \ref{it:2.2}. Lastly, we have from \eqref{e:loc.mg.L2.bounded} that 
    \[
    \E[\textup{(b)}^2]
    \le \int_{q_*}^1 s(t)^2
    \E[(\sigma_t)^2] \,dt
    \le 2\bigg\{ \frac{1+\epsilon}{\alpha}
        + \epsilon + 1\bigg\}\,.\]
Combining the above bounds gives altogether
    \[
    \E[X(1)^2]
    \le 3\bigg\{
    \E[X(q_*)^2]
    +\E[\textup{(a)}^2]
    +\E[\textup{(b)}^2]\bigg\}
    \le 3\bigg\{
    \frac{3(1+\epsilon)}{\alpha}
    +2 + 4\epsilon\bigg\}\,,
    \]
proving the claim.
\end{proof}
\end{lem}
We next give the \hyperlink{p:thm.control.problems.main}{proof of Theorem~\ref{thm:control.problems.main}}, assuming the results of \S\ref{ss:control.problems.Lip.approx}--\ref{ss:exact.endpoint.ising}.
Recall $\chi_{\cA_N}$ as defined in \eqref{e:p.corr.overlap}.
For $q_0 \in [0,1]$, define
\begin{align}
    \label{e:cM.Lip.q0}
    \cM^{\Lip}(\alpha,q_0;L,\epsilon) 
    &\equiv
    \bigg\{
        \begin{array}{ll}
        \mu \in \cP_2(\bbR):
        & \text{there exists a sequence of $L$-Lipschitz $(\cA_{N_j})_{j\ge 1}$ with} \\
        &\WERR(\cA_{N_j},\mu)
         \le \epsilon \text{ and }
         \lim_{j\to\infty} \chi_{\cA_{N_j}}(0) = q_0
        \end{array}
    \bigg\}\\    \label{e:sM.Lip.q0}
    \sM^{\Lip}(\alpha,q_0) 
    &\equiv \adjustlimits
     \bigcap_{\epsilon > 0} \bigcup_{L > 0} \cM^{\Lip}(\alpha,q_0;L,\epsilon)\,.
\end{align}
Compare this with our earlier definitions \eqref{e:cM.Lip} and \eqref{e:sM.Lip}.

\begin{proof}[\hypertarget{p:thm.control.problems.main}{Proof of Theorem~\ref{thm:control.problems.main}}, assuming \S\ref{ss:control.problems.Lip.approx}--\ref{ss:exact.endpoint.ising}] For all $q_0 \in [0,1)$, we have the inclusions:
    \begin{itemize}
        \item $\sM^{\Lip}(\alpha,q_0) \subseteq \sM^{\BOGP}(\alpha,q_0)$ (Lemma~\ref{l:Lip.in.BOGP});
        \item $\sM^{\BOGP}(\alpha,q_0) \subseteq \sM^{\si}(\alpha,q_0)$ (trivial);
        \item $\sM^{\si}(\alpha,q_0) \subseteq \sM^{\sii}(\alpha,q_0)$ (Lemma~\ref{l:fair.in.tmp1});
        \item $\sM^{\sii}(\alpha,q_0) \subseteq \sM^{\siii}(\alpha,q_0)$ (Lemma~\ref{l:tmp1.in.tmp2});
        \item $\sM^{\siii}(\alpha,q_0) \subseteq \sM^{\siv}(\alpha,q_0)$ (Lemma~\ref{l:tmp2.in.good});
        \item $\sM^{\siv}(\alpha,q_0) \subseteq \sM^{\sv}(\alpha,q_0)$ (Lemma~\ref{l:good.in.tmp3});
        \item $\sM^{\sv}(\alpha,q_0) \subseteq \sM^{\svi}(\alpha,q_0)$ (Proposition~\ref{p:continuity.in.auxilliary.budget});
        \item $\sM^{\svi}(\alpha,q_0) \subseteq \sM^{\svii}(\alpha,q_0)$ (Lemma~\ref{l:great.in.spr.prf4.in.prf2});
        \item $\sM^{\svii}(\alpha,q_0) \subseteq \sM^{\sviii}(\alpha,q_0)$ (Lemma~\ref{l:spr.in.prf1});
        \item $\sM^{\sviii}(\alpha,q_0) = \sM^{\six}(\alpha,q_0)$ (trivial: rename $(L_0,\epsilon_0)$ to $(L,\epsilon)$ and recall \eqref{e:measure.classes.simplification.without.L0});
        \item $\sM^{\six}(\alpha,q_0) \subseteq \sM^{\sx}(\alpha,q_0)$ (trivial);
        \item $\sM^{\sx}(\alpha,q_0) \subseteq \sM^{\sxi}(\alpha,q_0)$ (Lemma~\ref{l:prf2.in.prf3});
        \item $\sM^{\sxi}(\alpha,q_0) \subseteq \sM^{\IAMP}(\alpha,q_0)$ (Lemma~\ref{l:prf3.in.IAMP});
        \item $\sM^{\IAMP}(\alpha,q_0) \subseteq \sM^{\Lip}(\alpha,0)$ (Lemma~\ref{l:IAMP.in.Lip}).
    \end{itemize}
    If we set $q_0 = 0$, we deduce that for any $\star$ in \eqref{e:star.options} except $\sxii$, $\sideal$, and $\sidealast$,
    \beq
        \label{e:sM.star.coincide}
        \sM^\star(\alpha,0) = \sM^{\Lip}(\alpha,0)\,,
    \eeq
    and thus all of these $\sM^\star(\alpha,0)$ coincide.
    Furthermore, we have the inclusions and equalities:
    \begin{itemize}
        \item $\sM^{\sx}(\alpha,0) = \sM^{\sxii}(\alpha,0)$ (Lemmas~\ref{l:great.in.spr.prf4.in.prf2} and \ref{l:prf2.in.prf4.prf15.in.idealast});
        \item $\sM^{\sxii}(\alpha,0) \subseteq \sM^{\sideal}(\alpha,0)$ (trivial); 
        \item $\sM^{\sideal}(\alpha,0) \subseteq \sM^{\sxii}(\alpha,0)$ (Proposition~\ref{p:continuity.in.p.near.0}).
    \end{itemize}
    It follows that \eqref{e:sM.star.coincide} holds for $\star \in \{\sxii, \sideal\}$ as well.
    We also have the inclusions:
    \begin{itemize}
        \item $\sM^{\six}(\alpha,0) \subseteq \sM^{\sidealast}(\alpha,0)$ (Lemma~\ref{l:prf2.in.prf4.prf15.in.idealast});
        \item $\sM^{\sidealast}(\alpha,0) \subseteq \sM^{\sideal}(\alpha,0)$ (trivial).
    \end{itemize}
    Thus \eqref{e:sM.star.coincide} holds for $\star = \sidealast$ as well.    
    For general $q_0 \in [0,1)$, the above inclusions also show $\sM^{\Lip}(\alpha,q_0) \subseteq \sM^{\Lip}(\alpha,0)$, and thus
    \[
        \bigcup_{q_0 \in [0,1)}
        \sM^{\Lip}(\alpha,q_0)
        = \sM^{\Lip}(\alpha,0)
        = \sM^{\sideal}(\alpha,0)
        = \sM^{\sidealast}(\alpha,0)\,.
    \]
    Finally, we have:
    \begin{itemize}
        \item $\sM^{\Lip}(\alpha) = \bigcup_{q_0 \in [0,1]} \sM^{\Lip}(\alpha,q_0)$ (Lemma~\ref{l:sM.Lip.union.q0}); 
        \item $\sM^{\Lip}(\alpha,1) = \{\cN(0,1)\}$ (Lemma~\ref{l:sM.Lip.1.trivial}); 
        \item $\cN(0,1) \in \sM^{\sideal}(\alpha,0)$ (Lemma~\ref{l:ideal.contains.std.gaussian}).
    \end{itemize}
    Thus
    \begin{align*}
        \sM^{\Lip}(\alpha) &= \bigcup_{q_0 \in [0,1)} \sM^{\Lip}(\alpha,q_0) \cup \sM^{\Lip}(\alpha,1) \\
        &= \sM^{\sideal}(\alpha,0) \cup \{\cN(0,1)\}
        = \sM^{\sideal}(\alpha,0)
        = \sM^{\sidealast}(\alpha,0)\,. 
    \end{align*}
    As observed in Remark~\ref{r:measure.classes.simplification}, $\sM^{\sideal}(\alpha,0) = \ocM^{\Ising}(\alpha)$ and $\sM^{\sidealast}(\alpha,0) = \ocM^{\Ising, \concave}(\alpha)$.
    The result follows.
\end{proof}

\subsection{Approximations of Lipschitz algorithms}
\label{ss:control.problems.Lip.approx}

In this subsection, we use simple approximation arguments to show the following equalities.
\begin{itemize}
    \item Lemma~\ref{l:sM.Lip.union.q0} shows that $\sM^{\Lip}(\alpha)$ agrees with the union of $\sM^{\Lip}(\alpha,q_0)$ over $q_0\in[0,1]$.
    \item Lemma~\ref{l:sM.Lip.1.trivial} shows $\sM^{\Lip}(\alpha,1) = \{\cN(0,1)\}$.
\end{itemize}
The following lemma will be useful in the approximation arguments below.

\begin{lem}
\label{l:control.problems.approx}
Let $\cA_N$ and $\tilde\cA_N$ be $L$-Lipschitz algorithms with 
	\beq\label{e:ctrl.prbs.appx.assum}
	\frac{1}{N} \bbE\bigg[
		\Big\|\cA(\bG,\bg^\aux)-\tilde\cA(\bG,\bg^\aux)
		\Big\|_2^2
		\bigg] 
		\le \delta\,.
	\eeq
    There exist $C_1 = C_1(\alpha)$, $C_2 = C_2(\alpha,L)$, and an absolute constant $c_0$ such that 
    \[\WERR(\cA_N,\tilde{\cA}_N)
    \equiv
        \max\bigg\{
        \bbW_2\Big(\mu(\cA_N),
        \mu(\tilde \cA_N)
        	\Big)^2,
        \bbW_2\Big(
        	\mu^{\Ising}(\cA_N),
        \mu^{\Ising}(\tilde \cA_N)
        \Big)^2
        \bigg\}
        \le C_1 \delta + C_2e^{-c_0N}\,.
    \]
\end{lem}
\begin{proof}
    Recall from \eqref{e:mu.of.Alg} that, with the notation $\EmpDist(\cdot)$ defined in \eqref{eq:proj-pursuit-def},
    \[
        \mu(\cA_N) = \bbE \bigg[\EmpDist
            \bigg( \frac{\bG\cA(\bG,\bg^{\aux})}{N^{1/2}}\bigg)
        \bigg]\,,
    \]
Throughout this proof, we abbreviate $\bsig = \cA_N(\bG,\bg^{\aux})$ and $\tilde\bsig = \tilde \cA_N(\bG,\bg^{\aux})$. By convexity of the squared $\bbW_2$ distance, we can bound
    \begin{align*}
        \bbW_2\Big(
        \mu(\cA_N),\mu(\tilde \cA_N)
        \Big)^2
        &\le \bbE \bigg[\bbW_2\bigg(
            \EmpDist \bigg( \frac{\bG \bsig}{N^{1/2}}\bigg), 
            \EmpDist \bigg( \frac{\bG \tilde\bsig}{N^{1/2}}\bigg)
        \bigg)^2\bigg] \\
        &\le \bbE \bigg[\frac{1}{M} \bigg\|\frac{\bG (\bsig - \tilde\bsig)}{N^{1/2}}\bigg\|_2^2\bigg]
        \le \frac{\bbE[\|\bG\|_{\op}^2 \|\bsig - \tilde\bsig\|_2^2]}{MN}\,.
    \end{align*}
By Lemma~\ref{l:wishart}, there exists a finite constant $C_1 = C_1(\alpha)$ and a positive absolute constant $c_0$ such that 
the event $\cE \equiv \{\|\bG\|_{\op} \le C_1 N^{1/2}\}$ occurs with probability at least $1-\exp(-Nc_0)$.
In what follows we abbreviate $\E[A;\cE]\equiv \E[A\ind\{\cE\}]$. Then 
	\[
    \textup{(a)}
    \equiv
    \frac{\bbE[
	\|\bG\|_{\op}^2
	\|\bsig - \tilde\bsig\|_2^2
	;\cE]}{MN}
     \le \frac{(C_1)^2
     	\bbE[\|\bsig - \tilde\bsig\|_2^2]}{M} 
	\stackrel{\eqref{e:ctrl.prbs.appx.assum}}{\le}
	\frac{(C_1)^2 \delta}{\alpha}\,,
    \]
which is $\le C_1 \delta$ by adjusting $C_1$ (note that it still depends on $\alpha$ only). Meanwhile, by H\"older's inequality,
	\begin{align*}
	\textup{(b)} 
	&\equiv
	\frac{\bbE[
	\|\bG\|_{\op}^2 
	\|\bsig - \tilde\bsig\|_2^2;
	\cE^c]}{MN}
	\le
	 \frac{2\bbE[
	\|\bG\|_{\op}^2 
	(\|\bsig\|_2^2
	+\|\tilde\bsig\|_2^2);
	\cE^c]}{MN}\\
	&\le
	\frac{2
	\P(\cE^c)^{1/3}
	\E[\|\bG\|_{\op}^6]^{1/3}
	(\E[\|\bsig\|_2^6]^{1/3}
	+\E[\|\tilde\bsig\|_2^6]^{1/3})
	}{MN}\,.
	\end{align*}
The map $\bG \mapsto \|\bG\|_{\op}$ is $1$-Lipschitz, and the maps $(\bG,\bg^{\aux}) \mapsto \|\bsig\|_2$ and $(\bG,\bg^{\aux}) \mapsto \|\tilde\bsig\|_2$ are $L$-Lipschitz by assumption.
    Combining the gaussian concentration of these quantities 
    with the above bound on $\P(\cE^c)$ and the assumption $\bbE[\|\bsig\|^2] = \bbE[\|\tilde\bsig\|^2] = N$ (by Definition~\ref{d:Lip}), we conclude
    \[
        \max\bigg\{
            \bbE[\|\bG\|_{\op}^6]^{1/3}, 
            \bbE[\|\bsig\|_2^6]^{1/3}, 
            \bbE[\|\tilde\bsig\|_2^6]^{1/3}
        \bigg\} \le C_2N
    \]
for some $C_2 = C_2(\alpha,L)$. By adjusting $c_0$ and $C_2$ we conclude $\textup{(b)} \le C_2 e^{-c_0 N}$.
    Altogether this shows
    \[
        \bbW_2\Big(
        \mu(\cA_N),\mu(\tilde \cA_N)
        \Big)^2
        \le \textup{(a)} + \textup{(b)}
        \le C_1 \delta + C_2 e^{-c_0N}\,,
    \]
    proving the first estimate of the proposition.
    The second estimate is only simpler:
    \begin{align*}
        \bbW_2\Big(
        \mu^{\Ising}(\cA_N),\mu^{\Ising}(\tilde \cA_N)
        \Big)^2
        &\le \bbE \bigg[\bbW_2\Big(
            \EmpDist(\bsig),
            \EmpDist(\tilde\bsig)
        \Big)^2\bigg]
        \le 
        \frac{\bbE[\|\bsig - \tilde\bsig\|^2] }{N} 
        \stackrel{\eqref{e:ctrl.prbs.appx.assum}}{\le} 
        \delta\,.
    \end{align*}
This concludes the proof.
\end{proof}    

\begin{lem}\label{l:sM.Lip.union.q0}
The set $\sM^{\Lip}(\alpha)$ from \eqref{e:sM.Lip} can be expressed as
    \[\sM^{\Lip}(\alpha) = \bigcup_{q_0 \in [0,1]}\sM^{\Lip}(\alpha,q_0)
    \]
for $\sM^{\Lip}(\alpha,q_0)$
as defined by \eqref{e:sM.Lip.q0}.
\begin{proof}
    For any $q_0 \in [0,1]$, the inclusion $\sM^{\Lip}(\alpha,q_0) \subseteq \sM^{\Lip}(\alpha)$ is trivial, so it suffices to prove
    \beq\label{e:sM.Lip.union.q0.goal}
        \sM^{\Lip}(\alpha) \subseteq \bigcup_{q_0 \in [0,1]} \sM^{\Lip}(\alpha,q_0)\,.
    \eeq
Consider any sequence $\epsilon_n\downarrow0$. Then, for $\cM^{\Lip}(\alpha;L,\epsilon)$ defined in \eqref{e:cM.Lip}, we have
    \[ \sM^{\Lip}(\alpha) 
    =\adjustlimits
    \bigcap_n \bigcup_{L>0}
    \cM^{\Lip}(\alpha;L,\epsilon_n)
    \,.
    \]
Consider any measure $\mu \in \sM^{\Lip}(\alpha)$. For any $n\ge 1$, there exists $L_n$ (depending on $\mu,\epsilon_n$) and a sequence of $L_n$-Lipschitz algorithms $(\cA^n_{N_j})_{j\ge 1}$ with 
    \beq
        \label{e:sM.Lip.union.q0.W2.limit}
        \WERR(\cA^n_{N_j},\mu) =
           \max\bigg\{
            \bbW_2\Big(
            \mu(\cA^n_{N_j}),
            \mu\Big),
            \bbW_2\Big(
            \mu^{\Ising}(\cA^n_{N_j}), \cP(\{\pm 1\})\Big)
        \bigg\} \le \epsilon_n\,.
    \eeq
Abbreviate $\bsig_{n,j} 
\equiv \cA^n_{N_j}(\bG,\bg^\aux)$. 
By passing to a subsequence of the $N_j$, we can ensure that
    \beq
        \label{e:sM.Lip.union.q0.chi.limit}
	q_0(n,j)
	\equiv \chi_{\cA^n_{N_j}}(0)
	= \frac{\|\bbE\bsig_{n,j}\|^2}{N_j}	
	\stackrel{j\to\infty}
	{\longrightarrow} 
	q_0(n)
    \eeq
    for some $q_0(n) \in [0,1]$.
    Let $q_0 \in [0,1]$ be any subsequential limit of the $q_0(n)$ as $n\to\infty$.
    We will show that $\mu \in \sM^{\Lip}(\alpha,q_0)$, from which the result follows.

    By passing to a subsequence of the $\epsilon_n$, we can assume that $q_0(n) \to q_0$ as $n\to\infty$.
    We will show that for $C_1(\alpha)$ as in Lemma~\ref{l:control.problems.approx}, 
    $\tilde L_n \equiv 2L_n$, and
    \[\tilde\epsilon_n 
    = \epsilon_n + 2\bigg(
    	C_1(\alpha) \Big[
	|q_0 - q_0^n| + \epsilon_n
	\Big]
	\bigg)^{1/2}\]
    there exists a sequence of $\tilde L_n$-Lipschitz $(\tilde\cA^n_{N_j})_{j\ge 1}$ with
    \beq
        \label{e:sM.Lip.union.q0.W2.limit.perturbed}
        \WERR(\tilde{\cA}^n_{N_j},\mu)=
        \max\bigg\{
            \bbW_2\Big(\mu(\tilde \cA^n_{N_j}), \mu\Big),
            \bbW_2\Big(\mu^{\Ising}(\tilde \cA^n_{N_j}), \cP(\{\pm 1\})
            \Big)
        \bigg\}
        \le \tilde\epsilon_n
    \eeq
as well as (with the abbreviation $\tilde{\bsig}_{n,j}
\equiv \tilde{\cA}^n_{N_j}(\bG,\bg^\aux)$) 
    \beq
        \label{e:sM.Lip.union.q0.chi.limit.perturbed}
     \tilde{q}_0(n,j)
     \equiv 
     \chi_{\tilde\cA^n_{N_j}}(0) 
    = \frac{\|\bbE\tilde{\bsig}_{n,j}\|^2}{N_j} 
     = q_0
    \eeq
for every $n,j$. This implies that
    \[
        \mu \in \cM^{\Lip}(\alpha,q_0;\tilde L_n,\tilde\epsilon_n)
        \subseteq \bigcup_{L>0} \cM^{\Lip}(\alpha,q_0;L,\tilde\epsilon_n)\,.
    \]
    Since this set is decreasing in $\tilde\epsilon_n$, and $\lim_{n\to\infty} \tilde\epsilon_n = 0$, we infer that
    \[
        \mu \in 
        \adjustlimits
        \bigcap_n \bigcup_{L>0} \cM^{\Lip}(\alpha,q_0;L,\tilde\epsilon_n)
        = \adjustlimits
        \bigcap_{\epsilon>0} \bigcup_{L>0} \cM^{\Lip}(\alpha,q_0;L,\epsilon)
        = \sM^{\Lip}(\alpha,q_0)\,.
    \]
    As this is true for all $\mu \in \sM^{\Lip}(\alpha)$, the inclusion \eqref{e:sM.Lip.union.q0.goal} follows.

It remains to construct $\tilde L_n$-Lipschitz algorithms satisfying \eqref{e:sM.Lip.union.q0.W2.limit.perturbed} and \eqref{e:sM.Lip.union.q0.chi.limit.perturbed}. If $q_0<1$, 
then by omitting finitely many terms of the sequences, we can assume
$q_0(n,j) < 1$ and 
	\[\frac{1-q_0}{1-q_0(n,j)} \le 4\,.\]
If $q_0 > 0$, for the same reason we can assume $q_0(n,j) > 0$.
    Then define the perturbed algorithms $\tilde\cA^n_{N_j}$ by setting 
    \[
        \tilde{\bsig}_{n,j}
        \equiv \tilde\cA^n_{N_j}(\bG,\bg^\aux)
        \equiv \bigg(
        \frac{q_0}{q_0(n,j)}
        \bigg)^{1/2}
        (\bbE\bsig_{n,j})
        + \bigg(
        \frac{1-q_0}{1-q_0(n,j)}\bigg)^{1/2} (
        \bsig_{n,j}
         - \bbE\bsig_{n,j})\,.
    \]
If $q_0 = 0$, then we interpret the first term on the right-hand side to be zero. If $q_0=1$, then we interpret the second term on the right-hand side to be zero. If $q_0\in[0,1)$, then the Lipschitz constant of the perturbed algorithm $\tilde\cA^n_{N_j}$, as a function of $(\bG,\bg^\aux)$, is bounded by 
    \[
        \bigg(
        \frac{1-q_0}{1-q_0(n,j)}
        \bigg)^{1/2} L_n \le 2 L_n = \tilde L_n\,.
    \]
If $q_0=1$, then the perturbed algorithm is deterministic, and so $0$-Lipschitz. 
Since $\bbE[\|\bsig_{n,j}\|^2] = N_j$ by Definition~\ref{d:Lip}, we also have
    \[
        \frac{\bbE[\|\tilde{\bsig}_{n,j}\|^2]}{N_j} 
        = \frac{q_0}{q_0(n,j)} \cdot q_0(n,j) + \frac{1-q_0}{1-q_0(n,j)} \cdot (1-q_0(n,j)) = 1\,.
    \]
    Thus $(\tilde\cA^n_{N_j})_{j\ge 1}$ is a sequence of $\tilde L_n$-Lipschitz algorithms, again in the sense of Definition~\ref{d:Lip}. We have
    \[\tilde{q}_0(n,j)
        = \chi_{\tilde\cA^n_{N_j}}(0)
        = \frac{\|\bbE \tilde{\bsig}_{n,j}\|^2}{N_j} 
        = \frac{q_0 \|\bbE \bsig_{n,j}\|^2}{q_0(n,j)} 
        = q_0\,,
    \]
    which implies \eqref{e:sM.Lip.union.q0.chi.limit.perturbed}.
    Finally,
    \begin{align*}
        &\frac{\bbE[\|\tilde{\bsig}_{n,j} - \bsig_{n,j}\|^2]}{N_j} 
        = \frac{1}{N_j} \bbE \bigg[\bigg\|
            \bigg\{
            \bigg(
            \frac{q_0}{q_0(n,j)}
            \bigg)^{1/2} - 1\bigg\}\bbE[\bsig_{n,j}]
            + \bigg\{
            \bigg(\frac{1-q_0}{1-q_0(n,j)}
            \bigg)^{1/2} - 1\bigg\} (\bsig_{n,j} - \bbE[\bsig_{n,j}])
        \bigg\|^2\bigg] \\
        &\qquad= \bigg(q_0^{1/2} 
        	- q_0(n,j)^{1/2}\bigg)^2
        + \bigg((1 - q_0)^{1/2} - 
        (1 - q_0(n,j))^{1/2}\bigg)^2\\
        &\qquad\stackrel{(*)}{\le} 
        2|q_0 - q_0(n,j)|
        \stackrel{(**)}{\le} 2(|q_0 - q_0^n| + \epsilon_n)\,,
    \end{align*}
    where the step $(*)$ uses the inequality $|x^{1/2}-y^{1/2}| \le |x-y|^{1/2}$  for all $x,y\ge0$, and the step $(**)$ holds, after possibly omitting finitely many $j$, because $q_0(n,j)$ tends to $q_0(n)$ as $j\to\infty$.
    By Lemma~\ref{l:control.problems.approx}, we conclude
    \begin{align*}
        &\WERR(\cA^n_{N_j},\tilde\cA^n_{N_j})
        =\max\bigg\{
            \bbW_2\Big(\mu(\cA^n_{N_j}),\mu(\tilde \cA^n_{N_j})\Big),
            \bbW_2\Big(\mu^{\Ising}(\cA^n_{N_j}),\mu^{\Ising}(\tilde \cA^n_{N_j})\Big)
        \bigg\} \\
        &\qquad\le \bigg\{
            C_1(\alpha) \cdot 2(|q_0 - q_0^n| + \epsilon_n)
            + C_2(\alpha,\tilde L_n) e^{-c_0 N_j}
        \bigg\}^{1/2}
        \le 2\bigg\{
        C_1(\alpha) (|q_0 - q_0^n| + \epsilon_n)\bigg\}^{1/2}
    \end{align*}
again after possibly omitting finitely many $j$.
    Combined with \eqref{e:sM.Lip.union.q0.W2.limit} we conclude \eqref{e:sM.Lip.union.q0.W2.limit.perturbed}. 
\end{proof}
\end{lem}

\begin{lem}
    \label{l:sM.Lip.1.trivial}
The set $\sM^{\Lip}(\alpha,1)$, as defined by \eqref{e:sM.Lip.q0}, is a singleton set consisting of the standard gaussian measure:
$\sM^{\Lip}(\alpha,1)=\{\cN(0,1)\}$.
\begin{proof}
    Consider any $\mu \in \sM^{\Lip}(\alpha,1)$. By the definition \eqref{e:sM.Lip.q0}, this means that for any $\epsilon > 0$, there exists $L>0$ and a sequence of $L$-Lipschitz algorithms $(\cA_{N_j})_{j\ge 1}$ such that
    \begin{align*}
        &\WERR(\cA_{N_j}, \mu)
        =\max\bigg\{
            \bbW_2\Big( \mu(\cA_{N_j}), \mu \Big),
            \bbW_2\Big(\mu^{\Ising}(\cA_{N_j}), \cP(\{\pm 1\}) \Big)
        \bigg\} \le \epsilon\,, \\
        &q_0(j) \equiv \chi_{\cA_{N_j}}(0) 
        = \frac{\|\bbE \bsig_j\|_2^2}{N_j} 
        \stackrel{j\to\infty}{\longrightarrow} 1
        \,,
    \end{align*}
where we abbreviate $\bsig_j\equiv \cA_{N_j}(\bG,\bg^\aux)$. 
    Define the perturbed (deterministic) algorithm $\tilde\cA_{N_j}$ by setting 
    \[
        \tilde{\bsig}_j
        \equiv \tilde\cA_{N_j}(\bG,\bg^\aux) 
        = \frac{\bbE\bsig_j}{q_0(j)^{1/2}} \,.
    \]
This clearly satisfies $\bbE[\|\tilde{\bsig}_j\|_2^2]  = N_j$, so $\tilde\cA_{N_j}$ is a $0$-Lipschitz algorithm in the sense of Definition~\ref{d:Lip}.
Since $\tilde{\bsig}_j$ is a deterministic vector in $N_j$ dimensions of norm $(N_j)^{1/2}$, we have $\mu(\tilde \cA_{N_j}) = \cN(0,1)$.
    Note that
    \begin{align*}
        \frac{\bbE[
            \|\bsig_j - \tilde{\bsig}_j\|^2]}{N_j} 
        &= \frac{1}{N_j} \bbE\bigg[ \bigg\|
            \bigg(\frac{1}{q_0(j)^{1/2}} - 1 \bigg) 
            	(\bbE\bsig_j)
            - (\bsig_j - \bbE\bsig_j)
        \bigg\|^2 \bigg] \\
        &= \Big(1 - q_0(j)^{1/2}\Big)^2 + \Big(1 - q_0(j)\Big)
        \le 2\Big(1-q_0(j)\Big)\,,
    \end{align*}  
again using that inequality $|x^{1/2}-y^{1/2}| \le |x-y|^{1/2}$ for $x,y\ge0$.
    Then, by Lemma~\ref{l:control.problems.approx}, with constants $c_0$, $C_1(\alpha)$, $C_2(\alpha,L)$ as defined therein, we have
    \[
        \bbW_2\Big(\mu(\cA_{N_j}),
        	\mu(\tilde \cA_{N_j}) = \cN(0,1) 
	\Big)
        \le \bigg\{
        C_1(\alpha) 2(1-q_0(j)) + C_2(\alpha,L) e^{-c_0N_j} \bigg\}^{1/2}\,.
    \]
Combined with the assumption that $\bbW_2(\mu(\cA_{N_j}), \mu) \le \epsilon$ for all $j$, we find $\bbW_2(\mu,\cN(0,1)) \le \epsilon$.
    Since this holds for all $\epsilon > 0$, we conclude $\mu = \cN(0,1)$.
    We have thus shown $\sM^{\Lip}(\alpha,1) \subseteq \{\cN(0,1)\}$.
    The reverse inclusion $\cN(0,1) \in \sM^{\Lip}(\alpha,1)$ is witnessed by the (sequence of) deterministic algorithms
    \[
        \cA_N(\bG,\bg^{\aux}) = (1,\ldots,1) \in \bbR^N\,, 
    \]
    which are $0$-Lipschitz with $\mu(\cA_N) = \cN(0,1)$, $\mu^{\Ising}(\cA_N) = \delta_1 \in \cP(\{\pm 1\})$, and $\chi_{\cA_N}(0) = 1$.
\end{proof}
\end{lem}

\subsection{Consequences of earlier results}
\label{ss:control.problems.earlier.results}

In this subsection, we use results from \S\ref{s:sde}--\ref{sec:IAMP} to prove the following inclusions:
\begin{itemize}
    \item Lemma~\ref{l:IAMP.in.Lip} proves $\sM^{\IAMP}(\alpha,q_0) \subseteq \sM^{\Lip}(\alpha,0)$;
    \item Lemma~\ref{l:Lip.in.BOGP} proves $\sM^{\Lip}(\alpha,q_0) \subseteq \sM^{\BOGP}(\alpha,q_0)$;
    \item Lemma~\ref{l:prf3.in.IAMP} proves $\sM^{\sxi}(\alpha,q_0) \subseteq \sM^{\IAMP}(\alpha,q_0)$.
\end{itemize}
We recall that $\sM^{\Lip}(\alpha,q_0)$ is defined by \eqref{e:sM.Lip.q0}, while the remaining classes $\cM^\star(\alpha,q_0)$ are specified by Definitions \ref{d:SDE.conds}--\ref{d:measure.classes}.

\begin{lem}\label{l:IAMP.in.Lip}
    For all $q_0 \in [0,1)$, we have $\sM^{\IAMP}(\alpha,q_0) \subseteq \sM^{\Lip}(\alpha,0)$. 

\begin{proof}
    This is a consequence of Theorem~\ref{thm:IAMP-main}\ref{i:IAMP-main-centered}.
    Consider any $\mu \in \cM^{\IAMP}(\alpha,q_0;L,\epsilon)$. Then, there exists
    \[
        (q_*,b,\sigma,w,p,\zeta,\zeta^\Ising) \in \Adm^{\IAMP}(\alpha,q_0;L,\epsilon)
    \]
    such that $\mu(q_*,b,\sigma,p,\zeta) = \mu$.
    Let $\tilde p$ denote the restriction of $p$ to $[q_*,1] \subseteq [q_0,1]$.
    We will argue that $(b,\sigma,w,\tilde p,\zeta,\zeta^\Ising)$ satisfies the assumptions of Theorem~\ref{thm:IAMP-main}\ref{i:IAMP-main-centered} with $q_*$ in place of $q_0$.
    Indeed, conditions~\ref{it:1.1}, \ref{it:2.1}, \ref{it:3.1}, \ref{it:4.2}, \ref{it:5.1}, and \ref{it:6.1} exactly correspond to assumptions~(\ref{i:IAMP-main-coefs}), (\ref{i:IAMP-main-diffusivity}), (\ref{i:IAMP-main-p}), (\ref{i:IAMP-main-budget}), (\ref{i:IAMP-main-endpt}), and (\ref{i:IAMP-main-init}) of Theorem~\ref{thm:IAMP-main}, while condition~\ref{it:1.1} also provides the additional assumption in Theorem~\ref{thm:IAMP-main}\ref{i:IAMP-main-centered} that $w(t,\cdot)$ is even and $\zeta^\Ising$ is symmetric. By Theorem~\ref{thm:IAMP-main}\ref{i:IAMP-main-centered}, there exists $\iota = \iota(\epsilon)$ with $\lim_{\epsilon\rightarrow 0} \iota(\epsilon) = 0$ such that for all (sufficiently large) $N$, there exists a $C(L,\epsilon)$-Lipschitz algorithm $\cA_N$ with
    \[
        \max\Big(
            \bbW_2(\mu(\cA_N), \mu),
            \bbW_2(\mu^{\Ising}(\cA_N), \cP(\{\pm 1\}))
        \Big) \le \iota
    \]
and
    $\bbE[\cA_N(\bG,\bg^\aux)] = \bzero$.
    This implies $\chi_{\cA_N}(0) = 0$.
    Recalling the definition \eqref{e:cM.Lip.q0}, we have $\mu \in \cM^{\Lip}(\alpha,0;C(L,\epsilon),\iota(\epsilon))$.
    We have thus shown
    \[
        \cM^{\IAMP}(\alpha,q_0;L,\epsilon) \subseteq \cM^{\Lip}(\alpha,0;C(L,\epsilon),\iota(\epsilon))\,.
    \]
    Thus, for $\cB_\epsilon$ as defined in \eqref{e:W2.ball.notation},
    \[
        \cB_{\epsilon}\Big(\cM^{\IAMP}(\alpha,q_0;L,\epsilon)\Big)
        \subseteq \cB_{\epsilon}\Big(\cM^{\Lip}(\alpha,0;C(L,\epsilon),\iota(\epsilon))\Big)\,.
    \]
    Taking a union over $L>0$ shows
    \[
        \bigcup_{L>0} \cB_{\epsilon}\Big(\cM^{\IAMP}(\alpha,q_0;L,\epsilon)\Big)
        \subseteq \bigcup_{L>0} \cB_{\epsilon}\Big(\cM^{\Lip}(\alpha,0;L,\iota(\epsilon))\Big)
        \subseteq \bigcup_{L>0} \cB_{\tilde\epsilon}\Big(\cM^{\Lip}(\alpha,0;L,\tilde\epsilon)\Big)\,,
    \]
    where $\tilde\epsilon = \tilde\epsilon(\epsilon) = \max(\epsilon,\iota(\epsilon))$.
    Since $\lim_{\epsilon\rightarrow 0} \tilde\epsilon(\epsilon) = 0$, taking an intersection over $\epsilon$ yields
    \begin{align*}
        \sM^{\IAMP}(\alpha,q_0)
        &\stackrel{\eqref{e:measure.classes.simplification.without.L0}}{=}
        \adjustlimits
        \bigcap_{\epsilon>0} \bigcup_{L>0} \cB_{\epsilon}\Big(\cM^{\IAMP}(\alpha,q_0;L,\epsilon)\Big) \\
        &\,\,\subseteq
        \adjustlimits \bigcap_{\epsilon>0} \bigcup_{L>0} \cB_{\epsilon}\Big(\cM^{\Lip}(\alpha,0;L,\epsilon)\Big)
        \stackrel{\eqref{e:sM.Lip.q0}}{=} \sM^{\Lip}(\alpha,0)\,. \qedhere
    \end{align*}
  \end{proof}
\end{lem}

\begin{lem}\label{l:Lip.in.BOGP}
    For all $q_0 \in [0,1)$, we have $\sM^{\Lip}(\alpha,q_0) \subseteq \sM^{\BOGP}(\alpha,q_0)$.
\begin{proof}
    This is a consequence of Theorem~\ref{thm:BOGP-hardness-main}.
    Consider any $\mu \in \cM^{\Lip}(\alpha,q_0;L,\epsilon)$, where $L$ is sufficiently large depending on $\alpha,q_0,\epsilon$. Recalling \eqref{e:cM.Lip.q0}, this means there exists a sequence of $L$-Lipschitz algorithms $(\cA_j)^\circ$ (acting in $N_j$ dimensions) such that
    \begin{align} \nonumber 
    &\WERR( (\cA_j)^\circ,\mu)
 \le \epsilon\,,\\
        \label{e:Lip.in.BOGP.q0.limit}
   &     q_0(j)
        \equiv \chi_{ (\cA_j)^\circ}(0)
         \stackrel{j\to\infty}{\longrightarrow} q_0\,.
    \end{align}
    Let $(\cA_j)_{j\ge 1}$ be the perturbations of the $( (\cA_j)^\circ)_{j\ge 1}$ given by Proposition~\ref{p:wlogable}, which are also $L$-Lipschitz.
    Proposition~\ref{p:wlogable} further yields
    \[
        \frac{1}{N} \bbE \Big[\| (\cA_j)^\circ(\bG,\bg^{\aux}) 
        	- \cA_j(\bG,\bg^{\aux})\|^2\Big]
        \le \frac{2}{(1-q_0(j))L^2}\,.
    \]
Thus Lemma~\ref{l:control.problems.approx} implies that for $c_0$, $C_1(\alpha)$ and $C_2(\alpha,L)$ given therein,
	\[
	\WERR(\cA_j,(\cA_j)^\circ)
	\le 
	\bigg\{
        C_1(\alpha) \cdot \frac{2}{(1-q_0(j))L^2} + C_2(\alpha,L) e^{-c_0N_j}
        \bigg\}^{1/2} 
        \le \epsilon \]
    where the last inequality holds if $L$ is sufficiently large depending on $\alpha,q_0,\epsilon$ (consistent with our assumptions), and $j$ sufficiently large depending on $q_0,\epsilon,L$. Hence, by omitting finitely many $j$, we can assume it holds for all $j$ that
    \beq
        \label{e:Lip.in.BOGP.hardness.hypothesis}
        \WERR(\cA_j,\mu)=
        \max\bigg\{
            \bbW_2\Big(\mu(\cA_j), \mu\Big),
            \bbW_2\Big(\mu^{\Ising}(\cA_j), \cP(\{\pm 1\})\Big)
        \bigg\}
        \le 2\epsilon\,,
    \eeq
By passing to a further subsequence of the $j$'s, we can assume that for all $j$ we have 
    \[
         \max\bigg\{
            \bbW_2\Big(\mu(\cA_j), \mu\Big),
            \bbW_2\Big(\mu^{\Ising}(\cA_j), \mu^\Ising\Big)
        \bigg\}
         \le 3\epsilon
    \]
for some fixed $\mu^\Ising \in \cP(\{\pm 1\})$.
    By passing to a further subsequence, we can assume (as in Theorem~\ref{thm:BOGP-hardness-main}) that $\chi_{\cA_j}$ converges to a limiting function $\chi$, both pointwise and in $L^1$.
    Since Proposition~\ref{p:wlogable} ensures 
    \[
        \chi_{\cA_j}(0) = \chi_{\cA^\circ_j}(0) = q_0(j)\,,
    \]
    by \eqref{e:Lip.in.BOGP.q0.limit} we have $\chi(0) = q_0$.

We now apply Theorem~\ref{thm:BOGP-hardness-main} to the sequence $(\cA_j)_{j\ge 1}$ where, treating $\epsilon$ as a function of $L$ tending to $0$ sufficiently slowly, the hypothesis \eqref{e:BOGP-hardness-main-hypothesis} holds by \eqref{e:Lip.in.BOGP.hardness.hypothesis}.
    Then, for some $\tilde\epsilon = \tilde\epsilon(\epsilon)$ tending to $0$ sufficiently slowly as $\epsilon\rightarrow 0$, Theorem~\ref{thm:BOGP-hardness-main} provides controls $(q_*,b,\sigma,w,p\equiv \chi^{-1},\zeta,\zeta^\Ising)$ satisfying the conclusions of that theorem, where all $o_L(1)$ error terms are replaced by $\tilde\epsilon$. We assume without loss that $\tilde{\epsilon}\le1$. 
    We will now verify that for $C(L)$ given by Theorem~\ref{thm:BOGP-hardness-main},
    \beq
        \label{e:Lip.in.BOGP.controls.in.BOGP}
         (q_*,b,\sigma,w,p\equiv \chi^{-1},\zeta,\zeta^\Ising) 
         \in \Adm^{\BOGP}\Big(
         \alpha,q_0;
         \tilde C = \max\{C(L),L^2\}, \tilde\epsilon 
         \Big)\,.
    \eeq
    Indeed, we verify each of the conditions from Definition~\ref{d:SDE.classification} specifying the $\Adm^{\BOGP}$ class:
    \begin{itemize}
\item Theorem~\ref{thm:BOGP-hardness-main}\ref{it:BOGP-hardness-main-Lip}
    implies condition~\ref{it:1.1bogp};

\item Theorem~\ref{thm:BOGP-hardness-main}\ref{it:BOGP-hardness-main-diffusivity}
    implies condition~\ref{it:2.1}; 

\item We have $p(q_0) = 0$ because we defined $\chi(0) = q_0$. 
        Theorem~\ref{thm:BOGP-hardness-main} asserts $p(q_*) \le \tilde\epsilon$, $p$ is concave, and $\|p\|_{C^1([q_0,1])} \le L^2 \le \tilde C$, which implies condition~\ref{it:3.2}.
        
\item Theorem~\ref{thm:BOGP-hardness-main}\ref{it:BOGP-hardness-main-budget}
implies condition~\ref{it:4.1};

\item Let $\tilde\mu^\Ising = \mu^\Ising(q_*,w,\zeta^\Ising)$ and $Y(1) \sim \tilde\mu^\Ising$.
        Then, since $\mu^\Ising \in \cP(\{\pm 1\})$,
        \[
            \bbE[(|Y(1)|-1)^2]
            = \bbW_2(\tilde\mu^\Ising,\cP(\{\pm 1\}))^2
            \le \bbW_2(\tilde\mu^\Ising,\mu^\Ising)^2
            \le \tilde\epsilon^2
            \le \tilde\epsilon\,,
        \]
        where the second-last inequality is by the second inequality of \eqref{eq:BOGP-hardness-main-endpt} from Theorem~\ref{thm:BOGP-hardness-main},
        and the last inequality is by our assumption $\tilde\epsilon\le1$. This proves condition~\ref{it:5.1}.

\item
Theorem~\ref{thm:BOGP-hardness-main}\ref{it:BOGP-hardness-main-initial}
implies condition~\ref{it:6.1};
\item Condition~\ref{it:7.1} is trivial; 

        \item 
        Theorem~\ref{thm:BOGP-hardness-main}\ref{it:BOGP-hardness-main-qstar-condition} implies condition~\ref{it:8.1}.

\end{itemize}
The final assertion of Theorem~\ref{thm:BOGP-hardness-main} implies that $\tilde\mu = \mu(q_*,b,\sigma,p,\zeta)$ satisfies $\bbW_2(\tilde\mu,\mu) \le \tilde\epsilon$. Since \eqref{e:Lip.in.BOGP.controls.in.BOGP} implies $\mu\in\cM^{\BOGP}(\alpha,q_0;\tilde C,\tilde\epsilon)$, 
we conclude
    \[
        \mu \in \cB_{\tilde\epsilon}\Big(
            \cM^{\BOGP}(\alpha,q_0;\tilde C,\tilde\epsilon)\,.
        \Big)
    \]
    and thus (for $L$ sufficiently large depending on $\epsilon$, as we previously assumed)
    \[
        \cM^{\Lip}(\alpha,q_0;L,\epsilon)
        \subseteq \cB_{\tilde\epsilon(\epsilon)}\Big(
            \cM^{\BOGP}(\alpha,q_0;\tilde C(L),
            	\tilde\epsilon(\epsilon))
        \Big)\,.
    \]
Now recall \eqref{e:sM.Lip.q0} and \eqref{e:measure.classes.simplification.without.L0}: since $\cM^{\Lip}(\alpha,q_0;L,\epsilon)$ and $\cM^{\BOGP}(\alpha,q_0;L,\epsilon)$ are both increasing in $L$ and decreasing as $\epsilon \downarrow 0$, taking a union over $L$ and an intersection over $\epsilon$ concludes the proof.
\end{proof}
\end{lem}

\begin{lem}\label{l:prf3.in.IAMP}
For all $q_0 \in [0,1)$, we have $\sM^{\sxi}(\alpha,q_0) \subseteq \sM^{\IAMP}(\alpha,q_0)$. 

\begin{proof} We repeat the most relevant part of Definition~\ref{d:SDE.classification} for this proof:
	\begin{center}
	\begin{tabular}[h]{r|cccccccc}
	\MSRheader\\
	\hline
	\PRFTHRdefn\\
	\IAMPdefn
	\end{tabular}
	\end{center}
Consider any controls 
    \beq\label{e:prf3.in.IAMP.input.class}
        (\acute q_*=q_0, \acute b, \acute \sigma, \acute w, p, \acute \zeta=\delta_0, \acute \zeta^\Ising) \in \Adm^{\sxi}(\alpha,q_0;L)\,.
    \eeq
Let $\acute{X},\acute{Y}$ be defined by \eqref{e:formal.sde.main} and \eqref{e:formal.sde.ising} with these parameters. We verify that the $\sxi$ conditions imply that the assumptions of Theorem~\ref{thm:SDE-smoothing-general} are satisfied, with $\acute q_*=q_0$ by
condition~\ref{it:8.2}:
\begin{itemize}
\item Condition~\ref{it:1.2even}
implies that 
Theorem~\ref{thm:SDE-smoothing-general} assumption~\ref{i:SDE-smoothing-general.INPUT.1} is satisfied,
with $\mathscr{L}(Z^\Ising(t))$ even;
\item Condition~\ref{it:2.2}
implies that 
Theorem~\ref{thm:SDE-smoothing-general} assumption~\ref{i:SDE-smoothing-general.INPUT.2b}
is satisfied;
\item Condition~\ref{it:3.1}
implies that 
Theorem~\ref{thm:SDE-smoothing-general} assumption~\ref{i:SDE-smoothing-general.INPUT.implies.3a}
is satisfied;
\item Condition~\ref{it:4.3}
implies that 
Theorem~\ref{thm:SDE-smoothing-general} assumption~\ref{e:Z.budget.assumption}
is satisfied;

\item Condition~\ref{it:5.2}
implies that 
Theorem~\ref{thm:SDE-smoothing-general} assumption~\ref{i:SDE-smoothing-general.INPUT.5b}
is satisfied;

\item Condition~\ref{it:6.2}
implies that 
Theorem~\ref{thm:SDE-smoothing-general} assumption~\ref{i:SDE-smoothing-general.INPUT.6b} is satisfied;

\item Condition~\ref{it:7.2}
implies that 
Theorem~\ref{thm:SDE-smoothing-general} assumption~\ref{i:SDE-smoothing-general.INPUT.7b} is satisfied.
\end{itemize}
Consider any $\epsilon > 0$, and let $(q_*, b, \sigma, w, \zeta, \zeta^\Ising)$ and $C=C(L,\epsilon)$ be as given by Theorem~\ref{thm:SDE-smoothing-general}.
    We will argue that
    \beq\label{e:prf3.in.IAMP.output.class}
        (q_*, b, \sigma, w, p, \zeta, \zeta^\Ising) \in \Adm^{\IAMP}\Big(\alpha,q_0;
        \tilde C \equiv
        \max\{ 
        C(L,\epsilon),L\},
        \tilde\epsilon 
        \equiv
        3\epsilon + \epsilon^2\Big)\,.
    \eeq
Let $X,Y$ be defined by \eqref{e:formal.sde.main} and \eqref{e:formal.sde.ising} with these parameters. 
    We verify each of the conditions defining $\Adm^{\IAMP}$: 
\begin{itemize}
\item Theorem~\ref{thm:SDE-smoothing-general} conclusion \ref{it:SDE-smoothing-general.coefs}, combined with the final assertion of Theorem~\ref{thm:SDE-smoothing-general}, implies condition~\ref{it:1.1}.

\item Theorem~\ref{thm:SDE-smoothing-general} conclusion~\ref{it:SDE-smoothing-general.OUTPUT.diffus}
implies condition~\ref{it:2.1}, using that $\epsilon \le\tilde{\epsilon}$.

\item By assumption, $p$ satisfies condition~\ref{it:3.1} with respect to the parameters in \eqref{e:prf3.in.IAMP.input.class}, i.e. $1/L \le \acute p(q_0) \le \epsilon$ and $\|p\|_{C^2([q_0,1])}\le L$.
        Since $\tilde C \ge L$ and $\tilde \epsilon \ge \epsilon$, $p = \acute p$ also satisfies condition~\ref{it:3.1} with respect to the parameters in \eqref{e:prf3.in.IAMP.output.class}.

\item Theorem~\ref{thm:SDE-smoothing-general} conclusion~\ref{it:SDE-smoothing-general.budget}
        implies condition~\ref{it:4.2}, since $\epsilon \le \tilde\epsilon$.

\item Recalling \eqref{e:control.problems.endpoint.measures}, we have 
$\acute{Y}(1)\sim\acute \mu^\Ising = \mu^\Ising(\acute q_*, \acute w, \acute \zeta^\Ising)$ and 
$Y(1)\sim\mu^\Ising = \mu^\Ising(q_*,w,\zeta^\Ising)$. Recalling condition~\ref{it:5.2} of the $\sxi$ class, the second inequality of \eqref{eq:SDE-smoothing-general-endpt} from Theorem~\ref{thm:SDE-smoothing-general} implies
        \[
            \bbE\Big[
            (|Y(1)|-1)^2\Big]
            = \bbW_2(\mu^\Ising, \cP(\{\pm 1\}))^2
            \le \bbW_2(\mu^\Ising, \acute \mu^{\Ising})^2
            \le \epsilon^2
            \le \tilde\epsilon\,,
        \]
which implies condition~\ref{it:5.1}.

\item Again recalling \eqref{e:control.problems.endpoint.measures}, we have 
$\acute{X}(q_0)\sim\acute{\zeta}=\delta_0$,
$X(q_*)\sim\zeta$,
$\acute{Y}(q_0)\sim\acute{\zeta}^\Ising$, and
$Y(q_*)\sim\zeta^\Ising$. Thus, Theorem~\ref{thm:SDE-smoothing-general} conclusion~\ref{it:SDE-smoothing-general.initial} gives 
        \[
            \bbE[X(q_*)^2] = \bbW_2(\zeta,\acute\zeta)^2
            \le \epsilon^2
            \le \tilde\epsilon\,,
        \]
Recalling condition~\ref{it:6.2} of the $\sxi$ class, we also have
        \begin{align*}
           & \Big|\bbE [Y(q_*)^2]^{1/2} 
            - (q_0)^{1/2}
            \Big|
            =\Big|\bbE [Y(q_*)^2]^{1/2} - \bbE [\acute Y(q_0)^2]^{1/2}
            \Big|\\
            &\qquad= \Big|\bbW_2(\zeta^\Ising,\delta_0) - \bbW_2(\acute \zeta^\Ising,\delta_0)\Big|
            \le \bbW_2(\zeta^\Ising,\acute\zeta^\Ising)
            \le \epsilon\,.
        \end{align*}
Combining with Theorem~\ref{thm:SDE-smoothing-general} conclusion~\ref{it:SDE-smoothing-general.OUTPUT.qstar} 
gives
        \begin{align*}
            &\Big|\bbE [Y(q_*)^2] - q_*\Big| 
            \le \Big|\bbE [Y(q_*)^2] - q_0\Big| + |q_0 - q_*| \\
            &\qquad\le \Big(\bbE [Y(q_*)^2]^{1/2} + (q_0)^{1/2}\Big) \Big|\bbE [Y(q_*)^2]^{1/2} - (q_0)^{1/2}\Big| + \epsilon \\
            &\qquad\le (2q_0^{1/2} + \epsilon) \epsilon + \epsilon
            \le 3\epsilon + \epsilon^2
            \le \tilde\epsilon\,,
        \end{align*}
which implies condition~\ref{it:6.1}.
        \item Condition~\ref{it:7.1} is trivial.
        \item Theorem~\ref{thm:SDE-smoothing-general} conclusion~\ref{it:SDE-smoothing-general.OUTPUT.qstar} ensures $q_* \in [q_0,q_0+\epsilon] \subseteq [q_0,q_0+\tilde\epsilon]$, which implies condition~\ref{it:8.2}.
    \end{itemize}
    Finally, let
    \begin{align*}
        \acute \mu = \mu(\acute q_*, \acute b, \acute \sigma, \acute p, \acute \zeta) &\in \cM^{\sxi}(\alpha,q_0;L)\,, \\
        \mu = \mu(q_*,b,\sigma,p,\zeta) &\in \cM^{\IAMP}(\alpha,q_0;\tilde C,\tilde\epsilon)\,,
    \end{align*}
    be as defined by \eqref{e:control.problems.endpoint.measures}.
    Theorem~\ref{thm:SDE-smoothing-general} (see the assertion immediately following conclusion~\ref{it:SDE-smoothing-general.OUTPUT.qstar}) yields $\bbW_2(\acute \mu, \mu) \le \epsilon$.
    Since the above construction can be done for any $\acute \mu \in \cM^{\sxi}(\alpha,q_0;L)$, we conclude that
    \[
        \cB_\epsilon\Big(\cM^{\sxi}(\alpha,q_0;L)\Big)
        \subseteq \cB_{2\epsilon}\Big(\cM^{\IAMP}(\alpha,q_0;\tilde C(L,\epsilon),\tilde\epsilon(\epsilon))\Big)
        \subseteq \cB_{\tilde\epsilon}\Big(\cM^{\IAMP}(\alpha,q_0;\tilde C(L,\epsilon),\tilde\epsilon(\epsilon))\Big)\,.
    \]
    Taking a union over $L$ and an intersection over $\epsilon$ completes the proof (recalling again \eqref{e:measure.classes.simplification.without.L0}).
\end{proof}
\end{lem}

\subsection{Simple adjustments}
\label{ss:simple.adjustments}

In this subsection, we show how to make minor adjustments to the control processes to prove the following five inclusions: 
\begin{itemize} 
    \item Lemma~\ref{l:fair.in.tmp1} proves $\sM^{\si}(\alpha,q_0) \subseteq \sM^{\sii}(\alpha,q_0)$;
    \item Lemma~\ref{l:tmp1.in.tmp2} proves $\sM^{\sii}(\alpha,q_0) \subseteq \sM^{\siii}(\alpha,q_0)$;
    \item Lemma~\ref{l:tmp2.in.good} proves $\sM^{\siii}(\alpha,q_0) \subseteq \sM^{\siv}(\alpha,q_0)$;
    \item Lemma~\ref{l:good.in.tmp3} proves $\sM^{\siv}(\alpha,q_0) \subseteq \sM^{\sv}(\alpha,q_0)$;
    \item Lemma~\ref{l:prf2.in.prf3} proves $\sM^{\sx}(\alpha,q_0) \subseteq \sM^{\sxi}(\alpha,q_0)$.
\end{itemize}

\begin{lem}
\label{l:fair.in.tmp1}
For all $q_0 \in [0,1)$, we have
$\sM^{\si}(\alpha,q_0) \subseteq \sM^{\sii}(\alpha,q_0)$.

\begin{proof} We repeat the most relevant part of Definition~\ref{d:SDE.classification} for this proof:
	\begin{center}
	\begin{tabular}[h]{r|cccccccc}
	\MSRheader\\
	\hline
	\FAIRdefn\\
	\TMPONEdefn
	\end{tabular}
	\end{center}
We see that the $\si$ and $\sii$ classes differ only in the condition on the Ising initialization $Y(q_*)\sim\zeta^\Ising$: the $\si$ class satisfies the condition~\ref{it:7.1} (i.e., $\zeta^{\Ising} \in \cP_2(\R)$), while the $\sii$ class satisfies condition~\ref{it:7.2} (i.e., $\zeta^{\Ising} \in \cP([-1,1])$).

Therefore, let us suppose we have controls $(q_*,b,\sigma,w,p,\zeta,\zeta^\Ising)$ satisfying the $(\alpha,q_0,L,\epsilon)$-$\si$ conditions, which drive processes $X,Y$ as in \eqref{e:formal.sde.main}, \eqref{e:formal.sde.ising}. We will construct controls $(\tilde q_*,\tilde b,\tilde \sigma,\tilde w,\tilde p,\tilde \zeta,\tilde \zeta^\Ising)$ satisfying the $(\alpha,q_0,L,\tilde{\epsilon})$-$\sii$ conditions, for some $\tilde{\epsilon}$ that will tend to zero as $\epsilon\to0$. The latter controls will drive processes $\tilde X, \tilde Y$ that will approximate $X,Y$.

Let $(\tilde q_*,\tilde b,\tilde \sigma,\tilde p,\tilde \zeta) = (q_*,b,\sigma,p,\zeta)$.
    We then take 
    \[\tilde\zeta^{\Ising}
    =f_\# \zeta^{\Ising},\quad
        f(x) = \begin{cases}
            1 - \epsilon/x
            & 
            \textup{for $x \in (1,+\infty)$,} \\
            (1-\epsilon) x &
            \textup{for $[-1,1]$,} \\
            -1 + \epsilon/x &
            \textup{for $x \in (-\infty,1)$.}
        \end{cases}
    \]
Thus $\tilde\zeta^{\Ising} \in \cP([-1,1])$, so it satisfies condition~\ref{it:7.2}. We couple the initializations $Y(q_*)\sim\zeta^\Ising$ and 
	\[
	\tilde{Y}(q_*)
	=f(Y(q_*))
	\sim \tilde{\zeta}^\Ising\,.
	\]
Note that for $\epsilon\in(0,1)$, the function $f$ is an invertible mapping from $\bbR$ to $(-1,1)$.
    Recall from the definition \eqref{e:filt.X.Y} of $\cF_Y$ that we can write
    \[
        w_t = w_t\Big((W(s))_{s\in [q_*,t]},U',Y(q_*)\Big)\,,
    \]
where the $w_t$ on the right-hand side is a measurable function. We then define $\tilde w$ by setting
    \[
        \tilde w_t\Big((W(s))_{s\in [q_*,t]},U',\tilde Y(q_*)\Big)
        = w_t\Big(
        (W(s))_{s\in [q_*,t]},
        U',
        f^{-1}(\tilde Y(q_*))
        = Y(q_*) 
        \Big)\,.
    \]
Therefore $w=\tilde{w}$ as processes.
We let $\tilde X, \tilde Y$ be defined by \eqref{e:formal.sde.main} and \eqref{e:formal.sde.ising} with $(\tilde b,\tilde \sigma,\tilde w)$ in place of $(b,\sigma,w)$; and we couple them to
$X,Y$ by taking the initializations $\tilde X(q_*) = X(q_*)$ and $\tilde Y(q_*) = f(Y(q_*))$, as well as the same 
Brownian motions $B,W$. It is clear from these definitions that conditions \ref{it:1.2}, \ref{it:2.1}, \ref{it:3.2}, \ref{it:4.1}, \ref{it:8.1} remain satisfied for $(\tilde q_*,\tilde b,\tilde \sigma,\tilde w,\tilde p,\tilde \zeta,\tilde \zeta^\Ising)$; so it remains to verify that conditions \ref{it:5.1} and \ref{it:6.1} remain satisfied.

From the above definitions and \eqref{e:formal.sde.main}, it is clear that $X=\tilde{X}$ as processes. Since we already noted that $w=\tilde{w}$ as processes, it follows from \eqref{e:formal.sde.ising} that
	\beq\label{e:Y.tildeY.diffs.agree}
	Y(t)-Y(q_*)
	= \int_{q_*}^t w_u\,d W(u)
	=\tilde{Y}(t)-\tilde{Y}(q_*)
	\eeq
for all $t\in[q_*,1]$. Thus, conditions
\ref{it:5.1} and \ref{it:6.1} will essentially follow from showing that
$Y(q_*)$ and $\tilde{Y}(q_*)$ are close in $L^2$, so we turn to this next.

Let $g(x) = (|x|-1)_+^2$. Recall from Remark~\ref{r:true.mg} that $Y$ is an $L^2$-bounded continuous martingale. Since $g$ is nonnegative and convex, it follows that $g(Y(t))$ is a nonnegative submartingale. Thus
    \[
    \E g(Y(q_*))
    \le \E g(Y(1)) 
    \le \E\Big[(|Y(1)|-1)^2\Big]
    \le \epsilon\,,
    \]
where the last inequality follows from condition~\ref{it:5.1} for $Y(1)$. Thus we can bound
    \begin{align}\nonumber
        &\E\Big[(Y(q_*) - \tilde Y(q_*))^2\Big]
        = \E
        \Big[(Y(q_*) - f(Y(q_*)))^2
        \Big] \\
        &\qquad\le \E\Big[
        \ind\{|Y(q_*)| \le 1\}\epsilon^2 Y(q_*)^2\Big]
        + \E\Big[\ind\{|Y(q_*)| > 1\}(|Y(q_*)| - (1-\epsilon))^2\Big]
        \nonumber \\
        &\qquad\le \epsilon^2 + \E\Big[
            \ind\{|Y(q_*)| > 1\}(2(|Y(q_*)| - 1)^2 + 2\epsilon^2)
        \Big] 
        \le 3\epsilon^2 + 2\E g(Y(q_*))
        \le 4\epsilon\,,
    \label{eq:fair.implies.tmp1.wass}
    \end{align}
where the last inequality holds for $\epsilon$ small enough. It follows that
    \begin{align*}
        &\Big|\E[\tilde Y(q_*)^2] - \E[Y(q_*)^2]\Big|
        \le \E\bigg[
            \Big|
            \tilde Y(q_*) - Y(q_*)\Big|
            \Big|\tilde Y(q_*)
            +Y(q_*)\Big|
        \bigg] \\
        &\qquad \le
        \E \bigg[
            \Big(
            \tilde Y(q_*) - Y(q_*)
            \Big)^2
        \bigg]^{1/2}
        \E \bigg[
            \Big( \tilde Y(q_*)
            + Y(q_*)\Big)^2
        \bigg]^{1/2}
    \end{align*}
In the last line above, the second factor is at most $O(1)$ by condition~\ref{it:6.1} for $Y(q_*)$, together with the fact that $|\tilde Y(q_*)| \le |Y(q_*)|$ almost surely. Meanwhile, the first factor is at most $2\epsilon^{1/2}$ by \eqref{eq:fair.implies.tmp1.wass}.
    This proves condition~\ref{it:6.1} for $\tilde Y(q_*)$ by suitably adjusting $\epsilon$.
    Lastly, for the endpoint condition~\ref{it:5.1}, we bound
    \begin{align}
	\label{eq:endpoint-condition-stability}
        &\bigg|\bbE\Big[(|\tilde Y(1)| - 1)^2\Big]
        - \bbE\Big[(|Y(1)| - 1)^2\Big]\bigg|
        \le \bbE\bigg[
        \Big|\tilde Y(1) - Y(1) \Big|
        \Big(|\tilde Y(1)| + |Y(1)| + 2\Big)
        \bigg] \\
        \notag
        &\qquad\le \bbE\Big[
            (\tilde Y(1) - Y(1))^2\Big]^{1/2}
        \bbE\Big[(|\tilde Y(1)| + |Y(1)| + 2)^2\Big]^{1/2}\,.
    \end{align}
In the last line above, the second factor is at most $O(1)$ by  Remark~\ref{r:true.mg}, while the first factor is bounded by $2\epsilon^{1/2}$ by combining
\eqref{e:Y.tildeY.diffs.agree} with the estimate 
\eqref{eq:fair.implies.tmp1.wass}. This proves condition~\ref{it:5.1} for $\tilde Y(1)$, again after suitably adjusting $\epsilon$.

Finally, since we already noted that $X = \tilde X$ as processes, we have $\Law(\tilde X(1)) = \Law(X(1))$. Recalling the notations of Definition~\ref{d:measure.classes}, we have shown that for some $\tilde{\epsilon}$ tending to zero with $\epsilon$, we have
    \[
        \cM^{\si}(\alpha,q_0;L,\epsilon)
        \subseteq
        \cM^{\sii}(\alpha,q_0;L,\tilde \epsilon)\,.
    \]
It follows that $\sM^{\si}(\alpha,q_0) \subseteq  \sM^{\sii}(\alpha,q_0)$, concluding the proof.
\end{proof}
\end{lem}

\begin{lem}\label{l:tmp1.in.tmp2}
For all $q_0 \in [0,1)$, we have 
$\sM^{\sii}(\alpha,q_0) \subseteq \sM^{\siii}(\alpha,q_0)$.

\begin{proof}
We repeat the most relevant part of Definition~\ref{d:SDE.classification} for this proof:
	\begin{center}
	\begin{tabular}[h]{r|cccccccc}
	\MSRheader\\
	\hline
	\TMPONEdefn\\
	\TMPTWOdefn
	\end{tabular}
	\end{center}
We see that
 the $\sii$ and $\siii$ classes differ only in the condition on $q_*$: the $\sii$ class satisfies condition~\ref{it:8.1} (i.e., $q_* \in [q_0,q_0+\epsilon]$), while the $\siii$ class satisfies condition~\ref{it:8.2} (i.e., $q_*=q_0$).

Thus, let us suppose we have controls $(q_*,b,\sigma,w,p,\zeta,\zeta^\Ising)$ satisfying the $(\alpha,q_0,L,\epsilon)$-$\sii$ conditions, which drive processes $X,Y$ as in 
\eqref{e:formal.sde.main}, \eqref{e:formal.sde.ising}. We will construct controls $(\tilde q_*,\tilde b,\tilde \sigma,\tilde w,\tilde p,\tilde \zeta,\tilde \zeta^\Ising)$ satisfying the $(\alpha,q_0,L,\tilde \epsilon)$-$\siii$ conditions, for some $\tilde{\epsilon}$ that will tend to zero as $\epsilon\to0$. The latter controls will drive processes $\tilde X, \tilde Y$ that will approximate $X,Y$.

Let $\tilde q_* = q_0$, so that condition~\ref{it:8.2} is satisfied, and $(\tilde p,\tilde \zeta,\tilde \zeta^\Ising) = (p,\zeta,\zeta^\Ising)$. We will couple $\tilde{X}(q_0)=X(q_*)$ and 
$\tilde{Y}(q_0)=Y(q_*)$. 
From the
definition \eqref{e:filt.X.Y} of $\cF_X$ and $\cF_Y$, we can write
\begin{align*}
        b_t &= b_t( (B(s))_{s\in[q_*,t]},U,X(q_*)), \\
        \sigma_t &= \sigma_t( (B(s))_{s\in[q_*,t]},U,X(q_*)), \\
        w_t &= w_t( (W(s))_{s\in[q_*,t]},U',Y(q_*)),
    \end{align*}
where the $\sigma_t,b_t,w_t$ on the right-hand sides are measurable functions.
We then let $(\tilde{b},\tilde{\sigma},\tilde{w})$ be defined by
$(\tilde b_t,\tilde \sigma_t,\tilde w_t) = (0,1,1)$ for $t\in [q_0,q_*)$, and
    \begin{align*}
        \tilde b_t( (B(s))_{s\in[q_0,t]},U,
        \tilde{X}
        (q_0)
        = X(q_*)
        ) &= b_t( (B(s))_{s\in[q_*,t]},U,
        X(q_*))
        \,, \\
        \tilde \sigma_t( (B(s))_{s\in[q_0,t]},U,
        \tilde{X}
        (q_0)
        = X(q_*))
        &= \sigma_t( (B(s))_{s\in[q_*,t]},U,
        X(q_*))
        \,, \\
        \tilde w_t( (W(s))_{s\in[q_0,t]},U',
        \tilde{Y}(q_0)
        = X(q_*)) 
        &= w_t( (W(s))_{s\in[q_*,t]},U',Y(q_*))\,.
    \end{align*}
for $t\in [q_*,1]$. 
We then couple the resulting $\tilde X, \tilde Y$ to $X,Y$ using the initializations $\tilde X(q_0) = X(q_*)$ and $\tilde Y(q_0) = Y(q_*)$ as noted above, and the same Brownian motions $B,W$. From this it is clear that $(\tilde b_t,\tilde \sigma_t,\tilde w_t) = (b_t,\sigma_t,w_t)$ for all $t \in [q_*,1]$. It is clear from these definitions that conditions \ref{it:1.2}, \ref{it:2.1}, \ref{it:3.2}, \ref{it:6.1}, \ref{it:7.2} remain satisfied for the modified controls $(\tilde q_*,\tilde b,\tilde \sigma,\tilde w,\tilde p,\tilde \zeta,\tilde \zeta^\Ising)$; and we turn to verifying that conditions \ref{it:4.1} and \ref{it:5.1} also remain satisfied.

We first verify that the modified controls $(\tilde{b},\tilde{\sigma},\tilde{w},\tilde{p})$ also satisfy \ref{it:4.1}. Indeed, for the modified controls, for $t\in [q_0,q_*)$ we have $\Cost(t) = 0$ and $\budget(t;\alpha) = 1/\alpha$, so $(\Cost(t) - \budget(t;\alpha))_+ = 0$. For $t\in [q_0,q_*)$, the values of $\Cost(t)$ and $\budget(t;\alpha)$ are the same for $(\tilde b,\tilde \sigma,\tilde w,\tilde{p})$ as for $(b,\sigma,w,p)$. This verifies that the modified controls also satisfy condition~\ref{it:4.1}.

We next verify that $\tilde{Y}(1)$ satisfies condition~\ref{it:5.1} (with $\epsilon$ suitably adjusted). From the coupling described above, we have
\beq\label{e:tmp12.tildeY.Y.diffs.agree}
        \tilde Y(1) - Y(1)
        = \tilde Y(q_*) - Y(q_*)
        = \tilde Y(q_*) - \tilde Y(q_0)
        = \int_{q_0}^{q_*} 
        \,dW(t)\,,
\eeq
where the last step uses that $\tilde{\sigma}_t=1$ for $t\in[q_0,q_*]$.
This implies
    \beq\label{e:tmp12.tildeY.Y.L2.error}
	\E\Big[(\tilde Y(1) - Y(1))^2\Big]
	=\E\Big[(\tilde Y(q_*) - \tilde Y(q_0))^2\Big] 
	= q_* - q_0 \le \epsilon\,,
    \eeq
where the last step uses that $q_*$ satisfies condition~\ref{it:8.1}. Now recall from the proof of Lemma~\ref{l:fair.in.tmp1} that
\eqref{e:Y.tildeY.diffs.agree} and
\eqref{eq:fair.implies.tmp1.wass} together imply \eqref{eq:endpoint-condition-stability}: by the same argument, \eqref{e:tmp12.tildeY.Y.diffs.agree} and \eqref{e:tmp12.tildeY.Y.L2.error} together imply that $\tilde{Y}(1)$
satisfies condition~\ref{it:5.1}, after suitably adjusting $\epsilon$.

Similarly to \eqref{e:tmp12.tildeY.Y.diffs.agree}, we have
	\[
	\tilde{X}(1)-X(1)
	=\tilde{X}(q_*)-X(q_*)
	=\tilde{X}(q_*)-\tilde{X}(q_0)
	=\int_{q_0}^{q_*} s(t)\,dB(t)\,.
	\]
Recalling that $s(t) \equiv [(tp)'(t)]^{1/2}$, we obtain
	\[\bbW_2\Big(\Law(\tilde X(1)),\Law(X(1))\Big)^2
	=\E\bigg[
	\Big(\tilde{X}(1)-X(1)\Big)^2\bigg]
	\le \int_{q_0}^{q_*} s(t)^2\,dt
	= q_* p(q_*)
	\le \epsilon\,,
	\]
having used condition~\ref{it:3.2}.
Thus, for some $\tilde{\epsilon}$ tending to zero with $\epsilon$, 
    \[
        \cB_\epsilon\Big(
        \cM^{\sii}
        (\alpha,q_0;L,\epsilon
        )\Big)
        \subseteq 
        \cB_{\tilde\epsilon}
        \Big(\cM^{\siii}
        (\alpha,q_0;L,\tilde \epsilon)
        \Big)\,.
    \]
Recalling the notations of Definition~\ref{d:measure.classes}, the claim follows.
\end{proof}
\end{lem}

\begin{lem}\label{l:tmp2.in.good}
For all $q_0 \in [0,1)$, we have 
$\sM^{\siii}(\alpha,q_0) \subseteq \sM^{\siv}(\alpha,q_0)$.

\begin{proof}
We repeat the most relevant part of Definition~\ref{d:SDE.classification} for this proof:
	\begin{center}
	\begin{tabular}[h]{r|cccccccc}
	\MSRheader\\
	\hline
	\TMPTWOdefn\\
	\GOODdefn
	\end{tabular}
	\end{center}
We see that the $\siii$ and $\siv$ classes 
differ only in the conditions on the $L^2$ norms of the initializations $X(q_*)$ and $Y(q_*)$:
the $\siii$ class satisfies  condition~\ref{it:6.1}, while the $\siv$ class satisfies condition~\ref{it:6.2}. 

Therefore, let us suppose we have controls $({{q_*=q_0}},b,\sigma,w,p,\zeta,\zeta^\Ising)$ satisfying the $(\alpha,q_0,L,\epsilon)$-$\siii$ conditions, which drive processes $X,Y$ as in \eqref{e:formal.sde.main}, \eqref{e:formal.sde.ising}. We will construct controls $(\tilde q_*=q_0,\tilde b,\tilde \sigma,\tilde w,\tilde p,\tilde \zeta,\tilde \zeta^\Ising)$ satisfying the $(\alpha,q_0,L,\tilde \epsilon)$-$\siv$ conditions, for $\tilde{\epsilon}$ tending to zero as $\epsilon\to0$. The latter controls will drive processes $\tilde X, \tilde Y$ that will approximate $X,Y$. 

We take $\tilde p = p$ and $\tilde{X}(q_0)=0\sim\tilde \zeta = \delta_0$. We will define a measure $\tilde{\zeta}^\Ising$ such that $\tilde{Y}(q_0)\sim\tilde{\zeta}^\Ising$ satisfies the required condition~\ref{it:6.2}. We will then define $(\tilde{b},\tilde{\sigma},\tilde{w})$ to simulate the original controls $(b,\sigma,w)$: from the definition \eqref{e:filt.X.Y} of $\cF_X$ and $\cF_Y$, we can write
    \begin{align*}
        \sigma_t &= \sigma_t( (B(s))_{s\in[q_0,t]},U,X(q_0))
        \,, \\
        b_t &= b_t( (B(s))_{s\in[q_0,t]},U,X(q_0))
        \,, \\
        w_t &= w_t( (W(s))_{s\in[q_0,t]},U',Y(q_0))\,,
    \end{align*}
where the $\sigma_t,b_t,w_t$ on the right-hand sides are measurable functions. Then sample independent random variables $\hat{X}(q_0)\sim\zeta$ and $\tilde{U}\sim\unif([0,1])$. 
The laws of $\tilde{U}$ and $(U,\hat{X}(q_0))$ are atom-free Borel measures on Polish spaces, so by Maharam's theorem \cite[Chapter 33]{fremlin2000measure}, there exists a measure-preserving bijection $\rho$ such that $\rho(\tilde{U}) = (U,\hat{X}(q_0))$. We use this to define
     \begin{align*}
        \tilde\sigma_t( (B(s))_{s\in[q_0,t]},\tilde U, \tilde X(q_0))
        &= \sigma_t( (B(s))_{s\in[q_0,t]} ,U,\hat X(q_0))\,, \\
        \tilde b_t(
            (B(s))_{s\in[q_0,t]},\tilde U, \tilde X(q_0))
        &= b_t(
        (B(s))_{s\in[q_0,t]},U,\hat X(q_0))\,,
    \end{align*}
where $(U,\hat{X}(q_0))=\rho(\tilde{U})$. Note that as a result $\tilde{b}_t$ and $\tilde{\sigma}_t$ do not depend on $\tilde{X}(q_0)$, and in any case $\tilde{X}(q_0)=0$ by definition.

To define $\tilde{\zeta}^\Ising$ and $\tilde{w}$, we separate into cases according to $q_0$. If $q_0=0$, we let $\tilde Y(q_0)=0 \sim \tilde{\zeta}^\Ising=\delta_0$, and define $\tilde{w}$ analogously as $\tilde{b},\tilde{\sigma}$ above. Otherwise, suppose $q_0\in(0,1)$.
Since we take $\epsilon \downarrow 0$ in the definition of $\cM^\star(\alpha,q_0)$ (see Definition~\ref{d:measure.classes}), we may assume $\epsilon$ to be small relative to $q_0$, so that condition~\ref{it:6.1} for $Y(q_0)$ implies 
 $\E[Y(q_0)^2]\in(0,1)$. Consider the following family of functions $f_a : [-1,1] \to [-1,1]$ parametrized by $a \in [0,2]$: for $a \in (0,1]$, let $f_a(x) = ax$, and for $a\in (1,2]$, let
\[
    f_a(x) = \begin{cases}
        1 - (2-a)(1-x) & x \ge 0, \\
        -1 + (2-a)(1+x) & x < 0.
    \end{cases}
\]
Note that $f_1(x)=x$; $f_2(x)=\pm1$ gives the sign of $x$; and as a function of $a\in[1,2]$, the value of $f_a(x)$ interpolates linearly between $f_1(x)$ and $f_2(x)$. Thus, for any fixed $x\in [-1,1]$, $f_a(x)^2$ is a continuous, nondecreasing function of $a\in[0,2]$.  Moreover, it is strictly increasing over $a\in[0,1]$ provided $x\ne0$, and it is strictly increasing over $a\in[1,2]$ provided $x\notin\{-1,+1\}$. Condition~\ref{it:6.1} for $Y(q_0)$ implies that $Y(q_0)$ cannot be almost surely zero, and also cannot be almost surely $\pm1$. It follows that
\beq
    \label{eq:tmp2-implies-good-g}
    g(a)
    \equiv \E\Big[f_a(Y(q_0))^2\Big]
\eeq
is a continuous, strictly increasing function of $a\in[0,2]$, with $g(0)=0$ and $g(2)=1$. It follows that there is a unique $a_*\in(0,2)$ such that $g(a_*) = q_0$. Let $f_*\equiv f_{a_*}$, and
	\[
	\tilde{Y}(q_0)
	\equiv f_*(Y(q_0))\sim
	\tilde \zeta^\Ising
	\equiv f_\# \zeta^\Ising\,.
	\]
The choice $g(a_*)=q_0$ then implies
$\E[\tilde{Y}(q_0)^2]=q_0$, so that 
$\tilde{Y}(q_0)$ satisfies condition~\ref{it:6.2} as required. Since $f_*$ is a bijection, we can use it to define
\[
    \tilde w_t\Big( (W(s))_{s\in [q_0,t]},U',\tilde Y(q_0)\Big)
    = w_t\Big( (W(s))_{s\in [q_0,t]}, U', (f_*)^{-1}(\tilde Y(q_0))
    = Y(q_0) 
    \Big)\,.
\]
Let $\tilde{X}$ and $\tilde{Y}$ be defined by \eqref{e:formal.sde.main} and \eqref{e:formal.sde.ising} with $(\tilde{b},\tilde{\sigma},\tilde{w})$ in place of $(b,\sigma,w)$, with initializations $\tilde X(q_0) = 0$, $\tilde Y(q_0) = f_*(Y(q_0))$, and with the same Brownian motions $B,W$. From the above, it is clear that we have
    \[
    \Big(\tilde{b}_t,\tilde{\sigma}_t,\tilde{w}_t,
        \tilde{X}(t),  \tilde{Y}(t)-\tilde{Y}(q_0)
        \Big)
    =
    \Big(b_t,\sigma_t,w_t,X(t)-X(q_0),Y(t)-Y(q_0)\Big)\,,
    \]
and that the modified controls $(\tilde q_*,\tilde b,\tilde \sigma,\tilde w,\tilde p,\tilde \zeta,\tilde \zeta^\Ising)$ continue to satisfy conditions \ref{it:1.2}, \ref{it:2.1}, \ref{it:3.2}, \ref{it:4.1},  \ref{it:7.2}, and \ref{it:8.2}. For the endpoint condition~\ref{it:5.1}, we bound
\begin{align}\nonumber
    \bbE\Big[ (\tilde Y(1) - Y(1))^2\Big]
    &= \bbE\Big[ (\tilde Y(q_0) - Y(q_0))^2\Big]
    \stackrel{(*)}{=}
    \bbE\Big[ (|\tilde Y(q_0)| - |Y(q_0)|)^2\Big] \\
    &\le \bbE\Big[ \big|\tilde Y(q_0)^2 - Y(q_0)^2\big|\Big]
    \stackrel{(**)}{=} \Big|\bbE[\tilde Y(q_0)^2] - \bbE[Y(q_0)^2]\Big|
    \le \epsilon
    \label{e:tmp2.good.tildeY.Y.L2.error}
\end{align}
Here, the step $(*)$ is because $Y(q_0)$ and $\tilde Y(q_0)$ cannot take opposing signs. The step $(**)$ is because $\tilde Y(q_0)^2 - Y(q_0)^2$ either is almost surely nonnegative (if $a_*\in[1,2)$) or almost surely nonpositive (if $a_*\in(0,1]$). Then, similarly as in
\eqref{eq:endpoint-condition-stability} in the proof of Lemma~\ref{l:fair.in.tmp1}, we can bound
    \begin{align*}
        &\bigg|\bbE\bigg[
        \Big(|\tilde Y(1)| - 1\Big)^2\bigg]
        - \bbE\bigg[
        \Big(|Y(1)| - 1
        \Big)^2\bigg]\bigg|
        \le \bbE\bigg[
        \Big|\tilde Y(1) - Y(1) \Big|
        \Big(|\tilde Y(1)| + |Y(1)| + 2\Big)
        \bigg] \\
        &\qquad\le \bbE\bigg[
            \Big(\tilde Y(1) - Y(1)\Big)^2\bigg]^{1/2}
        \bbE\bigg[
        \Big(
        |\tilde Y(1)| + |Y(1)| + 2\Big)^2\bigg]^{1/2}.
    \end{align*}
In the last line, the second factor is at most $O(1)$ by Remark~\ref{r:true.mg}, while the first factor is at most $\epsilon^{1/2}$ by \eqref{e:tmp2.good.tildeY.Y.L2.error}. Thus $\tilde{Y}(1)$ also satisfies condition~\ref{it:5.1}, after adjusting $\epsilon$ as needed.

Finally, we recall from \eqref{e:control.problems.limiting.regime} that $\epsilon\ll\iota$, from which it follows that
\[
    \bbW_2\Big(\Law(\tilde X_1),\Law(X_1)\Big)^2
    = \bbE\bigg[ \Big(\tilde X(1) - X(1) \Big)^2\bigg]
        = \bbE[X(q_0)^2]
        \le \epsilon\,.
    \]
Thus for some $\tilde \epsilon$ tending to zero with $\epsilon$, we have shown  that
\[
    \cB_\epsilon\Big(\cM^{\siii}(\alpha,q_0;L,\epsilon,L_0
    )\Big)
    \subseteq 
    \cB_{\tilde\epsilon}
    \Big(\cM^{\siv}
    (\alpha,q_0;L,\tilde \epsilon,L_0
    )\Big),
\]
Recalling the notations of Definition~\ref{d:measure.classes}, the claim follows.
\end{proof}
\end{lem}

\begin{lem}\label{l:good.in.tmp3}
For all $q_0 \in [0,1)$, we have 
$\sM^{\siv}(\alpha,q_0) \subseteq \sM^{\sv}(\alpha,q_0)$.

\begin{proof}
We repeat the most relevant part of Definition~\ref{d:SDE.classification} for this proof:
	\begin{center}
	\begin{tabular}[h]{r|cccccccc}
	\MSRheader\\
	\hline
	\GOODdefn\\
	\TMPTHRdefn
	\end{tabular}
	\end{center}
We see that the
$\siv$ and $\sv$ classes differ only in the condition on 
$\E[(w_t)^2]$: the $\siv$ class satisfies condition~\ref{it:2.1} (i.e., $\E[(w_t)^2] \in [1-\epsilon,1+\epsilon]$), while the $\sv$ class satisfies condition~\ref{it:2.2} (i.e., $\E[(w_t)^2] = 1$). 

Therefore, let us suppose we have controls $(q_* = q_0,b,\sigma,w,p,\zeta,\zeta^\Ising)$ satisfying the $(\alpha,q_0,L,\epsilon)$-{{$\siv$}} conditions. We then define
$(\tilde{q}_*,\tilde{b},\tilde{\sigma},\tilde{p},\tilde{\zeta},\tilde{\zeta}^\Ising)
=(q_* ,b,\sigma,p,\zeta,\zeta^\Ising)$, and
	\[\tilde w_t=\frac{w_t}{\bbE[(w_t)^2]^{1/2}}\,,\]
so that $\tilde{w}$ satisfies
condition~\ref{it:2.2}. We let $\tilde{X}=X$, and let
$\tilde{Y}$ be given by \eqref{e:formal.sde.ising} with $\tilde{w}$ in place of $w$, but with
the same initialization $\tilde{Y}(q_0)=Y(q_0)$ and the same Brownian motion $W$. 
 It is clear from these definitions that the modified controls continue to satisfy
conditions~\ref{it:1.2}, \ref{it:3.2}, \ref{it:6.2}, \ref{it:7.2}, and \ref{it:8.2}. 
For condition~\ref{it:4.1}, we note that for the two choices of $w$, the difference in the resulting values of 
$\budget(t;\alpha)$ is
	\[
    \frac{|\E[w_t]^2 - \E[\tilde w_t]^2|}{\alpha} = \frac{\E[w_t]^2}{\alpha} \bigg|1 - \frac{1}{\bbE[(w_t)^2]}
    \bigg|
    \le \frac{1+\epsilon}{\alpha} \cdot \frac{\epsilon}{1-\epsilon}\,.
\]
Since this tends to zero as $\epsilon \to 0$, we see that $\tilde{w}$ also satisfies
condition~\ref{it:4.1}, provided we adjust $\epsilon$. Finally, it follows from \eqref{e:formal.sde.ising} that
    \[
        \bbE\Big[ (Y(1)-\tilde Y(1))^2\Big]
        \leq
        \int_{q_0}^1
        \bbE[(\tilde w_t-w_t)^2]
        \,dt
        =
        \int_{q_0}^1
        \Big(
        \bbE [(w_t)^2]^{1/2} - 1
        \Big)^2
        \,dt
        \le 
        \Big(1 - (1-\epsilon)^{1/2}
        \Big)^2\,.
    \]
By the same estimate as in \eqref{eq:endpoint-condition-stability} in the proof of Lemma~\ref{l:fair.in.tmp1}, $\tilde{Y}(1)$ also satisfies condition~\ref{it:5.1}, again provided we adjust $\epsilon$. 
Thus we have shown that for some $\tilde{\epsilon}$ tending to zero with $\epsilon$, we have
\[
    \cB_\epsilon(\cM^{\siv}(\alpha,q_0;L,\epsilon,L_0))
    \subseteq \cB_{\tilde\epsilon}(\cM^{\sv}(\alpha,q_0;L,\tilde \epsilon,L_0)),
\]
Recalling the notations of Definition~\ref{d:measure.classes}, the claim follows.\end{proof}\end{lem}

\begin{lem}\label{l:prf2.in.prf3}
For all $q_0 \in [0,1)$, we have $\sM^{\sx}(\alpha,q_0) \subseteq \sM^{\sxi}(\alpha,q_0)$. 

\begin{proof}
We repeat the most relevant part of Definition~\ref{d:SDE.classification} for this proof:
	\begin{center}
	\begin{tabular}[h]{r|cccccccc}
	\MSRheader\\
	\hline
	\PRFTWOdefn \\
	\PRFTHRdefn
	\end{tabular}
	\end{center}
Consider any controls
	\[(q_*=q_0, b, \sigma, w, p, \zeta=\delta_0, \zeta^\Ising) 
	\in \Adm^{\sx}(\alpha,q_0;L,\epsilon)\,.\]
	We will construct controls $(\tilde q_* = q_0, \tilde b, \tilde \sigma, \tilde w, \tilde p, \tilde \zeta = \delta_0, \tilde \zeta^\Ising) \in \Adm^{\sxi}(\alpha,q_0;L,\epsilon)$ so that the endpoint measures $\mu = \mu(q_*,b,\sigma,p,\zeta)$ and $\tilde\mu = \mu(\tilde q_*, \tilde b, \tilde \sigma, \tilde p, \tilde\zeta)$ coincide.
    This shows
    \[
        \Adm^{\sx}(\alpha,q_0;L,\epsilon)
        \subseteq 
        \Adm^{\sxi}(\alpha,q_0;L,\epsilon)
    \]
for all $L,\epsilon$, which implies the result by recalling Definition~\ref{d:measure.classes}. In fact, we will take \[(\tilde q_*, \tilde b, \tilde \sigma, \tilde p, \tilde\zeta) = (q_*,b,\sigma,p,\zeta)\,,\] so that $\mu = \tilde \mu$ is immediate.
    We will construct controls $(\tilde w, \tilde \zeta^\Ising)$, as well as $\tilde U' \sim \unif([0,1])$ and a Brownian motion $(\tilde W(t))_{t\in [q_0,1]}$, so that for $Y$ and $\tilde Y$ the processes defined by \eqref{e:formal.sde.ising} with (respectively) inputs 
    \begin{align*}
	&(q_*=q_0,w,\zeta^\Ising,
	(W(t))_{t\in [q_0,1]},U')\,,\\
	&(\tilde q_* = q_0,\tilde w,\tilde\zeta^\Ising,(\tilde W(t))_{t\in [q_0,1]},\tilde U')\,,\end{align*}
we have the equality of processes
    \beq
        \label{e:prf2.in.prf3.symmetrization}
        (\tilde Y(t))_{t\in [q_0,1]}
        = (\xi \cdot Y(t))_{t\in [q_0,1]}\,,
    \eeq
    for a random variable $\xi \sim \unif(\{\pm 1\})$ independent of the process $(Y(t))_{t\in [q_0,1]}$. To this end, we construct $\tilde\zeta^\Ising$ as the symmetrization of $\zeta^\Ising$:
    \[ 
    \tilde{\zeta}^\Ising(A)
    \equiv
    \frac{\zeta^\Ising(A)
    +\zeta^\Ising(-A)}{2}
	\] 
for all measurable $A\subseteq \R$. We then sample $\tilde{Y}(q_0)\in\tilde{\zeta}^\Ising$, so that clearly $|\tilde{Y}(q_0)|$ is equidistributed as $|Y(q_0)|$. This symmetrization can be reversed, in the following sense. Define $f : [0,1] \to [0,1]$ to be the conditional probability of the event $\{Y(q_0)\geq 0\}$ conditional on $|Y(q_0)|$, i.e.,
    \[
        f(x) = \bbP\Big(Y(q_0) \ge 0
        	\,\Big| \,|Y(q_0)| =x\Big).
    \]
For $x,U'' \in [0,1]$ define
    \[g_f(x,U'') \equiv \begin{cases}
            x & U'' \le f(x), \\
            -x & U'' > f(x).
        \end{cases}
    \]
Having sampled $\tilde{Y}(q_0)\in\tilde{\zeta}^\Ising$, we next define
	\[\hat{Y}(q_0) 
	\equiv g_f(|\tilde Y(q_0)|,U'')\,.
	\]
We then verify that $\hat{Y}(q_0)$ is equidistributed as $Y(q_0)\sim\zeta$: for $x>0$,
	\begin{align*}
	&\P(\hat{Y}(q_0) \in dx)
	= \P(|\tilde{Y}(q_0)| \in dx)
		\bbP\Big(Y(q_0) \ge 0
        	\,\Big| \,|Y(q_0)| =x\Big)\\
	&\qquad=\P(|Y(q_0)| \in dx)
		\bbP\Big(Y(q_0) \ge 0
        	\,\Big| \,|Y(q_0)| =x\Big)
	=\P(Y(q_0) \in dx)\,,
	\end{align*}
and similarly for $x\le0$. It follows that we can set $Y(q_0) = \hat{Y}(q_0) \sim\zeta$. We then define
	\[
	\xi = \sign\Big(
		\tilde Y(q_0)Y(q_0)\Big)\,,
	\]
and we claim that $\xi$ is a symmetric random sign that is independent of $Y(q_0)$: for $x>0$,
	\begin{align*}
	&\P\Big(\xi = +1, Y(q_0) \in dx\Big)
	=\P(\tilde{Y}(q_0)\in dx)
	\P\Big(Y(q_0)\ge0 
	\,\Big|\ 
	|Y(q_0)|=x\Big)\\
	&\qquad= \frac{\P(|Y(q_0)|\in dx)}{2}
	\P\Big(Y(q_0)\ge0 
	\,\Big|\ 
	|Y(q_0)|=x\Big)
	= \frac{\P(Y(q_0) \in dx)}{2}\,.
	\end{align*} 
Next let $(W(s))_{s\in [q_0,t]}$ be a Brownian motion independent of $U'',\tilde Y(q_0)$ (and thus $Y(q_0)$).
    Recalling the definition \eqref{e:filt.X.Y} of $\cF_Y$, $w_t$ is a measurable function
    \[
        w_t = w_t\Big( (W(s))_{s\in [q_0,t]},U',Y(q_0)\Big)\,,
    \]
    where $U' \sim \unif([0,1])$.
    We construct a Brownian motion $(\tilde W(s))_{s\in [q_0,t]}$ by $\tilde W(s) = \xi \cdot W(s)$.
    Note that $\tilde{W}$ is independent of $\xi$, and thus of $\tilde Y(q_0)$.
    We next construct $\tilde w$, which will be measurable with respect to 
    \[
        \tilde\cF_Y = \sigma\Big(
            (\tilde W(s) : q_0 \le s\le t),
            \tilde U',\tilde Y(q_0)
        \Big)
    \]
    for some $\tilde U' \sim \unif([0,1])$ independent of $\tilde Y(q_0)$ and $\tilde W$.
The laws of $\tilde U'$ and $(U',U'')$ are atom-free Borel measures on Polish spaces, so again by Maharam's theorem \cite[Chapter 33]{fremlin2000measure} there exists a measure-preserving bijection $\rho:[0,1]\to [0,1]^2$ such that $\rho(\tilde U') = (U',U'')$. We use this to define
    \[
        \tilde w_t\Big( (\tilde W(s))_{s\in [q_0,t]},\tilde U',\tilde Y(q_0)\Big)
        = 
        w_t \Big( (\xi \cdot W(s))_{s\in [q_0,t]},U', Y(q_0)\Big)\,.
    \]
    Note that the right-hand side is a function of the variables on the left-hand side because the random variables $Y(q_0) = g_f(|\tilde Y(q_0)|,U'')$ and $\xi = \sign(\tilde Y(q_0)Y(q_0))$ are both $(U'',\tilde Y(q_0))$-measurable.
    Then, for all $t\in [q_0,1]$, we have
    \[
        \tilde Y(t) = \tilde Y(q_0) + \int_{q_0}^t \tilde w_t\,d\tilde W(s)
        = \xi \cdot \bigg(Y(q_0) + \int_{q_0}^t w_t\,dW(s)\bigg)
        = \xi \cdot Y(t)\,.
    \]
This proves the claimed coupling \eqref{e:prf2.in.prf3.symmetrization}.

    Finally, we verify that our constructed controls satisfy $(\tilde q_*, \tilde b, \tilde \sigma, \tilde w, \tilde p, \tilde \zeta, \tilde \zeta^\Ising) \in \Adm^{\sxi}(\alpha,q_0;L,\epsilon)$.
    Since $(q_*, b, \sigma, w, p, \zeta, \zeta^\Ising) \in \Adm^{\sx}(\alpha,q_0;L,\epsilon)$, conditions \ref{it:3.1} and \ref{it:8.2} are immediate.
    Since $\tilde w_t = w_t$, conditions \ref{it:2.2} and \ref{it:4.3} also follow.
    Finally \eqref{e:prf2.in.prf3.symmetrization} implies conditions \ref{it:1.2even}, \ref{it:5.2}, \ref{it:6.2}, and \ref{it:7.2}. The claim follows.
\end{proof}
\end{lem}

\subsection{Continuity in the budget}
\label{ss:continuity.in.budget}

In this subsection we prove the inclusion $\cM^{\sv}(\alpha,q_0) \subseteq \cM^{\svi}(\alpha,q_0)$ (see Proposition~\ref{p:continuity.in.auxilliary.budget}), along with  Corollary~\ref{cor:feasible.distributions.continuity}
for the map $\alpha \mapsto  \ocM^{\Ising}(\alpha)$. For these arguments, it will be convenient to make the following definition (compare with conditions~\ref{it:4.1}--\ref{it:4.3} from Definition~\ref{d:SDE.conds}):

\begin{dfn}[generalized budget constraint]\label{d:X.adm}
Recall the set $\incr([q_0,1];[0,1])$ defined in Definition~\ref{d:incr.p}.
For $p\in \incr([q_0,1];[0,1])$ and a generalized ``budget'' function $f \in L^1([q_0,1];[0,+\infty))$, let
$\Adm(q_0;p,f)$
denote the set of progressively measurable control processes $(b,\sigma)$, with $\sigma \ge 0$,
satisfying
    \beq\label{e:budget-f}
    \bbE\bigg[
    (b_t)^2 + \frac{tp'(t)}{p(t)} (\sigma_t-1)^2
    \bigg]
    \le f(t)\eeq
for all $t \in [q_0,1]$. Then let
	\[
    \cM(q_0;p,f) = \bigg\{
        \mu(q_* = q_0,b,\sigma,p,\zeta = \delta_0) : (b,\sigma)
        \in \Adm(q_0;p,f)
    \bigg\}\,,
\]
where we recall from \eqref{e:control.problems.endpoint.measures} that $\mu(q_*,b,\sigma,p,\zeta)$ denotes the law of $X(1)$ where $X$ is given by \eqref{e:formal.sde.main} with parameters $(q_*,b,\sigma,p,\zeta)$. (Note that in order to accommodate controls $(b,\sigma)$ where the budget constraint is only satisfied in the averaged sense of condition~\ref{it:4.1}, we allow $f$ to map into $[0,+\infty)$. On the other hand, all classes of measures that we will represent using Definition~\ref{d:X.adm} will have $q_* = q_0$ and $\zeta = \delta_0$, so we will set these parameters in the definition.)
\end{dfn}

\begin{rmk}
We compare the above Definition~\ref{d:X.adm} with our earlier Definition~\ref{d:SDE.conds}. In Definition~\ref{d:X.adm}, we have the following:
\begin{itemize}
\item The control processes $b,\sigma$ are progressively measurable, in keeping with condition~\ref{it:1.2}.
\item The function $p$ is specified, allowing us to handle any of the conditions~\ref{it:3.1}--\ref{it:3.7}.
\item The budget $f$ is is specified, allowing us to handle any of the conditions~\ref{it:4.1}--\ref{it:4.3}.
\item The initial measure is $\zeta=\delta_0$, in keeping with the first requirement of condition~\ref{it:6.2}.
\item The value of $q_*$ is set to $q_0$, in keeping with condition~\ref{it:8.2}.
\end{itemize}
The remaining conditions of Definition~\ref{d:SDE.conds} (\ref{it:2.1} or \ref{it:2.2}; \ref{it:5.1} or \ref{it:5.2}; the second part of \ref{it:6.2}; \ref{it:7.1} or \ref{it:7.2}) all involve the control $w_t$ and the process $Y(t)$. In Definition~\ref{d:X.adm}, all these remaining conditions are summarized by the budget function $f$. 
\end{rmk}

\begin{ppn}
\label{p:continuity.in.auxilliary.budget}
For any $q_0 \in [0,1)$ and $p \in \incr([q_0,1];[0,1])$, 
the mapping $f\mapsto \cM(q_0;p,f)$,
as a function from
 $L^1([q_0,1];[0,+\infty))$ to $P(\cP_2(\bbR))$, is uniformly continuous:
    \[
        d_\cH\Big(\cM(q_0;p,f),\cM(q_0;p,\tilde f)\Big)^2
        \le 2 \|f-\tilde f\|_1\,,
    \]
where we recall that $P(\cP_2(\bbR))$ is metrized by the Hausdorff distance \eqref{e:hausdorff}. As a consequence, we have the inclusion  $\sM^{\sv}(\alpha,q_0)\subseteq\sM^{\svi}(\alpha,q_0)$.

\begin{proof}
Let $f,\tilde f \in L^1([q_0,1];[0,+\infty))$, and suppose $(b,\sigma)\in\Adm(q_0;p,f)$. Let
	\[
	a(t)
	\equiv \min\bigg\{
	\bigg(\frac{\tilde{f}(t)}
	{f(t)}\bigg)^{1/2}
	,1\bigg\} \in[0,1]\,,
	\]
so in particular if $f(t)=0$ then we take $a(t)=1$. Then define $\tilde{b}_t\equiv a(t) b_t$, and $(\tilde{\sigma}_t-1)\equiv a(t)(\sigma_t-1)$, and note that $\sigma\ge0$ implies
$\tilde{\sigma}_t \ge 1-a(t)\ge0$.  By definition, $(\tilde{b},\tilde{\sigma})$ satisfies condition~\eqref{e:budget-f}, hence it belongs to $\Adm(q_0;p,f)$ as specified by Definition~\ref{d:X.adm}. Let $X,\tilde{X}$ be the 
processes defined by 
\eqref{e:formal.sde.main} with controls $(b,\sigma)$ and $(\tilde{b}, \tilde{\sigma})$ respectively, coupled with the same initialization
$X(q_0) =\tilde{X}(q_0) = 0$ and the same Brownian motion $B$. Then, applying the Cauchy--Schwarz inequality similarly as in the proof of Lemma~\ref{l:L2-bound-on-stochastic-control}, we have
    \begin{align*}
       &\bbE\Big[(X(1) - \tilde{X}(1))^2\Big]
        \le 2 \bbE\bigg[ \bigg(
            \int_{q_0}^1
            p'(t)^{1/2} |b_t - \tilde b_t|\,dt
        \bigg)^2 \bigg]
        + 2 \int_{q_0}^1 
        s(t)^2
        \bbE\Big[ (\sigma_t
        -\tilde{\sigma}_t
        )^2\Big]\,dt \\
    &\qquad\le 2 \int_{q_0}^1
    \E\bigg[ (b_t - \tilde b_t)^2
    + s(t)^2 
    \Big( (\sigma_t-1)
    -(\tilde{\sigma}_t-1)\Big)^2\bigg]\,dt \\
    &\qquad= 2 \int_{q_0}^1
    (1-a(t) )^2
\E\Big[(b_t)^2 + s(t)^2 
(\sigma_t-1)^2\Big]\,dt \\
&\qquad\le
2 \int_{q_0}^1 (1-a(t) )^2
	 f(t)\,dt
	 \le 2\int_{q_0}^1 \Big(f(t)^{1/2}-\tilde{f}(t)^{1/2}\Big)^2\,dt
\le 2\|f-\tilde{f}\|_1\,.
\end{align*}
This proves the first assertion of the proposition.
For the second assertion, we repeat the most relevant part of Definition~\ref{d:SDE.classification}:
	\begin{center}
	\begin{tabular}[h]{r|cccccccc}
	\MSRheader\\
	\hline
	\TMPTHRdefn\\
	\GREATdefn
	\end{tabular}
	\end{center}
The $\sv$ and $\svi$ classes  differ only in the budget constraint: the $\sv$ class satisfies condition~\ref{it:4.1}, while the $\svi$ class satisfies condition~\ref{it:4.3}. Therefore, let us suppose we have 
	\begin{align*}
	&(q_*=q_0,b,\sigma,w,p,\zeta=\delta_0,\zeta^\Ising)
	\in
	\Adm^{\sv}(\alpha,q_0;L,\epsilon)\,,\\
	&\mu = \mu(q_*=q_0,b,\sigma,p,\zeta=\zeta_0)
 \in \cM^{\sv}(\alpha,q_0;L,\epsilon)\,.\end{align*}
Define the budget functions
        \[
        f(t) = \max\bigg\{
            \bbE\bigg[
            (b_t)^2 + \frac{s(t)^2}{p(t)}(\sigma_t-1)^2\bigg],
            \frac{(\bbE w_t)^2}{\alpha}
        \bigg\}, \quad
        \tilde f(t) =\frac{(\bbE w_t)^2}{\alpha}\,.
    \]
Then, using the notation of
Definition~\ref{d:X.adm}, we have
from condition~\ref{it:4.1}
 that
 $(b,\sigma)$ belongs to $\Adm(q_0;p,f)$, and $\|f - \tilde f\|_1 \le \epsilon$. Therefore, by the first assertion of this proposition, there exists $\tilde \mu \in \cM(q_0;p,\tilde f)$ such that
$\bbW_2(\mu,\tilde \mu)^2 \le 2\epsilon$. Moreover, following the above construction, we can take
	\[ \tilde \mu
=\mu(q_*=q_0,\tilde{b},\tilde{\sigma},p,\zeta=\delta_0)
\in \cM^{\svi}(\alpha,q_0;L,\epsilon)\,.\]
It follows that, for some $\tilde\epsilon=\tilde\epsilon(\epsilon)$ tending to $0$ with $\epsilon$, we have
\[
    \cB_{\epsilon}\Big(\cM^{\sv}(\alpha,q_0;L,\epsilon)\Big)
    \subseteq \cB_{\tilde\epsilon}\Big(\cM^{\svi}(\alpha,q_0;L,\tilde\epsilon)\Big)\,.
\]
Recalling \eqref{e:measure.classes.simplification.without.L0}, taking a union over $L$ and an intersection over $\epsilon$ proves the second assertion of this proposition.
\end{proof}

We conclude this section by presenting the
\hyperlink{p:cor.feasible.distributions.continuity.main}{proof of Corollary~\ref{cor:feasible.distributions.continuity} for $\ocM^{\Ising}(\alpha)$}. The \hyperlink{p:cor.feasible.distributions.continuity.sym}{proof of Corollary~\ref{cor:feasible.distributions.continuity} for $\ocM^{\Ising,\sym}(\alpha)$} appears in \S\ref{ss:sdes.symmetric.endpoints}.

\begin{proof}[\hypertarget{p:cor.feasible.distributions.continuity.main}{Proof of Corollary~\ref{cor:feasible.distributions.continuity} for $\alpha \mapsto \ocM^{\Ising}(\alpha)$}]
Recall from Remark~\ref{r:measure.classes.simplification} that $\ocM^{\Ising}(\alpha)=\sM^{\sideal}(\alpha,0)$. We repeat the most relevant part of Definition~\ref{d:SDE.classification}:
	\begin{center}
	\begin{tabular}[h]{r|cccccccc}
	\MSRheader\\
	\hline
	\SUPERPERFECTdefn
	\end{tabular}
	\end{center}
With this in mind, let $\BUDGETS(\alpha)$ denote the set of all  $f \in L^\infty([0,1];[0,1/\alpha])$ of the form
\[
    f(t) = \budget(t;\alpha) 
    = \frac{(\E w_t)^2}{\alpha}\,,
\]
for $\zeta^\Ising \in \cP([-1,1])$ and $w \in \ProgMsrbl(\cF_Y)$ satisfying the following conditions:
\begin{itemize}
\item $\bbE[(w_t)^2] = 1$ for all $t\in [0,1]$;
\item for $Y$ given by  \eqref{e:formal.sde.ising} with 
$Y(q_*=q_0=0)=0$,  we have 
$|Y(1)| = 1$ almost surely.
\end{itemize}
Recall from Definition~\ref{d:incr.p} that
$\incr([0,1];[0,1]) = \mathscr{P}$. 
Then, using the notation of Definition~\ref{d:X.adm}, we have
\[ \ocM^{\Ising}(\alpha)
    = \overline{
        \bigcup_{f\in \BUDGETS(\alpha)}
        \bigcup_{p\in\mathscr{P}}
        \cM\Big (0;p,f\Big)
    }\,,
\]
where the overline indicates closure with respect to the $\bbW_2$ metric. Similarly, for any $\alpha'$, we have
\[
    \ocM^{\Ising}(\alpha')
    =
    \overline{
        {{\bigcup_{f\in \BUDGETS(\alpha)}}}
        \bigcup_{p\in\mathscr{P}}
        \cM\Big (0;p, \alpha f/\alpha' \Big)
    }\,.
\]
The conclusion follows by applying Proposition~\ref{p:continuity.in.auxilliary.budget}.
\end{proof}
\end{ppn}

\subsection{Continuity in the function $p$}
\label{ss:continuity.in.p}

In this subsection we prove the following inclusions:
\begin{itemize}
    \item Lemma~\ref{l:prf2.in.prf4.prf15.in.idealast} proves $\sM^{\sx}(\alpha,q_0) \subseteq \sM^{\sxii}(\alpha,q_0)$ and $\sM^{\six}(\alpha,q_0) \subseteq \sM^{\sidealast}(\alpha,q_0)$;
    \item Lemma~\ref{l:great.in.spr.prf4.in.prf2} proves $\sM^{\svi}(\alpha,q_0) \subseteq \sM^{\svii}(\alpha,q_0)$ and $\sM^{\sxii}(\alpha,q_0) \subseteq \sM^{\sx}(\alpha,q_0)$;
    \item Proposition~\ref{p:continuity.in.p.near.0} proves $\sM^{\sideal}(\alpha,0) \subseteq \sM^{\sxii}(\alpha,0)$.
\end{itemize}
All these statements involve variations among a subset of the conditions \ref{it:3.1}--\ref{it:3.7}. We also remark that in the \hyperlink{p:thm.control.problems.main}{proof of Theorem~\ref{thm:control.problems.main}} above, we in fact only require the equality $\sM^{\sx}(\alpha,q_0) = \sM^{\sxii}(\alpha,q_0)$ and the inclusion $\sM^{\six}(\alpha,q_0) \subseteq \sM^{\sidealast}(\alpha,q_0)$ in the case $q_0=0$. 

For these inclusions, we will show that the sets of endpoint measures are suitably continuous with respect to the function $p$. Because the cases $\star \in \{\svi,\svii,\sx,\sxii,\sideal\}$ all use  conditions \ref{it:2.2} and \ref{it:4.3}, for each $t\in [q_0,1]$ we simply have
\[\Cost(t) \le
    \frac{(\E w_t)^2}{ \alpha} 
    \le \frac{\E[(w_t)^2]}{\alpha} =\frac1\alpha\,.
\]
Therefore we can work with elements of $\cM(q_0;p,f)$ (see Definition~\ref{d:X.adm}) where the budget function $f$ lies in $L^\infty([q_0,1];[0,1/\alpha])$.

\begin{ppn}
    \label{p:continuity.in.p}
    For fixed $q_0 \in [0,1)$ and $f \in L^\infty([q_0,1];[0,1/\alpha])$, the mapping $p \mapsto \cM(q_0;p,f)$, viewed as a function from $\incr([q_0,1];[0,1])$ to  $P(\cP_2(\bbR))$, is uniformly continuous:
    \[
        d_\cH\Big(\cM(q_0;p,f),\cM(q_0;\tilde p,f)\Big)^2
        \le 6(1+\|f\|_\infty) \|p'-\tilde p'\|_1\,.
    \]
Here $\incr([q_0,1];[0,1])$ is metrized by  $\mathsf{dist}(p,\tilde p) = \|p'-\tilde p'\|_1$, and we recall that  $P(\cP_2(\bbR))$ is metrized by the Hausdorff distance \eqref{e:hausdorff}. 

\begin{proof}
Let $p,\tilde p\in \incr([q_0,1];[0,1])$ with $\|p'-\tilde p'\|_1 = \delta$.
Since $p(1) = \tilde p(1) = 1$, it follows by integrating that $|p(t) - \tilde p(t)| \le \delta$ for all $t\in [q_0,1]$, and so certainly we have $\|p-\tilde{p}\|_1\le\delta$. Recalling the notation of Definition~\ref{d:X.adm}, let $(b,\sigma)\in \Adm(q_0;p,f)$. We will construct modified controls $(\tilde b,\tilde \sigma) \in \Adm(q_0;\tilde p,f)$ such that the corresponding solutions to \eqref{e:formal.sde.main} approximate each other. Recall $s(t)\equiv (tp)'(t)^{1/2}$, and similarly denote
$\tilde{s}(t)\equiv (t\tilde{p})'(t)^{1/2}$.
We set $\tilde{b}=b$, and let $\hat \sigma$ be defined by the relation
	\[\hat{\sigma}_t-1
	= a(t) 
	(\sigma_t-1)\,,\quad
	a(t)
	\equiv 
	\frac{s(t)}{p(t)^{1/2}}
	\frac{\tilde{p}(t)^{1/2}}
		{\tilde{s}(t)}
	\]
for all $t\in[q_0,1]$. We note that 
	\[
	\frac{s(t)}{p(t)^{1/2}}
	= \frac{[p(t) + tp'(t)]^{1/2}}
		{p(t)^{1/2}} \ge 1\,,
	\]
so the ratio $a(t)$ can only be large if $p'(t)$ is very large, which occurs on a set of times of small measure. It is clear from the construction that $(\tilde b,\hat{\sigma})$ satisfies the budget constraint \eqref{e:budget-f} defining $\Adm(q_0;\tilde p,f)$, with the exception that $\hat{\sigma}$ is not necessarily non-negative. Let $X$ be given by \eqref{e:formal.sde.main} with controls $(b,\sigma)$, and $\hat X$ be given by \eqref{e:formal.sde.main} with controls $(\tilde b, \hat \sigma)$. We couple these processes with the same initialization $X(q_0) = \hat{X}(q_0) = 0$ and the same Brownian motion $B$. We will first show the endpoint laws of $X(1)$ and $\hat{X}(1)$ are close, then address the issue of $\hat{\sigma}$ being possibly negative. To this end, we recall from \eqref{e:formal.sde.main} that
  \[
        X(1) = \int_{q_0}^1 p'(t)^{1/2} b_t\,dt
        + \int_{q_0}^1 s(t)\,dB(t)
        + \int_{q_0}^1 s(t) (\sigma_t-1)\,dB(t)
        \,,
    \]
and similarly for $\hat{X}(1)$. This gives us a decomposition
$X(1)-\hat{X}(1) = \textup{(a)}+\textup{(b)}+\textup{(c)}$, and we proceed to bound each of the three terms:
    \begin{align*}\E[\textup{(a)}^2]
    &= \E \bigg[
    \bigg(\int_{q_0}^1 \Big(
        p'(t)^{1/2}
        - \tilde{p}'(t)^{1/2}\Big) b_t\,dt\bigg)^2\bigg]
    \le \bigg(
     \int_{q_0}^1\Big(
        p'(t)^{1/2}
        - \tilde{p}'(t)^{1/2}\Big)^2\,dt\bigg)
    \bigg( \int_{q_0}^1 \E[(b_t)^2]\,dt\bigg) \\
    &\le\bigg(\int_{q_0}^1 |p'(t)-\tilde{p}'(t)|\,dt\bigg)
    \bigg(\int_{q_0}^1 f(t)\,dt\bigg) 
    =  \delta \|f\|_1 \,,
    \end{align*}
having used the Cauchy--Schwarz inequality along with the budget constraint~\eqref{e:budget-f}. Next,
    \[\E[\textup{(b)}^2]
    = \int_{q_0}^1 \Big( s(t)-\tilde{s}(t) \Big)^2\,dt
    \le \int_{q_0}^1 \Big|s(t)^2 - \tilde{s}(t)^2\Big|\,dt
    \le \|p-p\|_1+\|p'-p'\|_1 \le 2\delta\,.\]
Lastly, using the above definition of $\hat\sigma$, we have
	\begin{align*}
	\E[\textup{(c)}^2]
    &= \int_{q_0}^1
    \Big( s(t) (\sigma_t-1)-\tilde{s}(t)(\hat{\sigma}_t
    -1)\Big)^2\,dt\\
    &= \int_{q_0}^1
    \frac{s(t)^2}{p(t)}
    \E[(\sigma_t-1)^2]
    \cdot p(t)
    \bigg( 1- 
    	\frac{\tilde{s}(t)}{s(t)} a(t)
    	 \bigg)^2
	\,dt\\
	&\le
	\|f\|_\infty
	\int_{q_0}^1
	p(t)
    \bigg( 1- \frac{\tilde{p}(t)^{1/2}}
		{p(t)^{1/2}}
    	 \bigg)^2
	\,dt 
	\le \|f\|_\infty \delta\,,
	\end{align*}
having used the budget constraint~\eqref{e:budget-f} along with the calculation from the bound on (a) above. Collecting the above bounds gives
    \[
    \E\Big[( X(1)-\hat{X}(1))^2\Big]
    \le 
    3\bigg\{
    \E[\textup{({{a}})}^2]+\E[\textup{(b)}^2]+\E[\textup{(c)}^2]
    \bigg\}
    \le 6
    \delta
    \Big(1+\|f\|_\infty\Big)\,.
    \]
To conclude, we address the issue of $\hat{\sigma}$ being possibly negative: set $\tilde{\sigma}=|\hat{\sigma}|\ge0$. Then $(\tilde b,\tilde \sigma)$ continues to satisfy \eqref{e:budget-f}, so we have $(\tilde b,\tilde \sigma) \in \Adm(q_0;\tilde p,f)$. Define 
	\[
	\tilde B(t) = \int_{q_0}^t (\ind\{\hat \sigma \ge 0\} - \ind\{\hat \sigma < 0\}) \,\de B(t)\,,\]
so $\tilde{B}$ is marginally also a standard Brownian motion. Let $\tilde X$ be given by \eqref{e:formal.sde.main} with controls $(\tilde b,\tilde m)$, initialization $\tilde X(q_0) = 0$, and Brownian motion $\tilde B$ in place of $B$. Then $\tilde X = \hat X$ as processes, so it follows from the above analysis that $\E [(X(1) - \tilde X(1))^2] \le 6\delta (1+\|f\|_\infty)$. This proves the claim.
\end{proof}
\end{ppn}

\begin{lem}\label{l:prf2.in.prf4.prf15.in.idealast}
For all $q_0 \in [0,1)$, we have 
\begin{align*}
\sM^{\sx}(\alpha,q_0) 
&\subseteq \sM^{\sxii}(\alpha,q_0)\,,\\
\sM^{\six}(\alpha,q_0) 
&\subseteq \sM^{\sidealast}(\alpha,q_0)\,.
\end{align*}

\begin{proof}
We repeat the most relevant part of Definition~\ref{d:SDE.classification} for this proof:
	\begin{center}
	\begin{tabular}[h]{r|cccccccc}
	\MSRheader\\
	\hline
	\PRFTWOdefn\\
	\PRFFOURdefn\\
	\PRFONEHALFdefn\\
	\SUPERPERFECTASTdefn
	\end{tabular}
	\end{center}
The classes only differ in the conditions on the function $p$.

    Consider any $\mu \in \cM^{\sx}(\alpha,q_0;L,\epsilon)$ where $\epsilon < 1/2$. 
    Then $\mu$ is of the form $\mu = \mu(q_*=q_0,b,\sigma,p,\zeta = \delta_0)$ for some $(q_*=q_0,b,\sigma,w,p,\zeta = \delta_0,\zeta^\Ising) \in \Adm^{\sx}(\alpha,q_0;L,\epsilon)$.
    Since $p$ satisfies condition~\ref{it:3.1}, $p$ is twice differentiable and therefore continuous, and further satisfies $p(q_0) \le \epsilon$.
    So we can find $q' \in (q_0,1)$ such that $p(q') = 2\epsilon$, and define $\tilde p \in \incr([q_0,1];[0,1])$ by
    \[
        \tilde p(t) = \begin{cases}
            2\epsilon (t - q_0) / (q' - q_0) & t \in [q_0,q')\,, \\
            p(t) & t \in [q',1]\,.
        \end{cases}
    \]
    In particular $\tilde p(q_0) = 0$, so $\tilde p$ satisfies condition~\ref{it:3.5withzero}, as required for the $\sxii$ class.
    Since $p$ and $\tilde p$ agree on $[q',1]$, we have
    \[
        \|p'-\tilde p'\|_1
        = \int_{q_0}^{q'} |p'(t) - \tilde p'(t)|\,dt
        \le \int_{q_0}^{q'} (p'(t) + \tilde p'(t))\,dt
        \le p(q') + \tilde p(q')
        = 4\epsilon\,.
    \]
Recalling Definition~\ref{d:X.adm}, we have $\mu \in \cM(q_0;p,f)$ for budget function
$f(t) = (\E w_t)^2 / \alpha$. Applying
    Proposition~\ref{p:continuity.in.p} gives the existence of a measure
    \[
        \tilde\mu \in \cM(q_0;\tilde p,f)
        \subseteq \cM^{\sxii}(\alpha,q_0)
    \]
    with $\bbW_2(\mu,\tilde \mu)^2 \le 6(1+1/\alpha) \cdot 4\epsilon$.
    Since this holds for all $\mu \in \cM^{\sx}(\alpha,q_0;L,\epsilon)$ with $\epsilon < 1/2$, we conclude
    \[
        \cM^{\sx}(\alpha,q_0;L,\epsilon) \subseteq \cB_{\tilde\epsilon}\Big(\cM^{\sxii}(\alpha,q_0)\Big)
    \]
    for some $\tilde\epsilon$ tending to $0$ with $\epsilon$.
    By adjusting $\tilde\epsilon$ we can also ensure
    \[
        \cB_{\epsilon}\Big(\cM^{\sx}(\alpha,q_0;L,\epsilon)\Big) \subseteq \cB_{\tilde\epsilon}\Big(\cM^{\sxii}(\alpha,q_0)\Big)\,.
    \]
 Recall \eqref{e:measure.classes.simplification.without.L0}: taking a union over $L$ and an intersection over $\epsilon$ proves the first assertion, $\sM^{\sx}(\alpha,q_0) \subseteq \sM^{\sxii}(\alpha,q_0)$.

    The proof of the second assertion $\sM^{\six}(\alpha,q_0) \subseteq \sM^{\sidealast}(\alpha,q_0)$ is very similar, so we only describe the differences.
    Consider any $\mu \in \cM^{\six}(\alpha,q_0;L,\epsilon)$ where $\epsilon < 1/2$, realized by controls $(q_*=q_0,b,\sigma,w,p,\zeta = \delta_0,\zeta^\Ising) \in \Adm^{\six}(\alpha,q_0;L,\epsilon)$.
    In this case $p$ satisfies condition~\ref{it:3.1concave}, so $p$ is concave in addition to the properties of condition~\ref{it:3.1} used above.
    Let $q'$ and $\tilde p$ be defined as above.
    Since $p(q_0) \ge 0$ and $p(q') = 2\epsilon$, the concavity of $p$ implies $p(t) \ge \tilde p(t)$ for all $t\in [0,q']$.
    Thus $\tilde p$ is also concave, and in particular satisfies condition~\ref{it:3.3}, as required for the $\sidealast$ class.
Then, arguing identically as above completes the proof.
\end{proof}
\end{lem}

\begin{lem} \label{l:approximate-p-by-C2}
For any $\epsilon \in (0,1)$ and any function $p \in \incr([q_0,1];[0,1])$ with $p(q_0)=0$ (i.e., satisfying condition~\ref{it:3.5withzero}), there exists sufficiently large $L$
 (depending on $p,\epsilon$) and $\tilde p \in \incr([q_0,1];[0,1])$ satisfying $1/L \le \tilde{p}(q_0) \le \epsilon$ and $\|\tilde{p}\|_{C^2([q_0,1])} \le L$ (i.e., satisfying condition~\ref{it:3.1}) 
 such that $\|p' - \tilde p'\|_1 \le \epsilon$. 
If $p$ is furthermore concave (i.e., satisfying condition~\ref{it:3.3}), then the same conclusion holds with $\tilde p$ also concave (i.e., satisfying condition~\ref{it:3.1concave}). 
\begin{proof}
Since $p\in \incr([q_0,1];[0,1])$ is absolutely continuous, we have $p' \in L^1([q_0,1];[0,+\infty))$. Since $p(q_0) = 0$ and $p(1) = 1$, we have $\|p'\|_1 = 1$. 

Throughout this proof, we use the term ``step function'' to refer to a function on $[q_0,1]$ that is piecewise constant with finitely many discontinuities. Since step functions are dense in $L^1([q_0,1];[0,+\infty))$, there exists a step function $g_1 : [q_0,1] \to [0,+\infty)$ such that
	\[\|p' - g_1\|_1 \le \epsilon / 6\,.\]
For sufficiently large $L$, there exists a smooth function $g_2 : [q_0,1] \to [0,L]$ with $\|g_2\|_{C^2([q_0,1])} \le L/10$, such that $\|g_1 - g_2\|_1 \le \epsilon / 6$, which implies
	\[\|p' - g_2\|_1 \le \epsilon / 3\,.\]
Since $\|p'\|_1 = 1$, we have $\|g_2\|_1 \le 1 + \epsilon/3$. 
We then define
	\[\tilde p(t) = 1 
		- \bigg(1-\frac{\epsilon}{2}\bigg)
		\int_t^1 g_2(s)\,ds\,,\]
and we note that $\tilde{p}$ is smooth and nondecreasing with
$\|\tilde{p}\|_{C^2([q_0,1])} \le L$. For $\epsilon < 1$, we can bound
	\[
    \|p' - \tilde p'\|_1
    = \bigg\|p' - \bigg(1-\frac{\epsilon}{2}\bigg) g_2\bigg\|_1
    \le \|p' - g_2\|_1 + \frac{\epsilon}{2} \|g_2\|_1
    \le \frac{\epsilon}{3} + \bigg(1+\frac{\epsilon}{3}\bigg) \frac{\epsilon}{2}
    \le \epsilon\,.
\]
Since $p(1) = \tilde p(1) = 1$, it follows that
\[
    \tilde p(q_0) \le p(q_0) + \|p' - \tilde p'\|_1
    =\|p' - \tilde p'\|_1
     \le \epsilon\,.
\]
For sufficiently large $L$, we also have
\[
    \tilde p(q_0) = 1 - \bigg(1-\frac{\epsilon}{2}\bigg) \|g_2\|_1 
    \ge 1 - \bigg(1-\frac{\epsilon}{2}\bigg)
    	\bigg(1+\frac{\epsilon}{3}\bigg)
	\ge \frac{\epsilon}{6} 
    \ge \frac{1}{L}\,.
\]
Thus $\tilde p$ satisfies condition~\ref{it:3.1}, as desired.
If $p$ is concave, then $p'$ is nonincreasing, and we can further arrange for the approximations $g_1$ and $g_2$ to be nonincreasing: in this case the resulting $\tilde{p}$ is concave, and so it satisfies condition~\ref{it:3.1concave}. This concludes the proof.
\end{proof}
\end{lem}

\begin{lem}\label{l:great.in.spr.prf4.in.prf2}
For all $q_0 \in [0,1)$, we have 
\begin{align*}\sM^{\svi}(\alpha,q_0) 
&\subseteq \sM^{\svii}(\alpha,q_0)\,,\\
\sM^{\sxii}(\alpha,q_0) 
&\subseteq \sM^{\sx}(\alpha,q_0)\,.
\end{align*}
    
\begin{proof} 
We repeat the most relevant part of Definition~\ref{d:SDE.classification} for this proof:
	\begin{center}
	\begin{tabular}[h]{r|cccccccc}
	\MSRheader\\
	\hline
	\GREATdefn\\
	\SPRdefn\\
	\PRFFOURdefn\\
	\PRFTWOdefn
	\end{tabular}
	\end{center}
The $\svi$ and $\svii$ classes differ only in the condition on the function $p$: the $\svi$ class satisfies condition~\ref{it:3.2}, while the $\svii$ class satisfies condition~\ref{it:3.4}. The $\sxii$ 
and $\sx$ classes also differ only in the condition on the function $p$: the $\sxii$ 
class satisfies condition~\ref{it:3.5withzero}, while the $\sx$ class satisfies condition~\ref{it:3.1}.

    Consider any $\mu \in \cM^{\svi}(\alpha,q_0;L,\epsilon)$.
    Then $\mu$ is of the form $\mu = \mu(q_*=q_0,b,\sigma,p,\zeta = \delta_0)$ for some $(q_*=q_0,b,\sigma,w,p,\zeta = \delta_0,\zeta^\Ising) \in \Adm^{\svi}(\alpha,q_0;L,\epsilon)$.
    Since $p$ satisfies condition~\ref{it:3.2}, it is concave.

    Let $\epsilon_0 > 0$ be arbitrary. 
    By Lemma~\ref{l:approximate-p-by-C2} 
     (applied with $\epsilon_0,L_0$ in place of $\epsilon,L$), there exists sufficiently large $L_0 = L_0(\epsilon_0)$ and concave $\tilde p \in \incr([q_0,1];[0,1])$ such that $1/L_0 \le \tilde p(q_0) \le \epsilon_0$ and $\|\tilde p\|_{C^2([q_0,1])} \le L_0$ (i.e. satisfying condition~\ref{it:3.4}, as required of the $\svii$ class) and $\|p'-\tilde p'\|_1 \le \epsilon_0$.

By condition~\ref{it:4.3}, we have
$\mu \in \cM(q_0;p,f)$ for budget function
 $f(t) = (\E w_t)^2 / \alpha$.
    Proposition~\ref{p:continuity.in.p} then implies the existence of 
    \[
        \tilde\mu \in \cM(q_0;\tilde p,f)
        \subseteq \cM^{\svii}(\alpha,q_0;\epsilon,L_0,\epsilon_0)
    \]
    with $\bbW_2(\mu,\tilde \mu)^2 \le 6(1+1/\alpha) \cdot \epsilon_0 \equiv (\tilde\epsilon_0)^2$.
    We have thus shown that for any $\epsilon_0 > 0$ there exists $L_0 = L_0(\epsilon_0)$ such that
    \[
        \cM^{\svi}(\alpha,q_0;L,\epsilon) \subseteq \cB_{\tilde\epsilon_0}\Big(\cM^{\svii}(\alpha,q_0;\epsilon,L_0,\tilde \epsilon_0)\Big)\,.
    \]
Consequently, for $\hat\epsilon_0 = \epsilon_0 + \tilde\epsilon_0$, we have
    \[
        \cB_{\epsilon_0}\Big(\cM^{\svi}(\alpha,q_0;L,\epsilon)\Big) 
        \subseteq \cB_{\hat\epsilon_0}\Big(
        \cM^{\svii}(\alpha,q_0;\epsilon,L_0,\hat \epsilon_0)
        \Big)\,.
    \]
    Taking unions and intersections of the above, in the order specified by Definition~\ref{d:measure.classes}, gives the first assertion $\sM^{\svi}(\alpha,q_0) \subseteq \sM^{\svii}(\alpha,q_0)$.

    The proof of the remaining inclusion $\sM^{\sxii}(\alpha,q_0) \subseteq \sM^{\sx}(\alpha,q_0)$ is very similar. We consider any $\mu \in \cM^{\sxii}(\alpha,q_0)$, realized by controls $(q_*=q_0,b,\sigma,w,p,\zeta = \delta_0,\zeta^\Ising) \in \Adm^{\sxii}(\alpha,q_0)$, so that $p$ satisfies condition~\ref{it:3.5withzero}. 
    For arbitrary $\epsilon > 0$,  Lemma~\ref{l:approximate-p-by-C2} 
    implies the existence of sufficiently large $L = L(\epsilon)$ and $\tilde p \in \incr([q_0,1];[0,1])$ such that $1/L \le \tilde p(q_0) \le \epsilon$ and $\|\tilde p\|_{C^2([q_0,1])} \le L$ (i.e., satisfying condition~\ref{it:3.1}, as required of the $\sx$ class) and $\|p'-\tilde p'\|_1 \le \epsilon$.
    Arguing as above shows that for $\tilde\epsilon^2 \equiv 6(1+1/\alpha) \cdot \epsilon$,
    \[
        \cM^{\sxii}(\alpha,q_0) \subseteq \cB_{\tilde\epsilon}(\cM^{\sx}(\alpha,q_0;L,\tilde\epsilon))\,,
    \]
    which implies the result similarly to above.
\end{proof}
\end{lem}

For the next result, we will argue that if $q_0=0$, the sets of endpoint measures $\cM(0;\tilde p,f)$ are also continuous in $p$ with respect to arbitrary perturbations in an interval $[0,\epsilon)$. Note that $\|p' - \tilde p'\|_1$ is not necessarily small under such perturbations, so this notion of continuity is not captured by Proposition~\ref{p:continuity.in.p}. 

\begin{ppn}
    \label{p:continuity.in.p.near.0}
    Suppose $f \in L^\infty([0,1];[0,1/\alpha])$ and $p,\tilde p\in\sP=\incr([0,1];[0,1])$, and the restrictions of $p,\tilde p$ to $[\epsilon,1]$ are equal.
    Then
    \[
    	d_\cH\Big(\cM(0;p,f),\cM(0;\tilde p,f)\Big)^2
    	\le 9(1+\|f\|_\infty) \epsilon\,.
    \]
    As a consequence, we have the inclusion $\sM^{\sideal}(\alpha,0) \subseteq \sM^{\sxii}(\alpha,0)$.
\end{ppn}
\begin{proof}
	Let $(b,\sigma)\in \Adm(0;p,f)$. We construct the modified controls 
	\[
		(\tilde b_t,\tilde \sigma_t) = \begin{cases}
			(0,1) & t < [0,\epsilon)\,, \\
			(b_t,\sigma_t) & t \in [\epsilon,1]\,,
		\end{cases}
	\]
	so that clearly $(\tilde b,\tilde \sigma) \in \Adm(q_0;\tilde p,f)$.
	We will argue that the corresponding solutions to \eqref{e:formal.sde.main} approximate each other. 
	Let $X$ be given by \eqref{e:formal.sde.main} with controls $(b,\sigma)$, and $\tilde X$ be given by \eqref{e:formal.sde.main} with controls $(\tilde b, \tilde\sigma)$. We couple these processes with the same initialization $X(0) = \tilde X(0) = 0$ and the same Brownian motion $B$.
	Recall $s(t)\equiv (tp)'(t)^{1/2}$, and similarly denote $\tilde{s}(t)\equiv (t\tilde{p})'(t)^{1/2}$.
	Since $(p(t),s(t),b_t,\sigma_t)$ and $(\tilde p(t),\tilde s(t), \tilde b_t, \tilde\sigma_t)$ differ only on $t\in [0,\epsilon)$, we have
  	\begin{align*}
        X(1) - \tilde X(1) 
        &= \int_0^\epsilon (p'(t)^{1/2} b_t - \tilde p'(t)^{1/2} \tilde b_t)\,dt
        + \int_0^\epsilon (s(t) \sigma_t - \tilde s(t) \tilde\sigma_t)\,dB(t) \\
        &= \int_0^\epsilon p'(t)^{1/2} b_t \,dt
        + \int_0^\epsilon s(t) \sigma_t\,dB(t)
        + \int_0^\epsilon (-\tilde s(t))\,dB(t)
        \equiv \textup{(a)}+\textup{(b)}+\textup{(c)}\,.
    \end{align*}
    We proceed to bound each of the three terms.
    Note that since $s(t)^2 / p(t) \ge 1$, the budget constraint \eqref{e:budget-f} implies $\E[(\sigma_t-1)^2] \le f(t)$.
    Then, by the Cauchy--Schwarz inequality,
    \[
        \E[\textup{(a)}^2]
        = \E \bigg[\bigg(
            \int_0^\epsilon p'(t)^{1/2} b_t \,dt
        \bigg)^2\bigg]
        \le \bigg(\int_0^\epsilon p'(t)\,dt\bigg) \bigg(\int_0^\epsilon \E[(b_t)^2]\,dt\bigg)
        \le \epsilon \|f\|_\infty\,.
    \]
    Also,
    \begin{align*}
	\E[\textup{(b)}^2]
	&= \int_0^\epsilon s(t)^2 \E[(\sigma_t)^2]\,dt
	\le \int_0^\epsilon s(t)^2 \cdot 2\Big(\E[(\sigma_t-1)^2+1]\Big)\,dt \\
	&\le 2(\|f\|_\infty + 1) \int_0^\epsilon s(t)^2\,dt 
	= 2(\|f\|_\infty + 1) \epsilon p(\epsilon)
	\le 2(\|f\|_\infty + 1) \epsilon\,, \\
	\E[\textup{(c)}^2]
	&= \int_0^\epsilon \tilde s(t)^2\,dt
	= \epsilon \tilde p(\epsilon)
	\le \epsilon\,.
    \end{align*}
    Altogether we conclude
    \[
    	\E\Big[(X(1) - \tilde X(1))^2\Big]
    	\le 3 \Big(
    		\E[\textup{(a)}^2]
    		+ \E[\textup{(b)}^2]
    		+ \E[\textup{(c)}^2]
    	\Big)
    	\le 9(\|f\|_\infty + 1) \epsilon\,.
    \]
This proves the first assertion of the proposition.
For the second assertion, we first repeat the most relevant part of Definition~\ref{d:SDE.classification}:
	\begin{center}
	\begin{tabular}[h]{r|cccccccc}
	\MSRheader\\
	\hline
	\PRFFOURdefn\\
	\SUPERPERFECTdefn
	\end{tabular}
	\end{center}
The $\sxii$ and $\sideal$ classes differ only in the constraint on the function $p\in\sP$: the $\sxii$ class satisfies condition~\ref{it:3.5withzero} (i.e., $p(q_0 = 0) = 0$) while the $\sideal$ class satisfies condition~\ref{it:3.5} (i.e., no additional condition).
	Let us suppose we have controls $(q_*=0,b,\sigma,w,p,\zeta=\delta_0,\zeta^\Ising=\delta_0)$ satisfying the $(\alpha,0)$-$\sideal$ conditions, so that $p$ satisfies condition~\ref{it:3.5}. Let
	\[
		\mu = \mu(q_*=0,b,\sigma,p,\zeta=\zeta_0)
		\in \cM^{\sideal}(\alpha,0)
	\]
	be the associated endpoint measure.
	Define the budget function
	\[
		f(t) =\frac{(\bbE w_t)^2}{\alpha}\,.
	\]
By condition~\ref{it:4.3},  we have $(b,\sigma) \in \Adm(0;p,f)$.
	For an arbitrary $\epsilon > 0$, consider the function
	\[
		\tilde p(t) = \begin{cases}
			t p(\epsilon)/\epsilon & t\in [0,\epsilon)\,, \\\
			p(t) & t \in [\epsilon,1]\,,
		\end{cases}
	\]
for which $\tilde p( 0) = 0$ (meaning $\tilde{p}$ satisfies condition~\ref{it:3.5withzero}). 
	By the first assertion of this proposition, there exists $\tilde\mu \in \cM(0;\tilde p,f)$ such that $\bbW_2(\mu,\tilde \mu)^2 \le 9(1 + 1/\alpha) \epsilon$.
	Indeed, following the above construction we can take
	\[
		\tilde\mu = \mu(q_*=0,\tilde b, \tilde\sigma,\tilde p, \zeta=\delta_0)
		\in \cM^{\sxii}(\alpha,0)\,.
	\]
	Since $\epsilon$ was arbitrary, it follows that
	\[
		\cM^{\sideal}(\alpha,0)
		\subseteq \overline{\cM^{\sxii}(\alpha,0)}\,,
	\]
	where the overline denotes closure with respect to the $\bbW_2$ metric.
	However, recall from Remark~\ref{r:measure.classes.simplification} that $\sM^{\sxii}(\alpha,0)$ and $\sM^{\sideal}(\alpha,0)$ are 
	 themselves closed 
	with respect to the $\bbW_2$ metric.
	The second assertion of this proposition follows.
\end{proof}

\subsection{Ising endpoint condition}
\label{ss:exact.endpoint.ising}

The main result of this subsection is Lemma~\ref{l:spr.in.prf1}, which proves the inclusion $\sM^{\svii}(\alpha,q_0)\subseteq \sM^{\sviii}(\alpha,q_0)$. As a byproduct of our arguments, we will also obtain Lemma~\ref{l:ideal.contains.std.gaussian}, which shows that $\sM^{\sideal}(\alpha,0)$ contains the standard gaussian distribution.

Recalling Definitions~\ref{d:SDE.conds}--\ref{d:measure.classes}, proving $\sM^{\svii}(\alpha,q_0)\subseteq \sM^{\sviii}(\alpha,q_0)$ amounts to showing that the Ising exact endpoint condition $|Y(1)| = 1$ can be approximated by the relaxed condition $\bbE[(|Y(1)|-1)^2] \le \epsilon$.
We emphasize that, being purely a matter of stochastic control, this is different from the also important issue of whether \emph{approximate solutions to the random perceptron constraint satisfaction problem} (i.e., configurations $x\in\{-1,+1\}^N$ violating $o(1)$ fraction of the perceptron constraints) can be turned into exact solutions (satisfying all of the constraints). This latter task is essentially possible, and is a matter of suitably modifying the output of our main IAMP algorithm (from Section~\ref{sec:IAMP}): this will be handled in the companion work \cite{bogpinprogress}.

In preparation for proving Lemma~\ref{l:spr.in.prf1}, take control parameters
	\[(q_* = q_0,b,\sigma,w,p,\zeta=\delta_0,\zeta^\Ising)
	\in\Adm^{\svii}(\alpha,q_0;\epsilon,L_0,\epsilon_0)\,.\]
Note that this implies $\zeta^\Ising \in \cP([-1,1])$ and $\E[Y(q_0)^2] = q_0$ for $Y(q_0)\sim\zeta^\Ising$. Let $Y$ be given by \eqref{e:formal.sde.ising} with these parameters. Recalling Remark~\ref{r:true.mg}, we have that $Y$ is an $L^2$-bounded true martingale,
satisfying
    \[\E[Y(t)^2]=\bbE[ Y(q_0)^2 + \langle Y\rangle_t]
    = q_0+ \E\int_{q_0}^t (w_u)^2\,d\langle W\rangle_u
    = t\]
for all $t\in [q_0,1]$, and with $\bbE[(|Y(1)|-1)^2] \le \epsilon$. Define the stopping time
    \beq\label{e:tau}
    \tau = \inf\Big\{ t \ge q_0 : |Y_t| \ge 1 \Big\}\,,
    \eeq
and denote the stopped martingale $\hat{Y}(t)\equiv Y(t\wedge\tau)$.
For $t\in[q_0,1]$, we define
    \[\gamma(t)\equiv \E[\hat{Y}(t)^2]
    = \E[\hat{Y}(q_0)^2+ \langle \hat{Y}\rangle_t]
    = q_0 + \E\int_{q_0}^t (w_u)^2
        \ind\{u\le\tau\} \,du
    \,.
    \]
Then $\gamma(q_0)=q_0$, and for any $t\in [q_0,1]$ we have
    \beq
    \label{e:gamma.prime}
    \gamma'(t)
    =\E\Big[
    (w_t)^2 \ind\{t\le\tau\}
    \Big]
    \le \E[(w_t)^2]=1\,.
    \eeq
The following shows that stopping the martingale does not significantly affect the second moment:

\begin{ppn}
    \label{p:gamma-1}
    We have $\gamma(1) \ge 1 - 4\epsilon^{1/2}$.
\end{ppn}

We defer the
\hyperlink{proof:p.gamma-1}{proof of Proposition~\ref{p:gamma-1}}
to the end of the subsection, and proceed with the rest of the argument. We first provide an elementary lemma:

\begin{lem}\label{l:generalized.inverse.fn}
Let $\gamma:[q_0,1]\to[q_0,1]$ be a monotone and absolutely continuous function, and define its generalized inverse
	\[
	\lambda:[q_0,\gamma(1)]\to[q_0,1]\,,\quad
	\lambda(s)
	\equiv\inf\Big\{
	t : \gamma(t)>s
	\Big\}\,.
	\]
Then the set $\{s : \gamma'(\lambda(s))=0 \}$ has zero Lebesgue measure.

\begin{proof}
It is straightforward to verify that $\lambda$ is a right-continuous function with $\gamma(\lambda(s))=s$ for all $s$. For any measurable set $A\subseteq\R$, consider the preimage $\lambda^{-1}(A)$: we have $s\in \lambda^{-1}(A)$ if and only if $\lambda(s)\in A$, which implies $s=\gamma(\lambda(s)) \in \gamma(A)$. This proves that 
	\[\lambda^{-1}(A)
	\subseteq \gamma(A)\]
(the converse need not hold). Let $A=\{t\in[q_0,1] : \gamma'(t)=0\}$: by absolute continuity of $\gamma$, we have
	\[|\gamma(A)|
	= \int_{q_0}^1 \gamma'(t)
	\ind\{t\in A\}
	\,dt = 0\,.\]
Since $\lambda^{-1}(A)\subseteq \gamma(A)$, it follows that
	\[
	|\lambda^{-1}(A)|
	=\bigg| \Big\{s : \gamma'(\lambda(s))=0
	\Big\}\bigg|=0\,,
	\]
as claimed.
\end{proof}
\end{lem}

\begin{lem}\label{l:time.change}
For $\hat{Y}$ the stopped martingale defined above, we can construct a deterministic time change $\tilde{Y}$ that can be expressed as the stochastic integral
    \[
    \tilde{Y}(t)
    =\tilde{Y}(q_0) + \int_0^t \tilde{w}_t\,d\tilde{W}(t)
    \]
for $t\in[q_0,\gamma(1)]$,
where
$\tilde{Y}(q_0)=\hat{Y}(q_0)=Y(q_0)$,
$\tilde{W}$ is a standard Brownian motion, and $\tilde{w}$ satisfies
the constraint $\E[(\tilde{w}_t)^2]=1$ for all $t\in[q_0,\gamma(1)]$.

\begin{proof} We will construct a  martingale $\tilde{Y}$ to essentially be the deterministic time-change of $\hat{Y}$ such that $\E[\tilde{Y}(t)^2]=t$ for all $t\in[q_0,\gamma(1)]$. Then, for $t\in[\gamma(1),1]$, we will take $\tilde{Y}$ to be any martingale ending at $\pm1$ almost surely, which exists by Lemma~\ref{lem:Yt-process-always-is-completable} below. To formalize the time reparametrization, since $\gamma:[q_0,1]\to[q_0,\gamma(1)]$ is nondecreasing, we can define its right-continuous generalized inverse 
\[\lambda:[q_0,\gamma(1)]\to[q_0,1]\,.\]
As a consequence of \eqref{e:gamma.prime} we must have $\gamma(t)-\gamma(s)\le t-s$ and $\lambda(t)-\lambda(s)\ge t-s$ for all $q_0 \le s\le t\le1$. It  follows by combining these inequalities with Proposition~\ref{p:gamma-1} that
    \beq\label{e:gamma-estimate}
    \begin{aligned}
    s &\ge \gamma(s) \ge \gamma(1)+s-1 \ge s-4\epsilon^{1/2}\,,\\
    t &\le \lambda(t) \le \lambda(\gamma(1)) + t-\gamma(1)
    =1+t-\gamma(1) \le t+4\epsilon^{1/2}
    \end{aligned}
    \eeq
for all $s\in [q_0,1]$, and $t\in [q_0,\gamma(1)]$. It is clear from this that we must define $\tilde{Y}(t) = \hat{Y}(\lambda(t))$, so that
we will have $\E[\tilde{Y}(t)^2] = \E[\hat{Y}(\lambda(t))^2] = \gamma(\lambda(t))=t$. Since
    \[
    \hat{Y}(\lambda(t)) = Y(q_0) + \int_{q_0}^{\lambda(t)}
        w_u \ind\{u\le\tau\} \,dW(u)\,,
    \]
by considering the change of variables $u=\lambda(s)$ it is natural to define
    \beq\label{e:time-changed-W}
    \tilde{W}(t) = \int_{q_0}^{\lambda(t)} \gamma'(u)^{1/2} dW(u)\,.
    \eeq
Then $\tilde{W}$ is a local martingale with quadratic variation
    \[
    \langle \tilde{W}\rangle_t 
    = \int_{q_0}^{\lambda(t)} \gamma'(u)\,du
    = \gamma(\lambda(t))-\gamma(q_0) = t-q_0
    \]
for all $t\in[q_0,1]$, so by the L\'evy characterization it is a standard Brownian motion on $t\in[q_0,1]$, and is progressively measurable with respect to the time-changed filtration
$\tilde{\cF}_Y(t)\equiv\cF_Y(\lambda(t))$. We then define the modified control
    \[
    \tilde{w}_s
    \equiv
    \frac{w_{\lambda(s)}}{\gamma'(\lambda(s))^{1/2} }
      \ind\{\lambda(s)\le \tau\}\,,\]
where we define $\tilde{w}_s\equiv1$ in the case $\gamma'(\lambda(s))=0$ --- by Lemma~\ref{l:generalized.inverse.fn}, the set of such times $s$ is of measure zero. If $\gamma'(\lambda(s))=0$ then clearly we have $\E[(\tilde{w}_s)^2]=1$, and otherwise 
\eqref{e:gamma.prime} gives
	\[\E[(\tilde{w}_s)^2]
    =\frac{\E[(w_{\lambda(s)})^2 \ind\{
        \lambda(s)\le\tau\}]}
        {\gamma'(\lambda(s))}
    = 1\,.
    \]
We can then express the desired process as
    \[
    \tilde{Y}(t)
    = \hat{Y}(\lambda(t))
    =\tilde{Y}(q_0) + \int_{q_0}^t \tilde{w}_s\,d\tilde{W}(s)\,,
    \]
where $\tilde{Y}(q_0)=\hat{Y}(q_0)= Y(q_0)$. This concludes the proof.
\end{proof}
\end{lem}

The following lemma gives the extension to the full time interval $[q_0,1]$:

\begin{lem}\label{lem:Yt-process-always-is-completable}
On the time interval $[q_0,\gamma(1)]$ let $\tilde{Y}$, $\tilde{W}$, $\tilde{w}$, $\tilde{\cF}_Y$ be as given by Lemma~\ref{l:time.change}.
It is possible to extend these to the full time interval $[q_0,1]$ such that
$|\tilde Y(1)| = 1$ almost surely, and the same conditions as before are satisfied: $\tilde{W}$ is a standard Brownian motion,
$\tilde{w}$ is progressively measurable with respect to
    \[
    \tilde \cF_Y(t) = \sigma\Big(
        (\tilde W_s)_{s\in [q_0,t]}, U', \tilde Y(q_0)
    \Big)
\]
with $\bbE[(\tilde w_t)^2] = 1$ for all $t$,
and
 \[
        \tilde Y(t) = \tilde Y(\gamma(1))
            + \int_{\gamma(1)}^t \tilde w_s\,d\tilde{W}(s)\,.
    \]

\begin{proof} Abbreviate $t_1=\gamma(1)$.
Let $Z$ be a standard Brownian motion on the time interval $[t_1,{{\infty)}}$ with initialization $Z(t_1)=\tilde{Y}(t_1)$. Let
    \[
    \bar{\tau}
    \equiv\inf\Big\{t\ge t_1 : |Z(t)|\ge1\Big\}\,,
    \]
and define the stopped martingale
    \[\bar{Z}(t)\equiv Z(t\wedge\bar{\tau})
    =\tilde{Y}(t_1)+ \int_{q_0}^t \ind\{t\le\bar{\tau}\} \,dZ(t)\,.
    \]
Then, similarly as in the proof of Lemma~\ref{l:time.change}, we define $\tilde{Y}$ to be the time change of $\bar{Z}$ satisfying the desired properties. Explicitly, let 
    \[
    \bar{\gamma}(t)
    = \E[\bar{Z}(t_1+t)^2] -\E[\bar{Z}(t_1)^2] 
    = \E \int_{t_1}^{t_1+t} \ind\{t\le\bar{\tau}\}\,dt\,,
    \]
and denote its generalized inverse $\bar{\lambda}$.
We then take $\tilde{Y}(t)=\bar{Z}(t_1+\bar{\lambda}(t-t_1))$, so that
$\E[\tilde{Y}(t)^2]=t$. As in Lemma~\ref{l:time.change}, we can express
    \[
    \tilde{Y}(t)
    = \tilde{Y}(t_1) + \int_{t_1}^t \tilde{w}_t\,d\tilde{W}(t)
    \]
where we can define explicitly    \begin{align*}
    \tilde{W}(t) 
    &=\int_{t_1}^{t_1+\bar{\lambda}(t-t_1)}
	\bar{\gamma}'(u-t_1)^{1/2}\,dZ(u)\,,\\
    \tilde{w}_t
    &= \frac{\ind\{t_1+\bar{\lambda}(t-t_1)\le\tau\}}
        {\bar{\gamma}'(\bar{\lambda}(t-t_1))^{1/2}}\,.
    \end{align*}
The claim follows.
\end{proof}
\end{lem}

The next proposition implies that, averaged over time, the modified control $\tilde w_t$ provides nearly as much budget as the original control $w_t$:

\begin{ppn}\label{p:budget-loss-small}
On the time interval $[q_0,\gamma(1)]$,
for $\tilde{w}_t$ as defined by Lemma~\ref{l:time.change}, we have
    \[
        \int_{q_0}^{\gamma(1)}
        \Big(
            \bbE w_{\lambda(s)} - \bbE \tilde w_s
        \Big)_+ \,ds
        \le 2\epsilon^{1/4}
    \]
for all $\epsilon>0$ small enough.

\begin{proof} By Lemma~\ref{l:generalized.inverse.fn}, for bounding the integral we can ignore the measure-zero set of times $s$ for which $\gamma'(\lambda(s))=0$. Therefore, without loss, we restrict our attention to times $s$ for which
$\gamma'(\lambda(s))>0$. In this case, we have $\gamma'(\lambda(s)) \lambda'(s) = 1$, so the definition of $\tilde w_s$ for $s\in[q_0,\gamma(1)]$ 
can be written as $\tilde w_s = \lambda'(s)^{1/2} w_{\lambda(s)} \ind\{\lambda(s) \le \tau\}$. We therefore have
    \begin{align*}
        \int_{q_0}^{\gamma(1)}\lt(
            \bbE w_{\lambda(s)} - \bbE\tilde w_s
        \rt)_+ \,ds
        &\le
        \int_{q_0}^{\gamma(1)}
        \bbE\Big[
        w_{\lambda(s)} - \lambda'(s)^{1/2}
            w_{\lambda(s)} \ind\{ \lambda(s) \le \tau\}
    \Big]_+\,ds\\
        &\le
        \int_{q_0}^{\gamma(1)}
        \bbE\Big[ \lambda'(s)^{1/2} w_{\lambda(s)} \ind\{
        \lambda(s) > \tau\}\Big]\,ds\,,
    \end{align*}
    where the last inequality uses $\lambda' \ge 1$, as a consequence of \eqref{e:gamma.prime}. We then apply the Cauchy--Schwarz inequality to bound the above by
    \begin{align*}
      &  \bigg(
            \int_{q_0}^{\gamma(1)}
            \bbP(\lambda(s) > \tau)\,ds\bigg)^{1/2}
            \cdot\bigg(
            \int_{q_0}^{\gamma(1)}
            \lambda'(s)
            \bbE\Big[
            (w_{\lambda(s)})^2 \ind\{ \lambda(s) > \tau\}\Big]
            \,ds
        \bigg)^{1/2}\\
  &\quad      \le \bigg(
            \int_{q_0}^{\gamma(1)}
            \lambda'(s)
            \bbE\Big[
            (w_{\lambda(s)})^2 \ind\{
             \lambda(s) > \tau\}\Big]
            \,ds
        \bigg)^{1/2}.
    \end{align*}
Changing variables to $u = \lambda(s)$, we bound the last integral above by
    \begin{align*}
        \int_{q_0}^1
        \bbE[(w_u)^2 \ind\{u > \tau\}]\,du
        &= \bbE
            \int_\tau^1 (w_u)^2 \,du
       = \bbE[Y(1)^2 -\hat Y(1)^2]
        = 1 - \gamma(1)
        \le 4\epsilon^{1/2}\,,    \end{align*}
where the last step again uses
Proposition~\ref{p:gamma-1}. This proves the claim.
\end{proof}
\end{ppn}

We now define a set of ``good" times on which the controls $(\tilde b_t, \tilde \sigma_t)$ that we construct below can be successfully coupled to $(b_{\lambda(t)}, \sigma_{\lambda(t)})$. To be precise, let
\[
    I = \Big\{
        t \in [q_0,\gamma(1)] : \lambda'(t) \le 1 + \epsilon^{1/4} ~\text{and}~
        \bbE[w_{\lambda(t)}] \le \bbE[\tilde w_t] + \epsilon^{1/8}
    \Big\}.
\]
We next show this set has large measure.

\begin{lem}
    \label{l:I-large-measure}
    For sufficiently small $\epsilon$, we have $|\lambda(I)| \ge |I| \ge 1 - q_0 - \epsilon^{1/10}$.

\begin{proof}
Define the ``bad'' intervals
    \begin{align*}
    I_1 &\equiv |\{t \in [q_0,\gamma(1)] : \lambda'(t) > 1 + \epsilon^{1/4}\}|\,,\\
    I_2&\equiv   |\{t \in [q_0,\gamma(1)] : \bbE[w_{\lambda(t)}] > \bbE[\tilde w_t] + \epsilon^{1/8} \}|\,.
    \end{align*}
By Markov's inequality combined with Proposition~\ref{p:gamma-1},
    \[
    |I_1|
    \le \int_{q_0}^{\gamma(1)}
    \frac{\lambda'(s)-1}{\epsilon^{1/4}}\,ds
    =\frac{1-\gamma(1)}{\epsilon^{1/4}}
    \le \frac{4\epsilon^{1/2}}{\epsilon^{1/4}}
    = 4\epsilon^{1/4}\,.
    \]
Similarly, by Markov's inequality combined with Proposition~\ref{p:budget-loss-small},
    \[
    |I_2|
    \le
     \int_{q_0}^{\gamma(1)} \frac{\big(\E[w_{\lambda(s)}] - \E[w_s]\big)_+
     }
        {\epsilon^{1/8}}\,ds
    \le \frac{2\epsilon^{1/4}}{\epsilon^{1/8}}
    \le 2\epsilon^{1/8}\,.
    \]
Combining with Proposition~\ref{p:gamma-1} shows
    \begin{align*}
        |I| &\ge \gamma(1) - q_0 - (4\epsilon^{1/4} + 2\epsilon^{1/8}) \\
        &\ge 1 - q_0 - (4\epsilon^{1/2} + 4\epsilon^{1/4} + 2\epsilon^{1/8})
        \ge 1 - q_0 - \epsilon^{1/10}\,,
    \end{align*}
    where the last inequality holds
    for sufficiently small $\epsilon$.
    Finally, since $\lambda' \ge 1$, we have $|\lambda(I)| \ge |I|$.
\end{proof}
\end{lem}

\begin{lem}\label{l:spr.in.prf1}
For all $q_0 \in [0,1)$, we have $\sM^{\svii}(\alpha,q_0)\subseteq \sM^{\sviii}(\alpha,q_0)$.
\end{lem}

\begin{proof}
We repeat the most relevant part of Definition~\ref{d:SDE.classification} for this proof:
	\begin{center}
	\begin{tabular}[h]{r|cccccccc}
	\MSRheader\\
	\hline
	\SPRdefn\\
	\PRFONEdefn
	\end{tabular}
	\end{center}
We see that the $\svii$ and $\sviii$ classes differ only in the Ising endpoint condition: the $\svii$ class satisfies condition~\ref{it:5.1}, while the $\sviii$ class satisfies condition~\ref{it:5.2}.
Therefore, let us suppose we have controls $(q_* = q_0,b,\sigma,w,p,\zeta = \delta_0,\zeta^\Ising)$ belonging to $\Adm^{\svii}(\alpha,q_0;\epsilon,L_0,\epsilon_0)$. We will construct modified controls 
	\[(\tilde q_*=q_0, \tilde b, \tilde \sigma, \tilde w, \tilde p = p, \tilde \zeta = \delta_0, \tilde \zeta^\Ising = \zeta^\Ising) \in
	\Adm^{\sviii}(\alpha,q_0;L_0,\epsilon_0)\,.\] Note from the above choices that the modified controls will clearly satisfy conditions~\ref{it:3.4}, \ref{it:6.2}, \ref{it:7.2}, and \ref{it:8.2}.

Let $\tilde{w}$ be constructed as in Lemmas~\ref{l:time.change} and \ref{lem:Yt-process-always-is-completable}. Recall that $\tilde{w}$ is measurable with respect to the filtration $\tilde{\cF}_Y(t)$, and satisfies $\E[(\tilde{w}_t)^2]=1$ for all $t$. If $\tilde{Y}$ is defined by \eqref{e:formal.sde.ising} using $\tilde{w}$ and Brownian motion $\tilde{W}$ as above, then we will have $|\tilde{Y}(1)|=1$ almost surely. Thus the modified controls will satisfy
conditions \ref{it:2.2} and \ref{it:5.2}.

We will next construct the controls $\tilde{b},\tilde{\sigma}$, which will be
measurable with respect to an analogous time-changed filtration $\tilde \cF_X(t)$, so that condition~\ref{it:1.2} will be satisfied.  Analogously to \eqref{e:time-changed-W}, let $\tilde{B}$ be given by
    \[
        \tilde B(t) = \int_{q_0}^{\lambda(t)} \gamma'(s)^{1/2}\,dB(s)
    \]
for $t \in [q_0,\gamma(1)]$.
Then $\tilde{B}$ is a standard Brownian motion with respect to the time-changed filtration $\tilde{\cF}_X(t)=\cF_X(\lambda(t))$, for $t\in[q_0,\gamma(1)]$. We extend $\tilde B(t)$ to be a standard Brownian motion for all $t\in[q_0,1]$, and extend the filtration $\tilde \cF_X$ accordingly to all $t\in[q_0,1]$. Let $I\subseteq[q_0,\gamma(1)]$ be the interval of ``bad'' times from Lemma~\ref{l:I-large-measure}.
Recall that $s(t)\equiv (tp)'(t)^{1/2}$.
For $t\in I$, let
    \begin{align}\label{e:tilde.b}
        \tilde b_t &=
        b_{\lambda(t)} \cdot
        \min\bigg\{
            \frac{\bbE[\tilde w_t]}{\bbE[w_{\lambda(t)}]},
            1\bigg\}
 , \\
        (\tilde \sigma_t - 1) &=
        (\sigma_{\lambda(t)} - 1)
        \cdot \min\bigg\{
            \frac{\bbE[\tilde w_t]}{\bbE[w_{\lambda(t)}]}
            \cdot
            \frac{s(\lambda(t))/p(\lambda(t))^{1/2}}
                {s(t)/p(t)^{1/2}},
            1
        \bigg\}\,.\label{e:tilde.sigma}\end{align}
Then clearly $|\tilde{b}_t| \le |b_{\lambda(t)}|$, $|\tilde{\sigma}_t-1| \le|\sigma_{\lambda(t)}-1|$, and $\tilde{\sigma}_t\ge0$. Then
    $(\tilde{b}, \tilde{\sigma}, \tilde w)$ satisfies the budget constraint \ref{it:4.3} at time $t$ because $(b,\sigma,w)$ satisfies the constraint at time $\lambda(t)$. For $t\not\in I$, we simply take $(\tilde{b}_t, \tilde{\sigma}_t) = (0,1)$, which satisfies the budget constraint \ref{it:4.3} because the cost of this control is zero.

    Let $X$ be given by \eqref{e:formal.sde.main} with controls $(b,\sigma)$ and Brownian motion $B$; and let $\tilde X$ be given by \eqref{e:formal.sde.main}
    with controls
     $(\tilde{b}, \tilde{\sigma})$ and Brownian motion $\tilde{B}$.
    (Because we work under initialization condition~\ref{it:6.2}, both processes have the same initialization $X(q_0) = \tilde X(q_0) = 0$.)
    We will show that
 $X(1)$ and $\tilde{X}(1)$ are close in $L^2$, from which the conclusion will follow. To this end, recalling that $I\subseteq[q_0,\gamma(1)]$,  
let us denote $I^c \equiv [q_0,1] \setminus I$ and $\lambda(I)^c = [q_0,1] \setminus \lambda(I)$. We can then decompose
    \begin{align*}
        X(1) &= \int_{q_0}^1 p'(t)^{1/2} b_t\,dt
        + \int_{q_0}^1 s(t) \sigma_t \,dB(t) \\
        &= \int_I \lambda'(t) p'(\lambda(t))^{1/2} b_{\lambda(t)}\,dt
        + \int_I \lambda'(t)^{1/2} s(\lambda(t)) \,d\tilde{B}(t)
        + \int_I \lambda'(t)^{1/2}
        s(\lambda(t)) (\sigma_{\lambda(t)}-1) \,d\tilde{B}(t) \\
        &\qquad + \int_{\lambda(I)^c} p'(t)^{1/2} b_t\,dt
        + \int_{\lambda(I)^c} s(t) \sigma_t\,dB(t)\,,\\
        \tilde X(1) &= \int_I p'(t)^{1/2} \tilde b_t \,dt
        + \int_I s(t)\,d\tilde{B}(t)
        + \int_I s(t) (\tilde \sigma_t - 1) \,d\tilde{B}(t)
        + \int_{I^c} s(t) \tilde \sigma_t \,d\tilde{B}(t)\,,
    \end{align*}
where we recall that $\tilde b_t = 0$ for $t\in I^c$. Therefore
$\bbE[(X(1) - \tilde X(1))^2]
        \le O(1)(\Delta_1 + \ldots + \Delta_{{5}})$, where
    \begin{align*}
        \Delta_1 &= \bbE \bigg[
        \bigg(
            \int_I\Big[
            \lambda'(t) p'(\lambda(t))^{1/2} b_{\lambda(t)}
            - p'(t)^{1/2} \tilde b_t \Big] \,dt
        \bigg)^2
        \bigg] \,,\\
        \Delta_2 &=
        \int_I \lt(
            \lambda'(t)^{1/2} s(\lambda(t))
            - s(t)
        \rt)^2\,dt\,,\\
        \Delta_3 &= \bbE
            \int_I \lt(
                \lambda'(t)^{1/2}
                    s(\lambda(t))(\sigma_{\lambda(t)} - 1)
                -  s(t) (\tilde \sigma_t - 1)
            \rt)^2\,dt \,,\\
        \Delta_4 &= \bbE \bigg[\bigg(
            \int_{\lambda(I)^c} p'(t)^{1/2} b_t \,dt
        \bigg)^2\bigg] \,, \\
        \Delta_5 &= \bbE \lt[
            \int_{\lambda(I)^c} s(t) (\sigma_t)^2 \,dt
        \rt] + \bbE \lt[
            \int_{I^c} s(t) (\tilde \sigma_t)^2 \,dt
        \rt]
    \end{align*}
Recall from Lemma~\ref{l:I-large-measure} that $|\lambda(I)^c| \le |I^c| \le \epsilon^{1/10}$; we now use this to bound the last two quantities above.
    By the budget constraint \ref{it:4.3}, we have $\bbE[(b_t)^2] \le \alpha^{-1}$. Combining with the Cauchy--Schwarz inequality gives
    \[
        \Delta_4 \le
        \lt(
            \int_{\lambda(I)^c} p'(t) \,dt
        \rt)\lt(
            \int_{\lambda(I)^c} \bbE[(b_t)^2] \,dt
        \rt)
        \le \frac{|\lambda(I)^c|}{\alpha}
        \le \frac{\epsilon^{1/10}}{\alpha}\,.
    \]
Similarly, it follows from the budget constraint \ref{it:4.3} that
    \[
        \bbE[\sigma_t^2]
        \le 2\bbE\Big[(\sigma_t-1)^2 + 1\Big]
        \le 2\bigg( \frac{1}{\alpha}+1\bigg)\,.
    \]
    Moreover, since $\|p\|_{C^2([q_0,1])} \le L_0$ by condition~\ref{it:3.4}, we have ${{s(t)^2}} \le L_0+1$ for all $t$.
    Thus
    \[
        \Delta_5 \le
        4\bigg(\frac1\alpha+1\bigg) (L_0+1) \epsilon^{1/10}\,.
    \]
We now turn to bounding $\Delta_1,\Delta_2,\Delta_3$.
Recall from \eqref{e:gamma-estimate} that $\gamma$ is close to the identity function. Combining with the bound $\|p\|_{C^2([q_0,1])} \le L_0$
from condition~\ref{it:3.4}, we have
    \[
       \max\bigg\{ |p(s) - p(\gamma(s))|, |p'(s) - p'(\gamma(s))|\bigg\}
        \le L_0 |s - \gamma(s)|
        \le 4L_0\epsilon^{1/2}\,.
    \]
It follows that
    \beq\label{e:tpt-estimate}
    \begin{aligned}
        &|s(\lambda(t))^2 - s(t)^2|
        =|p(\lambda(t)) + \lambda(t) p'(\lambda(t))
        - p(t) - tp'(t) | \\
            &\qquad\le |p(\lambda(t)) - p(t)|
        + |\lambda(t) - t| p'(\lambda(t))
        + t |p'(\lambda(t)) - p'(t)|
        \le 12L_0\epsilon^{1/2}\,.
    \end{aligned}\eeq
Also recall from the definition of $I$ that for $t\in I$, $|\lambda'(t)-1| \le \epsilon^{1/4}$. It follows that
    \begin{align*}
        \Delta_2 &\le
        \int_I \Big| \lambda'(t) s(\lambda(t))^2 - s(t)^2\Big| \,dt
        \le
        \int_I
        |\lambda'(t) - 1| s(\lambda(t))^2
            \,dt
        + \int_I
        |s(\lambda(t))^2-s(t)^2|\,dt \\
        &
            \le \epsilon^{1/4} \int_I s(\lambda(t))^2\,dt
        + 12L_0\epsilon^{1/2}
        \le (L_0+1) \epsilon^{1/4} + 12L_0\epsilon^{1/2}\,.
    \end{align*}
    Similarly, for all $t\in I$, we can bound
    \begin{align*}
        &\Big|
            \lambda'(t) p'(\lambda(t))^{1/2} b_{\lambda(t)}
            - p'(t)^{1/2} \tilde b_t
        \Big| \\
        &\le
        |\lambda'(t) - 1| p'(\lambda(t))^{1/2} |b_{\lambda(t)}|
        + \Big|p'(\lambda(t))^{1/2} - p(t)^{1/2} \Big| |b_{\lambda(t)}|
        + p'(t)^{1/2} |b_{\lambda(t)} - \tilde b_t| \\
        &\le \epsilon^{1/4} (L_0)^{1/2} |b_{\lambda(t)}|
        + |p'(\lambda(t)) - p(t)|^{1/2} |b_{\lambda(t)}|
        + (L_0)^{1/2} |b_{\lambda(t)} - \tilde b_t|\,.\end{align*}
Combining with the definition \eqref{e:tilde.b} of $\tilde{b}$ gives that the above is upper bounded by
        \begin{align*}&= (L_0)^{1/2} \lt(
            3\epsilon^{1/4}
            + \bigg(1-\frac{\bbE [\tilde w_t]}{\bbE [w_{\lambda(t)}]}
                \bigg)_+
        \rt)
        |b_{\lambda(t)}|
        \le \frac{4(L_0)^{1/2}\epsilon^{1/8}}{\bbE [w_{\lambda(t)}]} |b_{\lambda(t)}|\,.
    \end{align*}
In the last inequality above, we used that $\bbE [w_{\lambda(t)}] \le 1$ and $\bbE[\tilde w_t] \ge \bbE [w_{\lambda(t)}] - \epsilon^{1/8}$ for all $t\in I$.
    Since the budget constraint \ref{it:4.3} implies $\bbE[(b_{\lambda(t)})^2] \le (\bbE w_{\lambda(t)})^2 / \alpha$, we obtain
    \[
        \Delta_1 \le
        \int_I \frac{16L_0\epsilon^{1/4} \bbE[ (b_{\lambda(t)})^2]}
            {(\bbE w_{\lambda(t)} )^2} \,dt
        \le \frac{16L_0\epsilon^{1/4}}{\alpha}.
    \]
It remains to bound $\Delta_3$, which we decompose further as
    \begin{align*}
        \Delta_3
        &\le 2(\Delta_{3,1} + \Delta_{3,2})
        \equiv 2 \bbE
            \int_I \Delta_{3,1}(t)\,dt
            +\E \int_I \Delta_{3,2}(t) \,dt\,,\\
        \Delta_{3,1}(t) &=
        \Big(\lambda'(t)^{1/2}-1\Big)^2 s(\lambda(t))^2
         (\sigma_{\lambda(t)} - 1)^2\,, \\
        \Delta_{3,2}(t) &= \Big(
            s(\lambda(t)) (\sigma_{\lambda(t)} - 1)
            - s(t) (\tilde \sigma_t - 1) \Big)^2\,.
    \end{align*}
Applying the budget constraint~\ref{it:4.3}, we can bound
    \[
    \E\Delta_{3,1}(t)
    \le
    \Big(\lambda'(t)^{1/2}-1\Big)^2
    p(\lambda(t)) \frac{ (\bbE w_{\lambda(t)})^2 }{ \alpha}
    \le
    \frac{(\lambda'(t)^{1/2}-1)^2}{\alpha}
    \le \frac{\epsilon^{1/4}}{\alpha}\,,
    \]
where the last bound is by recalling $t\in I$.  Similarly, recalling the definition~\eqref{e:tilde.sigma} of $\tilde{\sigma}$, we have
    \begin{align*}
     \E \Delta_{3,2}(t)
        &= \bigg(
        \min\bigg\{
            \frac{(\bbE \tilde w_t) p(t)^{1/2}}{ (\bbE w_{\lambda(t)}) p(\lambda(t))^{1/2}}
            , \frac{s(t)}{s(\lambda(t))}\bigg\}- 1
        \bigg)^2 s(\lambda(t))^2
    \E\Big[(\sigma_{\lambda(t)} - 1)^2\Big] \\
    &\le  \bigg(
        \min\bigg\{
            \frac{(\bbE \tilde w_t) p(t)^{1/2}}{ (\bbE w_{\lambda(t)}) p(\lambda(t))^{1/2}}
            , \frac{s(t)}{s(\lambda(t))}\bigg\}- 1
        \bigg)^2
            p(\lambda(t)) \frac{(\E w_{\lambda(t)})^2}{\alpha}\,,
    \end{align*}
having again used the budget constraint~\ref{it:4.3}.  Applying \eqref{e:tpt-estimate} gives
    \[
    \bigg(\frac{s(t)}{s(\lambda(t))}-1\bigg)^2 p(\lambda(t))
    \le \bigg| \frac{s(\lambda(t))^2-s(t)^2}{s(\lambda(t))^2}\bigg|
        p(\lambda(t))
    \le\frac{ 12 L_0 \epsilon^{1/2} p(\lambda(t))}{ p(\lambda(t)) +
        \lambda(t) p'(\lambda(t))}
    \le 12 L_0 \epsilon^{1/2}\,.
    \]
It follows that
    \[
    \E \Delta_{3,2}(t)
    \le
    \max
    \bigg\{
    \bar{\Delta}_{3,2}(t)
    \equiv
    \frac1\alpha
    \Big(
    (\E w_{\lambda(t)}) p(\lambda(t))^{1/2}
    - (\E\tilde{w}_t) p(t)^{1/2}
    \Big)_+^2
     ,
    \frac{12 L_0 \epsilon^{1/2}}{\alpha}
    \bigg\}\,.
    \]
It is straightforward to check that for any $x,y\in\R$, we can bound $(x+y)_+^2 \le 2[ x^2 + (y_+)^2]$. Thus
    \begin{align*}
        \bar{\Delta}_{3,2}(t)
        &\le \frac{2}{\alpha} \bigg\{ (\bbE\tilde w_t)^2
        \Big(
            p(t)^{1/2} - p(\lambda(t))^{1/2}
        \Big)^2
        +p(\lambda(t))
        \Big(
            \bbE [w_{\lambda(t)}] - \bbE [\tilde w_t]
        \Big)_+^2\bigg\}
        \le \frac{8L_0\epsilon^{1/2}}{\alpha} + \frac{2\epsilon^{1/4}}{\alpha}
        \,,
    \end{align*}
where we again use that $\bbE [w_{\lambda(t)}] \le \bbE[\tilde w_t] + \epsilon^{1/8}$ for $t\in I$. It follows that for fixed $L_0$, all the $\Delta_i$ tend to zero in the limit $\epsilon \to 0$. We therefore obtain
\[
    \bbW_2\Big(\Law(X(1)), \Law(\tilde X(1))\Big) \le o_{\epsilon \rightarrow 0}(1) \le \epsilon_0\,,
\]
where we recall the order of parameters \eqref{e:control.problems.limiting.regime}. This implies
\[
    \cB_{\epsilon_0}\Big(
    \cM^{\svii}(\alpha,q_0;\epsilon,L_0,\epsilon_0)\Big)
    \subseteq \cB_{2\epsilon_0}\Big(\cM^{\sviii}(\alpha,q_0;L_0,\epsilon_0)\Big)\,.
\]
Recalling the notations of Definition~\ref{d:measure.classes}, the claim follows.
\end{proof}

We turn to the deferred proof of Proposition~\ref{p:gamma-1}.
This relies on the following lemmas.

\begin{lem}\label{l:quadratic.ineq}
For $0<\epsilon < 1/2$ and any $x\in\R$,
    \[
    2\epsilon^{1/2}
    +\frac{2(|x|-1)^2}{\epsilon^{1/2}} \ge x^2-1\,.
    \]
\begin{proof}
We can suppose without loss $x\ge0$, and make the change of variables $y=x-1$. Then the desired inequality rearranges to
    \[
    \bigg(\frac{2}{\epsilon^{1/2}}-1\bigg) y^2
    - 2y+2\epsilon^{1/2}\ge 0\,.
    \]
The left-hand side is a quadratic function of $y$ which is strictly positive at $y=0$ and has negative discriminant for the range of $\epsilon$ given in the statement of the lemma. It follows that the left-hand side never reaches zero, and the claim follows.
\end{proof}
\end{lem}
\begin{lem}
    \label{lem:BM-time-outside}
Let $\acute{Y}$ be a continuous local martingale started from $\acute{Y}(0)=1$, and let $\theta$ be a stopping time with $\E\langle \acute{Y}\rangle_\theta<\infty$. If $\E [(|{{\acute{Y}}}_\theta|-1)^2] \le\epsilon$, then $\E\langle \acute{Y}\rangle_\theta \le 4\epsilon^{1/2}$.

\begin{proof}
The process
$\acute{M}(t)\equiv \acute{Y}^2(t)-\langle\acute{Y} \rangle_t$ is a continuous local martingale. It follows from e.g.\ \cite[Thm.~4.13]{MR3497465} that the stopped process $\acute{M}(t\wedge\theta)$ is a uniformly integrable true martingale, so the optional stopping theorem gives
    \[
    \E\langle \acute{Y}\rangle_\theta
    = \E\Big[ (\acute{Y}_\theta)^2-1\Big]\,.
    \]
Combining with Lemma~\ref{l:quadratic.ineq} and the hypothesis gives
    \[
    \E\langle \acute{Y}\rangle_\theta
    \le 2\epsilon^{1/2}
    + \frac{2\E[(|\acute{Y}_\theta|-1)^2]}{\epsilon^{1/2}}
    \le 4\epsilon^{1/2}\,,
    \]
as claimed.
\end{proof}
\end{lem}

\begin{proof}[\hypertarget{proof:p.gamma-1}{Proof of Proposition~\ref{p:gamma-1}}]
Recall from \eqref{e:tau}  that $\tau$ denotes the first time $t$ that $|Y_t|=1$, and we defined $\hat{Y}(1) = Y(\tau\wedge1)$. Thus, on the event $\tau>1$, we have $Y(1)=\hat{Y}(1)$. Thus, using the definition of $\gamma(t)$,
    \[
        1 - \gamma(1)
        = \bbE\Big[Y(1)^2 - \hat Y(1)^2\Big]
        = \bbE \lt[\int_\tau^1 (w_t)^2 \,dt
        ; \tau < 1 \rt]\,.\]
We will define a process $\acute{Y}$ which, on the event $\tau<1$, captures the evolution of $Y$ on the time interval $[\tau,1]$. First, let us express $\tau = \tau_+ \wedge \tau_-$ where
    \[\begin{aligned}
        \tau_+ = \inf\{ t \ge q_0 : Y_t \ge 1\}\,, \\
        \tau_- = \inf\{ t \ge q_0 : Y_t \le -1 \}\,.\end{aligned}
    \]
We then define, for $t\in[\tau,1]$ on the event $\tau<1$, 
    \[
        \acute{Y}(t) = 1 +
        \Big(\ind\{t \ge \tau = \tau_+\} - \ind\{t \ge \tau = \tau_-\}
            \Big)
        \int_\tau^t w_s \,dW(s)\,.
    \]
Then $\acute{Y}$ evolves exactly as $Y$ on the time interval $[\tau,1]$, up to a reflection in the case $Y_\tau=-1$. We have
    \[
        \bbE\Big[(|\acute{Y}(1)| - 1)^2\Big]
        = \bbE\Big[(|Y(1)| - 1)^2; \tau<1\Big]
        \le \bbE\Big[(|Y(1)| - 1)^2\Big]
        \le \epsilon\,,
    \]
and applying Lemma~\ref{lem:BM-time-outside} gives
    \[\E\langle \acute{Y}\rangle_1
    =
        \bbE\bigg[ \int_\tau^1 (w_t)^2 \,dt ; \tau<1 \bigg]
        \le 4\epsilon^{1/2}.
    \]
This proves the claim.
\end{proof}

Finally, we prove the last inclusion asserted at the start of this subsection:

\begin{lem}   \label{l:ideal.contains.std.gaussian}
We have $\cN(0,1) \in \sM^{\sideal}(\alpha,0)$.

\begin{proof}
We repeat the part of Definition~\ref{d:SDE.classification} specifying the $\sideal$ class:
	\begin{center}
	\begin{tabular}[h]{r|cccccccc}
	\MSRheader\\
	\hline
	\SUPERPERFECTdefn
	\end{tabular}
	\end{center}
We will exhibit controls $(q_*,b,\sigma,q,p,\zeta,\zeta^{\Ising}) \in \Adm^{\sideal}(\alpha,0)$ such that $\mu(q_*,b,\sigma,p,\zeta) = \cN(0,1)$. To this end, we let
    \[
        (q_*,b,\sigma,p,\zeta) = (0,0,1,1,\delta_0),
    \]
    so that the process $X(t)$ defined by \eqref{e:formal.sde.main} is simply
    \[
        X(t) = \int_0^t \,dB(u) = B(t)\,,
    \]
and clearly we have 
	\[
	\mu(q_*,b,\sigma,p,\zeta)
	=
	\Law(X(1)) = \cN(0,1)\,.\]
    We then set $\zeta^{\Ising} = \delta_0$, and $w$ to be any progressively measurable control with $\bbE[(w_t)^2] = 1$ for all $t\in [0,1]$ such that $Y(t)$ defined in \eqref{e:formal.sde.ising} satisfies $|Y(t)| = 1$ almost surely.
    The fact that such a $w$ exists follows from Lemma~\ref{lem:Yt-process-always-is-completable} (where we set $q_0 = \gamma(1) = 0$, which makes $\tilde Y, \tilde W, \tilde w, \tilde \cF_Y$ therein degenerate). 

    The above controls clearly satisfy conditions \ref{it:1.2}, \ref{it:2.2}, \ref{it:3.5}, \ref{it:5.2}, \ref{it:6.2}, \ref{it:7.2}, and \ref{it:8.2}.
    Moreover, for the above choices of $b,\sigma$, we have $\Cost(t) = 0$, so condition~\ref{it:4.3} is satisfied.
    Thus these controls belong to $\Adm^{\sideal}(\alpha,0)$, and the claim follows.
\end{proof}
\end{lem}

\subsection{Simplification to $p\equiv 1$ for symmetrized endpoint measures}
\label{ss:sdes.symmetric.endpoints}

In this subsection we give the 
\hyperlink{p:thm.control.problems.sym}{proof of Theorem~\ref{thm:control.problems.sym}}.
We adopt the notations of \S\ref{ss:control.problems.measure.classes}, and recall the conditions on control parameters from Definitions \ref{d:SDE.conds}--\ref{d:SDE.classification}.
We now introduce several more classes of controls and measures for the symmetrized setting, where we exclusively take $q_0 = 0$.

\begin{dfn}[symmetrized classification, cf.\ Defn.~\ref{d:SDE.classification}]
\label{d:SDE.classification.sym}
Let $q_0 = 0$. We classify the control parameters $(q_*,b,\sigma,w,p,\zeta,\zeta^\Ising)$ identically as in Definition~\ref{d:SDE.classification}, except for modified conditions on $p$:
\begin{center}
\begin{tabular}[h]{r|cccccccc}
\MSRheader\\
\hline
\IAMPdefn{sy}\\
\BOGPdefn{sy}\\
\FAIRdefn{sy}\\
\TMPONEdefn{sy}\\
\TMPTWOdefn{sy}\\
\GOODdefn{sy}\\
\TMPTHRdefn{sy}\\
\GREATdefn{sy}\\
\SPRdefn{sy}\\
\SUPERPERFECTdefn{sy}\\
\PRFTWOdefn{sy}\\
\PRFTHRdefn{sy}
\end{tabular}
\end{center}%
Compared to Definition~\ref{d:SDE.classification}, we note that the $\sviii$, $\six$, $\sxii$, and $\sidealast$ classes are no longer present, and neither are the parameters $L_0, \epsilon_0$. 
We denote the subsets of $(q_*,b,\sigma,w,p,\zeta,\zeta^\Ising)$ satisfying the above conditions by
\beq
    \label{e:star.options.sym}
    \Adm^{\star,\sym}(\alpha,0;L,\epsilon),\qquad
    \star \in \bigg\{
    \begin{array}{c}
    \IAMP,\BOGP,\si,\sii,\siii,\siv,\\
    \sv,\svi,\svii,\sideal,\sx,\sxi
    \end{array}\bigg\}.
\eeq
Similarly to Definition~\ref{d:SDE.classification}, if either of $L,\epsilon$ does not appear in the definition of a class $\star$, then we take $\Adm^{\star,\sym}(\alpha,0;L,\epsilon)$ to be the same for all values of that parameter and sometimes omit it from the notation, writing e.g. $\Adm^{\svii,\sym}(\alpha,0;L,\epsilon) = \Adm^{\svii,\sym}(\alpha,0;\epsilon)$.
\end{dfn}

For $\mu \in \cP_2(\bbR_{\ge 0})$ and $\cM \in \cP_2(\bbR_{\ge 0})$, we define the point-to-set distance $\bbW_2(\mu,\cM)$ as in \eqref{def:W2-to-set}, and we define the $\bbW_2$-ball (cf.\ \eqref{e:W2.ball.notation}) 
\beq\label{e:W2.ball.notation.sym}
    \cB_{\epsilon,\sym}(\cM) = \Big\{
        \mu \in \cP_2(\bbR_{\ge 0}) : \bbW_2(\mu,\cM) \le \epsilon
    \Big\}\,.
\eeq
We now use the classification of Definition~\ref{d:SDE.classification.sym} to define corresponding sets of achievable endpoint measures in the symmetrized setting, as follows:

\begin{dfn}[symmetrized measure classes, cf.\ Defn.~\ref{d:measure.classes}]
    \label{d:measure.classes.sym}
Analogously to \eqref{e:control.problems.endpoint.measures}, we let
    \[
        \mu^\sym(q_*,b,\sigma,p,\zeta) = \Law(|X(1)|)
    \]
    where $X$ is defined by \eqref{e:formal.sde.main}.
    For $\star$ as in \eqref{e:star.options.sym}, we denote
    \beq\label{e:measure.classes.sym.cM}
        \cM^{\star,\sym}(\alpha,0;L,\epsilon) = \bigg\{
            \mu^\sym(q_*,b,\sigma,p,\zeta) :
            (q_*,b,\sigma,w,p,\zeta,\zeta^\Ising) \in \Adm^{\star,\sym}(\alpha,0;L,\epsilon)
        \bigg\}\,,
    \eeq
We may omit the parameters $L,\epsilon$ if they do not appear in the definition of the class $\star$. We further denote
    \[
        \sM^{\star,\sym}(\alpha,0) =
        \adjustlimits
        \bigcap_{\epsilon>0}
        \bigcup_{L>0}
        \cB_{\epsilon,\sym}\Big(\cM^{\star,\sym}(\alpha,0;L,\epsilon)\Big)\,,
    \]
    for $\cB_{\epsilon,\sym}$ as defined by \eqref{e:W2.ball.notation.sym}.
\end{dfn}

\begin{rmk}\label{r:measure.classes.simplification.sym}
    Analogously to Remark~\ref{r:measure.classes.simplification}, we have 
    \[
    \sM^{\sideal,\sym}(\alpha,0) 
    = \ocM^{\Ising,\sym}(\alpha)\,\]
for $\ocM^{\Ising,\sym}(\alpha)$ specified in Definition~\ref{d:achievable-msrs-sym}.
\end{rmk}

Recall the notation \eqref{e:symmetrize.measure}.
For any subset of measures
$\cM \subseteq \cP_2(\bbR)$, define
\beq\label{e:symmetrize.set.of.measures}
    \sym(\cM) 
    \equiv
    \Big\{\sym(\mu) : \mu \in \cM
    \Big\}\,.
\eeq
Note that $\sym$ is a contraction on $\cP_2(\bbR)$ metrized by $\bbW_2$, and therefore also a contraction on the family of subsets of $\cP_2(\bbR)$, as metrized by the Hausdorff distance \eqref{e:hausdorff}.

\begin{proof}[\hypertarget{p:cor.feasible.distributions.continuity.sym}{Proof of Corollary~\ref{cor:feasible.distributions.continuity} for $\alpha \mapsto  \ocM^{\Ising,\sym}(\alpha)$}]
    Define $\BUDGETS(\alpha)$ and $\cM(0;p,f)$ as in the proof of this corollary for $\alpha \mapsto \ocM^{\Ising,\sym}(\alpha)$ 
    above.
    Then,
    \[
        \ocM^{\Ising,\sym}(\alpha)
        =
        \overline{
            \bigcup_{f\in \BUDGETS(\alpha)}
            \sym(\cM(0;p\equiv 1,f))
        }\,,
    \]
    where the overline denotes closure with respect to the $\bbW_2$ metric.
    Similarly, for any $\alpha'$, 
    \[
        \ocM^{\Ising,\sym}(\alpha')
        =
        \overline{
            \bigcup_{f\in \BUDGETS(\alpha)}
            \sym(\cM(0;p\equiv 1,\alpha f / \alpha'))
        }\,.
    \]
    The conclusion follows by Proposition~\ref{p:continuity.in.auxilliary.budget}, since $\sym$ is a contraction.
\end{proof}

Recalling the notation  \eqref{e:WERR.sym.notation}, we now define (cf.\ \eqref{e:cM.Lip.sym}--\eqref{e:sM.Lip.sym})
\begin{align}
    \label{e:cM.Lip.0.sym}
    \cM^{\Lip,\sym}(\alpha,0;L,\epsilon) &\equiv \left\{
        \begin{array}{ll}
        \mu \in \cP_2(\bbR_{\ge 0}): 
        & \textup{$\exists$ $L$-Lipschitz algorithms $(\cA_{N_j})_{j\ge 1}$ with} \\
        & \textup{$\WERR^\sym(\cA_{N_j},\mu)
         \le \epsilon$,
         $\chi_{\cA_{N_j}}(1-\epsilon) \le \epsilon$,}\\
         &\textup{and
        $\lim_{j\to\infty} \chi_{\cA_{N_j}}(0) = 0$}
        \end{array}
    \right\}\,,\\
    \label{e:sM.Lip.0.sym}
    \sM^{\Lip,\sym}(\alpha,0) 
    &\equiv \bigcap_{\epsilon > 0} \bigcup_{L > 0} \cM^{\Lip,\sym}(\alpha,0;L,\epsilon)\,.
\end{align}
We now give the proof of Theorem~\ref{thm:control.problems.sym}.
Like Theorem~\ref{thm:control.problems.main}, this is proved through a cycle of inclusions.
Most of these steps are proved by an identical argument to a step above, combined with the fact that $\sym$ is a contraction.
For these steps, we omit the proof and simply refer to the appropriate result.
The proofs of the remaining steps are presented below.
\begin{proof}[\hypertarget{p:thm.control.problems.sym}{Proof of Theorem~\ref{thm:control.problems.sym}}, assuming remainder of this section]
    We have the inclusions:
    \begin{itemize}
        \item $\sM^{\Lip,\sym}(\alpha,0) \subseteq \sM^{\BOGP,\sym}(\alpha,0)$ (Lemma~\ref{l:Lip.in.BOGP.sym} below);
        \item $\sM^{\BOGP,\sym}(\alpha,0) \subseteq \sM^{\si,\sym}(\alpha,0)$ (trivial);
        \item $\sM^{\si,\sym}(\alpha,0) \subseteq \sM^{\sii,\sym}(\alpha,0)$ (analogous to Lemma~\ref{l:fair.in.tmp1});
        \item $\sM^{\sii,\sym}(\alpha,0) \subseteq \sM^{\siii,\sym}(\alpha,0)$ (analogous to Lemma~\ref{l:tmp1.in.tmp2});
        \item $\sM^{\siii,\sym}(\alpha,0) \subseteq \sM^{\siv,\sym}(\alpha,0)$ (analogous to Lemma~\ref{l:tmp2.in.good});
        \item $\sM^{\siv,\sym}(\alpha,0) \subseteq \sM^{\sv,\sym}(\alpha,0)$ (analogous to Lemma~\ref{l:good.in.tmp3});
        \item $\sM^{\sv,\sym}(\alpha,0) \subseteq \sM^{\svi,\sym}(\alpha,0)$ (analogous to Proposition~\ref{p:continuity.in.auxilliary.budget});
        \item $\sM^{\svi,\sym}(\alpha,0) \subseteq \sM^{\svii,\sym}(\alpha,0)$ (Lemma~\ref{l:great.in.spr.ideal.in.prf2.sym} below);
        \item $\sM^{\svii,\sym}(\alpha,0) \subseteq \sM^{\sideal,\sym}(\alpha,0)$ (analogous to Lemma~\ref{l:spr.in.prf1});
                \item $\sM^{\sideal,\sym}(\alpha,0) \subseteq \sM^{\sx,\sym}(\alpha,0)$ (Lemma~\ref{l:great.in.spr.ideal.in.prf2.sym} below);
        \item $\sM^{\sx,\sym}(\alpha,0) \subseteq \sM^{\sxi,\sym}(\alpha,0)$ (analogous to Lemma~\ref{l:prf2.in.prf3});
        \item $\sM^{\sxi,\sym}(\alpha,0) \subseteq \sM^{\IAMP,\sym}(\alpha,0)$ (analogous to Lemma~\ref{l:prf3.in.IAMP});
        \item $\sM^{\IAMP,\sym}(\alpha,0) \subseteq \sM^{\Lip,\sym}(\alpha,0)$ (Lemma~\ref{l:IAMP.in.Lip.sym} below).
    \end{itemize}
    This implies that for any $\star$ in \eqref{e:star.options.sym},
    \beq
        \label{e:sM.star.coincide.sym}
        \sM^{\star,\sym}(\alpha,0) = \sM^{\Lip,\sym}(\alpha,0)\,,
    \eeq
    and thus all of these $\sM^{\star,\sym}(\alpha,0)$ coincide.
    For $\sM^{\Lip,\sym}(\alpha)$ defined in \eqref{e:sM.Lip.sym}, we further have the inclusions:
    \begin{itemize}
        \item $\sM^{\Lip,\sym}(\alpha) \subseteq \sM^{\IAMP,\sym}(\alpha,0)$ (Proposition~\ref{p:Lip.general.in.IAMP.sym} below);
        \item $\sM^{\Lip,\sym}(\alpha,0) \subseteq \sM^{\Lip,\sym}(\alpha)$ (trivial).
    \end{itemize}
    Thus $\sM^{\Lip,\sym}(\alpha)$ also coincides with the sets in \eqref{e:sM.star.coincide.sym}, and in particular equals $\sM^{\sideal,\sym}(\alpha,0)$.
    The result follows because $\sM^{\sideal,\sym}(\alpha,0) = \ocM^{\Ising,\sym}(\alpha)$, as was observed in Remark~\ref{r:measure.classes.simplification.sym}.
\end{proof}

\noindent We turn to the inclusions whose proofs are not analogous to above and begin with two preparatory lemmas.
\begin{lem}
    \label{l:control.problems.sym.preimage.compact}
    For any $\mu \in \cP_2(\bbR_{\ge 0})$, the preimage
    \[
        \sym^{-1}(\mu) \equiv \Big\{
            \tilde \mu \in \cP_2(\bbR) : \sym(\tilde\mu) = \mu
        \Big\}
    \]
    is compact with respect to the $\bbW_2$ metric.

\begin{proof}
The set $\sym^{-1}(\mu)$ is clearly closed. Therefore, it suffices to show that any sequence $(\tilde \mu_n)_{n\ge1}$ in $\sym^{-1}(\mu)$ has a $\bbW_2$-convergent subsequence. The collection of measures $\sym^{-1}(\mu)$ is clearly tight, since we have 
	\[
	\tilde{\mu}([-K,K])
	=\mu([0,K]) 
	\stackrel{K\to\infty}{\longrightarrow} 1\,,
	\]
uniformly over all $\tilde{\mu}\in \sym^{-1}(\mu)$. 
    By Prohorov's theorem, any sequence $(\tilde \mu_n)_{n\ge 1}$ in $\sym^{-1}(\mu)$ must have a subsequence $(\tilde \mu_{n_j})_{j\ge 1}$ that converges weakly to some measure $\tilde{\mu}$, which is also clearly in $\sym^{-1}(\mu)$. 
By Skorokhod's theorem, we have can find random variables $X_n,X$ such that $X_n$ is distributed according to $\tilde{\mu}_{n_j}$, $X$ is distributed according to $\tilde{\mu}$, and $X_n$ converges to $X$ almost surely as $n\to\infty$. Then, for any finite $M$, we have
	\[
	\E\Big[
	(X_n-X)^2; \max\{|X_n|,|X|\} \le M
	\Big]
	\le \E\Big[\max\Big\{(X_n-X)^2
	; (2M)^2
	\Big\}
	\Big] \stackrel{n\to\infty}{\longrightarrow} 0
	\]
by the bounded convergence theorem. On the other hand,	since $|X_n|$ and $|X|$ are distributed as samples from $\mu$, we have
	\begin{align*}
	&\E\Big[
	(X_n-X)^2; \max\{|X_n|,|X|\} 
		> M
	\Big]
	\le 2\Big\{
	\E[(X_n)^2;|X_n|>M]
	+\E[X^2;|X|>M]\Big\} \\
	&\qquad= 4 \E[X^2;|X|>M]
	\le o_M(1)\,.
	\end{align*}
By sending $M\to\infty$ we conclude $\tilde \mu_{n_j} \rightarrow \tilde \mu$ in $\bbW_2$, from which the compactness claim follows.
\end{proof}
\end{lem}

Similarly, define (cf.\ \eqref{e:cM.Lip.0.sym}--\eqref{e:sM.Lip.0.sym})
    \begin{align}
        \label{e:cM.Lip.0.chaos}
        \cM^{\Lip,\chaos}(\alpha,0;L,\epsilon) &\equiv \left\{
            \begin{array}{ll}
            \mu \in \cP_2(\bbR):
            & \textup{$\exists$ $L$-Lipschitz algorithms $(\cA_{N_j})_{j\ge 1}$ with} \\
             &  \textup{$\WERR(\cA_{N_j},\mu) \le \epsilon$,
             $\chi_{\cA_{N_j}}(1-\epsilon) \le \epsilon$,}\\
            &\textup{and $\lim_{j\to\infty} \chi_{\cA_{N_j}}(0) = 0$}
            \end{array}
        \right\}\\
        \label{e:sM.Lip.0.chaos}
        \sM^{\Lip,\chaos}(\alpha,0) 
        &\equiv \bigcap_{\epsilon > 0} \bigcup_{L > 0} \cM^{\Lip,\chaos}(\alpha,0;L,\epsilon)\,.
    \end{align}
We next prove some basic relations among these classes: 

\begin{lem}\label{l:control.problems.Lip.to.Lip.sym}
    For $\sM^{\Lip}(\alpha)$ and $\sM^{\Lip,\sym}(\alpha)$ as defined by \eqref{e:sM.Lip} and \eqref{e:sM.Lip.sym}, and the symmetrization operation $\sym$ as defined by \eqref{e:symmetrize.set.of.measures}, we have \[\sM^{\Lip,\sym}(\alpha) = \sym(\sM^{\Lip}(\alpha))\,.\]
Similarly, for $\sM^{\Lip,\sym}(\alpha,0)$ and $\sM^{\Lip,\chaos}(\alpha,0)$ as defined by \eqref{e:sM.Lip.0.sym} and \eqref{e:sM.Lip.0.chaos}, we have
	\[\sM^{\Lip,\sym}(\alpha,0)
	=\sym(\sM^{\Lip,\chaos}(\alpha,0))\,.\]

\begin{proof}
   By the definitions of these sets, we have
    \[
        \sym(\sM^{\Lip}(\alpha))
        \stackrel{\eqref{e:sM.Lip}}{=} 
        \sym\bigg(
            \bigcap_{\epsilon>0} \bigcup_{L>0}
            \cM^{\Lip}(\alpha;L,\epsilon)
        \bigg)
        \subseteq 
        \bigcap_{\epsilon>0} \bigcup_{L>0}
        \sym\Big(\cM^{\Lip}(\alpha;L,\epsilon)\Big)
        \stackrel{\eqref{e:sM.Lip.sym}}{=} \sM^{\Lip,\sym}(\alpha)\,.
    \]
Similarly, we have
	\begin{align*}
	&\sym(\sM^{\Lip,\chaos}(\alpha,0))
	\stackrel{\eqref{e:sM.Lip.0.chaos}}{=}\sym\bigg(
	\bigcap_{\epsilon > 0} \bigcup_{L > 0} \cM^{\Lip,\chaos}(\alpha,0;L,\epsilon)\bigg)\\
	&\qquad\subseteq
	\bigcap_{\epsilon > 0} \bigcup_{L > 0} \sym\Big(\cM^{\Lip,\chaos}(\alpha,0;L,\epsilon)\Big)
	\stackrel{\eqref{e:sM.Lip.0.sym}}{=}
	\sM^{\Lip,\sym}(\alpha,0)\,.
	\end{align*}
Therefore it suffices to show the reverse inclusions.

    Consider any $\mu \in \sM^{\Lip,\sym}(\alpha)$ and a sequence $\epsilon_n\downarrow0$.
    By the definition \eqref{e:sM.Lip.sym} of $\sM^{\Lip,\sym}(\alpha)$, for any $\epsilon_n$ there exists $L_n$ and a sequence of $L_n$-Lipschitz algorithms $\cA_{n,j}$ (acting in $N_j$ dimensions) such that
    \[\WERR^\sym(\cA_{n,j},\mu)
        =\max\bigg\{
            \bbW_2\Big(
            \sym(\mu(\cA_{n,j})),\mu
            \Big),
       \bbW_2\Big(\mu^\Ising(\cA_{n,j}),\cP(\{\pm 1\})\Big)
        \bigg\} \le \epsilon_n\,.
    \]
This implies the existence of $\tilde{\mu}_{n,j} \in \sym^{-1}(\mu)$ such that
	\[
	\WERR(\cA_{n,j},\tilde{\mu}_{n,j})
	=\max\bigg\{
            \bbW_2\Big(
            \mu(\cA_{n,j}),
            \tilde{\mu}_{n,j}
            \Big),
       \bbW_2\Big(\mu^\Ising(\cA_{n,j}),\cP(\{\pm 1\})\Big)
        \bigg\} \le \epsilon_n\,.
	\]
    In light of Lemma~\ref{l:control.problems.sym.preimage.compact}, by passing to a subsequence of the $j's$, we can ensure that
 $\tilde{\mu}_{n,j}$ converges in $\bbW_2$ to some $\tilde\mu^n \in \sym^{-1}(\mu)$ as $j\to\infty$.
    Then
    \[
        \bbW_2\Big(\mu(\cA_{n,j}),
        	\tilde{\mu}_n\Big) \le \epsilon_n + \bbW_2\Big(\tilde{\mu}_{n,j},\tilde{\mu}_n\Big) \le 2\epsilon_n\,,
    \]
    where the last inequality holds after omitting finitely many $j$'s.
    By passing to a subsequence of the $\epsilon_n$ we can ensure that as $n\to\infty$, $\tilde\mu^n$ converges in $\bbW_2$ to $\tilde\mu \in \sym^{-1}(\mu)$.
    Then
    \[
        \bbW_2\Big(\mu(\cA_{n,j}),\tilde \mu\Big) \le \tilde\epsilon_n \equiv 2\epsilon_n + \bbW_2(\tilde{\mu}_n, \tilde\mu)\,.
    \]
Recalling the definition \eqref{e:cM.Lip} of $\cM^{\Lip}(\alpha;L,\epsilon)$, it follows that
    \[
        \tilde\mu 
        \in \cM^{\Lip}(\alpha;L_n,\tilde\epsilon_n) 
        \subseteq \bigcup_{L>0} \cM^{\Lip}(\alpha;L,\tilde\epsilon_n)\,.
    \]
    Since the above holds for all $n$ and $\lim_{n\to\infty} \tilde\epsilon_n = 0$, taking an intersection over $n$ yields $\tilde\mu \in \sM^{\Lip}(\alpha)$, as defined by \eqref{e:sM.Lip}.
    This shows $\mu \in \sym(\sM^{\Lip}(\alpha))$, and thus $\sM^{\Lip,\sym}(\alpha) \subseteq \sym(\sM^{\Lip}(\alpha))$. This proves the first assertion of the lemma.

    For the remaining inclusion $\sM^{\Lip,\sym}(\alpha,0) \subseteq \sym(\sM^{\Lip,\chaos}(\alpha,0))$, consider any $\mu \in \sM^{\Lip,\sym}(\alpha,0)$. Repeat the above argument to obtain algorithms $\cA_{n,j}$ be as above; these satisfy $\chi_{\cA_{n,j}}(1-\epsilon) \le \epsilon$, so the above argument shows
    \[
        \tilde\mu 
        \in \cM^{\Lip,\chaos}(\alpha,0;L_n,\tilde\epsilon_n) 
        \subseteq \bigcup_{L>0} \cM^{\Lip,\chaos}(\alpha,0;L,\tilde\epsilon_n)\,.
    \]
    Taking an intersection over $n$ proves $\tilde\mu \in \sM^{\Lip,\chaos}(\alpha,0)$, which implies $\mu \in \sym(\sM^{\Lip,\chaos}(\alpha))$. This proves the second assertion of the lemma.
\end{proof}
\end{lem}

\begin{ppn}\label{p:Lip.general.in.IAMP.sym}
Recalling \eqref{e:sM.Lip.sym} and Definitions \ref{d:SDE.classification.sym}--\ref{d:measure.classes.sym}, we have
	\[\sM^{\Lip,\sym}(\alpha) \subseteq \sM^{\IAMP,\sym}(\alpha,0)\,.\]

\begin{proof}
Consider any $\mu \in \sM^{\Lip,\sym}(\alpha)$, as defined by \eqref{e:sM.Lip.sym}. 
    By Lemma~\ref{l:control.problems.Lip.to.Lip.sym}, we can find $\tilde\mu \in \sM^{\Lip}(\alpha)$ with $\sym(\tilde\mu) = \mu$.
It was shown in the \hyperlink{p:thm.control.problems.main}{proof of Theorem~\ref{thm:control.problems.main}} that $\sM^{\Lip}(\alpha) = \sM^{\IAMP}(\alpha,0)$, so $\tilde\mu \in \sM^{\IAMP}(\alpha,0)$.
Therefore, for any $\epsilon>0$ we can find $L>0$ and controls 
    \beq
        \label{e:Lip.general.in.IAMP.sym.controls}
        (q_*,b,\sigma,w,p,\zeta,\zeta^\Ising) \in \Adm^{\IAMP}(\alpha,0;L,\epsilon)
    \eeq
such that (recalling that $\sym$ is a contraction) we have
    \[
        \bbW_2\Big(\mu, \mu^{\sym}(q_*,b,\sigma,p,\zeta)\Big)
        \le\bbW_2\Big(\tilde \mu, \mu(q_*,b,\sigma,p,\zeta)\Big)
         \le \epsilon\,.
    \]
We see from Definitions~\ref{d:SDE.classification} and \ref{d:SDE.classification.sym}
that the sets of controls $\Adm^{\IAMP}(\alpha,0;L,\epsilon)$ and $\Adm^{\IAMP,\sym}(\alpha,0;L,\epsilon)$ are in fact identical, so the controls \eqref{e:Lip.general.in.IAMP.sym.controls} belong to $\Adm^{\IAMP,\sym}(\alpha,0;L,\epsilon)$, and 
    \[
        \mu \in \cB_\epsilon\Big(\cM^{\IAMP,\sym}(\alpha,0;L,\epsilon)\Big)\,.
    \]
    Since for all $\epsilon > 0$, there exists $L > 0$ so that this holds, we conclude $\mu \in \sM^{\IAMP,\sym}(\alpha,0)$.
    As this holds for all $\mu \in \sM^{\Lip,\sym}(\alpha)$, the claim follows.
\end{proof}
\end{ppn}

Compared to the inclusion $\sM^{\IAMP}(\alpha,q_0) \subseteq \sM^{\Lip}(\alpha,0)$ above (see Lemma~\ref{l:IAMP.in.Lip}), the main difference in the inclusion $\sM^{\IAMP,\sym}(\alpha,0) \subseteq \sM^{\Lip,\sym}(\alpha,0)$ is that the set $\cM^{\Lip,\sym}(\alpha,0;L,\epsilon)$ defined in \eqref{e:cM.Lip.0.sym} has the additional constraint \[\chi_{\cA_{N_j}}(1-\epsilon) \le \epsilon\,.\]
We show this inclusion using results from \S\ref{subsec:p=1-IAMP}, which imply that measures in $\sM^{\IAMP,\sym}(\alpha,0)$ can be attained by algorithms with this additional property.
This additional property passes through the proof of the inclusion $\sM^{\Lip,\sym}(\alpha,0) \subseteq \sM^{\BOGP,\sym}(\alpha,0)$; see Lemma~\ref{l:Lip.in.BOGP.sym} below. 
As a result, the controls in the $\Adm^{\BOGP,\sym}$ satisfy condition~\ref{it:3.6}, which is stronger than the condition~\ref{it:3.2} for the $\Adm^{\BOGP}$ class in that $p$ must satisfy $p(\epsilon) \ge 1-\epsilon$.
\begin{lem}\label{l:IAMP.in.Lip.sym}
We have $\sM^{\IAMP,\sym}(\alpha,0) \subseteq \sM^{\Lip,\sym}(\alpha,0)$.
\begin{proof} Consider any $\mu \in \cM^{\IAMP,\sym}(\alpha,0;L,\epsilon)$.
    Then, there exists
    \[
        (q_*,b,\sigma,w,p,\zeta,\zeta^\Ising) \in \Adm^{\IAMP,\sym}(\alpha,0;L,\epsilon)
    \]
    such that $\mu^{\sym}(q_*,b,\sigma,p,\zeta) = \mu$.
    Let $\tilde p$ denote the restriction of $p$ to $[q_*,1] \subseteq [q_0,1]$.
    Note that $(b,\sigma,w,\tilde p,\zeta,\zeta^\Ising)$ satisfies the assumptions of
    Theorem~\ref{thm:IAMP-main}\ref{i:IAMP-main-chaotic} with $q_*$ in place of $q_0$.
    Indeed, conditions~\ref{it:1.1}, \ref{it:2.1}, \ref{it:3.1}, \ref{it:4.2}, \ref{it:5.1}, and \ref{it:6.1} provide assumptions~(\ref{i:IAMP-main-coefs}), (\ref{i:IAMP-main-diffusivity}), (\ref{i:IAMP-main-p}), (\ref{i:IAMP-main-budget}), (\ref{i:IAMP-main-endpt}), and (\ref{i:IAMP-main-init}) of Theorem~\ref{thm:IAMP-main}.

    By Theorem~\ref{thm:IAMP-main}\ref{i:IAMP-main-chaotic}, there exists $\iota = \iota(\epsilon)$ with $\lim_{\epsilon\rightarrow 0} \iota(\epsilon) = 0$ such that for all (sufficiently large) $N$, there exists a $C(L,\epsilon)$-Lipschitz algorithm $\cA_N$ with
    \[
        \WERR^\sym(\cA_N,\mu)=
        \max\bigg\{
            \bbW_2\Big(\sym(\mu(\cA_N)), \mu\Big),
            \bbW_2\Big(\mu^{\Ising}(\cA_N), \cP(\{\pm 1\})\Big)
        \bigg\}
        \le \iota\,,
    \]
    which furthermore satisfies $\chi_{\cA_N}(0) = 0$ and $\chi_{\cA_N}(1-\iota) \le \iota$.
    This implies
    \[
        \cM^{\IAMP,\sym}(\alpha,0;L,\epsilon) \subseteq \cM^{\Lip,\sym}(\alpha,0;C(L,\epsilon),\iota(\epsilon))\,.
    \]
    The rest of the proof is identical to that of Lemma~\ref{l:IAMP.in.Lip}.
\end{proof}
\end{lem}

Recalling from \eqref{e:measure.classes.sym.cM}
the definition of $\cM^{\BOGP,\sym}(\alpha,0;L,\epsilon)$. Now 
define
    \[
        \cM^{\BOGP,\chaos}(\alpha,0;L,\epsilon) \equiv \bigg\{
            \mu(q_*,b,\sigma,p,\zeta) :
            (q_*,b,\sigma,w,p,\zeta,\zeta^\Ising) \in \Adm^{\BOGP,\sym}(\alpha,0;L,\epsilon)
        \bigg\}\,.
    \]
By comparison with \eqref{e:measure.classes.sym.cM} we see that
	\beq\label{e:BOGP.sym.relation}
	\sym\Big(\cM^{\BOGP,\chaos}(\alpha,0;L,\epsilon)\Big) = \cM^{\BOGP,\sym}(\alpha,0;L,\epsilon)\,.
	\eeq
We will use this in what follows. 

\begin{lem}\label{l:Lip.in.BOGP.sym}
We have $\sM^{\Lip,\sym}(\alpha,0) \subseteq \sM^{\BOGP,\sym}(\alpha,0)$.

\begin{proof}
We repeat the most relevant part of Definition~\ref{d:SDE.classification.sym}:
\begin{center}
\begin{tabular}[h]{r|cccccccc}
\MSRheader\\
\hline
\BOGPdefn{sy}\\
\end{tabular}
\end{center}
    By Lemma~\ref{l:control.problems.Lip.to.Lip.sym}, it suffices to show $\sym(\sM^{\Lip,\chaos}(\alpha,0)) \subseteq \sM^{\BOGP,\sym}(\alpha,0)$, where $\sM^{\Lip,\chaos}(\alpha,0)$ is defined in \eqref{e:sM.Lip.0.chaos}.
    Consider any $\mu \in \cM^{\Lip,\chaos}(\alpha,0;L,\epsilon)$ (recall \eqref{e:cM.Lip.0.chaos}), where $L$ is sufficiently large depending on $\alpha,\epsilon$.
    There exists a sequence of $L$-Lipschitz $((\cA_{N_j})^\circ)_{j\ge 1}$ such that
    \[
    \WERR((\cA_{N_j})^\circ,\mu)
    =
        \max\bigg\{
            \bbW_2\Big(\mu((\cA_{N_j})^\circ), \mu\Big),
            \bbW_2\Big(\mu^{\Ising}((\cA_{N_j})^\circ), \cP(\{\pm 1\})\Big)
        \bigg\} \le \epsilon\,,
    \]
    and the correlation functions $\chi_{(\cA_{N_j})^\circ}$ satisfy $\lim_{j\to\infty} \chi_{(\cA_{N_j})^\circ}(0) = 0$ and
    \beq
        \label{e:Lip.in.BOGP.sym.chaos}
        \chi_{\cA_{N_j}^\circ}(1-\epsilon) \le \epsilon\,.
    \eeq
    That is, we have the same information about the algorithms $((\cA_{N_j})^\circ)_{j\ge 1}$ as in the proof of Lemma~\ref{l:Lip.in.BOGP} (instantiated with $q_0 = 0$), and in addition the estimate \eqref{e:Lip.in.BOGP.sym.chaos}.
    For $\cA_{N_j}$ the perturbation of $(\cA_{N_j})^\circ$ given by Proposition~\ref{p:wlogable},  the estimate \eqref{e:Lip.in.BOGP.sym.chaos} implies
    \[
        \chi_{\cA_{N_j}}(1-\epsilon)
        = \bigg(1-\frac{1}{L^2}\bigg) \chi_{\cA_{N_j}^\circ}(1-\epsilon)
        + \frac{1-\epsilon}{L^2} 
        \le \bigg(1-\frac{1}{L^2}\bigg) \epsilon + \frac{1-\epsilon}{L^2} 
        \le 2\epsilon\,,
    \]
where the last inequality holds for $L$ sufficiently large depending on $\alpha,\epsilon$. By passing to a subsequence of the $j$'s, we may assume that $\chi_{\cA_{N_j}}$ converges both pointwise and in $L^1$ to a limiting function $\chi$, which must then satisfy 
    \[
        \chi(1-2\epsilon) \le \chi(1-\epsilon) \le 2\epsilon\,.
    \]
Therefore $p \equiv \chi^{-1}$ satisfies $p(2\epsilon) \ge 1-2\epsilon$, as required by condition~\ref{it:3.6} (with $2\epsilon$ in place of $\epsilon$).
Arguing identically to the proof of Lemma~\ref{l:Lip.in.BOGP} shows that for some $\hat\epsilon(\epsilon) \ge 2\epsilon$ tending to $0$ as $\epsilon \rightarrow 0$, we have 
    \[
        \cM^{\Lip,\chaos}(\alpha,0;L,\epsilon)
        \subseteq \cB_{\hat\epsilon}\Big(
            \cM^{\BOGP,\chaos}(\alpha,0;\tilde C(L),\hat\epsilon)
        \Big)\,.
    \]
    Taking a union in $L$ and intersection in $\epsilon$, and applying $\sym$ to both sides, shows
    \begin{align*}
       & \sym(\sM^{\Lip,\chaos}(\alpha,0))
        \subseteq \sym\bigg(
        \adjustlimits
            \bigcap_{\epsilon > 0} \bigcup_{L>0} \cB_{\epsilon}\Big(
                \cM^{\BOGP,\chaos}(\alpha,0;L,\epsilon)
            \Big)
        \bigg)\\
       &\qquad \subseteq 
        \adjustlimits
        \bigcap_{\epsilon > 0} \bigcup_{L>0}
        \cB_{\epsilon}\Big(
            \sym(\cM^{\BOGP,\chaos}(\alpha,0;L,\epsilon))
        \Big)\\
        &\qquad\stackrel{\eqref{e:BOGP.sym.relation}}{=} \adjustlimits
        \bigcap_{\epsilon > 0} \bigcup_{L>0}
        \cB_{\epsilon}\Big(
        \cM^{\BOGP,\sym}(\alpha,0;L,\epsilon)
        \Big)
        \stackrel{\eqref{e:measure.classes.sym.cM}}{=} \sM^{\BOGP,\sym}(\alpha,0)\,.
    \end{align*}
This concludes the proof.
\end{proof}

\end{lem}
\begin{lem}\label{l:great.in.spr.ideal.in.prf2.sym}
Recalling Definitions \ref{d:SDE.classification.sym}--\ref{d:measure.classes.sym}, we have
\begin{align*}
\sM^{\svi,\sym}(\alpha,0) 
&\subseteq \sM^{\svii,\sym}(\alpha,0)\,,\\
\sM^{\sideal,\sym}(\alpha,0) &\subseteq \sM^{\sx,\sym}(\alpha,0)\,.\end{align*}

\begin{proof}
We repeat the most relevant part of Definition~\ref{d:SDE.classification.sym} for the first assertion: 
\begin{center}
\begin{tabular}[h]{r|cccccccc}
\MSRheader\\
\hline
\GREATdefn{sy}\\
\SPRdefn{sy}
\end{tabular}
\end{center}
The $\svi$ and $\svii$ classes
 differ only in the constraint on $p$: the $\svi$ class satisfies condition~\ref{it:3.6} (i.e. $p(\epsilon) \ge 1-\epsilon$, among other conditions) while the $\svii$ class satisfies condition~\ref{it:3.7} (i.e. $p\equiv 1$). Consider any
	\[\mu = \mu^\sym(q_*=0,b,\sigma,p,\zeta=\zeta_0)
	 \in \cM^{\svi,\sym}(\alpha,0;L,\epsilon)\,,\]
corresponding to controls
	\[(q_*=0,b,\sigma,w,p,\zeta=\delta_0,\zeta^\Ising=\delta_0) \in \Adm^{\svi,\sym}(\alpha,0;L,\epsilon)\,.\]
In particular, $p$ satisfies condition~\ref{it:3.6}.
Recall the notation of Definition~\ref{d:X.adm}: 
we have $(b,\sigma) \in \Adm(0;p,f)$
for budget function $f(t) = (\bbE w_t)^2 / \alpha$, so
    \[
        \acute\mu = \mu(q_*=0,b,\sigma,p,\zeta=\zeta_0) \in \cM(0;p,f)\,.
    \]
Note that $\mu = \sym(\acute\mu)$.
Now consider the functions $\tilde p, \hat p \in \incr([0,1];[0,1])$ given by $\tilde p \equiv 1$ and
    \[
        \hat p(t) = \begin{cases}
            p(\epsilon) & t\in [0,\epsilon) \,, \\
            p(t) & t \in [\epsilon,1]\,.
        \end{cases}
    \]
    Since the restrictions of $p$ and $\hat p$ to $[\epsilon,1]$ agree, Proposition~\ref{p:continuity.in.p.near.0} implies
    \beq\label{e:grt.in.spr.hausdorff1}
        d_\cH\Big(\cM(0;p,f),\cM(0;\hat p,f)\Big) 
        \le \Big(9(1+1/\alpha) \epsilon\Big)^{1/2}\,,
    \eeq
    for $\cM(0;p,f)$ specified in Definition~\ref{d:X.adm}.
    Since $p$ satisfies condition~\ref{it:3.6}, we have $p(\epsilon) \ge 1-\epsilon$, and so
    \[
        \|\hat p' - \tilde p'\|_1
        =\|\hat p' \|_1 
        = \hat p(1) - \hat p(\epsilon)
        = p(1) - p(\epsilon)
        \le \epsilon\,.
    \]
    Then Proposition~\ref{p:continuity.in.p} implies
    \beq\label{e:grt.in.spr.hausdorff2}
        d_\cH\Big(\cM(0;\hat p,f),\cM(0;\tilde p,f)\Big) 
        \le \Big(6(1+1/\alpha) \epsilon\Big)^{1/2}\,.
    \eeq
    For some $\tilde\epsilon = \tilde\epsilon(\epsilon)$ tending to $0$ as $\epsilon \rightarrow 0$, \eqref{e:grt.in.spr.hausdorff1} and \eqref{e:grt.in.spr.hausdorff2} imply the existence of $\tilde \mu \in \cM(0;\tilde p,f)$ such that
    (since $\sym$ is a contraction)
    \beq\label{e:great.in.spr.sym.W2}
        \bbW_2\Big(
        \mu=\sym(\acute{\mu}), 
        \sym(\tilde \mu)
        \Big) 
        \le \bbW_2(\acute \mu,\tilde \mu)
        \le \tilde\epsilon\,.
    \eeq
    Let $\tilde\mu$ be attained by controls $(\tilde b,\tilde \sigma) \in \Adm(0;\tilde p,f)$, i.e.,
    \[
        \tilde\mu = \mu(q_*=0,\tilde b,\tilde \sigma,\tilde p,\zeta=\zeta_0)\,.
    \]
Since $\tilde{p}\equiv1$ satisfies condition~\ref{it:3.7}, we have 
    \[
        (q_*=0,\tilde b,\tilde \sigma,w,\tilde p,\zeta=\delta_0,\zeta^\Ising=\delta_0) 
        \in \Adm^{\svii,\sym}(\alpha,0;\epsilon)\,,
    \]
    and therefore $\sym(\tilde\mu) \in \cM^{\svii,\sym}(\alpha,0;\epsilon)$. 
    Since for all $\mu \in \cM^{\svi,\sym}(\alpha,0;L,\epsilon)$ we can find such $\sym(\tilde\mu)$ satisfying \eqref{e:great.in.spr.sym.W2}, we conclude
    \[
        \cM^{\svi,\sym}(\alpha,0;L,\epsilon) 
        \subseteq 
        \cB_{\tilde\epsilon}\Big( \sM^{\svii,\sym}(\alpha,0;\epsilon) \Big)\,.
    \]
Adjusting $\tilde\epsilon$ gives
    \[
        \cB_\epsilon\Big( \cM^{\svi,\sym}(\alpha,0;L,\epsilon) \Big)
        \subseteq 
        \cB_{\tilde\epsilon}\Big( \sM^{\svii,\sym}(\alpha,0;\tilde\epsilon) \Big)\,.
    \]
Recalling Definition~\ref{d:measure.classes.sym}, 
taking a union over $L$ and an intersection over $\epsilon$ proves the first assertion.

We now repeat the most relevant part of Definition~\ref{d:SDE.classification.sym} for the second assertion: 
\begin{center}
\begin{tabular}[h]{r|cccccccc}
\MSRheader\\
\hline
\SUPERPERFECTdefn{sy}\\
\PRFTWOdefn{sy}
\end{tabular}
\end{center}
The $\sideal$ and $\sx$ classes differ only in that the $\sideal$ class satisfies condition~\ref{it:3.7} (i.e. $p\equiv 1$) while the $\sx$ class satisfies condition~\ref{it:3.1} (i.e. $1/L \le p(0) \le \epsilon$ and $\|p\|_{C^2([0,1])} \le L$). We proceed by a similar argument as above. 
    Consider any 
    \[\mu 
    = \mu^\sym(q_*=0,b,\sigma,p\equiv 1,\zeta=\zeta_0)
    \in \cM^{\sideal,\sym}(\alpha,0)\,,\]
corresponding to controls
	\[(q_*=0,b,\sigma,w,p\equiv 1,\zeta=\delta_0,\zeta^\Ising=\delta_0) \in \Adm^{\sideal,\sym}(\alpha,0)\,.\]
    For any $\epsilon > 0$, consider
    \[
        \tilde p(t) = \begin{cases}
            1 - (2-\epsilon) (\epsilon - t)^2/(2\epsilon^2) & t\in [0,\epsilon)\,, \\ 
            1 & t \in [\epsilon,1]\,.
        \end{cases}
    \]
    Note that $\tilde p$ is twice differentiable, with $\tilde p(0) = \epsilon / 2$ and $\|\tilde p\|_{C^2([0,1])} \le 4/\epsilon^2$.
    This satisfies condition~\ref{it:3.1} if $L \ge 4/\epsilon^2$.
    Since the restrictions of $p$ and $\tilde p$ to $[\epsilon,1]$ agree, Proposition~\ref{p:continuity.in.p.near.0} implies
    \[
        d_\cH\Big(\cM(0;p,f),\cM(0;\tilde p,f)\Big) 
        \le \Big(9(1+1/\alpha) \epsilon\Big)^{1/2}\,.
    \]
    Arguing identically to above then shows that for some $\tilde\epsilon = \tilde\epsilon(\epsilon)$ tending to $0$ as $\epsilon \rightarrow 0$,
    \[
        \cB_\epsilon\Big( \cM^{\sideal,\sym}(\alpha,0) \Big)
        \subseteq 
        \cB_{\tilde\epsilon}\Big( \sM^{\sx,\sym}(\alpha,0;L,\epsilon) \Big)\,,
    \]
    provided $L \ge 4/\epsilon^2$.
Again recalling Definition~\ref{d:measure.classes.sym}, taking a union over $L$ and intersection over $\epsilon$ proves the second assertion of the lemma.
    \end{proof}
\end{lem}

\fi

\pagebreak\section{Confinement of Lipschitz achievable measures to ideal set}
\label{s:confinement}

\def\fm{{\mathfrak{m}}}
\def\fp{{\mathfrak{p}}}
\def\tbG{{\tilde\bG}}
\def\tbg{{\tilde\bg}}
\def\hbG{{\hat\bG}}

\iffull
% !TEX root = main.tex

In Section~\ref{sec:alternate-diffusions}, we showed that in the Ising setting, the sets of measures $\ocM^{\Ising}(\alpha)$ and $\ocM^{\Ising, \concave}(\alpha)$ specified in Definition~\ref{d:achievable-msrs} coincide and describe the asymptotically achievable $\bbW_2$-subsequential limits of $\mu(\cA_N)$ for $O(1)$-Lipschitz algorithms $\cA_N$.
We similarly showed that $\ocM^{\Ising,\sym}(\alpha)$ specified in Definition~\ref{d:achievable-msrs-sym} describes the asymptotically achievable symmetrized measures $\mu_{\sym}(\cA_N)$. The analogous results for the spherical setting also hold, as described in Remark~\ref{rmk:control.problems.spherical}.

However, a sequence $\mu(\cA_N)$ does not necessarily have a $\bbW_2$-convergent subsequence: in the spherical setting, we can simply consider the algorithm $\cA_N(\bG)=\bg^1$. Then $\cA_N$ is a Lipschitz algorithm (it satisfies the conditions of Definition~\ref{d:Lip}, with $L=1$). The measure $\mu(\cA_N)$ is approximately
	\[
	\bigg(1-\frac1N\bigg) \cN(0,1)
	+ \frac{\delta\{N^{1/2}\}}N\,,
	\]
which does not have any subsequence converging in $\bbW_2$. Thus, we cannot conclude that the point-to-set distance $\bbW_2(\mu(\cA_N), \ocM^{\Ising}(\alpha))$ tends to zero as $N\to\infty$.

In this section, we show that weakening the metric resolves this issue: for any $\fq\in [1,2)$, the $\mu(\cA_N)$ are asymptotically confined to $\ocM^{\Ising}(\alpha)$ in the $\bbW_\fq$ sense (and likewise for $\mu_{\sym}(\cA_N)$ and $\ocM^{\Ising,\sym}(\alpha)$).
Formally, recall the definitions of $\cM^{\Lip}(\alpha;L,\epsilon)$ and $\cM^{\Lip,\sym}(\alpha;L,\epsilon)$ from \eqref{e:cM.Lip} and \eqref{e:cM.Lip.sym}, and define
\begin{align}
	\label{e:cM.Lip.union.L}
	\cM^{\Lip}(\alpha;\epsilon) 
	&\equiv \bigcup_{L>0} \cM^{\Lip}(\alpha;L,\epsilon)\,, \\
	\label{e:cM.Lip.union.L.SYM}
	\cM^{\Lip,\sym}(\alpha;\epsilon) 
	&\equiv \bigcup_{L>0} \cM^{\Lip,\sym}(\alpha;L,\epsilon)\,.
\end{align}
As shown in Theorems~\ref{thm:control.problems.main} and \ref{thm:control.problems.sym}, as $\epsilon$ tends to zero the above sets converge respectively to $\sM^{\Lip}(\alpha) = \ocM^{\Ising}(\alpha)$ and $\sM^{\Lip,\sym}(\alpha) = \ocM^{\Ising,\sym}(\alpha)$.

\begin{dfn}\label{d:eps.approx.Lip}
	For $L,\epsilon> 0$, we say $\cA_N$ is an \textbf{$\epsilon$-approximate $L$-Lipschitz algorithm} if it satisfies
	\[\bbW_2(\mu^{\Ising}(\cA_N),\cP(\{\pm 1\})) \le \epsilon\]
in addition to the conditions from
Definition~\ref{d:Lip}.
\end{dfn}

\begin{ppn}
	\label{p:confinement}
	Fix $\alpha, L, \epsilon > 0$ and $\fq\in [1,2)$. 
	For any $\epsilon$-approximate $L$-Lipschitz algorithm $\cA_N$, we have
	\[
		\max\bigg\{
			\bbW_\fq\Big(\mu(\cA_N), \cM^{\Lip}(\alpha;\epsilon)\Big),
			\bbW_\fq\Big(\mu_{\sym}(\cA_N), \cM^{\Lip,\sym}(\alpha;\epsilon)\Big) 
		\bigg\} \le o_N(1)\,,
	\]
	where $o_N(1)$ is a term tending to $0$ as $N\to\infty$ (which may depend on $\alpha,L,\epsilon,\fq$).
\end{ppn}
\begin{rmk}
	Proposition~\ref{p:confinement} is stated for the Ising setting.
	In the corresponding result for the spherical setting, we let $(\cA_N)_{N\ge 1}$ be any sequence of $L$-Lipschitz algorithms, omitting the $\epsilon$-approximate hypothesis. 
We still define $\cM^{\Lip}(\alpha;\epsilon)$ and $\cM^{\Lip,\sym}(\alpha;\epsilon)$ as in \eqref{e:cM.Lip.union.L}, but with the constraint $\bbW_2(\mu^{\Ising}(\cA_{N_j}), \cP(\{\pm 1\})) \le \epsilon$ similarly omitted from the definitions \eqref{e:cM.Lip} and  \eqref{e:cM.Lip.sym} of $\cM^{\Lip}(\alpha;L,\epsilon)$ and $\cM^{\Lip,\sym}(\alpha;L,\epsilon)$, as discussed in Remark~\ref{rmk:control.problems.spherical}. The proof in the spherical setting is only simpler than the one below, using the spherical analogues of Theorems~\ref{thm:control.problems.main} and \ref{thm:control.problems.sym}.
\end{rmk}

The remainder of this section is organized as follows: 
\begin{itemize}
\item In \S\ref{ss:confinement.averaged} we give the proof of Proposition~\ref{p:confinement}.
\item In \S\ref{ss:confinement.main.achievability.results} we give the 
\hyperlink{proof:t.main.it:thm-main-IAMP}{proof of Theorem~\ref{thm:main}\ref{it:thm-main-IAMP}} (achievability in the general setting), as well as the 
\hyperlink{proof:t.symmetric.it:thm-symmetric-IAMP}{proof of Theorem~\ref{thm:symmetric}\ref{it:thm-symmetric-IAMP}} (achievability in the symmetric setting).

\item In \S\ref{ss:confinement.main.hardness.results} we give the 
\hyperlink{proof:t.main.hardness.results}{proofs of Theorem~\ref{thm:main}\ref{it:thm-main-BOGP} and Theorem~\ref{thm:symmetric}\ref{it:thm-symmetric-BOGP}}, our main hardness results for the general and symmetric settings.
\item In \S\ref{ss:confinement.perceptron.threshold.results}
we give the
proof of Corollary~\ref{cor:alg-for-optimization}, our result on the algorithmic threshold of perceptron optimization problems.
\end{itemize}

\subsection{Confinement of averaged inner product distributions}
\label{ss:confinement.averaged}

The main goal of this subsection is the \hyperlink{proof:p.confinement}{proof of Proposition~\ref{p:confinement}}. To this end, we first define a truncated version $\cA_{N[K]}$ of any Lipschitz algorithm $\cA_N$. We will see in Lemma~\ref{l:lip.trunc.effect.on.Wq} below that $\mu(\cA_N)$ and $\mu(\cA_{N[K]})$ are close in $\bbW_\fq$ for any $\fq \in [1,2)$ (though not necessarily in $\bbW_2$). 

\begin{dfn}[row-spliced matrix]
\label{d:row.splice}
	For $S\subseteq [M]$ and matrices $\bG, \hbG \in \R^{M\times N}$ with rows $(\bg^a : a\in [M])$ and $(\hat\bg^a : a\in [M])$, let $\tbG = \tbG(\bG,\hbG;S) \in \R^{M\times N}$ be the matrix with rows $(\tbg^a : a\in [M])$, where 
	\[
		\tbg^a = \begin{cases}
			\bg^a & a \not\in S\,, \\
			\hat\bg^a & a\in S\,.
		\end{cases}
	\]
That is to say, $\tbG$ uses the rows of $\hat{\bG}$ indexed by $S$, and the rows of $\bG$ indexed by $[M]\setminus S$.
\end{dfn}

\begin{dfn}[truncated Lipschitz algorithm]\label{d:truncated.lip.alg}
Consider an $L$-Lipschitz algorithm $\cA_N$ with input $(\bG, \bg^\aux)$. Denote its output $\bx\equiv\cA_N(\bG,\bg^\aux)$. For $a\in [M]$, define the random variable
	\[\fm_a\equiv
		\fm_a(\bG,\bg^\aux;\cA_N)
		\equiv
		\frac{(\bg^a,\bx)}{N^{1/2}}\,,
	\]
where we recall that $\bg^a$ is the $a$-th row of $\bG$.
	For $K > 0$, define
	\[
		S_K(\cA_N) = \bigg\{
			a \in [M] : 
			\Big|\E\fm_a(\bG,\bg^\aux;\cA_N)\Big|
			\ge K
		\bigg\}\,.
	\]
Let $\hbG \in \R^{M\times N}$ be a matrix of i.i.d. gaussians independent of $(\bG, \bg^\aux)$.
	Then define the \textbf{truncated algorithm} $\cA_{N[K]}$ as the function
	\[\cA_{N[K]}\Big(\bG, 
	\bg^{\aux+}
	\equiv(\bg^\aux, \hbG)\Big)
	\equiv
	\cA_N\Big(\tbG(\bG,\hbG;S_K(\cA_N)), \bg^\aux\Big)\,,
	\]
where $\tbG$ is the row-splicing operation from Definition~\ref{d:row.splice}. In words, if $|\E\fm_a| \ge K$, then $\cA_{N[K]}$ ignores the $a$-th row $\bg^a$ of $\bG$ and replaces this input with the independent gaussian vector $\hat\bg$. 	Clearly $\cA_{N[K]}$ is also an $L$-Lipschitz function, with output
	\[\tilde{\bx} \equiv 
		\cA_{N[K]}(\bG, (\bg^\aux, \hbG)) \stackrel{\textit{d}}{=}
		\bx
		= \cA_N(\bG, \bg^\aux)\,.
	\]
It follows that $\cA_{N[K]}$ is also an $L$-Lipschitz algorithm in the sense of Definition~\ref{d:Lip}. 
\end{dfn}

The following basic concentration and moment estimates will be useful. Throughout this section we abbreviate $\|\cdot\|$ for the euclidean norm. 

\begin{lem}\label{l:confinement.subgaussian.basic}
	Suppose $\cA_N$ is a $L$-Lipschitz algorithm.
	Then there exists $C(\alpha) > 0$ depending on only $\alpha$ and $c = c(\alpha,L) > 0$ such that
	\beq\label{e:confinement.subgaussian.basic.tail}
		\P\bigg(
			\frac{\|\bG\|_{\op}}{N^{1/2}}
			 \le C(\alpha)
			\,\,\text{and}\,\,
			\frac{\|\bx\|}{N^{1/2}} \le 2
		\bigg) \ge 1-e^{-cN}\,.
	\eeq
	Furthermore, for any $k\ge 1$, there exists $C(\alpha,k)$ such that for all sufficiently large $N$,
	\beq\label{e:confinement.subgaussian.basic.moment}
		\max\bigg\{
			\frac{\E[(\|\bG\|_{\op})^k]}
				{N^{k/2}}, 
			\frac{\E[\|\bx\|^k]}{N^{k/2}}
		\bigg\} \le C(\alpha,k)\,.
	\eeq
\begin{proof}
	By Lemma~\ref{l:wishart}, there exists $C(\alpha)$ such that 
	\beq
		\label{e:confinement.wishart.restated}
		\P\bigg(
		\frac{\|\bG\|_{\op}}{N^{1/2}}
		 \le C(\alpha)
		\bigg) \ge 1-e^{-cN}\,.
	\eeq
Since the mapping $(\bG,\bg^\aux) \mapsto \|\bx\|$ is $L$-Lipschitz, the random variable $\|\bx\|$ is $L$-subgaussian.
	Combined with the equality $\E[\|\bx\|^2] = N$ from Definition~\ref{d:Lip}, this readily implies
	\[
		\P\bigg(
		\frac{\|\bx\|}
			{N^{1/2}}
		 \le 2\bigg) 
		 \ge 1-e^{-cN}\,,
	\]
which gives \eqref{e:confinement.subgaussian.basic.tail}.
	By Jensen's inequality,
	\[
		\frac{\E\|\bx\|}{N^{1/2}}
		\le 
		\frac{\E[\|\bx\|^2]^{1/2}}
			{N^{1/2}}
		 = 1
		 \,.
	\]
Since $\|\bx\|$ is $L$-subgaussian, for any $u\ge1$ we have
	\[\P\bigg(
			\frac{\|\bx\|}{N^{1/2}} \ge u
		\bigg)
		\le \P\bigg(
			\frac{\|\bx\| - \E\|\bx\|]}{N^{1/2}} \ge u - 1
		\bigg) 
		\le \exp\bigg(
			-\frac{(u-1)^2 N}{2L^2}
		\bigg)\,.\]
Integrating this bound gives
	\begin{align*}
		\E\bigg[\frac{\|\bx\|^k}{N^{k/2}}\bigg]
		&\le 2^k + \int_{2^k}^\infty \P\bigg(
			\frac{\|\bx\|^k}{N^{k/2}} \ge t
		\bigg) \,dt
		= 2^k + \int_2^\infty \P\bigg(
			\frac{\|\bx\|}{N^{1/2}} \ge u
		\bigg) ku^{k-1}\,du \\
		&\le 2^k + \int_2^\infty \exp\bigg(-\frac{(u-1)^2 N}{2L^2}\bigg) ku^{k-1}\,du 
		\le 2^k + 1\,,
	\end{align*}
	where the last inequality holds for all sufficiently large $N$.
	Since the mapping $\bG \mapsto \|\bG\|_{\op}$ is $1$-Lipschitz, the random variable $\|\bG\|_{\op}$ is $1$-subgaussian.
	Then \eqref{e:confinement.wishart.restated} implies that $\E[\|\bG\|_{\op}] \le C(\alpha) N^{1/2}$. Then, by a similar calculation as above, we obtain
	\begin{align*}
		\E\bigg[\frac{\|\bG\|_{\op}^k}{N^{k/2}}\bigg]
		&\le (2C(\alpha))^k + \int_{2C(\alpha)^k}^\infty \P\bigg(
			\frac{\|\bG\|_{\op}^k}{N^{k/2}} \ge t
		\bigg) \,dt \\
		&= (2C(\alpha))^k + \int_{2C(\alpha)}^\infty \P\bigg(
			\frac{\|\bG\|_{\op}}{N^{1/2}} \ge u
		\bigg) ku^{k-1} \,du \\
		&\le (2C(\alpha))^k + \int_{2C(\alpha)}^\infty
		\exp\bigg(
			-\frac{(u-C(\alpha))^2 N}{2}
		\bigg) ku^{k-1}\,du
		\le (2C(\alpha))^k + 1\,,
	\end{align*}
	where the last inequality holds for all sufficiently large $N$.
	This proves \eqref{e:confinement.subgaussian.basic.moment}.
\end{proof}
\end{lem}

\begin{lem}\label{l:lip.alg.l2.l4.bound}
	There exist $C(\alpha) > 0$ and $C(\alpha,L,K) > 0$, such that the following holds for all sufficiently large $N$.
	If $\cA_N$ is an $L$-Lipschitz algorithm and $\cA_{N[K]}$ is as in Definition~\ref{d:truncated.lip.alg}, then 
	\begin{align}
		\label{e:lip.alg.l2.l4.bound.l2}
		\max\{ 
		\|\mu(\cA_N)\|_{L^2}^2
		,\|\mu(\cA_{N[K]})\|_{L^2}^2
		\}
		 &\le C(\alpha)\,, \\
		\label{e:lip.alg.l2.l4.bound.l4}
		\|\mu(\cA_{N[K]})\|_{L^4}^4 &\le C(\alpha,L,K)\,.
	\end{align}
	Furthermore, for $\mu_\bG(\bx)$ defined in \eqref{eq:proj-pursuit-def}, there exists $c = c(\alpha,L) > 0$ such that
	\beq\label{e:lip.alg.l2.l4.bound.tail}
		\P\bigg(\|\mu_\bG(\cA_N(\bG,\bg^\aux))\|_{L^2}^2 \le C(\alpha)
			\bigg) \ge 1 - e^{-cN}\,.
	\eeq

\begin{proof}
	Throughout this proof $c > 0$ is a constant depending on $\alpha,L$ which may change line by line. Recall that we abbreviate $\bx\equiv\cA_N(\bG,\bg^\aux)$. 
	We prove \eqref{e:lip.alg.l2.l4.bound.tail} first. Recalling the definition of $\fm$ from Definition~\ref{d:truncated.lip.alg}, note that
\beq\label{e:lip.alg.l2.l4.bound.step1}
		\|\mu_\bG(\cA_N(\bG,\bg^\aux))\|_{L^2}^2
		=\frac{\|\fm\|^2}{M} 
		= \frac{1}{M} \bigg\|\frac{\bG \bx}{N^{1/2}}\bigg\|^2
		\le \frac{(\|\bG\|_{\op} \|\bx\|)^2}{MN} \,.
	\eeq
	By \eqref{e:confinement.subgaussian.basic.tail} from Lemma~\ref{l:confinement.subgaussian.basic}, with probability $1-e^{-cN}$ the above is
	\[
	\le \frac{C(\alpha)^2 N \cdot 4N}{MN} 
		= \frac{4C(\alpha)^2}{\alpha}\,.
	\]
	This proves \eqref{e:lip.alg.l2.l4.bound.tail} after adjusting $C(\alpha)$.
	Next, note that
	\[
		\|\mu(\cA_N)\|_{L^2}^2
		= \E \|\mu_\bG(\cA_N(\bG,\bg^\aux))\|_{L^2}^2
		= 
		\frac{\E[\|\fm\|^2]}{M} 
		\stackrel{\eqref{e:lip.alg.l2.l4.bound.step1}}{\le} 
		\frac{\E[(\|\bG\|_{\op} \|\bx\|)^2]}{MN} \,.
	\]
By H\"older's inequality and the bound \eqref{e:confinement.subgaussian.basic.moment} from Lemma~\ref{l:confinement.subgaussian.basic}, the above is
	\[\le
		\frac{\E [(\|\bG\|_{\op})^4]^{1/2} \E[\|\bx\|^4]^{1/2}}{MN} 
		\le \frac{C(\alpha,4)N^2}{MN} 
		= \frac{C(\alpha,4)}{\alpha}\,,
	\]
	for $C(\alpha,4)$ defined in Lemma~\ref{l:confinement.subgaussian.basic}.
This proves the first part of \eqref{e:lip.alg.l2.l4.bound.l2} after adjusting $C(\alpha)$. For the second part of \eqref{e:lip.alg.l2.l4.bound.l2}, we split into cases depending on whether $a\in S_K(\cA_N)$. Recall that the truncated algorithm $\cA_{N[K]}$ resamples the rows indexed by $S_K(\cA_N)$. If $a\notin S_K(\cA_N)$, then $\tilde{\fm}_a$ is equidistributed as $\fm_a$, so they have the same second moment. If $a\in S_K(\cA_N)$, then 
	\[\tilde{\fm}_a
		= \frac{(\bg^a, \tilde{\bx})}{N^{1/2}}\,.
	\]
where $\bg^a$ is independent of $\tilde{\bx}$, as $\cA_{N[K]}$ does not use the input $\bg^a$.
	Thus, the distribution of $\tilde{\fm}_a$ conditional on $\tilde{\bx}$ is a (one-dimensional) centered gaussian random variable with standard deviation $\|\tilde{\bx}\| / N^{1/2}$. It follows in this case that
	\[
	\E[(\tilde{\fm}_a)^2]
	= \frac{\E[\|\tilde{\bx}\|^2]}{N}=1\,.
	\]
Altogether this gives
	\[
	\|\mu(\cA_{N[K]})\|^2
	=\frac{ \E[\|\tilde{\fm}\|^2]}{M}
	\le
	\frac{ \E[\|\fm\|^2 + M]}{M}
	\le C(\alpha)+1\,.
	\]
This proves the second part of \eqref{e:lip.alg.l2.l4.bound.l2} after adjusting $C(\alpha)$.

Finally we turn to the proof of  \eqref{e:lip.alg.l2.l4.bound.l4}. Recalling the notations of Definition~\ref{d:truncated.lip.alg}, let
	\[
	\tilde{\fm}_a
	\equiv
	\fm_a\Big(\bG,\bg^{\aux+}
		= (\bg^\aux, \hbG);\cA_{N[K]}
		\Big)
	= \frac{(\bg^a,\tilde{\bx})}{N^{1/2}}\,,
	\]
where we recall that $\tilde{\bx}\equiv \cA_{N[K]}(\bG,\bg^{\aux+}))$. 
We then aim to bound
	\[
		\|\mu(\cA_{N[K]})\|_{L^4}^4
		= \E \Big[\|\mu_\bG(\cA_{N[K]}(\bG,\bg^{\aux+}))\|_{L^4}^4\Big]
		= \frac{1}{M} \sum_{a=1}^M 
		\E\Big[
		(\tilde{\fm}_a)^4\Big]\,.
	\]
To this end, we will show that there exists $C(\alpha,L,K)$ such that
	\beq\label{e:lip.alg.l2.l4.bound.l4.goal}
		\E[(\tilde{\fm}_a)^4] \le C(\alpha,L,K)
	\eeq
	for all $a\in [M]$, which implies \eqref{e:lip.alg.l2.l4.bound.l4}.
	We again split into cases depending on whether $a\in S_K(\cA_N)$:\smallskip

\noindent\textbf{Case 1: $a\not\in S_K(\cA_N)$.} As noted above, in this case we have 
	\[
	\tilde{\fm}_a
		\stackrel{\textit{d}}{=}
		\fm_a
		= \fm_a(\bG,\bg^\aux;\cA_N)
		= \frac{(\bg^a, \bx)}{N^{1/2}}\,,
	\]
where we recall that that $\bx\equiv \cA_N(\bG,\bg^\aux)$.
We next define the truncated random variables 
	\begin{align*}
		\bar\bg^a 
		&\equiv \bg^a \min\bigg(1, \frac{2N^{1/2}}{\|\bg^a\|}\bigg)\,, 
		\\
		\bar{\bx}
		&\equiv\bar\cA_N(\bG,\bg^\aux) 
		\equiv \bx \min \bigg\{
			1, \frac{2N^{1/2}}{\|\bx\|}\bigg\}\,,\\
			\bar{\fm}_a&\equiv
		\bar{\fm}_a(\bG,\bg^\aux;\cA_N)
		= \frac{(\bar\bg^a, \bar{\bx}
		)}{N^{1/2}}\,.
	\end{align*}
	By \eqref{e:confinement.subgaussian.basic.tail} from Lemma~\ref{l:confinement.subgaussian.basic} and a standard bound on $\|\bg^a\|$, the event
	\[
		\cE \equiv \bigg\{
		\frac{\|\bg^a\|}{N^{1/2}}
		\le 2 \textup{ and }
		\frac{\|\bx\|}{N^{1/2}}
			\le 2
			\bigg\}\,.
	\]
occurs with probability at least $1-e^{-cN}$. We will next show that
	\beq\label{e:lip.alg.l2.l4.bound.l4.step2}
		\max_{\fp \in \{1,4\}}
		\bigg|
			\E \Big[(\fm_a)^\fp
			- (\bar\fm_a)^\fp\Big]
		\bigg|
		\le e^{-cN}\,.
	\eeq
	Note that on event $\cE$, we have
	$\fm_a=\bar{\fm}_a$.
	So, for $\fp \in \{1,4\}$, we can bound
	\[\bigg|
			\E \Big[(\fm_a)^\fp
			- (\bar\fm_a)^\fp\Big]
		\bigg|
		= 
		\bigg|
			\E \Big[
			\Big( (\fm_a)^\fp
			- (\bar\fm_a)^\fp\Big)
			; \cE^c\Big]
		\bigg|
		\le \E[
		|\fm_a|^\fp
		; \cE^c]\,,\]
	where the last inequality holds because $\fm_a$ and $\bar{\fm}_a$ always have the same sign, with $|\fm_a| \ge |\bar{\fm}_a|$. Recalling the definition of $\fm_a$, it follows using H\"older's inequality that the above is
	\begin{align*}
		&\le \E\bigg[
			\frac{(\|\bg^a\| \|\bx\|)^\fp}{N^{\fp/2}}; \cE^c
		\bigg] 
		\le \frac{\P(\cE^c)^{1/2} 
		\cdot \E[\|\bg^a\|^{4\fp}]^{1/4} 
		\cdot \E[\|\bx\|^{4\fp}]^{1/4}}{N^{\fp/2}} \,.
	\end{align*}
	The bound \eqref{e:confinement.subgaussian.basic.moment} from Lemma~\ref{l:confinement.subgaussian.basic} ensures that $\E[\|\bx\|^{4\fp}] \le C(\alpha,4\fp) N^{2\fp}$, while a standard gaussian bound shows $\E[\|\bg^a\|^{4\fp}] \le C(\fp) N^{2\fp}$ for some constant $C(\fp)$. Recalling that the event $\cE$ occurs with probability $1-e^{-cN}$, we conclude
	\[
		\bigg|
		\E \Big[
		(\fm_a)^\fp
		- (\bar{\fm}_a)^\fp
		\Big]
		\bigg|
		\le \frac{ e^{-cN/2} \cdot C(\fp)^{1/4} N^{\fp/2} \cdot C(\alpha,4\fp)^{1/4} N^{\fp/2}}{N^{\fp/2}} 
		\,,
	\]
	which is bounded by $e^{-cN}$ after adjusting $c$.
	This proves \eqref{e:lip.alg.l2.l4.bound.l4.step2}.

By the assumption $a\not\in S_K(\cA_N)$, we have $|\E\fm_a| \le K$.
	Together with \eqref{e:lip.alg.l2.l4.bound.l4.step2} (with $\fp=1$) this implies
	\beq\label{e:lip.alg.l2.l4.bound.l4.step3}
		|\E\bar\fm_a| \le K+1
	\eeq
	for sufficiently large $N$.
	Note that the mapping $\bx \mapsto \bx \min\{1,2N^{1/2}/\|\bx\|\}$ is a contraction, so $\bar\cA_N$ is also $L$-Lipschitz. Recall from above the notation
	$\bar{\bx}$, and abbreviate analogously
	$\bar{\bx}'\equiv \bar\cA_N(\bG',(\bg^\aux)')$.  We then bound 
	\begin{align*}
	&\Big|\bar\fm_a(\bG,\bg^\aux;\cA_N) - \bar\fm_a(\bG',(\bg^\aux)';\cA_N)\Big|
	=\bigg|
	\frac{(\bar{\bg}^a,\bar{\bx})}{N^{1/2}}
	-
	\frac{((\bar{\bg}^a)',\bar{\bx}')}{N^{1/2}}
	\bigg|
	 \\
		&\qquad\le \bigg|
			\frac{(\bar\bg^a - (\bar\bg^a)', \bar{\bx})}{N^{1/2}}
		\bigg| + \bigg|
			\frac{((\bar\bg^a)', \bar{\bx} - \bar{\bx}')}{N^{1/2}}
		\bigg| \le \frac{\|\bar\bg^a - (\bar\bg^a)'\| \|\bar{\bx}\|
			+ \|(\bar\bg^a)'\| \|\bar{\bx} - \bar{\bx}'\|}{N^{1/2}}\\ 		&\qquad\le 
			2\Big(
			\|\bar\bg^a - (\bar\bg^a)'\|
			+ \|\bar{\bx} - \bar{\bx}'\|
		\Big)
		\le 2(L+1) \Big\|(\bG,\bg^\aux) - (\bG',(\bg^\aux)')\Big\|_2\,.
	\end{align*}
	Thus the map $(\bG,\bg^\aux) \mapsto \bar\fm_a(\bG,\bg^\aux;\cA_N)$ is $2(L+1)$-Lipschitz.
	By gaussian concentration of measure, this implies that the random variable $\bar{\fm}_a$ is $2(L+1)$-subgaussian.
	Combined with \eqref{e:lip.alg.l2.l4.bound.l4.step3}, we deduce that
	\[
		\E[(\bar{\fm}_a)^4] \le C(\alpha,L,K)
	\]
	for some $C(\alpha,L,K)$.
	Combining with \eqref{e:lip.alg.l2.l4.bound.l4.step2} (with $\fp = 4$) yields \eqref{e:lip.alg.l2.l4.bound.l4.goal} in the case $a\notin S_K(\cA_N)$, after adjusting $C(\alpha,L,K)$.\smallskip

\noindent\textbf{Case 2:} $a\in S_K(\cA_N)$. As noted above, in this case, the distribution of $\tilde{\fm}_a$ conditional on $\tilde{\bx}$ is a (one-dimensional) centered gaussian random variable with standard deviation $\|\tilde{\bx}\| / N^{1/2}$.
	It follows that
	\[
		\E[(\tilde{\fm}_a)^4]
		= \E\bigg[
			\frac{3\|\tilde{\bx}\|^4}{N^2}
		\bigg]
		\le 3C(\alpha,4)\,,\]
	for $C(\alpha,4)$ given by Lemma~\ref{l:confinement.subgaussian.basic}. This proves \eqref{e:lip.alg.l2.l4.bound.l4.goal} in the case $a\in S_K(\cA_N)$. 
In conclusion, having proved 
\eqref{e:lip.alg.l2.l4.bound.l4.goal} in both cases $a\notin S_K(\cA_N)$
and $a\in S_K(\cA_N)$,
the claim \eqref{e:lip.alg.l2.l4.bound.l4} follows as noted above.
\end{proof}
\end{lem}

\begin{lem}\label{l:lip.trunc.effect.on.Wq}
	There exists $C(\alpha) > 0$ such that for all $\fq\in [1,2)$ and any sequence of $L$-Lipschitz algorithms $(\cA_N)_{N\ge 1}$, the following holds.
	For all sufficiently large $N$,
	\[
		\bbW_\fq\Big(\mu(\cA_N), \mu(\cA_{N[K]})\Big)
		\le \frac{C(\alpha)}{K^{2/\fq-1}}\,.
	\]
\begin{proof}
	Consider $N$ large enough that the conclusion of Lemma~\ref{l:lip.alg.l2.l4.bound} holds: this gives
	\[
	C(\alpha)
	\stackrel{\eqref{e:lip.alg.l2.l4.bound.l2} }{\ge}
	\|\mu(\cA_N)\|_{L^2}^2
	= \frac{\E[\|\fm\|^2]}{M}
	\ge \frac{\|\E\fm\|^2}{M}
	\ge \frac{|S_K(\cA_N)| K^2}{M}\,.
	\]	
We can thus couple $x\sim \mu(\cA_N)$ and $y\sim \mu(\cA_{N[K]})$ such that 
$\P(x\neq y) \le C(\alpha) / K^2$. It follows from  the first assertion of  Lemma~\ref{l:lip.alg.l2.l4.bound} that 
	\[
		\E[|x-y|^2]
		\le 2\E[x^2+y^2]
		= 2(\|\mu(\cA_N)\|_{L^2}^2 + \|\mu(\cA_{N[K]})\|_{L^2}^2)
		\stackrel{\eqref{e:lip.alg.l2.l4.bound.l2}}{\le} 4C(\alpha)\,.\]
Applying H\"older's inequality then gives
	\begin{align*}
		&\bbW_\fq\Big(\mu(\cA_N), \mu(\cA_{N[K]})\Big)^\fq
		\le \E[|x-y|^\fq]
		\le \P(x\neq y)^{1-\fq/2} \E[|x-y|^2]^{\fq/2}
	\\
	&\qquad
		\le \bigg(\frac{C(\alpha)}{K^2}\bigg)^{1-\fq/2} (4C(\alpha))^{\fq/2}
		= \frac{2^\fq C(\alpha)}{K^{2-q}}
		\le \bigg(\frac{2\max\{1,C(\alpha)\}}{K^{2/\fq-1}}\bigg)^\fq\,.
	\end{align*}
	The result follows by renaming $2\max\{1,C(\alpha)\}$ to $C(\alpha)$.
\end{proof}
\end{lem}

For $\fp \ge 1$ and $C>0$, let $\cP_\fp(\R,C)$ and $\cP_\fp(\R_{\ge0},C)$ be the spaces of Borel probability measures $\mu$ on $\R$ and $\R_{\ge0}$ satisfying the bound $\|\mu\|_{L^\fp}^\fp \le C$.

\begin{lem}\label{l:confinement.cPp.compact.Wq}
	For all $1\le \fq < \fp$ and $C > 0$, the spaces $\cP_\fp(\R,C)$ and $\cP_\fp(\R_{\ge0},C)$, endowed with the $\bbW_\fq$ metric,
are compact.

\begin{proof}
	We will show compactness of $\cP_\fp(\R,C)$, as compactness of $\cP_\fp(\R_{\ge0},C)$ is proved analogously.
	It suffices to show that every sequence $(\mu_n)_{n\ge 1}$ in $\cP_\fp(\R,C)$ has a subsequence converging in $\bbW_\fq$ to a measure in $\cP_\fp(\R,C)$.
	Consider any such sequence $(\mu_n)_{n\ge 1}$.
	The set of measures $\cP_\fp(\R,C)$ is clearly tight, so by Prohorov's theorem there exists a weakly convergent subsequence $(\mu_{n_j})_{j\ge 1}$. 
Let $\mu$ be its weak limit. By Skorokhod's theorem, we can find a coupling of random variables $X_j\sim\mu_{n_j}$ and $X\sim\mu$ such that
$X_j\to X$ almost surely. Then, by Fatou's lemma, we have
	\[
		\|\mu\|_{L^\fp}^\fp 
	=\E[|X|^p]^{1/p}
		\le
		\liminf_j \E[|X_j|^p]^{1/p}
		= \liminf_{j\to\infty} \|\mu_{n_j}\|_{L^\fp}^\fp
		 \le C\,,
	\]
	so $\mu \in \cP_\fp(\R,C)$. 
	It remains to show $(\mu_{n_j})_{j\ge 1}$ converges to $\mu$ in the $\bbW_\fq$ metric.
	 Note that for any $\nu \in \cP_\fp(\R_{\ge0},C)$, Markov's inequality implies
	\[
	\nu(|x|>K) 
	\le \frac{1}{K^\fp}\int
		|x|^\fp \,d\nu
		 \le \frac{C}{K^\fp}\,.
	\]
Combining this with H\"older's inequality gives
	\beq\label{e:confinement.cPp.compact.Wq.tail}
		\int_{|x|>K} |x|^\fq \,d\nu
		\le \nu(|x|>K)^{1 - \fq/\fp} \bigg(
			\int |x|^\fp \,d\nu
		\bigg)^{\fq/\fp}
		\le \frac{C}{K^{\fp - \fq}}\,.
	\eeq
	By \eqref{e:confinement.cPp.compact.Wq.tail} and the weak convergence $\mu_{n_j} \rightarrow \mu$,
	\begin{align*}
		\limsup_{j\to\infty} \bigg|
			\int |x|^\fq \,d\mu_{n_j}
			- \int |x|^\fq \,d\mu
		\bigg|
		&\le \limsup_{j\to\infty} \bigg|
			\int
			\max\{|x|,K\}^\fq
			\,d\mu_{n_j}
			- \int
			\max\{|x|,K\}^\fq
			 \,d\mu
		\bigg|\\
		&\qquad + \limsup_{j\to\infty} \int_{|x| > K} |x|^\fq \,d\mu_{n_j}
		+ \int_{|x| > K} |x|^\fq \,d\mu
		\le \frac{2C}{K^{\fp-\fq}}\,.
	\end{align*}
	Since $K$ was arbitrary, we infer the convergence of $\fq$-th moments
	\[
		\lim_{j\to\infty} 
		\E[|X_j|^p]
		=\lim_{j\to\infty} 
		\int |x|^\fq \,d\mu_{n_j}
		= \int |x|^\fq \,d\mu
		= \E[|X|^p]\,.\]
This implies that the random variables $\{X_j\}$ are uniformly integrable, so
	\begin{align*}
	&\bbW_\fq(\mu_{n_j}, \mu)^\fq
	\le \E[|X_j-X|^\fq]\\
	&\le\E\Big[
	\max\{|X_j-X|,2K\}^\fq\Big]
	+\E\Big[|X_j| ; |X_j|\ge K\Big]
	+ \E\Big[|X| ; |X|\ge K\Big]
	\\
	&\le o_j(1) + o_K(1)+ o_K(1)\,.
	\end{align*}
Since the left-hand side does not depend on $K$, we conclude that
$\mu_{n_j} \rightarrow \mu$ in the $\bbW_\fq$ metric.
\end{proof}
\end{lem}

\begin{lem}\label{l:confinement.lip.trunc.to.ideal}
	For any $\alpha,L,\epsilon,K>0$ and any sequence of $\epsilon$-approximate $L$-Lipschitz algorithms $(\cA_N)_{N\ge 1}$,
	\[
		\lim_{N\to\infty} \bbW_2\Big(\mu(\cA_{N[K]}), \cM^{\Lip}(\alpha;\epsilon)\Big) 
		= \lim_{N\to\infty} \bbW_2\Big(\mu_{\sym}(\cA_{N[K]}), \cM^{\Lip,\sym}(\alpha;\epsilon)\Big)
		= 0\,.
	\]

\begin{proof}
	We only prove that the first limit is zero, as the second limit is proved analogously.
	Suppose for contradiction that there exists $\delta > 0$ such that $\bbW_2(\mu(\cA_{N[K]}), \cM^{\Lip}(\alpha;\epsilon)) \ge \delta$ for infinitely many $N$.
	By \eqref{e:lip.alg.l2.l4.bound.l4} from Lemma~\ref{l:lip.alg.l2.l4.bound}, there exists $C = C(\alpha,L,K) > 0 $ such that for sufficiently large $N$, $\mu(\cA_{N[K]}) \in \cP_4(\R,C)$.
Then, for infinitely many $N$, we have
	\[\mu(\cA_{N[K]})\in
		S_\delta 
		\equiv \Big\{ \mu \in \cP_4(\R,C) : \bbW_2(\mu, \cM^{\Lip}(\alpha;\epsilon)) \ge \delta 
		\Big\}\,.
	\]
	By Lemma~\ref{l:confinement.cPp.compact.Wq}, the space $\cP_4(\R,C)$ endowed with the $\bbW_2$ metric is compact.
	The set $S_\delta$ is closed with respect to the $\bbW_2$ metric, hence also compact. Therefore
 we may find $\mu \in S_\delta$ and a subsequence $(N_j)_{j\ge 1}$ such that
	\[
		\lim_{j\to\infty} \bbW_2\Big(\mu(\cA_{N_j[K]}),
		\mu\Big) = 0\,.
	\]
	As noted in Definition~\ref{d:truncated.lip.alg}, we have $\mu^{\Ising}(\cA_{N_j[K]}) = \mu^{\Ising}(\cA_{N_j})$: thus, since each $\cA_{N_j}$ is an $\epsilon$-approximate $L$-Lipschitz algorithm, so is each $\cA_{N_j[K]}$. It follows by omitting finitely many $j$'s that we have
	\[
	\WERR(\cA_{N_j[K]},\mu)
	=\max\bigg\{
	\bbW_2\Big(\mu(\cA_{N_j[K]}),
		\mu\Big),
	\bbW_2\Big(
		\mu^\Ising(\cA_{N_j[K]}),
		\cP(\{\pm1\})
		\Big)
	\bigg\}\le\epsilon
	\]
for all $j$. Recalling the definition \eqref{e:cM.Lip} of  $\cM^{\Lip}(\alpha;L,\epsilon)$, this implies
	\[
		\mu \in \cM^{\Lip}(\alpha;L,\epsilon) \subseteq \cM^{\Lip}(\alpha;\epsilon)\,,
	\]
	which contradicts that $\mu \in S_\delta$.
\end{proof}
\end{lem}

\begin{proof}[\hypertarget{proof:p.confinement}{Proof of Proposition~\ref{p:confinement}}]

	We only prove that the first limit is zero, as the second limit is proved analogously.
	Fix any $K > 0$.
	Combining Lemmas~\ref{l:lip.trunc.effect.on.Wq} and \ref{l:confinement.lip.trunc.to.ideal} gives
	\begin{align*}
		&\limsup_{N\to\infty} \bbW_\fq\Big(\mu(\cA_N), \cM^{\Lip}(\alpha;\epsilon)\Big)\\ 
		&\qquad\le \limsup_{N\to\infty} \bbW_\fq\Big(\mu(\cA_N), \mu(\cA_{N[K]})\Big) 
		+ \limsup_{N\to\infty} \bbW_\fq\Big(\mu(\cA_{N[K]}), \cM^{\Lip}(\alpha;\epsilon)\Big) \\
		&\qquad\le \frac{C(\alpha)}{K^{2/\fq-1}}
		+ \limsup_{N\to\infty} \bbW_2(\mu(\cA_{N[K]}), \cM^{\Lip}(\alpha;\epsilon))
		= \frac{C(\alpha)}{K^{2/\fq-1}}\,,
	\end{align*}
	where $C(\alpha)$ is given by Lemma~\ref{l:lip.trunc.effect.on.Wq}.
	Since the left-hand side does not depend on $K$, we conclude that it in fact equals zero, and the claim follows.
\end{proof}

\subsection{Proofs of main achievability results}
\label{ss:confinement.main.achievability.results} In this subsection we present the 
\hyperlink{proof:t.main.it:thm-main-IAMP}{proof of Theorem~\ref{thm:main}\ref{it:thm-main-IAMP}} (achievability in the general setting), as well as the 
\hyperlink{proof:t.symmetric.it:thm-symmetric-IAMP}{proof of Theorem~\ref{thm:symmetric}\ref{it:thm-symmetric-IAMP}} (achievability in the symmetric setting). These
results essentially follow from results obtained in Sections~\ref{sec:IAMP}--\ref{sec:alternate-diffusions}.
We present these proofs below, after the following concentration estimate (which will also be used in \S\ref{ss:confinement.main.hardness.results} to prove our hardness results).

\begin{lem}\label{l:Wq.lip.to.fixed.msr.conc}
	For any $\alpha,L,\delta>0$ there exists $c = c(\alpha,L,\delta) > 0$ such that the following holds for all $\fq\in [1,2]$ and all sufficiently large $N$.
Recall $\mu_\bG(\bx)$ defined in \eqref{eq:proj-pursuit-def}. For any fixed $\mu \in \cP_2(\R)$ and $L$-Lipschitz algorithm $\cA_N$, abbreviate $\bx\equiv\bx_N\equiv \cA_N(\bG,\bg^\aux)$ as before, and define the random variables
	\begin{align*}
	X&\equiv X_N
	\equiv X(\bG,\bg^\aux)
	\equiv \bbW_\fq\Big(\mu_\bG(\bx_N),\mu\Big)\,,\\
	J&\equiv J_N
	\equiv J(\bG,\bg^\aux)
	\equiv 
	\bbW_\fq\Big(
	\EmpDist(\bx_N),\cP(\{\pm1\})
	\Big)\,.
	\end{align*}
 We then have the bounds
\begin{align}\label{e:Wq.lip.to.fixed.msr.conc.tail}
		\bbP\Big(
			|X_N - \E X_N|
			\le \delta
		\Big)
		& \ge 1-e^{-cN}\,,\\	\label{e:Wq.lip.to.fixed.msr.conc.var}
		\bbE\Big[
			(X_N - \E X_N)^2
		\Big] 
		&= o_N(1)\,,\\
		\label{e:Wq.lip.to.fixed.msr.conc.tail.ising}
		\bbP\Big(
			|J_N - \E J_N|
			\le \delta
		\Big)
		& \ge 1-e^{-cN}\\
\label{e:Wq.lip.to.fixed.msr.conc.var.ising}
		\bbE\Big[
			( J_N - \E J_N)^2
		\Big] 
		&= o_N(1)\,.
	\end{align}
where $o_N(1)$ denotes a term tending to $0$ as $N\to\infty$, which may depend on $\alpha,L,\mu$.
	Analogous inequalities to \eqref{e:Wq.lip.to.fixed.msr.conc.tail} and  \eqref{e:Wq.lip.to.fixed.msr.conc.var} hold with $\mu_{\bG,\sym}(\bx_N)$ in place of $\mu_\bG(\bx_N)$, and $\mu \in \cP_2(\R_{\ge 0})$.
\begin{proof}
	Throughout this proof, $c$ denotes a constant that depends on $\alpha,L,\delta$, which can change line by line.
	Let $C = C(\alpha)$ be given by Lemma~\ref{l:confinement.subgaussian.basic}. As in the proof of Lemma~\ref{l:lip.alg.l2.l4.bound}, let
	\[\bar{\bx}
	\equiv\bar\cA_N(\bG,\bg^\aux) 
		\equiv \bx \min \bigg\{
			1, \frac{2N^{1/2}}{\|\bx\|}\bigg\}\,.\] 
However, in contrast with the proof of Lemma~\ref{l:lip.alg.l2.l4.bound}, we now let $\bar{\bg}^a$ denote the $a$-th row of 
	\[\bar\bG = \bG \min\bigg(1,\frac{CN^{1/2}}{\|\bG\|_{\op}}\bigg)\,.\]
 We then define the measure
	\[
		\bar\mu(\bG,\bg^\aux) = \EmpDist\bigg(\frac{\bar\bG \bar{\bx}}{N^{1/2}}\bigg)
		\equiv \frac{1}{M} \sum_{a=1}^M \delta\bigg\{
		\frac{(\bar\bg^a, \bar{\bx})}{N^{1/2}}\bigg\}\,,
	\]
	By  the bound \eqref{e:confinement.subgaussian.basic.tail} from Lemma~\ref{l:confinement.subgaussian.basic}, the event
	\[
		\cE = \bigg\{
			\frac{\|\bG\|_{\op}}{N^{1/2}} \le C \textup{ and }
			\frac{\|\bx\|}{N^{1/2}} \le 2
		\Big\}\,.
	\]
occurs with probability at least $1-e^{-cN}$. We abbreviate
	\[Y \equiv Y_N
	\equiv Y(\bG,\bg^\aux) 
	\equiv \bbW_\fq\Big(\bar\mu(\bG,\bg^\aux),\mu\Big)\]
We will first show that
	\beq\label{e:Wq.lip.to.fixed.msr.conc.new2}
		\max_{\fp \in \{1,2\}}
		\Big|\E[X^\fp] - \E[Y^\fp]
		\Big| \le e^{-cN}\,.
	\eeq
		Note the trivial bound $\bbW_2(\nu,\mu) \le \|\nu\|_{L^2} + \|\mu\|_{L^2}$. Thus, for $\fp \in \{1,2\}$,
	\[
		\bbW_2(\nu,\mu)^\fp \le 2\Big\{
		(\|\nu\|_{L^2})^\fp + 
		(\|\mu\|_{L^2})^\fp
		\Big\}\,.
	\]
	Note also that on event $\cE$, $X=Y$. Thus, for each $\fp \in \{1,2\}$, we have
	\begin{align*}
		&|\E[X^\fp] - \E[Y^\fp]|
		\le \E\Big[
		(X^\fp + Y^\fp);\cE^c\Big] 
		\le \E\bigg[
			\Big\{
			\bbW_2\Big(\mu_\bG(\bx),\mu\Big)^\fp + \bbW_2\Big(\bar\mu(\bG,\bg^\aux),\mu\Big)^\fp
		\Big\}
		;\cE^c\bigg]\\
		&\qquad\le 2\E\bigg[
		\Big\{
		(\|\mu_\bG(\bx)\|_{L^2})^\fp 
		+ (\|\bar\mu(\bG,\bg^\aux)\|_{L^2})^\fp 
		+ 2(\|\mu\|_{L^2})^\fp
		\Big\}
		;\cE^c
		\bigg] \\
		&\qquad\stackrel{*}{\le} 4\E\bigg[ \Big\{
			(\|\mu_\bG(\bx)\|_{L^2})^\fp 
			+ (\|\mu\|_{L^2})^\fp\Big\}
		;\cE^c
		\bigg] 
		= 4\E\bigg[
		\bigg(
				\frac{\|\bG \bx\|^\fp}{(MN)^{\fp/2}}
				+ (\|\mu\|_{L^2})^\fp
			\bigg)
			;\cE^c
		\bigg] \\
		&\qquad\le 4\E\bigg[
			\ind\{\cE^c\}\bigg(
				\frac{(\|\bG\|_{\op})^\fp \|\bx\|^\fp}{(MN)^{\fp/2}}
				+ (\|\mu\|_{L^2})^\fp
			\bigg)
		\bigg]\,.
	\end{align*}
	The  inequality marked $*$ uses that $\|\mu_\bG(\cA_N(\bG,\bg^\aux))\|_{L^2}$ stochastically dominates $\|\bar\mu(\bG,\bg^\aux)\|_{L^2}$. Using H\"older's inequality, the above is
	\[\le
		\frac{4\P(\cE^c)^{1/2}
		\E[\|(\bG\|_{\op})^{4\fp}]^{1/4}
		\E[\|\bx\|^{4\fp}]^{1/4}}{(MN)^{\fp/2}}
		+ 4\P(\cE^c) 
		(\|\mu\|_{L^2})^\fp\,.
	\]
Recalling Lemma~\ref{l:confinement.subgaussian.basic}, from \eqref{e:confinement.subgaussian.basic.tail} we have $\P(\cE^c) \le e^{-cN}$,
while the moments in the above expression can be bounded by 
\eqref{e:confinement.subgaussian.basic.moment}. 
	This proves \eqref{e:Wq.lip.to.fixed.msr.conc.new2} after adjusting $c$.

	Note that $\bG \mapsto \bar\bG$ is a contraction, and is therefore $1$-Lipschitz.
	Similarly, $\bx \mapsto \bx \min\{1,2N^{1/2}/\|\bx\|_2\}$ is a contraction on $\R^N$, so $\bar\cA_N$ is an $L$-Lipschitz function.
	Recall from above the notation
	$\bar{\bx}$, and abbreviate analogously
	$\bar{\bx}'\equiv \bar\cA_N(\bG',(\bg^\aux)')$. We then have
	\begin{align*}
		&|Y(\bG,\bg^\aux) - Y(\bG',(\bg^\aux)')|
		\le \bbW_\fq(\bar\mu(\bG,\bg^\aux),\bar\mu(\bG',(\bg^\aux)')) \\
		&\qquad \le \bbW_2(\bar\mu(\bG,\bg^\aux),\bar\mu(\bG',(\bg^\aux)'))
		 \le \frac{\|\bar\bG\|_{\op} \|\bar{\bx} - \bar{\bx}'\|
			+ \|\bar\bG - \bar\bG'\|_{\op} \|\bar{\bx}'\|}{(MN)^{1/2}}\\ 		&\qquad \le 
			\frac{CN^{1/2} \cdot L \|(\bG,\bg^\aux) - (\bG',(\bg^\aux)')\|_2
			+ \|\bG - \bG'\|_2 \cdot 2N^{1/2}}{(MN)^{1/2}} 
			 \\
		&\qquad = \frac{CL+2}{M^{1/2}} \|(\bG,\bg^\aux) - (\bG',(\bg^\aux)')\|_2\,.
	\end{align*}
	Therefore, the mapping $(\bG,\bg^\aux) \mapsto Y(\bG,\bg^\aux)$ is $(CL+2) / M^{1/2}$-Lipschitz.
	By gaussian concentration of measure, the random variable $Y$ is $(CL+2) / M^{1/2}$-subgaussian, so we obtain
	\[
		\bbP\bigg( |Y - \E[Y]|  \le \frac{\delta}{2}\bigg) \ge 1- e^{-cN}\,.
	\]
	Combined with \eqref{e:Wq.lip.to.fixed.msr.conc.new2} (with $\fp = 1$) we deduce
	\[
		\bbP\Big(|Y - \E[X]| \le \delta \Big) \ge 1- e^{-cN}\,.
	\]
Since $X=Y$ on the event $\cE$, we conclude
	\[
	\bbP\Big(|X - \E[X]| > \delta \Big)
	\le \bbP(\cE^c)
	+ \bbP\Big(|Y - \E[X]| > \delta 
		\Big)
	\le e^{-cN}\,,
	\]
by adjusting $c$ appropriately. This
	 proves \eqref{e:Wq.lip.to.fixed.msr.conc.tail}.

	We next turn to the proof of  \eqref{e:Wq.lip.to.fixed.msr.conc.var}. 
	Since $Y$ is $(CL+2) / M^{1/2}$-subgaussian, we also have
	\beq\label{e:Wq.lip.to.fixed.msr.conc.new4}
		\E\Big[(Y - \E Y)^2
		\Big] = o_N(1)\,.
	\eeq
	We then bound
	\begin{align}
		\nonumber
		&\Big|
			\E\Big[ (X - \E X)^2\Big]
			- \E\Big[(Y - \E Y)^2\Big]
		\Big|
		\le \Big| \E(X^2) - \E(Y^2) 
			\Big|
		+ | \E X - \E Y | 
		( \E X + \E Y ) \\
		\label{e:Wq.lip.to.fixed.msr.conc.new5}
		&\qquad\stackrel{\eqref{e:Wq.lip.to.fixed.msr.conc.new2}}{\le} 
		e^{-cN} 
		\Big(1 + \E X + \E Y \Big)\,.
	\end{align}
	These remaining expectations are bounded by
	\begin{align*}
		\E X + \E Y
		&\le \E[\bbW_2(\mu_\bG(\cA(\bG,\bg^\aux)),\mu) + \bbW_2(\bar\mu(\bG,\bg^\aux),\mu)] \\
		&\le \E[\|\mu_\bG(\cA(\bG,\bg^\aux))\|_{L^2} + \|\bar\mu(\bG,\bg^\aux)\|_{L^2} + 2\|\mu\|_{L^2}] \\
		&\le 2\E[\|\mu_\bG(\cA(\bG,\bg^\aux))\|_{L^2}] + 2\E[\|\mu\|_{L^2}] \\
		&\le 2\E[(\|\mu_\bG(\cA(\bG,\bg^\aux))\|_{L^2})^2]^{1/2} + 2\E[\|\mu\|_{L^2}]
		=  2\|\mu(\cA_N)\|_{L^2} + 2\E[\|\mu\|_{L^2}]\,.
	\end{align*}
	By the estimate \eqref{e:lip.alg.l2.l4.bound.l2} from Lemma~\ref{l:lip.alg.l2.l4.bound}, the above is bounded independently of $N$.
	Combining with \eqref{e:Wq.lip.to.fixed.msr.conc.new4} and \eqref{e:Wq.lip.to.fixed.msr.conc.new5} gives 
	\[
		\E\Big[
		(X - \E X)^2\Big] = o_N(1)\,,
	\]
	which proves \eqref{e:Wq.lip.to.fixed.msr.conc.var}.
	The analogues of \eqref{e:Wq.lip.to.fixed.msr.conc.tail} and \eqref{e:Wq.lip.to.fixed.msr.conc.var} with $\mu_{\bG,\sym}(\bx)$ in place of $\mu_\bG(\bx)$ and $\mu \in \cP_2(\R_{\ge 0})$ are proved identically.

	The remaining estimates \eqref{e:Wq.lip.to.fixed.msr.conc.tail.ising} and \eqref{e:Wq.lip.to.fixed.msr.conc.var.ising} are similar but simpler.
	Denote $\bx'\equiv\cA_N(\bG',(\bg^\aux)')$. Then
	\begin{align*}
		&\bigg|\bbW_2\Big(\EmpDist(\bx),\cP(\{\pm 1\})\Big) 
		- \bbW_2\Big(\EmpDist(\bx'),\cP(\{\pm 1\})\Big)\bigg| \\
		&\le \bbW_2\Big(\EmpDist(\bx),\EmpDist(\bx')\Big)
		\le \frac{\|\bx - \bx'\|}{N^{1/2}} \le \frac{L}{N^{1/2}}\,,
	\end{align*}
since $\cA_N$ is $L$-Lipschitz. This implies that $(\bg,\bg^\aux) \mapsto \bbW_2(\EmpDist(\bx),\cP(\{\pm 1\}))$ is $L/N^{1/2}$-Lipschitz.
The estimates  \eqref{e:Wq.lip.to.fixed.msr.conc.tail.ising} and \eqref{e:Wq.lip.to.fixed.msr.conc.var.ising} then follow from gaussian concentration of measure.
\end{proof}
\end{lem}

\begin{proof}[\hypertarget{proof:t.main.it:thm-main-IAMP}{Proof of Theorem~\ref{thm:main}\ref{it:thm-main-IAMP}}]
	Consider any $\iota > 0$ and any $\mu \in \ocM^{\Ising}(\alpha)$.
	The \hyperlink{p:thm.control.problems.main}{proof of Theorem~\ref{thm:control.problems.main}} shows that 
$\ocM^{\Ising}(\alpha) = \sM^{\IAMP}(\alpha,0)$, where we recall from
	Definition~\ref{d:measure.classes} and Remark~\ref{r:measure.classes.simplification} that
	\[\sM^{\IAMP}(\alpha,0)
	\stackrel{\eqref{e:measure.classes.simplification.without.L0}}{=}
	\adjustlimits
	\bigcap_{\epsilon>0}
	\bigcup_{L>0}\cB_\epsilon\Big(
	\cM^{\IAMP}(\alpha,0;L,\epsilon)
	\Big)\,.
	\]
We set $\epsilon$ as the value provided by Theorem~\ref{thm:IAMP-main}\ref{i:IAMP-main-main}
(which depends only on $\alpha,\iota$).
	As the conclusion of Theorem~\ref{thm:IAMP-main}\ref{i:IAMP-main-main} is monotone in $\epsilon$, we may further assume $\epsilon \le \iota$.
	Then, there exists $L>0$ such that we have $\mu \in \cB_\epsilon(\cM^{\IAMP}(\alpha,0;L,\epsilon))$.
This means that  we can find $\tilde\mu \in \cM^{\IAMP}(\alpha,0;L,\epsilon)$ such that $\bbW_2(\mu,\tilde\mu) \le \epsilon$.
	Then Theorem~\ref{thm:IAMP-main}\ref{i:IAMP-main-main} provides a $C(L,\epsilon)$-Lipschitz algorithm $\cA_N$, which agrees with the output of an efficiently implementable IAMP algorithm with probability $1-e^{-cN}$, such that 
	\begin{align*}
		\P\Big(\cA_N(\bG,\bg^\aux) \in \Sigma(\iota)\Big) 
		&\ge 1-e^{-cN}\,,\\
		\E\bbW_2\Big(\mu_\bG(\cA_N(\bG,\bg^\aux)),\tilde\mu\Big)
		&\le \iota\,.
	\end{align*}
Combining the estimate \eqref{e:Wq.lip.to.fixed.msr.conc.tail} from Lemma~\ref{l:Wq.lip.to.fixed.msr.conc} yields that, with probability $1-e^{-cN}$,
	\[
		\bbW_2\Big(\mu_\bG(\cA_N(\bG,\bg^\aux)),\tilde\mu\Big) \le 2\iota\,.
	\]
On this event it follows that
	\[
		\bbW_2\Big(\mu_\bG(\cA_N(\bG,\bg^\aux)),\mu\Big) 
		\le \bbW_2\Big(\mu_\bG(\cA_N(\bG,\bg^\aux)),\tilde\mu\Big)
		+ \bbW_2\Big(\mu,\tilde\mu\Big)
		\le 2\iota + \epsilon \le 3\iota\,.
	\]
	We conclude that $\cA_N$ $(3\iota,1-e^{-cN})$-attains $\mu$.
	Since $\iota > 0$ was arbitrary, the result follows.
\end{proof}

\begin{proof}[\hypertarget{proof:t.symmetric.it:thm-symmetric-IAMP}{Proof of Theorem~\ref{thm:symmetric}\ref{it:thm-symmetric-IAMP}}]
	This proof proceeds similarly to above.
	Suppose $\iota > 0$ and $\mu \in \ocM^{\Ising,\sym}(\alpha)$.
	The \hyperlink{p:thm.control.problems.sym}{proof of Theorem~\ref{thm:control.problems.sym}} shows $\ocM^{\Ising,\sym}(\alpha) = \sM^{\IAMP,\sym}(\alpha,0)$, where this set is specified in Definition~\ref{d:measure.classes.sym}.
	So $\mu \in \sM^{\IAMP,\sym}(\alpha,0)$.
	Let $\epsilon$ be given by Theorem~\ref{thm:IAMP-main}\ref{i:IAMP-main-main}
	(depending on $\alpha,\iota$); as above we may assume $\epsilon \le \iota$.
	Arguing identically as above, we can find $L>0$ and $\tilde\mu \in \cM^{\IAMP,\sym}(\alpha,0;L,\epsilon)$ such that $\bbW_2(\mu,\tilde\mu) \le \epsilon$.
	Recall the set $\Adm^{\IAMP,\sym}(\alpha,0;L,\epsilon)$ from Definition~\ref{d:SDE.classification.sym}.
	Then we can find controls
	\[
		(q_*,b,\sigma,w,p,\zeta,\zeta^\Ising) \in \Adm^{\IAMP,\sym}(\alpha,0;L,\epsilon)
	\]
	such that 
	\[
		\tilde\mu = \mu_{\sym}(q_*,b,\sigma,p,\zeta)\,.
	\]
	Also recall the set $\Adm^{\IAMP}(\alpha,0;L,\epsilon)$ from Definition~\ref{d:SDE.classification}.
	By comparing Definitions~\ref{d:SDE.classification} and \ref{d:SDE.classification.sym}, we see that in fact 
		\[\Adm^{\IAMP}(\alpha,0;L,\epsilon) = \Adm^{\IAMP,\sym}(\alpha,0;L,\epsilon)\,.\]
Thus we also have
	\[
		(q_*,b,\sigma,w,p,\zeta,\zeta^\Ising) \in \Adm^{\IAMP}(\alpha,0;L,\epsilon)\,.
	\]
	Let $\hat\mu = \mu(q_*,b,\sigma,p,\zeta)$, so that $\sym(\hat\mu) = \tilde\mu$.
	Then, arguing as above,  Theorem~\ref{thm:IAMP-main}\ref{i:IAMP-main-main} 
	provides a $C(L,\epsilon)$-Lipschitz algorithm $\cA_N$ such that with probability $1-e^{-cN}$, $\cA_N(\bG,\bg^\aux) \in \Sigma(\iota)$ and
	\[\bbW_2\Big(\mu_{\bG,\sym}(\cA_N(\bG,\bg^\aux)),\tilde\mu\Big)
	\le
		\bbW_2\Big(\mu_\bG(\cA_N(\bG,\bg^\aux)),\hat\mu\Big) \le 2\iota
	\]
(having used that $\sym$ is a contraction).
This implies
	\[
		\bbW_2\Big(\mu_{\bG,\sym}(\cA_N(\bG,\bg^\aux)),\mu\Big) 
		\le \bbW_2\Big(\mu_{\bG,\sym}(\cA_N(\bG,\bg^\aux)),\tilde\mu\Big) 
		+ \bbW_2(\mu,\tilde\mu)
		\le 2\iota + \epsilon \le 3\iota\,.
	\]
	So $\cA_N$ $(3\iota,1-e^{-cN})$-symmetrically attains $\mu$, and the result follows.
\end{proof}

\subsection{Proofs of main hardness results}
\label{ss:confinement.main.hardness.results}

In this subsection we present the
\hyperlink{proof:t.main.hardness.results}{proofs of Theorem~\ref{thm:main}\ref{it:thm-main-BOGP} and Theorem~\ref{thm:symmetric}\ref{it:thm-symmetric-BOGP}}, our main hardness results for the general and symmetric settings.
We begin with some preparatory lemmas:

\begin{lem}\label{l:confinement.consequence.of.variance.bound}
	Fix any $\alpha,L > 0$, $\fq \in [1,2]$, and $\mu \in \cP_2(\R)$. 
	Then, for any $L$-Lipschitz algorithm $\cA_N$,
	\begin{align}\label{e:confinement.consequence.of.variance.bound}
		(\E \bbW_\fq(\mu_\bG(\cA_N(\bG,\bg^\aux)),\mu))^2
		&\ge \bbW_\fq(\mu(\cA_N),\mu)^2 - o_N(1)\,,\\
	\label{e:confinement.consequence.of.variance.bound.ising}
		(\E \bbW_2(\EmpDist(\cA_N(\bG,\bg^\aux)),\cP(\{\pm 1\})))^2
		&\ge \bbW_2(\mu^\Ising(\cA_N),\cP(\{\pm 1\}))^2 - o_N(1)\,.
	\end{align}
where $o_N(1)$ denotes a term tending to $0$ as $N\to\infty$, which may depend on $\alpha,L,\mu$.

\begin{proof}
	Abbreviate $X = \bbW_\fq(\mu_\bG(\cA_N(\bG,\bg^\aux)),\mu)$.
	By Jensen's inequality, we have
	\[
		\E[X^\fq]^{2/\fq} 
		\le \E[X^2] 
		= (\E X)^2 + \E[(X-\E X)^2] \stackrel{\eqref{e:Wq.lip.to.fixed.msr.conc.var}}{=} (\E X)^2 + o_N(1)\,.
	\]
By convexity of $\bbW_\fq(\cdot,\cdot)^\fq$, we have
	\[
		\E[X^\fq] = \E [\bbW_\fq(\mu_\bG(\cA_N(\bG,\bg^\aux)),\mu)^\fq] \ge \bbW_\fq(\mu(\cA_N),\mu)^\fq\,.
	\]
	Rearranging yields the first assertion \eqref{e:confinement.consequence.of.variance.bound}.
	The second assertion
	\eqref{e:confinement.consequence.of.variance.bound.ising} is proved in the same way, replacing $X$ with $J \equiv \bbW_2(\EmpDist(\cA_N(\bG,\bg^\aux)),\cP(\{\pm 1\}))$ in the above argument and using \eqref{e:Wq.lip.to.fixed.msr.conc.var.ising} from Lemma~\ref{l:Wq.lip.to.fixed.msr.conc} in place of \eqref{e:Wq.lip.to.fixed.msr.conc.var}.
\end{proof}
\end{lem}

\begin{lem}\label{l:confinement.Wq.point.to.set.whp}
	For any $\alpha,L,\epsilon,\delta>0$, and $\fq\in [1,2)$, there exists $c = c(\alpha,L,\delta) > 0$ such that the following holds for sufficiently large $N$. For any $L$-Lipschitz algorithm $\cA_N$, we have
	\[
		\P\bigg(
			\bbW_\fq\Big(\mu_\bG(\cA_N(\bG,\bg^\aux)), \cM^{\Lip}(\alpha;\epsilon)
			\Big) \le \delta
		\bigg) \ge 1 - e^{-cN}\,,
	\]
for $\cM^{\Lip}(\alpha;\epsilon)$ as defined by \eqref{e:cM.Lip.union.L}. The analogous estimate holds likewise for $\bbW_\fq(\mu_{\bG,\sym}(\cA_N(\bG,\bg^\aux))$ and $\cM^{\Lip,\sym}(\alpha;\epsilon))$, where the latter is defined by \eqref{e:cM.Lip.union.L.SYM}. 

\begin{proof}
	Throughout this proof, $c$ denotes a constant that depends on $\alpha,L,\delta$, which can change line by line.
	By \eqref{e:lip.alg.l2.l4.bound.tail} from Lemma~\ref{l:lip.alg.l2.l4.bound}, for $C = C(\alpha)$ as defined therein, we have
	\beq\label{e:confinement.mu.bG.to.P2}
		\P\bigg(
			\mu_\bG(\cA_N(\bG,\bg^\aux)) \in \cP_2(\R,C)
		\bigg) \ge 1 - e^{-cN}\,.
	\eeq
We now define
	\[
		S_\delta \equiv
		\Big\{ \mu \in \cP_2(\R,C) : \bbW_\fq(\mu, \cM^{\Lip}(\alpha;\epsilon)) \ge \delta \Big\}
	\]
(note this differs slightly from the set $S_\delta$ defined in the proof of Lemma~\ref{l:confinement.lip.trunc.to.ideal}). 
	By Lemma~\ref{l:confinement.cPp.compact.Wq}, the space $\cP_2(\R,C)$ endowed with the $\bbW_\fq$ metric is compact.
	The set $S_\delta$ is closed with respect to the $\bbW_\fq$ metric, and therefore also compact. We may thus cover $S_\delta$ with finitely many $\bbW_\fq$-balls,
	\[S_\delta = \bigcup_{\mu \in T_\delta} \cB_{\fq,\delta/5}(\mu)
		\,,\]
	where $T_\delta \subseteq S_\delta$ is a finite set depending on only $\fq, \delta$, and
		\[\cB_{q,\delta/5}(\mu) 
		\equiv\bigg\{
			\mu' \in S_\delta : \bbW_\fq(\mu,\mu') \le 
			\frac{\delta}{5}
		\bigg\}\,.\]
	Note that for each $\mu \in S_\delta$, we have
	\begin{align*}
		\bbW_\fq(\mu(\cA_N),\mu) 
		&\ge \bbW_\fq(\mu,\cM^{\Lip}(\alpha;\epsilon))
		- \bbW_\fq(\mu(\cA_N),\cM^{\Lip}(\alpha;\epsilon)) \\
		&\ge \delta - \bbW_\fq(\mu(\cA_N),\cM^{\Lip}(\alpha;\epsilon))\,,
	\end{align*}
where the last term on the right-hand side above does not depend on $
\mu$, and is $o_N(1)$ by
	Proposition~\ref{p:confinement}. It follows that for $N$ large enough, we have
	\[
		\inf_{\mu \in S_\delta} \bbW_\fq(\mu(\cA_N),\mu) \ge 
		\frac{4\delta}{5}\,.
	\]
	For each $\mu \in T_\delta$,  Lemma~\ref{l:confinement.consequence.of.variance.bound} implies
	\[
		\E \bbW_\fq\Big(\mu_\bG(\cA_N(\bG,\bg^\aux)),\mu\Big)
		\stackrel{\eqref{e:confinement.consequence.of.variance.bound}}{\ge} \bigg(\Big[\bbW_\fq(\mu(\cA_N),\mu)^2 - o_N(1)\Big]_+\bigg)^{1/2}
		\ge \frac{3\delta}{5}
	\]
	for sufficiently large $N$.
	Note that the $o_N(1)$ term above may depend on $\mu$; however, as $T_\delta$ is a finite set independent of $N$, the above estimate can be taken to hold simultaneously for all $\mu \in T_\delta$. By Lemma~\ref{l:Wq.lip.to.fixed.msr.conc} combined with a union bound over
$\mu \in T_\delta$, we obtain
	\[
		\P\bigg(
			\min_{\mu \in T_\delta} \bbW_\fq\Big(\mu_\bG(\cA_N(\bG,\bg^\aux),\mu\Big)
			\ge \frac{2\delta}{5}
		\bigg)
		\stackrel{\eqref{e:Wq.lip.to.fixed.msr.conc.tail}}{=} 1 - e^{-cN}\,.
	\]
	(This $c$ depends on $\alpha,L,\delta$ only, as $\epsilon,\fq$ affect only the $N$-independent number of terms $|T_\delta|$ in the union bound.)
	On this event, we further have
	\[\bbW_\fq\Big(\mu_\bG(\cA_N(\bG,\bg^\aux)), S_\delta\Big)
		\ge \min_{\mu \in T_\delta} \bbW_\fq\Big(\mu_\bG(\cA_N(\bG,\bg^\aux)),\mu\Big) - \frac{\delta}{5}
		\ge 
		\frac{\delta}{5}\,,\]
	and therefore $\mu_\bG(\cA_N(\bG,\bg^\aux)) \not\in S_\delta$.
	Combining with \eqref{e:confinement.mu.bG.to.P2} shows that 
	\[
		\mu_\bG(\cA_N(\bG,\bg^\aux)) \in \cP_2(\R,C) \setminus S_\delta
	\]
with probability $1-e^{-cN}$ (again adjusting $c$ appropriately). 
	On this event, we obtain the desired estimate
	\[
		\bbW_\fq\Big(\mu_\bG(\cA_N(\bG,\bg^\aux)), \cM^{\Lip}(\alpha;\epsilon)
		\Big) \le \delta\,.
	\]
	The analogous estimate for $\bbW_\fq(\mu_{\bG,\sym}(\cA_N(\bG,\bg^\aux))$ and $\cM^{\Lip,\sym}(\alpha;\epsilon))$ is proved in the same way.
\end{proof}
\end{lem}

\begin{lem}\label{l:confinement.cM.lip.to.nbd.of.ideal}
For any $\alpha,\iota > 0$ and $\fq\in [1,2)$ there exists $\epsilon = \epsilon(\alpha,\iota,\fq) > 0$ such that 
	\[
	\max\bigg\{
	\bbW_\fq(\mu,\ocM^{\Ising}(\alpha)) 
	: \mu\in \cM^{\Lip}(\alpha;\epsilon)
	\bigg\}
	\le \iota\,.
	\]
The same statement holds with
$\ocM^{\Ising,\sym}(\alpha)$ and
$\cM^{\Lip,\sym}(\alpha;\epsilon)$ in place of $\ocM^{\Ising}(\alpha)$
and $\cM^{\Lip}(\alpha;\epsilon)$. (The sets mentioned above specified by Definition~\ref{d:achievable-msrs}, \eqref{e:cM.Lip.union.L}, and \eqref{e:cM.Lip.union.L.SYM}.)

\begin{proof} Let
$\cA_N$ be an $L$-Lipschitz algorithm. 
If $\mu \in \cP_2(\R)$ satisfies $\bbW_2(\mu(\cA_N),\mu) \le \epsilon$, then Lemma~\ref{l:lip.alg.l2.l4.bound} gives, for $N$ large enough,
	\[
	\|\mu\|_{L^2} 
	\le 
	\|\mu(\cA_N)\|_{L^2} 
	+ \bbW_2(\mu(\cA_N),\mu)
	\stackrel{\eqref{e:lip.alg.l2.l4.bound.l2}}{\le} 
	C(\alpha)^{1/2} + \epsilon
	\le (2C(\alpha))^{1/2}\,,
	\]
where the last bound holds for $\epsilon$ small enough. For such $\epsilon$, it follows by recalling \eqref{e:cM.Lip.union.L} and \eqref{e:cM.Lip} that 
	\[\cM^{\Lip}(\alpha;\epsilon) 
	\stackrel{\eqref{e:cM.Lip.union.L}}{=}
	\bigcup_{L>0}
		\cM^{\Lip}(\alpha;L,\epsilon) \subseteq \cP_2(\R,2C(\alpha))\,.
	\]
For $\iota > 0$, define the set
	\[
		S_\iota = \Big\{\mu \in \cP_2(\R,2C(\alpha)) : \bbW_\fq(\mu,\ocM^{\Ising}(\alpha)) \ge \iota\Big\}\,.
	\]
	Now suppose for contradiction that $\cM^{\Lip}(\alpha;\epsilon)$ intersects $S_\iota$ for all $\epsilon > 0$.
	Recall from Lemma~\ref{l:confinement.cPp.compact.Wq} that the $\cP_2(\R,2C(\alpha))$
	endowed with the $\bbW_\fq$ metric
	 is compact.
	The set $S_\iota$ is closed with respect to this metric. Thus, for any sufficiently small $\epsilon > 0$, the set
	\[
		S_\iota \cap \overline{\cM^{\Lip}(\alpha;\epsilon)}
	\]
	(where the overline now denotes closure with respect to the $\bbW_\fq$ metric) is a nonempty closed subset of $\cP_2(\R,2C(\alpha))$ indexed by $\epsilon$, hence also compact.
	These sets form a nested family of nonempty compact sets, so their intersection is also nonempty.
	However, their intersection is
	\beq\label{e:confinement.cM.lip.to.nbd.of.ideal}
		S_\iota \cap \bigcap_{\epsilon > 0} \overline{\cM^{\Lip}(\alpha;\epsilon)}
		\subseteq  
		S_\iota \cap \bigcap_{\epsilon > 0} \cM^{\Lip}(\alpha;2\epsilon)
		\stackrel{\eqref{e:sM.Lip}}{=} 
		S_\iota \cap \cM^{\Lip}(\alpha)
		\stackrel{\textup{Thm.~\ref{thm:control.problems.main}}}{=} S_\iota \cap \ocM^{\Ising}(\alpha)
		= \emptyset\,,
	\eeq
which yields a contradiction.
This proves that there exists $\epsilon$ such that $\bbW_\fq(\mu,\ocM^{\Ising}(\alpha)) \le \iota$ for all $\mu \in \cM^{\Lip}(\alpha;\epsilon)$.
	The analogous claim for $\cM^{\Lip,\sym}(\alpha;\epsilon)$ and $\ocM^{\Ising,\sym}(\alpha)$ is proved similarly, using Theorem~\ref{thm:control.problems.sym} in place of Theorem~\ref{thm:control.problems.main} in \eqref{e:confinement.cM.lip.to.nbd.of.ideal}.
\end{proof}
\end{lem}

\begin{proof}[\hypertarget{proof:t.main.hardness.results}{Proofs of Theorem~\ref{thm:main}\ref{it:thm-main-BOGP} and Theorem~\ref{thm:symmetric}\ref{it:thm-symmetric-BOGP}}]
	Given parameters $\alpha,\iotamsr,\fq$, let $\epsilon = \epsilon(\alpha,\iotamsr/2,\fq) > 0$ be as given by Lemma~\ref{l:confinement.cM.lip.to.nbd.of.ideal}.
	We set $\iotasol = \epsilon / 3$.
Take $L>0$ and suppose $\cA_N$ is an $L$-Lipschitz algorithm.
	We will show that for sufficiently large $N$,
	\beq\label{e:main.BOGP.goal}
		\P\bigg(
			\cA_N(\bG,\bg^\aux) \in \Sigma(\iotasol)\,\text{and}\,
			\bbW_q\Big(\mu_\bG(\cA(\bG,\bg^\aux)),\ocM^{\Ising}(\alpha)
			\Big) \ge \iotamsr
		\bigg) \le e^{-cN}\,.
	\eeq
We divide the proof into two cases, depending on whether $\cA_N$ is an $\epsilon$-approximate $L$-Lipschitz algorithm (see Definition~\ref{d:eps.approx.Lip}):
\begin{itemize}
\item If $\cA_N$ is \emph{not} a $\epsilon$-approximate $L$-Lipschitz algorithm, then
	\[
		\bbW_2\Big(\mu^\Ising(\cA_N), \cP(\{\pm1\})\Big) \ge \epsilon\,.
	\]
	Then Lemma~\ref{l:confinement.consequence.of.variance.bound} gives
	\begin{align*}
		&\E \bbW_2\Big(\EmpDist(\cA_N(\bG,\bg^\aux)),\cP(\{\pm 1\}))\Big)\\
&\qquad		\stackrel{\eqref{e:confinement.consequence.of.variance.bound.ising}}{\ge} 
		\bigg(
			\Big[\bbW_2\Big(\mu^\Ising(\cA_N),\cP(\{\pm 1\})\Big)^2 - o_N(1)
			\Big]_+
		\bigg)^{1/2}
		\ge \frac{2\epsilon}{3}\,,
	\end{align*}
where the last bound holds for all $N$ large enough.
	By \eqref{e:Wq.lip.to.fixed.msr.conc.tail.ising} from Lemma~\ref{l:Wq.lip.to.fixed.msr.conc}, with probability $1-e^{-cN}$ we have
	\[\frac{\dist(\cA_N(\bG,\bg^\aux),\Sigma_N)}{N^{1/2}}
	=
		\bbW_2\Big(\EmpDist(\cA_N(\bG,\bg^\aux)), \cP(\{\pm 1\})\Big) 
		> \frac{\epsilon}{3} = \iotasol\,.
	\]
On this event $\cA_N(\bG,\bg^\aux) \not\in \Sigma(\iotasol)$, and \eqref{e:main.BOGP.goal} is proved.

\item If $\cA_N$ is an $\epsilon$-approximate $L$-Lipschitz algorithm, then Lemma~\ref{l:confinement.Wq.point.to.set.whp} gives
	\[
		\P\bigg(
			\bbW_\fq\Big(\mu_\bG(\cA_N(\bG,\bg^\aux)), \cM^{\Lip}(\alpha;\epsilon)
			\Big) \ge \frac{\iotamsr}{2}
		\bigg) \le e^{-cN}\,.
	\]
By the triangle inequality combined with Lemma~\ref{l:confinement.cM.lip.to.nbd.of.ideal}, we can bound
	\begin{align*}
	& \bbW_\fq\Big(\mu_\bG(\cA_N(\bG,\bg^\aux)), \ocM^{\Ising}(\alpha)\Big)\\
	&\qquad \le
	\min\bigg\{
	\bbW_\fq\Big(\mu_\bG(\cA_N(\bG,\bg^\aux)), \mu \Big)
	:\mu \in \cM^{\Lip}(\alpha;\epsilon)
	\bigg\}\\
	&\qquad\qquad + \max\bigg\{
	 \bbW_\fq\Big(\mu,
	\ocM^{\Ising}(\alpha)
	\Big)
	:\mu \in \cM^{\Lip}(\alpha;\epsilon)
	\bigg\}\\
	&\qquad \le
	\bbW_\fq
	\Big(\mu_\bG(\cA_N(\bG,\bg^\aux)), \cM^{\Lip}(\alpha;\epsilon) \Big)
	+ \frac{\iotamsr}{2}\,,
	\end{align*}
which implies the inclusion of events
	\begin{align*}
		&\bigg\{
		\bbW_\fq\Big(
		\mu_\bG(\cA_N(\bG,\bg^\aux)), \cM^{\Lip}(\alpha;\epsilon)
		\Big) \ge 
		\frac{\iotamsr}{2}
		\bigg\} \\
		&\qquad
		\supseteq 
		\bigg\{
		\bbW_\fq\Big(
		\mu_\bG(\cA_N(\bG,\bg^\aux)), \ocM^{\Ising}(\alpha)
		\Big)
		\ge \iotamsr
		\bigg\}\,,
	\end{align*}
	which proves \eqref{e:main.BOGP.goal}.
\end{itemize}
Theorem~\ref{thm:main}\ref{it:thm-main-BOGP} now follows.
	The proof of Theorem~\ref{thm:symmetric}\ref{it:thm-symmetric-BOGP} is analogous.
\end{proof}

\subsection{Algorithmic threshold for perceptron model}
\label{ss:confinement.perceptron.threshold.results}

In this subsection we give the
proof of Corollary~\ref{cor:alg-for-optimization}. We begin with a simple lemma: \begin{lem}\label{l:integration.against.phi.unif.cts}
	Fix any $\phi \in C_b(\R,\R)$ and define the domain
	\[
		\ocM^{\Ising+}(\alpha) = \Big\{
			\mu \in \cP_2(\R,C'(\alpha)) : \bbW_1(\mu, \ocM^{\Ising}(\alpha)) \le 1
		\Big\}\,.
	\]
	Then the function
	\[
		f_\phi(\mu) \equiv \int \phi\,d\mu
	\]
	is uniformly continuous on $\ocM^{\Ising+}(\alpha)$ equipped with the $\bbW_1$ metric.

\begin{proof}
The set $\ocM^{\Ising+}(\alpha)$ is clearly closed with respect to the $\bbW_1$ metric. By Lemma~\ref{l:confinement.cPp.compact.Wq}, the space $\cP_2(\R,C'(\alpha))$ endowed with the $\bbW_1$ metric is compact, so $\ocM^{\Ising+}(\alpha)$ is also compact. It thus suffices to show that $f_\phi$ is continuous. In fact, since we are on a metric space, it suffices to show that $f_\phi$ is sequentially continuous. To this end,  consider any sequence $\mu_n$ in $\ocM^{\Ising+}(\alpha)$ converging to $\mu$ in $\bbW_1$. Then $\mu_n$ also converges weakly to $\mu$. Since $\phi$ was assumed to be a bounded continuous function, we obtain 
	 \[\lim_{n\to\infty} f_\phi(\mu_n) = f(\mu)\,.\]
Therefore $f_\phi$ is indeed continuous, which concludes the proof.
\end{proof}
\end{lem}

Let $C(\alpha)$ be as given by Lemma~\ref{l:lip.alg.l2.l4.bound} and define
\[
	C'(\alpha) = \max\bigg\{
	C(\alpha), 
	\frac{10}{\alpha} + 
	7 \bigg\}\,.
\]
Recall $\ocM^\Ising(\alpha) = \sM^\sideal(\alpha,0)$ by the \hyperlink{p:thm.control.problems.main}{proof of Theorem~\ref{thm:control.problems.main}}. Then Lemma~\ref{l:L2-bound-on-stochastic-control} implies 
\beq\label{e:integration.against.phi.C.prime.1}
	\sup \bigg\{ 
		\|\mu\|_{L^2}^2 : 
		\mu \in 
		\ocM^\Ising(\alpha)
		= \sM^\sideal(\alpha,0)
	\bigg\} 
	\le
	\frac{10}{\alpha} +  6 + O(\epsilon)
	 \le C'(\alpha)\,.
\eeq
Lemma~\ref{l:lip.alg.l2.l4.bound} implies that for any $L>0$ there exists $c = c(\alpha,L) > 0$ such that
\beq\label{e:integration.against.phi.C.prime.2}
	\P\bigg( \|\mu_\bG(\cA_N(\bG,\bg^\aux))\|_{L^2}^2 \le C'(\alpha) \bigg) \stackrel{\eqref{e:lip.alg.l2.l4.bound.tail}}{\ge} 1 - e^{-cN}\,.
\eeq
for any $L$-Lipschitz algorithm $\cA_N$.

\begin{proof}[Proof of Corollary~\ref{cor:alg-for-optimization}\ref{it:cor-optimization-IAMP}]
	Consider any $\phi \in C_b(\R,\R)$ and $\iota > 0$.
Let $f_\phi$ and $\ocM^{\Ising+}(\alpha)$ be defined in Lemma~\ref{l:integration.against.phi.unif.cts}: by this lemma, we can find $\iota'>0$ such that
 \beq\label{e:integration.against.phi.unif.cts.consequence}\sup\bigg\{
 |f_\phi(\mu) - f_\phi(\mu')| 
 :
 \mu,\mu' \in \ocM^{\Ising+}(\alpha),
 \bbW_1(\mu,\mu') \le \iota'
 \bigg\}
 \le \frac{\iota}{2}
 \eeq
Without loss of generality we may assume $\iota' \le \min\{\iota, 1\}$.
	By the definition \eqref{eq:def-ALG-Ising} of $\ALG$, we can find $\mu \in \ocM^\Ising(\alpha)$ such that
	\[
		f_\phi(\mu) \ge \ALG - \iota / 2\,.
	\]
	By Theorem~\ref{thm:main}\ref{it:thm-main-IAMP}, there exists $L$ and an $L$-Lipschitz algorithm $\cA_N$ such that with probability $1-e^{-cN}$, the output $\bx = \cA_N(\bG,\bg^\aux)$ satisfies
	\begin{align}
	\label{e:cor-optimization-IAMP.coord}
	&\bx \in \Sigma_N(\iota') \subseteq \Sigma_N(\iota)\,,\\
	&\bbW_2(\mu_\bG(\bx),\mu) 
	\le \iota'\,.
	\label{e:cor-optimization-IAMP.mu.approx}
	\end{align}
Suppose \eqref{e:cor-optimization-IAMP.coord}, \eqref{e:cor-optimization-IAMP.mu.approx} and the event in \eqref{e:integration.against.phi.C.prime.2} all hold; this occurs with probability $1-e^{-cN}$.
It follows from
\eqref{e:integration.against.phi.C.prime.1} and \eqref{e:integration.against.phi.C.prime.2} that $\mu_\bG(\bx), \mu \in \cP_2(\R,C'(\alpha))$.
	Further, we have $\mu \in \ocM^\Ising(\alpha)$ by assumption while
	\[
		\bbW_1(\mu_\bG(\bx), \ocM^\Ising(\alpha))
		\le \bbW_1(\mu_\bG(\bx), \mu)
		\le \bbW_2(\mu_\bG(\bx), \mu)
		\stackrel{\eqref{e:cor-optimization-IAMP.mu.approx}}{\le}
		\iota' \le 1\,,
	\]
so both measures $\mu_\bG(\bx), \mu$ belong to $\ocM^{\Ising+}(\alpha)$. Recalling \eqref{e:integration.against.phi.unif.cts.consequence}, the bound $\bbW_1(\mu_\bG(\bx), \mu) \le \iota'$ implies
	\[
\Big|f_\phi(\mu_\bG(\bx)) - f_\phi(\mu)
\Big| \le \frac{\iota}{2}\,.
	\]
	We conclude that on this event, for $H_N(\bx)$ defined in \eqref{e:hamiltonian},
	\[
		\ALG - \iota \le f_\phi(\mu_\bG(\bx))
		= \int \phi\, d\mu_\bG(\bx)
		= \frac{H_N(\bx) }{M}\,.
	\]
	Together with \eqref{e:cor-optimization-IAMP.coord} this implies the result.
\end{proof}

\begin{proof}[Proof of Corollary~\ref{cor:alg-for-optimization}\ref{it:cor-optimization-BOGP}]
Consider any $\phi \in C_b(\R,\R)$ and $\iotaval > 0$.
	By Lemma~\ref{l:integration.against.phi.unif.cts}, we can find $\iotamsr > 0$ such that 
	\beq
	\label{e:integration.against.phi.unif.cts.second}
	\sup\bigg\{
	|f_\phi(\mu) - f_\phi(\mu')|
	: \mu,\mu' \in \ocM^{\Ising+}(\alpha),
	\bbW_1(\mu,\mu') \le \iotamsr
	\bigg\}
	 \le \frac{\iotaval}{2}
	 \eeq
Without loss of generality we may assume $\iotamsr \le 1$.
By Theorem~\ref{thm:main}\ref{it:thm-main-BOGP}, there exists $\iotasol > 0$ such that for any $L$, there exists $c>0$ such that all $L$-Lipschitz algorithms $\cA_N$, the output $\bx = \cA_N(\bG,\bg^\aux)$ satisfies
	\[
		\P\bigg(
			\bx \in \Sigma_N(\iotasol) \,\,\text{and}\,\,
			\bbW_1(\mu_\bG(\bx), \ocM^{\Ising}(\alpha)) \ge \iotamsr
		\bigg) \le e^{-cN}\,.
	\]
Taking a union bound with \eqref{e:integration.against.phi.C.prime.2} gives
	\[
		\P\bigg(
		\Big\{
		\bx \in \Sigma_N(\iotasol)
		\Big\}
		\cap \Big\{
				\bbW_1(\mu_\bG(\bx), \ocM^{\Ising}(\alpha)) \ge \iotamsr \,\,\text{or}\,\, \mu_\bG(\bx)\not\in \cP_2(\R,C'(\alpha))
			\Big\}
		\bigg) \le e^{-cN}\,.
	\]
	We will show the inclusion of events
	\beq\label{e:cor-optimization-BOGP.goal}
		\bigg\{
			\frac{H_N(\bx)}{M}
			\ge \ALG + \iotaval
		\bigg\} \subseteq \bigg\{
			\bbW_1(\mu_\bG(\bx), \ocM^{\Ising}(\alpha)) \ge \iotamsr \,\,\text{or}\,\, \mu_\bG(\bx)\not\in \cP_2(\R,C'(\alpha))
		\bigg\}\,,
	\eeq
	which directly implies the conclusion.
	This is equivalent to the contrapositive
	\beq\label{e:cor-optimization-BOGP.goal2}
		\bigg\{
			\bbW_1(\mu_\bG(\bx), \ocM^{\Ising}(\alpha)) < \iotamsr \,\,\text{and}\,\, \mu_\bG(\bx) \in \cP_2(\R,C'(\alpha))
		\bigg\} \subseteq 
		\bigg\{
		\frac{H_N(\bx)}{M} < \ALG + \iotaval
		\bigg\}\,.
	\eeq
Suppose that the event on the above left-hand side holds: then there exists
some $\mu \in \ocM^{\Ising}(\alpha)$ such that $\bbW_1(\mu_\bG(\bx), \mu) < \iotamsr$, so we have
	\[
		\bbW_1(\mu_\bG(\bx), \ocM^\Ising(\alpha))
		\le \bbW_1(\mu_\bG(\bx), \mu)
		\le \iotamsr \le 1\,.
	\]
We also have $\mu \in \cP_2(\R,C'(\alpha))$ by \eqref{e:integration.against.phi.C.prime.1}, so we conclude that the measures
$\mu_\bG(\bx), \mu$ both lie in $\ocM^{\Ising+}(\alpha)$. Recalling \eqref{e:integration.against.phi.unif.cts.second} gives the bound 
	\begin{align*}
		\frac{H_N(\bx)}{M}
		= \int \phi \,d\mu_\bG(\bx)
		= f_\phi(\mu_\bG(\bx))
		&\le f_\phi(\mu) + 
		\Big|f_\phi(\mu_\bG(\bx)) - f_\phi(\mu)\Big| \\
		&\le \sup\Big\{
		f_\phi(\mu') : \mu' \in \ocM^\Ising(\alpha)\Big\} 
		+ \frac{\iotaval}{2}
		= \ALG + \frac{\iotaval}{2}\,.
	\end{align*}
This proves \eqref{e:cor-optimization-BOGP.goal2}, and the result follows as noted above.
\end{proof}

To conclude, we recall that
Theorem~\ref{thm:main} part~\eqref{it:thm-main-IAMP} was proved in \S\ref{ss:confinement.main.achievability.results}, part~\eqref{it:thm-main-BOGP} was proved in \S\ref{ss:confinement.main.hardness.results}, and part~\eqref{it:thm-main-concave} was previously proved in Theorem~\ref{thm:control.problems.main}.
Theorem~\ref{thm:symmetric} part~\eqref{it:thm-symmetric-IAMP} was proved in \S\ref{ss:confinement.main.hardness.results}, and part~\eqref{it:thm-symmetric-BOGP} was proved in \S\ref{ss:confinement.main.hardness.results}.

\fi

\appendix

\pagebreak\section{Moment estimates and tightness}\label{s:kolmogorov}

\iffull
% !TEX root = main.tex

In this section we prove Theorem~\ref{t:tightness}. 
For the reader's convenience we start by reviewing key definitions from Section~\ref{s:rerand}. Recall from \eqref{e:process.X} the process $X=X^N$ of inner products revealed ``going down the tree.'' In the Ising case we additionally have the process $X^\Ising$ defined by \eqref{e:process.X.Ising}. We primarily focus on the Ising case since it is more difficult, and the spherical case is very similar but easier.
The processes are defined on discrete times $q_d$ where the time steps $q_{d+1}-q_d$ are of order $\delta$.
\begin{itemize}
\item In \eqref{e:X.decomp}--\eqref{e:vX} we defined the process
$\vX=(\vX^{\RomI:\RomIII},X^{\Ising})$. (Again, for the spherical perceptron we simply do not keep track of $X^{\Ising}$, but we will focus our discussion on the Ising perceptron.)

\item In Definition~\ref{d:trunc} we specified the truncated process $\vX^\trstar$, where $\vX^{\trstar,\RomI:\RomIII}$ is frozen upon leaving 
$[-\trK,\trK)^3$, and 
$\vX^{\trstar,\Ising}$ is frozen upon leaving $[-\trK,\trK)$.

\item In Definition~\ref{d:froze} we specified a ``smoothed'' variant $\vX^\frstar$, where
$\vX^{\frstar,\RomI:\RomIII}$ has an additional $1/\trK^6$ density
of frozen particles spaced evenly over $[-2\trK,2\trK]^3$; and similarly 
$\vX^{\frstar,\Ising}$ has an additional $1/\trK^3$ density
of frozen particles spaced evenly over $[-2\trK,2\trK]$.

\item In Definitions~\ref{d:blocks}--\ref{d:buckets.Ising}
we partitioned space into blocks of width $\eta$, and defined 
$\vX^\trunc$ and $\vX^\frozen$ to be variants of 
$\vX^\trstar$ and $\vX^\frstar$ with small spatial buckets frozen.
 
\item In Definition~\ref{d:rerand} we constructed
spatial rerandomizations $\vY$ and $\vY^\frozen$, for $\vX^\trunc$ and $\vX^\frozen$ respectively.

\end{itemize}
\textbf{In this section, we show that the processes $\vY$ and $\vY^\frozen$ have a well-defined scaling limit that has a semimartingale decomposition.}

The section is organized as follows:
\begin{itemize}
\item In \S\ref{ss:lip.conc} we record some preliminary estimates on Lipschitz bounds for some basic quantities that appear throughout this section.

\item In \S\ref{s:martingale.subgaussian.estimates} we show that the preprocessing steps of \S\ref{ss:preprocessing} incur small $\bbW_2$ error, provided the original process is close in $\bbW_2$ distance to a fixed measure (see assumption \eqref{e:w2.coupling.assumption}).

\item In \S\ref{ss:kolmogorov.Y} we prove Proposition~\ref{p:Y.kolmogorov}, which gives
annealed moment estimates for $\vY$ and $\vY^\frozen$.
\item In \S\ref{ss:kolmogorov.drift} we define discrete-time processes $\vD$ and $\vD^\frozen$ (see Definition~\ref{d:D}) that describe the ``drift'' of $\vY$ and $\vY^\frozen$, and prove annealed moment estimates for these processes (Proposition~\ref{p:drift.kolmogorov}).
\item In \S\ref{ss:kolmogorov.qv} we define discrete-time processes
$\vQ$ and $\vQ^\frozen$ (see Definition~\ref{d:Q}) that describe the ``quadratic variation'' and ``covariation'' of the components of $\vY$ and $\vY^\frozen$, and prove annealed moment estimates for these processes (Proposition~\ref{p:qv.kolmogorov}).

\item In \S\ref{ss:tightness} we apply the moment estimates of 
\S\ref{ss:kolmogorov.Y}--\ref{ss:kolmogorov.qv}, and argue similarly as in the standard proof of the Kolmogorov continuity lemma (see e.g.\ \cite[Thm.~2.9]{MR3497465}). This implies that,
as $N\to\infty$ with an appropriate set of parameters (see Assumption~\ref{a:params}),
the processes $\vY$, $\vY^\frozen$, $\vD$, $\vD^\frozen$, $\vQ$, and $\vQ^\frozen$ have subsequential scaling limits, which are continuous stochastic processes --- see Proposition~\ref{p:tightness} and Corollary~\ref{c:tightness}. 
\item We further show in Proposition~\ref{p:limiting.drift.qv} that the scaling limits of $\vD$ and $\vQ$ in fact correspond to the drift and covariation of the scaling limit of $\vY$, and similarly for the processes $\vY^\frozen$, $\vD^\frozen$, $\vQ^\frozen$.
\item The \hyperlink{proof:t.tightness}{proof of Theorem~\ref{t:tightness}} is given at the end of this section, and follows by collecting the above assertions.
\end{itemize}
In what follows it will be useful to recall the notation $p'(q_d)$ introduced in \eqref{e:deriv.p.bound}. Recall also the notation introduced in Remark~\ref{r:annealed.quenched} for quenched versus annealed probabilities.

\subsection{Preliminary estimates}
\label{ss:lip.conc}

In this subsection we prove that the quantities defined in \S\ref{ss:X.decomp} are suitably Lipschitz, or approximately Lipschitz, with respect to the gaussian disorder. The estimates written here will be used in later sections to obtain concentration of various quantities in the large-system limit. As above, let $\bXi\equiv (\bXi^\ell : 0\le\ell\le\dmax)$ denote a tuple of $M\times N$ matrices, with i.i.d.\ standard gaussian entries. 
Recall the notation introduced in
\eqref{e:gaus.decomp} and \eqref{e:normalized.vectors.g}--\eqref{e:normalized.vectors.y}, and denote
	\beq\label{e:bar.bG}
	\bar{\bG}
	\equiv
	\bar{\bG}(q_d)
	\equiv \frac{\bG(q_d)}
	{\sqrt{p_d}}
	=\frac1{\sqrt{p_d}}
	\sum_{\ell=0}^d
	(p_\ell-p_{\ell-1})^{1/2}
	\bXi^\ell
	\,,
	\eeq
so that $\bar{\bG}$ also has i.i.d.\ standard gaussian entries. When $p_d=0$, we use the convention from \S\ref{ss:X.decomp} that $\bar{\bG}(q_d)$ is an independent standard gaussian matrix. As above, we write $\bar{\bg}^a\equiv \bar{\bg}^a(q_d)$ for the rows of $\bar{\bG}(q_d)$, and $\bmeta^a \equiv \bmeta^a(q_{d+1})$ for the rows of $\bXi^{d+1}$. Recall the filtration $\cG(q_d)$ defined by \eqref{e:gaus.filt}. 

\begin{lem}[$\bar{\bx}$ is Lipschitz]\label{l:bar.x.subgaus}
Take $\bar{\bx}\equiv \bar{\bx}(q_d)$
 as in \eqref{e:normalized.vectors.x}, so it is measurable with respect to
 $\cG(q_d)$. Then
		\[
		\|\E\bar{\bx}\|^2
		\le \frac{Nq_0}{q_d} \le N\,.\]
and the centered random vector $\bar{\bx}  -\E\bar{\bx}$ is subgaussian with variance proxy $L^{O(1)}$.

\begin{proof} 
Recall $\chi(p)=\E R(p)$ with $R(p)$ as in \eqref{e:p.corr.overlap}, and $q_d=\chi(p_d)$. Note that
	\[
	\frac{\|\E\bar{\bx}\|^2}{N}
	= \frac{\|\E\cA(\bG)\|^2}{Nq_d}
	= \frac{\chi(0)}{q_d} 
	\le \frac{\chi(p_0)}{q_d} = \frac{q_0}{q_d}\le1\,,
	\]
which proves the first claim. Next note that we can write
	\[
	\sqrt{q_d}\bar{\bx}
	\stackrel{\eqref{e:normalized.vectors.x}}{=}
	\E\bigg( \cA
	\Big(\sqrt{p_d}\bar{\bG} +\sqrt{1-p_d}\bG'\Big)
	\,\bigg|\,\bar{\bG}\bigg)\,,
	\]
where $\bar{\bG} = \bG(q_d)/\sqrt{p_d}$ and $\bG'$ are independent $M\times N$ matrices, each with i.i.d.\ standard gaussian entries. Since $\cA$ is $L$-Lipschitz, and $p_d/q_d \le L^{O(1)}$ by \eqref{e:deriv.p.bound}, it follows that $\bar{\bx}$ is $L^{O(1)}$-Lipschitz with respect to $\bar{\bG}$. The claim then follows using Lemma~\ref{l:lip.subgaus}.
\end{proof}
\end{lem}

\begin{lem}[$\by$ is Lipschitz]\label{l:y.subgaus}
Take $\by\equiv \by(q_{d+1})$ as in \eqref{e:normalized.vectors.y}, so it is measurable with respect to $\cG(q_{d+1})$, as defined by \eqref{e:gaus.filt}. 
Conditional on $\cG(q_d)$,
the random vector $\by$ has mean zero, and is subgaussian with variance proxy $L^{O(1)}$.

\begin{proof} Let us abbreviate
	\[
	\bXi
	\equiv \bXi^{d+1}
	\equiv \frac{\bG(q_{d+1})
		-\bG(q_d)}{\sqrt{p_{d+1}-p_d}}\,,
	\]
so this is a gaussian independent of $\cG(q_d)$. Then we can express
	\begin{align*}
	\sqrt{\delta_d}\by
	&\stackrel{\eqref{e:normalized.vectors.y}}{=}
	\E\bigg[ \cA\Big(
		\bG(q_d)
		+(p_{d+1}-p_d)^{1/2}\bXi
		+ (1-p_{d+1})^{1/2} \bG'\Big)\,\bigg|\,\bG(q_d),\bXi
		\bigg]\\
	&\qquad-\E\bigg[ \cA\Big(
		\bG(q_d)
		+(p_{d+1}-p_d)^{1/2}\bXi
		+ (1-p_{d+1})^{1/2} \bG'\Big)\,\bigg|\,\bG(q_d)\bigg] \,.
	\end{align*}
Thus, conditional on $\cG(q_d)$,
the vector $\by$ has mean zero.
It follows by the assumption on $\cA$ together with assumption~\eqref{e:deriv.p.bound} that
 $\by$ is $L^{O(1)}$-Lipschitz in $\bXi$, so the subgaussian claim follows using Lemma~\ref{l:lip.subgaus}.
\end{proof}
\end{lem}

\begin{lem}[$\vX$ is approximately Lipschitz]
\label{l:X.lip}
For $\sigma\in\{\RomI,\RomII,\RomIII\}$,
let $h^{\sigma,d}$ denote the function which takes $\bXi$ as input, and outputs the resulting $M$-dimensional vector $(\Delta X^\sigma(q_d,a))_{a\le M}$.
Similarly, let $h^{\Ising,d}$ denote the function which takes $\bXi$ as input, and outputs the resulting $N$-dimensional vector $(\Delta X^\Ising(q_d,i))_{i\le N}$.
For $\sigma\in\{\RomI,\RomII,\RomIII,\Ising\}$ and for all $0\le d\le\dmax-1$, the functions $h^{\sigma,d}$ are $L^{O(1)}$-Lipschitz when restricted to the subset
	\[U_\MAX\equiv 
	\bigg\{\bXi:
	\MAX(\bXi)
	\le L^{O(1)}
	\bigg\}\,,
	\]
for $\MAX$ as defined by \eqref{e:MAX.N}.

\begin{proof}
We saw in the proof of Lemma~\ref{l:bar.x.subgaus} that $\bar{\bx}$ is $L^{O(1)}$-Lipschitz with respect to $\bar{\bG}$, which in turn is $L^{O(1)}$-Lipschitz with respect to $\bXi$. It then follows by recalling  \eqref{e:normalized.vectors.y} that $\delta^{1/2}\by$ is also $L^{O(1)}$-Lipschitz with respect to $\bXi$, which 
immediately gives the claim for $h^{\Ising,d}$ (even without the restriction to $U_\MAX$). For $\sigma\in\{\RomI,\RomII,\RomIII\}$,
it follows from
\eqref{e:X.decomp} that
	\begin{align*}
	h^{\RomI,d}(\bXi) 
	&= \delta_d p'(q_d)^{1/2}
		\frac{\bXi^{d+1}\by(\bXi)}{N^{1/2}}
	\\
	h^{\RomII,d}(\bXi) 
	&= [\delta_d p(q_d)]^{1/2}
		\frac{\bar{\bG}(\bXi) \by(\bXi)}{N^{1/2}}\,,\\
	h^{\RomIII,d}(\bXi)
	&= [\delta_d q_d p'(q_d)]^{1/2}
	\frac{\bXi^{d+1} \bar{\bx}(\bXi)}{N^{1/2}}\,,
	\end{align*}
where $\bar{\bG}(\bXi)$ denotes the matrix defined by \eqref{e:bar.bG} using disorder $\bXi$;
 and $\bar{\bx}(\bXi)$ and $\by(\bXi)$ denote the vectors defined by 
 \eqref{e:normalized.vectors.x} and \eqref{e:normalized.vectors.y} using disorder $\bXi$. Then, writing $\|\cdot\|_\textup{F}$ for Frobenius norm, we have
	\begin{align*}
	&\frac{\|h^{\RomI,d}(\bXi)-h^{\RomI,d}(\bTe)\|}
		{\delta_d p'(q_d)^{1/2}}
	\le
	\frac{1}{N^{1/2}}
	\bigg\{
	\|\bXi^{d+1}\|\|\by(\bXi)-\by(\bTe)\|
	+ \|\bXi^{d+1}-\bTe^{d+1}\|\|\by(\bTe\|
	\bigg\}\\
	&\qquad\le \frac{L^{O(1)}}{N^{1/2}}
	\bigg\{
	\frac{\|\bXi^{d+1}\|}{\delta^{1/2}}
	+\|\by(\bTe)\|
	\bigg\} \|\bXi-\bTe\|_\textup{F}
	\stackrel{\star}{\le}  \frac{L^{O(1)}}{\delta^{1/2}}
	\|\bXi-\bTe\|_\textup{F}\,,
	\end{align*}
where in the last line $\star$ indicates that the bound holds provided
$\bXi,\bTe\in U_\MAX$. Similarly,
	\begin{align*}
	&\frac{\|h^{\RomII,d}(\bXi)-h^{\RomII,d}(\bTe)\|}
		{[\delta_d p(q_d)]^{1/2}}
	\le \frac{1}{N^{1/2}}
	\bigg\{ \|\bar{\bG}(\bXi)\|
		\|\by(\bXi)-\by(\bTe)\|
	+\|
		\bar{\bG}(\bXi)-\bar{\bG}(\bTe)
		\|
	\|\by(\bTe)\|\bigg\} \\
	&\qquad
	\le \frac{L^{O(1)} }{N^{1/2}}
	\bigg\{
	\frac{\|\bar{\bG}(\bXi)\|}{\delta^{1/2}}	
	+\|\by(\bTe)\| \bigg\}
	\|\bXi-\bTe\|_\textup{F}
	\stackrel{\star}{\le} \frac{L^{O(1)}}{\delta^{1/2}}
	\|\bXi-\bTe \|_\textup{F}\,\\
	&\frac{\|h^{\RomIII,d}(\bXi)-h^{\RomIII,d}(\bTe)\|}
		{[\delta_d q_d p'(q_d)]^{1/2}}
	\le \frac{1}{N^{1/2}}
		\bigg\{ 
		\|\bXi^{d+1}\|
		\|
		\bar{\bx}(\bXi)-\bar{\bx}(\bTe)\|
	+ \|\bXi^{d+1}-\bTe^{d+1}\|
	\|\bar{\bx}(\bTe)\|\bigg\}\\
	&\qquad\le \frac{L^{O(1)}}{N^{1/2}}
	\bigg\{ \|\bXi^{d+1}\|
	+ \|\bar{\bx}(\bTe)\|
	\bigg\} 
	\|\bXi-\bTe\|_\textup{F}
	\stackrel{\star}{\le} L^{O(1)}
	\|\bXi-\bTe \|_\textup{F}\,.
	\,\end{align*}
This implies the claim.
\end{proof}
\end{lem}

\begin{cor}[empirical averages are approximately Lipschitz]\label{c:conc.empir.avg}
Let $f:\R\to\R$ be any $L^{O(1)}$-Lipschitz function. For all $0\le d\le \dmax$, and $\sigma\in\{\RomI,\RomII,\RomIII\}$, the empirical averages
	\[\frac1M
	\sum_{a\in[M]}
	f(X^\sigma(q_d,a))
	\]
are $(L^{O(1)}/(\delta N^{1/2}))$-Lipschitz on $U_\MAX$. Likewise the empirical averages
	\[\frac1N
	\sum_{i\in[N]}
	f(X^\Ising(q_d,i) )
	\]
are $(L^{O(1)}/(\delta N^{1/2}))$-Lipschitz on $U_\MAX$.

\begin{proof}
Let $\sigma\in\{\RomI,\RomII,\RomIII\}$.
Recall from Lemma~\ref{l:X.lip}
that the $\R^M$-valued functions $h^{\sigma,d}$ are $L^{O(1)}$-Lipschitz on the set $U_\MAX$. It follows that
the $\R^M$-valued functions 
	\[
	H^{\sigma,d}(\bXi)
	=(X^\sigma(q_d,a))_{a\le M}
	= \sum_{\ell=0}^{d-1}
	h^{\sigma,\ell}(\bXi)
	\]
are $(L^{O(1)}/\delta)$-Lipschitz on the set $U_\MAX$. Let $F:\R^M\to\R^M$ denote the map obtained by applying $f$ coordinatewise, and note that $F$ is $L^{O(1)}$-Lipschitz. Therefore,
	\[
	\frac1M
	\sum_{a\in[M]}
	f(X^\sigma(q_d,a))
	=
	\frac1M
	(\ind_M, F(H(\bXi))
	\]
is $(L^{O(1)}/(\delta M^{1/2}))$-Lipschitz on the set $U_\MAX$. A similar argument applies for $X^\Ising$ (in this case in fact the restriction to $U_\MAX$ is not necessary, since as we mentioned in the proof of Lemma~\ref{l:X.lip}, the functions $h^{\Ising,d}$ are Lipschitz on the full space).
\end{proof}
\end{cor}

The main application of Corollary~\ref{c:conc.empir.avg} will be in Corollary~\ref{c:trunc.small} below, where we will bound the probability that the trajectory of $\vX$ covers a large distance $\trK$. We will argue that the \emph{annealed} probability of this event (averaging over the gaussian disorder) is $o_{\trK}(1)$. We will then use Corollary~\ref{c:conc.empir.avg}, along with general results on concentration for Lipschitz functionals of gaussians (see Lemma~\ref{l:lip.subgaus}) to deduce that the \emph{quenched} probability (conditional on the disorder $\bXi$)
 is $o_{\trK(1)}$ with high probability.

% !TEX root = main.tex

\subsection{Estimates on preprocessing errors} 
\label{s:martingale.subgaussian.estimates}

In this subsection, we bound in $\bbW_2$ distance the errors introduced by the preprocessing steps of \S\ref{ss:preprocessing}. The main result of this subsection is the following, which is required for the \hyperlink{proof:t.BOGP-hardness-main}{proof of Theorem~\ref{thm:BOGP-hardness-main}}: 

\begin{ppn}\label{p:summarize.W2.errors.preprocessing} Assume that $L\ge1$ and that condition \eqref{e:deriv.p.bound} holds. Let $\vX\equiv \vX^N$ be the process defined by \eqref{e:X.decomp} and \eqref{e:Delta.X.Ising}. Let $X(q_d)$ denote the random variable $X(q_d,a(q_d))$ where $a(q_d)$ is sampled uniformly at random from $[M]$. Similarly denote $\vX(q_d)\equiv \vX(q_d,(a(q_d),i(q_d)))$, where $i(q_d)$ is sampled uniformly at random from $[N]$. Assume that for some fixed measure $\mu$ on $\R$ with finite second moment, we have
	\beq\label{e:w2.coupling.assumption}
	\bbW_2(\mathscr{L}(X(1)),\mu)
	\le \epsilon\,,\eeq
where the bound holds uniformly over $N$. Then each of the preprocessing steps defined in \S\ref{ss:preprocessing} incurs a negligible $\bbW_2$ error on the endpoint distribution: the processes $\vX, \vX^\trstar, \vX^\frstar,\vX^\trunc, \vX^\frozen$ can be coupled so that
	\[
	\E\bigg[ 
	\|\vX(1)-\vX^\trstar(1)\|^2
	+\|\vX^\trstar(1)-\vX^\frstar(1)\|^2
	+\|\vX^\trstar(1)-\vX^\trunc(1)\|^2
	+\|\vX^\frstar(1)-\vX^\frozen(1)\|^2
	\bigg]
	\le O(\epsilon) + o_{\trK}(1)\,,
	\]
uniformly in $N$.
\end{ppn}

The \hyperlink{proof:p.summarize.W2.errors.preprocessing}{proof of Proposition~\ref{p:summarize.W2.errors.preprocessing}} appears at the end of this subsection. In the course of the proof, we also obtain some moment estimates (Lemma~\ref{l:martingale.sup.bound.X.II.III.Ising}, Lemma~\ref{l:martingale.sup.bound.X.I.minus.mean}, Proposition~\ref{p:max.X.process.on.trstar.stopped.event}) which are similar to moment estimates derived in the following subsection (see Lemma~\ref{l:subgaus.X.II.III.Ising} ).

\begin{lem} \label{l:martingale.sup.bound.X.II.III.Ising}
Assume that $L\ge1$ and that condition \eqref{e:deriv.p.bound} holds. Let $\vX\equiv \vX^N$ be the process defined by \eqref{e:X.decomp} and \eqref{e:Delta.X.Ising}. We then have
	\begin{align*}
	\max_{\sigma\in\{\RomII,\RomIII\}}
	\max_{a\in[M]}
	\E\bigg[ \max_{0\le d\le\dmax}
	X^\sigma(q_d,a)^8\bigg] 
	&\le L^{O(1)}\,,\\
	\max_{i\in[N]}
	\E\bigg[ \max_{0\le d\le\dmax}
	X^\Ising(q_d,i)^8\bigg] 
	&\le L^{O(1)}\,,
	\end{align*}
with the estimates holding uniformly in $N$.

\begin{proof} First consider $X^{\RomII}(q_d,a)$ for fixed $a\in[M]$, as a process indexed by time $q_d$ with $0\le d\le \dmax$. For the purposes of the proof let us abbreviate
		\[
	\frac{\Delta X^{\RomII}(q_d,a)}{[\delta_d p(q_d)]^{1/2}}
	\equiv 
	\mathfrak{y}^{\RomII}(q_d,a)
	\equiv \frac{(\bar{\bg}^a(q_d), \by(q_{d+1}))}{N^{1/2}}
	\,.
	\]
 Note that for fixed $a\in[M]$, $X^{\RomII}(q_d,a)$ (indexed by time $q_d$) is a martingale with respect to the filtration $\cG(q_d)$. Define a modified version similarly to the argument of  Lemma~\ref{l:lip.alg.l2.l4.bound}: let
	\begin{align}\nonumber
	\tilde{\bg}^a(q_d)
	&\equiv
	\bar{\bg}^a(q_d)
	\min\bigg\{1,\frac{2N^{1/2}}{\|\bar{\bg}^a(q_d)\|}\bigg\}\,,\\ \nonumber
	\tilde{\by}(q_{d+1})
	&\equiv
	\by(q_{d+1})
	\min\bigg\{1,\frac{2N^{1/2}}{\| \by(q_{d+1})\|}\bigg\}\,, \\
	\frac{\Delta \tilde{X}^{\RomII}(q_d,a)}{[\delta_d p(q_d)]^{1/2}}
	&\equiv 
	\tilde{\mathfrak{y}}^{\RomII}(q_d,a)
	\equiv \frac{(\tilde{\bg}^a(q_d), \tilde{\by}(q_{d+1}))}{N^{1/2}}
	\,.\label{e:lipschitz.truncation.trick.II}
	\end{align}
We note that we have
$\tilde{\bg}^a(q_d)=\bar{\bg}^a(q_d)$ and 
$\tilde{\by}(q_{d+1})=\by(q_{d+1})$ on the event
	\[\cE\equiv
	\cE(q_d,a) 
	\equiv \bigg\{ 
	\frac{\|\bar{\bg}^a(q_d)\|}{2N^{1/2}} \le1 ,
	\frac{\| \by(q_{d+1})\|}{2N^{1/2}} \le1
	\bigg\} \,.
	\]
It follows from Lemma~\ref{l:y.subgaus} and standard estimates that the event $\cE(q_d,a)$ fails with probability $e^{-cN}$. We then argue similarly as in the proof of \eqref{e:lip.alg.l2.l4.bound.l4.step2}: the random variables $\mathfrak{y}^{\RomII}$ and $\tilde{\mathfrak{y}}^{\RomII}$ always have the same sign, with $|\tilde{\mathfrak{y}}^{\RomII}| \le |\mathfrak{y}^{\RomII}|$, and they agree on the event $\cE\equiv\cE(q_d,a)$: therefore,
	\begin{align}\nonumber
	&\bigg|\E\bigg[
	\Big(\mathfrak{y}^{\RomII}(q_d,a)
	- \tilde{\mathfrak{y}}^{\RomII}(q_d,a)
	\Big)^8
	\bigg]\bigg|
	\le \E[\mathfrak{y}^{\RomII}(q_d,a)^8;\cE^c]
	\le \E[\mathfrak{y}^{\RomII}(q_d,a)^{16}]^{1/2}
	\P(\cE^c)^{1/2}\\
	&\qquad\le
	\frac{\E[ \|\bar{\bg}^a(q_d)\|^{32}]^{1/4}
		\E[ \|\by(q_{d+1})\|^{32} ]^{1/4}}{N^4}
	\P(\cE^c)^{1/2}
	\le e^{-cN}\,.
	\label{e:RomII.tilde.discrepancy.second.mmt.bound}
	\end{align}
Essentially the same argument as for \eqref{e:RomII.tilde.discrepancy.second.mmt.bound} also gives
	\beq\label{e:RomII.by.tilde.discrepancy.second.mmt.bound}
	\E\bigg[ \Big\| \by(q_{d+1})-\tilde{\by}(q_{d+1})\Big\|^8\bigg]
	\le \E\Big[ \| \by(q_{d+1})\|^8 ; \cE^c \Big]
	\le e^{-cN}
	\eeq
Now let $\tilde{X}^{\RomII}(q_d,a)$ be the process started from zero with increments $\Delta\tilde{X}^{\RomII}(q_d,a)$. (Again, we consider this for fixed $a\in[M]$, as a process indexed by time $q_d$.) It follows from  \eqref{e:RomII.tilde.discrepancy.second.mmt.bound} above that
	\beq\label{e:martingale.X.II.to.tilde.X.II.error}
	\E\bigg[ \max_{0\le d\le \dmax}
	\Big(X^{\RomII}(q_d,a)
		-\tilde{X}^{\RomII}(q_d,a)\Big)^8\bigg]
	\le o_N(1)\,.
	\eeq
Note also that $\Delta \tilde{X}^{\RomII}(q_d,a)$ is a Lipschitz function of the gaussian input, with Lipschitz constant $L^{O(1)}(\delta_d)^{1/2}$.
Next define the re-centered increments
	\begin{align*}
	\frac{\Delta \bar{X}^{\RomII}(q_d,a)}{[\delta_d p(q_d)]^{1/2}}
	&\equiv
	\frac{\Delta \tilde{X}^{\RomII}(q_d,a)}{[\delta_d p(q_d)]^{1/2}}
	-\E\bigg[
	\frac{\Delta \tilde{X}^{\RomII}(q_d,a)}{[\delta_d p(q_d)]^{1/2}}
		\,\bigg|\, \cG(q_d)
		\bigg]\\
	&= \tilde{\mathfrak{y}}^{\RomII}(q_d,a)
	-\E\Big[\tilde{\mathfrak{y}}^{\RomII}(q_d,a)
		\,\Big|\,\cG(q_d)\Big]\,.
	\end{align*}
Then $\Delta \bar{X}^{\RomII}(q_d,a)$ is also Lipschitz in the gaussian noise, with Lipschitz constant $L^{O(1)} (\delta_d)^{1/2}$; moreover it has mean zero conditional on $\cG(q_d)$. Therefore, for fixed $a\in[M]$, if we define $\bar{X}^{\RomII}(q_d,a)$ to be the process started from zero with increments
$\Delta \bar{X}^{\RomII}(q_d,a)$, then we see that $\bar{X}^{\RomII}(q_d,a)$ is a \emph{martingale} with respect to $\cG(q_d)$. By the above discussion, the increments are subgaussian with variance proxy $L^{O(1)} \delta_d$, so by iterated expectations we can bound
	\begin{align}\nonumber
	&\E\Big[ \exp(\theta \bar{X}^{\RomII}(1,a))\Big]
	\le \E\bigg[ \exp(\theta \bar{X}^{\RomII}(q_{\dmax-1},a))
	\E\Big[
	\exp(\theta \Delta\bar{X}^{\RomII}(q_{\dmax-1},a))
	\,\Big|\, \cG(q_{\dmax-1})\Big]
	\bigg]\\
	&\qquad\le \exp\Big\{
		\theta^2 L^{O(1)}\delta_d
		\Big\}
	\E\Big[\exp(\theta \bar{X}^{\RomII}(q_{\dmax-1},a))\Big]
	\le \exp\Big\{
		\theta^2 L^{O(1)}
		\Big\}\label{e:bar.X.II.martingale.exponential.concentration}
	\end{align}
for all $\theta \in\R$, from which it follows that $\E[\bar{X}^{\RomII}(1,a)^8] \le L^{O(1)}$. Combining with Doob's $L^p$ martingale inequality gives
	\beq\label{e:martingale.bar.X.II.application.of.Doob}
	\E\bigg[ \max_{0\le d\le \dmax} \bar{X}^{\RomII}(q_d,a)^8
		\bigg]
	\le O(1)
	\E\Big[\bar{X}^{\RomII}(1,a)^8\Big]
	\le L^{O(1)}\,.\eeq
Finally, to control the discrepancy between $\tilde{X}^{\RomII}$ and $\bar{X}^{\RomII}$, let
	\[
	\hat{\mathfrak{y}}^{\RomII}(q_d,a)
	\equiv \frac{(\tilde{\bg}^a(q_d), \by(q_{d+1}))}{N^{1/2}}
	\]
and note that $\hat{\mathfrak{y}}^{\RomII}(q_d,a)$ has mean zero conditional on $\cG(q_d)$. We then bound
	\begin{align*}&\E\bigg[
	\E\Big[\tilde{\mathfrak{y}}^{\RomII}(q_d,a)
		\,\Big|\,\cG(q_d)\Big]^8\bigg]
	=\E\bigg[
	\E\Big[\tilde{\mathfrak{y}}^{\RomII}(q_d,a)
		-\hat{\mathfrak{y}}^{\RomII}(q_d,a)
		\,\Big|\,\cG(q_d)\Big]^8\bigg]\\
	&\qquad\le \E\bigg[ \Big(\tilde{\mathfrak{y}}^{\RomII}(q_d,a)
		-\hat{\mathfrak{y}}^{\RomII}(q_d,a)\Big)^8\bigg]
	\le \E\bigg[\frac{ \|\tilde{\bg}^a(q_d)\|^8}{N^2}
	\|\by-\tilde{\by}\|^8\bigg]
	\le O(1) \E\Big[ \|\by-\tilde{\by}\|^8\Big] 
	\stackrel{\eqref{e:RomII.by.tilde.discrepancy.second.mmt.bound}}
		{\le} e^{-cN}\,.
	\end{align*}
It follows by summing over increments that
	\beq\label{e:martingale.tilde.X.II.to.bar.X.II.error}
	\E\bigg[ \max_{0\le d\le \dmax}
	\Big(
	\tilde{X}^{\RomII}(q_d,a)-
	\bar{X}^{\RomII}(q_d,a)
	\Big)^8\bigg]
	\le o_N(1)\,.
	\eeq
It follows by combining the above estimates
\eqref{e:martingale.X.II.to.tilde.X.II.error}, \eqref{e:martingale.bar.X.II.application.of.Doob}, and \eqref{e:martingale.tilde.X.II.to.bar.X.II.error} that
	\[ \max_{a\in[M]}
	\E\bigg[ \max_{0\le d\le\dmax}
	X^\RomII(q_d,a)^8\bigg] \le L^{O(1)}\,,\]
which proves the claim for $X^\RomII$. We comment for later use that essentially the same argument also gives
	\beq\label{e:martingale.X.holder.estimate}
	\E\bigg[
	\max_{r\in[s,t]} 
	\bigg|
		\sum_{\ell=s}^{r-1}
	\Delta X^\RomII(q_\ell,a(q_{\ell+1})) \bigg|^8\bigg]
	\le L^{O(1)} (q_t-q_s)^4
	\eeq
The claim for $X^\RomIII$ follows by a very similar argument: recalling \eqref{e:X.decomp}, we define truncated versions of $\bmeta^a(q_{d+1})$, $\bar{\bx}(q_d)$, and $\Delta X^\RomIII (q_d,a)$, similarly to \eqref{e:lipschitz.truncation.trick.II}, and argue as above.
The claim for $X^\Ising$ follows by a similar but slightly simpler argument:  we use the truncation $\tilde{\by}(q_{d+1})$ from \eqref{e:lipschitz.truncation.trick.II}, and analogously to \eqref{e:Delta.X.Ising} we define
	\[\frac{\Delta \tilde{X}^\Ising(q_d,i)}{(\delta_d)^{1/2}}
	\equiv (\be_i,\tilde{\by}(q_{d+1}))\,,
	\]
where $\be_i$ is the $i$-th standard coordinate vector in $\R^N$. (Clearly, no truncation is needed for $\be_i$.) A simplification of the above argument then gives the claim for $X^\Ising$, which concludes the proof.
\end{proof}
\end{lem}

\begin{lem}\label{l:martingale.sup.bound.X.I.minus.mean} Assume that $L\ge1$ and that condition \eqref{e:deriv.p.bound} holds. Let $\vX\equiv \vX^N$ be the process defined by \eqref{e:X.decomp} and \eqref{e:Delta.X.Ising}. Then the process $X^\RomI$ satisfies the estimate
	\[ \max_{a\in[M]}
	\E\bigg[ \max_{0\le d\le\dmax}
	 \Big(X^\RomI(q_d,a)
	 -\E X^\RomI(q_d,a)\Big)
	 ^8\bigg] \le L^{O(1)}\,,\]
with the estimate holding uniformly in $N$.
\begin{proof} Similarly as in the proof of Lemma~\ref{l:martingale.sup.bound.X.II.III.Ising}, let us abbreviate
	\[
	\frac{\Delta X^{\RomI}(q_d,a)}{\delta_d p'(q_d)^{1/2}}
	\equiv \mathfrak{y}^{\RomI}(q_d,a)
	\equiv \frac{(\bmeta^a(q_{d+1}),\by(q_{d+1}))}{N^{1/2}}\,.
	\]
Define truncated $\tilde{\bmeta}$ and $\tilde{\by}$ similarly as
in the proof of Lemma~\ref{l:martingale.sup.bound.X.II.III.Ising}, and define
	\[
	\frac{\Delta \tilde{X}^{\RomI}(q_d,a)}{\delta_d p'(q_d)^{1/2}}
	\equiv \tilde{\mathfrak{y}}^{\RomI}(q_d,a)
	\equiv \frac{(\tilde{\bmeta}^a(q_{d+1}),
		\tilde{\by}(q_{d+1}))}{N^{1/2}}\,.
	\]
By the same argument as for \eqref{e:lip.alg.l2.l4.bound.l4.step2} or for \eqref{e:RomII.tilde.discrepancy.second.mmt.bound}, we have
	\beq\label{e:lip.alg.l2.l4.bound.l4.step2.RomI}
	\bigg|\E\bigg[
	\Big(
	\mathfrak{y}^{\RomI}(q_d,a)
	- \tilde{\mathfrak{y}}^{\RomI}(q_d,a)
	\Big)^8
	\bigg]\bigg|
	\le e^{-cN}\,.
	\eeq
Let $\tilde{X}^{\RomI}(q_d,a)$ be the process started from zero with increments
$\Delta \tilde{X}^{\RomI}(q_d,a)$. (As before, we regard this as a process indexed by time $q_d$, for fixed $a\in[M]$.) Applying \eqref{e:lip.alg.l2.l4.bound.l4.step2.RomI} gives
	\beq\label{e:martingale.X.I.to.tilde.X.I.error}
	\E\bigg[ \max_{0\le d \le \dmax}\Big(
	X^{\RomI}(q_d,a)- \tilde{X}^{\RomI}(q_d,a)\Big)^8\bigg]
	\le o_N(1)\,.
	\eeq
Then consider the centered increments
	\begin{align*}
	\frac{\Delta \bar{X}^{\RomI}(q_d,a)}{\delta_d p'(q_d)^{1/2}}
	&\equiv
	\frac{\Delta \tilde{X}^{\RomI}(q_d,a)}{\delta_d p'(q_d)^{1/2}}
	-\E\bigg[
	\frac{\Delta \tilde{X}^{\RomI}(q_d,a)}{\delta_d p'(q_d)^{1/2}}
		\,\bigg|\, \cG(q_d)
		\bigg]\\
	&= \tilde{\mathfrak{y}}^{\RomI}(q_d,a)
	- \E\Big[\tilde{\mathfrak{y}}^{\RomI}(q_d,a)
		\,\Big|\,\cG(q_d)\Big]\,.
	\end{align*}
Each increment above is a Lipschitz function of the gaussian noise, with Lipschitz constant $L^{O(1)}\delta_d$, and has mean zero conditional on $\cG(q_d)$.
Let $\bar{X}^{\RomI}(q_d,a)$ be the process started from zero with increments
$\Delta \bar{X}^{\RomI}(q_d,a)$. (Again, we regard this as a process indexed by time $q_d$, for fixed $a\in[M]$.)
Then $\bar{X}^{\RomI}(q_d,a)$ is a \emph{martingale}, whose increments are subgaussian with variance proxy $L^{O(1)} (\delta_d)^2$. Similar arguments as in the proof of Lemma~\ref{l:martingale.sup.bound.X.II.III.Ising} give (cf.\ \eqref{e:bar.X.II.martingale.exponential.concentration})
	\[
	\E \Big[\exp\{ \theta \bar{X}^{\RomI}(1,a)\}\Big]
	\le \exp( \theta^2 L^{O(1)} \delta)\,,\]
and combining with Doob's $L^p$ martingale inequality gives (cf.\ \eqref{e:martingale.bar.X.II.application.of.Doob})
	\beq\label{e:martingale.bar.X.I.application.of.Doob}
	\E\Big[ \max_{0\le d\le \dmax}
	\bar{X}^{\RomI}(q_d,a)^8\Big] 
	\le L^{O(1)} \delta^4\,.
	\eeq
It remains to address the discrepancy between $X^{\RomI}$ and $\bar{X}^{\RomI}$. Since $\tilde{X}^{\RomI}$ has subgaussian increments, we have
	\[
	\E\bigg[
	\Big( \Delta \tilde{X}^{\RomI}(q_d,a) - \E \Delta \tilde{X}^{\RomI}(q_d,a)
	\Big)^8\bigg] \le  L^{O(1)}\delta^8\,,
	\]
and summing over increments gives
	\beq\label{e:martingale.tilde.X.I.close.to.mean}
	\E\bigg[ \max_{0\le d\le \dmax}
	\Big(
	 \tilde{X}^{\RomI}(q_d,a) - \E \tilde{X}^{\RomI}(q_d,a)\Big)^8\bigg]
	 \le L^{O(1)}\,.
	\eeq
It also follows from  
\eqref{e:lip.alg.l2.l4.bound.l4.step2.RomI}
or \eqref{e:martingale.X.I.to.tilde.X.I.error}
that
 \beq\label{e:martingale.tilde.X.I.close.to.X.I}
 \max_{0\le d\le \dmax}
 \bigg|\E \Big[ \tilde{X}^{\RomI}(q_d,a)
 - X^{\RomI}(q_d,a) \Big] \bigg| \le o_N(1)
 \eeq
 where we recall $\E X^{\RomI}(q_d,a)
 = \E X(q_d,a)$.  Combining \eqref{e:martingale.X.I.to.tilde.X.I.error}, \eqref{e:martingale.tilde.X.I.close.to.mean}, and \eqref{e:martingale.tilde.X.I.close.to.X.I} gives the conclusion.
\end{proof}
\end{lem}

Recall from Definition~\ref{d:trunc} the truncated process $\vX^\trstar$. Consider the random time
	\beq\label{e:trstar.stopping.time}
	\tau\equiv \tau(a)
	= \min\bigg\{1,
	\inf\Big\{q_d : \vX(q_d,a) \notin [-\trK,\trK)^3\Big\}
	\bigg\}\,.\eeq
Then the processes $\vX(q_d,a)$ and $\vX^\trstar(q_d,a)$ (for fixed $a\in[M]$, indexed by time $q_d$) agree up to time $\tau(a)$.
Note that if we fix any index $a\in[M]$, then $\tau(a)$ is a stopping time with respect to the filtration $\cG(t)$. Note that if for $\sigma\in\{\RomI,\RomII,\RomIII\}$ we define the stopping times
	\beq\label{e:trstar.stopping.time.I.II.III}
	\tau^\sigma
	\equiv \tau^\sigma(a)
	= \min\bigg\{1,
	\inf\Big\{q_d :
	X^\sigma(q_d,a) \notin [-\trK,\trK)\Big\}
	\bigg\}\,,
	\eeq
then we can also express $\tau$ as the minimum of $\tau^{\RomI}$, $\tau^{\RomII}$, $\tau^{\RomIII}$. We analogously define
	\beq\label{e:trstar.stopping.time.ISING}
	\tau^\Ising
	\equiv \tau^\Ising(i)
	= \min\bigg\{1,
	\inf\Big\{q_d :
	X^\Ising(q_d,i) \notin [-\trK,\trK)\Big\}
	\bigg\}\,,
	\eeq
The following controls the maximum of the $\vX$ process on the events
$\tau<1$ and $\tau^\Ising<1$.

\begin{ppn}\label{p:max.X.process.on.trstar.stopped.event}
Assume that $L\ge1$ and that condition \eqref{e:deriv.p.bound} holds. Let $\vX\equiv \vX^N$ be the process defined by \eqref{e:X.decomp}  and \eqref{e:Delta.X.Ising}, and assume that condition~\eqref{e:w2.coupling.assumption} holds. Then, for the stopping times defined by \eqref{e:trstar.stopping.time}--\eqref{e:trstar.stopping.time.ISING}, we have
	\begin{align*}
	\max_{\sigma\in\{\RomI,\RomII,\RomIII\}}
	\frac1M \sum_{a=1}^M
	\E\bigg[
	\max_{0\le d\le \dmax}
	X^\sigma(q_d,a)^2; \tau(a)<1
	\bigg]
	&\le O(\epsilon) + o_{\trK}(1)\,,\\
	\frac1N \sum_{i=1}^N
	\E\bigg[
	\max_{0\le d\le \dmax}
	X^\Ising(q_d,i)^2; 
	\tau^\Ising(i)<1
	\bigg]
	&\le o_{\trK}(1)
	\end{align*}
uniformly over $N$. 

\begin{proof} 
Throughout this proof, we write $\bE_a$ for expectation over a uniformly random index $a\in[M]$, we write $\E$ for expectation over the gaussian disorder, and we write $\E'\equiv \bE_a\E$ for expectation over both. The assumption \eqref{e:w2.coupling.assumption} says that there exists a coupling of $X(1)$ with $Z\sim\mu$ such that
	\[
	\E'\Big[(X(1)-Z)^2\Big]
	\le \epsilon\,.
	\]
We can use this to bound
	\begin{align}\nonumber
	&\frac1M  \sum_{a=1}^M (\E X^\RomI(1,a))^2
	=\frac1M  \sum_{a=1}^M (\E X(1,a))^2
	\le
	\frac1M  \sum_{a=1}^M \E[ X(1,a)^2 ]
	=\E'[ X(1)^2 ] \\
	&\qquad\le 
	2 \E'\Big[ (X(1)-Z)^2  + Z^2\Big]
	\le 2\Big( \epsilon + \E'(Z^2)\Big)
	\le O(1)\,,
	\label{e:X.I.squared.expectation.avg.bound}
	\end{align}
where the $O(1)$ indicates a large constant depending only on $\mu$. 
Recall that the $X^{\RomII}$ and $X^{\RomIII}$ processes have mean zero, so $\E X^{\RomI}=\E X$. Next, a straightforward gaussian computation gives the identity
	\[
	\E X^{\RomI}(q_d,a) = p(q_d) \E X^{\RomI}(1,a)\,.
	\]
It follows by combining with the result of Lemma~\ref{l:martingale.sup.bound.X.I.minus.mean} that
	\begin{align*}
	&\frac{1}{M}\sum_{a=1}^M
	\E \bigg[\max_{0\le d\le\dmax} X^{\RomI}(q_d,a)^2\bigg]\\
	&\qquad\le\frac{2}{M}\sum_{a=1}^M
		\E \bigg[\max_{0\le d\le\dmax}
			\Big( X^{\RomI}(q_d,a) - \E X^{\RomI}(q_d,a) \Big)^2
			\bigg]
	+ \frac{2}{M}\sum_{a=1}^M \max_{0\le d\le\dmax} 
		(\E X^{\RomI}(q_d,a))^2 \\
	&\qquad
	\le L^{O(1)}
	+ \frac{1}{M}\sum_{a=1}^M 
		(\E X^{\RomI}(1,a))^2
	\le L^{O(1)} +O(1) \le L^{O(1)}\,.
	\end{align*}
As a consequence (recalling $L\ll \trK$) we conclude
	\[\frac1M\sum_{a=1}^M \P(\tau^{\RomI}(a) < 1)
	\le
	\frac{1}{M\trK^2} \sum_{a=1}^M 
	\E\bigg[ \max_{0\le d\le \dmax} X^{\RomI}(q_d,a)^2
		\bigg] \le o_{\trK}(1)\,,
	\]
uniformly in $N$. We also have from Lemma~\ref{l:martingale.sup.bound.X.II.III.Ising} that 
	\beq\label{e:truncation.K.stopping.time.small.prob.ISING}
	\bigg\{
	\max_{a\in[M]}
	\max_{\sigma\in\{\RomII,\RomIII\}}
	\P(\tau^\sigma(a) < 1)\bigg\}
	+ \bigg\{
	\max_{i\in[N]}
	\P(\tau^\Ising(i) < 1) \bigg\}
	\le o_{\trK}(1)\,,\eeq
so altogether we obtain
	\beq\label{e:truncation.K.stopping.time.small.prob}
	\frac1M\sum_{a=1}^M \P(\tau(a) < 1) 
	\le o_{\trK}(1)\,,
	\eeq
uniformly in $N$. We can also estimate
	\begin{align*}
	(*)&\equiv \frac1M
	\sum_{a=1}^M
	\E\bigg[
	\max_{0\le d\le \dmax}  X^{\RomI}(q_d,a)^2; \tau(a)<1
	\bigg]
	\le 2\textup{(A)} + 2\textup{(B)}\,,\\
	\textup{(A)}
	&\equiv \frac1M \sum_{a=1}^M
	\E\bigg[
	\max_{0\le d\le \dmax}
	\Big( X^{\RomI}(q_d,a) - \E  X^{\RomI}(q_d,a) \Big)^2;
		\tau(a)<1
	\bigg]\,,\\
	\textup{(B)}
	&\equiv  \frac1M \sum_{a=1}^M
	\max_{0\le d\le \dmax}
	(\E  X^{\RomI}(q_d,a))^2 \P(\tau(a)<1)
	= \frac1M \sum_{a=1}^M
	(\E  X^{\RomI}(1,a))^2 \P(\tau(a)<1)\,.
	\end{align*}
We then use Lemma~\ref{l:martingale.sup.bound.X.I.minus.mean} to bound
	\begin{align*}\textup{(A)}
	&\le\frac1M
	\sum_{a=1}^M
	\E\bigg[\max_{0\le d\le \dmax}
	\Big( X^{\RomI}(q_d,a) - \E  X^{\RomI}(q_d,a) \Big)^4\bigg]^{1/2}
		\P(\tau(a)<1)^{1/2} \\
	&\le L^{O(1)}
	\bE_a [\P(\tau(a)<1)^{1/2}]
	\le L^{O(1)} \Big( \bE_a\P(\tau(a)<1)\Big)^{1/2}
	\stackrel{\eqref{e:truncation.K.stopping.time.small.prob}}{\le}
	o_{\trK}(1)\,.
	\end{align*}
Recall again the coupling of $X(1)$ with $Z\sim\mu$. Under this coupling, let $Z(a)$ denote the random variable $Z$ conditioned on the index $a\in[M]$. We then bound
	\begin{align*}\textup{(B)}
	&= \bE_a\bigg[
	\E[  X^{\RomI}(1,a)]^2 \P(\tau(a)<1)\bigg]
	\le 2\textup{(B1)} + 2\textup{(B2)}\,,\\
	\textup{(B1)}
	&\equiv \bE_a\bigg[
	\Big(\E[X^{\RomI}(1,a) - Z(a)]\Big)^2 \P(\tau(a)<1)\bigg] \\
	&= \bE_a\bigg[
	\Big(\E[X(1,a) - Z(a)]\Big)^2 \P(\tau(a)<1)\bigg]
	\le \E'\bigg[\Big(X(1) - Z \Big)^2\bigg] \le \epsilon\,,\\
	\textup{(B2)}
	&\equiv 
	\bE_a\bigg[
	\E[ Z(a)^2] \P(\tau(a)<1)\bigg]
	\equiv
	\bE_a\bigg[
	\E[ Z(a)^2] f(a) \bigg]
	= \E'(Z^2 f)\,.
	\end{align*}
Here $f(a)\equiv \P(\tau(a)<1)$. We know that $Z\sim\mu$ is a single random variable (i.e., not depending on $N$) with $\E'(Z^2) \le O(1)$ (a finite constant depending only on $\mu$). Meanwhile, the random variable $f$ depends on $N$, but $\E'f \le o_{\trK}(1)$ uniformly in $N$ by \eqref{e:truncation.K.stopping.time.small.prob}. Since $0\le f\le 1$ almost surely, it follows by the dominated convergence theorem that $\textup{(B2)} \le \E'(Z^2 f) \le o_{\trK}(1)$.

Combining the above estimates gives
	\[(*)=\frac1M \sum_{a=1}^M
	\E\bigg[
	\max_{0\le d\le \dmax}  X^{\RomI}(q_d,a)^2; \tau(a)<1
	\bigg]
	\le 4\epsilon + o_{\trK}(1)\,.
	\]
For $\sigma\in\{\RomII,\RomIII\}$, applying Lemma~\ref{l:martingale.sup.bound.X.II.III.Ising} gives the simpler estimates
	\begin{align*}&\frac1M\sum_{a=1}^M
	\E\bigg[
	\max_{0\le d\le \dmax}  X^{\sigma}(q_d,a)^2; \tau(a)<1
	\bigg]\\
	&\qquad\le\frac1M\sum_{a=1}^M
	 \E\bigg[
	\max_{0\le d\le \dmax}
	 X^{\sigma}(q_d,a)^4
	\bigg]^{1/2} \P(\tau(a) <1)^{1/2} 
	\le L^{O(1)}
	\frac1M\sum_{a=1}^M  \P(\tau(a) <1)^{1/2} \\
	&\qquad\le
	L^{O(1)}\Big(
	\bE_a \P(\tau(a) <1)\Big)^{1/2}
	\stackrel{\eqref{e:truncation.K.stopping.time.small.prob}}{\le}
	 o_{\trK}(1)\,,
	\end{align*}
as well as the analogous estimate for $X^\Ising$ (using \eqref{e:truncation.K.stopping.time.small.prob.ISING} in place of \eqref{e:truncation.K.stopping.time.small.prob}). This concludes the proof.
\end{proof}
\end{ppn}

\begin{cor}\label{c:X.versus.trstar.X.W2.error.estimate}
Assume that $L\ge1$ and that condition \eqref{e:deriv.p.bound} holds. Let $\vX\equiv \vX^N$ be the process defined by \eqref{e:X.decomp}
 and \eqref{e:Delta.X.Ising}, and assume that condition~\eqref{e:w2.coupling.assumption} holds. Then the process $\vX^\trstar$ of Definition~\ref{d:trunc} satisfies
	\[
	\E'\Big[\|\vX(1)-\vX^\trstar(1)\|^2\Big]
	\le O(\epsilon) + o_{\trK}(1)\,,
	\]
uniformly in $N$.

\begin{proof} 
For $\sigma\in\{\RomI,\RomII,\RomIII\}$ we have $X^{\trstar,\sigma}(1,a)=X^\sigma(\tau(a),a)$, with $\tau$ the stopping time defined by \eqref{e:trstar.stopping.time}. It follows that, for $\sigma\in\{\RomI,\RomII,\RomIII\}$, we have
	\[
	\E'\Big[( X^\sigma(1)-X^{\trstar,\sigma}(1))^2\Big]
	\le O(1)
	\E'\bigg[
	\max_{0\le d\le \dmax}
	X^\sigma(q_d)^2;\tau<1\bigg] 
	\le O(\epsilon) +o_{\trK}(1)\,,
	\]
where the last bound is by Proposition~\ref{p:max.X.process.on.trstar.stopped.event}. Similarly,
	\[
	\E'\Big[( X^\Ising(1)-X^{\trstar,\Ising}(1))^2\Big]
	\le O(1)
	\E'\bigg[
	\max_{0\le d\le \dmax}
	X^\Ising(q_d)^2;\tau^\Ising<1\bigg] 
	\le o_{\trK}(1)\,,
	\]
again by Proposition~\ref{p:max.X.process.on.trstar.stopped.event}. The conclusion follows.
\end{proof}
\end{cor}

Recall the process $\vX^\frstar(q_d,a)$ from Definition~\ref{d:froze}, where we added a small positive density of frozen particles. We denote
$\vX^\frstar(q_d)\equiv\vX^\frstar(q_d,a^\frozen(q_d))$, where $a^\frozen(q_d)$ is a uniformly random index from $[M^\frozen]$. The following lemma says that the addition of frozen particles incurs a negligible $\bbW_2$ error in the endpoint distribution of the process:

\begin{lem}\label{l:W2.error.from.adding.frozen.particles}
Assume that $L\ge1$ and that condition \eqref{e:deriv.p.bound} holds. Let $\vX\equiv \vX^N$ be the process defined by \eqref{e:X.decomp}
and \eqref{e:Delta.X.Ising}, and assume that condition~\eqref{e:w2.coupling.assumption} holds. Then the processes $\vX^\trstar$ and $\vX^\frstar$ of Definitions~\ref{d:trunc} and \ref{d:froze} can be coupled so that
	\[
	\E\Big[\|\vX^\trstar(1)
		-\vX^\frstar(1)\|^2\Big]
	\le o_{\trK}(1)\,,
	\]
uniformly in $N$.

\begin{proof}
Consider $\sigma\in\{\RomI,\RomII,\RomIII\}$. Abbreviate 
	\[\epsilon(\trK)\equiv 
		\frac{1/\trK^6}
		{1+1/\trK^6}\,.\]
If we let $U$ denote a random variable sampled uniformly at random from the interval $[0,1]$, then 
$X^{\frstar,\sigma}(1)$ is equidistributed as
	\[
	X^{\frstar,\sigma}(1,U)
	\equiv
	\begin{cases} \displaystyle
	X^{\trstar,\sigma}\bigg(1,
		\bigg\lceil
		\frac{M\cdot U}{1-\epsilon(\trK)}
		\bigg\rceil\bigg)
	& \textup{if }U\in[0,1-\epsilon(\trK)]\,,\\
	T(U) & 
		\textup{if }U\in(1-\epsilon(\trK),1]\,,
	\end{cases}
	\]	
where $T$ is a random variable confined to $[-2\trK,2\trK]$. On the other hand, $X^{\trstar,\sigma}(1)$ is equidistributed as
	\[
	X^{\trstar,\sigma}(1,U)
	\equiv 
	\begin{cases} \displaystyle
	X^{\trstar,\sigma}\bigg(1,
		\bigg\lceil
		\frac{M\cdot U}{1-\epsilon(\trK)}
		\bigg\rceil\bigg)
	& \textup{if }U\in[0,1-\epsilon(\trK)]\,,\\ \displaystyle
	X^{\trstar,\sigma}\bigg(1,
		\bigg\lceil
		\frac{M\cdot 
			[U-(1-\epsilon(\trK))]}{\epsilon(\trK)}
		\bigg\rceil\bigg)
	&\textup{if } U\in(1-\epsilon(\trK),1]
	\end{cases}\,.
	\]
Under this coupling, we have
	\[\E\Big[(
	X^{\frstar,\sigma}(1,U)
	-X^{\trstar,\sigma}(1,U))^2
	\Big]
	\le O(1) \epsilon(\trK)
	\E\Big[
	T(U)^2 + X^{\trstar,\sigma}(1)^2\Big]\,.
	\]
Since $T(U)$ is confined to $[-2\trK,2\trK]$, we clearly have $\E[T(U)^2] \le 4\trK^2$. It follows from Lemma~\ref{l:martingale.sup.bound.X.II.III.Ising} that $X^\sigma(1)$ has bounded second moment for $\sigma\in\{\RomII,\RomIII,\Ising\}$. Under condition~\eqref{e:w2.coupling.assumption}, it follows from Lemma~\ref{l:martingale.sup.bound.X.I.minus.mean} combined with \eqref{e:X.I.squared.expectation.avg.bound} that  $X^\RomI(1)$ also has bounded second moment. Combining with
Corollary~\ref{c:X.versus.trstar.X.W2.error.estimate} gives that $\vX^\trstar(1)$ has bounded second moment. Substituting into the above bound, and recalling the definition of $\epsilon(\trK)$, we conclude
	\[\E\Big[(
	X^{\frstar,\sigma}(1,U)
	-X^{\trstar,\sigma}(1,U))^2
	\Big]
	\le O(1) 
	\frac{\trK^2 + O(1)}{\trK^6}
	\le o_{\trK}(1)\,,
	\]
which proves the claim for $\sigma\in\{\RomI,\RomII,\RomIII\}$.
A similar argument proves the analogous claim for $X^\Ising$.
\end{proof}
\end{lem}

\begin{lem}\label{l:W2.error.from.freezing.small.buckets}
Assume that $L\ge1$ and that condition \eqref{e:deriv.p.bound} holds. Let $\vX\equiv \vX^N$ be the process defined by \eqref{e:X.decomp}
and \eqref{e:Delta.X.Ising}, and assume that condition~\eqref{e:w2.coupling.assumption} holds. Then the processes $\vX^\trstar,\vX^\frstar$ 
and $\vX^\trunc,\vX^\frozen$ 
of Definitions \ref{d:trunc}--\ref{d:buckets.Ising} can be coupled so that
	\begin{align*}
	\E\Big[\|\vX^\trstar(1)
		-\vX^\trunc(1)\|^2\Big]
	&\le O(\epsilon) + o_{\trK}(1)\,,\\
	\E\Big[\|\vX^\frstar(1)
		-\vX^\frozen(1)\|^2\Big]
	&\le O(\epsilon) + o_{\trK}(1)\,,
	\end{align*}
uniformly in $N$.
\begin{proof}
For fixed $a\in[M]$, consider $\vX^\trstar(q_d,a)$ as a process indexed by time $q_d$. Recalling Definitions~\ref{d:buckets} and \ref{d:buckets.Ising}, let $\tau^\trunc$ be the first time that index $a$ lands in a small bucket, and analogously $\tau^{\trunc,\Ising}$:
	\begin{align*}
	\tau^\trunc(a)
	&\equiv \min\bigg\{ 1,
	\min\bigg\{
		q_d : |B(q_d,a)| < \frac{M}{\BUCKETS^2}\bigg\} 
	\bigg\}\,,\\
	\tau^{\trunc,\Ising}(i)
	&\equiv \min\bigg\{ 1,
	\min\bigg\{
		q_d : |B^\Ising(q_d,i)| < \frac{N}{\BUCKETS^2}\bigg\} 
	\bigg\}\,.
	\end{align*}
It follows from the discussion around \eqref{e:total.number.of.particles.in.small.buckets} that
	\[\frac1M \sum_{a=1}^M
	\ind\{ \tau^\trunc(a) < 1 \}
	+\frac1N \sum_{i=1}^N
	\ind\{ \tau^{\trunc,\Ising}(i) < 1 \}
	\le o_\eta(1)\,.
	\]
This is similar to but stronger than our earlier claim 
\eqref{e:truncation.K.stopping.time.small.prob} for the stopping time $\tau$ from \eqref{e:trstar.stopping.time}. Therefore, repeating the argument around \eqref{e:truncation.K.stopping.time.small.prob} gives the result of
Proposition~\ref{p:max.X.process.on.trstar.stopped.event} with $\tau^\trunc$ in place of $\tau$. Repeating the argument of Corollary~\ref{c:X.versus.trstar.X.W2.error.estimate} then gives the conclusion for the $\bbW_2$ error between $\vX^\trstar(1)$ and $\vX^\trunc(1)$. Essentially the same argument applies to bound the 
$\bbW_2$ error between $\vX^\frstar(1)$ and $\vX^\frozen(1)$.
\end{proof}
\end{lem}

\begin{proof}[\hypertarget{proof:p.summarize.W2.errors.preprocessing}{Proof of Proposition~\ref{p:summarize.W2.errors.preprocessing}}] Follows by combining Corollary~\ref{c:X.versus.trstar.X.W2.error.estimate}, Lemma~\ref{l:W2.error.from.adding.frozen.particles}, and Lemma~\ref{l:W2.error.from.freezing.small.buckets}.
\end{proof}

\subsection{Moment estimates for process}
\label{ss:kolmogorov.Y}

Recall from \eqref{e:X.decomp}--\eqref{e:vX} the decomposition $\vX=(\vX^{\RomI:\RomIII},X^\Ising)$. 
We have also the variants $\vX^\trunc$, $\vX^\frozen$ from Definitions \ref{d:buckets} and \ref{d:buckets.Ising}, and their spatial rerandomizations $\vY$, $\vY^\frozen$ from Definitions~\ref{d:rerand} and \ref{d:rerand.Ising}. \textbf{The main result of this subsection is Proposition~\ref{p:Y.kolmogorov}, \hyperlink{proof:p.Y.kolmogorov}{whose proof} is given at the end of this subsection.
We also prove a quenched estimate, 
Proposition~\ref{p:X.second.mmt.quenched}, which will be used later.} It follows from a few intermediate results which we turn to next.
We start by comparing
increments of $\vX^\trunc$ with increments of the spatial rerandomization $\vY$:

\begin{lem}[errors from spatial rerandomization]
\label{l:X.Y.discrep}
Let $\vY$ be the spatial rerandomization of $\vX^\trunc$, defined by \eqref{e:def.Y}
with indices $(a(q_0), \ldots, a(\dmax))$ sampled according to Definition~\ref{d:rerand}. Then we have
	\begin{align*} 
	\max\bigg\{
	\Big|\Delta Y^{\sigma}(q_d)
		-\Delta X^{\trunc,\sigma}(q_d, a(q_{d+1}))\Big|
		: \sigma\in\{
		\RomI,\RomII,\RomIII
		\}
		\bigg\}
		&\le 2\eta\,,\\
	\bigg|
	\Delta Y^{\Ising}(q_d)
	-\Delta X^{\trunc,\Ising}(q_d,i(q_{d+1}))
	\bigg|
	&\le2\eta\,.
	\end{align*}
for all $0\le d\le \dmax-1$. The same holds with $(\vY^\frozen,\vX^\frozen)$ in place of $(\vY,\vX^\trunc)$.

\begin{proof}
We prove the estimate for
$Y^{\RomI}$ versus $X^{\trunc,\RomI}$; the other bounds follow by the same argument. It follows from the definition~\eqref{e:def.Y} that
	\[
	\Delta Y^{\RomI}(q_d)
	= \sum_{a\in B(q_{d+1},a(q_{d+1}))}
	\frac{X^{\trunc,\RomI}(q_{d+1},a)}{|B(q_{d+1},a(q_{d+1}))|}
	-\sum_{a\in B(q_d,a(q_d))}
	\frac{X^{\trunc,\RomI}(q_d,a)}{|B(q_d,a(q_d))|}\,.
	\]
Marginally, $a(q_d)$ is a uniformly random index from $[M]$. Conditional on $\cG(q_d)$, the index $a(q_{d+1})$ is sampled uniformly at random from  $B=B(q_d,a(q_d))$, the bucket at time $q_d$ that contains $a(q_d)$. 
It follows from the definitions of the buckets (see~\eqref{e:bucket}) that
	\begin{align*}
	\max\bigg\{ \Big|X^{\trunc,\RomI}(q_{d+1},a)
	-X^{\trunc,\RomI}(q_{d+1},a(q_{d+1}))\Big|
	: a\in B(q_{d+1},a(q_{d+1}))
	\bigg\} &\le \eta\,,\\
	\max\bigg\{ \Big|X^{\trunc,\RomI}(q_d,a)
	-X^{\trunc,\RomI}(q_d,a(q_{d+1}))\Big|
	: a\in B(q_d,a(q_{d+1}))
	\bigg\} &\le \eta\,,
	\end{align*}
where the second bound uses that
$a(q_{d+1})$ and $a(q_d)$ are in the same bucket $B$ at time $q_d$.
Since $\vY(q_d)$ is obtained by averaging $\vX^\trunc(q_d,a)$ over $a\in B$, we immediately obtain
	\[
	\max\bigg\{
	\Big|Y^{\RomI}(q_d)
	- X^{\trunc,\RomI}(q_d,a(q_d))\Big|,
	\Big|Y^{\RomI}(q_d)
	- X^{\trunc,\RomI}(q_d,a(q_{d+1}))\Big|
	\bigg\} \le\eta\,.
	\]
Likewise, since $\vY(q_{d+1})$ is obtained by averaging 
$\vX^\trunc(q_{d+1},a')$ over $a'\in B'=B(q_{d+1},a(q_{d+1}))$, we have
	\beq\label{e:Y.next}
	\Big|Y^{\RomI}(q_{d+1})
	- X^{\trunc,\RomI}(q_{d+1},a(q_{d+1}))\Big| \le\eta\,.
	\eeq
Combining these gives the claimed bound.
\end{proof}
\end{lem}

The following quenched estimate will be used numerous times in the rest of this section: 

\begin{ppn}[quenched second moment bounds for $\vX$, $\vX^\trunc$, and $\vX^\frozen$]
\label{p:X.second.mmt.quenched}
Let $\hat{a}$, $\hat{a}^\frozen$,
$\hat{\iota}$, and $\hat{\iota}^\frozen$
be sampled uniformly at random from 
$[M]$, $[M^\frozen]$,
$[N]$, and $[N^\frozen]$, independently of the gaussian disorder and of one another. 
Recall that $\bE$ denotes expectation over $\hat{a}$, $\hat{a}^\frozen$, $\hat{\iota}$, and $\hat{\iota}^\frozen$, as applicable.
Then for 
all $0\le d\le \dmax-1$ we have
	\[\max\bigg\{
	\frac{\bE[\Delta X^{\RomI}(q_d,\hat{a})^2]}
		{(\delta_d)^2 p'(q_d)},
	\frac{\bE[\Delta X^{\RomII}(q_d,\hat{a})^2]}
		{\delta_d  p(q_d)},
	\frac{\bE[\Delta X^{\RomIII}(q_d,\hat{a})^2]}
		{\delta_d  q_d p'(q_d)},
	\frac{\bE[\Delta X^{\Ising}(q_d,{\hat{\iota}})^2]}
		{\delta_d}
	\bigg\}
    \le (\MAX_N)^2\,,\]
where $\MAX_N$ is defined by \eqref{e:MAX.N}. The same bound holds with $\vX^\trunc$ or $\vX^\frozen$ in place of $\vX$
(using $(\hat{a}^\frozen,\hat{\iota}^\frozen)$ in place of $(\hat{a},\hat{\iota})$ for the case of $\vX^\frozen$).

\begin{proof}
Recall the notation introduced in \eqref{e:normalized.vectors.g}--\eqref{e:normalized.vectors.y}. From \eqref{e:X.decomp}, we have
	\begin{align*}
	\frac{\Delta X^{\RomI}(q_d,a)}{\delta_dp'(q_d)^{1/2}}
	&= 
	\frac{(\bmeta^a(q_{d+1}),
	\by(q_{d+1}))}{N^{1/2}}\,,\\
	\frac{\Delta X^{\RomII}(q_d,a)}{(\delta_d p_d)^{1/2}}
	&= \frac{(\bar{\bg}^a(q_d),\by(q_{d+1}))}{N^{1/2}}\,,\\
	\frac{\Delta X^{\RomIII}(q_d,a)}{[\delta_d q_d p'(q_d)]^{1/2}}
	&= 
	\frac{(\bmeta^a(q_{d+1}),\bar{\bx}(q_d))}{N^{1/2}}\,,\\
	\frac{\Delta X^{\Ising}(q_d,
		{i})}{(\delta_d)^{1/2}}
	&=(\be_i,\by(q_{d+1}))\,,
	\end{align*}
where $a\in[M]$ and $i\in[N]$.
We therefore have the (quenched) second moment bounds
	\begin{align*}
	\frac{\bE[\Delta X^{\RomI}(q_d,
		\hat{a})^2]}
		{(\delta_d)^2 p'(q_d)}
	&=\frac{1}{M}
	\sum_{a=1}^M \frac{(
	\bmeta^a(q_{d+1}),
		\by(q_{d+1}))^2}{N}
	\le \frac{
	\|\bXi(q_{d+1})\|^2
	\|\by(q_{d+1})\|^2}
		{MN}
	\le (\MAX_N)^2\,,\\
	\frac{\bE[\Delta X^{\RomII}(q_d,\hat{a})^2]}
		{\delta_d p_d}
	&=
	\frac{1}{M}
	\sum_{a=1}^M 
	\frac{(\bar{\bg}^a(q_d),\by(q_{d+1}))^2}{N}
	\le
	\frac{\|\bar{\bG}(q_d)\|^2
	\|\by(q_{d+1})\|^2}{MN}
	\le (\MAX_N)^2\,,\\
	\frac{\bE[\Delta X^{\RomIII}(q_d,
		\hat{a})^2]}{\delta_d q_d p'(q_d)}
	&=\frac{1}{M}
	\sum_{a=1}^M
	\frac{(\bmeta^a(q_{d+1}),\bar{\bx}(q_d))^2}{N}
	\le \frac{\|\bXi(q_{d+1})\|^2 \|\bar{\bx}(q_d)\|^2}{MN}
	\le (\MAX_N)^2\,,\\
	\frac{\bE[\Delta X^{\Ising}
	(q_d,{\hat{\iota}})^2]}{\delta_d}
	&=
	\sum_{i=1}^N \frac{(\be_i,\by(q_{d+1}))^2}{N}
	= \frac{\|\by(q_{d+1})\|^2}{N}
	\le \MAX_N\,.
	\end{align*}
This proves the claim for $\vX$.
Then note that $\Delta X^{\trunc,\RomI}(q_d,a)$ equals either $\Delta X^{\RomI}(q_d,a)$ or zero, where the latter case corresponds to scenarios where the trajectory has exited $[-\trK,\trK)^3$ (Definition~\ref{d:trunc}) or has landed in a small bucket (Definition~\ref{d:buckets}). 
Likewise, $\Delta X^{\frozen,\RomI}(q_d,a)$ equals either $\Delta X^{\trunc,\RomI}(q_d,a)$ or zero, 
where the latter case corresponds to the same scenarios mentioned above for $X^{\trunc,\RomI}$, or to the scenario where $a$ is the index of one of the frozen particles added in Definition~\ref{d:froze}. It follows that
	\[
	\bE\Big[\Delta X^{\frozen,\RomI}(q_d,\hat{a}^\frozen)^2\Big]
	\le 
	\bE\Big[\Delta X^{\trunc,\RomI}(q_d,\hat{a})^2\Big]
	\le \bE\Big[
		\Delta X^{\RomI}(q_d,\hat{a})^2\Big]\,,\]
and similarly for $X^{\RomII}$, $X^{\RomIII}$, and $X^{\Ising}$. The claim then follows for $\vX^\trunc$ and $\vX^\frozen$.
\end{proof}
\end{ppn}

\begin{cor}[quenched second moment bounds for $Y^{\RomI}$ and $Y^{\RomI,\frozen}$]
\label{c:Y.I.second.mmt}
For $N$ large enough, we have
	\beq\label{e:Y.I.second.quenched}
	\max\bigg\{
	\bE\Big[\Delta Y^{\RomI}(q_d)^2\Big]
	: 0 \le d\le \dmax-1
	\bigg\}
    \le L^{O(1)} \delta^2 (\MAX_N)^2
    \,,
    \eeq
with $\MAX_N$ as in Proposition~\ref{p:Y.kolmogorov}. The same bound holds with $\vY^\frozen$ in place of $\vY$. 

\begin{proof} It follows using
Lemma~\ref{l:X.Y.discrep} that
	\[\bE\Big(
	\Delta Y^{\RomI}(q_d)^2\Big)
	\le \bE\Big( \Delta X^{\trunc,\RomI}(q_d,a(q_{d+1}))^2\Big)
	+ O(\eta)\bigg\{ \bE\Big| \Delta X^{\trunc,\RomI}(q_d,a(q_{d+1}))\Big|
	+ \eta\bigg\}\,.\]
From the spatial rerandomization (Definition~\ref{d:rerand}), the index $a(q_{d+1})$ is a uniformly random sample from $[M]$, even if we condition on the gaussian disorder. Thus we can apply Proposition~\ref{p:X.second.mmt.quenched} with $a(q_{d+1})$ in place of $\hat{a}$. Together with Jensen's inequality and \eqref{e:deriv.p.bound}, it implies
	\[
	\bE\Big[\Delta Y^{\RomI}(q_d)^2\Big]
    \le L^{O(1)}\Big\{ \delta^2 (\MAX_N)^2
    + \eta^2\Big\}
    \le L^{O(1)}\delta^2 (\MAX_N)^2\,,
	\]
using that $\MAX_N\ge1$ and $\eta\ll\delta$ from Assumption~\ref{a:params}. This proves \eqref{e:Y.I.second.quenched}, and the analogous bound for $Y^{\frozen,\RomI}$ follows by the same argument.
\end{proof}
\end{cor}

Recall $\cG(q_d)$ from \eqref{e:gaus.filt}. We also define the expanded $\sigma$-algebras
	\begin{align}
	\label{e:gaus.filt.plus}
	\cFt(q_d)
	&\equiv \sigma\Big(
	\cG(q_d),
	a(q_0),i(q_0)
	\ldots,
	a(q_{d+1}),i(q_{d+1})
	\Big)\,,\\
	\cFtplus(q_d)
	&\equiv \sigma\Big(
	\cG(q_d),
	a^\frozen(q_0),i^\frozen(q_0)
	\ldots,
	a^\frozen(q_{d+1}),i^\frozen(q_{d+1})
	\Big)\,.
	\end{align}
Recall that $\cG(q_d)$ determines the processes $\vX$, $\vX^\trunc$, and $\vX^\frozen$ up to time $q_d$, so it also determines the partitions
$\mathcal{B}$, $\mathcal{B}^\frozen$, 
$\mathcal{B}^\Ising$, and $\mathcal{B}^{\frozen,\Ising}$ up to time $q_d$. According to the spatial rerandomization (Definition~\ref{d:rerand}), these partitions fully determine the law of 
the indices 
$(a(q_\ell): \ell \le d+1)$ and
$(i(q_\ell): \ell \le d+1)$,
so revealing these indices gives no additional information about $\bXi^\ell$ for $\ell>d$. Likewise, revealing 
the indices
$(a^\frozen(q_\ell): \ell \le d+1)$ and
$(i^\frozen(q_\ell): \ell \le d+1)$ gives no additional information about $\bXi^\ell$ for $\ell>d$. 

\begin{lem}\label{l:subgaus.X.II.III.Ising} 
Let indices $a(q_d),a^\frozen(q_d)$
be sampled as in Definition~\ref{d:rerand}, and let
indices $i(q_d),i^\frozen(q_d)$
be sampled as in Definition~\ref{d:rerand.Ising}. 
For $N$ large enough, we have
	\[\max_{\sigma\in\{\RomII,\RomIII\}}
	\E\bigg[
	\max_{r\in[s,t]} 
	\bigg|
	\sum_{\ell=s}^{r-1}
	\Delta X^\sigma(q_\ell,a(q_{\ell+1}))
	\bigg|^8
	\bigg]
	\le L^{O(1)} (q_t-q_s)^4
	\]
for all $0\le s\le t\le \dmax$.
The same bound holds
	with $\Delta X^{\trunc,\sigma}(q_d,a(q_{d+1}))$
	or $\Delta X^{\frozen,\sigma}(q_d,a^\frozen(q_{d+1}))$
	in place of $\Delta X^\sigma(q_d,a(q_{d+1}))$.
We also have
	\[
	\E \left[
	\max_{r\in[s,t]}
	\left|
	\sum_{\ell=s}^{r-1}
	\Delta X^{\Ising}(q_\ell,i(q_{\ell+1}))
	\right|^8
	\right]
	\le L^{O(1)}(q_t-q_s)^4\,,
	\]
and likewise with $\Delta X^{\trunc,\Ising}(q_\ell,i(q_{\ell+1}))$ or $\Delta X^{\Isplus}(q_\ell,i^\frozen(q_{\ell+1}))$ in place of $\Delta X^{\Ising}(q_\ell,i(q_{\ell+1}))$.

\begin{proof}
This follows by a very similar argument as for Lemma~\ref{l:martingale.sup.bound.X.II.III.Ising}. We will describe the proof for $X^\RomII$. A key difference is that in Lemma~\ref{l:martingale.sup.bound.X.II.III.Ising} we considered $X^\RomII(q_d,a)$ as a process indexed by time $q_d$ for a fixed index $a\in[M]$, with respect to the filtration $\cG(t)$. Here, we instead consider the process $X^\RomII(q_d,a(q_d))$, with respect to the filtration $\cFt(t)$ defined by \eqref{e:gaus.filt.plus}. With this change, repeating the argument of
Lemma~\ref{l:martingale.sup.bound.X.II.III.Ising} gives the conclusion for $X^\RomII$ (cf.\ \eqref{e:martingale.X.holder.estimate}). The conclusions for $X^\RomIII$ and $X^\Ising$ follow similarly. It is straightforward to verify that the same bounds hold for the modified processes 
$X^{\trunc,\sigma}$ and $X^{\frozen,\sigma}$.
\end{proof}
\end{lem}

\begin{proof}[\hypertarget{proof:p.Y.kolmogorov}{Proof of Proposition~\ref{p:Y.kolmogorov}}]
We start with the estimate
\eqref{e:Y.I.Kolmogorov} on $Y^{\RomI}$. In the simplest case, we have $q_d \le s \le t \le q_{d+1}$, in which case we have
	\[
	\bE\Big[\Big( Y^{\RomI}(t)-Y^{\RomI}(s)\Big)^2\Big]
	=\bigg(\frac{t-s}{q_{d+1}-q_d}\bigg)^2
	\bE\Big[\Big( Y^{\RomI}(q_{d+1})-Y^{\RomI}(q_d)\Big)^2\Big] 
	\le L^{O(1)} (\MAX_N)^2 (t-s)^2\,,
	\]
by Corollary~\ref{c:Y.I.second.mmt}. The other possibility is that
$s \le q_a \le \ldots \le q_b \le t$, in which case we can bound
	\begin{align*}
	\bE\Big[\Big( Y^{\RomI}(t)-Y^{\RomI}(s)\Big)^2\Big]
	&\le
	9 \bigg\{
	\bE\Big[\Big( Y^{\RomI}(q_a)-Y^{\RomI}(s)\Big)^2\Big]
	+\bE\Big[\Big( Y^{\RomI}(t)-Y^{\RomI}(q_b)\Big)^2\Big]\\
	&\qquad
	+\bE\Big[\Big( Y^{\RomI}(q_b)-Y^{\RomI}(q_a)\Big)^2\Big]
	\bigg\}\,.
	\end{align*}
The first two terms on the right-hand side can be bounded by the previous calculation. For the third term, the Cauchy--Schwarz inequality gives
	\beq\label{e.CS.short.to.long}
    \bE\Big[\Big( Y^{\RomI}(q_b)-Y^{\RomI}(q_a)\Big)^2\Big]
	\le
	(b-a) \sum_{\ell=a}^{b-1}
	\bE\Big[ \Delta Y(q_\ell)^2
	\Big]\,.
    \eeq
Applying Corollary~\ref{c:Y.I.second.mmt} again gives
	\[\frac{\bE[( Y^{\RomI}(q_b)-Y^{\RomI}(q_a))^2]}
		{L^{O(1)} (\MAX_N)^2}
	\le (b-a)
	\sum_{\ell=a}^{b-1}
	(\delta_\ell)^2
	\le (b-a)\delta
	\sum_{\ell=a}^{b-1} \delta_\ell
	=(b-a)\delta (q_b-q_a)
	\le L^{O(1)}
	(q_b-q_a)^2\,,
	\]
where the last bound uses \eqref{e:delta.cutoffs} and \eqref{e:deriv.p.bound}. Altogether we conclude
	\[
	\frac{\bE[( Y^{\RomI}(t)-Y^{\RomI}(s))^2]}
		{ (\MAX_N)^2}
	\le L^{O(1)} \bigg\{
	(q_a-s)^2 + (t-q_b)^2
	+ (q_b-q_a)^2
	\bigg\}
	\le L^{O(1)} (t-s)^2\,,
	\]
which proves the quenched second moment bound \eqref{e:Y.I.Kolmogorov} on $Y^{\RomI}$. The estimates on the random variable $\MAX_N$ from \eqref{e:MAX.N} follow by combining
Lemmas \ref{l:bar.x.subgaus},
\ref{l:y.subgaus},
\ref{l:sg.norm.bound}, and \ref{l:wishart}, and making use of Assumption~\ref{a:params}. The eighth moment bounds \eqref{e:Y.II.III.Is.Kolmogorov} on $Y^\sigma$ ($\sigma\in\{\RomII,\RomIII,\Ising\}$) follow more straightforwardly by combining Lemma~\ref{l:X.Y.discrep} with Lemma~\ref{l:subgaus.X.II.III.Ising},  using that  $\eta \ll\delta$ from Assumption~\ref{a:params}.
The bounds for $\vY^\frozen$ follow by repeating the same argument.
\end{proof}

\subsection{Moment estimates for drift}
\label{ss:kolmogorov.drift}

In this subsection we define and analyze the ``drift'' of $\vY$ and $\vY^\frozen$. \textbf{The main result of this subsection is Proposition~\ref{p:drift.kolmogorov} below. In addition we prove several quenched estimates, in particular Proposition~\ref{p:drift.quenched} and Lemma~\ref{l:barD.I.second.quenched}, which will be used in Section~\ref{s:sde}.}

\begin{dfn}[discrete drift processes] \label{d:D}
We let $\vD(s)$ and $\vD^\frozen(s)$ be the processes which start from 
 $\vD(q_0)= \vY(q_0)$
and $\vD^\frozen(q_0)= \vY^\frozen(q_0)$, 
are linear on each interval $q_d \le s \le q_{d+1}$, and have increments defined by
	\begin{align}
	\label{e:def.drift}
	\Delta\vD(q_d) 
	&= \bE\Big[ \vY(q_{d+1})
		\,\Big|\,\cG(q_d)\Big] 
		- \vY(q_d)
	=\bE\Big[ \Delta\vY(q_d) 
	\,\Big|\,\cG(q_d)\Big]
	\,,\\
	\Delta\vD^\frozen(q_d) 
	&= \bE\Big[ \vY^\frozen(q_{d+1})
		\,\Big|\,\cG(q_d)\Big] 
		- \vY^\frozen(q_d)
	=\bE\Big[ \Delta\vY^\frozen(q_d) 
	\,\Big|\,\cG(q_d)\Big]
	\,,
	\label{e:def.drift.plus}
	\end{align}
with $\cG(q_d)$ as in \eqref{e:gaus.filt}.
\end{dfn}

Note that for any $s=q_a < q_b=t$, 
the bound \eqref{e:Y.II.III.Is.Kolmogorov} from Proposition~\ref{p:Y.kolmogorov} immediately gives
	\begin{align}\nonumber
	&\E\bigg[\Big( D^\sigma(t)-D^\sigma(s)\Big)^8\bigg]
	= \E\bigg[
	\Big(
	\bE[ Y^\sigma(t)-Y^\sigma(s)
		\,|\, \cG(s) ] \Big)^8\bigg] \\
	&\qquad \le \E\bigg[ \Big(Y^\sigma(t)-Y^\sigma(s) \Big)^8\bigg]
	\le L^{O(1)} (t-s)^4
	\label{e:D.naive.kolmogorov.bound}
	\end{align}
for $\sigma\in\{\RomII,\RomIII,\Ising\}$. Next, for further analysis, we define an approximate version of $\vD$ that will be easier to work with:

\begin{dfn}[approximate drift processes]\label{d:barD}
Define
$\barvD$ to be the process which starts from 
$\barvD(q_0)=\vY(q_0)$, is linear on each interval $q_d \le s \le q_{d+1}$, and has increments given by
	\beq\label{e:barD}
	\Delta\bar{D}^\sigma(q_d)
	\equiv \begin{cases}
	\displaystyle
	\Delta\bar{D}^\sigma\Big(q_d,B(q_d,a(q_d))\Big)
	\equiv
	\sum_{a\in B(q_d,a(q_d))}
	\frac{\Delta X^{\trunc,\sigma}(q_d,a)}
		{|B(q_d,a(q_d))|} 
	& \textup{for $\sigma\in\{\RomI,\RomII,\RomIII\}$,}\\
	\displaystyle
	\Delta\bar{D}^\sigma\Big(q_d,B^\Ising(q_d,i(q_d))\Big)
	\equiv \sum_{i\in B^\Ising(q_d,i(q_d))}
	\frac{\Delta X^{\trunc,\Ising}(q_d,i)}
		{|B^\Ising(q_d,i(q_d))|} 
	& \textup{for $\sigma=\Ising$.}\\
		\end{cases}
	\eeq
Similarly, we define $\barvD^\frozen$ to be the process which starts from 
$\barvD^\frozen(q_0)=\vY^\frozen(q_0)$, is linear on each interval $q_d \le s \le q_{d+1}$, and has increments given by
	\beq\label{e:barD.plus}
	\Delta\bar{D}^{\frozen,\sigma}(q_d)
	\equiv \begin{cases}
	\displaystyle
	\sum_{a\in B^\frozen(q_d,a^\frozen(q_d))}
	\frac{\Delta X^{\frozen,\sigma}(q_d,a)}
		{|B^\frozen(q_d,a^\frozen(q_d))|}
	& \textup{for $\sigma\in\{\RomI,\RomII,\RomIII\}$,}\\
	\displaystyle
	\sum_{i\in B^{\Isplus}(q_d,i^\frozen(q_d))}
	\frac{\Delta X^{\Isplus}(q_d,i)}
		{|B^\Isplus(q_d,i^\frozen(q_d))|} 
	& \textup{for $\sigma=\Ising$.}\\
		\end{cases}
	\eeq
Write $B^\frozen=B^\frozen(q_d,a^\frozen(q_d))$, and note that 
if $B^\frozen\cap[M]=\emptyset$ then $
\Delta\bar{D}^{\frozen,\sigma}(q_d)=0$ for all $\sigma\in\{\RomI,\RomII,\RomIII\}$. Otherwise,
if $B^\frozen\cap[M]=B\ne\emptyset$, then we have
	\beq\label{e:barD.vs.plus}
	\Delta\bar{D}^{\frozen,\sigma}(q_d,B^\frozen)
	=
	\frac{|B|}{|B^\frozen|}
	\Delta\bar{D}^\sigma(q_d,B)
	\eeq
for all $\sigma\in\{\RomI,\RomII,\RomIII\}$. A similar relation holds for $\sigma=\Ising$.
\end{dfn}

We will show in Lemma~\ref{l:D.barD.discrep} that $\vD$ and $\barvD$ are very close, so it suffices to analyze $\barvD$. Likewise, it suffices to analyze $\barvD^\frozen$ in place of $\vD^\frozen$. 
Recall the stopping times 
$q_T$, $q_{\Ising,T}$, $q_{\frozen,T}$, $q_{\frozen,\Ising,T}$ 
from \S\ref{ss:tightness.statements}. Recall from \eqref{e:X.decomp} that $\Delta X^\RomI(q_d)$ has a factor $p'(q_d)^{1/2}$. Recalling 
\eqref{e:def.b} and \eqref{e:review.b}, we have
	\begin{align} 
	\nonumber 
	\bbb^\trunc(q_d)
	&\equiv
	\frac{\Delta\bar{D}^{\RomI}(q_d)}{\delta_d p'(q_d)^{1/2}}
	=\ind\{ q_d < q_T\} 
	\frac{1}{|B|}
	\sum_{a\in B} 
	\frac{(\bmeta^a,\by)}{N^{1/2}}\,, \\
	\label{e:def.b.plus}
	\bbb^\frozen(q_d)
	&\equiv
	\frac{\Delta\bar{D}^{\frozen,\RomI}(q_d)}{\delta_d p'(q_d)^{1/2}}
	= \ind\{ q_d < q_{\frozen,T}\} 
	\frac{1}{|B^\frozen|}
	\sum_{a\in B^\frozen \cap [M]}
	\frac{(\bmeta^a,\by)}{N^{1/2}}\,,
	\end{align}
where $B=B(q_d,a(q_d))\in\mathcal{B}(q_d)$ and $B^\frozen=B^\frozen(q_d,a^\frozen(q_d))\in\mathcal{B}^\frozen(q_d)$.
For $\sigma\in\{\RomII,\RomIII,\Ising\}$, we have
	\begin{align} 
	\label{e:drift.error.II}
	\frac{\Delta\bar{D}^{\RomII}(q_d)}{\delta_d}
	&= \ind\{q_d< q_T\}
	\bigg(\frac{p_d}{\delta_d}\bigg)^{1/2}
	\frac{1}{|B|}
	\sum_{a\in B} \frac{(\bar{\bg}^a,\by)}{N^{1/2}}
	\,,\\
	\frac{\Delta\bar{D}^{\RomIII}
		(q_d)}{\delta_d}
	&= \ind\{q_d< q_T\}
	\frac{[q_d(p_{d+1}-p_d)]^{1/2}}{\delta_d}
	\frac{1}{|B|}
	\sum_{a\in B} 
	\frac{(\bmeta^a,\bar{\bx})}{N^{1/2}}
	\,,\label{e:drift.error.III}\\
	\label{e:drift.error.Ising}
	\frac{\Delta \bar{D}^{\Ising}(q_d)}{\delta_d}
	&=\ind\{q_d< q_{\Ising,T}\}
	(\delta_d)^{-1/2}
	\frac{1}{|B^{\Ising}|}\sum_{i\in B^{\Ising}}
	(\be_i,\by)\,,
	\end{align}
where $B^{\Ising}=B^{\Ising}(q_d,i(q_d))\in\mathcal{B}^{\Ising}(q_d)$. Note that
$\Delta\bar{D}^\sigma(q_d)$ is determined by the buckets $B$ and $B^\Ising$; we shall sometimes emphasize this by writing
$\Delta\bar{D}^\sigma(q_d,B)
$ for $\sigma\in\{\RomI,\RomII,\RomIII\}$ and 
$\Delta\bar{D}^\Ising(q_d,B^\Ising)$. Similarly we will also denote
$\bbb^\trunc(q_d)\equiv\bbb^\trunc(q_d,B)$ and $\bbb^\frozen(q_d)\equiv\bbb^\frozen(q_d,B^\frozen)$. 
\textbf{We will show that $\bar{D}^{\RomI}$ and $\bar{D}^{\frozen,\RomI}$ are the only asymptotically nontrivial processes. More precisely, the main results of this subsection are the following:}

\begin{ppn}[estimates for $\barvD$ and $\barvD^\frozen$]
\label{p:drift.kolmogorov}
Let $\barvD$ and $\barvD^\frozen$ be as specified by Definition~\ref{d:barD}.
 For $N$ large enough, it holds for all
 $q_0\le s < t\le1$ that
	\beq\label{e:bar.D.I.second}
	\frac{\E[( \bar{D}^{\RomI}(t)-\bar{D}^{\RomI}(s))^2]}
		{(t-s)^2}
	\le L^{O(1)}\,,\eeq
as well as (with $\BUCKETS$ as in \eqref{e:buckets.scB})
	\beq\label{e:bar.D.notI.unif}
	\sum_{\sigma\in\{\RomII,\RomIII,\Ising\}}
	\E\sup\bigg\{
	\Big|\bar{D}^\sigma(t)\Big| : 0\le t\le 1\bigg\}
	\le\bigg( \frac{L^{O(1)} \BUCKETS}{N\delta}
	\bigg)^{1/2}\,.
	\eeq
The same estimates hold with $\barvD^\frozen$ in place of
$\barvD$.
\end{ppn}

\begin{ppn}[quenched estimates for $\vD$ and $\vD^\frozen$]\label{p:drift.quenched} 
We have
	\[\frac{\bE[( D^\RomI(t)-D^\RomI(s))^2]}{(t-s)^2}
	\le L^{O(1)} (\MAX_N)^2\]
where $\MAX_N$ is as in \eqref{e:MAX.N} and Proposition~\ref{p:X.second.mmt.quenched}; and
	\[\sum_{\sigma\in\{\RomII,\RomIII,\Ising\}}
	\frac{\bE[( D^\sigma(t)- D^\sigma(s))^2]}{(t-s)^2}
	\le L^{O(1)} 
	\bigg\{ 
	\frac{\eta^2}{\delta^2}+
	\frac{(\SMAX_N)^2\BUCKETS}{N\delta}\bigg\}\]
where $\SMAX_N$ is a random variable 
(defined explicitly by \eqref{e:SMAX.N})
satisfying \eqref{e.SMAX.bound}.
The same bound holds with $\vD^\frozen$ in place of
$\vD$.
\end{ppn}

The \hyperlink{proof:p.drift.kolmogorov}{proof of Proposition~\ref{p:drift.kolmogorov}} 
and \hyperlink{proof:p.drift.quenched}{proof of Proposition~\ref{p:drift.quenched}} 
appear at the end of this subsection. They are obtained from a few intermediate estimates which we turn to next.

\begin{proof}[\hypertarget{proof:l.barD.I.second.quenched}{Proof of Lemma~\ref{l:barD.I.second.quenched}}]
It follows from Jensen's inequality and Proposition~\ref{p:X.second.mmt.quenched} that
	\beq\label{e.jensen.example}
	\bE[\bbb^\trunc(q_d)^2]
	=\frac{\bE[\Delta \bar{D}^{\RomI}(q_d)^2]}
		{(\delta_d)^2  p'(q_d)}
	\le
	\frac{\bE[\Delta X^{\trunc,\RomI}(q_d,\hat{a})^2]}
		{(\delta_d)^2  p'(q_d)}
	\le
	(\MAX_N)^2\,.
	\eeq
The bound for $\bar{D}^{\frozen,\RomI}$ follows by the same argument, since Proposition~\ref{p:X.second.mmt.quenched} applies to $X^{\frozen,\RomI}$ as well as to $X^{\trunc,\RomI}$.
\end{proof}

We next prove annealed second moment bounds (Lemmas \ref{l:barD.II.second}--\ref{l:barD.Ising.second}) on increments of $\bar{D}^\sigma$ for $\sigma\in\{\RomII,\RomIII,\Ising\}$, which will be used in the proof of 
\eqref{e:bar.D.notI.unif} in
Proposition~\ref{p:drift.kolmogorov}.

\begin{lem}[annealed second moment bounds for 
$\bar{D}^{\RomII}$ and 
$\bar{D}^{\frozen,\RomII}$]
\label{l:barD.II.second}
Let $\bar{D}^{\RomII}$ and $\bar{D}^{\frozen,\RomII}$ 
be defined by \eqref{e:barD} and \eqref{e:barD.plus} (see also \eqref{e:drift.error.II}). We have 
the second moment bound 
	\[\frac{\E[\Delta\bar{D}^{\RomII}(q_d)^2]}
		{p(q_d)}
	\le \frac{L^{O(1)}\BUCKETS \delta}{N}\,,
	\]
for all $0\le d\le \dmax-1$. The same bound holds with $\bar{D}^{\frozen,\RomII}$ in place of $\bar{D}^{\RomII}$.
\begin{proof}
Recalling the notation from
\eqref{e:normalized.vectors.g} and \eqref{e:normalized.vectors.y},
we can rewrite \eqref{e:drift.error.II} as
	\[
	\frac{\Delta\bar{D}^{\RomII}(q_d)}
		{[\delta_d p(q_d)]^{1/2}}
	= 
	\frac1{|B|^{1/2}}
	\bigg(\frac{1}{N^{1/2}}\sum_{a\in B}
	\frac{\bar{\bg}^a(q_d)}{|B|^{1/2}},\by(q_{d+1})\bigg)
	\equiv
	\frac1{|B|^{1/2}}
	\bigg( \frac{\bar{\bg}^B(q_d)}{N^{1/2}},\by(q_{d+1}) \bigg)\,.
	\]
Note that $\bar{\bg}^B(q_d)$ is measurable with respect to $\cG(q_d)$. 
It follows by applying Lemma~\ref{l:y.subgaus} that 
	\[
	\frac{\E[\Delta\bar{D}^{\RomII}(q_d)^2]}{\delta_d p(q_d)}
	=\frac{1}{M}\E\bigg[\sum_{B\in\mathcal{B}(q_d)} 
	\bigg( \frac{\bar{\bg}^B(q_d)}{N^{1/2}},\by(q_{d+1}) \bigg)^2
		\bigg]
	\le \frac{L^{O(1)}}{M}
	\E\sum_{B\in\mathcal{B}(q_d)} 
	\frac{\|\bar{\bg}^B(q_d)\|^2}{N}\,.
	\]
We then note that if $\ind_B\in\R^M$ is the indicator of the coordinates in $B$, then
	\[
	\|\bar{\bg}^B(q_d)\|
	=\bigg\| \bar{\bG}(q_d) \frac{\ind_B}{|B|^{1/2}}\bigg\|
	\le \|\bar{\bG}(q_d)\|\,,
	\]
from which it follows that
	\[\frac{\E[\Delta\bar{D}^{\RomII}(q_d)^2]}{\delta_d p(q_d)}
	\le
	\frac{L^{O(1)}\BUCKETS 
		\E[\|\bar{\bG}(q_d)\|^2]}{MN}
	\le \frac{L^{O(1)}\BUCKETS}{N}\,,
	\]
where the last bound follows from Lemma~\ref{l:wishart}.
This proves the claim for $\bar{D}^{\RomII}$, and the claim for
$\bar{D}^{\frozen,\RomII}$ follows by using \eqref{e:barD.vs.plus}.
\end{proof}
\end{lem}

\begin{lem} [annealed second moment bounds for 
$\bar{D}^{\RomIII}$ and 
$\bar{D}^{\frozen,\RomIII}$]
\label{l:barD.III.second}
Let $\bar{D}^{\RomIII}$ and $\bar{D}^{\frozen,\RomIII}$ 
be defined by \eqref{e:barD} and \eqref{e:barD.plus} (see also \eqref{e:drift.error.III}).  We then have
	\[
	\frac{\E[\Delta\bar{D}^{\RomIII}(q_d)^2]}
		{q_d p'(q_d)}
	\le \frac{L^{O(1)}\BUCKETS \delta}{N}
	\]
for all $0\le d\le \dmax-1$. 
The same bound holds with $\bar{D}^{\frozen,\RomIII}$ in place of $\bar{D}^{\RomIII}$.
\begin{proof}
Recalling the notation from \eqref{e:normalized.vectors.x},
we can rewrite \eqref{e:drift.error.III} as
	\[
	\frac{\Delta\bar{D}^\RomIII(q_d)}
		{[\delta_d q_d p'(q_d)]^{1/2}}
	= \frac1{|B|^{1/2}}
	\bigg(
	\sum_{a\in B}
	\frac{\bmeta^a(q_{d+1})}{|B|^{1/2}},
		\frac{\bar{\bx}(q_d)}{N^{1/2}}
		\bigg)
	\equiv
	\frac1{|B|^{1/2}}
	\bigg(\bmeta^B(q_{d+1}),
	\frac{\bar{\bx}(q_d)}{N^{1/2}}
		\bigg)\,.
	\]
Note that $\bar{\bx}(q_d)$ is measurable with respect to $\cG(q_d)$, and $\bmeta^B(q_{d+1})$ is a standard gaussian in $\R^N$ conditional on $\cG(q_d)$. It follows that
	\[\frac{\E[\Delta\bar{D}^\RomIII(q_d)^2]}
		{\delta_d q_d p'(q_d)}
	\le 
	\frac{1}{M}\E\bigg[ \sum_{B\in\mathcal{B}(q_d)} 
	\bigg(\bmeta^B(q_{d+1}),
	\frac{\bar{\bx}(q_d)}{N^{1/2}}
		\bigg)^2\bigg]
	\le \frac{O(1)\BUCKETS
		\E[\|\bar{\bx}(q_d)\|^2]}{MN}\,,
	\]
and the claim for $\bar{D}^\RomIII$ follows by combining Lemma~\ref{l:bar.x.subgaus} with Lemma~\ref{l:sg.norm.bound}. The claim for
$\bar{D}^{\frozen,\RomIII}$ then follows by using \eqref{e:barD.vs.plus}.
\end{proof}
\end{lem}

\begin{lem}[annealed second moment bounds for $\bar{D}^{\Ising}$ and $\bar{D}^{\frozen,\Ising}$]
\label{l:barD.Ising.second}
For $\bar{D}^{\Ising}$ as in \eqref{e:drift.error.Ising}, we have
	\[\E\Big[\Delta\bar{D}^{\Ising}(q_d)^2\Big]
	\le\frac{ L^{O(1)}\BUCKETS\delta}
		{N}
	 \]
for all $0\le d\le \dmax-1$, with $\BUCKETS$ as in 
\eqref{e:buckets.scB}. The same bound holds with $\bar{D}^{\frozen,\Ising}$ in place of $\bar{D}^\Ising$.
	
\begin{proof} Recalling the notation from \eqref{e:normalized.vectors.y}, we can rewrite \eqref{e:drift.error.Ising} as
	\[
	\frac{\Delta \bar{D}^{\Ising}(q_d)}{(\delta_d)^{1/2}}
	=\frac{1}{|B|^{1/2}}
	\bigg(\sum_{i\in B}
	\frac{\be_i}{|B|^{1/2}},
	\by(q_{d+1})\bigg)
	\equiv \frac{(\be_B,\by(q_{d+1}))}{|B|^{1/2}}\,,
	\]
where we abbreviate $B\equiv B^\Ising$.
The vector $\be_B$ is measurable with respect to $\cG(q_d)$. It follows by applying 
Lemma~\ref{l:y.subgaus} that
	\[
	\frac{\E[\Delta \bar{D}^{\Ising}(q_d)^2]}{\delta_d}
	= \frac{1}{N}\E\bigg[
	\sum_{B\in\mathcal{B}^{\Ising}(q_d)}
	(\be_B,\by)^2\bigg]
	\le \frac{L^{O(1)}\BUCKETS}{N}\,.
	\]
This proves the claim for $\bar{D}^{\Ising}$, and the claim for
$\bar{D}^{\frozen,\Ising}$ follows by using \eqref{e:barD.vs.plus}.
\end{proof}
\end{lem}

In addition to Lemmas \ref{l:barD.II.second}--\ref{l:barD.Ising.second}, we also prove quenched second moment bounds for $\bar{D}^\sigma$ and $\bar{D}^{\frozen,\sigma}$ for $\sigma\in\{\RomII,\RomIII,\Ising\}$ in Lemmas \ref{l:barD.II.second.quenched}--\ref{l:barD.Ising.second.quenched} below, which will be used in the proof of Proposition~\ref{p:drift.quenched}.

\begin{lem}[quenched second moment bounds for 
$\bar{D}^{\RomII}$ and 
$\bar{D}^{\frozen,\RomII}$]
\label{l:barD.II.second.quenched}
In addition to Lemma~\ref{l:barD.II.second} we have the quenched second moment bound
	\[\max\bigg\{
	\frac{\bE[\Delta\bar{D}^{\RomII}(q_d)^2]}{p(q_d)}
	: 0\le d\le \dmax-1\bigg\}
	\le \frac{\SMAX_N(\RomII)^2 \BUCKETS \delta}{M}\,,\]
where $\SMAX_N(\RomII)$ is a random variable satisfying the bound
	\[\P\Big(\SMAX_N(\RomII) \ge \log N\Big)
	\le\exp\bigg(-\frac{(\log N)^2}{L^{O(1)}}\bigg)\,.\]
The same bound holds with $\bar{D}^{\frozen,\RomII}$ in place of
$\bar{D}^{\RomII}$.

\begin{proof} Recall that
we can rewrite \eqref{e:drift.error.II} as
	\[
	\frac{\Delta\bar{D}^{\RomII}(q_d)}{(p_d)^{1/2}}
	= 
	\frac{(\delta_d)^{1/2}}{|B|^{1/2}}
	\bigg(\frac{1}{N^{1/2}}\sum_{a\in B}
	\frac{\bar{\bg}^a(q_d)}{|B|^{1/2}},\by(q_{d+1})\bigg)
	\equiv
	\frac{(\delta_d)^{1/2}}{|B|^{1/2}}
	\bigg( \frac{\bar{\bg}^B(q_d)}{N^{1/2}},\by(q_{d+1}) \bigg)\,.
	\]
If we define the random variable
	\[
	\SMAX_N(\RomII)
	\equiv\max
	\bigg\{
	\bigg|
	\bigg( \frac{\bar{\bg}^B(q_d)}{N^{1/2}},\by(q_{d+1}) \bigg)
	\bigg|
	: 
	0\le d\le \dmax-1,
	B\in \mathcal{B}(q_d)
	\bigg\}\,,
	\]
then we have the quenched second moment bound
	\[
	\frac{\bE[\Delta\bar{D}^{\RomII}(q_d)^2]}{p_d}
	= \frac{\delta_d}{M}\sum_{B\in\mathcal{B}(q_d)}
	\bigg( \frac{\bar{\bg}^B(q_d)}{N^{1/2}},\by \bigg)^2
	\le \frac{\SMAX_N(\RomII)^2\BUCKETS\delta}{M}\,,
	\]
Note that for any $B\subseteq[M]$ we have
$\bar{\bg}^B= \bar{\bG}^\st \be_B$, so $\|\bar{\bg}^B\| \le\|\bar{\bG}\|_\textup{op}$. It then follows from Lemma~\ref{l:wishart} that for $C(\alpha)$ large enough, we have 
	\beq\label{e:II.union.bound}
	\P\Big(\|\bar{\bg}^B\|^2 \ge 
	C(\alpha) N \textup{ for any }B\subseteq[M]\Big)
	\le \frac{1}{2^N}\,.
	\eeq
The vectors $(\bar{\bg}^B : B\in\mathcal{B}(q_d))$ are measurable with respect to $\cG(q_d)$, so applying Lemma~\ref{l:y.subgaus} gives
	\begin{align*}
	&\P\Big(\SMAX_N(\RomII) \ge \log N\Big)
	=\P\bigg(
	\bigg|
	\bigg( \frac{\bar{\bg}^B(q_d)}{N^{1/2}},\by(q_{d+1}) \bigg)
	\bigg|
		\ge \log N
	\textup{ for any } 0\le d\le \dmax- 1, B \in \mathcal{B}(q_d)
	\bigg)\\
	&\qquad\le
	\dmax\BUCKETS
	\bigg\{ \frac{1}{2^N}
	+ \exp\bigg(-\frac{(\log N)^2}{L^{O(1)}}\bigg)
	\bigg\}
	\le\exp\bigg(-\frac{(\log N)^2}{L^{O(1)}}\bigg)
	\,,
	\end{align*}
where the $1/2^N$ in the second line comes from \eqref{e:II.union.bound}, and the final simplification uses Assumption~\ref{a:params}. This proves the claim for $\bar{D}^{\RomII}$, and the claim for
$\bar{D}^{\frozen,\RomII}$ follows by using \eqref{e:barD.vs.plus}.
\end{proof}
\end{lem}

\begin{lem}[quenched second moment bounds for 
$\bar{D}^{\RomIII}$ and 
$\bar{D}^{\frozen,\RomIII}$]
\label{l:barD.III.second.quenched}
In addition to Lemma~\ref{l:barD.III.second} we have the quenched second moment bound
	\[\max\bigg\{
	\frac{\bE[\Delta\bar{D}^{\RomIII}(q_d)^2]}{q_d p'(q_d)}
	: 0\le d\le\dmax-1
	\bigg\}
	\le \frac{\SMAX_N(\RomIII)^2\BUCKETS\delta}{M}\]
where $\SMAX_N(\RomIII)$ is a random variable satisfying the bound
	\[\P\Big(\SMAX_N(\RomIII) \ge \log N\Big)
	\le \exp
	\bigg(
	-\frac{(\log N)^2}{L^{O(1)}}
	\bigg)
	\,.\]
The same bound holds with $\bar{D}^{\frozen,\RomIII}$ in place of
$\bar{D}^{\RomIII}$.
\begin{proof}
Recall that we can rewrite \eqref{e:drift.error.III} as
	\[
	\frac{\Delta\bar{D}^\RomIII(q_d)}
		{[q_d p'(q_d)]^{1/2}}
	= \frac{(\delta_d)^{1/2}}{|B|^{1/2}}
	\bigg(
	\sum_{a\in B}
	\frac{\bmeta^a(q_{d+1})}{|B|^{1/2}},
		\frac{\bar{\bx}(q_d)}{N^{1/2}}
		\bigg)
	\equiv
	\frac{(\delta_d)^{1/2}}{|B|^{1/2}}
	\bigg(\bmeta^B(q_{d+1}),
	\frac{\bar{\bx}(q_d)}{N^{1/2}}
		\bigg)\,.
	\]
If we define the random variable
	\[
	\SMAX_N(\RomIII)
	\equiv\max\bigg\{
	\bigg|\bigg(\bmeta^B(q_{d+1}),
	\frac{\bar{\bx}(q_d)}{N^{1/2}}
		\bigg)\bigg|
	: 
	0\le d\le \dmax-1,
	B\in \mathcal{B}(q_d)
	\bigg\}\,,
	\]
then we have the quenched second moment bound
	\[
	\frac{\bE[\Delta\bar{D}^{\RomIII}(q_d)^2]}{q_d p'(q_d)}
	\le \frac{\delta_d}{M}
	\sum_{B\in\mathcal{B}(q_d)}
	\bigg(\bmeta^B(q_{d+1}),
	\frac{\bar{\bx}(q_d)}{N^{1/2}}
		\bigg)^2
	\le \frac{\SMAX_N(\RomIII)^2\BUCKETS\delta}{M}\,.
	\]
The vector $\bar{\bx}(q_d)$ is measurable with respect to $\cG(q_d)$; and it follows from Lemmas~\ref{l:bar.x.subgaus} and \ref{l:sg.norm.bound} that
	\beq\label{e:III.norm.bound}
	\P\bigg(\frac{\|\bar{\bx}(q_d)\|^2}{N} \ge L^{O(1)}\bigg)
	\le \frac{1}{2^N}
	\eeq
provided we take the $L^{O(1)}$ large enough. Conditional on $\cG(q_d)$, the vector $\bmeta^B(q_{d+1})$ is distributed as a standard gaussian in $\R^N$. Therefore,
	\begin{align*}
	&\P\Big(
	\SMAX_N(\RomIII)\ge \log N
	\Big)
	\le \P\bigg(
	\bigg|\bigg(\bmeta^B(q_{d+1}),
	\frac{\bar{\bx}(q_d)}{N^{1/2}}
		\bigg)\bigg|
	\ge\log N
	\textup{ for any }
	0\le d\le \dmax-1,
	B\in\mathcal{B}(q_d)
	\bigg) \\ 
	&\qquad \le \dmax\BUCKETS
	\bigg\{ \frac{1}{2^N}
	+ \exp
	\bigg(
	-\frac{(\log N)^2}{L^{O(1)}}
	\bigg) \bigg\}
	\le \exp
	\bigg(
	-\frac{(\log N)^2}{L^{O(1)}}
	\bigg)
	\end{align*}
where the $1/2^N$ in the second line comes from \eqref{e:III.norm.bound}, and the final simplification uses Assumption~\ref{a:params}. This proves the claim for $\bar{D}^{\RomIII}$, and the claim for
$\bar{D}^{\frozen,\RomIII}$ follows by using \eqref{e:barD.vs.plus}.
\end{proof}
\end{lem}

\begin{lem}[quenched second moment bounds for 
$\bar{D}^{\Ising}$ and 
$\bar{D}^{\frozen,\Ising}$]
\label{l:barD.Ising.second.quenched}
In addition to Lemma~\ref{l:barD.Ising.second} we have the quenched second moment bound
	\[\max\bigg\{
	\bE\Big[\Delta\bar{D}^{\Ising}(q_d)^2\Big]
	: 0\le d\le\dmax-1
	\bigg\}
	\le \frac{\SMAX_N(\Ising)^2\BUCKETS\delta}{N}\]
where $\SMAX_N(\Ising)$ is a random variable satisfying the bound
	\[\P\Big(\SMAX^N(\Ising) \ge \log N\Big)
	\le \exp
	\bigg(
	-\frac{(\log N)^2}{L^{O(1)}}
	\bigg)
	\,.\]
The same bound holds with $\bar{D}^{\frozen,\Ising}$ in place of
$\bar{D}^{\Ising}$.

\begin{proof}
Recall that we can rewrite \eqref{e:drift.error.Ising} as
	\[
	\frac{\Delta \bar{D}^{\Ising}(q_d)}{(\delta_d)^{1/2}}
	=\frac{1}{|B|^{1/2}}
	\bigg(\sum_{i\in B}
	\frac{\be_i}{|B|^{1/2}},
	\by\bigg)
	\equiv \frac{(\be_B,\by(q_{d+1}))}{|B|^{1/2}}\,,
	\]
where we abbreviate $B\equiv B^\Ising$.
If we define the random variable
	\[
	\SMAX_N(\Ising)
	\equiv\max
	\bigg\{
	\Big|( \be_B,\by(q_{d+1}))\Big|
	: 
	0\le d\le \dmax-1, B\in \mathcal{B}^\Ising(q_d)
	\bigg\}\,,
	\]
then we have the quenched second moment bound
	\[
	\bE[\Delta \bar{D}^{\Ising}(q_d)^2]
	= \frac{\delta_d}{N}
	\sum_{B\in \mathcal{B}^\Ising(q_d)}
	(\be_B,\by(q_{d+1}))^2
	\le  \frac{\SMAX_N(\Ising)^2\BUCKETS\delta}{N}\,.
	\]
The bound on $\SMAX_N(\Ising)$ follows by a similar argument as the bound for $\SMAX_N(\RomII)$ from Lemma~\ref{l:barD.II.second.quenched}, with the simplification that $\|\be_B\|=1$ is trivially bounded. 
This proves the claim for $\bar{D}^{\Ising}$, and the claim for
$\bar{D}^{\frozen,\Ising}$ follows by using \eqref{e:barD.vs.plus}.
\end{proof}
\end{lem}

\begin{proof}[\hypertarget{proof:p.drift.kolmogorov}{Proof of Proposition~\ref{p:drift.kolmogorov}}] 
Combining 
Lemmas \ref{l:barD.II.second}--\ref{l:barD.Ising.second}
with the Cauchy--Schwarz inequality gives
	\[\sum_{\sigma\in\{\RomII,\RomIII,\Ising\}}
	\E|\Delta \bar{D}^\sigma(q_d)|
	\le
	\sum_{\sigma\in\{\RomII,\RomIII,\Ising\}}
	\E\Big[ \Delta \bar{D}^\sigma(q_d)^2
	\Big]^{1/2}
	\le\bigg( \frac{L^{O(1)} \BUCKETS\delta}{N}\bigg)^{1/2}\,.
	\]
Summing over $0\le d\le \dmax-1 \le L^{O(1)}/\delta$ gives the estimate \eqref{e:bar.D.notI.unif} on $\bar{D}^\sigma$ for $\sigma\ne\RomI$.
The second moment estimate \eqref{e:bar.D.I.second} on increments of $\bar{D}^{\RomI}$
follows by the same argument as in the
\hyperlink{proof:p.Y.kolmogorov}{proof of
\eqref{e:Y.I.Kolmogorov} in Proposition~\ref{p:Y.kolmogorov}}, using 
the estimate provided
by Lemma~\ref{l:barD.I.second.quenched}
for the individual intervals $[q_d,q_{d+1}]$, and recalling $\E[(\MAX_N)^2] \le L^{O(1)}$ from Proposition~\ref{p:Y.kolmogorov}. The corresponding bounds for $\barvD^\frozen$ follow by the same argument.
\end{proof}

The next lemma verifies that the discrepancy between $\vD$ and $\barvD$ is small:

\begin{lem}\label{l:D.barD.discrep}
Let $\vD$ and $\barvD$ be as specified by Definitions~\ref{d:D} and \ref{d:barD}. Then
	\[\max\bigg\{
	\Big\|\Delta\vD(q_d) 
	-\Delta \barvD(q_d) \Big\|_\infty,
	\Big\|\Delta \vD^\frozen(q_d) 
	-\Delta\barvD^\frozen(q_d) \Big\|_\infty
	\bigg\} \le \eta\,.
    \]
\end{lem}

\begin{proof}
Consider $\sigma\in\{\RomI,\RomII,\RomIII\}$.
Recall from the definition \eqref{e:def.Y} that 
$Y^\sigma(q_d)$ is the average of $X^{\trunc,\sigma}(q_d,a)$ over $a\in B=B(q_d,a(q_d))$. Meanwhile, from the proof 
of Lemma~\ref{l:X.Y.discrep} (see \eqref{e:Y.next}), $Y^\sigma(q_{d+1})$ is $\eta$-close to $X^{\trunc,\sigma}(q_{d+1},a(q_{d+1}))$, where $a(q_{d+1})$ conditional on $\cG(q_d)$ is sampled uniformly at random from $B$. By comparing \eqref{e:def.drift} with \eqref{e:barD}, this implies that $\Delta D^\sigma$ and $\Delta \bar{D}^\sigma$ are $\eta$-close for $\sigma\in\{\RomI,\RomII,\RomIII\}$. The remaining bounds follow by the same argument.
\end{proof}

\begin{proof}
[\hypertarget{proof:p.drift.quenched}{Proof of Proposition~\ref{p:drift.quenched}}] 
First, it follows from
\eqref{e:deriv.p.bound} and Lemma~\ref{l:barD.I.second.quenched} that
	$\bE[\Delta \bar{D}^{\RomI}(q_d)^2]
	\le L^{O(1)} \delta^2 (\MAX_N)^2$.
Combining with Lemma~\ref{l:D.barD.discrep} gives
	\[
	\bE[\Delta D^{\RomI}(q_d)^2]
	\le L^{O(1)}\bigg\{\eta^2 + \delta^2 (\MAX_N)^2 \bigg\}
	\le L^{O(1)} \delta^2 (\MAX_N)^2\,,
	\]
since $\MAX_N\ge1$ and $\eta\ll\delta$ by Assumption~\ref{a:params}.
Similarly, it follows from 
\eqref{e:deriv.p.bound} and 
Lemmas \ref{l:barD.II.second.quenched}--\ref{l:barD.Ising.second.quenched} that
	\[\sum_{\sigma\in\{\RomII,\RomIII,\Ising\}}
	\max\bigg\{
	\bE[\Delta\bar{D}^\sigma(q_d)^2]
	: 0 \le d \le \dmax-1\bigg\}
	\le L^{O(1)}
	\frac{(\SMAX_N)^2\BUCKETS\delta}{N}\,,\]
where the random variable $\SMAX_N$ is defined by
	\beq\label{e:SMAX.N}
 	\SMAX_N\equiv\max\bigg\{\SMAX_N(\RomII),\SMAX_N(\RomIII),\SMAX_N(\Ising)\bigg\}
		\,.\eeq
Combining with Lemma~\ref{l:D.barD.discrep} gives
	\[\sum_{\sigma\in\{\RomII,\RomIII,\Ising\}}
	\max\bigg\{
	\bE[\Delta D^\sigma(q_d)^2]
	: 0 \le d \le \dmax-1\bigg\}
	\le L^{O(1)}
	\bigg\{
	\eta^2 + \frac{(\SMAX_N)^2\BUCKETS\delta}{N}
	\bigg\}\,.
	\]
The claim then follows by
the same argument as in the
\hyperlink{proof:p.Y.kolmogorov}{proof of
\eqref{e:Y.I.Kolmogorov} in Proposition~\ref{p:Y.kolmogorov}},
using the estimates above
 for the individual intervals $[q_d,q_{d+1}]$.
\end{proof}

\subsection{Moment estimates for quadratic variation}
\label{ss:kolmogorov.qv} 
In this subsection we define and analyze the ``quadratic variation'' of $\vY$ and $\vY^\frozen$. \textbf{The main result of this subsection is Proposition~\ref{p:qv.kolmogorov} below. In addition we prove some quenched estimates, in particular Lemma~\ref{l:tQ.first.quenched}
and Proposition~\ref{p:Y.notI.second.quenched}, which will be used in Section~\ref{s:sde}.}

\begin{dfn}[discrete quadratic variation processes]
\label{d:Q}
Let $\vY$ and $\vY^\frozen$ be the spatial rerandomization processes from Definitions~\ref{d:rerand} and \ref{d:rerand.Ising}. Let $\vD$ and $\vD^\frozen$ be the associated drift processes from Definition~\ref{d:D}.
We now let 
	\begin{align*}
	\vZ &\equiv \vY-\vD\,,\\
	\vZ^\frozen 
		&\equiv \vY^\frozen-\vD^\frozen\,,
	\end{align*}
so these are martingales with respect to the filtration $\cG(q_d)$ from \eqref{e:gaus.filt}. In this subsection we consider the ``quadratic variations'' (and ``covariations'') for the $Z$ processes. To this end, for $\sigma,\tau\in\{\RomI,\RomII,\RomIII,\Ising\}$, define 
$Q^{\sigma,\tau}$ to be the piecewise linear process with increments
	\beq\label{e:qv.covar}
	\Delta Q^{\sigma,\tau}(q_d) 
	=  \bE \bigg[
		\Big( \Delta Z^\sigma(q_d)\Big)
		\Big(\Delta Z^\tau(q_d)\Big)
		\,\bigg|\, \cG(q_d)
		\bigg]\,,\eeq
started from $Q^{\sigma,\tau}(q_0)=0$. Clearly $Q^{\sigma,\tau}=Q^{\tau,\sigma}$; and we will also abbreviate $Q^{\sigma}\equiv Q^{\sigma,\sigma}$. We define likewise $Q^{\frozen,\sigma,\tau}$.
\end{dfn}

Note that for any $s=q_a < q_b=t$, 
the bound \eqref{e:Y.II.III.Is.Kolmogorov} from Proposition~\ref{p:Y.kolmogorov},
combined with the bound
\eqref{e:D.naive.kolmogorov.bound}, readily gives
\begin{align} \nonumber
&\E\bigg[\Big( Q^{\tau,\tau'}(t)-Q^{\tau,\tau'}(s)\Big)^4\bigg]
=\E\bigg[ \Big( 
	\bE[ 
	(Z^\tau(t)-Z^\tau(s))
	(Z^{\tau'}(t)-Z^{\tau'}(s))
	\,|\,\cG(s)]
	\Big)^4\bigg]
	\\
&\qquad \le O(1) 
	\max_{\sigma\in\{\RomII,\RomIII,\Ising\}}
\E\bigg[ \Big(Z^\sigma(t)-Z^\sigma(s) \Big)^8 \bigg]
\nonumber \\
&\qquad\le  O(1)\max_{\sigma\in\{\RomII,\RomIII,\Ising\}}
	\bigg\{
	\E\bigg[ \Big(Y^\sigma(t)-Y^\sigma(s) \Big)^8 \bigg]
	+\E\bigg[ \Big(D^\sigma(t)-
	D^\sigma(s) \Big)^8 \bigg]\bigg\}
\le L^{O(1)} (t-s)^4 
\label{e:Q.naive.kolmogorov}
\end{align} 
for all $\tau,\tau'\in\{\RomII,\RomIII,\Ising\}$. Next, for further analysis, we define an approximate version of $Q$ that will be easier to work with:

\begin{dfn}[approximate quadratic variation processes]
\label{d:tQ}
For $\sigma,\tau\in\{\RomI,\RomII,\RomIII\}$, define $\tilde{Q}^{\sigma,\tau}$ to be the piecewise linear process with increments
	\beq\label{e:tQ.notIs}
	\Delta\tilde{Q}^{\sigma,\tau}(q_d)
	\equiv
	\Delta\tilde{Q}^{\sigma,\tau}(q_d,B)
	\equiv
	\frac{1}{|B|}
	\sum_{a\in B}
	\Big(\Delta X^{\trunc,\sigma}(q_d,a)\Big)
	\Big(\Delta X^{\trunc,\tau}(q_d,a)\Big)\,.
	\eeq
For the Ising case, we also define
	\beq\label{e:tQ.notIs.Is}
	\Delta\tilde{Q}^{\sigma,\Ising}(q_d)
	\equiv \frac{1}{|B||B^{\Ising}|}
	\sum_{a\in B,i\in B^{\Ising}}
	\Big(\Delta X^{\trunc,\sigma}(q_d,a)\Big)
	\Big(\Delta X^{\trunc,\Ising}(q_d,i)\Big)
	=
	\Big( \Delta \bar{D}^\sigma(q_d)\Big)
	\Big( \Delta \bar{D}^{\Ising}(q_d)\Big)
	\eeq
for $\sigma\in\{\RomI,\RomII,\RomIII\}$; 
as well as
	\beq\label{e:tQ.Is}
	\Delta\tilde{Q}^{\Ising}(q_d)
	\equiv
	\Delta\tilde{Q}^{\Ising}(q_d,B^\Ising)
	\equiv 
	\frac{1}{|B^{\Ising}|}
	\sum_{i\in B^{\Ising}}
	\Big(\Delta X^{\trunc,\Ising}(q_d,i)\Big)^2\,.
	\eeq
Denote also $E^{\sigma,\tau}\equiv Q^{\sigma,\tau}-\tilde{Q}^{\sigma,\tau}$. We define likewise $\Delta\tilde{Q}^\frozen$ using increments of $X^\frozen$. Then, 
analogously to \eqref{e:barD.vs.plus},
for all $\sigma , \tau \in\{\RomI,\RomII,\RomIII\}$ we have for $B^\frozen \cap [M] = B \neq \emptyset$
	\beq\label{e:tQ.vs.plus}
	\Delta\tilde{Q}^{\frozen,\sigma,\tau}(q_d,B^\frozen)
	=\frac{|B|}{|B^\frozen|}
	\Delta\tilde{Q}^{\sigma,\tau}(q_d,B)
	\,.\eeq
A similar relation holds between $\Delta\tilde{Q}^{\sigma,\Ising}$ and $\Delta\tilde{Q}^{\frozen,\sigma,\Ising}$, and likewise between
$\Delta\tilde{Q}^{\frozen,\Ising}$
and $\Delta\tilde{Q}^{\Ising}$.
\end{dfn}

Recall the definition of the stopping times $q_T$, $q_{\Ising,T}$, $q_{\frozen,T}$, $q_{\frozen,\Ising,T}$ appearing in \S\ref{ss:tightness.statements}. For $q_d\ge q_T$ we can take $\Delta\tilde{Q}^{\sigma,\tau}(q_d)=\Delta E^{\sigma,\tau}(q_d)=0$ for all $\sigma,\tau$. Recall from \eqref{e:def.v} and \eqref{e:def.u} that
	\begin{align*}
	\vvv^\trunc(q_d)
	&\equiv\frac{\Delta\tilde{Q}^{\RomII}(q_d)}{\delta_d p(q_d)}
	= \ind\{q_d<q_T\} 
	\frac{1}{|B|}\sum_{a\in B}
	\frac{(\bar{\bg}^a,\by)^2}{N}\,,\\
	\uuu^\trunc(q_d)
	&\equiv
	\frac{\Delta\tilde{Q}^{\RomII,\RomIII}(q_d)}{\delta_d
	[p(q_d) q_d p'(q_d)]^{1/2}}
	= \ind\{q_d<q_T\} 
	\frac{1}{|B|}\sum_{a\in B}
	\frac{(\bmeta^a,\bar{\bx})
		(\bar{\bg}^a,\by)}{N} \,.
	\end{align*}
Similarly, recall from \eqref{e:def.w} that 
	\[
	\www^\trunc(q_d)
	\equiv \frac{\Delta\tilde{Q}^{\Ising}(q_d)}
		{\delta_d}
	= \ind\{q_d<q_{\Ising,T}\} 
	\frac{1}{|B^{\Ising}|}
	\sum_{i\in B^{\Ising}}
	(\be_i,\by)^2\,.\]
Recalling \eqref{e:def.r}, we further decompose $\tilde{Q}^{\RomIII} = C^{\RomIII} + \tilde{E}^{\RomIII}$, according to the definition 
	\begin{align}
	\label{e:APPX.def.r}
	\rrr^\trunc(q_d)
	&\equiv
	\frac{\Delta C^{\RomIII}(q_d)}{\delta_d}
	\equiv \ind\{q_d<q_T\} 
	q_d p'(q_d) \,,
	\\
	\frac{\Delta \tilde{E}^{\RomIII}(q_d)}{\delta_d}
	&= \ind\{q_d<q_T\} 
	 q_d p'(q_d)
	\bigg\{
	\frac{1}{|B|}
	\sum_{a\in B}\bigg(
	\frac{(\bmeta^a,\bar{\bx})^2}{N}
	-\frac{\|\bar{\bx}\|^2}{N}\bigg)
	+\bigg( \frac{\|\bar{\bx}\|^2}{N}-1\bigg)
	\bigg\}\,.
	\label{e:qv.ezero}
	\end{align}
Recalling \eqref{e:tQ.vs.plus}, we also decompose
$\tilde{Q}^{\frozen,\RomIII} = C^{\frozen,\RomIII} + \tilde{E}^{\frozen,\RomIII}$ where 
	\[
	\rrr^{\frozen}(q_d)
	\equiv
	\frac{\Delta C^{\frozen,\RomIII}}{\delta_d}
	\equiv \frac{|B|}{|B^\frozen|}
		\rrr(q_d)\,.
	\]
\textbf{The quantities 
$(\bbb,\vvv,\uuu,\rrr,\www)^\frozen$
will play an important role in the SDE limit; see Section~\ref{s:sde}.}

We will show in Lemma~\ref{l:qv.simplification} that the processes $Q$ and $\tilde{Q}$ are close. We will show that all the processes $Q^{\RomI,\sigma}$ are asymptotically negligible, as are the processes $Q^{\sigma,\Ising}$ for $\sigma\in\{\RomI,\RomII,\RomIII\}$.  Thus the only processes that survive in the limit are $Q^{\RomII}$, $Q^{\RomII,\RomIII}$, $Q^{\RomIII}$, and $Q^{\Ising}$. We will prove an analogous statement for the $Q^\frozen$ processes. \textbf{The main result of this subsection is the following estimate (analogous to 
Proposition~\ref{p:Y.kolmogorov}
from \S\ref{ss:apriori}, and Proposition~\ref{p:drift.kolmogorov}
from \S\ref{ss:kolmogorov.drift}):} 

\begin{ppn}\label{p:qv.kolmogorov}
The following estimates hold for $N$ large enough. We have
	\beq\label{e:Q.I.small}
	\E Q^{\RomI}(1) \le L^{O(1)} \delta\,.
	\eeq
In the Ising case, we also have
	\beq\label{e:Q.Ising.cross.small}
	\sum_{\sigma\in\{\RomII,\RomIII\}}
	\E \sup\bigg\{ |\tilde{Q}^{\sigma,\Ising}(t)|
	: 0\le t\le 1\bigg\}
	\le
	\frac{L^{O(1)} \BUCKETS}{N}\,,
	\eeq
where we recall $\BUCKETS$ from \eqref{e:buckets.scB}. The errors $E^{\sigma,\tau}$ satisfy 
	\beq\label{e:qv.error.unif.bound}
	\sum_{\sigma,\tau\in\{\RomII,\RomIII,\Ising\}}
	\E \bigg[ \sup\bigg\{ |E^{\sigma,\tau}(t)|
	: 0\le t\le 1\bigg\}\bigg]
	\le
	L^{O(1)} 
		\bigg\{\frac{\eta}{\delta^{1/2}} 
		+ \frac{\BUCKETS}{N}
	\bigg\}
	\,.
	\eeq
Finally, the error $\tilde{E}^{\RomIII}$ satisfies
	\beq\label{e:tildeE.III.unif}
	\E
	\sup\bigg\{
	|\tilde{E}^{\RomIII}(t)| : 0\le t\le 1\bigg\}
	\le
	L^{O(1)} 
	\bigg\{
	\frac{1}{N\delta}
	+ \frac{\BUCKETS}{N}
	\bigg\}^{1/2}
	\eeq
The ``$\frozen$'' variant of this proposition also holds.
\end{ppn}

The \hyperlink{proof:p.qv.kolmogorov}{proof of Proposition~\ref{p:qv.kolmogorov}}
appears at the end of this subsection. It follows from several intermediate lemmas which we turn to next.

\begin{proof}[\hypertarget{proof:l.tQ.first.quenched}{Proof of Lemma~\ref{l:tQ.first.quenched}}]
For $\sigma,\tau\in\{\RomI,\RomII,\RomIII\}$, it follows from the definition \eqref{e:tQ.notIs} that
	\begin{align*}
	&\bE|\Delta\tilde{Q}^{\sigma,\tau}(q_d)|
	=\frac1{M} \sum_{B\in\mathcal{B}(q_d)} 
	\bigg|\sum_{a\in B}
	\Big(\Delta X^{\trunc,\sigma}(q_d,a)\Big)
	\Big(\Delta X^{\trunc,\tau}(q_d,a)\Big)
	\bigg| \\
	&\qquad\le
	\frac{1}{M}\sum_{a=1}^M\bigg|
	\Big(\Delta X^{\trunc,\sigma}(q_d,a)\Big)
	\Big(\Delta X^{\trunc,\tau}(q_d,a)\Big)
	\bigg|\le\bigg\{ \bE\Big[\Delta X^{\trunc,\sigma}(q_d,\hat{a})^2\Big]
	\bE\Big[\Delta X^{\trunc,\tau}(q_d,\hat{a})^2\Big]\bigg\}^{1/2}\,,
	\end{align*}
where $\hat{a}$ is a uniformly random sample from $[M]$ and $\bE$ includes expectation over $\hat{a}$.
It follows by applying Proposition~\ref{p:X.second.mmt.quenched} that
	\[
	\frac{\bE\Delta\tilde{Q}^{\RomI}(q_d)}{(\delta_d)^2 p'(q_d)}
	\le
	\frac{\bE\Delta X^{\trunc,\RomI}(q_d,\hat{a})^2}
		{(\delta_d)^2 p'(q_d)}
	\le (\MAX_N)^2\,.
	\]
We similarly obtain the bounds
	\begin{align*}	 
	\bE\vvv^\trunc(q_d)
	&\stackrel{\eqref{e:def.v}}{=}
	\frac{\bE\Delta\tilde{Q}^{\RomII}(q_d)}{\delta_d p(q_d)}
	\le
	\frac{\bE[ \Delta X^{\trunc,\RomII}(q_d,\hat{a})^2]}
	{\delta_d p(q_d)}
	\le (\MAX_N)^2\,,\\
	\bE|\uuu^\trunc(q_d)|
	&\stackrel{\eqref{e:def.u}}{=}
	\frac{\bE|\Delta\tilde{Q}^{\RomII,\RomIII}|}
	{\delta_d [p(q_d) q_d p'(q_d)]^{1/2}}
	\le
	\bigg\{
	\frac{\bE[\Delta X^{\trunc,\RomII}(q_d,\hat{a}^2)]}{p(q_d) }
	\frac{\bE[\Delta X^{\trunc,\RomIII}(q_d,\hat{a}^2)]}{q_d  p'(q_d)}
	\bigg\}^{1/2}
	\le (\MAX_N)^2\,,\\
	\bE\www^\trunc(q_d)
	&\stackrel{\eqref{e:def.w}}{=}
	\frac{\bE\Delta\tilde{Q}^{\Ising}(q_d)}{\delta_d}
	\le\frac{\bE[\Delta X^{\trunc,\Ising}(q_d,\hat{\iota})^2]}{\delta_d}
	\le(\MAX_N)^2\,.
	\end{align*}
The same bound applies for the $\frozen$ variant of the processes, since Proposition~\ref{p:X.second.mmt.quenched} also applies to $\vX^\frozen$.
\end{proof}

We have the following comparison between $Q$ and $\tilde{Q}$:

\begin{lem}[error between $Q$ and $\tilde{Q}$]
\label{l:qv.simplification}
Let $Q^{\sigma,\tau}=E^{\sigma,\tau}+\tilde{Q}^{\sigma,\tau}$ as given by \eqref{e:qv.covar} and \eqref{e:tQ.notIs}. Then
	\[\Big|
	\Delta E^{\sigma,\tau}(q_d)\Big| 
	\le O(1)\bigg\{
	\eta^2
	+
	\sum_{\rho\in\{\sigma,\tau\}}
	\bigg[
	\eta\bigg(
	\min\Big\{
	\Delta Q^\rho(q_d),
	\Delta\tilde{Q}^\rho(q_d)\Big\}
	\bigg)^{1/2}
	+ \Big( \Delta \bar{D}^\rho(q_d)\Big)^2
	\bigg]
	\bigg\}\]
for all $\sigma,\tau\in\{\RomI,\RomII,\RomIII,\Ising\}$.
The same bound holds with $Q^\frozen$, $\tilde{Q}^\frozen$, $E^\frozen$
in place of $Q$, $\tilde{Q}$, $E$.

\begin{proof} Recall that $\vZ\equiv \vY-\vD$, so we can rewrite \eqref{e:qv.covar} as
	\[
	\Delta Q^{\sigma,\tau}(q_d)
	= \bE \bigg[
		\Big( \Delta Y^\sigma(q_d)
			-\Delta D^\sigma(q_d)
			\Big)
		\Big(
		\Delta Y^\tau(q_d)
		-\Delta D^\tau(q_d)
		\Big)
		\,\bigg|\, \cG(q_d)
		\bigg]\,.
	\]
For simplicity of exposition, let us consider the case $\sigma=\tau=\RomI$.
Similarly to the argument for
Lemma~\ref{l:D.barD.discrep},
let $B=B(q_d,a(q_d))$ denote the bucket at time $q_d$ that contains $a(q_d)$. Lemma~\ref{l:X.Y.discrep} then tells us that $\Delta Y^\RomI(q_d)$ is within $2\eta$ of $\Delta X^{\trunc,\RomI}(q_d,a(q_{d+1}))$, where
$a(q_{d+1})$ is a uniformly random sample from
$B$. Thus define $\bar{Q}^{\RomI}$ to be the process with increments
	\beq\label{e:qv.process.simplified}
	\Delta\bar{Q}^{\RomI}(q_d)
	= \frac{1}{|B|}
	\sum_{a\in B}
	\Big(
	\Delta X^{\trunc,\RomI}(q_d,a)-\Delta D^\RomI(q_d)
	\Big)^2
\,.
	\eeq
 It follows from Lemma~\ref{l:X.Y.discrep} that
	\[
	\Delta Q^{\RomI}(q_d)
	=\frac{1}{|B|}
	\sum_{a\in B}
	\bigg\{
	O(\eta)
	+\Delta X^{\trunc,\RomI}(q_d,a) 
	- \Delta D^{\RomI}(q_d)
	\bigg\}^2\,,\]
which allows us to bound
	\begin{align}\nonumber
	&\Big|\Delta Q^{\RomI}(q_d)
	-\Delta\bar{Q}^{\RomI}(q_d)\Big| 
	\le O(\eta^2)
	+ O(\eta)
	\frac{1}{|B|}
	\sum_{a\in B}
	\bigg|
	\Delta X^{\trunc,\RomI}(q_d,a) 
	- \Delta D^{\RomI}(q_d)\bigg|
	\\ \nonumber
	&\qquad\le O(\eta^2)
	+ O(\eta) 
	\bigg[\frac{1}{|B|}
	\sum_{a\in B}
	\Big(
	\Delta X^{\trunc,\RomI}(q_d,a) 
	- \Delta D^{\RomI}(q_d) \Big)^2\bigg]^{1/2}\\
	&\qquad\le
	O(1)\bigg\{ \eta^2+\eta 
		\Big( \Delta\bar{Q}^{\RomI}(q_d) \Big)^{1/2}\bigg\}\,.
	\label{e:Q.to.barQ}
	\end{align}
This shows that increments of
$Q^{\RomI}$ and $\bar{Q}^{\RomI}$ are close. We can further bound
	\begin{align}\nonumber
	&\Big|
	\Delta\bar{Q}^{\RomI}(q_d)
	-\Delta\tilde{Q}^{\RomI}(q_d)\Big| 
	=\bigg| \frac{1}{|B|}
	\sum_{a\in B}
	\bigg\{
	\Delta X^{\trunc,\RomI}(q_d,a)-\Delta D^{\RomI}(q_d)
	\bigg\}^2
	- \frac{1}{|B|}\Big(\Delta X^{\trunc,\RomI}(q_d,a)\Big)^2\bigg|
	\\ \nonumber
	&\qquad=\bigg|
	\Big(\Delta D^{\RomI}(q_d)\Big)^2
	- 2
	\Big( \Delta D^{\RomI}(q_d)\Big)
	\frac{1}{|B|}
	\sum_{a\in B}
	\Delta X^{\trunc,\RomI}(q_d,a)
	\Big)
	\bigg|\\ \nonumber
	&\qquad=\bigg|
	\Big(\Delta D^{\RomI}(q_d)\Big)^2
	-2\Big(\Delta D^{\RomI}(q_d)\Big)
	\Big(\Delta \bar{D}^{\RomI}(q_d)\Big)
	\bigg|\\
	&\qquad\le O(1)  
	\bigg\{
	\Delta \bar{D}^{\RomI}(q_d)^2
	+ \eta\Big|\Delta \bar{D}^{\RomI}(q_d)\Big| + \eta^2\bigg\}
	\le O(1)
	\max\bigg\{\Delta \bar{D}^{\RomI}(q_d)^2,\eta^2\bigg\}
	\label{e:barQ.to.tildeQ}
	\,,\end{align}
having used the definition of $\bar{D}$ from \eqref{e:barD}, together with Lemma~\ref{l:D.barD.discrep}. Combining the bounds gives
	\[
	\Big|
	\Delta E^{\RomI}(q_d)\Big| 
	\le  O(1)
	\bigg\{ \eta^2 
	+\eta\Big( \Delta\tilde{Q}^{\RomI}(q_d) \Big)^{1/2}
	+\Big(\Delta \bar{D}^{\RomI}(q_d)\Big)^2
	\bigg\}\,.
	\]
Next note that \eqref{e:barQ.to.tildeQ} implies
	\[\Delta\tilde{Q}^{\RomI}(q_d)
	\le O(1)\max\bigg\{ 
	\Delta \bar{Q}^{\RomI}(q_d),
	\Delta \bar{D}^{\RomI}(q_d)^2,
	\eta^2
	\bigg\}\,.\]
Meanwhile \eqref{e:Q.to.barQ} implies that either 
$\Delta\bar{Q}^{\RomI}(q_d) \le O(\eta^2)$, in which case
$\Delta Q^{\RomI}(q_d) \le O(\eta^2)$ also,
or $\Delta\bar{Q}^{\RomI}(q_d)$ is much larger than $\eta^2$, in which case it is of the same order as
 $\Delta Q^{\RomI}(q_d)$. It follows that
	\[\Delta\tilde{Q}^{\RomI}(q_d)
	\le O(1)\max\bigg\{ 
	\Delta Q^{\RomI}(q_d),
	\Delta \bar{D}^{\RomI}(q_d)^2,
	\eta^2
	\bigg\}\,.\]
Combining these inequalities gives
	\[
	\Big|
	\Delta E^{\RomI}(q_d)\Big| 
	\le  O(1)
	\bigg\{ \eta^2 
	+\eta\bigg(
	\min\Big\{ 
	\Delta Q^{\RomI}(q_d) ,
	\Delta\tilde{Q}^{\RomI}(q_d) 
	\Big\}
	\bigg)^{1/2}
	+\Big(\Delta \bar{D}^{\RomI}(q_d)\Big)^2
	\bigg\}\,,
	\]
which proves the claimed bound for $\Delta Q^{\RomI}(q_d)$.
The general bound for
$\Delta Q^{\sigma,\tau}(q_d)$ follows by a very similar calculation.
\end{proof}
\end{lem}

\begin{cor}[quenched first moment bounds on $Q^\sigma$ and $Q^{\sigma,\frozen}$]\label{c:Q.first.quenched}
For $Q$ as in Definition~\ref{d:Q}, we have
	$\bE[\Delta Q^\RomI(q_d)]
	\le L^{O(1)}(\MAX_N)^2\delta^2$
with $\MAX_N$ as in Proposition~\ref{p:Y.kolmogorov}. We also have
	\[
	\sum_{\sigma\in\{\RomII,\RomIII,\Ising\}}
	\bE\Delta Q^\sigma(q_d)
	\le L^{O(1)}\bigg\{
	(\MAX_N)^2
	+ \frac{(\SMAX_N)^2 \BUCKETS}{N}
	\bigg\} \delta
	\]
with $\SMAX_N$ as in Proposition~\ref{p:drift.quenched}. The same bound holds for the $Q^{\frozen,\sigma}$ processes.
	
\begin{proof}
Recall that $Q^\sigma=E^\sigma+\tilde{Q}^\sigma$. It follows from Lemma~\ref{l:qv.simplification} that
	\[
	\Big|\Delta E^\sigma(q_d)\Big| 
	\le O\bigg(
	\eta^2
	+\eta\Big(\Delta\tilde{Q}^\sigma(q_d)\Big)^{1/2}
	+\Big( \Delta \bar{D}^\sigma(q_d)\Big)^2\bigg)\,.
	\]
For $\sigma=\RomI$, it follows by combining Lemma~\ref{l:barD.I.second.quenched} with Lemma~\ref{l:tQ.first.quenched} that
	\[
	\bE\Delta Q^\RomI(q_d)
	\le \bE\Delta \tilde{Q}^\RomI(q_d)
		+\bE\Big|\Delta E^\RomI(q_d)\Big|
	\le
	L^{O(1)} \bigg\{
	\eta^2 + (\MAX_N)^2\delta^2 \bigg\}\,.
	\]
For $\sigma\in\{\RomII,\RomIII,\Ising\}$,
it follows by combining Lemmas~\ref{l:barD.II.second.quenched}--\ref{l:barD.Ising.second.quenched} with Lemma~\ref{l:tQ.first.quenched} that
	\[
	\bE\Delta Q^\sigma(q_d)
	\le\bE\Delta \tilde{Q}^\sigma(q_d)+ \bE\Big|\Delta E^\sigma(q_d)\Big| 
	\le L^{O(1)}\bigg\{
	\eta^2
	+(\MAX_N)^2\delta
	+ \frac{(\SMAX_N)^2 \BUCKETS\delta}{N}
	\bigg\}\,.\]
This proves the claim for the $Q^\sigma$ processes,
using that $\MAX_N\ge1$ and $\eta\ll\delta$. The analogous bound for the $Q^{\frozen,\sigma}$ processes follows by the same argument.
\end{proof}
\end{cor}

We next prove Proposition~\ref{p:Y.notI.second.quenched}, which should be compared with 
the bound \eqref{e:Y.I.Kolmogorov} from Proposition~\ref{p:Y.kolmogorov}.

\begin{proof}[\hypertarget{proof:p.Y.notI.second.quenched}{Proof of Proposition~\ref{p:Y.notI.second.quenched}}]
For simplicity of exposition we assume $s=q_a \le q_b =t$; the general case follows by arguing as in the
\hyperlink{proof:p.Y.kolmogorov}{proof of Proposition~\ref{p:Y.kolmogorov}}.
It follows from Definitions~\ref{d:D} and \ref{d:Q} that for $q_a \le q_b$, we have 
    \begin{align*}
        \bE\bigg[\Big(Y^\sigma(q_b)-Y^\sigma(q_a)\Big)^2\bigg]
        &\le 2\bigg\{
    	   \bE\bigg[\Big(Z^\sigma(q_b)-Z^\sigma(q_a)\Big)^2\bigg]
    	   + \bE\bigg[\Big(D^\sigma(q_b)-D^\sigma(q_a)\Big)^2\bigg]
    	\bigg\} \\
        &= 2\bigg\{
    	   \bE\Big[ Q^\sigma(q_b)-Q^\sigma(q_a)\Big]
    	   + \bE\bigg[\Big(D^\sigma(q_b)-D^\sigma(q_a)\Big)^2\bigg]
    	\bigg\}.
    \end{align*}
where the second step uses that $Z^\sigma$ is a martingale. 
For $\sigma\in\{\RomII,\RomIII,\Ising\}$, recall from
Proposition~\ref{p:drift.quenched} that
	\[\bE\Big[( D^\sigma(t)- D^\sigma(s))^2\Big]
	\le L^{O(1)} (t-s)^2
	\bigg\{ 
	\frac{\eta^2}{\delta^2}+
	\frac{(\SMAX_N)^2\BUCKETS}{N\delta}\bigg\}\,.\]
On the other hand, Corollary~\ref{c:Q.first.quenched} gives
	\[
	\bE\Big[ Q^\sigma(t)-Q^\sigma(s)\Big]
	\le L^{O(1)} (t-s)
	\bigg\{
	(\MAX_N)^2
	+ \frac{(\SMAX_N)^2 \BUCKETS}{N}
	\bigg\}\,.
	\]
Combining these bounds proves the claim, recalling that $\MAX_N\ge1$ and $\eta\ll\delta$ from Assumption~\ref{a:params}. The analogous bound for $\vY^\frozen$ follows by the same argument.
\end{proof}

\begin{lem}[annealed first moment bounds for 
$E^{\sigma,\tau}$ and $E^{\frozen,\sigma,\tau}$
with $\sigma,\tau\in\{\RomII,\RomIII,\Ising\}$]
\label{l:Q.vs.tilde.Q}
Let
$Q^{\sigma,\tau}=E^{\sigma,\tau}+\tilde{Q}^{\sigma,\tau}$ as in Definitions~\ref{d:Q} and \ref{d:tQ}. Then, for $N$ large enough, we have
	\[\bbE
	\Big|\Delta E^{\sigma,\tau}(q_d)\Big| 
	\le L^{O(1)}\bigg\{\eta\delta^{1/2} + \frac{\BUCKETS\delta}{N}
	\bigg\}
	\]
for all $\sigma,\tau\in\{\RomII,\RomIII,\Ising\}$. The same bound holds with $E^\frozen$ in place of $E$.
\begin{proof}
This is similar
to the proof of Corollary~\ref{c:Q.first.quenched}. Recall 
that Lemma~\ref{l:qv.simplification} implies
	\[
	\Big|\Delta E^{\sigma,\tau}(q_d)\Big| 
	\le O\bigg(
	\eta^2
	+\sum_{\rho\in\{\sigma,\tau\}}\bigg[
	\eta\Big(\Delta\tilde{Q}^\rho(q_d)\Big)^{1/2}
	+\Big( \Delta \bar{D}^\rho(q_d)\Big)^2
	\bigg]\bigg)\,.
	\]
By Lemmas~\ref{l:barD.II.second}--\ref{l:barD.Ising.second}, for all $\sigma\in\{\RomII,\RomIII,\Ising\}$ we have
	\[\E\Big[\Delta\bar{D}^\sigma(q_d)^2\Big]
	\le\frac{ L^{O(1)}\BUCKETS\delta}{N}
	\,.\]
By Lemma~\ref{l:tQ.first.quenched} and the bound $\E[(\MAX_N)^2]\le L^{O(1)}$ from Proposition~\ref{p:Y.kolmogorov}, for all $\sigma\in\{\RomII,\RomIII,\Ising\}$ we have
	\[
	\E\bigg[\Big(\Delta\tilde{Q}^\sigma(q_d)\Big)^{1/2}\bigg]
	\le L^{O(1)} \delta^{1/2}\,.
	\]
Combining the bounds gives the claim, recalling that $\eta\ll\delta$ from Assumption~\ref{a:params}.
\end{proof}
\end{lem}

\begin{lem}[annealed first moment bounds for 
$\tilde{Q}^{\sigma,\Ising}$ and $\tilde{Q}^{\frozen,\sigma,\Ising}$
with $\sigma\in\{\RomII,\RomIII\}$]
\label{l:tilde.Q.sigma.Ising}
Recalling \eqref{e:tQ.notIs.Is}, we have
	\[\E \Big|\Delta\tilde{Q}^{\sigma,\Ising}(q_d)\Big|
	\le\frac{ L^{O(1)} \BUCKETS \delta}
		{N}\,.
	\]
for all $\sigma\in\{\RomII,\RomIII\}$.
The same bound holds with $\tilde{Q}^{\frozen,\sigma,\Ising}$
in place of $\tilde{Q}^{\sigma,\Ising}$.

\begin{proof}
Combining \eqref{e:tQ.notIs.Is} with the Cauchy--Schwarz inequality gives
	\[
	\E \Big|\Delta\tilde{Q}^{\sigma,\Ising}(q_d)\Big|
	\le
	\E\bigg[ \Big(\Delta\bar{D}^{\sigma}(q_d)\Big)^2\bigg]^{1/2}
	\E\bigg[ \Big(\Delta\bar{D}^{\Ising}(q_d)\Big)^2\bigg]^{1/2}\,.
	\]
The claim follows by applying Lemmas~\ref{l:barD.II.second}--\ref{l:barD.Ising.second}.
\end{proof}
\end{lem}	

\begin{lem}[annealed second moment bound for $\tilde{E}^{\RomIII}$]
\label{l:mmt.ezero}
For $\tilde{E}^{\RomIII}$ as defined by \eqref{e:qv.ezero},
for $N$ large enough we have
	\[
	\E\bigg[
	\bigg(\frac{\Delta\tilde{E}^{\RomIII}(q_d)}{\delta_d}
	\bigg)^2\bigg]
	\le
	L^{O(1)} 
	\bigg\{
	\frac{1}{N\delta}
	+ \frac{\BUCKETS}{N}
	\bigg\}
	\]
for all $0\le d\le \dmax-1$.
 The same bound holds with $\tilde{E}^{\frozen,\RomIII}$ in place of $\tilde{E}^{\RomIII}$.

\begin{proof}
For $q_d < q_T$ we rewrite \eqref{e:qv.ezero} as
	\beq\label{e:qv.ezero.rewrite}
	\frac{\tilde{E}^{\RomIII}(q_{d+1})-\tilde{E}^{\RomIII}(q_d)}{\delta_d}
	=
	\frac{q_d(p_{d+1}-p_d)}{\delta_d}
	\Big\{
	F+G
	\Big\}\,,\eeq
for $F$ and $G$ defined by 
	\[
	F
	\equiv \frac{1}{|B|}
	\sum_{a\in B}\bigg(
	\frac{(\bmeta^a,\bar{\bx})^2}{N}
	-\frac{\|\bar{\bx}\|^2}{N}\bigg)
	\equiv\frac{1}{|B|}
	\sum_{a\in B} \xi^a
	\,,\quad
	G
	\equiv \frac{\|\bar{\bx}\|^2}{N}-1\,.\]
Note that $\xi_a$ has mean zero conditional on $\cG(q_d)$. It follows that
	\begin{align*}
	\E[F^2]
	&=\E \bigg[\sum_B \frac{|B|}{M}
	\frac{1}{|B|^2}
	\bigg(\sum_{a\in B} \xi_a\bigg)^2\bigg]
	=\E \bigg[\sum_B \frac{1}{{M}|B|}
	\sum_{a\in B}
	\E(\xi_a^2\,|\,\cG(q_d))
	\bigg] \\
	&\le O(1)\,
	\E\bigg[ \frac{|\mathcal{B}(q_d)|}{M}
	{\frac{\|\bar{\bx}\|^4}{N^2}}\bigg]
	\le \frac{ L^{O(1)} \E[ |\mathcal{B}(q_d)|^2]^{1/2}}{N\alpha}
    \le \frac{\BUCKETS L^{O(1)}}{N}
	 \,,
	\end{align*}
where the second-to-last bound uses the Cauchy--Schwarz inequality together with Lemmas~\ref{l:sg.norm.bound} and \ref{l:bar.x.subgaus}. To bound the second moment of $G$, recall that in the case $q_d=0$, we define $\bar{\bx}=N^{1/2}\ind\in\R^N$, so we have simply $G=0$ in this case. Otherwise, if $q_d\ne0$, then we have
	\[
	\|\bar{\bx}\|^2
	\stackrel{\eqref{e:normalized.vectors.x}}{=}
	\frac{\|\bx(q_d)\|^2}{q_d}
	\le \frac{N}{q_d}
	\le \frac{N L^{O(1)}}{\delta}\,,
	\]
where the last bound is by \eqref{e:delta.cutoffs} and \eqref{e:deriv.p.bound}. Thus, if we regard $G$ as a function of the gaussian disorder $\bg\equiv\bar{\bG}(q_d)
\equiv \bG(q_d)/(p_d)^{1/2}$, we can bound
	\begin{align*}
	&|G(\bg)-G(\bg')|
	\le \bigg|
	\frac{(\bar{\bx}(\bg),\bar{\bx}(\bg)-\bar{\bx}(\bg') )}{N}
	\bigg|
	+\bigg|\frac{(\bar{\bx}(\bg)-\bar{\bx}(\bg'),\bar{\bx}(\bg'))}{N}
	\bigg| \\
	&\qquad\le
	\frac{\|\bar{\bx}(\bg)\| \|\bar{\bx}(\bg)-\bar{\bx}(\bg')\|}{N}
	+\frac{\|\bar{\bx}(\bg)-\bar{\bx}(\bg')\|
		\|\bar{\bx}(\bg')\|}{N}
	\le \frac{L^{O(1)}}{N} \bigg(\frac{N}{\delta}\bigg)^{1/2}
		\,.
	\end{align*}
Therefore $G$ is subgaussian with variance proxy $L^{O(1)} /(N\delta)$, and it follows that
	\[\E[G^2]
	\le \frac{L^{O(1)}}{N\delta}\,.\]
Substituting the bounds into \eqref{e:qv.ezero.rewrite}, and combining with \eqref{e:deriv.p.bound}, gives the claim.
\end{proof}
\end{lem}

\begin{proof}[\hypertarget{proof:p.qv.kolmogorov}{Proof of Proposition~\ref{p:qv.kolmogorov}}] 
We collect the estimates provided above for the individual intervals $[q_d,q_{d+1}]$:
\begin{itemize}
\item \eqref{e:Q.I.small} follows by the first moment bound of
Corollary~\ref{c:Q.first.quenched} for $Q^{\RomI}$, combined with the annealed moment bound $\E[(\MAX_N)^2]\le L^{O(1)}$
from Proposition~\ref{p:Y.kolmogorov};
\item \eqref{e:Q.Ising.cross.small} follows by
 the first moment bounds of Lemma~\ref{l:tilde.Q.sigma.Ising}
for $\tilde{Q}^{\sigma,\Ising}$ with $\sigma\in\{\RomII,\RomIII\}$;
\item \eqref{e:qv.error.unif.bound} follows by
the first moment bounds of Lemma~\ref{l:Q.vs.tilde.Q}
for $E^{\sigma,\tau}$ with $\sigma,\tau\in\{\RomII,\RomIII,\Ising\}$;
\item \eqref{e:tildeE.III.unif} follows by
the second moment bound
from 
Lemma~\ref{l:mmt.ezero} for $\tilde{E}^{\RomIII}$.
\end{itemize}
The ``$\frozen$'' variant of the statement follows by the same argument.
\end{proof}

\subsection{Tightness}\label{ss:tightness}

In this subsection we combine the estimates from the previous subsections
\S\ref{ss:kolmogorov.Y}--\S\ref{ss:kolmogorov.qv}
to prove tightness of the $Y$, $D$, and $Q$ processes in the limiting regime where $N\to\infty$, $\eta\to0$, and $\delta\to0$ with
$1/N \ll \eta \ll \delta$, as discussed in Assumption~\ref{a:params}.
\textbf{We do not yet send $\trK\to\infty$.} 
We give the proof for the Ising case, since the proof for the spherical case is only simpler. Recall that we have parameters $p_d,q_d$ as in  \eqref{e:p.q}; these parameters can depend on $N$. As in \eqref{e:gaus.decomp} we have
	\[
	\bG^N(p_d)
	\equiv (p_0)^{1/2}\bXi^{N,0} 
	+\sum_{\ell=1}^d (p_\ell-p_{\ell-1})^{1/2} \bXi^{N,\ell}\,,
	\]
where $\bXi^{N,\ell}$ are i.i.d.\ $M\times N$ matrices, each with i.i.d.\ standard gaussian entries.  Using the gaussian disorder $\bG^N$, define the processes 
	\[\Gamma^N
	\equiv \Big(
	\bar{D}^{\RomI},
	\tilde{Q}^{\RomII},
	\tilde{Q}^{\RomII,\RomIII},
	C^{\RomIII}, \tilde{Q}^{\Ising}
	\Big)\,,\quad
	\Upsilon^N
	\equiv \Big(Z^{\RomII},Z^{\RomIII}, Z^{\Ising}\Big)\]
as well as the error process
	\[\mathcal{E}^N
	\equiv
	|Z^{\RomI}|
	+\sum_{\sigma\in\{\RomI,\RomII,\RomIII,\Ising\}}
		|D^{\sigma}-\bar{D}^{\sigma}|
	+\sum_{\sigma\in\{\RomII,\RomIII,\Ising\}}
		|\bar{D}^\sigma|
	+\sum_{\sigma\in\{\RomII,\RomIII\}}
		|\tilde{Q}^{\sigma,\Ising}|
	+\sum_{\sigma,\tau\in\{\RomII,\RomIII,\Ising\}}
		|E^{\sigma,\tau}|
	+ |\tilde{E}^{\RomIII}|\,.
	\]
These are continuous processes defined for $q_0 \le q \le q_{\dmax}=1$,
and we now extend them to be simply constant on $0\le q\le q_0$ so that they are continuous for all $0\le q \le 1$. 

As before, we let $\bP\equiv \bP_{\bG^N}$ denote
the probability measure conditional on the disorder matrices
$\bG^N\equiv(\bG^N(q_d) : 0\le d\le \dmax)$.
We denote the law of the process $(\Gamma^N,\Upsilon^N,\mathcal{E}^N)$ under $\bP_{\bG^N}$ as
	\beq\label{e:law.of.path.Gamma}
	\mu_{\bG^N}(\cdot)
	\equiv
	\bP_{\bG^N}\bigg((\Gamma^N,\Upsilon^N,\mathcal{E}^N)
	\in \cdot\bigg)\,,\eeq
so this is a $\bG^N$-measurable random probability measure over the space of continuous paths $C([0,1],\R^{9})$ (equipped with the uniform topology). We then denote the law of this random measure as
	\beq\label{e:pm.on.pm.over.paths.Gamma}
	\Q_N(\cdot)
	\equiv
	\P\Big(
	\mu_{\bG^N}
	\in\cdot\Big)\,,
	\eeq
so that $\Q_N$ is a deterministic element of $\mathscr{P}(\mathscr{P}(C([0,1],\R^{9})))$.

\textbf{We now prove
Proposition~\ref{p:tightness}, which shows tightness for the family of measures $\{\Q_N : N\ge0\}$.} This implies weak convergence along subsequences (as $N\to\infty$) to limiting probability measures.
Recalling Definition~\ref{d:froze}, we also have the ``$\frozen$'' variants of all of the above objects, corresponding to adding a positive $1/\trK^8$ density of frozen particles on $[-2\trK,2\trK]^4$. We will show that all the same tightness statements hold for the ``$\frozen$'' variant. \textbf{We show in Section~\ref{s:sde} that the limit of the ``$\frozen$'' variant can be characterized by a well-behaved SDE.}

\begin{ppn}[subsequential limits]\label{p:tightness}
Let $N\to\infty$ with parameters satisfying Assumption~\ref{a:params}. Then the measures $\Q_N$, as defined by \eqref{e:pm.on.pm.over.paths.Gamma}, are tight. Moreover, for any subsequential limit $\Q$, if $\mu\sim\Q$ and $(\Gamma,\Upsilon,\mathcal{E})\sim\mu$, then $\mathcal{E}\equiv0$. The ``$\frozen$'' variant of the statement also holds.

\begin{proof}
Write $\E_N$ for expectation with respect to $\Q_N$. We use Assumption~\ref{a:params} throughout. We will argue that $\Gamma^N$ and $\Upsilon^N$ are asymptotically equivalent to
	\[
	\bar{\Gamma}^N
	\equiv \Big(
	\bar{D}^{\RomI},
	Q^{\RomII},Q^{\RomII,\RomIII},C^{\RomIII}, Q^{\Ising}
	\Big)\,,\quad
	\bar{\Upsilon}^N
	\equiv\Big( 
	Y^{\RomII},Y^{\RomIII},Y^{\Ising}
	\Big)\,.
	\]
For the purposes of this proof, it is more convenient to treat $\bar{D}^{\RomI}$ separately, and group together the remaining processes as
	\[
	\bar{f}^N
	\equiv
	\Big(
	Q^{\RomII},Q^{\RomII,\RomIII},C^{\RomIII}, Q^{\Ising},
		Y^{\RomII},Y^{\RomIII},Y^{\Ising}
	\Big)
	\,.\]
By Proposition~\ref{p:drift.kolmogorov}, for $N$ large enough we have
	\beq\label{e:kolmogorov.bar.Gamma}
	\frac{\E_N[\|
		\bar{D}^{\RomI,N}(t) -
		\bar{D}^{\RomI,N}(s)
		\|^2]}{(t-s)^2}
	\le L^{O(1)}
	\le \mathscr{C}\,.
	\eeq
Meanwhile, by Proposition~\ref{p:Y.kolmogorov}, together with \eqref{e:Q.naive.kolmogorov} and the definition~\eqref{e:def.r} of $C^{\RomIII}$, for $N$ large enough we also have 
	\beq\label{e:kolmogorov.bar.Upsilon}
	\frac{\E_N[\|\bar{f}^N(t)-
		\bar{f}^N
		(s)\|^4]}{(t-s)^2}
	\le L^{O(1)}
	\le \mathscr{C}
	\eeq
for $N$ large enough. We next address the error terms, using Assumption~\ref{a:params}. For any function $f$ on $[0,1]$
we use $\|f\|_\infty$ to denote the supremum norm.
\begin{itemize}
\item It follows by the bound \eqref{e:Q.I.small} from Proposition~\ref{p:qv.kolmogorov}, combined with Doob's $L^2$
martingale inequality, that we have
	\[
	\E_N\bigg[ \Big( \| Z^{\RomI}\|_\infty \Big)^2\bigg]
	\le 4 \E_N[Z^{\RomI}(1)^2]
	= 4\E_N Q^{\RomI}(1)
	\le L^{O(1)} \delta
	\ll1\,,
	\]
recalling that we assumed $\delta\to0$ as $N\to\infty$.
\item It follows from Lemma~\ref{l:D.barD.discrep} that
	\[
	\sum_{\sigma\in\{\RomI,\RomII,\RomIII,\Ising\}}
	\Big\|D^\sigma-\bar{D}^\sigma\Big\|_\infty
	\le \frac{L^{O(1)}\eta}{\delta}\ll1\,,
	\]
recalling that we assumed $\eta\ll\delta$.
\item It follows by the bound \eqref{e:bar.D.notI.unif} from Proposition~\ref{p:drift.kolmogorov}
that 
	\[
	\sum_{\sigma\in\{\RomII,\RomIII,\Ising\}}
	\E_N \| \bar{D}^\sigma\|_\infty
	\ll 1\,,
	\]
using that $\BUCKETS \ll N\delta$.
\item It follows by the bounds \eqref{e:Q.Ising.cross.small},
\eqref{e:qv.error.unif.bound}, and \eqref{e:tildeE.III.unif}
from Proposition~\ref{p:qv.kolmogorov}
that
	\[
	\sum_{\sigma\in\{\RomII,\RomIII\}}
		\E_N\|\tilde{Q}^{\sigma,\Ising}\|_\infty
	+\sum_{\sigma,\tau\in\{\RomII,\RomIII,\Ising\}}
		\E_N\|E^{\sigma,\tau}\|_\infty
	+ \E_N\|\tilde{E}^{\RomIII}\|_\infty\ll1\,,
	\]
using that $\eta\ll \delta\ll \delta^{1/2}$ and 
$\BUCKETS \ll N\delta \ll N$.
\end{itemize}
Altogether this proves that $\mathcal{E}^N$ converges in probability to the zero path. By passing to a subsequence we can assume that $\mathcal{E}^N$ converges almost surely to the zero path.
 
Now let $J_n(H)$
be the event that $\bar{D}^{\RomI}$
fails to be $\gamma$-H\"older with constant $H$ at dyadic scale $2^{-n}$, and similarly $F_n(H)$ for the process $\bar{f}$:
	\begin{align*}
	J_n(H)
	&\equiv
	\bigg\{
	\bigg\|
	\bar{D}^{\RomI}
	\bigg(\frac{i}{2^n}\bigg)
	-\bar{D}^{\RomI}
	\bigg(\frac{i-1}{2^n}\bigg)
	\bigg\|
	> \frac{H}{(2^n)^\gamma}
	\textup{ for some }
	1\le i\le 2^n
	\bigg\}\,,\\
	F_n(H)
	&\equiv 
	\bigg\{
	\bigg\|
	\bar{f}
	\bigg(\frac{i}{2^n}\bigg)
	-\bar{f}
	\bigg(\frac{i-1}{2^n}\bigg)
	\bigg\|
	> \frac{H}{(2^n)^\gamma}
	\textup{ for some }
	1\le i\le 2^n
	\bigg\}\,.
	\end{align*}
Finally, the moment estimates above give
\[
	\E_N\big[|Y^{\Ising}(0)|^2\big]\le \mathscr{C}\,.
\]
Let $I(H):=\{|Y^{\Ising}(0)|>H\}$. Then $\P_N(I(H))\le \mathscr{C}/H^2$.
It follows from the estimates 
\eqref{e:kolmogorov.bar.Gamma} and
 \eqref{e:kolmogorov.bar.Upsilon}
(taking a union bound over $1\le i\le 2^n$) that 
for $N$ large enough we have
	\begin{align*}
	\P_N(J_n(H) )
	&\le 2^n 
	\frac{\mathscr{C}
	 (1/2^n)^2}
	 {H^2 / (2^n)^{2\gamma}}
	= \frac{\mathscr{C}
		 (2^n)^{2\gamma}}
		{H^2 2^n}
	\le \frac{\mathscr{C}}
	{H^2 2^{n/2}}\,,\\
	\P_N(F_n(H))
	&\le 
	2^n \frac{\mathscr{C}
		(1/2^n)^2}
	{H^4/(2^n)^{4\gamma}}
	= \frac{\mathscr{C}
	(2^n)^{4\gamma}
	}{H^4 2^n}
	= \frac{\mathscr{C}
	}{H^4 2^{n/5}}\,,
	\end{align*}
where in each line the last step holds for $\gamma=1/5$. Consequently, if $A(H)$ denotes the union of $I(H)$ and the events $J_n(H) \cup F_n(H)$ over all $n\ge0$, for $N$ large enough we have
	\beq\label{e:A.H}
	\P_N (A(H))
	\le \sum_{n=0}^\infty
	\bigg\{
	\P_N(J_n(H))
	+\P_N(F_n(H))\bigg\}
	+\P_N(I(H))
	\le 
	\frac{8 \mathscr{C}}{H^2}
	+\frac{8 \mathscr{C}}{H^4}
	+\frac{\mathscr C}{H^2}
	\le \frac{ 17\mathscr{C}}{H^2}
	\,.
	\eeq
It follows by combining the above with Markov's inequality that, for $N$ large enough,
	\[
	\P_N\bigg(
	\bP_{\bG^N}(A(H))
	\ge 
	\frac{ 17 \mathscr{C}}{H}
	\bigg) 
	\le \frac{H \P_N(A(H))}{ 17\mathscr{C}}
	\le \frac{1}{H}\,.
	\]
Consequently, for any sequence $(H_k)_{k\ge0}$ that grows quickly enough, for all $N\ge N_0$ we have
	\beq\label{e:pm.pm.tightness}
	\P_N\bigg(
	\bP_{\bG^N}( A(H_k))
	>
		\frac{ 17 \mathscr{C}}{H_k}
	\textup{ for any }k\ge0
	\bigg)
	\le \sum_{k\ge0}\frac{1}{H_k}\,,
	\eeq
where the right-hand side can be made arbitrarily small by the choice of the sequence $(H_k)_{k\ge0}$. 

Analogously to \eqref{e:law.of.path.Gamma}, let $\bar{\mu}_{\bG^N}$ denote the law of the process $(\bar{\Gamma}^N,\bar{\Upsilon}^N)$ conditional on the disorder $\bG^N$.
Likewise, analogously to 
\eqref{e:pm.on.pm.over.paths.Gamma},
let $\bar{\Q}_N$ denote the law of the random measure $\bar{\mu}_{\bG^N}$.
Then the left-hand side of \eqref{e:pm.pm.tightness} above can be rewritten as $\bar{\Q}_N(\mathcal{K}^c)$, where
	\beq\label{e:compact.K}
	\mathcal{K}
	=\mathcal{K}((H_k)_{k\ge0})
	=\bigg\{\mu'\in
	\mathscr{P}(C([0,1],\R^8))
	: \mu'( A(H_k) ) \le 
	\frac{ 17 \mathscr{C}}{H_k}
	\textup{ for all }k\bigg\}\,.
	\eeq
For every finite $H$, the complement of $A(H)$ is compact (in the space 
$C([0,1],\R^8)$) by the Arzel\`a--Ascoli theorem. It follows by Prohorov's theorem (applied to probability measures over
$C([0,1],\R^8)$)
 that $\mathcal{K}$ is a compact set. Then
\eqref{e:pm.pm.tightness} can be rewritten as
	\[
	\bar{\Q}_N\bigg(\Big(\mathcal{K}(H_k)_{k\ge0}\Big)^c\bigg)
	\le \sum_{k\ge0}\frac{1}{H_k}\,,
	\]
which means that for every $\epsilon>0$ there is a compact set $\mathcal{K}$ (depending on $\epsilon$) that has measure at least $1-\epsilon$ under $\bar{\Q}_N$, for all $N\ge N_0$. It then follows by another application of Prohorov's theorem (now applied to probability measures over
$\mathscr{P}(C([0,1],\R^8))$) that the measures $\bar{\Q}_N$ are tight. 

Now suppose  $\bar{\Q}_N\Rightarrow \bar{\Q}$ along some subsequence of integers $N\to\infty$, where we use ``$\Rightarrow$'' to denote weak convergence. This means there exists a Skorohod coupling of the disorder matrices $\bG^N$ for which, with probability one over $(\bG^N)_{N\ge1}$, we have
$\bar{\mu}_{\bG^N}\Rightarrow\bar{\mu}$ with $\bar{\mu}\sim\bar{\Q}$. Therefore there exists a Skorohod coupling of random paths $(\bar{\Gamma}^N,\bar{\Upsilon}^N)\sim\bar{\mu}_{\bG^N}$ such that 
	\[\Big(\bar{\Gamma}^N,\bar{\Upsilon}^N\Big)
	\stackrel{N\to\infty}{\longrightarrow}
	\Big(\bar{\Gamma},\bar{\Upsilon}\Big)
	\sim\bar{\mu}
	\]
uniformly (almost surely). Under the same coupling, 
$\mathcal{E}^N$ converges uniformly to zero. 
This implies that $\Gamma^N-\bar{\Gamma}^N$ 
and $\Upsilon^N-\bar{\Upsilon}^N$ also converge uniformly to zero, which gives the non-``$\frozen$'' variant of the result.  The ``$\frozen$'' variant follows by the same argument, since all the preceding results of this section were also proved also for this setting. 
\end{proof}
\end{ppn}

\begin{cor}\label{c:trunc.small}
For $\sigma\in\{\RomI,\RomII,\RomIII\}$, let $p^\sigma\equiv p^\sigma(\bXi)$ denote the empirical averages
	\[
	\frac1M\sum_{a\le M}
	\max_{0\le d\le \dmax}
	\ind \Big\{
	X^{\trunc,\sigma}(q_d,a)
	\notin  [-\trK,\trK)\Big\}\,.
	\]
Similarly, let $p^\Ising\equiv p^\Ising(\bXi)$ denote the empirical average
	\[
	\frac1N\sum_{i\le N}
	\max_{0\le d\le \dmax}
	\ind \Big\{
	X^{\trunc,\Ising}(q_d,i)
	\notin  [-\trK,\trK)\Big\}\,.
	\]
Then $p^\sigma(\bXi) \le o_{\trK}(1)$
for all $\sigma\in\{\RomI,\RomII,\RomIII,\Ising\}$, with probability
at least $1- \exp(-N^{1/2})$.

\begin{proof}
Let $\bar{p}^\sigma$ be defined as $p^\sigma$, but with $[-\trK+1,\trK-1)$ in place of $[-\trK,\trK)$.
Let $f:\R\to\R$ be the Lipschitz function defined by
	\[f(x)=\begin{cases}
	0 & \textup{if $|x|\le \trK-1$,}\\
	1 & \textup{if $|x|\ge \trK$,}\\
		|x|-(\trK-1)
		&\textup{if $\trK-1
			\le |x| \le \trK$}\,.
	\end{cases}\]
Since $\ind\{x\notin[-\trK,\trK)\} \le f(x)
\le \ind\{x\notin[-\trK+1,\trK-1)\} 
$, for $\sigma\in\{\RomI,\RomII,\RomIII\}$ we have
	\[p^\sigma(\bXi)
	\le
	F^\sigma(\bXi)
	\equiv
	\frac1M
	\sum_{a\le M}
	\max_{0\le d\le \dmax}
	f(X^{\trunc,\sigma}(q_d,a))
	\le \bar{p}^\sigma(\bXi)
	\,.
	\]
We first argue that the \emph{annealed} expectation of $\bar{p}^\sigma$ is small for each $\sigma$. Note that since $\vX^{\trunc,\RomI:\RomIII}$ is frozen upon leaving $[-\trK,\trK)^3$ (Definition~\ref{d:trunc}), the probability that $\vX^{\trunc,\RomI:\RomIII}(q_d,a)$ lies outside $[-\trK,\trK)^3$ is maximized at time $d=\dmax$. From the spatial rerandomization (Definition~\ref{d:rerand}), this is the same as the probability that $\vY^{\RomI:\RomIII}(q_d)$ lies outside $[-\trK,\trK)^3$ at time $d=\dmax$. This (annealed) probability can be bounded using the proof of  Proposition~\ref{p:tightness}: recall from \eqref{e:A.H} that the event $A(H)$ has probability tending to zero as $H\to\infty$. On the complement of $A(H)$, the paths
$\bar{D}^\RomI$,
$Y^\RomII$, $Y^\RomIII$, $Y^\Ising$
are all $\gamma$-H\"older with H\"older norm $O(H)$, and in particular the maximum deviation of each path must be $O(H)$.
Moreover, since $\mathcal{E}^N$ converges uniformly to zero, the discrepancy between $\bar{D}^\RomI$ and $Y^\RomI$ converges uniformly to zero, so we can conclude that the maximum deviation of $Y^\RomI$ is also $O(H)$.
If we take $H$ to be a small constant multiple of $\trK$, then we conclude that $\E[\bar{p}^\sigma(\bXi)] \le o_{\trK}(1)$ for $\sigma\in\{\RomI,\RomII,\RomIII\}$. A similar argument gives 
$\E[\bar{p}^\Ising(\bXi)] \le o_{\trK}(1)$. 

The above implies that the annealed expectation of each $F^\sigma$ is also small. It remains to show that $F^\sigma$ is highly concentrated about this expectation. For $\sigma\in\{\RomI,\RomII,\RomIII\}$,
it follows by Corollary~\ref{c:conc.empir.avg} 
that $F^\sigma$ is $(L^{O(1)}/(\delta N^{1/2}))$-Lipschitz on the set $U_\MAX$. It follows by 
Assumption~\ref{a:params} that this Lipschitz constant is very small, for instance,
	\[
	\frac{L^{O(1)}}{\delta N^{1/2}}
	\ll \frac{1}{N^{1/3}}\,.
	\]
By the Kirszbraun extension theorem, there exists a function $\tilde{F}^\sigma(\bXi)$ which agrees with $F^\sigma$ for $\bXi\in U_\MAX$, but such that $\tilde{F}^\sigma$ is $(1/N^{1/3})$-Lipschitz on the entire space. 
It follows by standard concentration bounds (see e.g.\ Lemma~\ref{l:lip.subgaus}) that $\tilde{F}^\sigma$ is highly concentrated about its mean $\E[\tilde{F}^\sigma(\bXi)]$. We finally bound
	\[
	\E\Big[\tilde{F}^\sigma(\bXi)\Big]
	\le \E\Big[F^\sigma(\bXi)\Big]
	+ \E\Big[\tilde{F}^\sigma(\bXi);(U_\MAX)^c\Big]
	\le 
	\E\Big[F^\sigma(\bXi)\Big]
	+\bigg\{ \E\Big[\tilde{F}^\sigma(\bXi)^2\Big]
	\P\Big((U_\MAX)^c \Big)
	\bigg\}^{1/2}
	\le o_{\trK}(1)\,.
	\]
In the last step, we used that $\tilde{F}^\sigma$ is Lipschitz so its second moment must be bounded, while  \eqref{e.MAX.bound} implies that $\P((U_\MAX)^c)$ is very small. This proves the claim for $\sigma\in\{\RomI,\RomII,\RomIII\}$, and a similar argument gives the claim for $\sigma=\Ising$.
\end{proof}
\end{cor}

Note that for the proof of Proposition~\ref{p:tightness}, it was more convenient to work with discretized processes that were piecewise linear, since these can easily be shown to satisfy Kolmogorov estimates (Propositions \ref{p:Y.kolmogorov}, \ref{p:drift.kolmogorov}, and \ref{p:qv.kolmogorov}). However, it will often be more convenient to work with the piecewise \emph{constant} versions of the processes. Recall that we use $\vY^{\bullet,N}$ to denote the piecewise constant version of $\vY^N$: this means $\vY^{\bullet,N}(q_d)=\vY^N(q_d)$ for all $d$, and $\vY^{\bullet,N}$ is piecewise constant on each interval $[q_d,q_{d+1})$. We next argue that we still have weak convergence for the piecewise constant processes:

\begin{cor}[subsequential limits for piecewise constant processes]
\label{c:tightness}
Analogously to $\mu_{\bG^N}$ and $\Q_N$, let $\mu_{\bullet,\bG^N}$ denote the law of the piecewise constant processes
	\[\Big(
	\Gamma^{\bullet,N},\Upsilon^{\bullet,N},\mathcal{E}^{\bullet,N}
	\Big)\]
conditional on $\bG^N$, and let $\Q_{\bullet,N}$ denote the law of $\mu_{\bullet,\bG^N}$. Suppose $\Q_N\to\Q$ along an integer subsequence $N\to\infty$, as guaranteed by Proposition~\ref{p:tightness}. Then $\Q_{\bullet,N}\Rightarrow\Q$ along the same subsequence. The ``$\frozen$'' variant of the statement also holds.

\begin{proof}
By Proposition~\ref{p:tightness}, for any sequence of integers $N\to\infty$, there exists a further subsequence along which $\Q_N$ converges weakly to $\Q$, with $\mathcal{E}^N$ converging almost surely to the zero path. Take any such subsequence:
to simplify notation, we can assume (by re-indexing) that it is the entire sequence of integers $N\ge1$. 
This means there exists a Skorohod coupling of the disorder matrices $\bG^N$ for which, with probability one over $(\bG^N)_{N\ge0}$, we have $\mu_{\bG^N}\Rightarrow\mu$ along the same subsequence, where $\mu_{\bG^N}$ is the law of the process $(\Gamma^N,\Upsilon^N,\mathcal{E}^N)$ conditional on $\bG^N$. Moreover $\mathcal{E}^N$ converges almost surely to the zero path, which implies
that $\mathcal{E}^{\bullet,N}$ does also.

Now abbreviate $\Xi^N\equiv(\Gamma^N,\Upsilon^N)$
and $\Xi\equiv(\Gamma,\Upsilon)$.
 Let $L^\infty([0,1],\R^8)$ denote the space of (not necessarily continuous) functions $[0,1]\to\R^8$, equipped with the uniform metric. If $\Phi$ is a bounded Lipschitz functional on this space, then we have
	\beq\label{e:observable.a.s.conv}
	\lim_{N\to\infty} 
	\bE_{\bG^N} \Phi(\Xi^N)
	=\bE \Phi(\Xi)
	\eeq
almost surely (along the subsequence, under the Skorohod coupling).
Let $\Xi^{\bullet,N}$ be the piecewise constant variant of $\Xi^N$, and let $\mu_{\bG_\bullet,N}$ be the law of $\Xi^{\bullet,N}$. To conclude the result, it suffices to show that
	\[
	\lim_{N\to\infty} 
	\bE_{\bG^N} \Phi(\Xi^{\bullet,N})
	=\bE \Phi(\Xi)
	\]
almost surely for all bounded Lipschitz functionals $\Phi$, possibly after passing to a further subsequence of integers $N_k\to\infty$, as long as the subsequence does not depend on $\Phi$. In light of \eqref{e:observable.a.s.conv}, it suffices to show that
	\[\lim_{N\to\infty} 
	\bE_{\bG^N} \Big[
	\Phi(\Xi^N)
	-\Phi(\Xi^{\bullet,N})
	\Big]=0
	\]
almost surely along the subsequence $N_k$, for all bounded Lipschitz functionals $\Phi$.

Recall the proof of Proposition~\ref{p:tightness}. For any sequence $\epsilon_N\to0$, for each $N$ we have a compact set
 $\mathcal{K}_N=\mathcal{K}((H_{N,k})_{k\ge0})$ defined by \eqref{e:compact.K}, such that $\Q_N(\mathcal{K}_N^c) \le \epsilon_N$ for all $N \ge N_0$. We can take $\epsilon_N\to0$ slowly enough so that $H_{N,1} \delta^\gamma \to0$, where $\gamma=1/5$ is the H\"older exponent from the proof of Proposition~\ref{p:tightness}. We then have
 	\begin{align*}
	&\bigg|\bE_{\bG^N} \Big[
	\Phi(\Xi^N)
	-\Phi(\Xi^{\bullet,N})
	\Big] \bigg|
	\le
	2\|\Phi\|_\infty \ind\{ (\mathcal{K}_N)^c\}
	+\ind\{\mathcal{K}_N\}
	\bigg|
	\bE_{\bG^N} \Big[
	\Phi(\Xi^N)
	-\Phi(\Xi^{\bullet,N})
	\Big] 
	\bigg| \\
	&\qquad\le
	2\|\Phi\|_\infty \ind\{ (\mathcal{K}_N)^c\}
	+\ind\{\mathcal{K}_N\}\bigg\{
	2\|\Phi\|_\infty \mu( A(H_{N,1}))
	+ \|\Phi\|_{\textup{Lip}}
	\bE_{\bG^N}\Big[
	\|\Xi^N-\Xi^{\bullet,N}\|_\infty
	; A(H_{N,1})^c
	\Big]
	\bigg\} \\
	&\qquad\le
	2\|\Phi\|_\infty \ind\{ (\mathcal{K}_N)^c\}
	+\bigg\{
	2\|\Phi\|_\infty \frac{17\mathscr{C}}{H_{N,1}}
	+ \|\Phi\|_{\textup{Lip}} O(1) H_{N,1} \delta^\gamma
	\bigg\}\,.
	\end{align*}
If we now take a subsequence $N_k\to\infty$ such that $\epsilon_{N,k}$ is summable, then the above converges to zero almost surely along this subsequence. This concludes the proof of the non-``$\frozen$'' version of the statement. The ``$\frozen$'' variant follows by the same argument, since all the preceding results of this section were also proved for this setting. 
\end{proof}
\end{cor}

\begin{ppn}\label{p:limiting.drift.qv}
For any subsequential limit $\Q$ as obtained by Proposition~\ref{p:tightness}, let $\mu\sim\Q$ and $(\Gamma,\Upsilon,\mathcal{E})\sim\mu$ (so $\mathcal{E}\equiv0$). Denote the limiting processes as 
	\begin{align}\nonumber
	\Gamma
	&\equiv \Big(
	D^{\RomI}=D=Y^{\RomI},
	Q^{\RomII},Q^{\RomII,\RomIII},Q^{\RomIII},Q^{\Ising}
	\Big)\,,\\
	\Upsilon
	&\equiv\Big( 
	Z^{\RomII}=Y^{\RomII},
	Z^{\RomIII}=Y^{\RomIII},
	Z^{\Ising}=Y^{\Ising}
	\Big)\,.
	\label{e:limiting.process}
	\end{align}
Then $\Gamma$ is a finite-variation process, while $\Upsilon$ is a martingale. For each $\sigma,\tau\in\{\RomII,\RomIII,\Ising\}$, the covariation of $Z^\sigma$ with $Z^\tau$ is given by $Q^{\sigma,\tau}$, where
$Q^{\sigma,\Ising}\equiv0$ for $\sigma\in\{\RomII,\RomIII\}$. 
 The ``$\frozen$'' variant of the statement also holds.

\begin{proof}
Suppose $\Q_N\to\Q$ along some subsequence of integers $N\to\infty$.
By Corollary~\ref{c:tightness}, this means there exists a Skorohod coupling of the disorder matrices $\bG^N$ for which, with probability one over $(\bG^N)_{N\ge1}$, we have
$\mu_{\bullet,\bG^N}\Rightarrow \mu$ with $\mu\sim\bar{\Q}$. Let $\cG(s)$ be the filtration generated by the limiting process $(\Gamma,\Upsilon)$. It follows from \eqref{e:def.drift} that if 
$0\le s\le t\le 1$, then for all $\sigma\in\{\RomI,\RomII,\RomIII,\Ising\}$ and $A\in\cG(s)$, we have 
	\beq\label{e:discrete.mg.eq}
	\bE_{\bG^N}\bigg[ \Big( Y^{\bullet,\sigma}(t)
		-Y^{\bullet,\sigma}(s)\Big)
		\ind\{A\}\bigg]
	=\bE_{\bG^N}\bigg[ \Big( D^{\bullet,\sigma}(t)
		-D^{\bullet,\sigma}(s)\Big)
		\ind\{A\}\bigg]\,,
	\eeq
where $\bE_{\bG^N}$ denotes expectation with respect to $\mu_{\bullet,\bG^N}$.
Now fix an event $A\in\cG(s)$. Let $L_R,U_R$ be bounded continuous $\cG(s)$-measurable functions such that $L_R \le \ind\{A\} \le U_R$, and both $L_R,U_R$ converge to $\ind\{A\}$ in the limit $R\to\infty$. Moreover, we can arrange for the convergence to happen fast enough so that
	\beq\label{e:ind.A.approx}
	\E_N\Big(U_R-L_R\Big) \le \frac{o_R(1)}{R}\,,
	\eeq
where we recall that $\E_N$ denotes the overall expectation.
Define also the function
	\[
	I_R(x) = \max\Big\{\min\{x,R\},-R\Big\} \in[-R,R]\,.
	\]
Let $\bE$ denote expectation with respect to $\mu$.
Since $\mu_{\bullet,\bG^N}\Rightarrow\mu$, for fixed $R$ we have
	\begin{align*}
	\lim_{N\to\infty}
	\bE_{\bG^N}\bigg[ I_R\Big( Y^{\bullet,\sigma}(t)
		-Y^{\bullet,\sigma}(s)\Big)
	L_R \bigg]
	&=\bE\bigg[ I_R\Big(Y^\sigma(t)-Y^\sigma(s)\Big) L_R\bigg]\,,\\
	\lim_{N\to\infty}
	\bE_{\bG^N}\bigg[ I_R\Big( 
		D^{\bullet,\sigma}(t)-D^{\bullet,\sigma}(s)\Big)
	L_R\bigg]
	&=\bE\bigg[ I_R\Big(D^\sigma(t)-D^\sigma(s)\Big) L_R\bigg]
	\end{align*}
(and likewise with $U_R$ in place of $L_R$).
We then bound
	\begin{align*}
	&\bigg| \bE_{\bG^N}\Big[ I_R
		\Big( Y^{\bullet,\sigma}(t)-Y^{\bullet,\sigma}(s)\Big)
	L_R \Big] - \bE_{\bG^N}\Big[ 
		\Big( Y^{\bullet,\sigma}(t)-Y^{\bullet,\sigma}(s)\Big)
		;A\Big]\bigg| \\
	&\qquad\le R\bE_{\bG^N} \Big(U_R-L_R\Big)
	+ \bE_{\bG^N} \bigg[
	\Big| Y^{\bullet,\sigma}(t)-Y^{\bullet,\sigma}(s)\Big|
	; \Big| Y^{\bullet,\sigma}(t)-Y^{\bullet,\sigma}(s)\Big| \ge R\bigg]
	\end{align*}
Since $\mu_{\bG^N}\Rightarrow\mu$, for fixed $R$ we have
	\[
	\lim_{N\to\infty}
	R\bE_{\bG^N} \Big(U_R-L_R\Big) 
	= R\bE\Big(U_R-L_R\Big)\,,
	\]
and the overall expectation of this is small by
\eqref{e:ind.A.approx}.
It also follows using the moment bounds 
\eqref{e:Y.I.Kolmogorov}--\eqref{e:Y.II.III.Is.Kolmogorov} 
from Proposition~\ref{p:Y.kolmogorov} that
	\begin{align*}
	&\E_N \bigg[
	\Big| Y^{\bullet,\sigma}(t)-Y^{\bullet,\sigma}(s)\Big|
	; \Big| Y^{\bullet,\sigma}(t)-Y^{\bullet,\sigma}(s)\Big| \ge R\bigg] \\
	&\qquad\le \frac{\E_N[(Y^{\bullet,\sigma}(t)-Y^{\bullet,\sigma}(s))^2]}{R}
	\le \frac{L^{O(1)}}{R}\,,
	\end{align*}
where the bound holds uniformly over all $N$ large enough. Therefore, we have altogether
	\beq\label{e:Y.N.remove.truncation}
	\limsup_{N\to\infty} \E_N
	\bigg| \bE_{\bG^N}\Big[ I_R
		\Big( Y^{\bullet,\sigma}(t)-Y^{\bullet,\sigma}(s)\Big)
	L_R \Big] - \bE_{\bG^N}\Big[ 
		\Big( Y^{\bullet,\sigma}(t)-Y^{\bullet,\sigma}(s)\Big)
		;A\Big]\bigg|
	\le o_R(1)\,.
	\eeq
A similar but simpler argument gives, in the limit,
	\beq
	\label{e:Y.limit.remove.truncation}
	\E
	\bigg|
	\bE\Big[ \Big( Y^\sigma(t)-Y^\sigma(s)\Big)
		\ind\{A\}\Big]
	-
	\bE\Big[ I_R\Big(Y^\sigma(t)-Y^\sigma(s)\Big) L_R\Big]
	\bigg| \le o_R(1)\,.
	\eeq
Indeed, similarly as above, the left-hand side of \eqref{e:Y.limit.remove.truncation} has expectation upper bounded by
	\begin{align*}
	&\le R\E \Big(U_R-L_R\Big)
	+ \E\bigg[ \Big|Y^\sigma(t)-Y^\sigma(s)\Big|; 
	\Big|Y^\sigma(t)-Y^\sigma(s)\Big|
	\ge R\bigg] \\
	&\qquad\le o_R(1)
	+ \frac{\E[(Y^\sigma(t)-Y^\sigma(s))^2]}{R}
	\le o_R(1)
	+\frac{L^{O(1)}}{R}\,,
	\end{align*}
having again used \eqref{e:ind.A.approx}, as well as the moment bounds \eqref{e:Y.I.Kolmogorov}--\eqref{e:Y.II.III.Is.Kolmogorov} from Proposition~\ref{p:Y.kolmogorov} together with Fatou's lemma. This justifies \eqref{e:Y.limit.remove.truncation}. 
It follows by combining 
\eqref{e:Y.N.remove.truncation} and \eqref{e:Y.limit.remove.truncation} that
	\[\lim_{N\to\infty}
	\bE_{\bG^N}\bigg[ \Big( Y^{\bullet,\sigma}(t)-Y^{\bullet,\sigma}(s)\Big)
		;A\bigg]
	= \bE\bigg[ \Big( Y^\sigma(t)-Y^\sigma(s)\Big)
		;A\bigg]
	\]
for any $A\in\cG(s)$, where the limit holds in probability. A similar argument (using Proposition~\ref{p:drift.kolmogorov} and Lemma~\ref{l:D.barD.discrep} in place of Proposition~\ref{p:Y.kolmogorov}) gives
	\[
	\lim_{N\to\infty}
	\bE_{\bG^N}\bigg[ \Big( D^{\bullet,\sigma}(t)-D^{\bullet,\sigma}(s)\Big)
		;A\bigg]
	= \bE\bigg[ \Big( D^\sigma(t)-D^\sigma(s)\Big)
		;A\bigg]
	\]
in probability. It follows that the martingale equation \eqref{e:discrete.mg.eq} passes to the limit. Thus we conclude that $D^{\RomI}$ is a finite-variation process. For $\sigma\in\{\RomII,\RomIII,\Ising\}$, we already saw in Proposition~\ref{p:tightness} that $D^\sigma\equiv0$ in the limit, so we can now conclude that $Z^\sigma$ is a martingale.

It remains to prove that for each $\sigma,\tau\in\{\RomII,\RomIII,\Ising\}$, the covariation of $Z^\sigma$ with $Z^\tau$ is given by $Q^{\sigma,\tau}$. Analogously to \eqref{e:discrete.mg.eq}, we have the discrete martingale equations 
	\[
	\bE_{\bG^N}\bigg[ \Big( Z^{\bullet,\sigma}(t)Z^{\bullet,\tau}(t)
		-Z^{\bullet,\sigma}(s)Z^{\bullet,\tau}(s)
		\Big)
		;A\bigg]
	=\bE_{\bG^N}\bigg[ \Big( Q^{\bullet,\sigma,\tau}(t)
	-Q^{\bullet,\sigma,\tau}(s)\Big)
		;A\bigg]
	\]
for $0\le s\le t\le 1$ and $A\in\cG(s)$. These equations can be passed to the limit using the bounds from Propositions \ref{p:Y.kolmogorov} and \ref{p:qv.kolmogorov}, and this implies the non-``$\frozen$'' version of the statement. The ``$\frozen$'' variant follows by the same argument, since all the preceding results of this section were also proved for this setting.
\end{proof}
\end{ppn}

\begin{proof}[\hypertarget{proof:t.tightness}{Proof of Theorem~\ref{t:tightness}}]
This follows by combining the statements of Proposition~\ref{p:tightness},
Corollary~\ref{c:tightness}, and Proposition~\ref{p:limiting.drift.qv}.
\end{proof}

\fi

\pagebreak\section{Derivation of budget constraint}\label{a:freeprob}

\iffull
% !TEX root = main.tex

\newcommand{\BIsing}{K}

\newcommand{\cD}{\mathcal{D}}
\newcommand{\cQ}{\mathcal{Q}}
\newcommand{\hcD}{\hat\cD}

\newcommand{\Budget}{\textup{\textsf{budget}}}
\newcommand{\Surplus}{\textup{\textsf{surplus}}}
\newcommand{\Extr}{\textup{\textsf{Extr}}}
\newcommand{\Viol}{\textup{\textsf{viol}}}
\newcommand{\Surf}{\textup{\textsf{surf}}}
\newcommand{\Proj}{\textup{\textsf{P}}}

\newcommand{\ytrunc}{\hat{\by}_{\textup{tr}}}

\newcommand{\pol}{\textup{\textup{pol}}}

\newcommand{\tbmeta}{\tilde{\bmeta}}

\def\baln#1\ealn{%
    \begin{align*}%
    #1%
    \end{align*}%
}

\newcommand{\Lagr}{\mathscr{L}}

\newcommand{\YY}{\mathfrak{Y}}
\newcommand{\hY}{\hat{Y}}

\newcommand{\BBb}{\mathfrak{B}}
\newcommand{\UUu}{\mathfrak{U}}
\newcommand{\xxxx}{\mathfrak{x}}
\newcommand{\yyy}{\mathfrak{y}}

% removed underlines
\newcommand{\vBBb}{{\BBb}}
\newcommand{\vUUu}{{\UUu}}
\newcommand{\vbbb}{{\bbb}}
\newcommand{\vvvv}{{\vvv}}
\newcommand{\vuuu}{{\uuu}}
\newcommand{\vwww}{{\www}}
\newcommand{\vrrr}{{\rrr}}

\newcommand{\hBBb}{\hat \BBb}
\newcommand{\hUUu}{\hat \UUu}
\newcommand{\hxxx}{\hat \xxxx}
\newcommand{\hvvv}{\hat \vvv}
\newcommand{\hwww}{\hat \www}
\newcommand{\hvBBb}{\hat \vBBb}
\newcommand{\hvUUu}{\hat \vUUu}
\newcommand{\hvxxx}{\hat \vec \xxxx}
\newcommand{\hvvvv}{\hat \vvvv}
\newcommand{\hvwww}{\hat \vwww}
\newcommand{\budg}{\textup{bd}}
\newcommand{\cW}{{\mathcal{W}}}
\newcommand{\bM}{{\boldsymbol{M}}}
\newcommand{\eig}{{\textup{eig}}}
\newcommand{\bone}{{\boldsymbol 1}}
\newcommand{\vecB}{{\vec B}}
\newcommand{\vecv}{{\vec v}}
\newcommand{\vecU}{{\vec U}}
% \newcommand{\vecw}{{w}}

% removed underlines
\newcommand{\hvB}{\hat B}
\newcommand{\hvv}{\hat v}
\newcommand{\hvU}{\hat U}
\newcommand{\hvw}{\hat w}

\newcommand{\tv}{\tilde v}
\newcommand{\tvv}{\tilde v}
\newcommand{\tx}{\tilde x}
\newcommand{\hB}{\hat B}
\newcommand{\hv}{\hat v}
\newcommand{\hU}{\hat U}
\newcommand{\hw}{\hat w}
\newcommand{\hx}{\hat x}
\newcommand{\N}{\mathbb{N}}

% versions of F
\newcommand{\ri}{\textup{i}}
\newcommand{\rii}{\textup{ii}}
\newcommand{\riii}{\textup{iii}}
\newcommand{\riv}{\textup{iv}}
\newcommand{\rv}{\textup{v}}
\newcommand{\rvi}{\textup{vi}}
\newcommand{\rvii}{\textup{vii}}
\newcommand{\rviii}{\textup{viii}}

\newcommand{\scB}{\mathscr{B}}
\newcommand{\gd}{\textup{good}}
\newcommand{\btheta}{\boldsymbol \theta}
\newcommand{\bTheta}{\boldsymbol \Theta}
\newcommand{\dbg}{\dot \bg}
\newcommand{\hbg}{\hat \bg}
\newcommand{\bD}{\boldsymbol D}
\newcommand{\dbD}{\dot \bD}
\newcommand{\dbX}{\dot{\boldsymbol{X}}}
\newcommand{\hbD}{\hat \bD}
\newcommand{\dbA}{\dot \bA}
\newcommand{\hbA}{\hat \bA}
\newcommand{\dA}{\dot A}
\newcommand{\dg}{\dot g}
\newcommand{\bI}{\boldsymbol I}
\newcommand{\vc}{\vec c}

\newcommand{\barbg}{\bg}

In this section we prove the budget constraints as stated in Theorem~\ref{t:free.prob}. We begin by recalling the setting of Theorem~\ref{t:free.prob}, which we assume throughout this section, and introducing some additional notations.
Let $\omega$ denote a random variable.
Let $\bG, \bXi$ be $M\times N$ matrices with i.i.d.\ standard gaussian entries, with rows $\barbg^1,\ldots,\barbg^M$ and $\bmeta^1,\ldots,\bmeta^M$ respectively, such that $\bXi$ is independent of $(\bG,\omega)$. Let $\bar{\bx} = \bar{\bx}(\bG)$ and $\hat{\by} = \hat{\by}(\bG, \bXi)$ be $L$-Lipschitz, $\R^N$-valued functions such that 
\begin{align*}
    &\E \|\bar{\bx}\|^2 = \E \|\hat{\by}\|^2 = N\,, \\ 
    &\E(\hat{\by}\,|\,\bG) = \bzero\,\,\text{almost surely}\,.
\end{align*}
Let $r,s \ge 1$ be fixed. 
Suppose we have partitions, which are allowed to depend on $(\bG,\omega)$ but not $\bXi$:
\begin{align*}
    [M] = \bigsqcup_{j=1}^r B_j \textup{ with } 
    \frac{|B_j|}{M}
    = \lambda_j(\bG,\omega)\,,\\
    [N] = \bigsqcup_{j=1}^s 
    K_j \textup{ with } 
    \frac{|K_j|}{N}
    = \mu_j(\bG,\omega)\,.
\end{align*}
We write 
\beq
    \label{e:fp-partition}
    \hat{\cB} = \hat{\cB}(\bG,\omega)
    = \Big(B_1,\ldots,B_r,K_1,\ldots,K_s\Big)\,.
\eeq
Further define
\begin{align}
    \label{e:lambda-mu-random}
    \lambda(\bG,\omega) &= (\lambda_1(\bG,\omega),\ldots,\lambda_r(\bG,\omega))\,, &
    \mu(\bG,\omega) &= (\mu_1(\bG,\omega),\ldots,\mu_s(\bG,\omega))\,.
\end{align}
We assume it holds $(\bG,\omega)$-almost surely that
\beq
    \label{e:def-Lambda}
    (\lambda(\bG,\omega),\mu(\bG,\omega)) \in
    \Lambda 
    \equiv \bigg\{
        \begin{array}{l}
        (\lambda,\mu) \in [0,1]^r \times [0,1]^s: 
        \|\lambda\|_1 = \|\mu\|_1 = 1, \textup{ and} \\
       \qquad\qquad\qquad\qquad\qquad
       \min \{
            \min_{j\in [r]} \lambda_j, 
            \min_{j\in [s]} \mu_j \} \ge \iota
        \end{array} \bigg\}\,.\eeq
Note that this assumption implies
	\beq\label{e:r-s-bd}\max\{r,s\} \le \frac1\iota\,,\eeq
which we also assume throughout this section. Let $\be_i$ denote the $i$-th standard basis vector in $\R^N$. Then, recalling the notation \eqref{e:freeprob.b}--\eqref{e:freeprob.w}, we define the random variables
\begin{align}
    \label{e:fp.b}
    \bbb_j
    &\equiv
    \frac{1}{|B_j|} \sum_{a\in B_j}
    \frac{(\bmeta^a,\hat{\by})}{N^{1/2}}
    \,, \\
    \label{e:fp.v}
    \vvv_j
    &\equiv
    \frac{1}{|B_j|} \sum_{a\in B_j}
    \frac{(\barbg^a,\hat{\by})^2}{N}
    \,, \\
    \label{e:fp.u}
    \uuu_j
    &\equiv
    \frac{1}{|B_j|} \sum_{a\in B_j}
    \frac{(\bmeta^a,\bar{\bx})(\barbg^a,\hat{\by})}{N}
    \,, \\
    \label{e:fp.w}
    \www_j
    &\equiv
    \frac{1}{|\BIsing_j|} \sum_{i\in \BIsing_j}
    (\be_i,\hat{\by})^2 \,.
\end{align}
Recall from \eqref{e:cost} the function $\Cost\equiv \Cost_0 : \R \times \R \times [0,+\infty) \to \R$ defined by
\[
    \Cost(b,v,u)
    = b^2 + u^2 + \bigg(
        \Big[(v - u^2)_+\Big]^{1/2} - 1\bigg)^2\,.
\]
The main result of this appendix is the following, which restates Theorem~\ref{t:free.prob}:

\begin{ppn}\label{p:fp}
Fix parameters $\alpha, L, r,s, \iota, \epsilon$. As above, let $\omega$ denote a random variable, and let $\bG, \bXi$ be $M\times N$ matrices with i.i.d.\ standard gaussian entries, with rows $\barbg^1,\ldots,\barbg^M$ and $\bmeta^1,\ldots,\bmeta^M$ respectively, such that $\bXi$ is independent of $(\bG,\omega)$. Let $\bar{\bx} = \bar{\bx}(\bG)$ and $\hat{\by} = \hat{\by}(\bG, \bXi)$ be $L$-Lipschitz, $\R^N$-valued functions such that 
\begin{align*}
    &\E \|\bar{\bx}\|^2 = \E \|\hat{\by}\|^2 = N\,, \\ 
    &\E(\hat{\by}\,|\,\bG) = \bzero\,\,\text{almost surely}\,.
\end{align*} As in \eqref{e:fp-partition}, let $ \hat{\cB} = \hat{\cB}(\bG,\omega)$ be a random partition of $[M]$ and $[N]$ into $r$ and $s$ parts respectively, with proportions $\lambda(\bG,\omega),\mu(\bG,\omega)$ as in \eqref{e:lambda-mu-random}. Assume these proportions satisfy condition~\eqref{e:def-Lambda}, so that $r,s$ must satisfy condition~\eqref{e:r-s-bd}. Then, for all sufficiently large $N$,  we have
    \beq
        \label{e:fp}
        \sum_{j=1}^r
        \lambda_j(\bG,\omega) \,
        \Cost(\bbb_j,\uuu_j,\vvv_j)
        \le \frac1\alpha
         \bigg(
            \sum_{j=1}^s
            \mu_j(\bG,\omega)
            (\www_j)^{1/2}
        \bigg)^2
        + \epsilon
    \eeq
with probability at least $1-e^{-cN}$, where $c$ is a positive constant depending on the parameters $\alpha, L, \iota, \epsilon$.
\end{ppn}

Throughout this section we assume without loss of generality that $\epsilon \le 1$. We also assume $\iota \le 1$, as otherwise $\Lambda$ defined in \eqref{e:def-Lambda} is empty and the result is vacuous. Recalling the notation \eqref{e:fp.b}--\eqref{e:fp.w}, we let 
    \beq\label{e:freeprob.random.YY}
    \YY \equiv 
    (\vBBb\equiv\bbb^2,\vvvv,
    \vUUu\equiv\uuu^2,\vwww)\,,
    \eeq
where we square $\bbb$ coordinatewise to obtain the vector $\vBBb\equiv\bbb^2$ of the same dimension, and similarly for $\vUUu$. We must show that, with high probability, the random vector $\YY$ satisfies the constraints \eqref{e:fp}. This section is organized as follows:
\begin{itemize}
    \item \S\ref{ss:fp-linear-proj} reduces Proposition~\ref{p:fp} to Propositions~\ref{p:fp-projection} and \ref{p:convex.body.polytope.bound.new}:
\begin{itemize}
\item
The proof of Proposition~\ref{p:fp-projection} occupies most of this section, and we summarize its result as follows. For the random vector $\YY$ defined by
 \eqref{e:fp.b}--\eqref{e:fp.w} and \eqref{e:freeprob.random.YY}, we define $F\equiv F_{\lambda,\mu}(\hY;\YY)$ to be the (scalar-valued) projection of $\YY$ in the direction $\hY$; the precise definition is given in \eqref{e:1dproj} below. For comparison, we define $F_\star \equiv F_{\star,\lambda,\mu}(\hY)$ to be the maximum of this linear projection over \emph{all} vectors $\YY$ satisfying the constraint \eqref{e:fp}; see \eqref{e:Fstar} below. The result of Proposition~\ref{p:fp-projection} is that for any fixed $\lambda,\mu\in\Lambda$ and fixed $\hY$,  the random quantity 
$F$
is upper bounded by the deterministic quantity $F_\star$, except with probability $\exp(-cN)$.
\item The basic idea of the reduction is that intersecting the bound of Proposition~\ref{p:fp-projection} in $O(1)$ directions should imply that the random vector $\YY$ lies (with high probability) in or near the convex body of all vectors $\YY$ satisfying \eqref{e:fp}; and this should yield the desired result Proposition~\ref{p:fp}. To make this intuition rigorous,   Proposition~\ref{p:convex.body.polytope.bound.new} provides,
for the convex body defined by \eqref{e:fp}, an approximation via a polytope whose number of facets is appropriately uniform in the key parameters of the problem. This allows us to take a union bound with the conclusion of Proposition~\ref{p:fp-projection} summed over an appropriately bounded number of terms, thereby completing the reduction. 
\end{itemize}
Separately, we also show  that the optimization problem $F_\star$ is equivalent to a simplified optimization problem $F_\diamond \equiv F_{\diamond,\lambda,\mu}$; see Lemma~\ref{l:fp-Fst-to-Fdiam}.

\item \S\ref{ss:polytope.appx.new} is devoted to the 
\hyperlink{proof:p.convex.body.polytope.bound.new}{proof of Proposition~\ref{p:convex.body.polytope.bound.new}}. The rest of this appendix (\S\ref{ss:k-orth-replica}--\ref{ss:fp-eval}) is then devoted to Proposition~\ref{p:fp-projection}. This will be done by a succession $F^{\rii},F^{\riii},\ldots$ of approximate stochastic upper bounds for the original quantity $F^{\ri} \equiv F = F_{\lambda,\mu}(\hY;\YY)$, which gradually reduce the desired estimate to more tractable problems.
\item The purpose of \S\ref{ss:k-orth-replica}--\ref{ss:fp-uc} is to reduce from the setting of Proposition~\ref{p:fp-projection}, where we have a random partition \eqref{e:fp-partition} depending on $(\bG,\omega)$, to the setting of a deterministic partition (independent of $(\bG,\omega)$). This proceeds by a \textbf{uniform concentration} idea due to \cite{subag2018free}.
    In \S\ref{ss:k-orth-replica}, we bound $F$ by an optimization problem $F^{\riii}(k)$ over $k$ orthogonal points in $S_N$; see Propositions~\ref{p:fp-f-to-i} and \ref{p:fp-i-to-ii}. In \S\ref{ss:fp-uc} we show that for each fixed partition $\cB$, this quantity concentrates with probability $1-\exp(-ckN)$. For $k$ sufficiently large, this is enough to overcome a union bound over all possible $\cB$ (Corollary~\ref{c:fp-uc}). The main result of \S\ref{ss:fp-uc}, given by Propositions \ref{p:fp-ii-to-iii} and \ref{p:fp-iii-to-iv}, is to reduce $F^{\riii}(k)$ to an optimization problem $F^{\rv}$, which can be viewed as the top eigenvalue of a certain random matrix: a sum of a deterministic diagonal matrix, a deformed Wishart matrix, and a random $O(1)$-rank spike.
    
\item
    The remaining goal is to show $F^{\rv}$ is bounded by $F_\diamond$, which amounts to controlling the top eigenvalue of the random matrix we just described. As outlined informally in \S\ref{ss:heuristic.budget}, one way to do this is with free probability theory. In our formal proof, we take a different approach: in \S\ref{ss:fp-gaussian-comparison}, we bound $F^{\rv}$ using Slepian and Gordon's gaussian comparison inequalities. This leads to a minimax problem $F^{\rvi}$ on $\R^N \times \R^M$ (Proposition~\ref{p:fp-vi-to-v}) in which the deformed Wishart matrix is replaced by several \emph{linear} terms, which are each the inner product of the optimization variable with a gaussian vector. This allows us to reduce the main result Proposition~\ref{p:fp-projection} to a simpler statement, Proposition~\ref{p:fp-v-estimate}, which bounds $F^{\rvi}$
 in terms of $F_\diamond$. 
 
\item In \S\ref{ss:fp-F-diamond} 
we evaluate $F_\diamond$ using the method of Lagrange multipliers, resulting in a non-variational formula $V(z_\star)$ (Proposition~\ref{p:fp-Fdiam-to-V}). Finally, in \S\ref{ss:fp-eval}
we bound $F^{\rvi}$ by $V(z_\star)$, leading to the proof of Proposition~\ref{p:fp-v-estimate} (and hence Proposition~\ref{p:fp-projection}). 
\end{itemize}

\subsection{Linear projections of statistics}
\label{ss:fp-linear-proj}

In what follows, we will uniformly upper bound the probability that \eqref{e:fp} fails to hold on any permissible realization of $\lambda(\bG,\omega)$, $\mu(\bG,\omega)$. For any fixed $(\lambda,\mu) \in \Lambda$, with $\Lambda$ as defined in \eqref{e:def-Lambda}, we define the event
\beq
    \label{e:E-lambda-mu}
    \cE(\lambda,\mu) = \{\lambda(\bG,\omega) = \lambda, \mu(\bG,\omega) = \mu\}\,.
\eeq
Define the domain
\beq\label{e:cW}
    \cW_\mu \equiv \Big\{
        \vwww \in [0,+\infty)^s :
        (\mu, \vwww) = 1
    \Big\}\,,
\eeq
as well as
\begin{align}
    \label{e:cD-plus}
    \cD_+ &\equiv [0,+\infty)^r \times [0,+\infty)^r \times [0,+\infty)^r \times [0,+\infty)^s\,, \\
    \nonumber
    \cD_\mu &\equiv \{
        \YY\equiv(\vBBb,\vvvv,\vUUu,\vwww) \in \cD_+ : 
        \vwww \in \cW_\mu, \vvvv \ge \vUUu
    \}\,,
\end{align}
where the notation means that $\vvvv \ge \vUUu$ coordinatewise. We also let
	\beq\label{e:D.mu.affine.closure}
	\mathcal{F}_\mu
	\equiv\Big\{
	\YY\equiv(\vBBb,\vvvv,\vUUu,\vwww) \in \R^{3r+s}
	: (\mu,\vwww)=1\Big\}
	\,,\eeq
and note that $\mathcal{F}_\mu$ is the affine closure of $\cD_\mu$. 
Recall the above notations \eqref{e:fp.b}--\eqref{e:fp.w}, and also recall that we denote
$\BBb_j = (\bbb_j)^2$ and $\UUu_j = (\uuu_j)^2$. 
The random vector $\YY \equiv (\vBBb,\vvvv,\vUUu,\vwww)$ of 
\eqref{e:freeprob.random.YY} clearly takes values in $\cD_+$, and we will see in Lemma~\ref{l:random-YY-in-cD} below that it takes values approximately in $\cD_{\mu(\bG,\omega)}$.
Define the weighted inner products
\begin{align}
    \label{e:weighted.inner.prod}
    (\hvv,\vvvv)_{\lambda} &= \sum_{j=1}^r \lambda_j \hv_j \vvv_j \,, &
    (\hvw,\vwww)_{\mu}
    &= \sum_{j=1}^s \mu_j \hw_j \www_j \,.
\end{align}
We also abbreviate $\hY \equiv (\hvB,\hvv,\hvU,\hvw)$, which takes values in $\hcD \equiv
    \bbR^{3r+s}$.
We further assume that $\hY$ satisfies the normalization condition 
\beq
    \label{e:hY-normalization}
    \Big\|(\lambda\hvB, \lambda\hvv, \lambda\hvU, \mu\hvw)
    	\Big\|_2 = 1\,.
\eeq
In the above, $\lambda\hvB$ denotes the entrywise product of the $r$-dimensional vectors $\lambda$ and $\hvB$, and similarly for $\lambda\hvv, \lambda\hvU, \mu\hvw$. We note for use later that for any $(\lambda,\mu) \in \Lambda$, the normalization condition~\eqref{e:hY-normalization} implies 
\beq \label{e:fp-normalization-bd}
    \|\hB,\hv,\hU,\hw\|_\infty
     \le \frac1\iota\,.
\eeq
The vector $\hY$ specifies the direction of a one-dimensional projection of $\YY$ (see \eqref{e:1dproj} just below), so the normalization \eqref{e:hY-normalization} is without loss of generality. We then write
	\beq\label{e:1dproj}
    F_{\lambda,\mu}(\hY;\YY)
    \equiv (\YY,\hY)_{\lambda,\mu}
    \equiv (\hvB,\vBBb)_{\lambda}
    + (\hvv,\vvvv)_{\lambda}
    + (\hvU,\vUUu)_{\lambda}
    + (\hvw,\vwww)_{\mu}\,,
    \eeq
using the weighted inner product notation of \eqref{e:weighted.inner.prod}. Our proof will be based on bounding one-dimensional projections of the form \eqref{e:1dproj}. 
Define the cost and budget quantities 
    \begin{align}
    \label{e:appx.freeprob.cost}
    \Cost(\YY;\lambda)
    &\equiv
    \sum_{j=1}^r
        \lambda_j \,
        \bigg(
        \BBb_j + \UUu_j
        + \Big([(\vvv_j-\UUu_j)_+]^{1/2}
        -1 \Big)^2\bigg)\,,\\
    \Budget(\YY;\mu)
    &\equiv \frac{1}{\alpha}
        \bigg(
            \sum_{j=1}^s
            \mu_j (\www_j)^{1/2}
        \bigg)^2\,.
    \label{e:appx.freeprob.budget}
    \end{align}
This is defined for all $\YY \in \cD_+$; note that for $\YY \in \cD_\mu$ the $(\cdot)_+$ in \eqref{e:appx.freeprob.cost} can be omitted.
We then define the \textbf{ideal feasible set} as 
\beq\label{e:fp.ideal.feasible.set}
    S_\star(\lambda,\mu)
    = \Big\{ 
        \YY\in \cD_\mu : 
        \Cost(\YY;\lambda)
        \le
        \Budget(\YY;\mu)
    \Big\}\,.
\eeq 
We let $F_{\star,\lambda,\mu}$ denote the maximum value of the linear projection \eqref{e:1dproj} over the ideal feasible set, that is,
\beq
    \label{e:Fstar}
    F_{\star,\lambda,\mu}(\hY)
    = \sup\Big\{
    F_{\lambda,\mu}(\hY, \YY)
        :
        \YY \in S_\star(\lambda,\mu)
    \Big\}.
\eeq
Let us clarify that in the definition \eqref{e:1dproj} of $F_{\lambda,\mu}$, the notation $\YY\equiv (\vBBb,\vvvv,\vUUu,\vwww)$ indicates the random variable defined by \eqref{e:fp.b}--\eqref{e:fp.w} and \eqref{e:freeprob.random.YY}. By contrast, in the definition \eqref{e:Fstar} of $F_{\star,\lambda,\mu}$, the notation $\YY$ indicates any element of $S_\star(\lambda,\mu)$, and we optimize over all such $\YY$. As an example of the meaning of $F_{\star,\lambda,\mu}$, we remark on a simple special case of $\hY=(\hvB,\hvv,\hvU,\hvw)$.
Let $t\in \bbR$. If we take $\hvB,\hvv,\hvU\equiv \vec 0 \in \bbR^r$ and $\hvw\equiv t\cdot \vec 1 \in \bbR^s$, then $F_{\star,\lambda,\mu}(\hY)$ simplifies to
\[
    F_{\star,\lambda,\mu}(\hY)
    = \sup\Big\{
      t(\mu,\www)
      : \YY \in S_\star(\lambda,\mu)\Big\}
    =t\,,
  \]
where the last equality holds because $\YY \in \cD_\mu$ implies $\www\in\cW_{\mu}$.
\textbf{The following estimate is the main input to the \hyperlink{proof:p.fp.given.projection}{proof of Proposition~\ref{p:fp}:}}

\begin{ppn}[proved in \S\ref{ss:k-orth-replica}--\ref{ss:fp-eval}]
\label{p:fp-projection}
Fix parameters $\alpha, L, \iota, \epsilon'$. 
Recall the definition of $\Lambda$ from \eqref{e:def-Lambda}; and recall the definitions of $F_{\lambda,\mu},F_{\star,\lambda,\mu}$ from \eqref{e:1dproj}, \eqref{e:Fstar}.
There exists a positive constant $c=c(\alpha,L,\iota,\epsilon')$ such that the following holds for any fixed 
    $(\lambda,\mu) \in \Lambda$, and any fixed $\hY \in \hcD$ satisfying the normalization condition~\eqref{e:hY-normalization}:
    \beq
        \label{e:fp-projection}
        \P\bigg( 
            F_{\lambda,\mu}(\hY;\YY) 
            > F_{\star,\lambda,\mu}(\hY) + \epsilon'
            ; \cE(\lambda,\mu)
        \bigg)
        \le e^{-cN}\,,
    \eeq
    where $\P(A;E)\equiv \P(A\cap E)$ and the event $\cE(\lambda,\mu)$ is defined in \eqref{e:E-lambda-mu}.
\end{ppn}

Our proof of Proposition~\ref{p:fp} relies in addition on the following result.
For any subset $S \subseteq \cD_\mu$ and any $\epsilon' > 0$, let $S[\epsilon']$ denote the closed $\epsilon'$-neighborhood of $S$ within $\cD_\mu$.
That is,
\beq\label{e:closed.nbd.inside.D.mu}
    S[\epsilon'] = \bigg\{
        \YY \in \cD_\mu : \inf\Big\{
            \|\tilde\YY - \YY\|_2 : \tilde\YY \in S
        \Big\} \le \epsilon'
    \bigg\}\,.\eeq
We say $P\subseteq \cD_\mu$ is a \textbf{closed $\cD_\mu$-polytope of $n$ facets} if there exist closed half-spaces $H_j \subseteq \bbR^{3r+s}$, indexed by $j\in[n]$, such that $P$ can be expressed as
\beq
    \label{e:def.polytope}
    P = \cD_\mu \cap \bigcap_{j=1}^n H_j\,.
\eeq
The following proposition shows that for any $(\lambda,\mu) \in \Lambda$, the convex body $S_\star(\lambda,\mu)$ of \eqref{e:fp.ideal.feasible.set} can be approximated by a polytope with an appropriately bounded number of facets: 

\begin{ppn}[proved in \S\ref{ss:polytope.appx.new}]\label{p:convex.body.polytope.bound.new}
Fix parameters $\alpha,\iota,\epsilon$. There exist $\epsilon',C_\pol > 0$, depending only on $\alpha,\iota,\epsilon$, such that the following holds for any $(\lambda,\mu) \in \Lambda$: for $S_\star(\lambda,\mu)$ defined by \eqref{e:fp.ideal.feasible.set}, the closed $\epsilon'$-neighborhood $S_\star(\lambda,\mu)[\epsilon']$
(as defined by \eqref{e:closed.nbd.inside.D.mu}) is contained in a closed $\cD_\mu$-polytope $P$ of $n$ facets for some $n\le C_{\pol}$, such that every $\YY\in P$ satisfies the bound
\[
    \Surplus(\YY;\lambda,\mu)
    \equiv\Cost(\YY;\lambda)
    - \Budget(\YY;\mu) 
    \le \frac{\epsilon}{2}\,,
\]
for $\Cost$ and $\Budget$ as defined by \eqref{e:appx.freeprob.cost} and \eqref{e:appx.freeprob.budget}.
\end{ppn}

We now work towards the \hyperlink{proof:p.fp.given.projection}{proof of Proposition~\ref{p:fp}, assuming Propositions~\ref{p:fp-projection} and \ref{p:convex.body.polytope.bound.new}}.
We begin with the following simple estimates, which will be used throughout this appendix.

\begin{lem}\label{l:fp-gc}
For any constant $\delta > 0$, there exists $c = c(\alpha, L, \delta) > 0$ such that with probability $1-e^{-cN}$,
    \begin{align}
        \label{e:fp-x-conc}
        \bigg|\frac{\|\bar{\bx}\|^2}{N} - 1\bigg| &\le \delta\,,\\
    \label{e:fp-y-cond-conc}
        \bigg|\frac{\E[ \|\hat{\by}\|^2 | \bG]}{N} - 1\bigg| &\le
        \frac{\delta}{2}\,.
    \end{align} 
Moreover, conditional on any realization of $\bG$ satisfying \eqref{e:fp-y-cond-conc} the following holds. Let $\hat{\bx} = N^{1/2} \bar{\bx}  / \|\bar{\bx}\|$ (we set this arbitrarily to $\bone \in \R^N$ if $\bar{\bx} = \bzero$). Let $\bXi'$ be an independent copy of $\bXi$ and $\hat{\by}' = \hat{\by}(\bG,\bXi')$. Then
    \beq
        \label{e:fp-y-conc}
        \lt|\frac{\|\hat{\by}\|^2}{N} - 1\rt|\,,
        \frac{|(\hat{\bx},\hat{\by})|}{N},
        \frac{|(\hat{\by},\hat{\by}')|}{N} \le \delta
    \eeq
   with probability at least $1-e^{-cN}$.

\begin{proof}
    Since $\bG \mapsto \|\bar{\bx}\|$ is $L$-Lipschitz, $\|\bar{\bx}\|$ is $O(1)$-subgaussian by gaussian concentration of measure.
    Combined with $\E \|\bar{\bx}\|^2 = N$ this implies \eqref{e:fp-x-conc}.
    We next argue that the function
    \[
        f(\bG) = \E[ \|\hat{\by}(\bG,\bXi)\|^2 | \bG]^{1/2}
    \]
    is also $L$-Lipschitz. Indeed,
    \[
        \nabla f(\bG)
        = \frac{\E[ (\nabla_\bG \hat{\by}(\bG,\bXi))(\bG,\bXi) \hat{\by}(\bG,\bXi) | \bG]}{\E[ \|\hat{\by}(\bG,\bXi)\|^2 | \bG]^{1/2}}\,,
    \]
    and thus
    \[
        \|\nabla f(\bG)\|^2
        = \frac{\|\E[ (\nabla_\bG \hat{\by}(\bG,\bXi)) \hat{\by}(\bG,\bXi) | \bG]\|^2}{\E[ \|\hat{\by}(\bG,\bXi)\|^2 | \bG]}
        \le \frac{\E[ \|\nabla_\bG \hat{\by}(\bG,\bXi)\|_\op^2 \|\hat{\by}(\bG,\bXi)\|^2 | \bG]}{\E[ \|\hat{\by}(\bG,\bXi)\|^2 | \bG]}
        \le L^2\,.
    \]
    Thus $\E[ \|\hat{\by}\|^2 | \bG]^{1/2}$ is also $O(1)$-subgaussian.
    Combined with $\E \|\hat{\by}\|^2 = N$ this proves \eqref{e:fp-y-cond-conc}. 
    The first inequality of \eqref{e:fp-y-conc} follows similarly to \eqref{e:fp-x-conc}, applying gaussian concentration on $\bXi$ conditional on $\bG$.
    The second inequality of \eqref{e:fp-y-conc} holds by gaussian concentration because $\bXi \mapsto (\hat{\bx},\hat{\by}) / N$ is $L/N^{1/2}$-Lipschitz and $\E[\hat{\by} | \bG] = \bzero$. 
    The third inequality of \eqref{e:fp-y-conc} follows similarly, conditioning further on a realization of $\hat{\by}$ such that $\|\hat{\by}\| \le 2N^{1/2}$ (which holds with probability $1-e^{-cN}$).
\end{proof}
\end{lem}
Also recall from Lemma~\ref{l:wishart} that we can choose constants $C,c$ depending only on $\alpha$ such that
    \beq\label{e:wishart.bound.repeated}
    \max\bigg\{
    \frac{\|\bG\|_\op}{N^{1/2}},
    \frac{\|\bXi\|_\op}{N^{1/2}}
    \bigg\} \le C
    \eeq
with probability at least $1-e^{-cN}$.

\begin{lem}
    \label{l:random-YY-in-cD}
    Let $\YY = (\vBBb,\vvvv,\vUUu,\vwww)$ be the random vector defined in \eqref{e:freeprob.random.YY}, and $\mu(\bG,\omega)$ be defined in \eqref{e:lambda-mu-random}.
    For any $\epsilon'' > 0$ there exists $c = c(\alpha,\iota,\epsilon'') > 0$ such that
    \[
        \P\bigg(
            \inf\big\{
                \|\tilde \YY - \YY\|_2 : \tilde \YY \in \cD_{\mu(\bG,\omega)}
            \big\} \le \epsilon''
        \bigg)
        \ge 1 - e^{-cN}\,.
    \]
\begin{proof}
Throughout the proof, we abbreviate $(\lambda,\mu) = (\lambda(\bG,\omega),\mu(\bG,\omega))$, although we keep in mind that they are permitted to depend on $(\bG,\omega)$. We will construct $\tilde\YY = (\tilde\vBBb,\tilde\vvvv,\tilde\vUUu,\tilde\vwww) \in \cD_\mu$ such that $\|\tilde \YY - \YY\|_2 \le \epsilon''$ on an event with probability $1-e^{-cN}$.
    We set $(\tilde\vBBb,\tilde\vUUu) = (\vBBb,\vUUu)$, and $\tilde\vvvv = \max(\vvvv,\vUUu)$ where the maximum is taken entrywise.
    Note that for all $j\in [r]$, the Cauchy--Schwarz inequality gives
    \baln
        \vUUu_j - \vvvv_j
        &=
        \bigg(
            \frac{1}{|B_j|} \sum_{a\in B_j}
            \frac{(\bmeta^a,\bar{\bx})(\barbg^a,\hat{\by})}{N}
        \bigg)^2
        - \frac{1}{|B_j|} \sum_{a\in B_j}
        \frac{(\barbg^a,\hat{\by})^2}{N} \\
        &\le 
        \bigg(
            \frac{1}{|B_j|} \sum_{a\in B_j}
            \frac{(\barbg^a,\hat{\by})^2}{N}
        \bigg)
        \bigg(
            \frac{1}{|B_j|} \sum_{a\in B_j} \frac{(\bmeta^a,\bar{\bx})^2}{N} - 1
        \bigg)\,,
    \ealn
from which it follows that
    \[
        |\tilde\vvvv_j - \vvvv_j|
        = (\vUUu_j - \vvvv_j)_+
        \le \bigg(
            \frac{1}{|B_j|} \sum_{a\in B_j}
            \frac{(\barbg^a,\hat{\by})^2}{N}
        \bigg)
        \bigg|
            \frac{1}{|B_j|} \sum_{a\in B_j} \frac{(\bmeta^a,\bar{\bx})^2}{N} - 1
        \bigg|\,.
    \]
    By Lemma~\ref{l:fp-gc} and \eqref{e:wishart.bound.repeated}, with probability $1-e^{-cN}$,
    \[
        \frac{1}{|B_j|} \sum_{a\in B_j}
        \frac{(\barbg^a,\hat{\by})^2}{N}
        \le \frac{\|\hat{\by}\|^2 \|\bG\|_\op^2}{N|B_j|}
        \le \frac{(2N^{1/2})^2 (CN^{1/2})^2}{NM 
		\lambda_j}
        \le \frac{4C^2}{\alpha \iota}\,.
    \]
Set $\delta > 0$ small enough that
    \[
        [(1-\delta)^2,(1+\delta)^2] \subseteq \bigg[
            1 - \frac{\epsilon''\alpha\iota^2}{8C^2},
            1 + \frac{\epsilon''\alpha\iota^2}{8C^2}
        \bigg]\,.
    \]
    Since $(\bmeta^a : 1\le a\le M)$ are standard gaussians in $\bbR^N$ independent of $\bar{\bx}$ and $B_j$, we have
    \[
        \frac{1}{|B_j|} \sum_{a\in B_j} \frac{(\bmeta^a,\bar{\bx})^2}{N}
        \stackrel{d}{=}
        \frac{\|\bar{\bx}\|^2}{N} \cdot \frac{\|\bz\|^2}{|B_j|}
    \]
    where $\bz$ is a standard gaussian in $\bbR^{|B_j|}$.
    By Lemma~\ref{l:fp-gc} and a standard gaussian tail bound,
    \[
        \max\bigg\{
            \bigg|\frac{\|\bar{\bx}\|^2}{N} - 1\bigg|, 
            \bigg|\frac{\|\bz\|^2}{|B_j|} - 1\bigg|
        \bigg\} \le \delta
    \]
    with probability $1-e^{-cN}$.
It follows from the assumption on $\delta$ that
    \[
        \bigg| \frac{\|\bar{\bx}\|^2}{N} \cdot \frac{\|\bz\|^2}{|B_j|} - 1 \bigg|
        \le \frac{\epsilon''\alpha\iota^2}{8C^2}\,.
    \]
Combining the above bounds gives
    \beq\label{l:random-YY-in-cD-step1}
        |\tilde\vvvv_j - \vvvv_j|
        \le \frac{4C^2}{\alpha \iota}
        \cdot \frac{\epsilon''\alpha\iota^2}{8C^2}
        = \frac{\epsilon'' \iota}{2}\,.
    \eeq
    By a union bound, the above holds for all $j\in [r]$ with probability $1-e^{-cN}$.

    Next we set $\tilde\vwww$.
    If $(\mu,\vwww) \ge 1$ (resp. $(\mu,\vwww) < 1$), let $\tilde\vwww \in [0,+\infty)^s$ be any vector satisfying $\tilde\vwww \le \vwww$ (resp. $\tilde\vwww\ge \vwww$) with $(\mu,\tilde\vwww)=1$.
    Note that
    \[
        (\mu,\vwww)
        = \sum_{j=1}^s \mu_j \cdot 
        \frac{1}{|\BIsing_j|} \sum_{i\in \BIsing_j}
        (\be_i,\hat{\by})^2
        = \frac{1}{N} 
        \sum_{i=1}^N 
        (\be_i,\hat{\by})^2
        = \frac{\|\hat{\by}\|^2}{N}\,.
    \]
    Lemma~\ref{l:fp-gc} then ensures that, with probability $1-e^{-cN}$,
$|(\mu,\vwww) - 1| \le \epsilon'' \iota / 2$.
    On this event, using the assumption that $\tilde\vwww - \vwww$ is either all non-negative or all non-positive, we can bound 
    \[
        \epsilon'' \iota / 2
        \ge \Big|(\mu,\tilde\vwww - \vwww)\Big|
        \stackrel{\eqref{e:def-Lambda}}{\ge} \iota \|\tilde\vwww - \vwww\|_1
        \ge \iota \|\tilde\vwww - \vwww\|_2\,.
    \]
It follows that $\|\tilde\vwww - \vwww\|_2 \le \epsilon''/2$ with probability $1-e^{-cN}$. Combining with \eqref{l:random-YY-in-cD-step1} gives
    \[
        \|\tilde\YY - \YY\|_2
        \le \sum_{j=1}^r |\tilde\vvvv_j - \vvvv_j| + \|\tilde\vwww - \vwww\|_2
        \stackrel{\eqref{e:r-s-bd}}{\le}
        \frac{1}{\iota} \cdot \frac{\epsilon'' \iota}{2} + \frac{\epsilon''}{2}
        \le \epsilon''\,.
    \]
Since $\tilde\YY\in\cD_\mu$ by construction, this concludes the proof.
\end{proof}
\end{lem}

\begin{lem}
    \label{l:fp.surplus.approximation}
Recall the definition of $\Surplus(\YY;\lambda,\mu)$ from the statement of Proposition~\ref{p:convex.body.polytope.bound.new}. 
Recall $\cD_+$ from \eqref{e:cD-plus}.
Given any $\mu$, define 
    \[
        \cD_{+,\mu}(2) \equiv \Big\{
            \YY = (\BBb,\vvv,\UUu,\www) \in \cD_+ : (\mu, \www) \le 2
        \Big\}\,,
    \]
For $(\lambda,\mu) \in \Lambda$ and $\dot{\YY}, \ddot{\YY} \in \cD_{+,\mu}(2)$, we have
    \[\Big|\Surplus(\dot{\YY};\lambda,\mu) - \Surplus(\ddot{\YY};\lambda,\mu)\Big|
        \le
        \frac{4\|\Delta\YY\|_2
            +  6\|\Delta\YY\|_2^{1/2}}{\iota \min\{1,\alpha\}}\,,
    \]
where we abbreviate $\Delta\YY \equiv \dot{\YY} - \ddot{\YY}$. 
That is, $\YY \mapsto \Surplus(\YY;\lambda,\mu)$ is uniformly continuous in $\cD_{+,\mu}(2)$, uniformly in $(\lambda,\mu) \in \Lambda$.

\begin{proof}
We will bound separately the differences 
	\baln
	\Delta\Cost
	& \equiv\Cost(\dot{\YY};\lambda) - \Cost(\ddot{\YY};\lambda)\,,\\
	\Delta\Budget
	& \equiv\Budget(\dot{\YY};\mu) - \Budget(\ddot{\YY};\mu)\,.
	\ealn
Note that for any $\xxxx \in \bbR^r$, the Cauchy--Schwarz inequality gives
    \[
        \sum_{j=1}^r \lambda_j |\xxxx_j|
        \le \bigg(\sum_{j=1}^r \lambda_j |\xxxx_j|^2\bigg)^{1/2}
        \le \|\xxxx\|_2\,.
    \]
    Thus, abbreviating $\Delta\YY
    \equiv(\Delta\vBBb,\Delta\vvvv,\Delta\vUUu,\Delta\vwww)
    \equiv\dot{\YY}-\ddot{\YY}$, we have
    \baln
        |\Delta\Cost|
        &\le \sum_{j=1}^r \lambda_j \bigg\{
            |\Delta\vBBb_j|
            + |\Delta\vUUu_j|
            + \Big|
            (\dot{\vvvv}_j-\dot{\vUUu}_j)_+ - (\ddot{\vvvv}_j-\ddot{\vUUu}_j)_+
            \Big| 
            + 2 \Big|(\dot{\vvvv}_j-\dot{\vUUu}_j)_+^{1/2} - (\ddot{\vvvv}_j-\ddot{\vUUu}_j)_+^{1/2} 
            \Big|
        \bigg\} \\
        &\le \|\Delta\vBBb\|_2
        + \|\Delta\vUUu\|_2
        + \Big\|(\dot{\vvvv}-\dot{\vUUu})_+ - (\ddot{\vvvv}-\ddot{\vUUu})_+
        	\Big\|_2 
        + 2\Big\|(\dot{\vvvv}-\dot{\vUUu})_+^{1/2} - (\ddot{\vvvv}-\ddot{\vUUu})_+^{1/2}
        \Big\|_2\,,
    \ealn
    where the operations $x\mapsto x_+$ and $x\mapsto (x_+)^{1/2}$ are applied entrywise.
    Since $x\mapsto x_+$ is a contraction,
    \[
        \|(\dot{\vvvv}-\dot{\vUUu})_+ - (\ddot{\vvvv}-\ddot{\vUUu})_+\|_2
        \le \|(\dot{\vvvv}-\dot{\vUUu}) - (\ddot{\vvvv}-\ddot{\vUUu})\|_2
        \le \|\Delta\vvvv\|_2 + \|\Delta\vUUu\|_2\,.
    \]
Furthermore, since $|a^{1/2}-b^{1/2}| \le |a-b|^{1/2}$ for all $a,b\ge0$, for any vectors $\vec{a},\vec{b}\in[0,\infty)^r$ we can bound
	\[
	\Big(\|\vec{a}^{1/2}-\vec{b}^{1/2}\|_2\Big)^2
	=\sum_{j=1}^r \Big((a_j)^{1/2}-(b_j)^{1/2}\Big)^2
	\le \sum_{j=1}^r |a_j-b_j| = \|\vec{a}-\vec{b}\|_1\,.
	\]
Applying this estimate gives 
    \baln
        &\bigg(\Big\|(\dot{\vvvv}-\dot{\vUUu})_+^{1/2} - (\ddot{\vvvv}-\ddot{\vUUu})_+^{1/2}\Big\|_2\bigg)^2
        \le 
        \Big\|(\dot{\vvvv}-\dot{\vUUu})_+ - (\ddot{\vvvv}-\ddot{\vUUu})_+
        \Big\|_1 \\
        &\qquad \le 
            \|\Delta\vvvv\|_1
            + \|\Delta\vUUu\|_1 
            \le r^{1/2} \Big(
            \|\Delta\vvvv\|_2
            + \|\Delta\vUUu\|_2
        \Big)\,.
    \ealn
    Finally, crudely upper bounding 
$\|\Delta\vBBb\|_2,
        \|\Delta\vvvv\|_2,
        \|\Delta\vUUu\|_2$ by $\|\Delta\YY\|_2$,
    we obtain
    \[
        |\Delta\Cost|
        \le 4 \|\Delta\YY\|_2
        + 2
        \Big(2 r^{1/2} \|\Delta\YY\|_2 \Big)^{1/2} \,.
    \]
This gives a continuity estimate for the $\Cost$ function.

We next turn to the $\Budget$ function.
    For all $\YY  \in \cD_{+,\mu}(2)$, we have
    \[
        0\le \alpha \cdot \Budget(\YY;\mu) \le (\mu,\vwww) \le 2\,,
    \]
    where the intermediate step is by the Cauchy--Schwarz inequality.
    Therefore,
    \baln
        \alpha |\Delta\Budget| 
        &= \alpha^{1/2} \Big|
        \Budget(\dot{\YY};\mu)^{1/2} + 
        \Budget(\ddot{\YY};\mu)^{1/2}\Big| 
        \cdot 
        \alpha^{1/2} \Big|
        \Budget(\dot{\YY};\mu)^{1/2} 
        - \Budget(\ddot{\YY};\mu)^{1/2}
        \Big| \\
        &\le 2\sqrt{2} \cdot
        \alpha^{1/2} \Big|\Budget(\dot{\YY};\mu)^{1/2} 
        - \Budget(\ddot{\YY};\mu)^{1/2}\Big|\,.
    \ealn
    Using again that $|a^{1/2}-b^{1/2}| \le |a-b|^{1/2}$ for all $a,b\ge0$, we can bound
    \baln
        &\alpha^{1/2} \Big|\Budget(\dot{\YY};\mu)^{1/2} 
        - \Budget(\ddot{\YY};\mu)^{1/2}\Big|
        \le \sum_{j=1}^s \mu_j \Big|
        	(\dot{\vwww}_j)^{1/2} - (\ddot{\vwww}_j)^{1/2}\Big|
        \le \sum_{j=1}^s \mu_j |\dot{\vwww}_j - \ddot{\vwww}_j|^{1/2} \\
        &\qquad\le \Big(\sum_{j=1}^s \mu_j |\dot{\vwww}_j - \ddot{\vwww}_j|^2\Big)^{1/4}
        \le (\|\Delta\vwww\|_2)^{1/2} 
        \le  (\|\Delta\YY\|_2)^{1/2} \,.
    \ealn
Combining the above bounds gives
    \[
        \alpha |\Delta\Budget|
        \le 2\Big( 2\|\Delta\YY\|_2\Big)^{1/2}\,.
    \]
Altogether we obtain
    \baln
        &|\Surplus(\dot{\YY};\lambda,\mu) 
        - \Surplus(\ddot{\YY};\lambda,\mu)| 
        \le
        	\bigg\{
	4 \|\Delta\YY\|_2
        + 2
        \Big(2 r^{1/2} \|\Delta\YY\|_2 \Big)^{1/2} \bigg\}
        + \frac{2}{\alpha} 
        \Big( 2\|\Delta\YY\|_2\Big)^{1/2}\\
        &\qquad\stackrel{\eqref{e:r-s-bd}}{\le}
        \frac{4\|\Delta\YY\|_2
            + 6\|\Delta\YY\|_2^{1/2}
            }{\iota \min\{1,\alpha\}} \,.
    \ealn
This proves the claim. 
\end{proof}
\end{lem}

\begin{proof}[\hypertarget{proof:p.fp.given.projection}{Proof of Proposition~\ref{p:fp}, assuming Propositions~\ref{p:fp-projection}  and \ref{p:convex.body.polytope.bound.new}}]
 Let $\epsilon', C_\pol$ be the constants given by Proposition~\ref{p:convex.body.polytope.bound.new}, depending only on $\alpha,\iota,\epsilon$. Set $\epsilon'' > 0$ small enough that $\epsilon'' \le \min\{\epsilon'/2,1\}$ and
    \beq
        \label{e:fp-epsilon''}
        \frac{4\epsilon'' +  6(\epsilon'')^{1/2}}{\iota \min\{1,\alpha\}} 
        \le \frac{\epsilon}{2}\,.
    \eeq
Let $\YY$ be the random variable defined in \eqref{e:freeprob.random.YY}, and let
    \[
        \tilde\cE = \bigg\{
            \inf \Big\{
                \|\tilde \YY - \YY\| : \tilde \YY \in \cD_{\mu(\bG,\omega)}
            \Big\} \le \epsilon''
        \bigg\}
    \]
be the event from Lemma~\ref{l:random-YY-in-cD}. 
Recall from 
\eqref{e:E-lambda-mu} the definition of the event $\cE(\lambda,\mu)$.
We will argue that for any fixed $(\lambda,\mu) \in \Lambda$, we have
    \beq
        \label{e:p.fp.given.projection.goal}
        \bbP\bigg(
            \Surplus(\YY;\lambda,\mu)
            > \epsilon;
            \cE(\lambda,\mu)  \cap \tilde\cE 
        \bigg)
        \le e^{-cN}\,,
    \eeq
Proposition~\ref{p:convex.body.polytope.bound.new} implies the existence of closed half-spaces $H_j$, indexed by $1\le j\le n \le C_\pol$, such that
    \beq
        \label{e:p.fp.given.projection.polytope.approximation}
        S_\star(\lambda,\mu)[\epsilon'] \subseteq
        P = \cD_\mu \cap \bigcap_{j=1}^n H_j
        \subseteq 
        \bigg\{
            \tilde\YY \in \cD_\mu : 
            \Surplus(\tilde\YY;\lambda,\mu)
            \le  \frac\epsilon2
        \bigg\}\,.
     \eeq 
Consider any half-space $H_j$, which can be written as
    \begin{align*}
    H_j &= \Big\{
            \tilde\YY \in \tilde\cD : (\tilde\YY,\hY^j)_{\lambda,\mu} 
            \le F_{\star,\lambda,\mu}(\hY^j) + t^j\Big\}\\
   & = \Big\{
            \tilde\YY \in \tilde\cD : (\tilde\YY,\acute{Y}^j) 
            \le
            F_{\star,\lambda,\mu}(\hY^j) + t^j
    \Big\}\,,
    \end{align*}
where $t^j \in \bbR$, 
$\hY^j$ satisfies the normalization \eqref{e:hY-normalization}, and $\acute{Y}^j$ is the reweighting of 
$\hY^j$ defined by
    \[
        \acute Y^j = (\lambda \hB^j,\lambda\hv^j, 
        	\lambda \hU^j, 
        \mu\hw^j)\,,
    \]
which we note is a unit vector by \eqref{e:hY-normalization}.
Since $S_\star(\lambda,\mu)[\epsilon']$ is contained in $H_j$, it is straightforward to argue that we must have $t^j\ge\epsilon'$. Now define the smaller half-space
    \begin{align*}
    \acute{H}_j &= \bigg\{
            \tilde\YY \in \tilde\cD : (\tilde\YY,\hY^j)_{\lambda,\mu} \le F_{\star,\lambda,\mu}(\hY^j) + \frac{t^j}{2}
    \bigg\}
    = \bigg\{
            \tilde\YY \in \tilde\cD : (\tilde\YY,\acute{Y}^j) \le F_{\star,\lambda,\mu}(\hY^j) + \frac{t^j}{2}
    \bigg\}\,.
    \end{align*} 
Applying Proposition~\ref{p:fp-projection} with $\epsilon'/2$ in place of $\epsilon'$ gives
    \beq
        \label{e:p.fp.given.projection.proof.one.direction}
        \bbP\Big(
            \YY  \not\in 
            \acute{H}_j; 
            \cE(\lambda,\mu)
        \Big)
      \le
      \bbP\bigg(
            (\YY, \hY^j)_{\lambda,\mu} > F_{\star,\lambda,\mu}(\hY^j) + \frac{\epsilon'}{2};
            \cE(\lambda,\mu)
        \bigg)
        \le e^{-cN}\,.
    \eeq
Taking a union bound of \eqref{e:p.fp.given.projection.proof.one.direction} over $1\le j\le n$ (recalling $n\le C_\pol$ is bounded independently of $N$) shows
    \beq
        \label{e:p.fp.given.projection.proof.many.directions}
        \P\bigg(
            \YY \not\in
            \bigcap_{j=1}^n \acute{H}_j; \cE(\lambda,\mu) \cap \tilde\cE
        \bigg)
        \le 
        \sum_{j=1}^n
        \bbP\big(
            \YY \not\in \acute{H}_j;
            \cE(\lambda,\mu)
        \big)
        \le
        e^{-cN}\,.
    \eeq
    We will next show the inclusion of events
    \beq
        \label{e:p.fp.given.projection.event.inclusion}
        \bigg\{
            \Surplus(\YY;\lambda,\mu)
            > \epsilon;
            \cE(\lambda,\mu) \cap \tilde\cE
        \bigg\}
        \subseteq
        \bigg\{
            \YY \not\in
            \bigcap_{j=1}^n \acute{H}_j; \cE(\lambda,\mu) \cap \tilde\cE
        \bigg\}\,,
    \eeq
which is equivalent to showing that
    \beq
        \label{e:p.fp.given.projection.event.inclusion.goal}
        	\bigg\{
            \YY \in
            \bigcap_{j=1}^n \acute{H}_j; \cE(\lambda,\mu) \cap \tilde\cE
        \bigg\}
        \subseteq
        	\bigg\{
            \Surplus(\YY;\lambda,\mu)
            \le \epsilon;
            \cE(\lambda,\mu) \cap \tilde\cE
        \bigg\}\,.
	\eeq
Since $\cE(\lambda,\mu) \cap \tilde\cE$ holds, there exists $\tilde\YY \in \cD_\mu$ such that
$\|\tilde\YY - \YY\|_2 \le \epsilon'' \le \epsilon'/2$.
    Since $H_j$ is the closed $t^j/2$-neighborhood of $\acute{H}_j$ and $t^j/2 \ge \epsilon'/2$ for all $j$, it follows that
    \beq\label{e:p.fp.given.projection.event.inclusion.step1}
        \tilde\YY \in \cD_\mu \cap \bigcap_{j=1}^n H_j = P
        \stackrel{\eqref{e:p.fp.given.projection.polytope.approximation}}{\subseteq }
        \bigg\{\YY' \in \cD_\mu : 
            \Surplus(\YY';\lambda,\mu)
            \le  \frac\epsilon2
        \bigg\}\,.\eeq
Furthermore, writing $\tilde\YY = (\tilde\vBBb,\tilde\vvvv,\tilde\vUUu,\tilde\vwww)$, we have
    \[
        |(\mu,\tilde\vwww - \vwww)|
        \le \bigg(\sum_{j=1}^s \mu_j (\tilde\vwww_j - \vwww_j)^2
        	\bigg)^{1/2}
        \le \|\tilde\vwww - \vwww\|_2
        \le \epsilon'' \le 1\,.
    \]
    Since $\tilde\YY \in \cD_\mu$ implies $(\mu,\tilde\vwww) = 1$, we conclude that $\YY, \tilde\YY \in \cD_{+,\mu}(2)$.
    Then Lemma~\ref{l:fp.surplus.approximation} and \eqref{e:fp-epsilon''} imply
    \[
        \Big|\Surplus(\YY;\lambda,\mu) - \Surplus(\tilde\YY;\lambda,\mu)\Big| 
        \le \frac{4\epsilon'' +  6(\epsilon'')^{1/2}}
        	{\iota \min\{1,\alpha\}} 
        \le\frac\epsilon2\,.
    \]
Combining with \eqref{e:p.fp.given.projection.event.inclusion.step1} proves \eqref{e:p.fp.given.projection.event.inclusion.goal}, and \eqref{e:p.fp.given.projection.event.inclusion} follows. Finally, \eqref{e:p.fp.given.projection.proof.many.directions} and \eqref{e:p.fp.given.projection.event.inclusion} readily imply \eqref{e:p.fp.given.projection.goal} for each fixed $(\lambda,\mu) \in \Lambda$. Since $(\lambda(\bG,\omega),\mu(\bG,\omega)) \in \Lambda$ almost surely, and the total number of $(\lambda,\mu) \in \Lambda$ can be crudely bounded by $M^rN^s \le (MN)^{1/\iota}$, taking a union bound of \eqref{e:p.fp.given.projection.goal} over all 
$(\lambda,\mu) \in \Lambda$ gives
	\[
    \bbP\bigg(
        \Surplus(\YY;\lambda(\bG,\omega),\mu(\bG,\omega))
        > \epsilon;
        \tilde\cE
    \bigg) \le e^{-cN}\,.
\]
Taking a final union bound with Lemma~\ref{l:random-YY-in-cD} concludes the proof.
\end{proof}

We next note that we can reduce the function $F_{\star,\lambda,\mu}$ of \eqref{e:Fstar} to a simpler optimization problem: let
\[
    S_\diamond(\lambda,\mu) = \lt\{
        (\vvvv,\xxxx,\vwww) \in [0,+\infty)^r \times [0,+\infty) \times \cW_\mu :
        \sum_{j=1}^r \lambda_j ((\vvv_j)^{1/2} - 1)^2 + \xxxx
        \le \frac1\alpha \bigg(
            \sum_{j=1}^s
            \mu_j (\www_j)^{1/2}
        \bigg)^2
    \rt\}.
\]
This is a convex body (relative to the topology of $[0,+\infty)^r \times [0,+\infty) \times \cW_\mu$), because $\vvv \mapsto (\vvv^{1/2}-1)^2$ is convex while $\cW_\mu \ni \vwww \mapsto (\sum_{j=1}^s \mu_j (\www_j)^{1/2})^2$ is concave.
Further define 
\beq
    \label{eq:fp-hx}
    \hx 
    = \max_{j\in [r]} \bigg\{
        \max\Big\{\hv_j + \hU_j, \hB_j\Big\}
    \bigg\}\,.
\eeq
We then define the function
\beq
    \label{e:fp-Fdiam}
    F_{\diamond,\lambda,\mu}
    \equiv F_{\diamond,\lambda,\mu} (\hv,\hx,\hvw)
    = \sup\bigg\{
      (\hvv,\vvvv)_{\lambda} + \hx\xxxx + (\hvw,\vwww)_{\mu} :
      (\vvvv,\xxxx,\vwww) \in S_\diamond(\lambda,\mu)
    \bigg\}\,.
\eeq

\begin{lem}\label{l:fp-Fst-to-Fdiam}
    We have $F_{\star,\lambda,\mu}(\hY) = F_{\diamond,\lambda,\mu} (\hv,\hx,\hvw)$, where the correspondence between $\hY=(\hvB,\hvv,\hvU,\hvw)$ and $(\hv,\hx,\hvw)$ is given by \eqref{eq:fp-hx}.
\end{lem}
\begin{proof}
Recall that the notation $\lambda\hB$ denotes the coordinatewise product of the vectors $\lambda$ and $\hB$. Let $\circ$ denote vector concatenation. With this notation, we can express
    \baln
        F_{\star,\lambda,\mu}
        &=\sup\lt\{
            \begin{array}{l}
            (\hvB,\vBBb)_{\lambda}
            + (\hvv,\vvvv-\vUUu)_{\lambda}
            + (\hvv+\hvU,\vUUu)_{\lambda}
            + (\hvw,\vwww)_{\mu} : \\
            \qquad \YY\in \cD,\vwww\in\cW_\mu,
            \vvvv \ge \vUUu, \textup{ and } \\
            \qquad \displaystyle \sum_{i=1}^r \lambda_j \bigg(
                \BBb_j + \UUu_j
                + \Big((\vvv_j-\UUu_j)^{1/2}-1\Big)^2
            \bigg)
            \le \frac1\alpha
            \bigg(\sum_{j=1}^s \mu_j (\www_j)^{1/2}
                \bigg)^2
            \end{array}
        \rt\} \\
        &= \sup\lt\{
            \begin{array}{l}
            (\hvv,\vvvv')_{\lambda}
            + (\hvB \circ (\hvv+\hvU),\lambda \vBBb \circ \lambda \vUUu)
            + (\hvw,\vwww)_{\mu} : \\
            \qquad
            (\vBBb,\vvvv',\vUUu,\vwww) \in \cD,
           \vwww\in\cW_\mu, \textup{ and }\\
            \qquad \displaystyle
            \sum_{i=1}^r \Big((\vvv'_j)^{1/2}-1\Big)^2
            + \Big\|\lambda \vBBb \circ \lambda \vUUu\Big\|_1
            \le \frac1\alpha
            \bigg(\sum_{j=1}^s \mu_j (\www_j)^{1/2} \bigg)^2
            \end{array}
        \rt\} \,.
    \ealn
    In the final maximization, the objective and constraint are both linear in the vector $\lambda \vBBb \circ \lambda \vUUu$.
    Therefore, there is a maximizer where this vector has at most one nonzero entry, whose corresponding entry in $\hvB \circ (\hvv+\hvU)$ is maximal.
    This is equivalent to the maximization problem defining $F_{\diamond,\lambda,\mu}$.
\end{proof}

\textbf{The remainder of this section is devoted to the proofs of Propositions~\ref{p:fp-projection} and \ref{p:convex.body.polytope.bound.new}}, and is organized as follows.  In \S\ref{ss:polytope.appx.new} we present the
\hyperlink{proof:p.convex.body.polytope.bound.new}{proof of Proposition~\ref{p:convex.body.polytope.bound.new}}, which gives the polytope approximation of $S_\star(\lambda,\mu)$. In \S\ref{ss:k-orth-replica} through \S\ref{ss:fp-gaussian-comparison}, we prove that $F^{\ri}\equiv F_{\lambda,\mu}(\hY;\YY)$ can be approximately stochastically upper bounded by a sequence of modifications
$F^{\rii},\ldots,F^{\rvi}$. 
At the end of \S\ref{ss:fp-gaussian-comparison}, we \hyperlink{proof:p.fp-projection}{reduce Proposition~\ref{p:fp-projection} to a simplified statement, Proposition~\ref{p:fp-v-estimate},} which states that $F_{\diamond,\lambda,\mu}$ is an approximate high-probability upper bound for $F^{\rvi}$. The remainder of the section is then devoted to the \hyperlink{proof:p.fp-v-estimate}{proof of Proposition~\ref{p:fp-v-estimate}}, which implies Proposition~\ref{p:fp-projection}.

\subsection{Upper bound on number of facets in polytope approximation}
\label{ss:polytope.appx.new}

This subsection is devoted to the proof of Proposition~\ref{p:convex.body.polytope.bound.new}.
Throughout this subsection we fix $\alpha,\iota,\epsilon$; all parameters below can implicitly depend on these constants. We denote
	\beq\label{e:lam.mu.delta.ball}
        B_\delta(\lambda,\mu)
        \equiv \bigg\{
            (\lambda',\mu') \in \Lambda
            : \|(\lambda',\mu') - (\lambda,\mu)\|_\infty \le \delta
        \bigg\}\,.\eeq
 For $\epsilon_1,\epsilon_2 \ge 0$ define
\beq\label{e:relaxed.S.star}
    S_\star(\lambda,\mu;\epsilon_1,\epsilon_2) = \bigg\{
        \YY \in \cD_\mu :
        \Surplus(\YY;\lambda,\mu) + \epsilon_1 \alpha^2 \iota^3 \sum_{i=1}^r (\BBb_i)^2
        \le\epsilon_2
    \bigg\}\,.\eeq
We first provide a crude bound on a set of the form \eqref{e:relaxed.S.star}.

\begin{lem} 
\label{l:p.polytope.bd.local.aux}
For any $(\lambda,\mu) \in \Lambda$, we have 
	\[ S_\star\bigg(\lambda,\mu;
		\frac\epsilon5,
		\frac\epsilon2\bigg) 
	\subseteq \bigg[0,
		\frac8{\iota_0}
		\bigg]^{3r+s}\,.\]
where $\iota_0\equiv\iota \min\{1,\alpha\}$.
\begin{proof} Consider any $\YY$ belonging to the set on the left-hand side of the above display. By the Cauchy--Schwarz inequality,
    \[
        \Budget(\YY;\mu) \le \frac{1}{\alpha} (\mu,\vwww) = \frac{1}{\alpha}\,,
    \]
 Recall from the discussion after Proposition~\ref{p:fp}
 that  $\epsilon \le 1$. Then, for any $j\in [r]$, 
assumption~\eqref{e:def-Lambda} gives
    \[\max\bigg\{
            \vBBb_j,\vUUu_j,
            \Big((\vvvv_j-\vUUu_j)^{1/2} - 1\Big)^2
        \bigg\}
        \le \frac{\Cost(\YY;\lambda)}{\iota}
        \le \frac{\Budget(\YY;\mu) + \epsilon/2}{\iota}
        \le \frac{2}{\iota\min\{1,\alpha\}}\,.
    \]
    This further implies (recalling $\iota,\epsilon \le 1$)
    \[
        \vvvv_j \le \vUUu_j + \bigg(1 + 
        	\bigg( \frac{2}{\iota \min\{1,\alpha\}}\bigg)^{1/2}
	\bigg)^2
        \le \frac{2
        	+(1 + \sqrt{2})^2}{\iota \min\{1,\alpha\}} 
        \le \frac{8}{\iota \min\{1,\alpha\}}\,.
    \]
Finally, we can use the constraint $(\mu,\vwww) = 1$ to crudely bound
    \[
        \vwww_j \le \frac{1}{\mu_j} \stackrel{\eqref{e:def-Lambda}}{\le} \frac{1}{\iota} \le \frac{8}{\iota \min\{1,\alpha\}}\,.
    \] This concludes the proof.
\end{proof}
\end{lem}

The next two lemmas compare the set defined by
\eqref{e:relaxed.S.star}
 with the set $S_\star(\lambda,\mu)$ from \eqref{e:fp.ideal.feasible.set}.

\begin{lem}
    \label{l:S.star.in.S.star.eps}
    For any $\epsilon_1 > 0$ and $(\lambda,\mu) \in \Lambda$, we have $S_\star(\lambda,\mu) \subseteq S_\star(\lambda,\mu;\epsilon_1,\epsilon_1)$.
\begin{proof}
    Let $\YY \in S_\star(\lambda,\mu)$.
    Then for any $i \in [r]$,
    \[
        \lambda_i \BBb_i
        \le \Cost(\YY;\lambda)
        \le \Budget(\YY;\mu)
        = \frac{1}{\alpha}
        \bigg(
            \sum_{j=1}^s
            \mu_j
            (\www_j)^{1/2}
        \bigg)^2
        \le \frac{1}{\alpha}
        (\mu,\www)
        = \frac{1}{\alpha}\,,
    \]
where the last equality uses that $\www \in \cW_\mu$. Since $(\lambda,\mu) \in \Lambda$, this implies $\BBb_i \le 1/(\alpha \iota)$. It follows that
    \[
        \epsilon_1 \alpha^2 \iota^3 \sum_{i=1}^r (\BBb_i)^2
        \le \epsilon_1 \iota r
        \stackrel{\eqref{e:r-s-bd}}{\le} \epsilon_1\,,
    \]
and the conclusion follows.
\end{proof}
\end{lem}

\begin{lem}\label{l:D.mu.continuity.w} 
For $\mu,\mu'$ satisfying condition~\eqref{e:def-Lambda},
if we have $\vwww'\ge0$ with $(\mu',\vwww')=1$,
then we can find $\vwww\ge0$ with $(\mu,\vwww)=1$ 
such that
	\[
	\|\vwww-\vwww'\|_1
	\le \frac{\|\mu-\mu'\|_\infty}{\iota^2}
	\,.\]
Note if $\YY' \equiv (\vBBb,\vvvv,\vUUu,\vwww)'\in\cD_{\mu'}$, and we find $\vwww$ as above, then
    $\YY=(\vBBb',\vvvv',
    \vUUu',\vwww)\in\cD_\mu$.

\begin{proof}
We start with $\vwww'\ge0$ satisfying $(\mu',\vwww')=1$, which by condition~\eqref{e:def-Lambda} implies that
	\[\|\vwww'\|_1
	\le \frac{(\mu',\vwww')}{\iota}
	=\frac{1}{\iota}\,.
	\]
Abbreviating $\delta\equiv \|\mu-\mu'\|_\infty$, it follows that
	\[
	\Big|(\mu,\vwww')-1\Big|
	=  \Big|(\mu-\mu',\vwww')\Big|
	\le \|\mu-\mu'\|_\infty \|\vwww'\|_1
	\le \frac{\delta}{\iota}\,.
	\]
We then choose $\vwww$ to satisfy $(\mu,\vwww)=1$, which we can do with a similar construction as in the proof of Lemma~\ref{l:random-YY-in-cD}: if $(\mu,\vwww')\ge1$
(resp.\ $(\mu,\vwww')<1$), let $\vwww$ be any vector with $\vwww\le\vwww'$ (resp.\ $\vwww\ge\vwww'$) such that $(\mu,\vwww)=1$. Then, similarly as in the proof of Lemma~\ref{l:random-YY-in-cD}, we can use the above to bound
	\[\frac{\delta}{\iota}
	\ge
	\Big|(\mu,\vwww')-1\Big|
	=
	\Big|(\mu,\vwww'-\vwww)\Big|
	\ge \iota\|\vwww'-\vwww\|_1\,,
	\]
which concludes the proof.
\end{proof}
\end{lem}

\begin{lem}\label{l:S.star.cts.dH} 
For any $\epsilon_1,\epsilon_2\ge0$, 
and for $(\lambda,\mu) \in \Lambda$, 
the mapping
$(\lambda,\mu)\mapsto S_\star(\lambda,\mu;\epsilon_1,\epsilon_2)$,
is continuous with respect to the Hausdorff metric $d_H$.

\begin{proof} Fix $\epsilon_1,\epsilon_2\ge0$.
Recall the notation \eqref{e:lam.mu.delta.ball}, and suppose
$(\lambda',\mu')\in\Lambda$ and 
 $(\lambda,\mu)\in B_\delta(\lambda',\mu')$. Given $\YY'\in S_\star(\lambda',\mu';\epsilon_1,\epsilon_2)$, we will construct $\YY\in S_\star(\lambda,\mu;\epsilon_1,\epsilon_2)$ which is $o_\delta(1)$-close to $\YY'$. We set $\vwww$ as in Lemma~\ref{l:D.mu.continuity.w}, so that we have the bound $\|\vwww-\vwww'\|_1 \le \delta/\iota^2$. Recalling Lemma~\ref{l:p.polytope.bd.local.aux}, we see that
 $\vwww,\vwww'$ are confined to a compact domain $[0,9/\iota_0]^s$.
The $\Budget$ function is uniformly continuous on the compact domain
$[0,9/\iota_0]^{3r+s}\times\Lambda$, so
	\[\sup\bigg\{
	\begin{array}{l}
	|\Budget(\YY;\mu)
	-\Budget(\YY';\mu')|
	:\YY,\YY'\in[0,9/\iota_0]^{3r+s},
	(\lambda,\mu),(\lambda',\mu')\in\Lambda,\\
	\qquad\qquad\qquad\qquad\qquad\qquad\qquad
	\|\mu-\mu'\|_\infty \le \delta,
	\|\YY-\YY'\|_1 \le \delta/\iota^2
	\end{array}
	\bigg\} \le o_\delta(1)\,,
	\]
where $o_\delta(1)$ denotes an error tending to zero as $\delta\to0$. 
For the purposes of this proof, let us abbreviate
	\begin{align*}
	c(\YY;\lambda)
	&\equiv \Cost(\YY;\lambda)
	+ \epsilon_1\alpha^2\iota^2(\|\vBBb\|_2)^2\,,\\
	g(\YY;\lambda,\mu)
	&\equiv c(\YY;\lambda) - \Budget(\YY;\mu) - \epsilon_2\,.
	\end{align*}
Then we can similarly bound
	\[\sup\bigg\{\begin{array}{l}
	|c(\YY';\lambda)
	-c(\YY';\lambda')|
	: \YY'\in[0,9/\iota_0]^{3r+s}, 
	(\lambda,\mu),(\lambda',\mu')\in\Lambda,\\
	\qquad\qquad\qquad\qquad\qquad\qquad
	\|\lambda-\lambda'\|_\infty \le \delta,
	\end{array}
	\bigg\}
	 \le o_\delta(1)\,.
	\]
Combining the above bounds gives
	\beq\label{e:surplus.small.error.bound}
	c(\YY';\lambda)
	-\Budget(\YY;\mu) - \epsilon_2
	\le g (\YY';\lambda',\mu') + o_\delta(1)
	\le o_\delta(1)\,,
	\eeq
uniformly over 
$(\lambda,\mu),(\lambda',\mu')$ in $\Lambda$ within distance $\delta$,
and
$\YY,\YY'\in [0,9/\iota_0]^{3r+s}$ within distance $\delta/\iota^2$.
If the above left-hand side is in fact $\le0$ then we are done, since we can then set
$\YY=(\vBBb',\vvvv',\vUUu',\vwww)$
 to achieve
$g(\YY;\lambda,\mu)\le 0$. If the above left-hand side is
 strictly positive, then it remains to choose 
 $\YY=(\vBBb,\vvvv,\vUUu,\vwww)$
 close to 
 $\YY(0)\equiv(\vBBb',\vvvv',\vUUu',\vwww)$ to make $g(\YY;\lambda,\mu)\le 0$.
In this case, consider the path in $\cD_\mu$ given by the line segment joining $\YY(0)\equiv(\vBBb',\vvvv',\vUUu',\vwww)$
to $\YY(1)\equiv(\vec{0},\vvvv'-\vUUu',\vec{0},\vwww)$. If we parametrize this segment as a linear function $\YY(t)$ of $t\in[0,1]$, then along the segment we have
	\[\Cost(\YY(t);\lambda)
	-\Cost(\YY(0);\lambda)
	=-(\lambda,\vBBb'+\vUUu')t
	\equiv -A(\YY')t\,.
	\]
Recall from \eqref{e:surplus.small.error.bound} that at the start of the line segment we have $g(\YY(0);\lambda,\mu) \le o_\delta(1)$.  Let $\eta = g(\YY(0);\lambda,\mu)$. We then divide into several cases:
\begin{itemize}
\item If  $A(\YY')\ge\eta^{1/2}$, then we can choose $t_0= \eta / A(\YY') \le \eta^{1/2}$ such that
	\[
	g(\YY(t_0);\lambda,\mu) 
	\le  \eta^{1/2} - A(\YY') t_0 
	\le 0\,,
	\]
and (since $\YY,\YY' \in [0,9/\iota_0]^{3r+s}$) $\|\YY(t_0)-\YY(0)\|_2 \lesssim \eta^{1/2}$. In this case, we can conclude by setting $\YY = \YY(t_0)$.
\item Otherwise, if 
$A(\YY')\le \eta^{1/2}$, then
assumption~\eqref{e:def-Lambda} implies that 
\[
    \max\Big(\|\vBBb'\|_1,\|\vUUu'\|_1\Big) \le \frac{\eta^{1/2}}{\iota}\,,
\]
so that we have 
\[
    \|\YY(1)-\YY(0)\|_2
    \le \|\vBBb'\|_2 + 2\|\vUUu\|_2
    \le r^{1/2} \Big(\|\vBBb'\|_1 + 2\|\vUUu\|_1\Big)
    \stackrel{\eqref{e:r-s-bd}}{\le} \frac{\eta^{1/2}}{\iota^{3/2}}\,.
\]
If
$g(\YY(1);\lambda,\mu) \le0$
then we can conclude by setting $\YY = \YY(1)$.
\item If the above bound does not hold, then we still have
	\[
	g(\YY(1 );\lambda,\mu)
	\le
	g(\YY(0 );\lambda,\mu)
	= \eta 
	\,.
	\]
In this case, we consider the path in $\cD_\mu$ given by the line segment joining $\YY(1)=(\vec{0},\vvvv'-\vUUu',\vec{0},\vwww)$ to $\YY_{\min}=(\vec{0},\vec{1},\vec{0},\vwww)$, which we note satisfies
 	\[g(\YY_{\min};\lambda,\mu)
	\le
	\Cost(\YY_{\min};\lambda) = 0\,.
	\]
If $\|\YY(0)-\YY_{\min}\| \le  \eta^{1/3}$, then we can conclude by setting $\YY=\YY_{\min}$.

\item Otherwise, the vector $x\equiv \vvvv'-\vUUu'$ in $\YY(1)$, must have some coordinate $j$ such that 
\[
    |x_j-1| \ge \frac{\eta^{1/3}}{\sqrt{r}} \stackrel{\eqref{e:r-s-bd}}{\ge} \eta^{1/3} \iota^{1/2}\,.
\]
In this last case, we take $\YY(1)$ and decrease $|x_j-1|$ by $\eta^{1/3} \iota^{1/2}$ to obtain $\YY$. This results in
 	\[g(\YY;\lambda,\mu)
	\le 
	g(\YY(0);\lambda,\mu)
	- \iota \cdot (\eta^{1/3} \iota^{1/2})^2 
	\le \eta - \eta^{2/3} \iota^2 
	\le0\,,
	\]
where the last inequality holds for sufficiently small $\eta$.
\end{itemize}
Recall $\eta = o_\delta(1)$. It follows that for sufficiently small $\delta$, in all cases we can find
 $\YY\in\cD_\mu$ within $O(\eta^{1/3}) o_\delta(1)^{1/3}$ of $\YY(0)=(\vBBb',\vvvv',\vUUu',\vwww)$ with 
 	\[g(\YY;\lambda,\mu)
	\le 0\,.\]
Thus $\YY\in S_\star(\lambda,\mu;\epsilon_1,\epsilon_2)$ 
 with $\|\YY-\YY'\|_2 \le \delta'$, where the last inequality can be guaranteed by choosing $\delta$ sufficiently small. 
The claim follows.
\end{proof}
\end{lem}

\begin{cor}\label{c:S.star.cts.nbd.dH}
 For any $\epsilon_1,\epsilon_2,\epsilon\ge0$, 
and for $(\lambda,\mu) \in \Lambda$, 
the mapping
$(\lambda,\mu)\mapsto S_\star(\lambda,\mu;\epsilon_1,\epsilon_2)[\epsilon]$,
is continuous with respect to the Hausdorff metric $d_H$.

\begin{proof}
Fix $\epsilon_1,\epsilon_2,\epsilon\ge0$ and $\delta'>0$. By Lemma~\ref{l:S.star.cts.dH}, there exists $\delta>0$ such that if 
$(\lambda,\mu)\in\Lambda$ and 
$(\lambda',\mu')\in B_\delta(\lambda,\mu)$, then the sets
\begin{align*}
K&\equiv S_\star(\lambda,\mu;\epsilon_1,\epsilon_2)\,,\\
K'&\equiv S_\star(\lambda',\mu';\epsilon_1,\epsilon_2)
\end{align*}
lie within Hausdorff distance $\delta'$ of one another.
For any $\hat{\YY}\in K[\epsilon]$, there exists $\YY\in K$ with
$\|\hat{\YY}-\YY\|_2\le\epsilon$, and then there exists $\YY'\in K'$ with
$\|\YY-\YY'\|_2\le\delta'$. By Lemma~\ref{l:D.mu.continuity.w}, 
since $\hat{\YY}\in\cD_\mu$,
there exists $\tilde{\YY}\in\cD_{\mu'}$
with $\|\hat{\YY}-\tilde{\YY}\|_1\le\delta/\iota^2$. Then
	\[
	\|\tilde{\YY}-\YY'\|_2
	\le\|\tilde{\YY}-\hat{\YY}\|_2
	+ \|\hat{\YY} - \YY\|_2
	+ \|\YY - \YY'\|_2
	\le \epsilon + \delta' + \frac{\delta}{\iota^2}
	\le \epsilon + 2\delta'\,,
	\]
where the last inequality can be arranged by taking $\delta$ small enough (depending on $\delta'$).
Since $\YY'$ and $\tilde{\YY}$ both lie in $\cD_{\mu'}$, it follows by convexity of $\cD_{\mu'}$ that $\tilde{\YY}$ lies within distance
$2\delta'$ of a vector $\hat{\YY}' \in K'[\epsilon]$. We then finally have
	\[
	\|\hat{\YY}-\hat{\YY}'\|_2
	\le \|\hat{\YY}-\tilde{\YY}\|_2
		+ \|\tilde{\YY}-\hat{\YY}'\|_2
	\le\frac{\delta}{\iota^2} + 2\delta'
	\le 3\delta'\,.
	\]
This implies that $K[\epsilon]$ and $K'[\epsilon]$ lie within Hausdorff distance $3\delta'$ of one another, and the claim follows.
\end{proof}
\end{cor}

\begin{lem}
    \label{l:S.star.nbd.in.S.star.eps}
    There exists $\epsilon' > 0$ such that 
    \[S_\star(\lambda,\mu)[\epsilon'] 
    \subseteq S_\star\bigg(\lambda,\mu;
    	\frac\epsilon5, \frac\epsilon3
    	\bigg)\]
for all $(\lambda,\mu) \in \Lambda$.

\begin{proof}
    We first prove the claim for a fixed $(\lambda,\mu) \in \Lambda$, where $\epsilon'$ for now is permitted to depend on $(\lambda,\mu)$.
Suppose for contradiction that for all $\epsilon'>0$ we have
    \[ \emptyset\ne
    S_\star(\lambda,\mu)[\epsilon'] \cap 
	S_\star\bigg(\lambda,\mu;\frac\epsilon5, \frac\epsilon3\bigg)^c 
	\subseteq
    S_\star(\lambda,\mu)[\epsilon'] \cap 
\overline{S_\star\bigg(\lambda,\mu;\frac\epsilon5, \frac\epsilon{3}\bigg)^c}\,,
    \]
    where $\overline{S}$ denotes the closure of $S$ relative to the topology of $\cD_\mu$. This forms a nested family of nonempty compact sets indexed by $\epsilon' > 0$, so the intersection
    \[
        \bigcap_{\epsilon'>0} \bigg\{
        S_\star(\lambda,\mu)[\epsilon'] \cap \overline{S_\star\bigg(\lambda,\mu;
        	\frac\epsilon5, \frac\epsilon{3}\bigg)^c}\bigg\}
        = S_\star(\lambda,\mu) \cap
        	\overline{S_\star\bigg(\lambda,\mu;
			\frac\epsilon5, \frac\epsilon{3}\bigg)^c}
    \]
is also nonempty. That is, we can find $\YY\in S_\star(\lambda,\mu)$ with
	\[
	\Surplus(\YY;\lambda,\mu)
	+ \frac\epsilon5 \alpha^2\iota^3 (\|\vBBb\|_2)^2
	\ge \frac{\epsilon}{3}\,.
	\]
This contradicts Lemma~\ref{l:S.star.in.S.star.eps}, which implies 
	\[
	\Surplus(\YY;\lambda,\mu)
	+ \frac\epsilon5 \alpha^2\iota^3 (\|\vBBb\|_2)^2
	\le \frac{\epsilon}{5}\,.
	\]
This proves that for any fixed $(\lambda,\mu)\in\Lambda$, the claim holds for some $\epsilon'>0$ that may depend on $(\lambda,\mu)$. 
Next, for any $(\lambda,\mu)\in\Lambda$, let 
\[\epsilon'(\lambda,\mu)
\equiv \max\bigg\{
\epsilon':
S_\star(\lambda,\mu)[\epsilon'] \subseteq 
S_\star\bigg(\lambda,\mu;\frac\epsilon5,\frac\epsilon3\bigg)\bigg\}\,,\] 
where the maximum exists because both sets are closed. We will prove below that $\epsilon'(\lambda,\mu)$ is a continuous function of $(\lambda,\mu)\in\Lambda$, and we have just shown that it is strictly positive on $\Lambda$. It follows that
\[
        \inf\bigg\{\epsilon'(\lambda,\mu) : (\lambda,\mu) \in \Lambda\bigg\} > 0\,,
    \]
and we may take $\epsilon'$ equal to this infimum to conclude.
 
We now justify the above claim that $\epsilon'(\lambda,\mu)$ is a  continuous function of $(\lambda,\mu)\in\Lambda$. (Note for our purposes it suffices to show that $\epsilon'(\lambda,\mu)$ is lower semicontinuous in $(\lambda,\mu)$, but it is not much more difficult to verify that it is in fact continuous.) To this end, 
fix $(\lambda,\mu)\in\Lambda$ and abbreviate $\epsilon'=\epsilon'(\lambda,\mu)$, so that we have
	\[S_\star(\lambda,\mu)[\epsilon']
	\subseteq
	S_\star\bigg(\lambda,\mu;\frac\epsilon5,\frac\epsilon3\bigg)\,.
	\]
To this end, note that the function
    \beq\label{e:f.strictly.convex.surplus.relaxation}
        f(\YY;\lambda,\mu) 
        \equiv \Surplus(\YY;\lambda,\mu) 
        	+ \frac{\epsilon}{5} \alpha^2 \iota^3 (\|\vBBb\|_2)^2\,,
    \eeq
is \emph{strictly} convex on $\cD_\mu$, as the $(\|\BBb\|_2)^2$ term ensures strict convexity in $\vBBb$. As a result, for any sufficiently small $\eta'>0$, there exists $\delta'>0$ such that
	\[S_\star(\lambda,\mu)[\epsilon'-\delta']
	\subseteq
	S_\star\bigg(\lambda,\mu;\frac\epsilon5,
		\frac\epsilon3-\eta'\bigg)\,.
	\]
Let $(\lambda',\mu')\in B_\delta(\lambda,\mu)$. By Corollary~\ref{c:S.star.cts.nbd.dH}, we have
	\[
	d_H\bigg(
	S_\star(\lambda,\mu)[\epsilon'-\delta'],
	S_\star(\lambda',\mu')[\epsilon'-\delta']
	\bigg)
	\le o_\delta(1)\,.
	\]
Since the function $f$ is uniformly continuous on the compact domain
$[0,9/\iota_0]^{3r+s}\times \Lambda$,
we conclude 
	\[
	\sup\bigg\{
	 f(\YY';\lambda,\mu)
	 : \YY' \in S_\star(\lambda',\mu')[\epsilon'-\delta']
	 \bigg\}
	\le \frac\epsilon3-\eta' + o_\delta(1)
	\le\frac\epsilon3\,,
	\]
where the last inequality can be arranged by choosing $\delta$ small enough. This implies
	\[
	S_\star(\lambda',\mu')[\epsilon'-\delta']
	\subseteq
	S_\star\bigg(\lambda',\mu';\frac\epsilon5,
		\frac\epsilon3
		\bigg)\,,
	\]
and so $\epsilon'(\lambda',\mu)\ge\epsilon'(\lambda,\mu)-\delta'$. Since $\delta$ does not depend on $\mu$, by symmetry we will also have $\epsilon'(\lambda,\mu)\ge\epsilon'(\lambda',\mu')-\delta'$.
 Since $\delta'$ was an arbitrary positive constant, shows that $\epsilon'(\lambda,\mu)$ is a  continuous function of $(\lambda,\mu)\in\Lambda$. As explained above, the lemma follows.
\end{proof}
\end{lem}

For the rest of this subsection, let $\epsilon'$ be given by Lemma~\ref{l:S.star.nbd.in.S.star.eps}. Say $K \subseteq \cD_\mu$ is a \textbf{$\cD_\mu$-convex body} if it is convex and compact with nonempty interior relative to the topology of $\cD_\mu$.
\begin{lem}
\label{l:approximating-polytope-exists}
For any $(\lambda,\mu) \in \Lambda$, 
there exists a closed $\cD_\mu$-polytope $P$ such that
    \beq\label{e:approximating-polytope}
       K_1 \equiv S_\star(\lambda,\mu)[\epsilon']
        \subseteq P
        \subseteq S_\star\bigg(\lambda,\mu;
        	\frac\epsilon5,\frac\epsilon2\bigg) 
	\equiv K_2\,,
    \eeq
where $K_1,K_2$ are both
 $\cD_\mu$-convex bodies.

\begin{proof}
Note that $\Cost(\YY;\lambda)$, $(\|\BBb\|_2)^2$, and $-\Budget(\YY;\mu)$ are all convex functions of $\YY\in\cD_\mu$. (This uses the fact that for $\YY \in \cD_\mu$, the input to the $(\cdot)_+$ in $\Cost(\YY;\lambda)$ is always nonnegative.) This immediately implies that $K_2$ is a $\cD_\mu$-convex body. It also implies that $S_\star(\lambda,\mu)$ is a $\cD_\mu$-convex body, and therefore so is its $\epsilon'$-neighborhood $K_1$.

We now turn to the construction of the polytope $P$. Lemma~\ref{l:S.star.nbd.in.S.star.eps} gives
    \beq\label{e:S.star.nbd.in.S.star.eps.consequence} 
        K_1 \subseteq S_\star\bigg(\lambda,\mu;
        		\frac\epsilon5,\frac{\epsilon}{3}
		\bigg)
        \subseteq (K_2)^\circ\,,
    \eeq
    where $(K_2)^\circ$ denotes the interior of $K_2$ relative to the topology of $\cD_\mu$. We let $C$ denote the boundary of $K_2$ relative to the topology of $\cD_\mu$, that is,
	\[
	C\equiv \partial K_2
	\equiv K_2 \setminus (K_2)^\circ\,.
	\]
Note that $C$ is a compact set in the topology of $\cD_\mu$. It follows from \eqref{e:S.star.nbd.in.S.star.eps.consequence} that $K_1$ and $C$ are disjoint.

Recall from \eqref{e:D.mu.affine.closure} that we use $\mathcal{F}_\mu$ to denote the affine closure of $\cD_\mu$.
For each $\YY\in C$, by the (strict) hyperplane separating theorem, there exists an open half-space $U_\YY$ of $\mathcal{F}_\mu$ such that $\YY\in U_\YY$ and $U_\YY\cap K_1=\emptyset$. Then $V_\YY\equiv \cD_\mu \cap U_\YY$ is open in the $\cD_\mu$-topology, and we have an open cover
	\[
	C\subseteq
	\bigcup_{\YY\in C} V_\YY\,.
	\]
Since $C$ is compact, this must admit a finite subcover, so there exists finite $S\subseteq C$ such that
	\[
	C\subseteq
	\bigcup_{\YY\in S} V_\YY\,.
	\]
Let $A_\YY \equiv \mathcal{F}_\mu\setminus U_\YY$, a closed half-space of $\mathcal{F}_\mu$. It follows that we can express $A_\YY = \mathcal{F}_\mu \cap H_\YY$ where $H_\YY$ is a closed half-space of the full space $\R^{3r+s}$. We then define the closed $\cD_\mu$-polytope 
	\[
	P\equiv
	\cD_\mu\cap
	\bigg( \bigcap_{\YY\in S} H_\YY
	\bigg)\,.
	\]
Since each $U_\YY$ is disjoint from $K_1$, each $H_\YY$ must contain $K_1$, so we have $K_1\subseteq P$.
On the other hand, since 
	\[
	\cD_\mu \setminus P
	= \cD_\mu \cap 
	\bigg( \bigcup_{\YY\in S} U_\YY
	\bigg)
	\supseteq C\,,
	\]
we see that $P$ does not intersect $C$. Since $P$ is convex, it must be fully contained in either $(K_2)^\circ$ or
$(\cD_\mu\setminus K_2)^\circ$.
Since we already saw that $K_1\subseteq P$, the former must hold, so $P\subseteq K_2$ as claimed.\end{proof}
\end{lem}

We now define $C(\lambda,\mu)$ to be the minimum number of facets of a closed $\cD_\mu$-polytope $P$ such that \eqref{e:approximating-polytope} holds; so that Lemma~\ref{l:approximating-polytope-exists} gives  $C(\lambda,\mu)<\infty$ for every $(\lambda,\mu) \in \Lambda$. \textbf{The main technical result of this subsection is the following:}

\begin{ppn}\label{p:polytope.bd.local} For any $(\lambda,\mu) \in \Lambda$, there exist $\delta,\tilde C > 0$  such that
    \[
        \sup\bigg\{
            C(\lambda',\mu') : 
            (\lambda',\mu') \in B_\delta(\lambda,\mu) 
        \bigg\} \le \tilde C\,.
    \]
Here the parameters $\delta,\tilde C$ are permitted to depend on $\lambda,\mu$
   in addition to $\alpha,\iota,\epsilon$.
   
\end{ppn}
The \hyperlink{proof:p.polytope.bd.local}{proof of Proposition~\ref{p:polytope.bd.local}} appears below, after some preparatory lemmas. 
In what follows, we will rely on a classical result on the stability of linear programs under perturbation. The following is a simplified statement:
\begin{thm}[\cite{MR51275}; see also \cite{MR317760}] \label{t:hoffman}
 Let $A$ and $C$ be real matrices of dimensions $m\times n$ and $q\times n$ respectively. Let $b$ and $d$ denote vectors in $\R^m$ and $\R^q$ respectively. Define
	\begin{align*}
	S_{A,C}(b,d)
	&\equiv \Big\{x\in\R^n: Ax\le b, Cx=d\Big\}\,,\\
	F(A,C)
	&\equiv\Big\{(b,d) \in\R^{m+q}:
	S_{A,C}(b,d) \ne\emptyset\Big\}\,.\end{align*}
Let $(b,d)\in F(A,C)$, so that 
$S_{A,C}(b,d)$ is nonempty. If $\hat{x}$ is any point in $\R^n$, then there is a point $\bar{x}\in S_{A,C}(b,d)$ with
	\[
	\|\bar{x}-\hat{x}\|_2
	\le \sigma(A,C)\left\|
	\begin{pmatrix}
	(A\hat{x}-b)_+\\ C\hat{x}-d
	\end{pmatrix}\right\|\,,
	\]
where $\sigma(A,C)$ is a finite number depending only on $A,C$.
\end{thm}

This has the following immediate consequences:
\begin{cor}[\cite{MR51275}; see also \cite{MR317760}] 
\label{c:hoffman}
Take $A,b,C,d$ as in Theorem~\ref{t:hoffman} with $(b,d)\in F(A,C)$.
Let $B(M_0)$ denote the ball of radius $M_0$ in $\R^n$, centered at the origin.
\begin{enumerate}[a.]
\item \label{i:hoffman.rhs} Suppose $(b',d')\in F(A,C)$, and
suppose further that the (nonempty) sets
$S_{A,C}(b,d)$ and $S_{A,C}(b',d')$ are both contained in $B(M_0)$.
Then
	\[d_H\Big(S_{A,C}(b,d),S_{A,C}(b',d')\Big)
	\le c\left\| \begin{pmatrix}
	b-b'\\d-d'\end{pmatrix}\right\|_2
	\]
where $c$ is a finite constant depending only on $A,C,M_0$.
\item \label{i:hoffman.semicty}
Now suppose we have $(A,b,C,d)'$ of the same dimensions,
with $(b',d')\in F(A',C')$. 
Suppose further that the (nonempty) sets
$S\equiv S_{A,C}(b,d)$ and $S'\equiv S_{A',C'}(b',d')$ are both contained in $B(M_0)$.
Then, for any $x'\in S'$ there exists $x\in S$ with
	\[\|x-x'\|_2
	\le c \left\| \begin{pmatrix} A-A'\\ C-C'\\
	b-b'\\d-d'\end{pmatrix}\right\|_2\,,
	\]
where $c$ is a finite constant depending only on $A,C,M_0$, and on the right-hand side we treat $A$ and $C$ as vectors of length $mn$ and $qn$ respectively.
\end{enumerate}

\begin{proof}
Suppose $x'\in S'$. Then
	\[
	Ax'-b
	=[A'+(A-A')] x'- [b' + (b-b')]
	\le (A-A')x'-(b-b')\,,
	\]
from which it follows that
	\[
	\|(Ax'-b)_+\|
	\le \| [(A-A')x'-(b-b')]_+\|
	\le \|A-A'\|_\textup{op} M_0 + \|b-b'\|\,.
	\]
Similarly, we can bound
	\[
	\|Cx'-d\|
	\le \|C-C'\|_\textup{op} M_0 + \|d-d'\|\,.
	\]
The operator norm is bounded by the Frobenius norm, so the conclusion of part~\eqref{i:hoffman.semicty} follows.
Part \eqref{i:hoffman.rhs} then follows from part~\eqref{i:hoffman.semicty}.
\end{proof}
\end{cor}

\begin{lem}\label{p:polytope.bd.local.step.one}
In the setting of Proposition~\ref{p:polytope.bd.local}, 
we can construct a family of half-spaces $(\tilde{H}_j : j\in[\tilde{n}])$
such that for $\delta_0>0$ small enough, we have
	 \begin{align}
        \label{e:polytope-tilde1}
        P\supseteq 
        \tilde{P}_\mu &\equiv \cD_\mu \cap \bigcap_{j=1}^{\tilde n} \tilde H_j 
        \supseteq S_\star(\lambda,\mu)[\epsilon']\,, \quad \text{and} \\
        \label{e:polytope-tilde2}
        	\emptyset \ne 
        \tilde{P}_{\mu'} &\equiv \cD_{\mu'} \cap \bigcap_{j=1}^{\tilde n} \tilde H_j 
        \subseteq 
        S_\star\bigg(\lambda',\mu';
        	\frac\epsilon5,\frac\epsilon2-\delta_0
        	\bigg)\,.
    \end{align}
for all $(\lambda',\mu') \in B_{\delta_0}(\lambda,\mu)$. 

\begin{proof}Let $n = C(\lambda,\mu)$. Let $(H_j:j\in[n])$ be a collection of closed half-spaces of $\bbR^{3r+s}$ such that
	\[
        K_1=S_\star(\lambda,\mu)[\epsilon']
        \subseteq P
        \equiv \cD_\mu \cap \bigcap_{j=1}^n H_j
        \subseteq S_\star\bigg(\lambda,\mu;
        	\frac\epsilon5,\frac\epsilon2\bigg)
        =K_2\,.\]
Abbreviate $u_\mu\equiv (0^{3r},\mu)$, and recall from \eqref{e:D.mu.affine.closure} that the affine space $\mathcal{F}_\mu$ is defined by the constraint $(\YY,u_\mu)=1$. In the above definition of $P$, we can assume without loss that none of the $H_j$ have normal vector parallel to $u_\mu$: such $H_j$ are redundant because of the intersection with $\cD_\mu$.

Recall that the function $f$ from \eqref{e:f.strictly.convex.surplus.relaxation}
is strictly convex on $\cD_\mu$.
We also have
	\[
	P\subseteq K_2
	=\bigg\{\YY \in \cD_\mu : f(\YY;\lambda,\mu) \le 	
	\frac\epsilon2
		\bigg\}\,.
	\]
The strict convexity of $f$ implies that the set 
    \[
        \cQ \equiv \bigg\{
        	\YY \in P : f(\YY;\lambda,\mu) = \frac\epsilon2
		\bigg\}
    \]
must be a subset of the extreme points of $P$ (the points that cannot be written as a convex combination of two other points in $P$). Therefore, $\cQ$ is a discrete set.
    Moreover, Lemma~\ref{l:S.star.nbd.in.S.star.eps} implies
    \[
        S_\star(\lambda,\mu)[\epsilon'] \subseteq
        \bigg\{\YY \in P : f(\YY;\lambda,\mu) \le 
        	\frac\epsilon3
        	\bigg\}\,,
    \]
and thus $\cQ \cap S_\star(\lambda,\mu)[\epsilon'] = \emptyset$. Enumerate the points of $\cQ$ as $\YY_1,\ldots,\YY_{\hat{n}}$, and for each $j\in[\hat{n}]$ let $\hat{H}_j$ be a closed half-space in $\R^{3r+s}$ whose interior contains $K_1$, and which does not contain $\YY_j$. We then construct the collection $(\tilde H_j : j\in[\tilde{n}])$ by appending these half-spaces to $(H_j:j\in[n])$, so that $\tilde{n}=n+\hat{n}$. Note that $\hat{H}_j$ strictly separates $\YY_j\in\cD_\mu$ from $K_1\subseteq\cD_\mu$, so the normal vector of $\hat{H}_j$ cannot be parallel to $u_\mu$. 

Recall from \eqref{e:polytope-tilde1} and \eqref{e:polytope-tilde2}
that we denote 
	\[
	\tilde{P}_{\mu'}
	\equiv \cD_{\mu'} \cap \bigcap_{j=1}^{\tilde{n}} \tilde{H}_j\,.
	\]
From the above construction combined with Lemma~\ref{l:p.polytope.bd.local.aux}, we have
	\[
	K_1\subseteq \tilde{P}_\mu \subseteq P\setminus\cQ
	\subseteq P
	\subseteq K_2
	\subseteq \cD_\mu
	\cap \bigg[0,\frac{8}{\iota_0}\bigg]^{3r+s}
	 \subseteq \mathcal{F}_\mu
	 \cap \bigg[0,\frac{8}{\iota_0}\bigg]^{3r+s}\,.
	\]
In particular, we know that $K_1$ has nonempty interior in the $\mathcal{F}_\mu$-topology, so $\tilde{P}_\mu$ also has nonempty interior in the $\mathcal{F}_\mu$-topology; we denote this set $\mathrm{int}(\tilde{P}_\mu;\mathcal{F}_\mu)$. If $\YY\in \mathrm{int}(\tilde{P}_\mu;\mathcal{F}_\mu)$, then for all $\YY'$ orthogonal to $u_\mu$ with $\|\YY'\|$ small enough, we must have $\YY+\YY' \in\tilde{P}_\mu$. We claim that any such $\YY$ must also belong to the interior of $\tilde{H}_j$ in the full space $\R^{3r+s}$. Indeed, in the full space $\R^{3r+s}$, we can express
	\[\tilde{H}_j
	= \Big\{\YY''\in\R^{3r+s}
		: (\YY'',v_j) \le c_j\Big\}\,.
	\]
Suppose for contradiction that $\YY$ lies on the boundary of $\tilde{H}_j$ relative to the full space $\R^{3r+s}$. Then we must have $(\YY,v_j)=c_j$, as well as
$(\YY+\YY',v_j)\le c_j$ for all sufficiently small $\YY'$ orthogonal to $u_\mu$. This is only possible if $v_j$ is parallel to $u_\mu$, yielding a contradiction. This proves that $\YY$ is indeed in the interior of $\tilde{H}_j$ in the full space $\R^{3r+s}$. Since $\YY\in\tilde{P}_\mu$, we conclude that there must exist a set $W_\YY$, a small open ball in the topology of $\R^{3r+s}$, such that 
	\beq
	\label{e:tilde.P.mu.interior.claim}
	\YY \in W_\YY \subseteq 
	\tilde{A}\equiv
	\bigg[0,\frac{9}{\iota_0}\bigg]^{3r+s}
	\cap
	\bigcap_{j=1}^{\tilde{n}} \tilde{H}_j\,.\eeq
It follows that there exists $\delta^+>0$ small enough  such that for all $(\lambda',\mu')\in \Lambda_+\equiv B_{\delta^+}(\lambda,\mu)$, the set $\cD_{\mu'}$ also intersects the small open ball $W_\YY\subseteq\tilde{A}$. For the remainder of the proof, we restrict attention to $(\lambda',\mu')\in\Lambda_+$. The above implies that for all $(\lambda',\mu')\in\Lambda_+$, the set
	\[
	\tilde{A}_{\mu'}
	\equiv \cD_{\mu'} \cap \tilde{A}
	= \cD_{\mu'} \cap 
	\bigg[0,\frac{9}{\iota_0}\bigg]^{3r+s}
	\cap
	\bigcap_{j=1}^{\tilde{n}} \tilde{H}_j
	= \bigg[0,\frac{9}{\iota_0}\bigg]^{3r+s}
	\cap \tilde{P}_{\mu'}
	\]
is nonempty. 

Finally we will argue that for sufficiently small $\delta \in(0,\delta^+)$, it holds for all $(\lambda',\mu')\in B_\delta(\lambda,\mu)$ that
	\[
	\tilde{P}_{\mu'}
	\subseteq S_\star\bigg(
		\lambda',\mu';
		\frac\epsilon5,
		\frac\epsilon2-\delta\bigg)
	=\bigg\{\YY'\in \cD_{\mu'}
	: f(\YY';\lambda',\mu') \le \frac\epsilon2-\delta
	\bigg\}\,.
	\]
The result will follow by renaming $\delta$ to $\delta_0$. To show this, we start by recalling that $\tilde{P}_\mu\subseteq P\setminus\cQ \subseteq[0,8/\iota_0]^{3r+s}$, so
	\[\max\bigg\{f(\YY;\lambda,\mu)
	:\YY\in \tilde{A}_\mu
	= \bigg[0,\frac{9}{\iota_0}\bigg]^{3r+s}
	\cap \tilde{P}_{\mu}
	=\tilde{P}_{\mu}
	\bigg\} < \frac{\epsilon}{2}\,.
	\]
The function $f$ is uniformly continuous on the compact domain
	\[\bigg[0,\frac{9}{\iota_0}\bigg]^{3r+s}\times \Lambda_+\,.\]
If $(\lambda',\mu')\in\Lambda_+$, then we argued above that $\tilde{A}_{\mu'}$ is nonempty. For any
$\YY'\in\tilde{A}_{\mu'}$, applying Corollary~\ref{c:hoffman}\ref{i:hoffman.semicty} gives the existence of
$\YY\in\tilde{A}_\mu$ with
	\[\|\YY-\YY'\|_2
	\le c\|\mu-\mu'\|_\infty\,,
	\]
where $c$ is a finite constant depending on $\mu,\tilde{A}_\mu=\tilde{P}_\mu,\iota_0$. We then use the uniform continuity of $f$ to conclude that there exists $\delta>0$ small enough so that
	\[
	\max\bigg\{
	f(\YY';\lambda',\mu')
	:(\lambda',\mu')\in B_\delta(\lambda,\mu),
	\YY' \in \tilde{A}_{\mu'}
	\bigg\} \le \frac\epsilon2-\delta\,.
	\]
This implies
	\[
	\tilde{A}_{\mu'}
	= \bigg[0,\frac{9}{\iota_0}\bigg]^{3r+s}
	\cap \tilde{P}_{\mu'}
	\subseteq
	S_\star\bigg(
		\lambda',\mu';
		\frac\epsilon5,
		\frac\epsilon2-\delta\bigg)
	\subseteq \bigg[0,\frac{8}{\iota_0}\bigg]^{3r+s}\,,
	\]
where the last containment is again by Lemma~\ref{l:p.polytope.bd.local.aux}. It follows that $\tilde{P}_{\mu'}$ does not intersect 
	\[\bigg[0,\frac{9}{\iota_0}\bigg]^{3r+s}
	\setminus \bigg[0,\frac{8}{\iota_0}\bigg]^{3r+s}\,.
	\]
Since we argued above that $\tilde{A}_{\mu'}$ is nonempty, we conclude that $\tilde{P}_{\mu'}$ must intersect 
$[0,8/\iota_0]^{3r+s}$. Since $\tilde{P}_{\mu'}$ is a convex subset of $\cD_{\mu'}$, we conclude that it must be fully contained in 
$[0,8/\iota_0]^{3r+s}$. This concludes the proof.
\end{proof}
\end{lem}

\begin{proof}[\hypertarget{proof:p.polytope.bd.local}{Proof of Proposition~\ref{p:polytope.bd.local}}]
Lemma~\ref{p:polytope.bd.local.step.one} gives that for some $\delta_0>0$, there exist $(\tilde H_j:j\in[\tilde{n}])$, a family of closed half-spaces in $\R^{3r+s}$, 
such that for all $(\lambda',\mu') \in B_{\delta_0}(\lambda,\mu)$, the sets $\tilde{P}_\mu$, $\tilde{P}_{\mu'}$ defined in \eqref{e:polytope-tilde1}, \eqref{e:polytope-tilde2} satisfy the conditions
    \begin{align}
        \label{e:polytope-local-step1-guarantee.first}
        \tilde{P}_\mu &\supseteq S_\star(\lambda,\mu)[\epsilon']\,, \\
        \emptyset \ne
        \tilde{P}_{\mu'} &\subseteq S_\star(\lambda',\mu';\epsilon/5,\epsilon/2-\delta_0)
        \label{e:polytope-local-step1-guarantee}\,.
    \end{align}
We will now construct closed half-spaces 
$\acute{H}_j$, indexed by $j\in[n']$, such that
    \beq
        \label{e:polytope-prime-goal}
        S_\star(\lambda',\mu')[\epsilon']
        \subseteq P_{\mu'}
        \equiv \cD_{\mu'} \cap \bigcap_{j=1}^{n'} \acute{H}_j
        \subseteq S_\star\bigg(\lambda',\mu';
        	\frac\epsilon5,\frac\epsilon2
        	\bigg)\,.
    \eeq
for some $\delta > 0$
and for all $(\lambda',\mu') \in B_\delta(\lambda,\mu)$. 
This implies the result with the above choice of $\delta$ and $\tilde C = n'$. To this end, for $\delta_1 \ge 0$, let $\tilde H_j[\delta_1]$ denote the closed $(\delta_1)$-neighborhood of $\tilde{H}_j$ in $\R^{3r+s}$. Then, for each $(\lambda',\mu') \in B_{\delta_0}(\lambda,\mu)$, let $\delta_1(\lambda',\mu')$ be the largest $\delta_1$ such that
    \beq\label{e:def.delta.one.of.vlam.vmu}
    \tilde{P}_{\mu'}\{\delta_1\}
    \equiv
        \cD_{\mu'} \cap \bigcap_{j=1}^{\tilde n} \tilde H_j[\delta_1] \subseteq 
	S_\star\bigg(\lambda',\mu';
		\frac\epsilon5,
		\frac\epsilon2
		\bigg)
	=\bigg\{
	\YY'\in\cD_{\mu'}
	: f(\YY;\lambda',\mu')
	\le\frac\epsilon2
	\bigg\}\,.
    \eeq
Note that $(\tilde{P}_{\mu'}\{\delta_1\}:\delta_1\ge0)$ defines a nested family of compact sets, with $\tilde{P}_{\mu'}\{0\}=\tilde{P}_{\mu'}\ne\emptyset$. It follows from Corollary~\ref{c:hoffman}\ref{i:hoffman.rhs} that for $\delta_1\ge0$, the mapping
$\delta_1\mapsto\tilde{P}_{\mu'}\{\delta_1\}$
is continuous with respect to Hausdorff distance $d_H$. This implies that
the largest value $\delta_1(\lambda',\mu')$
for \eqref{e:def.delta.one.of.vlam.vmu} to hold is well-defined. Since \eqref{e:polytope-local-step1-guarantee}
holds, and $f$ is a uniformly continuous function as noted above, we must have $\delta_1(\lambda',\mu') > 0$ for each $(\lambda',\mu') \in B_{\delta_0}(\lambda,\mu)$. Moreover, we argue that $\delta_1(\lambda',\mu')$ is a lower semicontinuous function of $(\lambda',\mu')$: indeed, suppose $\bar{\delta}_1=\delta_1(\lambda',\mu')$, so that \eqref{e:def.delta.one.of.vlam.vmu} holds. Suppose $\YY$ is any maximizer of $f$ on the set $\tilde{P}_{\mu'}\{\bar{\delta}_1\}$, and suppose $\YY$ lies in the interior of every $\tilde{H}_j[\bar{\delta}_1]$. Then, as long as $\vBBb'$ is small enough, the vector
$\YY'=\YY+(\vBBb',0,0,0)$ also lies in $\tilde{P}_{\mu'}\{\bar{\delta}_1\}$.
But we have $f(\YY')>f(\YY)$, contradicting the assumption that $\YY$ was a maximizer of $f$ on $\tilde{P}_{\mu'}\{\bar{\delta}_1\}$. Therefore any maximizer of $f$ on $\tilde{P}_{\mu'}\{\bar{\delta}_1\}$ must lie on the boundary of at least one of the half-spaces $\tilde{H}_j[\bar{\delta}_1]$. As a result, for any $\eta_1>0$ there must exist $\eta_2>0$ such that 
	\[\tilde{P}_{\mu'}\{\bar{\delta}_1-\eta_1\}=
	\cD_{\mu'}
	\cap \bigcap_{j=1}^{\tilde n}
	\tilde H_j[\bar{\delta}_1-\eta_1]
	\subseteq
	 S_\star\bigg(\lambda',\mu';
		\frac\epsilon5,
		\frac\epsilon2
		-\eta_2
		\bigg)\,.
	\]
Then, if $(\lambda'',\mu'')$ is sufficiently close to $(\lambda',\mu')$, we claim that
	\beq\label{e:lampp.mupp.containmnet}
	\tilde{P}_{\mu''}\{\bar{\delta}_1-\eta_1\}=
	\cD_{\mu''}
	\cap \bigcap_{j=1}^{\tilde n}
	\tilde H_j[\bar{\delta}_1-\eta_1]
	\subseteq
	S_\star\bigg(\lambda'',\mu'';
		\frac\epsilon5,
		\frac\epsilon2
		\bigg)\,.\eeq
To see that \eqref{e:lampp.mupp.containmnet} holds, recall that
 $\tilde{P}_{\mu''}\{\bar{\delta}_1-\eta_1\}
 \supseteq\tilde{P}_{\mu''}\ne\emptyset$. It follows from Corollary~\ref{c:hoffman}\ref{i:hoffman.semicty} that for $\YY''\in \tilde{P}_{\mu''}\{\bar{\delta}_1-\eta_1\}$, we can find $\YY'\in \tilde{P}_{\mu'}\{\bar{\delta}_1-\eta_1\}$ with
 	\[\|\YY'-\YY''\|_2
	\le c\|\mu'-\mu''\|_\infty\,,\]
where $c$ is a finite constant depending on $\mu',\tilde{P}_{\mu'}\{\bar{\delta}_1-\eta_1\},\iota_0$.
Since $f$ is uniformly continuous, we can then take $(\lambda'',\mu'')$ sufficiently close to $(\lambda',\mu')$ to guarantee
$f(\YY'';\lambda'',\mu'')\le\epsilon/2$, so that
$\YY''$ belongs to the right-hand side of \eqref{e:lampp.mupp.containmnet} as claimed.
Altogether this implies $\delta_1(\lambda'',\mu'')\ge\delta_1(\lambda',\mu') -\eta_1$. Since $\eta_1>0$ was arbitrary, this proves that $\delta_1(\lambda',\mu')$ is indeed a lower semicontinuous function of $(\lambda',\mu')$, as claimed. The set $B_{\delta_0}(\lambda,\mu)$ is compact, so we therefore conclude
    \[\inf\bigg\{\delta_1(\lambda',\mu') : (\lambda',\mu') \in B_{\delta_0}(\lambda,\mu)\bigg\} > 0\,.
    \]
Let $\delta_1$ denote this infimum.
    We then finally let $\acute{n} = \tilde n$, and define $(\acute{H}_j:j\in[\acute{n}])$ by setting $\acute{H}_j = \tilde H_j[\delta_1]$.
    Thus for all $(\lambda',\mu') \in B_{\delta_0}(\lambda,\mu)$, for $P_{\mu'}$ as defined in \eqref{e:polytope-prime-goal}, we have
    \[
        P_{\mu'} = \cD_{\mu'} \cap \bigcap_{j=1}^{\acute{n}} \acute{H}_j \subseteq S_\star\bigg(
        	\lambda',\mu';
	\frac\epsilon5,\frac\epsilon2
	\bigg)\,.
    \]
This proves the right-most inclusion of \eqref{e:polytope-prime-goal}.
    Next, consider
    \baln
        \tilde S &= \bigcap_{j=1}^{\tilde n} \tilde H_j 
        \supseteq \cD_\mu \cap \tilde S
        	=
        	\tilde{P}_\mu 
	\stackrel{\eqref{e:polytope-local-step1-guarantee.first}}{\supseteq}
	S_\star(\lambda,\mu)[\epsilon']
	\,, \\
        \acute{S}
        &= \bigcap_{j=1}^{\acute{n}} \acute{H}_j
        = \bigcap_{j=1}^{\tilde n} \tilde H_j[\delta_1]
    \ealn
as subsets of $\bbR^{3r+s}$, and note that they satisfy
$d(\tilde S, \bbR^{3r+s} \setminus \acute{S})\ge\delta_1$. 
Now recall from Corollary~\ref{c:S.star.cts.nbd.dH} that if $(\lambda',\mu')\in B_\delta(\lambda,\mu)$, then 
	\[
	d_H\bigg(
	S_\star(\lambda,\mu)[\epsilon'],
	S_\star(\lambda',\mu')[\epsilon']
	\bigg) \le o_\delta(1)\,.
	\]
This implies that we can take $\delta$ small enough so that 
$S_\star(\lambda',\mu')[\epsilon'] \subseteq \acute{S}$. 
Since $S_\star(\lambda',\mu')[\epsilon'] \subseteq \cD_{\mu'}$ by definition, we finally conclude
    \[
        S_\star(\lambda',\mu')[\epsilon'] \subseteq \cD_{\mu'} \cap \acute{S} = P_{\mu'}\,,
    \]
    proving the left-most inclusion of \eqref{e:polytope-prime-goal}. This concludes the proof.
\end{proof}

\begin{proof}[\hypertarget{proof:p.convex.body.polytope.bound.new}{Proof of Proposition~\ref{p:convex.body.polytope.bound.new}}]
Again let $\epsilon'$ be given by Lemma~\ref{l:S.star.nbd.in.S.star.eps}. A locally bounded function on a compact set is uniformly bounded, so we have
	\beq\label{e:polytope.facets.uniform.bd}
	C_\pol
	\equiv\sup\Big\{
            C(\lambda,\mu) :
            (\lambda,\mu) \in \Lambda
        \Big\}
        < \infty\,.
	\eeq
Note that
 and $\epsilon',C_\pol$ depend on only $\alpha,\iota,\epsilon$.
    We have thus shown that for all $(\lambda,\mu) \in \Lambda$, there exists a closed $\cD_\mu$-polytope $P$ with at most $C_\pol$ facets such that
    \[
        S_\star(\lambda,\mu)[\epsilon']
        \subseteq P
        \subseteq S_\star\bigg(
        \lambda,\mu;\frac\epsilon5,\frac\epsilon2
        \bigg)\,.
    \]
    The second inclusion further implies that for all $\YY \in P$,
    \[
        \Cost(\YY;\lambda)
        \le \Cost(\YY;\lambda) + \frac{\epsilon}{5} \alpha^2 \iota^3 
        	(\|\vBBb\|_2)^2
        \le \Budget(\YY;\mu) + \frac{\epsilon}{2}\,,
    \]
which proves the claim.
\end{proof}

\subsection{Reduction to optimization over many orthogonal replicas}
\label{ss:k-orth-replica}

The remainder of this section builds towards the 
\hyperlink{proof:p.fp-projection}{proof of Proposition~\ref{p:fp-projection}}.
From here onwards, we observe the following convention on implicit constants:

\begin{rmk}\label{r:fp-implicit-constants} 
Throughout what follows, the parameters $\alpha, L, \iota, \epsilon'$ will be considered fixed. The notations $O(\cdot)$ and $\lesssim$ will hide constants that depend on these parameters, which will be uniform over $r,s$ satisfying \eqref{e:r-s-bd}, $\lambda,\mu \in \Lambda$
as defined by \eqref{e:def-Lambda}, and $\hY$ satisfying the normalization condition \eqref{e:hY-normalization}. For example, by \eqref{e:r-s-bd} and \eqref{e:fp-normalization-bd}, all factors of $r,s$ and $\hB_j, \hv_j, \hU_j, \hw_j$ can be absorbed into uniform implicit constants. 
As a result, we shall fix 
 $r,s,\lambda,\mu,\hY$ and often drop these parameters from our notation.  Thus, we write $F$, $F_\star$, and $F_\diamond$ for the quantities \eqref{e:1dproj}, \eqref{e:Fstar}, and \eqref{e:fp-Fdiam} respectively; and $\cD,\cW$ for $\cD_\mu,\cW_\mu$. We will often state that events hold with probability at least $1-e^{-cN}$; when we do, the constant $c$ can depend on the above parameters $\alpha,L,\iota,\epsilon'$, and is similarly uniform over $r,s,\lambda,\mu,\hY$. The dependence of $c$ on any other parameter will always be made explicit.
\end{rmk}

 Throughout what follows, we also assume that the event $\cE(\lambda,\mu)$ of \eqref{e:E-lambda-mu} holds almost surely, which is without loss of generality: if $\cE(\lambda,\mu)$ does not hold almost surely, then we can modify the random partition outside the event $\cE(\lambda,\mu)$, such that the modified partition satisfies $\cE(\lambda,\mu)$ almost surely. Towards the \hyperlink{proof:p.fp-projection}{proof of Proposition~\ref{p:fp-projection}}, this subsection has two main results:
\begin{itemize}
\item Proposition~\ref{p:fp-f-to-i} shows that the quantity $F$ of \eqref{e:1dproj} 
is approximately stochastically dominated by a quantity $F^{\rii}(k)$, defined by \eqref{e:F.rii} below, in which the random vector $\hat{\by} = \hat{\by}(\bG, \bXi)$ is replaced by optimizing an average over $k$ orthogonal replicas $\by^{1:k}$. 
\item Proposition~\ref{p:fp-i-to-ii} then shows that the quantity $F^{\rii}(k)$ is in turn stochastically dominated by a similar but simpler quantity $F^{\riii}(k)$, defined by \eqref{e:F.iii} below. \end{itemize} 
Recalling \eqref{e:1dproj}, let $F^{\ri}(\bG,\bXi,\by)$ be the sum of the functions
    \begin{align}
    \label{e:F.B}
    F_B(\bG,\bXi,\by)
    &\equiv
    (\hvB,\vBBb)_{\lambda}
    =\sum_{j=1}^r
    \frac{\hB_j}{MN|B_j|}
    \bigg(\sum_{a\in B_j} \bmeta^a, \by \bigg)^2\,,\\
    \label{e:F.v}
    F_v(\bG,\by)
    &\equiv
    (\hvv,\vvvv)_{\lambda}
    =
    \sum_{j=1}^r
    \frac{\hv_j}{M}
    \sum_{a\in B_j}
    \frac{(\barbg^a,\by)^2}{N}\,,\\
    \label{e:F.U}
    F_U(\bG,\bXi,\by)
    &\equiv (\hvU,\vUUu)_{\lambda}
    =\sum_{j=1}^r
    \frac{\hU_j}{M|B_j|}
    \bigg(
        \sum_{a\in B_j}
        \frac{(\bmeta^a,\bar{\bx})(\barbg^a,\by)}{N}
    \bigg)^2 \,, \\
    \label{e:F.w}
    F_w(\by)
    &\equiv
    (\hvw,\vwww)_{\mu}
    = \sum_{j=1}^s
    \frac{\hw_j}{N}
    (\Pi_j, \by^{\otimes 2})\,,
\end{align}
where $\Pi_j$ denotes orthogonal projection onto the coordinates in $\BIsing_j$. The sets $B_j$ and $\BIsing_j$ depend on the random partition $\hat{\cB}(\bG,\omega)$  defined in \eqref{e:fp-partition}, but we for now suppress this dependence. Then note that the quantity $F$ of \eqref{e:1dproj} can be written as
    \beq\label{e:freeprob.F.equals.F.ri}
    F = F^{\ri}(\bG,\bXi,\hat{\by}(\bG,\bXi))\,.
    \eeq
Let $k$ be a large constant we will set later, which will depend on only the parameters $\alpha,L,\iota,\epsilon'$ indicated in Remark~\ref{r:fp-implicit-constants}. Define
    \[
    S_\perp(k) = \left\{
    \begin{array}{l}
    \by^{1:k}
    \equiv 
    (\by^1,\ldots,\by^k) \in (\R^N)^k : \,
        \|\by^\ell\|^2 = N, \,
        (\by^\ell, \bar{\bx}) = 0, \,\\
    \qquad\qquad\qquad\qquad\qquad
    (\by^\ell,\by^{\ell'}) = 0 \,
        \text{for all $\ell \neq \ell'$}
        \end{array}
    \right\}\,.
\]
Let $\bXi^1,\ldots,\bXi^k$ be $k$ i.i.d. samples of $\bXi$, and define the $(\bG,\bXi^1,\ldots,\bXi^k)$-measurable random variable
    \beq\label{e:F.rii}
    F^{\rii}(k)
    =
    \sup\bigg\{
    F^{\ri}(\bG,
        \bXi^{1:k},\by^{1:k}) 
    \equiv
    \frac{1}{k}
    \sum_{\ell=1}^k
    F^{\ri}(\bG,\bXi^\ell,\by^\ell)
    :\by^{1:k} \in S_\perp(k)
    \bigg\} \,.
    \eeq
We remark again that the quantities $F,F^{\ri},F^{\rii}$ all depend on the random partition $\hat{\cB}(\bG,\omega)$. For now we have suppressed this dependence from the notation, but this dependence will become relevant in \S\ref{ss:fp-uc}. 

\begin{ppn}\label{p:fp-f-to-i}
Recall $F=F^{\ri}$ from \eqref{e:freeprob.F.equals.F.ri} and $F^{\rii}$ from \eqref{e:F.rii}. 
    For sufficiently large $N$, for all $t\in \R$,
    \[
        \P(F \ge t) \le
        \P\bigg(F^{\rii}(k) +
        \frac{\epsilon'}{5} \ge t
        \bigg) + e^{-cN} \,.
    \]
That is, $F^{\rii}(k)$ approximately stochastically dominates $F$.
\end{ppn}

The \hyperlink{proof:l.fp-f-to-i}{proof of Proposition~\ref{p:fp-f-to-i}} appears below, after the next two lemmas. Define the truncation
\[
    \ytrunc
    = \ytrunc(\bG,\bXi)
    = \hat{\by} \min\bigg\{
    1, \frac{2N^{1/2}}{\|\hat{\by}\|}
    \bigg\}\,.
\]
Note that $\ytrunc$ is also an $L$-Lipschitz function of $(\bG,\bXi)$.

\begin{lem}\label{l:fp-F-kirszbraun}
    Let $C$ be as in \eqref{e:wishart.bound.repeated}.
    Fix any realization of $\bG$ satisfying
    \beq
    \label{e:fp-F-kirszbraun-bG-event}
    \max\bigg\{
       \frac{ \|\bG\|_{\op}}{C},
       \frac{\|\bar{\bx}\|}{2}
      \bigg\} \le N^{1/2}\,.
    \eeq
For such $\bG$, we have the following:
    \begin{enumerate}[(a)]
        \item \label{i:fp-F-kirszbraun-yLip} Under the additional constraint that $\bXi$ satisfies
        \beq
            \label{e:fp-F-kirszbraun-bXi-event}
            \frac{\|\bXi\|_{\op}}{N^{1/2}}
             \le C\,,
        \eeq
        the function
$\by \mapsto F^{\ri}(\bG,\bXi,\by)$,
        restricted to the domain $\{\|\by\| \le 2 N^{1/2}\}$, is $O(N^{-1/2})$-Lipschitz.
        \item \label{i:fp-F-kirszbraun-bXiLip}
Restricted to the (convex) domain \eqref{e:fp-F-kirszbraun-bXi-event},
the function
        \beq
        \label{e:fp-F-kirszbraun-fn-bXi}
            \bXi \mapsto F^{\ri}(\bG,\bXi,\ytrunc(\bG,\bXi))\,,
        \eeq
         is $O(N^{-1/2})$-Lipschitz.
    \end{enumerate}
Consequently, for any realization of $\bG$ satisfying \eqref{e:fp-F-kirszbraun-bG-event}, there exists a function $\tilde{F}(\bXi;\bG)$ which is $O(N^{-1/2})$-Lipschitz in $\bXi$, and agrees with the function
\eqref{e:fp-F-kirszbraun-fn-bXi} on the set of $\bXi$ satisfying \eqref{e:fp-F-kirszbraun-bXi-event}.

\begin{proof}
Fix $\bG$, $\bXi$ satisfying \eqref{e:fp-F-kirszbraun-bG-event}, \eqref{e:fp-F-kirszbraun-bXi-event}. Recalling \eqref{e:F.B}--\eqref{e:F.w}, we calculate
    \begin{align*}
    \nabla_\by F_B(\bG,\bXi,\by)
    &=
    \sum_{j=1}^r
        \frac{2\hB_j}{MN|B_j|}
        \bigg(\sum_{a\in B_j} \bmeta^a, \by \bigg)
        \bigg(\sum_{a\in B_j} \bmeta^a \bigg)\,,\\
            \nabla_\by F_v(\bG,\bXi,\by)
    &=\sum_{j=1}^r
        \frac{2\hv_j}{M}
        \sum_{a\in B_j}
        \frac{(\barbg^a,\by) \barbg^a}{N}
    \\
            \nabla_\by F_U(\bG,\bXi,\by)
    &=
    \sum_{j=1}^r
        \frac{2\hU_j}{M|B_j|}
        \bigg(
            \sum_{a\in B_j}
            \frac{(\bmeta^a,\bar{\bx})(\barbg^a,\by)}{N}
        \bigg)\bigg(
            \sum_{a\in B_j}
            \frac{(\bmeta^a,\bar{\bx})\barbg^a}{N}
        \bigg)\\
            \nabla_\by F_w(\by)
    &=\sum_{j=1}^s\frac{2\hw_j}{N} \Pi_j \by\,.
    \end{align*}
Let $\bG^{B_j} \in \R^{B_j\times N}$ denote the restriction of $\bG$ to the rows indexed by $B_j$, and similarly denote $\bXi^{B_j}$.
Let $\ind_{B_j}\in\R^M$ denote the indicator of the $B_j$ coordinates.
Then we can bound
    \begin{align*}
    \Big| \nabla_\by F_B(\bG,\bXi,\by) \Big|
    &= \bigg| \sum_{j=1}^r
        \frac{2\hB_j}{MN|B_j|}
    (\ind_{B_j})^\st \bXi \by
    \cdot \bXi^\st\ind_{B_j} \bigg|
    \le  \bigg| \sum_{j=1}^r
        \frac{2\hB_j}{MN}
            (\|\bXi\|_\op)^2 \|\by\|  \bigg|
    \lesssim \frac{1}{N^{1/2}}\,,\\
    \Big| \nabla_\by F_v(\bG,\bXi,\by) \Big|
    &=\bigg| \sum_{j=1}^r
        \frac{2\hv_j}{MN}
        (\bG^{B_j})^\st\bG^{B_j} \by \bigg|
            \le
    \bigg|\sum_{j=1}^r
        \frac{2\hv_j}{MN}
    (\|\bG\|_\op)^2 \|\by\|
    \bigg|\lesssim \frac{1}{N^{1/2}}\,,\\
    \Big| \nabla_\by F_U(\bG,\bXi,\by) \Big|
    &=\bigg|\sum_{j=1}^r
        \frac{2\hU_j}{MN^2|B_j|}
        (\bXi^{B_j}\bar{\bx},\bG^{B_j}\by)
        \cdot
            (\bG^{B_j})^\st \bXi^{B_j}\bar{\bx}\bigg| \\
    &\le \bigg|\sum_{j=1}^r
        \frac{2\hU_j}{MN^2|B_j|}
        (\|\bXi\|_\op\|\bG\|_\op)^2
        \|\bar{\bx}\|^2 \|\by\|
        \bigg|
        \lesssim \frac{1}{N^{1/2}}\,,\\
        \bigg|\nabla_\by F_w(\by)\bigg|
    &\le \bigg|
    \sum_{j=1}^s\frac{2\hw_j}{N} \|\by\|\bigg|
           \lesssim \frac{1}{N^{1/2}}\,.
    \end{align*}
Here, all factors of $\hB_j, \hv_j, \hU_j, \hw_j$ can be absorbed in the $\lesssim$ notation, as indicated by Remark~\ref{r:fp-implicit-constants}. Combining these estimates proves part \eqref{i:fp-F-kirszbraun-yLip}.

For the proof of part \eqref{i:fp-F-kirszbraun-bXiLip}, we abbreviate the quantity of interest as $\nabla_{\bXi}\equiv \nabla_{\bXi} F^{\ri}(\bG,\bXi,\ytrunc(\bG,\bXi))$, and decompose it as $\nabla_\bXi =\nabla_\circ+\nabla_\bullet$ where $\nabla_\circ$ denotes the contribution from differentiating through $\ytrunc$, and $\nabla_\bullet$ denotes the remainder. Explicitly,
    \[\nabla_\circ
    =
    (\nabla_\bXi \ytrunc)
        \nabla_{\by} F^{\ri}(\bG,\bXi,\ytrunc(\bG,\bXi))
        \]
Because $\ytrunc$ is $L$-Lipschitz, we have $\|\nabla_{\bXi} \ytrunc\|_{\op} \le L$, and combining with the result of part \eqref{i:fp-F-kirszbraun-yLip} gives $\|\nabla_\circ\|\lesssim1/N^{1/2}$. To bound $\nabla_\bullet$, we calculate that if $a\in B_j$, then
    \begin{align*}
    \nabla_{\bmeta^a} F_B(\bG,\bXi,\ytrunc)
    &= \frac{2\hB_j}{MN|B_j|}
    \bigg(\sum_{b\in B_j} \bmeta^b, \ytrunc \bigg)
    \ytrunc
    =\frac{2\hB_j}{MN|B_j|}
    (\ind_{B_j})^\st\bXi\ytrunc
    \cdot\ytrunc\,,\\
    \nabla_{\bmeta^a} F_U(\bG,\bXi,\ytrunc)
    &=\frac{2\hU_j}{M|B_j|}
       \bigg(
        \sum_{b\in B_j}
        \frac{(\bmeta^b,\bar{\bx})(\barbg^b,\ytrunc)}{N}
    \bigg)
    \frac{(\bg^a,\ytrunc)}{N} \bar{\bx} \\
 &=\frac{2\hU_j}{MN|B_j|}
    (\bXi^{B_j}\bar{\bx},\bG^{B_j}\ytrunc)
    \cdot (\bg^a,\ytrunc) \bar{\bx}\,.
    \end{align*}
We can straightforwardly bound
    \[\Big\|\nabla_{\bmeta^a} F_B(\bG,\bXi,\ytrunc)
        \Big\|
    \lesssim\frac{ \|\bXi\|_\op \|\ytrunc\|^2}{MN |B_j|^{1/2}}
    \lesssim \frac{1}{N}\,,\]
from which it follows that
    \[\Big\|\nabla_{\bXi} F_B(\bG,\bXi,\ytrunc)\Big\|^2
    = \sum_{a=1}^M
    \Big\|\nabla_{\bmeta^a} F_B(\bG,\bXi,\ytrunc)
        \Big\|^2
    \lesssim \frac{1}{N}\,.
    \]
Lastly, we can bound
    \[
    \Big\|\nabla_{\bmeta^a} F_U(\bG,\bXi,\ytrunc)
        \Big\|
    \lesssim \frac{
        \|\bXi\|_\op\|\bG\|_\op \|\bar{\bx}\|^2\|\ytrunc\|
    }{MN^2 |B_j|}
    |(\bg^a,\ytrunc)|
    \lesssim \frac{|(\bg^a,\ytrunc)|}{N^{3/2}}\,,
    \]
from which it follows that
    \[\Big\|\nabla_{\bXi} F_U(\bG,\bXi,\ytrunc)\Big\|^2
    = \sum_{a=1}^M
    \Big\|\nabla_{\bmeta^a} F_U(\bG,\bXi,\ytrunc)
        \Big\|^2
    \lesssim \frac{1}{N^3}
    \sum_{a=1}^M (\bg^a,\ytrunc)^2
    \le \frac{( \|\bG\|_\op\|\ytrunc\|)^2}{N^3}\lesssim \frac{1}{N}\,.
    \]
Combining the above estimates proves $\|\nabla_\bullet\|\lesssim1/N^{1/2}$, and this concludes the proof of part~\eqref{i:fp-F-kirszbraun-bXiLip}.

The final assertion follows using the Kirszbraun extension theorem.
\end{proof}
\end{lem}

\begin{lem}\label{l:F-one-replica-conc}
Let $\tilde{F}$ be the function from the statement of Lemma~\ref{l:fp-F-kirszbraun}. We then have
    \[\Big|F^{\ri}(\bG,\bXi,\hat{\by}(\bG,\bXi)) - \E [\tilde{F}(\bXi;\bG) | \bG]
    \Big| \le
    \frac{\epsilon'}{15}\]
with probability at least $1-e^{-cN}$.

\begin{proof}
    By Lemma~\ref{l:fp-gc} and
    \eqref{e:wishart.bound.repeated}, it holds with probability at least $1-e^{-cN}$ that $\bG$ satisifes \eqref{e:fp-F-kirszbraun-bG-event}, $\bXi$ satisfies \eqref{e:fp-F-kirszbraun-bXi-event}, and $\hat{\by} = \ytrunc$.
    On this event, $F^{\ri}(\bG,\bXi,\hat{\by}(\bG,\bXi)) = \tilde{F}(\bXi;\bG)$. It follows by gaussian concentration of measure (Lemma~\ref{l:lip.subgaus}) that
        \[\Big|\tilde{F}(\bXi;\bG) - \E [\tilde{F}(\bXi;\bG)
    \,|\,\bG] \Big| \le \frac{\epsilon'}{15}\]
with probability at least $1-e^{-cN}$, so the claim follows.
\end{proof}
\end{lem}

\begin{proof}[\hypertarget{proof:l.fp-f-to-i}{Proof of Proposition~\ref{p:fp-f-to-i}}] Let $\bXi^1,\ldots,\bXi^k$ be as above, and $\acute{\by}^\ell = \hat{\by}(\bG,\bXi^\ell)$.
    Let $\delta>0$ be a small constant we will set later (depending on only the parameters in Remark~\ref{r:fp-implicit-constants}) and recall $\hat{\bx} = N^{1/2}\bar{\bx} / \|\bar{\bx}\|$.
    By Lemma~\ref{l:fp-gc}, with probability $1-e^{-cN}$, we have
    \[\bigg\{
        \bigg|\frac{\|\acute{\by}^\ell\|^2}{N}-1\bigg|,
        \frac{|(\acute{\by}^\ell,\hat{\bx})|}{N},
        \frac{|(\acute{\by}^\ell,\acute{\by}^{\ell'})|}{N}
        : \ell,\ell'\in[k], \ell\neq \ell'
 \bigg\}       \le \delta\,.
    \]
Let $\hat{\bx}, \hat{\by}^1, \ldots, \hat{\by}^k$ be the output of the Gram--Schmidt algorithm on input $\hat{\bx}, \acute{\by}^1, \ldots, \acute{\by}^k$, where we normalize the outputs to have norm $N^{1/2}$. A standard estimate shows that for sufficiently small $\delta$, we have an absolute constant $C_6$ such that
    \[
    \max\bigg\{
    \|\hat{\by}^\ell - \acute{\by}^\ell\| : \ell\in[k]
    \bigg\} \le C_6 k\delta N^{1/2}
    \,.\]
    Let $F^\ell = F^{\ri}(\bG,\bXi^\ell,\hat{\by}^\ell)$ and $\acute{F}^\ell = F^{\ri}(\bG,\bXi^\ell,\acute{\by}^\ell)$.
    By
    \eqref{e:wishart.bound.repeated}, it holds with probability at least $1-e^{-cN}$ that $\bG$ satisfies \eqref{e:fp-F-kirszbraun-bG-event}, while each $\bXi^\ell$ satisfies \eqref{e:fp-F-kirszbraun-bXi-event}.
    We also have $\|\hat{\by}^\ell\|, \|\acute{\by}^\ell\| \le 2\sqrt{N}$.
    Therefore, writing $C_7/N^{1/2}$ for the Lipschitz constant in Lemma~\ref{l:fp-F-kirszbraun}\eqref{i:fp-F-kirszbraun-yLip}, we can bound
    \[
        |F^\ell - \acute{F}^\ell|
        \le \frac{C_7}{N^{1/2}}
         \|\hat{\by}^\ell - \acute{\by}^\ell\|
        \le C_6C_7k\delta
        \le \frac{\epsilon'}{15}\,,
    \]
where the last bound follows by setting $\delta$ small enough. Note that we will set $k$ in the next subsection depending on only the parameters in Remark~\ref{r:fp-implicit-constants}, and thus $\delta$ depends on only these parameters as well. 

    Recall $F = F^{\ri}(\bG,\bXi,\hat{\by}(\bG,\bXi))$.
    By Lemma~\ref{l:F-one-replica-conc}, with probability $1-e^{-cN}$, $|F - \E [\tilde{F}(\bXi;\bG) | \bG]| \le \epsilon'/15$ and $|\acute{F}^\ell - \E [\tilde{F}(\bXi;\bG) | \bG]| \le \epsilon'/15$ for all $\ell \in [k]$.
    On the intersection of these events,
    we have
    \[
        F
        \le \E [\tilde{F}(\bXi;\bG) | \bG] + \frac{\epsilon'}{15}
        \le \frac{1}{k} \sum_{\ell=1}^{\ell} \acute{F}^\ell + \frac{2\epsilon'}{15}
        \le \frac{1}{k} \sum_{\ell=1}^{\ell} F^\ell + \frac{\epsilon'}{5}
        \le F^{\rii}(k) + \frac{\epsilon'}{5} \,,
    \]
which proves the claim.
\end{proof}
Next, we bound $F^{\rii}(k)$ by a slightly simpler $F^{\riii}(k)$, which removes the dependence on $\bar{\bx}$ and simplifies the dependence on the $\bXi^\ell$. Unlike for $F, F^{\ri}, F^{\rii}(k)$, we also make explicit the dependence of $F^{\riii}(k)$ on the underlying partition of $[M]$ and $[N]$, which we now allow to be an arbitrary partition $\cB$, rather than the partition $\hat{\cB}(\bG,\omega)$. This is in preparation for the uniform concentration argument in the next subsection, which will allow us to reduce to the case of a deterministic partition. With this in mind, let us denote
    \[
    \scB(\lambda,\mu)
    \equiv \bigg\{
    \cB = (B_1,\ldots,B_r,\BIsing_1,\ldots,\BIsing_s)
    : \frac{|B_j|}{M} = \lambda_j
    \textup{ and }
    \frac{|\BIsing_j|}{N} = \mu_j
    \textup{ for all $j$}
    \bigg\}\,.
    \]
Recall $\hat{\bx} =N^{1/2} \bar{\bx}  / \|\bar{\bx}\|$.
Recalling \eqref{e:F.B}--\eqref{e:F.w}, let us
write for example $F_B(\bG,\bXi,\by,\cB)$ for the expression
\eqref{e:F.B} for the given partition $\cB$. We further rewrite
 \eqref{e:F.B} and \eqref{e:F.U} as
    \begin{align}
    \label{e:F.B.rewrite}
    F_B(\bG,\bXi,\by,\cB)
    &=\sum_{j=1}^r
            \frac{\hB_j}{MN}
            \bigg(
                \frac{1}{|B_j|^{1/2}}
                \sum_{a\in B_j}
                \bmeta^a, \by
            \bigg)^2\,,\\
    F_U(\bG,\bXi,\by,\cB)
    &=
    \frac{\|\bar{\bx}\|^2}{N}
            \sum_{j=1}^r
            \frac{\hU_j}{MN|B_j|}
            \bigg(
                \sum_{a\in B_j}
                \frac{
                (\bXi\hat{\bx})_a
                }{N^{1/2}} \barbg^a,
                \by
            \bigg)^2\,.
    \label{e:F.U.rewrite}
    \end{align}
Then, in comparison to the function $F^{\ri}$ defined by summing \eqref{e:F.B}--\eqref{e:F.w},
we let
$\bar{F}$ be the function defined by summing  $F_v(\bG,\by,\cB)$ and $F_w(\by,\cB)$, together with the simplified functions 
    \begin{align}
    \label{e:bar.F.B}
    \bar{F}_B(\tbmeta,\by,\cB)
    &\equiv
    \sum_{j=1}^r
            \frac{\hB_j}{MN}
            (\tbmeta^j, \by
            )^2\,,\\
           \label{e:bar.F.U}
    \bar{F}_U(\bG,\bz,\by,\cB)
    &\equiv
    \sum_{j=1}^r
            \frac{\hU_j}{MN|B_j|}
            \bigg(
                \sum_{a\in B_j}
                z_a
                \barbg^a,
                \by
            \bigg)^2\,,
    \end{align}
where $\tbmeta\equiv (\tbmeta^1,\ldots,\tbmeta^r)$ is a list of $r$ i.i.d.\ standard gaussian vectors in $\R^N$, and $\bz$ is a standard gaussian vector in $\R^M$. Similarly as before, suppose we have i.i.d.\ replicas
$(\tbmeta^1,\ldots,\tbmeta^k)$ and $(\bz^1,\ldots,\bz^k)$.
Define
\[
    S_+(k)
    = \Big\{
    \by^{1:k}
    \in (\R^N)^k : \,
        \|\by^\ell\|^2 = N, \,
        (\by^\ell,\by^{\ell'}) = 0 \,
        \text{for all $\ell \neq \ell'$}
    \Big\}\,,
\]
and note this corresponds to $S_\perp(k)$ without the condition $(\by^\ell, \bar{\bx}) = \bzero$. We then define
    \beq\label{e:F.iii}
    F^{\riii}(k,\cB)
    \equiv
    \sup\bigg\{
    \bar{F}(\bG,\tbmeta^{1:k},
        \bz^{1:k},\by^{1:k},\cB)
    \equiv \frac{1}{k}
    \sum_{\ell=1}^k
    \bar{F}(\bG,\tbmeta^\ell,
        \bz^\ell,\by^\ell,\cB) 
        : \by^{1:k}
    \in S_+(k)
    \bigg\}\,,
    \eeq
and we conclude this subsection by proving the following:

\begin{ppn}
    \label{p:fp-i-to-ii}
    Recall $F^{\rii}$ from \eqref{e:F.rii}
    and $F^{\riii}$ from \eqref{e:F.iii}. 
    For sufficiently large $N$, for all $t\in \R$,
    \[
        \P(F^{\rii}(k) \ge t)
        \le \P
        \bigg(
        F^{\riii}(k,\hat{\cB}(\bG,\omega))
         + \frac{\epsilon'}{5} \ge t
        \bigg) + e^{-cN} \,.
    \]
That is, $F^{\riii}(k)$ approximately stochastically dominates
$F^{\rii}(k)$.

\begin{proof} Condition on any realization of $\bG$. Recall that
$F^\rii(k)$ is defined by \eqref{e:F.rii}
as the average of $F^\ri$ over $k$ orthogonal replicas, maximized over $S_\perp(k)$. Recall moreover that $F^\ri$ is defined by summing \eqref{e:F.B}--\eqref{e:F.w}, where we rewrote \eqref{e:F.B} and \eqref{e:F.U} in \eqref{e:F.B.rewrite} and \eqref{e:F.U.rewrite} above.
\begin{itemize}
\item Considering \eqref{e:F.U.rewrite}, 
note that $\bXi^\ell \hat{\bx} / \sqrt{N}$ is a standard gaussian in $\R^M$, independent of the gaussian process $\{\bXi^\ell \by^\ell : (\hat{\bx}, \by^\ell) = 0\}$.
Thus, the quantity inside
the maximum defining $F^\rii(k)$
(see \eqref{e:F.rii}),
as a process indexed by $\by^{1:k}$, remains distributionally unchanged if we replace each $\bXi^\ell \hat{\bx} / \sqrt{N}$ with $\bz^\ell$, which corresponds to replacing $F_U$ with $\bar{F}_U$. 
\item Considering \eqref{e:F.B.rewrite}, we can set
     \[
        \tbmeta^j
        = \frac{1}{|B_j|^{1/2}}
        \sum_{a\in B_j}
        \bmeta^{\ell,a}\,,
    \]
and note that the $\tbmeta^j$ are i.i.d.\ standard gaussian vectors in $\R^N$.
\end{itemize}
This shows that
    \begin{align*}
    F^{\rii}(k)&
       \stackrel{d}{=}
       \sup\bigg\{
       \frac{1}{k} \sum_{\ell=1}^k
     \bigg[
       \bar{F}_B(\tbmeta^\ell,\by^\ell,\hat{\cB}(\bG,\omega))
       +F_v(\bG,\by^\ell,\hat{\cB}(\bG,\omega)) \\
 &\qquad\qquad\qquad\qquad 
       +\frac{\|\bar{\bx}\|^2}{N}
       \bar{F}_U(\bG,\bz^\ell,\by^\ell,\hat{\cB}(\bG,\omega))
       +F_w(\by^\ell,\hat{\cB}(\bG,\omega)) \bigg]
        : \by^{1:k}
    \in S_\perp(k)
       \bigg\} \,.
    \end{align*}
Enlarging the domain from $S_\perp(k)$ to $S_+(k)$ can only decrease the above supremum. Finally, the difference between the quantity inside the above supremum, versus the quantity inside the supremum \eqref{e:F.iii} (with $\cB = \hat{\cB}(\bG,\omega)$),
is upper bounded by
    \baln
        &\lt|\frac{\|\bar{\bx}\|^2}{N} - 1 \rt|
        \frac{1}{k}
        \sum_{\ell=1}^k
        \sum_{j=1}^r
        \frac{\hU_j}{MN|B_j|}
        \bigg(
            \sum_{a\in B_j}
            (\bz^\ell)_a \barbg^a,
            \by^\ell
        \bigg)^2 \\
        &\le \bigg|\frac{\|\bar{\bx}\|^2}{N} - 1 \bigg|
        \frac{1}{k}
        \sum_{\ell=1}^k
        \sum_{j=1}^r
        \frac{\hU_j}{MN|B_j|}
        \|\bz^\ell(B_j)\|^2
        \|\bG^{B_j}\|_{\op}^2
        \|\by^\ell\|^2\,.
    \ealn
    By standard estimates, we have $\|{\bz^\ell(B_j)}\|^2 \le 2|B_j|$ with probability at least $1-e^{-cN}$.
    Combined with the estimate $\|\bG^{B_j}\|_{\op} \le C\sqrt{N}$ from \eqref{e:wishart.bound.repeated}, and $\|\by^\ell\| = \sqrt{N}$, we find that the last display is
        \[
    \le O(1)
    \bigg|\frac{\|\bar{\bx}\|^2}{N} - 1 \bigg|
    \le \frac{\epsilon'}{5}\,,
    \]
where the last inequality holds with probability at least $1-e^{-cN}$ by
 Lemma~\ref{l:fp-gc}.
\end{proof}
\end{ppn}

\subsection{Application of uniform concentration}
\label{ss:fp-uc}

In the previous subsection, we showed that the quantity $F$ of \eqref{e:1dproj} is approximately stochastically dominated by the quantity $F^{\rii}(k)$ of \eqref{e:F.rii} (Proposition~\ref{p:fp-f-to-i}), which in turn is approximately stochastically dominated by the quantity $F^{\riii}(k)$ of \eqref{e:F.iii} (Proposition~\ref{p:fp-i-to-ii}). At this point we recall that all of these quantities depend on the random partition $\hat{\cB}(\bG,\omega)$  defined in \eqref{e:fp-partition}, which we can indicate more explicitly by writing 
\begin{align*}
F^{\rii}(k)&\equiv F^{\rii}(k,\hat{\cB}(\bG,\omega))\,,\\
F^{\riii}(k)&\equiv F^{\riii}(k,\hat{\cB}(\bG,\omega))\,,
\end{align*} and so on. We further recall from
\eqref{e:F.iii}
that $F^{\riii}(k,\cB)$ is the maximum over $k$ orthogonal replicas $\by^{1:k}$ of the average value of $\bar{F}$, where $\bar{F}$ is defined
by summing $F_v$ from \eqref{e:F.v}, $F_w$ from \eqref{e:F.w},
$\bar{F}_B$ from \eqref{e:bar.F.B},
and $\bar{F}_U$ from \eqref{e:bar.F.U}.

\textbf{The main result of this subsection is Proposition~\ref{p:fp-ii-to-iii} below, which shows that for fixed partition $\cB \in \scB(\lambda,\mu)$ (independent of the random matrix $\bG$), the quantity}
    \beq\label{e:F.iv}
    F^{\riv}(\cB)
    \equiv
    \sup\bigg\{
    \bar{F}(\bG,\tbmeta,\bz,\by,\cB)
    : \|\by\|^2=N
    \bigg\}
    \eeq
\textbf{is an approximate stochastic upper bound for the maximum value of $F^{\riii}(k,\cB')$ over all partitions $\cB' \in \scB(\lambda,\mu)$, and thus for $F^{\riii}(k,\hat{\cB}(\bG,\omega))$.}
This is achieved by applying a \textbf{uniform concentration} idea, previously introduced in \cite{subag2018free} and used to determine algorithmic thresholds in \cite{huang2023algorithmic}: due to the optimization over $k$ orthogonal replicas, for a fixed partition $\cB'$, the quantity $F^{\riii}(k,\cB')$ concentrates with failure probability $e^{-ckN}$. For large enough $k$, we can take a union bound over all partitions $\cB'$ to obtain the desired conclusion. The main result of this subsection is formally as follows: 

\begin{ppn}
    \label{p:fp-ii-to-iii}
   Recall $F^{\riii}$ from \eqref{e:F.iii}
    and $F^{\riv}$ from \eqref{e:F.iv}.
    Fix any $\cB \in \scB(\lambda,\mu)$ (independent of $\bG$).
    For sufficiently large $N$, we have for all $t\in \R$
    \[
        \P\Big(F^{\riii}(k,\hat{\cB}(\bG,\omega)) \ge t\Big)
        \le \P\bigg(F^{\riv}(\cB) + \frac{\epsilon'}{5} \ge t\bigg)
         + e^{-cN} \,.
    \]
That is, $F^{\riv}(\cB)$ approximately stochastically dominates $F^{\riii}(k,\hat{\cB}(\bG,\omega))$.
\end{ppn}

A further simplification is given by Proposition~\ref{p:fp-iii-to-iv}, which appears at the end of this subsection. 
The \hyperlink{proof:p.fp-ii-to-iii}{proof of Proposition~\ref{p:fp-ii-to-iii}} occupies most of this subsection, and relies on the following intermediate claims. For $C$ as in \eqref{e:wishart.bound.repeated}, define the event
\[
    \cE_{\gd}
    = \bigg\{
        \frac{\|\bG\|_{\op}}{N^{1/2}}
         \le C,
        \frac{\|\bz^\ell\|}{M^{1/2}}
            \le 2\,,\,
        \frac{\|\tbmeta^{j,\ell}\|}
            {N^{1/2}} \le 2
        \,\textup{for all
        $\ell \in [k]$,
        $j\in [r]$}
    \bigg\}\,.
\]
It follows by \eqref{e:wishart.bound.repeated} and standard estimates that
this event occurs with probability at least $1-e^{-cN}$.
The following lemma and corollary formalize the claim that the average over $k$ orthogonal replicas concentrates with failure probability $e^{-ckN}$.

\begin{lem}\label{l:F.iii.sqrtk.Lip}
    For any fixed $\cB \in \scB(\lambda,\mu)$, the function
    \beq
        \label{e:fp-Fii-kirszbraun-fn}
        (\bG,
        \tbmeta^{1:k},
        \bz^{1:k}
        )
        \mapsto F^{\riii}(k,\cB)\,,
    \eeq
restricted to domain $\cE_{\gd}$, is $O((kN)^{-1/2})$-Lipschitz. Consequently, for any $\cB \in \scB(\lambda,\mu)$, there exists an $O((kN)^{-1/2})$-Lipschitz function
    $\tilde{F}^{\riii}(k,\cB)$
    which agrees with the function \eqref{e:fp-Fii-kirszbraun-fn} on the set $\cE_{\gd}$. Moreover, we can choose the functions $\tilde{F}^{\riii}(k,\cB)$ such that $\bbE \tilde{F}^{\riii}(k,\cB)$ is constant over all $\cB\in\scB(\lambda,\mu)$.

\begin{proof}
Recall that $\bar{F}(\bG,\tbmeta^{1:k},\bz^{1:k},\by^{1:k},\cB)$ denotes the function appearing inside the supremum \eqref{e:F.iii}.
Since a supremum of Lipschitz functions remains Lipschitz, it
suffices to show that for any fixed $\by^{1:k} \in S_+(k)$, the function
    \[(\bG,\tbmeta^{1:k},
        \bz^{1:k})
        \mapsto
       \bar{F}
       (\bG,\tbmeta^{1:k},
        \bz^{1:k},\by^{1:k},\cB)\,,
    \]
restricted to domain $\cE_{\gd}$, is $O((kN)^{-1/2})$-Lipschitz. Let us abbreviate $\bar{F}_k$ for the above function, and note that it can be expressed as the sum of functions
$\bar{F}_{B,k}$, $F_{v,k}$, $\bar{F}_{U,k}$, and $F_{w,k}$, where for example $\bar{F}_{B,k}$ denotes the average of the function $\bar{F}_B$ from \eqref{e:bar.F.B} over $k$ orthogonal replicas. If $a\in B_j$, then we calculate
    \[
    \nabla_{\bg^a}
    F_{v,k}
    = \frac{2\hv_j}{kM}
    \sum_{\ell=1}^k
    \frac{(\bg^a,\by^\ell)}{N}\by^\ell\,.
    \]
By the assumption that $\by^{1:k}\in S_+(k)$, we can bound
    \[
    \|\nabla_{\bg^a}F_{v,k}\|^2
    = \bigg(
    \frac{2\hv_j}{kM}\bigg)^2
    \sum_{\ell=1}^k
    \frac{(\bg^a,\by^\ell)^2}{N}
    \lesssim
    \frac{1}{k^2N^3}
    \sum_{\ell=1}^k
    (\bg^a,\by^\ell)^2\,.
    \]
Summing over $a\in[M]$ gives
    \[\|\nabla_\bG F_{v,k}\|^2
    =\sum_{a=1}^M
    \|\nabla_{\bg^a}F_{v,k}\|^2
    \lesssim
    \frac{1}{k^2N^3}
    \sum_{\ell=1}^k
    \|\bG\by^\ell\|^2
    \le \frac{(\|\bG\|_\op)^2 }{k N^2}
    \lesssim \frac{1}{kN}\,.
    \]
For $a\in B_j$, we also calculate
    \[
    \nabla_{\bg^a} \bar{F}_{U,k}
    = \frac{2\hU_j}{kMN|B_j|}
    \sum_{\ell=1}^k
    \bigg(
                \sum_{b\in B_j} (\bz^\ell)_b \barbg^b, \by^\ell
            \bigg)
            (\bz^\ell)_a \by^\ell\,.
    \]
Again by the assumption $\by^{1:k}\in S_+(k)$, we have
    \begin{align*}
    \|\nabla_{\bg^a} \bar{F}_{U,k}\|^2
    &= \bigg(
    \frac{2\hU_j}{kMN|B_j|}\bigg)^2
    \sum_{\ell=1}^k
    [\bz^\ell(B_j)\bG^{B_j}\by^\ell]^2
    ((\bz^\ell)_a)^2 N \\
    &\lesssim
    \frac{1}{k^2 N^5}
    \sum_{\ell=1}^k
    (\|\bG\|_\op)^2 \|\bz^\ell\|^2 \|\by^\ell\|^2
    ((\bz^\ell)_a)^2
    \lesssim
    \frac{1}{k^2 N^2}
    \sum_{\ell=1}^k
    ((\bz^\ell)_a)^2\,.
    \end{align*}
Summing over $a\in[M]$ gives
    \[
    \|\nabla_\bG \bar{F}_{U,k}\|^2
    =\sum_{a=1}^M
    \|\nabla_{\bg^a} \bar{F}_{U,k}\|^2
    \lesssim
    \frac{1}{k^2 N^2}
    \sum_{\ell=1}^k
    \|\bz^\ell\|^2
    \lesssim\frac{1}{kN}\,.
    \]
The remaining two functions $\bar{F}_{B,k}$ and $F_{w,k}$ do not depend on $\bG$, so altogether this proves the desired bound on $\nabla_{\bG}\bar{F}_k$. We next turn to the derivative with respect to $\tbmeta^{1:k}$:
for $a\in B_j$,
    \[
    \nabla_{\tbmeta^{j,\ell}}
    \bar{F}_k
    =
    \nabla_{\tbmeta^{j,\ell}}
    \bar{F}_{B,k}
    =\frac{2\hB_j}{kMN}
    (\tbmeta^{j,\ell},\by^\ell)
    \by^\ell\,,
    \]
from which it follows that
    \[
    \|\nabla_{\tbmeta^{1:k}}
    \bar{F}_k\|^2 \le
    \sum_{\ell=1}^k \sum_{j=1}^r
    \bigg(\frac{2\hB_j}{kMN}\bigg)^2
    \|\tbmeta^{j,\ell}\|^2
    \|\by^\ell\|^4
    \lesssim
    \frac{1}{kN}\,.
    \]
Finally we consider the derivative with respect to $\bz^{1:k}$:
    \[
    \nabla_{\bz^\ell(B_j)}\bar{F}_k
    =\nabla_{\bz^\ell(B_j)}\bar{F}_{U,k}
    = \frac{2\hU_j}{kMN|B_j|}
    \bigg(\sum_{a \in B_j}
        (\bz^\ell)_a \bg^a,
        \by^\ell\bigg)
    \bG^{B_j}\by
    \,,
    \]
from which it follows that
    \[
    \|\nabla_{\bz^\ell(B_j)}\bar{F}_k\|^2
    \le
    \bigg( \frac{2\hU_j}{kMN|B_j|}\bigg)^2
    (\|\bG\|_\op)^4 \|\bz^\ell\|^2
    \|\by^\ell\|^4
    \lesssim
    \frac{1}{k^2 N}\,.
    \]
Summing over $\ell\in[k]$ and $j\in [r]$  gives
    \[
    \|\nabla_{\bz^{1:k}}
        \bar{F}_k\|^2
    = \sum_{\ell=1}^k \sum_{j=1}^r
    \|\nabla_{\bz^\ell(B_j)}\bar{F}_k
        \|^2
    \lesssim \frac{1}{kN}\,.
    \]
Combining the above estimates gives $\|\nabla\bar{F}_k\|^2 \lesssim 1/(kN)$, as claimed.

We now construct the functions $\tilde F^{\riii}(k,\cB)$.
First, choose an arbitrary representative $\cB \in \scB(\lambda,\mu)$.
By the Kirszbraun extension theorem, there exists a $O((kN)^{-1/2})$-Lipschitz function $\tilde F^{\riii}(k,\cB)$ which agrees with \eqref{e:fp-Fii-kirszbraun-fn} on $\cE_{\gd}$. Next, consider any $\acute{\cB} \in \scB(\lambda,\mu)$, and write
\baln
    \cB &= \Big(B_1,\ldots,B_r,\BIsing_1,\ldots,\BIsing_s\Big)\,,\\
    \acute\cB &= \Big(\acute{B}_1,\ldots,\acute{B}_r,\acute{\BIsing}_1,\ldots,\acute{\BIsing}_s\Big)\,.
\ealn
Since $\cB$ and $\acute{\cB}$ both belong to $\scB(\lambda,\mu)$, there must exist permutations $\sigma : [M] \to [M]$ and $\pi : [N] \to [N]$ which map $\cB$ to $\acute{\cB}$, in the sense that $\sigma$ maps $B_j$ to $\acute{B}_j$ for all $j\in[r]$, while $\pi$ maps each $\BIsing_j$ to $\acute{\BIsing}_j$ for all $j\in[s]$.
Let $\bG[\sigma,\pi] \in \bbR^{M\times N}$ be the matrix with entries $\bG[\sigma,\pi]_{i,j} = \bG_{\sigma^{-1}(i),\pi^{-1}(j)}$.
Analogously define $\bz^\ell[\sigma] \in \bbR^M$ and $\tbmeta^\ell[\pi] \in \bbR^N$ for $\ell\in[k]$.
Then define $\tilde F^{\riii}(k,\acute{\cB})$ by
\beq
    \label{e:tilde-F-iii-symmetric-defn}
    \tilde F^{\riii}(k,\acute{\cB})(\bG,\tbmeta^{1:k},\bz^{1:k})
    \equiv \tilde F^{\riii}(k,\cB)(
        \bG[\sigma,\pi],
        \tbmeta^{1:k}[\sigma],
        \bz^{1:k}[\pi]
    )\,.
\eeq
Note that $(\bG,\tbmeta^{1:k},\bz^{1:k}) \in \cE_{\gd}$ if and only if $(\bG[\sigma,\pi],\tbmeta^{1:k}[\sigma],\bz^{1:k}[\pi]) \in \cE_{\gd}$, and on this event the quantity  \eqref{e:tilde-F-iii-symmetric-defn} equals
\[
    F^{\riii}(k,\cB)(
        \bG[\sigma,\pi],
        \tbmeta^{1:k}[\sigma],
        \bz^{1:j}[\pi]
    ) = F^{\riii}(k,\acute{\cB})(\bG,\tbmeta^{1:k},\bz^{1:k})\,,
\]
where the last equality is by the obvious symmetries of the function $F^{\riii}(k,\cB)$ defined by \eqref{e:F.iii}. This verifies that $\tilde F^{\riii}(k,\acute{\cB})$ agrees with $F^{\riii}(k,\acute{\cB})$ on the event  $\cE_{\gd}$.
From \eqref{e:tilde-F-iii-symmetric-defn}, it is also clear that $\tilde F^{\riii}(k,\cB)$ and $\tilde F^{\riii}(k,\cB')$ have the same Lipschitz constant as functions of $(\bG,\tbmeta^{1:k},\bz^{1:k})$, and the same distribution as random variables.
The conclusion follows.
\end{proof}
\end{lem}

\begin{cor}\label{c:fp-uc}
Let $\tilde{F}^{\riii}(k,\cB)$ be the function from Lemma~\ref{l:F.iii.sqrtk.Lip}. As stated therein, $\E\tilde{F}^{\riii}(k,\cB)$ does not depend on the choice of $\cB \in \scB(\lambda,\mu)$, so we denote this quantity simply $\E\tilde{F}^{\riii}(k)$. We then have
        \[\P\bigg(
    \Big|F^{\riii}(k,\cB) - \E\tilde{F}^{\riii}(k)\Big|
    \ge \frac{\epsilon'}{10}
    ;\cE_{\gd}
    \bigg) \le e^{-ckN}\]
for any fixed partition $\cB \in \scB(\lambda,\mu)$.

\begin{proof}
From Lemma~\ref{l:F.iii.sqrtk.Lip},
the function $\tilde{F}^{\riii}(k,\cB)$ is $O((kN)^{-1/2})$-Lipschitz, and therefore concentrates around its mean $\E\tilde{F}^{\riii}(k)$ with failure probability $e^{-ckN}$ (by Lemma~\ref{l:lip.subgaus}). On the event $\cE_\gd$, the functions
$\tilde{F}^{\riii}(k,\cB)$
and $F^{\riii}(k,\cB)$ coincide,
so the claim follows.
\end{proof}
\end{cor}

\begin{proof}[\hypertarget{proof:p.fp-ii-to-iii}{Proof of Proposition~\ref{p:fp-ii-to-iii}}]
Fix $\cB \in \scB(\lambda,\mu)$ (independent of $\bG$).
By Corollary~\ref{c:fp-uc} and a union bound over all $\cB'\in \scB(\lambda,\mu)$, we have
    \baln
       &\P\bigg(F^{\riii}(k,\hat{\cB}(\bG,\omega)) \ge \E\tilde{F}^{\riii}(k) + \frac{\epsilon'}{10}\bigg)\\
        &\qquad\le \P\bigg( \max
        \Big\{
         F^{\riii}(k,\cB')
        :\cB' \in \scB(\lambda,\mu)
        \Big\} \ge \E\tilde{F}^{\riii}(k)+
            \frac{\epsilon'}{10}
    \bigg) \\
        &\qquad\le \P((\cE_{\gd})^c)
        + \sum_{\cB' \in \scB(\lambda,\mu)}
        \P\bigg( F^{\riii}(k,\cB') \ge \E\tilde{F}^{\riii}(k) +
        \frac{\epsilon'}{10};
        \cE_{\gd}\bigg) \\
        &\qquad\le \frac{1}{e^{cN}}
        + \frac{|\scB(\lambda,\mu)|}{e^{ckN}}
        \le e^{-cN} \,.
    \ealn
The final bound follows by taking $k$ large enough (depending on only the parameters in Remark~\ref{r:fp-implicit-constants}), and noting that
$|\scB(\lambda,\mu)| \le e^{O(N)}$. Applying Corollary~\ref{c:fp-uc} again, we have
    \[
    F^{\riii}(k,\hat{\cB}(\bG,\omega))
    \le \E\tilde{F}^{\riii}(k) + \frac{\epsilon'}{10}
    \le F^{\riii}(k,\cB)
        +\frac{\epsilon'}{5}
    \]
with probability at least $1-e^{-cN}$. Finally we note that $F^{\riii}(k,\cB) \le F^\riv(\cB)$ almost surely, so the claim follows.
\end{proof}

 \textbf{For the remainder of this entire section, we will fix a deterministic partition $\cB \in \scB(\lambda,\mu)$ and suppress $\cB$ from our notation. All objective functions below will be defined using this fixed partition $\cB$.} 
We now conclude the current subsection with yet another simplification of the function.
Recall from
\eqref{e:F.iv} that $F^{\riv}$ is the supremum over the sphere of the function $\bar{F}$.
The function $\bar{F}$
is defined
by summing $F_v$ from \eqref{e:F.v}, $F_w$ from \eqref{e:F.w},
$\bar{F}_B$ from \eqref{e:bar.F.B},
and $\bar{F}_U$ from \eqref{e:bar.F.B}.
Abbreviate $r'\equiv2r$, and let $\tbmeta^{r+1},\ldots,\tbmeta^{r'}$ be i.i.d.\ standard gaussian vectors in $\R^N$ (independent of all else).
Abbreviate also
    \beq
    \label{e:hx-j-defn}
    \hx_j
    \equiv\begin{cases}
    \hB_j&\textup{for $1\le j\le r$,}\\
    \hv_{j-r}+\hU_{j-r}
        &\textup{for
        $r+1 \le j \le r'$.}
    \end{cases}
    \eeq
Let $\acute{F}(\bG,\tbmeta,\by)$ denote the sum of $F_w(\by)$ (as in \eqref{e:F.w}), $F_v(\bG,\by)$ (as in \eqref{e:F.v}),
and
    \beq\label{e:F.x}
    F_x(
        \tbmeta,\by)
    \equiv\sum_{j=1}^{r'} \frac{\hx_j}{MN}
        (\tbmeta^j,\by)^2\,.
    \eeq
Lastly, in comparison with
\eqref {e:F.iv}, define
    \beq\label{e:F.rv}
    F^{\rv}
    \equiv\sup\bigg\{
    \acute{F}(\bG,\tbmeta,\by)
    : \|\by\|^2=N\bigg\}\,,\eeq
where now $\tbmeta = (\tbmeta^1,\ldots,\tbmeta^{r'})$. 
We then have the following:

\begin{ppn}
\label{p:fp-iii-to-iv}
Recall $F^{\riv}$ from \eqref{e:F.iv}, and
$F^{\rv}$ from \eqref{e:F.rv}. 
    For sufficiently large $N$, we have for all $t\in \R$
    \[
    \P(F^{\riv} \ge t)
    \le \P\bigg(
    F^{\rv}+ \frac{\epsilon'}{5}
    \ge t
    \bigg) + e^{-cN} \,.
    \]
That is, $F^{\rv}$
approximately stochastically dominates
$F^{\riv}$.

\begin{proof}
For each $j\in[r]$, let
    \[\bv_1(B_j)\equiv
    \frac{\bz(B_j)}{\|\bz(B_j)\|}\,,\]
and complete this to a (Haar random)
orthonormal basis
$(\bv_a(B_j) : 1\le a \le |B_j|)$
of $\R^{|B_j|}$. Let $\tau_j \equiv \|\bz^{B_j}\|^2 / |B_j|$. Then we can rewrite \eqref{e:F.v} and \eqref{e:bar.F.U} as
    \begin{align*}
    F_v(\bG,\by)
    &= \sum_{j=1}^r
        \frac{\hv_j}{MN}
        \sum_{a=1}^{|B_j|}
        (
         \bv_a(B_j),\bG^{B_j} \by)^2
         \,,\\
    \bar{F}_U(\bG,\bz,\by)
    &= \sum_{j=1}^r
        \frac{\tau_j \hU_j}{MN}
        ( \bv_1(B_j),\bG^{B_j} \by)^2\,.
    \end{align*}
If we split off the $a=1$ term of $F_v$ and combine it with $\bar{F}_U$, then we see that $F_v+\bar{F}_U$ is the same as the sum of the two terms
    \begin{align*}
    F_{v\setminus1}(\bG,\by)
    &= \sum_{j=1}^r
        \frac{\hv_j}{MN}
        \sum_{a=2}^{|B_j|}
        (
         \bv_a(B_j),\bG^{B_j} \by)^2
         \,,\\
    \bar{F}_{U+}(\bG,\bz,\by)
    &= \sum_{j=1}^r
        \frac{\hv_j+
            \tau_j \hU_j}{MN}
        ( \bv_1(B_j),\bG^{B_j} \by)^2\,.
    \end{align*}
Fix arbitrary representatives $a_j \in B_j$, and denote
    \begin{align*}
    \acute{F}_{v\setminus}
    (\bG,\by)
    &= \sum_{j=1}^r
        \frac{\hv_j}{MN}
        \sum_{a\in B_j\setminus a_j}
        (\bg^a,\by)^2
    \le F_v(\bG,\by)\,, \\
    \acute{F}_{U+}
        (\tbmeta,\by)
    &= \sum_{j=1}^r
        \frac{\hv_j+
            \tau_j \hU_j}{MN}
    (\tbmeta^{r+j},\by)^2\,.
    \end{align*}
Recalling $F_w$ from \eqref{e:F.w} and
$\bar{F}_B$ from \eqref{e:bar.F.B}, define
    \begin{align*}
    \acute{F}_-(\bG,\tbmeta,\by)
    &\equiv F_w(\by)
    +\acute{F}_{v\setminus}
        (\bG,\by)
    +\bar{F}_B(\tbmeta,\by)
    +\acute{F}_{U+}
        (\tbmeta,\by)
    \\
    &\le F_w(\by)
    +\acute{F}_v(\bG,\by)
        +\bar{F}_B(\tbmeta,\by)
        +\acute{F}_{U+}
        (\tbmeta,\by)
    \equiv\ddot{F}(\bG,\tbmeta,\by)
    \end{align*}
By gaussian rotational invariance,
$\bar{F}$ has the same distribution (as processes indexed by $\by$) as $\acute{F}_-$, which in term is dominated by $\ddot{F}$. Finally, the difference between $\ddot{F}$ and $\acute{F}$
is bounded, for all $\|\by\|^2=N$, by
\[
        \bigg|\sum_{j=1}^r
        \frac{(\tau_j-1) \hU_j}{MN}
        (\tbmeta^{r+j}, \by)^2\bigg|
        \le \sum_{j=1}^r
        \frac{|\tau_j-1| |\hU_j|}{M}
        \|\tbmeta^{r+j}\|^2\,.
    \]
It holds with probability at least $1-e^{-cN}$ that
    \[|\tau_j-1| \le \frac{\alpha \iota \epsilon'}{20r},\quad
    \max_{j\in[r]}
    \frac{\|\tbmeta^{r+j}\|}{N^{1/2}}
     \le 2\,.
    \]
On this event, the above quantity is bounded by 
    \[
        r \cdot \frac{\alpha \iota \epsilon'}{20} \cdot \frac{\|\hU\|_\infty}{M} \cdot 4N 
        \stackrel{\eqref{e:fp-normalization-bd}}{\le} \frac{\epsilon'}{5}\,,
    \]
which concludes the proof.\end{proof}
\end{ppn}

\subsection{Application of gaussian comparison inequalities}
\label{ss:fp-gaussian-comparison}

It follows from Propositions~\ref{p:fp-f-to-i}, \ref{p:fp-i-to-ii}, \ref{p:fp-ii-to-iii}, and \ref{p:fp-iii-to-iv}
that the random variable $F$ of interest is approximately stochastically dominated by the random variable
$F^{\rv}$ from \eqref{e:F.rv}. As heuristically indicated in \S\ref{ss:heuristic.budget}, one way to estimate $F^{\rv}$ is to note that it is the largest eigenvalue of a certain random matrix, which can be computed using free probability theory.  In this section, for the formal proof, we will take a different approach based on gaussian comparison inequalities, which will give the same result. This is summarized by Proposition~\ref{p:fp-vi-to-v} below, which bounds $F^{\rv}$ in terms of a simpler quantity $F^{\rvi}$. This allows us to reduce the main result Proposition~\ref{p:fp-projection} to a simpler statement, Proposition~\ref{p:fp-v-estimate}, which bounds $F^{\rvi}$ in terms of the quantity $F_\diamond$ from \eqref{e:fp-Fdiam}. At the end of this subsection we present the \hyperlink{proof:p.fp-projection}{proof of Proposition~\ref{p:fp-projection}, assuming Proposition~\ref{p:fp-v-estimate}}. The remainder of this appendix is then devoted to the proof of Proposition~\ref{p:fp-v-estimate}. 

Let $J_+$, $J_-$ be the disjoint subsets of $[r]$ defined by
    \begin{align}
    \label{e:Jplus}
    J_+ &\equiv
        \{j\in [r]: \hv_j > 0\}\,, \\
    \label{e:Jminus}
    J_- &\equiv
    \{j\in [r]: \hv_j < 0\}\,.
    \end{align}
Then let $B_+$, $B_-$ be the disjoint subsets of $[M]$ defined by
    \[B_\pm
    \equiv \bigcup_{j \in J_\pm} B_j\,.\]
Further, let $\hbg_+,\hbg_-$ be standard gaussian vectors in $\R^{B_+}, \R^{B_-}$
while $\dbg_+,\dbg_-$ are i.i.d.\ standard gaussian vectors in $\R^N$. Recall the function $F_w$ from \eqref{e:F.w},
and the function $F_v$ from \eqref{e:F.v}, and
the function $F_x$ from \eqref{e:F.x}.
Now define
    \begin{align}
    \label{e:bar.F.v.plus}
    \bar{F}_{v+}(\hbg_+,
        \dbg_+,\by)
    &\equiv\sup\bigg\{
    \frac{2(\hbg_+,\btheta_+)}{M}
        + \frac{2\|\btheta_+\|(\dbg_+,\by)}{MN^{1/2}}
         - \sum_{j\in J_+} \frac{\|\btheta_+(B_j)\|^2}{M\hv_j}
         :\btheta_+\in\R^{B_+}
    \bigg\}\,,\\
    \label{e:bar.F.v.minus}
    \bar{F}_{v-}(\hbg_-,
        \dbg_-,\by)
    &\equiv\inf 
        \bigg\{
    \frac{2(\hbg_-,\btheta_-)}{M}
        + \frac{2\|\btheta_-\|
        (\dbg_-,\by)}{MN^{1/2}}
         - \sum_{j\in J_-} \frac{\|\btheta_-(B_j)\|^2}{M\hv_j}
         :\btheta_-\in\R^{B_-}
    \bigg\}\,.
    \end{align}
(Note that $\bar{F}_{v+}$ is defined as a supremum of a concave function of $\btheta_+$, while $\bar{F}_{v-}$ is defined as an infimum of a convex function of $\btheta_-$.) Define
    \beq\label{e:F.rvi}
    F^{\rvi}
    \equiv
    \sup_{\|\by\|^2=N}
    \bigg\{
    F_w(\by)
    +F_x(\tbmeta,\by)
    +\bar{F}_{v+}(
    \hbg_+,\dbg_+, \by)
    +\bar{F}_{v-}(\hbg_-,
    \dbg_-,\by)
    \bigg\}\,.\eeq
We then have the following:

\begin{ppn}\label{p:fp-vi-to-v}
Recall $F^{\rv}$ from \eqref{e:F.rv}
and $F^{\rvi}$ from \eqref{e:F.rvi}. For all $t\in \R$, we have
    \[
        \P(F^{\rv} \ge t) \le 4\P(F^{\rvi} \ge t)\,.
    \]

\begin{proof}
Recall from \eqref{e:F.rv} that the function $F^{\rv}$ is defined by the supremum of the function $\acute{F}$ over $\|\by\|^2=N$. The function $\acute{F}$ is given by the sum of $F_w$ from \eqref{e:F.w}, $F_v$ from \eqref{e:F.v}, and $F_x$ from \eqref{e:F.x}. Note that we can decompose $F_v = F_{v+}+ F_{v-}$ where         \[
    F_{v+}(\bG,\by)
    \equiv\sum_{j\in J_+}
        \frac{\hv_j}{M}
        \sum_{a\in B_j}
        \frac{(\bg^a,\by)^2}{N}
    = \sum_{j\in J_+}
        \frac{\hv_j}{MN}
        \|\bG^{B_j}\by\|^2\,,
    \]
and $F_{v-}$ is defined analogously. By calculus, we can rewrite
    \begin{align*}
    F_{v+}(\bG,\by)
    &= \sup\bigg\{
    \frac{2(\btheta_+,\bG^{B_+} \by)}{M\sqrt{N}}
            - \sum_{j\in J_+} \frac{\|\btheta_+(B_j)\|^2}{M\hv_j}
            :\btheta_+
            \in\R^{B_+}
    \bigg\}\,,\\
    F_{v-}(\bG,\by)
    &=\inf\bigg\{
    \frac{2(\btheta_-,\bG^{B_-} \by)}{M\sqrt{N}}
            - \sum_{j\in J_-} \frac{\|\btheta_-(B_j)\|^2}{M\hv_j}
    :\btheta_-
            \in\R^{B_-}
    \bigg\}\,.
    \end{align*}
For standard gaussians $Z_+, Z_-$ independent of all else, let
    \begin{align*}
    \tilde{F}_{v+}(\bG,Z_+,\by)
    &\equiv
    \sup\bigg\{
    \frac{2(\btheta_+,\bG^{B_+} \by) + 2\|\btheta_+\|\|\by\| Z_+}{MN^{1/2}}
            - \sum_{j\in J_+} \frac{\|\btheta_+(B_j)\|^2}{M\hv_j}
    :\btheta_+\in\R^{B_+}
    \bigg\}\,,\\
    \tilde{F}_{v-}
    (\bG,Z_-,\by)
    &\equiv
    \inf\bigg\{
    \frac{2(\btheta_-,\bG^{B_-} \by) + 2\|\btheta_-\|\|\by\| Z_-}{MN^{1/2}}
            - \sum_{j\in J_-} \frac{\|\btheta_-(B_j)\|^2}{M\hv_j}
            :\btheta_-\in\R^{B_-}
    \bigg\}\,.
    \end{align*}
Note that we have $Z_+,Z_-\ge0$ with probability $1/4$; on this event, the quantities $F_{v\pm}$ are dominated by the quantities $\tilde{F}_{v\pm}$. Let
    \begin{align*}
    \tilde{F}^{\rv}
    &\equiv
    \sup_{\|\by\|^2=N}
    \bigg\{
    F_w(\by)
    +F_x(\tbmeta,\by)
    +
    \tilde{F}_{v+}(\bG,Z_+,\by)
    +
    \tilde{F}_{v-} (\bG,Z_-,\by)
    \bigg\}\,.
    \end{align*}
Finally, in order to compare $\tilde{F}_{v\pm}$ with $\bar{F}_{v\pm}$ as defined by \eqref{e:bar.F.v.plus} and \eqref{e:bar.F.v.minus}, let $Z$ denote a standard gaussian random variable, and let us compare the gaussian processes
    \begin{align*}
    A(\btheta,\by)
    &\equiv (\btheta,\bG\by)
        + \|\btheta\|\|\by\| Z\,,\\
    B(\btheta,\by)
    &\equiv\|\by\| (\hbg,\btheta) + \|\btheta\| (\dbg,\by)\,.
    \end{align*}
We then calculate
    \begin{align*}
    \E[A(\btheta,\by)
        A(\btheta',\by')]
    &= (\btheta,\btheta')
        (\by,\by')
        +\|\btheta\|\|\btheta'\|
        \|\by\|\|\by'\|\,,\\
    \E[B(\btheta,\by)
        B(\btheta',\by')]
    &=(\btheta,\btheta')
    \|\by\|\|\by'\|
    +\|\btheta\|\|\btheta'\|
    (\by,\by')\,.
    \end{align*}
These quantities are equal whenever $\btheta=\btheta'$ or $\by=\by'$. In general, the $A$-covariances upper bound the $B$-covariances: by scaling, it suffices to show this in the case $\btheta,\btheta',\by,\by'$ are all unit vectors, in which case it reduces to the inequality
    \[
    1-(\btheta,\btheta')
        -(\by,\by')
        +(\btheta,\btheta')
        (\by,\by')
    = \Big(1-(\btheta,\btheta')\Big)
    \Big(1-(\by,\by')\Big) \ge0\,.
    \]
Therefore, we can use Slepian's inequality to compare $\tilde{F}_{v+}$ with $\bar{F}_{v+}$; and we can use Gordon's inequality to compare $\tilde{F}_{v-}$ with $\bar{F}_{v-}$. This implies that $\tilde{F}^{\rv}$ is stochastically dominated by
$\tilde{F}^{\rvi}$. The claim follows by recalling the above observation that
$F^{\rv} \le \tilde{F}^{\rv}$
with probability at least $1/4$.
\end{proof}
\end{ppn}

We have thus reduced our task to proving the following, from which Proposition~\ref{p:fp-projection} follows easily.

\begin{ppn}
    \label{p:fp-v-estimate}
For $F_\diamond$ as in \eqref{e:fp-Fdiam} and $F^{\rvi}$ as in \eqref{e:F.rvi}, we have
    \[\P\bigg(F^{\rvi} 
    \ge F_\diamond +
        \frac{\epsilon'}{5}\bigg) \le e^{-cN}\,.\]
\end{ppn}

The
\hyperlink{proof:p.fp-v-estimate}{proof of Proposition~\ref{p:fp-v-estimate}}
appears at the end of this section.

\begin{rmk}
    \label{r:fp-rpr} Although in the above discussion $r' = 2r$, Proposition~\ref{p:fp-v-estimate} will not require this fact. In the proofs below, we will allow $r'$ and the sequence $(\hx_1,\ldots,\hx_{r'})$ appearing in $F^{\rvi}$ to be arbitrary. We hereafter assume that $r',\hv,\hw,(\hx_1,\ldots,\hx_{r'})$ satisfy
    \beq
        \label{e:reduced-normalization}
        \max\bigg\{ r', \|\hv\|_\infty, \|\hw\|_\infty, \|(\hx_1,\ldots,\hx_{r'})\|_\infty
        \bigg\} \le \frac{2}{\iota}\,.
    \eeq
    Note that for the choice \eqref{e:hx-j-defn} of $(\hx_1,\ldots,\hx_{r'})$
   with $r'=2r$, this follows from \eqref{e:r-s-bd} and \eqref{e:fp-normalization-bd}.
    The $\hx$ in $F_\diamond$ will be $\max_{j\in [r']} \hx_j$; in light of \eqref{e:hx-j-defn}, this is consistent with the choice \eqref{eq:fp-hx} we made when originally defining $F_\diamond$. 
    All implicit constants in the following proofs can depend on the parameters $\alpha,L, \iota,\epsilon'$ and will be uniform in $r,s$ satisfying \eqref{e:r-s-bd}, $\lambda,\mu \in \Lambda$ (as defined by \eqref{e:def-Lambda}), and $r',\hv,\hw,(\hx_1,\ldots,\hx_{r'})$ satisfying \eqref{e:reduced-normalization}. Taking the choice \eqref{e:hx-j-defn}, we obtain the uniformity discussed in Remark~\ref{r:fp-implicit-constants}.
\end{rmk}

\begin{proof}[\hypertarget{proof:p.fp-projection}{Proof of Proposition~\ref{p:fp-projection}, assuming Proposition~\ref{p:fp-v-estimate}}]
    Combining the results of Propositions~\ref{p:fp-f-to-i}, \ref{p:fp-i-to-ii}, \ref{p:fp-ii-to-iii}, \ref{p:fp-iii-to-iv}, \ref{p:fp-vi-to-v}, and \ref{l:fp-Fst-to-Fdiam}, we obtain
    \[
    \P(F \ge F_\star + \epsilon')
    \le 4\P\bigg(F^{\rvi} 
    \ge F_\star + \frac{\epsilon'}5
        \bigg) + e^{-cN}
        = 4\P\bigg(F^{\rvi}
    \ge F_\diamond +
        \frac{\epsilon'}5
    \bigg) + e^{-cN}\,.
    \]
    By Proposition~\ref{p:fp-v-estimate}, this final bound is at most $e^{-cN}$, adjusting $c$ as needed.
\end{proof}

\textbf{The rest of this section is devoted to the proof of Proposition~\ref{p:fp-v-estimate}.} We also make the following simplifying remark which will be useful in what follows:

\begin{rmk}\label{l:freeprob.simplifications}
For the proof of Proposition~\ref{p:fp-v-estimate}, we may assume the following without loss of generality:
\begin{itemize}
    \item Each of $\hv_1,\ldots,\hv_r$ is nonzero.
    \item The entries of $(\hv_1,\ldots,\hv_r)$ and $(\hw_1,\ldots,\hw_s)$ are distinct.
    \item $\hv_1 = \max_{j\in [r]} \hv_j$ and $\hw_1 = \max_{j\in [s]} \hw_j$.
\end{itemize}
Indeed, any $\hv_j = 0$ does not appear in $F^{\rvi}$, while in the optimization problem defining $F_\diamond$, the optimal setting of $\vvv_j$ for such $j$ is clearly $\vvv_j = 1$.
Thus, if $J_{\neq0} = \{j : \hv_j \neq 0\}$ and $\lambda_{\neq0} = \sum_{j \in J_{\neq0}} \lambda_j$, both $F^{\rvi}$ and $F_\diamond$ reduce to a smaller problem (i.e. with smaller $r+s$) with $\lambda_{\neq0} \alpha$ in place of $\alpha$ and (for $j \in J_{\neq0}$) $\lambda_j / \lambda_{\neq0}$ in place of $\lambda_j$.
These problems are self-reducible (i.e. always take the same form as the original) thanks to Remark~\ref{r:fp-rpr}.
Similarly, if $\hv_j = \hv_{j'}$ (resp. $\hw_j = \hw_{j'}$), by convexity of $S_\diamond(\lambda,\mu)$, in the optimization problem defining $F_\diamond$ optimality is attained at $\vvv_j=\vvv_{j'}$ (resp. $\www_j = \www_{j'}$). Thus both $F^{\rvi}$ and $F_\diamond$ reduce to a smaller problem with $\lambda_j+\lambda_{j'}$ in place of $\lambda_j,\lambda_{j'}$ (and  $\mu_j+\mu_{j'}$ in place of $\mu_j,\mu_{j'}$). Note that these reductions all preserve the assumptions $(\lambda,\mu) \in \Lambda$, \eqref{e:r-s-bd}, and \eqref{e:reduced-normalization}.
\end{rmk}

\subsection{Non-variational description of $F_\diamond$}
\label{ss:fp-F-diamond}

The main result of this subsection is Proposition~\ref{p:fp-Fdiam-to-V}, which gives an equivalent non-variational formula for the $F_\diamond$ defined by \eqref{e:fp-Fdiam}. We also provide some technical estimates (Lemmas \ref{l:fp-g2}--\ref{l:fp-uniform-bound}) which will be used in the analysis of the next subsection. Throughout this subsection we assume the simplifications of Remark~\ref{l:freeprob.simplifications} --- in particular, recall for now that $\hw_1$ is maximal among the $\hw_j$. Define $m : (\hw_1,+\infty) \to (0,+\infty)$ by
\beq\label{e:freeprob.m.z}
    m(z) \equiv m_1(z) \equiv \sum_{j=1}^s \frac{\mu_j}{z - \hw_j}\,.
\eeq
By continuity we define $m(\hw_1)\equiv +\infty$.
We will see below that the function $m(z)$ arises from the Lagrangian calculation of $F_\diamond$. For now, we note that $m(z)$ is a continuous, strictly decreasing, strictly convex  function with inverse $m^{-1}(z) : (0,+\infty) \to (\hw_1,+\infty)$. We also define for convenience
\begin{align}\label{e:freeprob.m2.z}
    m_2(z) &\equiv -m'(z) = \sum_{j=1}^s \frac{\mu_j}{(z - \hw_j)^2}\,, \\
    m_3(z) &\equiv \frac{m''(z)}{2}
        = \sum_{j=1}^s \frac{\mu_j}{(z - \hw_j)^3}\,.\label{e:freeprob.m3.z}
\end{align}
Now define the value
\beq\label{e:freeprob.def.z0}
    z_0 =
    \begin{cases}
        \hw_1 & \hv_1 < 0\,, \\
        m^{-1}(\alpha/\hv_1) & \hv_1 > 0\,,
    \end{cases}\eeq
and consider the function
\beq\label{e:freeprob.V}
    V(z) \equiv z + \alpha \sum_{j=1}^r \frac{\lambda_j \hv_j}
        {\alpha - \hv_j m(z)}\,.\eeq
We will see below that the function $V$ also arises in the Lagrangian calculation.
Recalling that $\hv_1$ is maximal among the $\hv_j$, the minimum denominator in the above sum is $\alpha-\hv_1 m(z)$. If $\hv_1<0$, then the denominator is strictly positive as long as $z>\hw_1$. If instead $\hv_1>0$, then the denominator is strictly positive as long as $z>m^{-1}(\alpha/\hv_1)$. This explains the above definition of $z_0$, and shows that we can regard $V$ as a function $(z_0,+\infty) \to \R$. We can extend $V$ to $z_0$ as follows:  if $\hv_1 < 0$, then $m(z)\uparrow+\infty$ as $z\downarrow z_0$, and $V$ extends continuously to $V(z_0) \equiv z_0$. If instead $\hv_1>0$, then $m(z)\uparrow \alpha/\hv_1$ as $z\downarrow z_0$, and by continuity we define $V(z_0)\equiv +\infty$.
The next two lemmas establish some key properties of $V$.

\begin{lem}\label{l:freeprob.V.convex}
The function $V$ defined by \eqref{e:freeprob.V} is strictly convex on the interval $(z_0,+\infty)$, with $V(z)\to\infty$ as $z\to\infty$.
As a consequence, $V$ has a unique minimizer on the interval $[z_0,+\infty)$, which we denote $z_{\min}$.

\begin{proof} For $z>z_0$, we calculate
    \beq
        \label{e:fp-V-pr}
        V'(z)
        = 1 - \alpha \frac{m_2(z)}{m(z)^2} \sum_{j=1}^r \lambda_j
            \bigg(\frac{\hv_j m(z)}{\alpha - \hv_j m(z)}\bigg)^2\,.
    \eeq
Recall that $m:(\hw_1,\infty)\to(0,\infty)$ is strictly decreasing, and we have
$\alpha - \hv_j m(z) > 0$ for all $z>z_0$. It follows that, on  the interval
$(z_0,\infty)$, the function
    \[
    \frac{m(z)}{\alpha-\hv_j m(z)}
    = \frac{1}{\alpha/m(z)-\hv_j}
    \]
is strictly decreasing, \emph{regardless of the sign of $\hv_j$}. Consequently, in \eqref{e:fp-V-pr}, each term
    \[      \bigg(\frac{\hv_j m(z)}{\alpha - \hv_j m(z)}\bigg)^2\]
is strictly decreasing on $(z_0,\infty)$, again \emph{regardless of the sign of $\hv_j$}. In the special case $s=1$, we have $m_2(z)=m(z)^2$, so we conclude in this case that $V'$ is strictly increasing on the interval $(z_0,\infty)$. For $s\ge2$, we note that $m_2(z)/m(z)^2$ is positive, and
    \[
        \frac{d}{dz}\bigg[ \frac{m_2(z)}{m(z)^2}\bigg]
        = \frac{2}{m(z)^3} \Big(m_2(z)^2 - m(z)m_3(z)\Big)
        < 0
    \]
by the Cauchy--Schwarz inequality. It therefore follows in this case as well that $V'$ is strictly increasing, and so $V$ is strictly convex, on the interval $(z_0,\infty)$.
Finally, as $z\to\infty$, $m(z) \to 0$. Therefore the first term of $V(z)$ tends to $\infty$ while the second term remains bounded, and so $V(z) \to \infty$.
\end{proof}
\end{lem}

\begin{lem}
    \label{l:fp-zmin-equal-z0}
Recall from Lemma~\ref{l:freeprob.V.convex} that $z_{\min}$ is the unique minimizer of $V$ on the interval $[z_0,\infty)$.
    We have $z_{\min} = z_0$ if and only if
    \beq
        \label{e:zmin-equal-z0}
        \alpha \le \mu_1 \textup{ and } \hv_1 < 0\,.
    \eeq
    If \eqref{e:zmin-equal-z0} does not hold, $z_{\min}$ is the unique solution to the equation $V'(z) = 0$ for $z\in(z_0,\infty)$.

\begin{proof}
Recall that if $\hv_1 > 0$, then $V(z)\uparrow\infty$ as $z\downarrow z_0$, so in this case clearly $z_{\min} > z_0$.
If instead $\hv_1 < 0$, then $z_0 = \hw_1$, so we have $z_{\min}=z_0$ if and only if the derivative $V'$ stays nonnegative as $z\downarrow \hw_1$. To this end, recall that as $z\downarrow \hw_1$ we have $m(z)\uparrow\infty$, so
    \[\sum_{j=1}^r \lambda_j \bigg(
        \frac{\hv_j m(z)}{\alpha - \hv_j m(z)}\bigg)^2
        \longrightarrow \sum_{j=1}^r \lambda_j = 1\,.\]
We also have $m_2(z)/m(z)^2 \to 1/\mu_1$, and substituting into
\eqref{e:fp-V-pr} gives
    \[
        V'(z)
        \longrightarrow 1 - \frac{\alpha}{\mu_1}
    \]
   as $z\downarrow\hw_1$. This is nonnegative if and only if $\alpha \le \mu_1$, which yields the condition \eqref{e:zmin-equal-z0}. Otherwise, we must have $z_{\min} > z_0$; in this case, since $V$ is strictly convex and twice differentiable, the value $z_{\min}$ must be the unique solution to the equation $V'(z) = 0$ for $z\in(z_0,\infty)$.
\end{proof}
\end{lem}

Let us now define
\beq\label{e:fp.z.star}
    z_\star
    = \begin{cases}
        m^{-1}(\alpha/\hx) & \textup{if $\hx > \alpha / m(z_{\min})$,} \\
        z_{\min} & \textup{otherwise.}
    \end{cases}\eeq
We take the convention that $\alpha / m(z_{\min}) = 0$ if $z_{\min} = \hw_1$.
\textbf{The main result of this subsection is the following:}

\begin{ppn}
    \label{p:fp-Fdiam-to-V}
    We have $F_\diamond = V(z_\star)$, where we recall that $F_\diamond$ is defined by \eqref{e:fp-Fdiam}, and $V(z_\star)$ is defined by \eqref{e:freeprob.V} and \eqref{e:fp.z.star}.
\end{ppn}

\begin{rmk}\label{r:fp.cases} The behavior of $z_\star$ can be delineated by the following three cases:
\begin{enumerate}
    \item \label{i:fp-case1} \eqref{e:zmin-equal-z0} holds and $\hx \le 0$.
    Note that \eqref{e:zmin-equal-z0} implies $z_{\min} = z_0 = \hw_1$, and as a result $\alpha / m(z_{\min}) = 0 \ge \hx$.
    In this case we have
    $z_\star=z_{\min}$, and $V(z_\star)=z_\star$.
    \item \label{i:fp-case2} \eqref{e:zmin-equal-z0} does not hold and $\hx \le \alpha / m(z_{\min})$.
    Then $z_\star = z_{\min} > z_0 \ge \hw_1$ is the unique solution to $V'(z_\star) = 0$ on $(z_0,+\infty)$.
    \item \label{i:fp-case3} Otherwise, we must have $\hx > \alpha / m(z_{\min})$, and so $z_\star = m^{-1}(\alpha/\hx) > z_{\min} \ge z_0 \ge \hw_1$.
\end{enumerate}
In most of what follows, we separate the analysis according to these cases.
\end{rmk}

\begin{proof}[Proof of Proposition~\ref{p:fp-Fdiam-to-V} under case \eqref{i:fp-case1} of Remark~\ref{r:fp.cases}]
Recall from Remark~\ref{l:freeprob.simplifications} that we assume $\hv_1$ is maximal among the $\hv_j$. Therefore, under case \eqref{i:fp-case1}, all the $\hv_j$ are negative. Recall also that $\hw_1$ is maximal among the $\hw_j$. It follows that
    \[
        (\hvv,\vvvv)_{\lambda} + \hx\xxxx + (\hvw,\vwww)_{\mu}
        \le (\hvw,\vwww)_{\mu}
            \le \hw_1 (\mu,\vwww) = \hw_1 = V(z_\star)\,,
    \]
where the last step uses the constraint $(\mu, \vwww) = 1$.
    Equality in the above is attained at $\xxxx = 0$, $\vvvv \equiv 0$, and $\vwww = (1/\mu_1,0,\ldots,0)$, which lies in $S_\diamond(\lambda,\mu)$ because
    \[
        \sum_{j=1}^r \lambda_j \Big((\vvv_j)^{1/2} - 1\Big)^2 + \xxxx
        = 1
        \le \frac{\mu_1}{\alpha}
        = \frac1\alpha \bigg(
            \sum_{j=1}^s
            \mu_j (\www_j)^{1/2}
        \bigg)^2\,.
    \]
This proves in case \eqref{i:fp-case1} that $F_\diamond = V(z_\star)$, as claimed.
\end{proof}

\begin{lem}
    \label{l:fp-case23-prep}
    If case \eqref{i:fp-case1} of Remark~\ref{r:fp.cases} does not hold, then in any local
   maximizer $(\vvvv,\xxxx,\vwww)$ of the function $F_\diamond$ defined by \eqref{e:fp-Fdiam}, the vector $\vwww$ must be entrywise positive.
\begin{proof}
    If case \eqref{i:fp-case1} does not hold, we have $\hv_1 > 0$, $\hx > 0$, or $\mu_1 < \alpha$. We assume for contradiction that $\vwww$ is not entrywise positive, and divide the analysis into a few cases:
\begin{itemize}
\item     \emph{Case 1:} $\hv_1 > 0$ or $\hx > 0$. In this case, let  $\www_{j_1}$ be a positive entry of $\vwww$, which must exist due to the constraint $(\mu,\vwww) = 1$. Now suppose for contradiction that we have $\www_{j_2} = 0$. For small positive $\delta$, the replacement
\[(\www_{j_1},\www_{j_2})
\mapsto \bigg(\www_{j_1} - \frac{\mu_{j_2}}{\mu_{j_1}} \delta, \delta
    \bigg)\]
increases the value of the budget
    \[\frac1\alpha
    \bigg( \sum_{j=1}^s \mu_j (\www_j)^{1/2}\bigg)^2\]
by $\Theta(\delta^{1/2})$. It follows that we may replace $\vvv_1$ with $\vvv_1 + \Theta(\delta^{1/2})$, or $\xxxx$ with $\xxxx + \Theta(\delta^{1/2})$, while remaining inside the set $S_\diamond(\lambda,\mu)$. Recalling
\eqref{e:fp-Fdiam}, we can do this in a way which increases the value of $F_\diamond$ by $\Theta(\delta^{1/2})$, contradicting the assumption that we were at a local maximizer.

\item \emph{Case 2:} $\hv_1 < 0$, $\hx \le 0$, and $\mu_1 < \alpha$.
We divide further into two cases:
\begin{itemize}
\item     \emph{Case 2a:} $\vvvv$ is not identically zero.
    As above, let $\www_{j_1} > 0$ and $\www_{j_2} = 0$.
    Further let $\vvv_{j_3} > 0$, and note that $\hv_{j_3} \le \hv_1 < 0$.
    Therefore the replacement
    \[
        (\www_{j_1},\www_{j_2},\vvv_{j_3}) \mapsto
        \lt(\www_{j_1} - \frac{\mu_{j_2}}{\mu_{j_1}} \delta, \delta,
        \vvv_{j_3} - \Theta(\delta^{1/2})\rt)
    \]
    improves \eqref{e:fp-Fdiam} by $\Theta(\delta^{1/2})$ while remaining inside $S_\diamond(\lambda,\mu)$.
    This gives a contradiction.

\item     \emph{Case 2b:} $\vvvv$ is identically zero. Then the constraints defining $S_\diamond(\lambda,\mu)$ reduce to
 \beq
        \label{e:fp-S-diam}
        1 + \xxxx \le \frac1\alpha \bigg(\sum_{j=1}^s \mu_j (\www_j)^{1/2}\bigg)^2\,, \quad
        (\mu,\vwww) = 1\,.\eeq
    Since $\mu_1 < \alpha$, we have $\vwww \neq (1/\mu_1, 0, \ldots, 0)$.
   Thus there must exist at least one index $j_1\ge2$ for which $\www_{j_1}$ is positive.
    If $\www_1 = 0$, then the replacement
    \[(\www_1,\www_{j_1}) \mapsto
        \bigg(\delta,\www_{j_1} - \frac{\mu_1}{\mu_{j_1}} \delta\bigg)\] improves the value of \eqref{e:fp-Fdiam}
while remaining inside $S_\diamond$.
    Therefore we must have $\www_1 > 0$.
    Finally, if $\www_{j_2} = 0$ for any $j_2 \ge 2$, the replacement
    \[
        (\www_1,\www_{j_1},\www_{j_2}) \rightarrow
        \lt(\www_1 + \delta^{2/3}, \www_{j_1} - \frac{\mu_1}{\mu_{j_1}} \delta^{2/3} - \frac{\mu_{j_2}}{\mu_{j_1}} \delta, \delta \rt)
    \]
    improves the value of \eqref{e:fp-Fdiam} by $\Theta(\delta^{2/3})$ while remaining inside $S_\diamond$.
\end{itemize}
\end{itemize}
This concludes the proof.
\end{proof}
\end{lem}

\begin{proof}[Proof of Proposition~\ref{p:fp-Fdiam-to-V} under cases \eqref{i:fp-case2} and \eqref{i:fp-case3} of  Remark~\ref{r:fp.cases}] We will calculate $F_\diamond$ using the method of Lagrange multipliers.
    Introduce dual variables $\tvv \in [0,+\infty)^r$, $\tx \in [0,+\infty)$, $z\in \R$, $\zeta \in [0,+\infty)$ and consider the Lagrangian
    \baln
        \Lagr(\vvvv,\xxxx,\vwww;\tvv,\tx,z,\zeta)
        &= (\hvv+\tvv,\vvvv)_{\lambda} + (\tx + \hx) \xxxx + (\hvw,\vwww)_{\mu} + z \Big(1 - (\mu, \vwww)\Big) \\
        &\qquad + \zeta \bigg\{
            \alpha^{-1} \bigg(\sum_{j=1}^s \mu_j (\www_j)^{1/2} \bigg)^2
            - \sum_{j=1}^r \lambda_j \Big((\vvv_j)^{1/2} - 1\Big)^2
            - \xxxx
        \bigg\}\,.
    \ealn
    (By Lemma~\ref{l:fp-case23-prep}, any maximizer $\vwww$ is entrywise positive, so we do not need a dual variable for the constraint $\vwww \succeq 0$.) Any maximizer must satisfy $\nabla\Lagr= 0$, as well as the complementary slackness conditions: $\tx \xxxx = 0$, $\tv_j \vvv_j = 0$ for all $j$, and
    \[
        \zeta \bigg\{
            \alpha^{-1} \bigg(\sum_{j=1}^s \mu_j (\www_j)^{1/2}\bigg)^2
            - \sum_{j=1}^r \lambda_j \Big(
                (\vvv_j)^{1/2} - 1\Big)^2
            - \xxxx
        \bigg\} = 0\,.
    \]
Note that if $s=1$, then $\www=\www_1$ and $\mu=\mu_1=1$, and the constraint $\www\in\cW$ implies $\www_1=1$. Therefore let us assume instead that $s\ge2$. 
Stationarity of $\Lagr$ in $\vwww$ yields
    \[
        \mu_j \bigg\{
        \hw_j - z + \frac{\zeta}{\alpha (\www_j)^{1/2}}
        \sum_{k=1}^s \mu_k (\www_k)^{1/2}
        \bigg\} = 0
    \]
    for all $j \in [s]$. There are no division by zero issues because $\vwww$ is entrywise positive. Rearranging gives
    \beq
        \label{e:fp-stationarity-w}
        (z - \hw_j) (\www_j)^{1/2}
        = \frac{\zeta}{\alpha}
        \sum_{k=1}^s \mu_k (\www_k)^{1/2}\,.
    \eeq
Recall that $\zeta\ge0$. Since $\vwww$ is entrywise positive, if $\zeta=0$ then we would have $z=\hw_j$ for all $j$, contradicting the assumption from Remark~\ref{l:freeprob.simplifications} that the $\hw_j$ are distinct. It follows that we must have $\zeta > 0$ and $z > \hw_j$ for all $j$. Moreover, since the right-hand side of \eqref{e:fp-stationarity-w} is independent of $j$, we must have
    \beq
        \label{e:fp-stationarity-w2}
     \frac{  (z - \hw_j) (\www_j)^{1/2} }{ (z - \hw_k) (\www_k)^{1/2}}
     =1\,.
    \eeq
for all $j,k\in[s]$. Substituting \eqref{e:fp-stationarity-w2} into \eqref{e:fp-stationarity-w} yields
    \[
        \frac{\alpha}{\zeta} = \sum_{j=1}^s \frac{\mu_j}{z - \hw_j} = m(z)\,.
    \]
    Similarly, substituting \eqref{e:fp-stationarity-w2} into the constraint $(\mu,\vwww) = 1$ yields
    \[
        \www_j
        = \frac{(z-\hw_j)^{-2}}{\sum_{k=1}^s \mu_k (z-\hw_k)^{-2}}
        = \frac{1}{(z - \hw_j)^2 m_2(z)}\,.
    \]
Note that the above equation is also valid in the case $s=1$, in which case it simply gives $\www=\www_1=1$. We can then derive that
    \[
        \bigg(\sum_{j=1}^s \mu_j (\www_j)^{1/2}\bigg)^2
        = \frac{m(z)^2}{m_2(z)}\,.
    \]
Next, stationarity of $\Lagr$ in $\vvvv$ yields the equation
    \[\lambda_j\bigg\{
        \hv_j + \tv_j
     - \zeta\bigg(1 - \frac1{(\vvv_j)^{1/2}} \bigg)
        \bigg\}= 0
    \]
    for all $j\in [r]$. Thus $\vvv_j > 0$, which implies $\tv_j = 0$ by complementary slackness. Recalling from above that $\alpha/\zeta=m(z)$, we obtain
    \[(\vvv_j)^{1/2} = \frac{1}{1-\hv_j/\zeta}
        = \frac{\alpha}{\alpha - \hv_j m(z)}\,.\]
Note this implies $\alpha>\hv_1 m(z)$, which combined with $z > \hw_1$ implies $z \ge z_0$ (by the definition of $z_0$).
Furthermore, the above expression for $\vvv_j$ implies
    \[
        \sum_{j=1}^r \lambda_j \Big(
            (\vvv_j)^{1/2} - 1\Big)^2
        = m(z)^2 \sum_{j=1}^r
        \lambda_j
        \bigg(\frac{\hv_j}{\alpha - \hv_j m(z)}\bigg)^2\,.
    \]
Since we found above that $\zeta$ is strictly positive, the complementary slackness condition says that the budget constraint is saturated, that is,
    \begin{align} \nonumber
        \xxxx
        &= \frac1\alpha \bigg(
            \sum_{j=1}^s
            \mu_j (\www_j)^{1/2}
        \bigg)^2
        - \sum_{j=1}^r \lambda_j \Big((\vvv_j)^{1/2} - 1\Big)^2\\
    &= \frac{m(z)^2}{\alpha m_2(z)}
    \bigg\{ 1
            - \alpha m_2(z)
            \sum_{j=1}^r
        \lambda_j
        \bigg(\frac{\hv_j}{\alpha - \hv_j m(z)}\bigg)^2
        \bigg\}
        =
        \frac{m(z)^2}{\alpha m_2(z)} V'(z)\,.    \label{e:fp-budget-saturated}
    \end{align}
Since $\xxxx \ge 0$, the above implies $V'(z)\ge0$.
Combined with our above observation that $z\ge z_0$, and again recalling Lemma~\ref{l:freeprob.V.convex}, this implies $z \ge z_{\min}$.
Lastly, stationarity of $\Lagr$ in $\xxxx$ yields the equation
    \beq
        \label{e:fp-stationarity-x}
        \tx + \hx = \zeta = \frac{\alpha }{ m(z)}\,.
    \eeq
We now turn to evaluating $F_\diamond$. To this end, note that the above equation for $\vvv$ implies 
    \begin{align*}
    (\hvv,\vvvv)_{\lambda}
    &= \sum_{j=1}^r \lambda_j\frac{\alpha^2  \hv_j}{(\alpha - \hv_j m(z))^2}
    =\alpha
    \sum_{j=1}^r \lambda_j \hv_j
        \frac{ (\alpha
        - \hv_j m(z) )+  \hv_j m(z)}{(\alpha - \hv_j m(z))^2}\\
    &= \alpha
    \sum_{j=1}^r
        \frac{\lambda_j \hv_j}{ \alpha
        - \hv_j m(z) }
    + \alpha m(z) \sum_{j=1}^r \lambda_j \bigg(  \frac{\hv_j}
        {\alpha - \hv_j m(z)}\bigg)^2\,.
    \end{align*}
Similarly, the above equation for $\www$ implies
    \[
    (\hvw,\vwww)_{\mu}
    =\frac{1}{m_2(z)} \sum_{j=1}^s\mu_j \frac{ \hw_j}{(z - \hw_j)^2}
    = \frac{1}{m_2(z)} \sum_{j=1}^s\mu_j \frac{ z-(z-\hw_j)}{(z - \hw_j)^2}
    = z - \frac{m(z)}{m_2(z)}\,.\]
Combining and rearranging gives
    \begin{align*}
    (\hvv,\vvvv)_{\lambda}+
    (\hvw,\vwww)_{\mu}
    &= \bigg\{ z + \alpha
    \sum_{j=1}^r
        \frac{\lambda_j \hv_j}{ \alpha
        - \hv_j m(z) }\bigg\}
    - \frac{m(z)}{m_2(z)} \bigg\{ 1 -
    \alpha m_2(z)\sum_{j=1}^r \lambda_j \bigg(  \frac{\hv_j}
        {\alpha - \hv_j m(z)}\bigg)^2
    \bigg\}\\
    &= V(z) - \frac{m(z)}{m_2(z)} V'(z)\,.
    \end{align*}
On the other hand, recalling the above equation for $\xxxx$, combined with \eqref{e:fp-stationarity-x} and the complementary slackness condition $\tx \xxxx = 0$, we have 
    \[\hx \xxxx
    = (\tx+\hx) \xxxx
    = \frac{\alpha}{m(z)}\xxxx = \frac{m(z)}{m_2(z)} V'(z)\,.
    \]
As a result we can simplify
    \[
    F_\diamond
    = (\hvv,\vvvv)_{\lambda} + \hx \xxxx + (\hvw,\vwww)_{\mu}
    = V(z)\,,
    \]
and it remains to determine $z$. We now separate into cases: 
\begin{itemize}
\item In case \eqref{i:fp-case2} of Remark~\ref{r:fp.cases}, we have $\hx \le \alpha / m(z_{\min})$.
   If $z > z_{\min}$, then we have $\hx < \alpha / m(z)$, and so \eqref{e:fp-stationarity-x} implies $\tx > 0$, hence $\xxxx = 0$ by complementary slackness.
    Combined with \eqref{e:fp-budget-saturated}, this implies $z = z_{\min}$, a contradiction.     Therefore $z = z_{\min}$, which coincides with $z_\star$ in this case.
\item In case \eqref{i:fp-case3} of Remark~\ref{r:fp.cases}, we have $\hx > \alpha / m(z_{\min})$, so
    \[
        \frac{\alpha}{m(z)}=\zeta
        = \tx + \hx
        \ge \hx
        = \frac{\alpha}{m(z_\star)}
        > \frac{\alpha }{ m(z_{\min})}\,,
    \]
where the last equality is by the definition of $z_\star$. This implies $z>z_{\min}$, so $V'(z)>0$ and $\xxxx>0$, hence $\tx=0$ by complementary slackness. It follows that the first inequality above is an equality, therefore $z=z_\star$.
\end{itemize}
This proves in both cases \eqref{i:fp-case2} and \eqref{i:fp-case3} that
$F_\diamond=V(z_\star)$, as claimed.
\end{proof}

We conclude this subsection by introducing several more scalar functions that will be important in the below arguments.
For $\gamma \in(0,\infty)$, the function
\[
    g_\gamma(\beta)
    \equiv g(\beta,\gamma)= \beta + \gamma m(\beta)
\]
is strictly convex on $(\hw_1,\infty)$ and tends to $\infty$ at both $\hw_1$ and $\infty$. Therefore it has a unique minimizer $\beta_{\min} = \beta_{\min}(\gamma) \in (\hw_1,\infty)$, which solves the equation
\beq
    \label{e:fp-betamin}
    m_2(\beta_{\min}) = -m'(\beta_{\min})
    = \frac1\gamma\,.
\eeq
Note as a result that $\beta_{\min}$ is strictly increasing in $\gamma$. Define
    \beq
    \label{e:fp-betast}
    \beta_\star
    \equiv \beta_\star(\gamma)
    \equiv \begin{cases}
        m^{-1}(\alpha/\hx)
        &\textup{if $\hx > \alpha / m(\beta_{\min})$,} \\
        \beta_{\min}
        & \textup{otherwise.}
    \end{cases}
\eeq
Finally, for $\gamma\in(0,\infty)$ define the function
    \beq\label{e:freeprob.g.of.gamma}
    g(\gamma)
    \equiv
    g_\gamma(\beta_\star(\gamma))
    \,.
    \eeq
We extend the above definitions continuously to $\gamma = 0$ by defining $\beta_{\min}(0) = \hw_1$,
    \[
    \beta_\star(0)
    \equiv\begin{cases}
    m^{-1}(\alpha/\hx)
    &\textup{if $\hx > \alpha/m(\beta_{\min}(0))
    = 0$,}\\
    \beta_{\min}(0)=\hw_1
    &\textup{otherwise,}
    \end{cases}
    \]
and $g(0) = \beta_\star(0)$.

\begin{lem}
\label{l:fp-g2}The functions 
$\beta_\star(\gamma)$ and $g(\gamma)$, as given by \eqref{e:fp-betast} and \eqref{e:freeprob.g.of.gamma},
 have the following properties:
\begin{enumerate}[(a)]
        \item \label{i:fp-g-diff} $g$ is differentiable on $(0,+\infty)$, with $g'(\gamma) = m(\beta_\star(\gamma))$. If $\hx>0$ the statement extends to the right derivative of $g$ at
         $\gamma=0$. If $\hx\le0$, then $g'(\gamma)\uparrow\infty$ as $\gamma\downarrow0$, so $g$ does not have a finite right derivative at $\gamma=0$. 

        \item \label{i:fp-beta-incr} $\beta_\star$ is nondecreasing in $\gamma$.
        \item \label{i:fp-g-growth} As $\gamma \to \infty$ for fixed $\hvv, \hx$, we have $g(\gamma) = \Theta(\gamma^{1/2})$.
    \end{enumerate}

\begin{proof} Note that if $\hx\le0$ then it is never the case that $\hx\ge \alpha/m(\beta_{\min})$, so in this case we always have $\beta=\beta_{\min}$. With this in mind, let us define 
    \[
    \tilde{\beta}
    \equiv
    \begin{cases}
    m^{-1}(\alpha/\hx)
        &\textup{if $\hx>0$,}\\
    \hw_1&\textup{otherwise.}
    \end{cases}
    \]
It is then straightforward to verify that
    \[
        g(\gamma) =
         \min\bigg\{
         \beta + \gamma m(\beta)
         : \beta \ge \tilde \beta
         \bigg\}\,,
    \]
and the minimizing $\beta$ is precisely $\beta_\star(\gamma)$. We then separate into two cases:
\begin{itemize}
\item In the case $\hx\le0$, we always have $\beta_\star(\gamma)=\beta_{\min}(\gamma)$, and the conclusion of part \eqref{i:fp-g-diff} for $\gamma>0$ is an immediate consequence of the standard ``envelope theorem'':
    \begin{align}\nonumber
    g'(\gamma)
    &= \frac{d}{d\gamma}
    g( \beta_{\min}(\gamma), \gamma)
    =
    \bigg\{
    \frac{\partial g}{\partial \beta}
    \frac{d\beta_{\min}}{d\gamma}
    +\frac{\partial g}{\partial \gamma}
    \bigg\}
    ( \beta_{\min}(\gamma), \gamma)\\
    &= \frac{\partial g}{\partial \gamma}
    ( \beta_{\min}(\gamma), \gamma)
    = m(\beta_{\min}(\gamma))\,,
    \label{e:fp.envelope}
    \end{align}
having used that $\partial g/\partial\beta$ must be zero at the point $(\beta_{\min}(\gamma),\gamma)$. As $\gamma\downarrow0$, we note that expanding the equation $m_2(\beta_{\min}(\gamma))=1/\gamma$ implies
    \[\beta_{\min}(\gamma)
    =\hw_1 + (\mu_1\gamma)^{1/2} \Big[1+O(\gamma)\Big]
    \,.\]
It follows that as $\gamma\downarrow0$ we have
    \[g'(\gamma)
    = m(\beta_{\min}(\gamma))
    =\frac{\mu_1}{\beta_{\min}(\gamma)-\hw_1}+O(1)
    =\bigg( \frac{\mu_1}{\gamma}\bigg)^{1/2}+O(1)\,.
    \]
It follows that $g'(\gamma)\to\infty$ as $\gamma\downarrow0$.

\item In the case $\hx>0$, we have $\tilde{\beta}=m^{-1}(\alpha/\hx)$, and \eqref{e:fp-betast} can be rewritten as
    \beq
    \label{e:beta-star-transition}
    \beta_\star(\gamma)
    = \max\Big\{
    \beta_{\min}(\gamma),\tilde{\beta}
    \Big\}
    = \begin{cases}
    \tilde{\beta}
        &\textup{if }
        \gamma <
        \tilde{\gamma}
        \equiv 1/m_2(\tilde{\beta})\\
    \beta_{\min}(\gamma)
        &\textup{otherwise.}
    \end{cases}
    \eeq
At any $\gamma \in(\tilde{\gamma},\infty)$, the function $\beta_\star  = \beta_{\min}$ is differentiable in $\gamma$, and the same calculation as \eqref{e:fp.envelope} shows that $g$ is differentiable at $\gamma$, with derivative $g'(\gamma)=m(\beta_\star(\gamma))$. Similarly, at any $\gamma \in [0,\tilde{\gamma})$, the function $\beta_\star(\gamma) = \tilde{\beta}$ is constant in $\gamma$, so we again conclude that $g$ is differentiable at $\gamma$, with derivative
\[
    g'(\gamma)
    = \frac{\partial g}{\partial \gamma}(\beta_\star(\gamma),\gamma)
    = m(\beta_\star(\gamma))\,.
\]
At $\gamma=\tilde{\gamma}$, note that $\beta_\star(\gamma)$ has well-defined left and right derivatives, equal to $0$ and $(\beta_{\min})'(\tilde{\gamma})$ respectively. The above arguments imply that $g$ has well-defined left and right derivatives, which both equal $m(\beta_\star(\tilde{\gamma}))$, so $g$ is in fact differentiable at $\gamma=\tilde{\gamma}$.
\end{itemize}
The above concludes the proof of \eqref{i:fp-g-diff}. 
The conclusion of \eqref{i:fp-beta-incr} follows from the above representation of $\beta_\star$ in terms of $\beta_{\min}$, and our earlier observation that $\beta_{\min}$ is strictly increasing in $\gamma$.     Finally, note that for $\gamma$ large enough, we always have $\beta_\star(\gamma)=\beta_{\min}(\gamma)$. As $\gamma \to \infty$, it follows from \eqref{e:fp-betamin} that $\beta_{\min}\asymp \gamma^{1/2}
    \asymp 1/m(\beta_{\min})$, and the conclusion of \eqref{i:fp-g-growth} follows.
\end{proof}
\end{lem}

\begin{lem} \label{l:fp-vii-estimate-prep}
Recall the definition of $z_\star$ from \eqref{e:fp.z.star}.
Recalling the expression for $V'(z)$ from \eqref{e:fp-V-pr}, define
    \beq\label{e:fp.gamma.function}
    \gamma(z)
    \equiv
    \alpha \sum_{j=1}^r \lambda_j
    \bigg(
    \frac{ \hv_j}{\alpha - \hv_j m(z)}
    \bigg)^2\,.\eeq
Then, writing $\gamma_\star\equiv\gamma(z_\star)$, we have $\beta_\star(\gamma_\star) = z_\star$.

\begin{proof}
We divide the proof according to the cases from the definition
\eqref{e:fp.z.star} of $z_\star$:
\begin{itemize}
\item   First suppose $z_\star = z_{\min}$.
    Then, recalling \eqref{e:fp-V-pr}, we see that $z_\star$ solves the equation
    \[
    0 = V'(z_\star)
    = 1- m_2(z_\star) \gamma(z_\star)\,,
    \]
which rearranges to
$m_2(z_\star) = 1/\gamma_\star$.
This is precisely the characterization of $\beta_{\min}(\gamma_\star)$, so we conclude $ \beta_{\min}(\gamma_\star)=z_\star=z_{\min}$.
    Finally, since $z_\star = z_{\min}$, we have \[\hx \le
    \frac{\alpha}{m(z_{\min})}
    = \frac{\alpha}{m(\beta_{\min}(\gamma_\star))}
    \,,\] and therefore $\beta_\star(\gamma_\star) = \beta_{\min}(\gamma_\star) = z_\star$.
\item In the other case, we must have $\hx>0$, and
    \[
    z_\star = m^{-1}\bigg(\frac{\alpha}{\hx}\bigg) > z_{\min}\,.
    \]
Abbreviate $\beta_\star\equiv\beta_\star(\gamma_\star)$ and $\beta_{\min}\equiv\beta_{\min}(\gamma_\star)$.
We want to again conclude that $\beta_\star=z_\star$, so suppose for contradiction that this does not hold, in which case it follows from \eqref{e:fp-betast} combined with the above that
    \[
    \beta_\star
    =\beta_{\min}
    > m^{-1}\bigg(\frac{\alpha}{\hx}\bigg)
     = z_\star > z_{\min}\,.
    \]
Now we recall that $z_{\min}$ satisfies $V'(z_{\min})=0$, which by \eqref{e:fp-V-pr} rearranges to
    \beq\label{e:fp.char.of.zmin}
    1 = m_2(z_{\min})
    \gamma(z_{\min})
    =
    \frac{m_2(z_{\min})}{m(z_{\min})^2}
    \Big[
    m(z_{\min})^2\gamma(z_{\min})
    \Big]
    \,.\eeq
Meanwhile, $\beta_{\min}$ is characterized by the equation
    \beq\label{e:fp.char.of.betamin}
    1
    = m_2(\beta_{\min})
    \gamma(z_\star)
    = \frac{m_2(\beta_{\min})}
        {m(z_\star)^2}
    \Big[
    m(z_\star)^2\gamma(z_\star)
    \Big]
    \,.\eeq
We saw in the proof of Lemma~\ref{l:freeprob.V.convex} that
the function $m_2(z)/m(z)^2$ is nonincreasing in $z$. Since $m_2$ is strictly decreasing in $z$, it follows that
\[
    \frac{m_2(\beta_{\min})}
    {m(z_\star)^2}
    <\frac{m_2(z_\star)}
    {m(z_\star)^2}
    \le \frac{m_2(z_{\min})}
    {m(z_{\min})^2}
    \]
We also saw in the proof of Lemma~\ref{l:freeprob.V.convex} that
the function $m(z)^2\gamma(z)$ is strictly decreasing in $z$, which implies
$m(z_\star)^2\gamma(z_\star)<
    m(z_{\min})^2\gamma(z_{\min})$. This shows that the above equations \eqref{e:fp.char.of.zmin} and \eqref{e:fp.char.of.betamin} cannot hold simultaneously, a contradiction.
\end{itemize}
Thus we have shown in all cases that $\beta_\star(\gamma_\star)=z_\star$, as claimed.
\end{proof}
\end{lem}

The following technical estimates will be useful in the sequel for obtaining high-probability estimates on $F^{\rvi}$ that are uniform in the sense of Remark~\ref{r:fp-rpr}.

\begin{lem}\label{l:fp-m-fns-conversion}
Let $j,\ell\in[3]$ with $j\neq\ell$. Recall the functions $m_1,m_2,m_3$ from \eqref{e:freeprob.m.z}, \eqref{e:freeprob.m2.z}, and \eqref{e:freeprob.m3.z}; and let $s_\ell$ denote the inverse of $m_\ell$. For any $z>0$, we have the bounds
    \[
        \iota z^{j/\ell} 
        \le m_j(s_\ell(z))
        \le \bigg(\frac{z}{\iota}\bigg
            )^{j/\ell}\,.
    \]

\begin{proof}
From the definitions \eqref{e:freeprob.m.z}--\eqref{e:freeprob.m3.z}, 
for each $\ell\in[3]$ and for any $z>\hw_1$ we have
    \[
    \frac{\mu_1}{(z-\hw_1)^\ell}
    \le m_\ell(z)=\sum_{j=1}^s
    \frac{\mu_j}{(z-\hw_j)^\ell}
    \le \frac{1}{(z-\hw_1)^\ell}\,.
    \]
Thus, for $t_- = \hw_1 + (\mu_1/z)^{1/\ell}$ and $t_+ = \hw_1 + z^{-1/\ell}$, we have
    \[
        m_\ell(t_+) 
        \le \frac{1}{(t_+ - \hw_1)^\ell}
        = z 
        = \frac{\mu_1}{(t_- - \hw_1)^\ell}
        \le m_\ell(t_-)\,.
    \]
    Since $m_\ell$ is strictly decreasing, this implies $t_- \le s_\ell (z) \le t_+$. Then, since $m_j$ is also strictly decreasing, it follows that
    \[
        \mu_1 z^{j/\ell}
        = \frac{\mu_1}{(t_+ - \hw_1)^j}
        \le m_j(t_+)
        \le m_j(s_\ell (z))
        \le m_j(t_-)
        \le \frac{1}{(t_- - \hw_1)^j}
        = (z/\mu_1)^{j/\ell}\,.
    \]
The conclusion follows by recalling that $\mu_1 \ge \iota$ because $(\lambda,\mu) \in \Lambda$.
\end{proof}
\end{lem}

\begin{lem}
\label{l:fp-scalar-fn-lipschitz}
The following functions are $C_6$-Lipschitz, with $C_6 = 2/\iota$:
    \begin{enumerate}[(a)]
        \item \label{i:fp-m-inv-lipschitz} $z \mapsto m^{-1}(1/z)$, on $(0,+\infty)$.
        \item \label{i:fp-g-sq-lipschitz} $\gamma \mapsto g(\gamma^2)$, on $[0,+\infty)$.
        \item \label{i:fp-betast-sq-lipschitz} $\gamma \mapsto \beta_\star(\gamma^2)$, on $[0,+\infty)$.
\end{enumerate}

\begin{proof}
    A routine calculation and application of Lemma~\ref{l:fp-m-fns-conversion} implies
    \[
        \frac{d}{dz} m^{-1}(1/z) = \frac{1/z^2}{m_2(m^{-1}(1/z))}
        \le \frac{1/z^2}{\iota / z^2} = \frac{1}{\iota}\,.
    \]
    This proves part~\eqref{i:fp-m-inv-lipschitz}.
    Since $m$ is strictly decreasing, and $\beta_\star \ge \beta_{
    \min}$, we also have, for $\gamma > 0$,
    \[
        \frac{d}{d\gamma} g(\gamma^2) 
        = 2\gamma m(\beta_\star(\gamma^2))
        \le 2\gamma m(\beta_{\min}(\gamma^2))
        = 2\gamma m(s_2(1/\gamma^2))
        \le \frac{2}{\iota^{1/2}}\,,
    \]
where the final inequality is by Lemma~\ref{l:fp-m-fns-conversion}.
This shows that $\gamma \mapsto g(\gamma^2)$ is $2/\iota^{1/2}$-Lipschitz on $(0,+\infty)$, and the estimate extends to $[0,+\infty)$ since $g$ was defined to be continuous on $[0,\infty)$ (see the comments below \eqref{e:freeprob.g.of.gamma}). 
    This proves part~\eqref{i:fp-g-sq-lipschitz}. For part~\eqref{i:fp-betast-sq-lipschitz}, it suffices to prove that the function
    \[
        \gamma \mapsto \beta_{\min}(\gamma^2)
        = s_2\bigg(\frac1{\gamma^2}\bigg)
    \]
    is $O(1)$-Lipschitz on $[0,+\infty)$. Similarly as above, a routine calculation and application of Lemma~\ref{l:fp-m-fns-conversion} shows that at any $\gamma > 0$, 
    \[
        \frac{d}{d\gamma} 
        s_2\bigg(\frac1{\gamma^2}\bigg)
        = \frac{1/\gamma^3}
            {m_3(s_2(1/\gamma^2))}
        \le \frac{1/\gamma^3}{\iota / \gamma^3}
        = \frac{1}{\iota}\,.
    \]
This shows that $\gamma \mapsto \beta_{\min}(\gamma^2)$ is $1/\iota$-Lipschitz on $(0,+\infty)$, and the estimate extends to $[0,+\infty)$ since the function $\beta_{\min}$ was defined to be continuous on $[0,\infty)$ (again we refer to the comments below \eqref{e:freeprob.g.of.gamma}). 
This proves part~\eqref{i:fp-betast-sq-lipschitz}.
\end{proof}
\end{lem}

\begin{lem}
    \label{l:fp-uniform-bound}
There exists a constant $C_7=C_7(\alpha,\iota)$ such that for all parameters $\hw,\hv,(\hx_1,\ldots,\hx_{r'})$ satisfying \eqref{e:reduced-normalization}, we have
    \[
        \max\bigg\{
            \gamma_\star, \frac{1}{\alpha - \hv_1 m(z_\star)}
        \bigg\} \le C_7\,,
    \]
for $z_\star$ as defined by \eqref{e:fp.z.star} and $\gamma_\star$ as defined by Lemma~\ref{l:fp-vii-estimate-prep}.

\begin{proof} Recall the sets $J_+,J_- \subseteq [r]$ defined in \eqref{e:Jplus}, \eqref{e:Jminus}, which index the positive and negative entries of $\hv$. Let $\gamma_+(z)$ denote the contribution to \eqref{e:fp.gamma.function} from indices $j\in J_+$, and denote similarly $\gamma_-(z)$. Then, for any $z>\hw_1$, we have 
    \beq\label{e:fp-uniform-bound-Jminus-contribution}
    \gamma_-(z)
    = \alpha \sum_{j\in J_-}
            \lambda_j \bigg(
            \frac{\hv_j}{\alpha - \hv_j m(z)}
        \bigg)^2
        \le \alpha 
        \sum_{j\in J_-}
            \lambda_j
    \bigg(
            \frac{\hv_j}{\alpha}
        \bigg)^2
        \stackrel{\eqref{e:reduced-normalization}}{\le} \frac{4}{\alpha\iota^2}\,.
    \eeq
If $J_+ = \emptyset$, this upper bounds $\gamma(z)=\gamma_-(z)$ for any $z>\hw_1$. 
    In this case, we also have \[
    \frac1{\alpha - \hv_1 m(z)} \le \frac1\alpha\,,\]
for any $z>\hw_1$, so the result follows. 

Otherwise, assume $J_+ \neq \emptyset$, so $\hv_1 > 0$. By Lemma~\ref{l:fp-zmin-equal-z0}, this implies $z_{\min} > z_0$. Recall that \eqref{e:fp-uniform-bound-Jminus-contribution} above already bounds $\gamma_-(z)$ for all $z>\hw_1$. Meanwhile, both $m(z)$ and $\gamma_+(z)$ are strictly decreasing on the interval $z \in (z_0,+\infty)$. Since $z_\star \ge z_{\min} > z_0$ and $\gamma(z_\star) = \gamma_\star$, it suffices to show 
\beq
        \label{e:fp-uniform-bound-goal}
        \max\bigg\{
        \gamma_+(z_{\min})
        \le \gamma(z_{\min}),
            \frac{1}{\alpha - \hv_1 m(z_{\min})}
        \bigg\} \le C_7\,.
    \eeq
We divide the analysis into two cases:
    \begin{itemize}
        \item If $m(z_{\min}) \le \alpha / (2\hv_1)$, then 
        \begin{align*}
    \gamma_+(z_{\min})\
    &=
            \alpha \sum_{j\in J_+} \lambda_j \bigg(
            \frac{\hv_j}
            {\alpha - \hv_j m(z_{\min})}
            \bigg)^2\\
            &\le \alpha \bigg(
                \frac{\hv_1}{\alpha - \hv_1 m(z_{\min})}
            \bigg)^2 
            \le \frac{4(\hv_1)^2}{\alpha}
            \stackrel{\eqref{e:reduced-normalization}}{\le}
            \frac{16}{\alpha \iota^2}\,,\\
       \frac{1}{\alpha - \hv_1 m(z_{\min})}
           &\le \frac{2}{\alpha}\,.
        \end{align*}

        \item If $m(z_{\min}) > \alpha / (2\hv_1)$, then Lemma~\ref{l:fp-m-fns-conversion} gives
        \[
            m_2(z_{\min}) >
            m_2\bigg(m^{-1}\bigg(\frac{\alpha}{2\hv_1}\bigg)\bigg)
            \ge \frac{\alpha^2 \iota }{4(\hv_1)^2}\,.
        \]
        Since $\hv_1 > 0$, it follows from Lemma~\ref{l:fp-zmin-equal-z0} that $V'(z_{\min}) = 0$.
        Recalling the formula for $V'$ from \eqref{e:fp-V-pr}, this implies
        \[
        \frac{\alpha\lambda_1(\hv_1)^2}
            {(\alpha-\hv_1m(z_{\min}))^2}
        \le
        \gamma(z_{\min}) 
            = \frac{1}{m_2(z_{\min})}
            \le \frac{4(\hv_1)^2}{\alpha^2\iota}
            \stackrel{\eqref{e:reduced-normalization}}{\le}
            \frac{16}{\alpha^2\iota^3}\,.
        \]
Rearranging the above and recalling $(\lambda,\mu) \in \Lambda$ gives
        \[
            \frac{1}{\alpha - \hv_1 m(z_{\min})}
            \le \frac{2}{\alpha^{3/2}\iota}\,.
        \]
    \end{itemize}
In all cases, the conclusion follows.
\end{proof}
\end{lem}

\subsection{Conclusion of proof}
\label{ss:fp-eval}

In this subsection we conclude the \hyperlink{proof:p.fp-v-estimate}{proof of Proposition~\ref{p:fp-v-estimate}}, which says that with very high probability, the quantity $F^{\rvi}$ defined by \eqref{e:F.rvi} is approximately upper bounded by the value $F_\diamond$ defined by \eqref{e:fp-Fdiam}. We will use the characterization $F_\diamond=V(z_\star)$ from Proposition~\ref{p:fp-Fdiam-to-V}, where $V(z_\star)$ is defined by \eqref{e:freeprob.V} and \eqref{e:fp.z.star}. In what follows, it will be convenient to abbreviate
    \begin{align}\label{e:fp.dbD}
    \dbD &\equiv \sum_{j=1}^s \hw_j \Pi_j \in \R^{N\times N}\,,\\
    \dbX &\equiv
        \sum_{j=1}^{r'} \frac{\hx_j}{M}
            \tbmeta^j(\tbmeta^j)^\st
            \in \R^{N\times N}\,.\label{e:fp.dbX}
    \end{align}
Similarly, for each $j\in[r]$, we let $\hat{\Pi}_j$ denote orthogonal projection onto coordinates $a\in B_j$.  Recall the sets $J_+,J_- \subseteq [r]$ defined in \eqref{e:Jplus} and  \eqref{e:Jminus}, and let 
    \begin{align*}
    \hbD_+ &\equiv \sum_{j\in J_+} \hv_j \hat{\Pi}_j \in \R^{B_+ \times B_+}\,,\\
    \hbD_- &\equiv \sum_{j\in J_-} \hv_j \hat{\Pi}_j \in \R^{B_- \times B_-}\,.
    \end{align*}
Note that $\dbD$, $\hbD_+$, and $\hbD_-$ are all diagonal matrices. Then, for the reader's convenience, we recall here that
    \beq\label{e:fp-F-v}
    F^{\rvi}
    \stackrel{\eqref{e:F.rvi}}{=}
    \sup_{\|\by\|^2=N}
    \bigg\{
    F_w(\by)
    +F_x(\tbmeta,\by)
    +\bar{F}_{v+}(
    \hbg_+,\dbg_+,\by)
    +\bar{F}_{v-}(\hbg_-,\dbg_-,\by)
    \bigg\}
    \eeq
where the terms are defined by
    \begin{align}
    \label{e:fp.recall.Fw}
    F_w(\by)
    &\stackrel{\eqref{e:F.w}}{=}
    (\hvw,\vwww)_{\mu}
    = \sum_{j=1}^s
    \frac{\hw_j}{N}
    (\Pi_j, \by^{\otimes 2})
    = \frac{(\dbD\by,\by)}{N}
    \,,\\
    \label{e:fp.recall.Fx}
    F_x(\tbmeta,\by)
    &\stackrel{\eqref{e:F.x}}{=}
    \sum_{j=1}^{r'} \frac{\hx_j}{MN}
        (\tbmeta^j,\by)^2
    = \frac{(\dbX\by,\by)}{N}\,,\\
    \label{e:fp.recall.barFvplus}
    \bar{F}_{v+}(
    \hbg_+,\dbg_+,\by)
    &\stackrel{\eqref{e:bar.F.v.plus}}{=}\sup\bigg\{
    \frac{2(\hbg_+,\btheta_+)}{M}
        + \frac{2\|\btheta_+\|(\dbg_+,\by)}{MN^{1/2}}
         -\frac{(\btheta_+,
         (\hbD_+)^{-1}\btheta_+)}
            {M}
    :\btheta_+\in\R^{B_+}
    \bigg\}\,,\\
        \label{e:fp.recall.barFvminus}
    \bar{F}_{v-}(
    \hbg_-,\dbg_-,\by)
    &\stackrel{\eqref{e:bar.F.v.minus}}{=}\inf
        \bigg\{
    \frac{2(\hbg_-,\btheta_-)}{M}
        + \frac{2\|\btheta_-\|
        (\dbg_-,\by)}{MN^{1/2}}
         - \frac{(\btheta_-,
         (\hbD_-)^{-1}\btheta_-)}
            {M}
        :\btheta_-\in\R^{B_-}
    \bigg\}\,.
    \end{align}

\begin{proof}[\hypertarget{proof:p.fp-v-estimate.case.one}{Proof of Proposition~\ref{p:fp-v-estimate} under case \eqref{i:fp-case1}
of Remark~\ref{r:fp.cases}}]
Since $\hv_1 < 0$, the set $B_+$ is empty, so we trivially have $\bar{F}_{v+}=0$. Setting $\btheta_- = \bzero$ gives the upper bound $\bar{F}_{v-}\le0$. Since $\hx = \max_{j\in [r']} \hx_j \le 0$, we also conclude $F_x\le0$. Therefore
    \[
        F^{\rvi} 
        \le \sup_{\|\by\|^2=N}
        \bigg\{
            \sum_{j=1}^s \frac{\hw_j}{N} (\Pi_j, \by^{\otimes 2})
        \bigg\}
        \le \hw_1\,,
    \]
recalling that $\hw_1$ is maximal among the $\hw_j$.
\end{proof}

Thus we can assume from now on that
 we are in case \eqref{i:fp-case2} or \eqref{i:fp-case3} of Remark~\ref{r:fp.cases}.
Then $z_\star > z_0$, so $m(z_\star) < \infty$ and $\alpha - \hv_1 m(z_\star) > 0$. \textbf{We begin with a heuristic derivation for these cases.} Let
\[
    \gamma \equiv \frac{\|\btheta_+\|^2 + \|\btheta_-\|^2}{M\alpha}\,,
\]so that we can express
\[
    \frac{\|\btheta_+\| \dbg_+ + \|\btheta_-\| \dbg_-}{M\sqrt{N}}
    = \frac{\gamma^{1/2} \dbg}{N}\,,
\]
where $\dbg$ is a standard gaussian vector in $\R^N$.
\textbf{We now non-rigorously exchange the order of optimizations in \eqref{e:fp-F-v}, first maximizing over $\by$ for fixed $(\btheta_+,\btheta_-)$ in \eqref{e:fp.recall.barFvplus} and \eqref{e:fp.recall.barFvminus}.} We then proceed by the method of Lagrange multipliers: omitting terms that do not depend on $\by$, the Lagrangian for the new optimization problem is
    \beq\label{e:fp-y-lagrangian}
    \Lagr(\by;\beta)
    = \frac{(\dbD \by, \by)}{N}
    + \frac{(\dbX \by, \by)}{N} 
    + \frac{2\gamma^{1/2}(\dbg,\by)}{N}
    + \beta \lt(1 - \frac{\|\by\|^2}{N}\rt)\,.
\eeq
Define $\dbA_\beta \equiv (\beta \bI_N - \dbD)^{-1}$. Stationarity in $\by$ implies
\beq
    \label{e:fp-g-stationarity}
    \by = \dbA_\beta \Big(
        \gamma^{1/2} \dbg
        + \dbX \by
    \Big)\,.\eeq
Substituting this back into the Lagrangian gives
\beq
    \label{e:fp-y-lagrangian-bd}
    \Lagr(\by;\beta)
    = \beta
    - \frac{((\beta \bI_N - \dbD) \by, \by)}{N}
    + \frac{(\dbX\by,\by)}{N}
    + \frac{2\gamma^{1/2}(\dbg,\by)}{N}
    = \beta + \frac{\gamma^{1/2}(\dbg,\by)}{N}\,.
\eeq
In this heuristic derivation we will use $\approx$ to mean equal up to additive $o_N(1)$ error with high probability.
Note that $(\dbA_\beta \dbg, \dbg) / N \approx m(\beta)$ and $(\dbA_\beta \dbg, \tbmeta^j) \approx 0$, so we conclude
\beq\label{e:fp.heuristic.lagrangian}
    \Lagr(\by;\beta)
    \approx \beta + \gamma m(\beta)
    = g_\gamma(\beta)\,.\eeq
We next determine the possible values of $\beta$; and we claim that the optimum will be achieved by taking $\beta\approx\beta_\star(\gamma)$. To this end, note that the constraint on the norm of $\by$ implies
\beq
    \label{e:fp-y-norm}
    1 = \frac{\|\by\|^2}{N}
    \stackrel{\eqref{e:fp-g-stationarity}}{=}
    \frac1N \bigg\|\dbA_\beta \Big(
        \gamma^{1/2} \dbg
        +  \dbX \by \Big)
        \bigg\|^2 
    \approx m_2(\beta) \bigg\{\gamma + \sum_{j=1}^{r'}
    \bigg(
    \frac{\hx_j (\tbmeta^j, \by)}{N\alpha} \bigg)^2
    \bigg\}
    \ge \gamma m_2(\beta)\,.
\eeq
In the rest of this heuristic discussion we will treat approximate equalities as equalities. Since $\beta_{\min}(\gamma)$ is characterized by the relation
\eqref{e:fp-betamin}, and $m_2$ is decreasing in $\beta$, this implies
 $\beta \ge \beta_{\min}(\gamma)$. Recall from the discussion around \eqref{e:fp-betamin} that for any positive $\gamma$, $g_\gamma$ is a strictly convex function of $\beta\in(\hw_1,\infty)$, with global minimum at $\beta_{\min}$. Thus it is strictly increasing on the range $\beta\ge\beta_{\min}$. Consequently, for the optimization problem \eqref{e:fp.heuristic.lagrangian}, we seek the largest possible $\beta$ for which constraint \eqref{e:fp-y-norm} holds. We shall argue that in all cases, the optimal choice is given by $\beta=\beta_\star(\gamma)$, as defined by \eqref{e:fp-betast}. To this end, note that for each $j\in [r']$, we have
\[
    \frac{(\tbmeta^j,\by)}{N}
    \stackrel{\eqref{e:fp-g-stationarity}}{=}
    \frac1N \bigg(\tbmeta^j, \dbA_\beta \Big(
        \gamma^{1/2} \dbg
        +  \dbX \by \Big) \bigg) 
    \approx
    \frac{\hx_j m(\beta)}{\alpha}
    \cdot \frac{(\tbmeta^j, \by)}{N}\,,
\]
which implies either $(\tbmeta^j,\by) / N \approx 0$ or $m(\beta) \hx_j / \alpha \approx 1$. We now divide into cases:
\begin{itemize}
\item If $\hx < \alpha/m(\beta_{\min})$, then for all $j\in[r']$ we must have
    \[
    \hx_j \le \hx < \frac{\alpha}{m(\beta_{\min})}
     \le \frac{\alpha}{m(\beta)}\,,
    \]
where the last inequality follows from our above observation that $\beta\ge \beta_{\min}$, together with the fact that $m$ is a decreasing function. Therefore we cannot have $\hx_j m(\beta) / \alpha \approx 1$ for any $j$, and we must instead have $(\tbmeta^j,\by) / N \approx 0$ for all $j\in[r']$. This implies that \eqref{e:fp-y-norm} holds with (approximate) equality, and therefore $\beta=\beta_{\min}$, which in this case also agrees with the definition of $\beta_\star$.

\item If $\hx\ge \alpha/m(\beta_{\min})$, then our options are to let $\beta=\beta_{\min}$ as above, or to let $m(\beta)\hx_j/\alpha=1$ for some $j\in[r']$. The largest such $\beta$ is $\beta=m^{-1}(\alpha/\hx)$, which in this case also agrees with the definition of $\beta_\star$.
\end{itemize} This proves in all cases that $\beta=\beta_\star(\gamma)$, resulting in
\[
\Lagr(\by;\beta) \approx g_\gamma(\beta_\star(\gamma))
    = g(\gamma)\,,
\] as claimed. Now define $\btheta\in\R^M$ as the concatenation of $\btheta_+$ and $\btheta_-$, and similarly $\hbg \in \R^M$.
Let $\hbD \in \R^{M\times M}$ be the concatenation of $\hbD_+$ and $\hbD_-$. Altogether, the above yields
\[
    F^{\rvi} \approx \adjustlimits
    \sup_{\btheta_+ \in \R^{B_+}}
    \inf_{\btheta_- \in \R^{B_-}} \lt\{
        - \frac{(\hbD^{-1} \btheta, \btheta)}{M}
        + \frac{2(\hbg,\btheta)}{M}
        + g\lt(\frac{\|\btheta\|^2}{M\alpha}\rt)
    \rt\}\,.
\]
Stationarity in $\btheta$ yields
\beq\label{e:fp.heuristic.btheta}
    \btheta = \lt(
        \hbD^{-1} - \frac{g'(\gamma)}{\alpha} \bI_M
    \rt)^{-1} \hbg
    = \lt(
        \hbD^{-1} - \frac{m(\beta_\star(\gamma))}{\alpha} \bI_M
    \rt)^{-1} \hbg\,,
\eeq
where the last equality is by Lemma~\ref{l:fp-g2}\eqref{i:fp-g-diff}. It follows that
    \begin{align*}\nonumber
    \gamma &= \frac{\|\btheta\|^2}{M\alpha}
    = \frac{1}{M\alpha} \bigg\|
        \bigg(
            \hbD^{-1} - \frac{m(\beta_\star(\gamma))}{\alpha} \bI
        \bigg)^{-1} \hbg
    \bigg\|^2
    \approx \frac1\alpha \sum_{j=1}^r \lambda_j \bigg(
    \frac{1}{\hv_j}-\frac{m(\beta_\star(\gamma))}{\alpha}
    \bigg)^{-2}
    \\
& = \alpha \sum_{j=1}^r \lambda_j
\bigg(
\frac{ \hv_j}{\alpha - \hv_j m(\beta_\star(\gamma))}
\bigg)^2 = \gamma(\beta_\star(\gamma))\,,
\end{align*}
where in the last step we use the definition of the function $\gamma(z)$ appearing in  Lemma~\ref{l:fp-vii-estimate-prep}. This gives a fixed-point equation in $\gamma$, which we rewrite as
    \beq\label{e:why-z-star}
    m(\beta_\star(\gamma))^2 \gamma
    = m(\beta_\star(\gamma))^2 \gamma(\beta_\star(\gamma))\,.
    \eeq
We will show this equation can have at most one solution $\gamma$. To this end, recall from Lemma~\ref{l:fp-g2}\eqref{i:fp-beta-incr} that $\beta_\star(\gamma)$ is nondecreasing in $\gamma$. Also recall that it was noted in the proof of Lemma~\ref{l:fp-vii-estimate-prep} that $m(z)^2\gamma(z)$ is strictly decreasing in $z$. It follows that the right-hand side of \eqref{e:why-z-star} is nonincreasing in $\gamma$. Abbreviating $f(\gamma)$ for the left-hand side of \eqref{e:why-z-star}, we calculate
    \[f'(\gamma)
    =
    m(\beta_\star(\gamma))
    \bigg\{ m(\beta_\star(\gamma))+ \gamma 2
        m'(\beta_\star(\gamma)) (\beta_\star)'(\gamma)
        \bigg\}\,.
    \]
Recalling the definition of $\beta_\star$ from \eqref{e:fp-betast}, we see that if $\hx>\alpha/m(\beta_{\min})$ then $(\beta_\star)'(\gamma)=0$, so in this case clearly $f'(\gamma)>0$. Otherwise, $\beta_\star=\beta_{\min}$, in which case we have $m'(\beta_\star)=-1/\gamma$ and consequently
    \[
     (\beta_\star)'(\gamma) = \frac{1}{m''(\beta_\star(\gamma)) \gamma^2 }
     = \frac{m'(\beta_\star(\gamma))^2}{m''(\beta_\star(\gamma))}
     = \frac{m_2(\beta_\star(\gamma))^2}{2m_3(\beta_\star(\gamma))}\,.
    \]
It follows in this case that
    \[
    f'(\gamma)
    = \frac{m(\beta_\star(\gamma))}{ m_3(\beta_\star(\gamma))}
    \bigg\{ m(\beta_\star(\gamma)) m_3(\beta_\star(\gamma))
    - m_2(\beta_\star(\gamma))^2
        \bigg\}\ge0\,,
    \]
using the Cauchy--Schwarz inequality.
In the case $s\ge2$, then the Cauchy--Schwarz inequality in fact implies $f'(\gamma)>0$: thus, the left-hand side of \eqref{e:why-z-star} is strictly increasing in $\gamma$ while the right-hand side is nonincreasing, so the equation has at most one solution $\gamma$. In the case $s=1$, recalling \eqref{e:beta-star-transition}, there are two scenarios: 
\begin{itemize}
\item For $\gamma\le\tilde{\gamma}$, the function $\beta_\star(\gamma)$ is constant in $\gamma$, so the right-hand side of \eqref{e:why-z-star} is constant in $\gamma$ while the left-hand side is increasing.
\item For $\gamma\ge\tilde{\gamma}$, the function $\beta_\star(\gamma)=\beta_{\min}(\gamma)$ is strictly increasing in $\gamma$, so the right-hand side of \eqref{e:why-z-star} is strictly decreasing in $\gamma$. Meanwhile the left-hand side of \eqref{e:why-z-star}
    \[
    m(\beta_\star(\gamma))^2\gamma
    =m(\beta_{\min}(\gamma))^2\gamma
    =m_2(\beta_{\min}(\gamma))\gamma
    =1/\gamma\,,
    \]
having used that $m_2(z)=m(z)^2$ for $s=1$.
\end{itemize}
Thus, in all cases, the difference between the right-hand side and left-hand side of \eqref{e:why-z-star} is strictly decreasing in $\gamma$, and so \eqref{e:why-z-star} has at most one solution.

If we take $\gamma=\gamma_\star\equiv \gamma(z_\star)$, then
Lemma~\ref{l:fp-vii-estimate-prep} tells us that $\beta_\star(\gamma) = z_\star$: this tells us that $\gamma=\gamma_\star$ satisfies \eqref{e:why-z-star}, and thus it must be the unique solution of \eqref{e:why-z-star}. It follows that
\baln
    F^{\rvi}
    &=-\frac{(\hbD^{-1} \btheta, \btheta)}{M}
    + \frac{2(\hbg,\btheta)}{M}
    + g\lt(\frac{\|\btheta\|^2}{M\alpha}\rt) \\
    &= -\frac{((\hbD^{-1} - m(z_\star) \bI_M / \alpha) \btheta, \btheta)}{M}
    - \frac{m(z_\star) \|\btheta\|^2}{M\alpha}
    + \frac{2(\hbg,\btheta)}{M}
    + g\lt(\frac{\|\btheta\|^2}{M\alpha}\rt) \\
    &\approx - \gamma_\star m(z_\star) + g(\gamma_\star) + \frac{(\hbg,\btheta)}{M}
    \approx z_\star + \alpha \sum_{j=1}^r \frac{\lambda_j \hv_j}{\alpha - \hv_j m(z_\star)}
    = V(z_\star) = F_\diamond\,.
\ealn
\textbf{This concludes our heuristic derivation of the assertion of Proposition~\ref{p:fp-v-estimate}.} We next turn to the rigorous proof, beginning with a few preparatory lemmas.

\begin{lem}\label{l:fp-theta-plus-bdd}
There exists a constant $C_8=C_8(\alpha,\iota)$ such that the following holds with probability at least $1-e^{-cN}$ over $\dot{\bg}_+, \hat{\bg}_+$, where $c$ is an absolute constant. 
    For any $\|\by\|^2=N$, the function
    \[
        \btheta_+ \mapsto
        \frac{2(\hbg_+,\btheta_+)}{M}
        + \frac{2\|\btheta_+\|(\dbg_+,\by)}{M\sqrt{N}}
        - \frac{((\hbD_+)^{-1} \btheta_+, \btheta_+)}{M}
    \]
    attains its supremum over $\R^{B_+}$ at some $\btheta_+$ satisfying
    $\|\btheta_+\| \le C_8N^{1/2}$.
    (We take the assertion to be vacuously true if $B_+ = \emptyset$.)

\begin{proof}
    Since the function tends to $-\infty$ as $\|\btheta_+\| \to +\infty$, the supremum is attained at a local maximizer $\btheta_+$.
    The function is differentiable everywhere except possibly at $\btheta_+=\bzero$, so the only possible local maximizer are $\btheta_+=\bzero$ and stationary points $\btheta_+$. 
    The stationarity condition gives
    \[
      \btheta_+
     =\hbD_+\bigg(   \hbg_+
        + \frac{(\dbg_+,\by) \btheta_+}{N^{1/2} \|\btheta_+\|}
        \bigg) \,.
    \]
With probability $1-e^{-cN}$, we have $\|\hbg_+\| \le 2 M^{1/2}$, $\|\dbg_+\| \le 2 N^{1/2}$, and consequently $|(\dbg_+,\by)| \le 2N$. Since the operator norm of $\hbD_+$ is at most $\hv_1$, it follows that any stationary point $\btheta_+$ must satisfy
    \[
        \|\btheta_+\|
        \le\hv_1 \bigg(
            \|\hbg_+\|
            + \bigg|
            \frac{(\dbg_+,\by) }{N^{1/2} }
            \bigg|
        \bigg)
        \le
        2\hv_1\Big( M^{1/2} + N^{1/2}
         \Big)
        \stackrel{\eqref{e:reduced-normalization}}{\le}
        \frac{4(\alpha^{1/2}+1)}{\iota} N^{1/2}\,.
    \]
This proves the claim.
\end{proof}
\end{lem}

In light of Lemma~\ref{l:fp-theta-plus-bdd}, the supremum in $\btheta_+$ in \eqref{e:fp-F-v} can be restricted to the domain
    \beq\label{e:fp.restricted.Theta}
    \bTheta_+ \equiv
    \Big\{\btheta_+
: \|\btheta_+\| \le C_8N^{1/2}
\Big\}\eeq
with probability $1-e^{-cN}$, where $c$ is an absolute constant. 

\begin{lem}\label{l:fp-gaussian-matrix-prod} For any constants $C,\delta > 0$, there exists a positive constant $c(C,\delta)$ such that the following holds: let $\bA$ be any $N\times N$ diagonal matrix with $\|\bA\|_{\op} \le C$, and suppose $\bg,\acute{\bg}$ are independent standard gaussian vectors in $\R^N$. Then,
    \[
        \frac{|(\bA \bg, \bg) - \Tr(\bA)|}{N}\,,
        \frac{|(\bA \bg, \acute{\bg})|}{N}
        \le \delta
    \]
    with probability at least $1-e^{- c(C,\delta) N}$.

\begin{proof}
We rewrite the random variables of interest as
    \begin{align*}
    (\bA \bg, \bg) - \Tr(\bA)
    &=
    \sum_{i=1}^N a_{i,i}
    \Big((g_i)^2-1\Big)\,,\\
    (\bA\bg,\bg')
    &=
    \sum_{i=1}^N a_{i,i}
        g_i\acute{g}_i\,.
    \end{align*}
The claim then follows by standard Chernoff bounds.
\end{proof}
\end{lem}

\noindent
In contrast with the heuristic derivation above, we will rigorously bound $F^{\rvi}$ by setting the value of $\btheta_-$ in the infimum, with the choice of $\btheta_-$ coming from the above heuristic calculation. The details of this are carried out in Proposition~\ref{p:fp-v-to-vi} below, and will result in the following
bound: recall the definition of $z_\star$ from \eqref{e:fp.z.star}, and let
\beq\label{e:fp.gamma.minus}
    \gamma_- 
    \equiv \alpha
    \sum_{j\in J_-}
    \lambda_j
    \bigg(
    \frac{ \hv_j}{\alpha
        - \hv_j m(z_\star)}\bigg)^2
    \le
    \gamma(z_\star)=\gamma_\star\,.\eeq
Recalling the definition of $F^{\rvi}$ from \eqref{e:fp-F-v}, we let
    \beq\label{e:fp.F.vii}
    F^{\rvii}
    =\sup_{\|\by\|^2=N}
    \bigg\{
    F_w(\by)
    +F_x(\tbmeta,\by)
    +\ddot{F}_{v+}(
    \hbg_+,\dbg_+,
    \by)
    +\ddot{F}_{v-}(\hbg_-,
    \dbg_-,\by)
    \bigg\}
    \eeq
where $\ddot{F}_{v+}$ is defined in the same way as $\bar{F}_{v+}$ from
\eqref{e:bar.F.v.plus} or \eqref{e:fp.recall.barFvplus}, but with the supremum taken over the restricted domain $\bTheta_+$ from \eqref{e:fp.restricted.Theta}. Meanwhile, $\ddot{F}_{v-}$ is defined as
    \beq\label{e:fp.ddot.F.vminus}
    \ddot{F}_{v-}
    (\dbg_-,\by)
    \equiv
    \alpha \sum_{j\in J_-} \frac{\lambda_j \hv_j}{\alpha - \hv_j m(z_\star)}
    +
    2(\gamma_-)^{1/2}
    \frac{(\dbg_-,\by)}{N}
        - m(z_\star) \gamma_-\,.
    \eeq
(Note that the first term on the right-hand side of \eqref{e:fp.ddot.F.vminus} corresponds to one of the terms in the definition \eqref{e:freeprob.V} of $V(z_\star)$.)

\begin{ppn}
\label{p:fp-v-to-vi}
Let $F^{\rvi}$ and $F^{\rvii}$ be as defined by \eqref{e:fp-F-v} and \eqref{e:fp.F.vii} respectively. There exists a positive constant $c$, depending on all the parameters in Remark~\ref{r:fp-rpr}, such that we have
    \[
        \P(F^{\rvi} \ge t) \le \P\bigg(F^{\rvii} \ge t + \frac{\epsilon'}{20}\bigg)
         + e^{-cN}\,.
    \]
for all $t\in\R$. That is, $F^{\rvi}$ is approximately stochastically dominated by $F^{\rvii}$.

\begin{proof}
It follows from Lemma~\ref{l:fp-theta-plus-bdd} that
$\bar{F}_{v+}$ and $\ddot{F}_{v+}$
coincide with probability at least $1-e^{-cN}$, for $c$ an absolute constant. Thus, the main part of the proof is to compare $\bar{F}_{v-}$ and $\ddot{F}_{v-}$, as defined by \eqref{e:fp.recall.barFvminus} and \eqref{e:fp.ddot.F.vminus} respectively.
If $B_-$ is empty then there is nothing to show, so assume that it is nonempty. Then, since $(\lambda,\mu) \in \Lambda$, we have $|B_-| \ge \iota N$. The definition \eqref{e:fp.recall.barFvminus}
of $\bar{F}_{v-}$ involves an infimum
 over $\btheta_-$, so we can upper bound its value by choosing any value of $\btheta_-$. We choose it according to the heuristic calculation \eqref{e:fp.heuristic.btheta}: recalling from the discussion following \eqref{e:why-z-star} that we take $\gamma$ in \eqref{e:fp.heuristic.btheta} such that $\beta_\star(\gamma)=z_\star$, we now take
    \[
    \btheta_-
    = \hbA_- \hbg_-\,,\quad
    \textup{where}\quad
    \hbA_- = \bigg((\hbD_-)^{-1} - \frac{m(z_\star)}{\alpha}\bI_{B_-}\bigg)^{-1}\,.
    \]
Recall that $\hbD_-$ is a diagonal matrix with (negative) diagonal entries $\hv_j$. It follows that $\hbA_-$ is a diagonal matrix with (negative) diagonal entries
    \[\bigg(
    \frac{1}{\hv_j}-\frac{m(z_\star)}{\alpha}\bigg)^{-1}
    = \frac{\alpha\hv_j}{\alpha-\hv_j m(z_\star)} \ge \hv_j\,,
    \]
from which we see that the operator norm of $\hbA_-$ satisifes
\[\|\hbA_-\|_\textup{op}\le\max\{ -\hv_j : j\in J_-\}  \le \frac2\iota\,,\] where the last estimate is by \eqref{e:reduced-normalization}.
    Let  $\delta$ be a sufficiently small constant depending on the parameters in Remark~\ref{r:fp-rpr}, which we will set below, and let $\approx_\delta$ denote equality up to additive $O(\delta)$ error.
    By Lemma~\ref{l:fp-gaussian-matrix-prod} (with $|B_-|$ in place of $N$, and $2\iota^{-1}$ in place of $C$), we have
    \begin{align*}
        \frac{(\hbg_-,\btheta_-)}{M}
        &\approx_\delta
        \frac{\Tr(\hbA_-)}{M}
        = \alpha \sum_{j\in J_-} \frac{\lambda_j \hv_j}{\alpha
        -\hv_j m(z_\star) }\,, \\
        \frac{\|\btheta_-\|^2}{M}
        &\approx_\delta
        \frac{\Tr((\hbA_-)^2)}{M}
        = \alpha \gamma_-\,,
    \end{align*}
with probability at least $1-e^{-c(\delta) N}$, where $c(\delta)$ depends on $\delta$ and implicitly on all the parameters in Remark~\ref{r:fp-rpr}. 
Returning to the definition
\eqref{e:fp.recall.barFvminus} of $\bar{F}_{v-}$, and taking the specific choice of $\btheta_-$ described above, we conclude that
    \begin{align*}
    &\bar{F}_{v-}(
    \hbg_-,\dbg_-,\by)
    \le
    \frac{2(\hbg_-,\btheta_-)}{M}
        + \frac{2\|\btheta_-\|
        (\dbg_-,\by)}{MN^{1/2}}
         - \frac{(\btheta_-,
         (\hbD_-)^{-1}\btheta_-)}
            {M} \\
    &\qquad
    =    \frac{2(\hbg_-,\btheta_-)}{M}
        + \frac{2\|\btheta_-\|
        (\dbg_-,\by)}
        {MN^{1/2}}
    - \frac{(\btheta_-,
         [(\hbD_-)^{-1}
         -m(z_\star)\bI_{B_-} /\alpha
         ]
         \btheta_-)}
            {M}
    - m(z_\star)
     \frac{\|\btheta_-\|^2}{M\alpha}\\
    &\qquad
    = \frac{(\hbg_-,\btheta_-)}{M}
    +2 \frac{\|\btheta_-\|
        }{(M\alpha)^{1/2} }
        \frac{(\dbg_-,\by)}{N}
        - m(z_\star)
     \frac{\|\btheta_-\|^2}{M\alpha}\\
    &\qquad\approx_\delta
    \alpha \sum_{j\in J_-} \frac{\lambda_j \hv_j}{\alpha
        -\hv_j m(z_\star) }
        + 2(\gamma_-)^{1/2}
        \frac{(\dbg_-,\by)}{N}
        -m(z_\star)\gamma_-
        =\ddot{F}_{v-}
    (\hbg_-,
    \dbg_-,\by)
    \end{align*}
with probability at least $1-e^{- c(\delta) N}$. For the last step, we used the above calculations, together with the observation that $|(\dbg_-,\by)| \le 2N$ with probability at least $1-e^{-cN}$.
The claim follows by setting $\delta$ small enough (depending on the parameters in Remark~\ref{r:fp-rpr}) that the errors in the $\approx_{\delta}$ sum to less than $\epsilon'/20$.
\end{proof}
\end{ppn}

We next make rigorous the optimization over $\|\by\|^2 = N$ via the method of Lagrange multipliers, as done heuristically following \eqref{e:fp-y-lagrangian} above. Let $C_9>0$ be a sufficiently large constant, depending on only the parameters in Remark~\ref{r:fp-rpr}, which we will set below.

\begin{ppn}\label{p:fp-g}
Let $\dbg$ be a standard gaussian in $\R^N$ and 
\beq
    \label{e:def-V-gamma-y}
    V_\gamma(\by)\equiv \frac{(\dbD \by, \by)}{N}
            +\frac{(\dbX \by, \by)}{N}
            + \frac{2\gamma^{1/2}(\dbg,\by)}{N}\,.
\eeq
Then, with probability $1-e^{-cN}$, 
we have
    \beq
        \label{e:fp-g}
        \sup \Big\{ V_\gamma(\by) :\|\by\|^2=N\Big\}
        \le g(\gamma)
        + \frac{\epsilon'}{10}
    \eeq
simultaneously for all $\gamma \in [0, C_9]$. The constant $c$ depends on $C_9$ and the parameters of Remark~\ref{r:fp-rpr}. 
\end{ppn}

The \hyperlink{proof:p.fp-g}{proof of Proposition~\ref{p:fp-g}} appears below,  after several preparatory lemmas. 

\begin{lem}\label{l:fp-g-prep}
For any constant $C,\delta > 0$, there exists $c(C,\delta) > 0$ such that the following holds. Let $\tbmeta^j$ be i.i.d.\ standard gaussian vectors in $\R^N$.
For $\vc = (c_0,\ldots,c_{r'}) \in \R^{r'+1}$, define
    \beq
        \label{e:fp-bybeta}
        \by(\beta)
        \equiv \dbA_\beta  \sum_{j=0}^{r'}
            c_j \tbmeta^j\,,
    \eeq
where we recall $\dbA_\beta$ is defined below \eqref{e:fp-g-stationarity}. Then, with probability $1-e^{-c(C,\delta)N}$, we have
the bounds
    \begin{align}
        \label{e:fp-g-prep-conclusion1}
        \bigg|\frac{\|\by(\beta)\|^2}{N} - \|\vc\|^2 m_2(\beta)\bigg|
        &\le \delta\,,\\
        \label{e:fp-g-prep-conclusion2}
        \max\bigg\{
        \bigg|\frac{(\tbmeta^j,\by(\beta))}{N} - c_j m(\beta) \bigg|
        : 0 \le j \le r' 
        \bigg\}
      &\le \delta\,,
    \end{align}
simultaneously
for all $\beta \in [\hw_1+ C^{-1},+\infty)$ and all  $\|\vc\|_\infty \le C$.

\begin{proof}
We first consider a fixed $\beta \in [\hw_1+1/C,+\infty)$: it follows from the definition of $\dbA_\beta$ that for any such $\beta$ we have $\|\dbA_\beta\|_\textup{op}\le C$. Thus Lemma~\ref{l:fp-gaussian-matrix-prod} can be applied, which gives that 
with probability at least $1-e^{- c(C,\delta)N}$, for all $j,k \in \{0,\ldots,r\}$ with $j\ne k$ we have
    \beq
        \label{e:fp-g-prep-event1}
        \max\bigg\{
        \bigg|\frac{((\dbA_\beta)^2 \tbmeta^j, \tbmeta^j)}{N}
            - m_2(\beta)\bigg|,
        \frac{|((\dbA_\beta)^2 \tbmeta^j, \tbmeta^k)|}{N}
        \bigg\}
        \le \delta_1 \equiv \frac{\delta}{C^2(r'+1)^2}\,.
    \eeq
Taking a linear combination of the above quantities implies the first conclusion \eqref{e:fp-g-prep-conclusion1} 
for all $\|\vc\|_\infty \le C$.
    Similarly, for all $j,k \in \{0,\ldots,r\}$ with $j\ne k$ we have
    \beq
        \label{e:fp-g-prep-event2}
        \lt|\frac{(\dbA_\beta \tbmeta^j, \tbmeta^j)}{N} - m(\beta)\rt|,
        \frac{|(\dbA_\beta \tbmeta^j, \tbmeta^k)|}{N}
        \le \delta_1\,,
    \eeq
and taking a linear combination implies the second conclusion \eqref{e:fp-g-prep-conclusion2}, again
for all $\|\vc\|_\infty \le C$.

Finally, we show that this event holds simultaneously over all $\beta \in [\hw_1+C^{-1},+\infty)$ with probability at least $1-e^{-c(C,\delta)N}$.
Note that there exists $T=T(C,\delta)$ and $\beta_1,\ldots,\beta_T \in [\hw_1+C^{-1},+\infty)$ such that for any $\beta \in [\hw_1+C^{-1},+\infty)$, there exists $t\in[T]$ such that
\begin{align*}
    \max\Big\{
        \|\dbA_\beta - \dbA_{\beta_t}\|_{\op}\,,
        \|(\dbA_\beta)^2 - (\dbA_{\beta_t})^2\|_{\op}
    \Big\} &\le \delta_1\,,\\
    \max\Big\{
    |m(\beta) - m(\beta_t)|,
    |m_2(\beta) - m_2(\beta_t)|\Big\}
    &\le \delta_1\,.
\end{align*}
Consider the event that the estimates \eqref{e:fp-g-prep-event1} and \eqref{e:fp-g-prep-event2} hold for $\beta_t$ for all $t\in [T]$, and moreover that $\|\tbmeta^j\| \le 2N^{1/2}$ for all $j$. This event occurs with probability $1-e^{-c(C,\delta) N}$. On this event, for any $j,k\in \{0,\ldots,r\}$, we can bound
\[
    \frac{|((\dbA_\beta)^2 \tbmeta^j, \tbmeta^k)-((\dbA_{\beta_t})^2 \tbmeta^j, \tbmeta^k)|}{N}
    \le 4\|(\dbA_\beta)^2 - (\dbA_{\beta_t})^2\|_{\op} \le 4\delta_1\,.
\]
It follows from the triangle inequality that
\baln
    \bigg|\frac{((\dbA_\beta)^2 \tbmeta^j, \tbmeta^j)}{N} - m_2(\beta)\bigg|
    &\le 
    \frac{|((\dbA_\beta)^2 \tbmeta^j, \tbmeta^j)-((\dbA_{\beta_t})^2 \tbmeta^j, \tbmeta^j)|}{N}\\
    &\qquad+
    \bigg|\frac{((\dbA_{\beta_t})^2 \tbmeta^j, \tbmeta^j)}{N} - m_2(\beta_t)\bigg|
    + |m_2(\beta_t)-m_2(\beta) | 
    \le 6\delta_1\,,
\ealn
and for $j\neq k$, similarly
\[
    \frac{|((\dbA_\beta)^2 \tbmeta^j, \tbmeta^k)|}{N}
    \le \frac{|((\dbA_\beta)^2 \tbmeta^j, \tbmeta^k)-((\dbA_{\beta_t})^2 \tbmeta^j, \tbmeta^k)|}{N} 
    + \frac{|((\dbA_{\beta_t})^2 \tbmeta^j, \tbmeta^k)|}{N}
    \le 6\delta_1\,.
\]
Taking a linear combination of the above quantities shows that with probability at least $1-e^{-c(C,\delta)N}$, 
 the conclusions \eqref{e:fp-g-prep-conclusion1}
and \eqref{e:fp-g-prep-conclusion2} hold for all $\beta\in[\hw_1+C^{-1},\infty)$ and all $\|\vc\|_\infty \le C$,
with $6\delta$ in place of $\delta$.
The conclusion follows by adjusting $\delta$.
\end{proof}
\end{lem}

\begin{lem}\label{l:fp-g-prep2}
Recall the function $V_\gamma(\by)$ from \eqref{e:def-V-gamma-y}.
For any $\gamma \ge 0$, any local maximizer of $V_\gamma(\by)$ on the sphere $\{\|\by\|^2=N\}$ is of the form
    \beq
        \label{e:fp-g-prep2-stationarity}
        \by =\dbA_\beta\Big(\gamma^{1/2}\dbg+\dbX\by\Big)\,,
    \eeq
    where we recall $\dbA_\beta \equiv (\beta \bI_N - \dbD)^{-1}$. Furthermore, for any such $\by$ we must have $\beta \ge \hw_1$.
\begin{proof}
We apply the method of Lagrange multipliers with the Lagrangian $\Lagr(\by;\beta)$ as defined in \eqref{e:fp-y-lagrangian}. Note that the Lagrangian Hessian (with respect to $\by$) is given by
    \[
        \nabla^2 \Lagr(\by;\beta)
        = -\frac{2}{N} \Big(
            \beta \bI_N - \dbD - \dbX
        \Big)\,.
    \]
At any local maximum, this Hessian must be negative semi-definite. Recall the definitions of $\dbD$ and $\dbX$ from \eqref{e:fp.dbD} and \eqref{e:fp.dbX} respectively. Since $|B_1|\ge\iota N$ by the assumption $(\lambda,\mu) \in \Lambda$, we can find a non-zero vector $\bv$ that is zero outside the coordinates of $B_1$, and is additionally orthogonal to each of the vectors $\tbmeta^j$ ($1\le j\le r'$) in the definition \eqref{e:fp.dbX} of $\dbX$. It follows that
    \[
    0 \le \bv^\st \Big(\beta -\dbD-\dbX \Big)\bv
    =\bv^\st \Big(\beta -\dbD \Big)\bv
    =(\beta-\hw_1)\|\bv\|^2\,,
    \]
so we conclude $\beta\ge\hw_1$. 
For any local maximizer $\by$ satisfying $\|\by\|^2=N$, stationarity of $\Lagr$ at $\by$ yields the expression~\eqref{e:fp-g-prep2-stationarity}.
\end{proof}
\end{lem}

We next prove Proposition~\ref{p:fp-g} holds for a \emph{fixed} $\gamma$.  In the \hyperlink{proof:p.fp-g}{proof of Proposition~\ref{p:fp-g}} below, we will extend the result to hold simultaneously for all $\gamma \in [0,C_9]$. 

\begin{ppn}
\label{p:fp-g-one-gamma}
For any fixed $C_9$, and any fixed $\gamma \in [1/C_9,C_9]$, we have
    \[
        \sup \Big\{ V_\gamma(\by) :\|\by\|^2=N\Big\}
         \le g(\gamma) + \frac{\epsilon'}{20}\,.
    \]
with probability at least $1-e^{-cN}$, where $c$ depends on $C_9$ and the parameters in Remark~\ref{r:fp-rpr}. 

\begin{proof}In this proof, we will set constants $\upsilon, \upsilon_1, \delta > 0$, where $\upsilon$ is sufficiently small depending on $C_9$, $\upsilon_1$ sufficiently small depending on $\upsilon,C_9$, and $\delta$ sufficiently small depending on $\upsilon_1,\upsilon,C_9$. That is,
\beq\label{e:fp-g-one-gamma-parameters-order}
    \delta \ll \upsilon_1 \ll \upsilon \ll \frac1{C_9}\,.
\eeq
Note that $C_9, \upsilon, \upsilon_1, \delta$ nonetheless depend on only the parameters in Remark~\ref{r:fp-rpr}.

Consider any local maximizer $\by$, and recall that $\dbA_\beta \equiv (\beta I_N-\dbD)^{-1}$.
By Lemma~\ref{l:fp-g-prep2}, $\by$ is of the form \eqref{e:fp-g-prep2-stationarity}, which we expand as
    \[
    \by
    =\dbA_\beta\bigg(\gamma^{1/2}\dbg
        +\sum_{j=1}^{r'}
        \frac{\hx_j}{M}(\tbmeta^j,\by)
        \tbmeta^j
        \bigg)
    =\dbA_\beta\sum_{j=0}^{r'} c_j \tbmeta^j\,.
    \]
For the last equality, we have defined $c_0\equiv \gamma^{1/2}$, $\tbmeta^0\equiv\dbg$, and 
    \[c_j\equiv \frac{\hx_j(\tbmeta^j, \by)}{M}\]
    for $1\le j\le r'$. With probability at least  $1-e^{-cN}$ (with $c$ an absolute constant), we have $\|\tbmeta^j\| \le 2\sqrt{N}$ for all $j \in 1[r']$, which implies 
        \[|c_j| \le \frac{2 |\hx_j|}{\alpha}
        \stackrel{\eqref{e:reduced-normalization}}{\le} \frac{4\alpha}{\iota}
        \]
for all $j\in[r']$. On this event, we have
\[
    \|\vc\|_\infty \le \max\bigg\{
        \gamma^{1/2}, \max\bigg\{
            \frac{2 |\hx_j|}{\alpha} : j\in [r']
        \bigg\}
    \bigg\} \le
    \max\bigg\{
        (C_9)^{1/2}, \frac{4}{\alpha \iota}
    \bigg\}
    \equiv C_{10}\,.
\]
Because $\beta_{\min} = \beta_{\min}(\gamma)$ solves $m_2(\beta_{\min}) = 1/\gamma \in [1/C_9,C_9]$, we have
    \begin{align*}
    C_9 &\ge \frac1\gamma
    = m_2(\beta_{\min}) 
    \ge \frac{\mu_1}{(\beta_{\min}-\hw_1)^2} \,,\\
    \frac1{C_9} &\le \frac1\gamma = 
    m_2(\beta_{\min}) \le \frac1{(\beta_{\min}-\hw_1)^2} \,,\end{align*}
Rearranging these inequalities gives (recalling $\mu_1 \ge \iota$ because $(\lambda,\mu) \in \Lambda$) 
\beq
    \label{e:fp-g-beta-min-estimate}
    (C_9)^{1/2} \ge
    \beta_{\min} - \hw_1
    \ge \bigg( \frac{\mu_1 }{ C_9}\bigg)^{1/2} 
    \ge
    \bigg(  \frac{\iota}{ C_9}\bigg)^{1/2}
    \equiv \upsilon_0\,.
\eeq
Therefore, for all $t\in \bbR$ with  $\hw_1-\beta_{\min} < t \le (C_9)^{1/2}$, we have
\beq
    \label{e:fp-m2-deriv-lb}
    - (m_2)'(\beta_{\min}+t) 
    = 2 m_3(\beta_{\min}+t)
    \ge \frac{2\mu_1}{(\beta_{\min}+t-\hw_1)^3}
    \ge
    \frac{\iota}{4(C_9)^{3/2}}\,.
\eeq
Given $\by$, let $\by(\beta)$ be defined as above for any value of $\beta$: that is,
    \[
    \by(\beta)
    \equiv \dbA_\beta
    \sum_{j=0}^{r'} c_j\tbmeta^j\,,
    \]
where $\vc$ is defined above (and depends on $\by$).  Thus, for any fixed $\by$, we may regard $\by(\beta)$ as function of $\beta$. At a maximizer of the constrained optimization problem, we must have $\by=\by(\beta)$. We now take
    \[\upsilon \in \bigg(0,
        \min\bigg\{ \frac{\upsilon_0}{2}, (C_9)^{1/2}\bigg\}\bigg)\,.\]
Note that \eqref{e:fp-g-beta-min-estimate} implies $\beta_{\min} - \upsilon \ge \hw_1 + \upsilon_0 / 2$. It follows that
\beq
    \label{e:fp-g-dbA-bd}
    \|\dbA_{\beta_{\min} - \upsilon}\|_{\op} 
    = \frac{1}{\beta_{\min} - \upsilon - \hw_1} 
    \le \frac{1}{\upsilon_0/2}
    = 2 \bigg(\frac{C_9}{\iota}\bigg)^{1/2}
        \equiv C_{11}\,.
\eeq
By Lemma~\ref{l:fp-g-prep} with $\max(C_{10},C_{11})$ in place of $C$, we have
    \begin{align*}
      &\frac{\|\by(\beta_{\min} - \upsilon)\|^2}{N}
        \approx_\delta
        \|\vc\|^2 m_2(\beta_{\min} - \upsilon)
        \ge \gamma m_2(\beta_{\min} - \upsilon) \\
       &\qquad\stackrel{\eqref{e:fp-m2-deriv-lb}}{\ge}
        \gamma \bigg(m_2(\beta_{\min}) + \frac{\upsilon \iota}{4C_9^{3/2}}\bigg)
        = 1 + \frac{\gamma \upsilon \iota}{4C_9^{3/2}}\,. 
    \end{align*}
with probability at least $1-e^{-cN}$, where now $c$ depends on $C_9, \alpha,\iota,\delta$. We can take $\delta$ small enough depending on $\upsilon$ (and on $C_9,\iota$) 
to conclude 
        \[\frac{\|\by(\beta_{\min} - \upsilon)\|^2}{N} > 1\,.\]
Since $\beta \mapsto \|\by(\beta)\|^2$ is decreasing on $\beta>\hw_1$, this implies
$\|\by(\beta)\|^2 / N > 1$ for all $\beta \in (\hw_1,\beta_{\min} - \upsilon]$.
However, if $\by$ is a maximizer of the constrained optimization problem, then $\by=\by(\beta)$ must have norm $N^{1/2}$. 
So we conclude that for any such maximizer $\by=\by(\beta)$, we must have $\beta > \beta_{\min} - \upsilon$.
If $\by=\by(\beta)$, then, as in \eqref{e:fp-y-lagrangian-bd}, we have
    \[
        V_\gamma(\by)
        = \Lagr(\by;\beta)
        = \beta + \frac{\gamma^{1/2}(\dbg,\by)}{N}\,.
    \]
    Since $\beta > \beta_{\min} - \upsilon \ge \hw_1+ \upsilon_0 / 2 = \hw_1 + 1/C_{11}$, Lemma~\ref{l:fp-g-prep} further implies
    \[
        V_\gamma(\by)
        \approx_\delta
        \beta + \gamma m(\beta)
        = g_\gamma(\beta)\,.
    \]
To finish the proof, we will bound $g_\gamma(\beta)$, which will result from restricting the possible values of $\beta$. 
By Lemma~\ref{l:fp-g-prep}, for any constant $\delta > 0$ and any $j\in [r']$, we have
    \[
        \frac{(\tbmeta^j,\by)}{N}
        \approx_\delta c_j m(\beta) = 
        \frac{(\tbmeta^j,\by)}{N}
        \cdot \frac{\hx_j}{\alpha} m(\beta)\,,
    \]
using the above definition of $c_j$. Rearranging gives
\[
    \frac{|(\tbmeta^j,\by)|}{N} \cdot 
    \bigg|\frac{\hx_j}{\alpha} m(\beta) - 1\bigg|
    \le O(\delta)\,.
\]
If $\delta$ is set sufficiently small depending on $\upsilon_1$, this implies 
    \[
    \min\bigg\{
    \frac{|(\tbmeta^j,\by)|}{N} ,
    \bigg|\frac{\hx_j}{\alpha} m(\beta) - 1\bigg|
    \bigg\} < \upsilon_1
    \]
for all $j\in[r']$. We divide into cases:
\begin{itemize}
\item First consider the case that $|(\tbmeta^j,\by)|/N<\upsilon_1$ for all $j\in[r']$. Using Lemma~\ref{l:fp-g-prep}, we have
    \baln
        \frac{\|\by(\beta_{\min} + \upsilon)\|^2}{N}
        &\approx_\delta
        \|\vc\|^2 m_2(\beta_{\min} + \upsilon)
        \le \bigg(\gamma +  \sum_{j=1}^{r'}
            \bigg(\frac{\upsilon_1\hx_j}{\alpha}\bigg)^2
    \bigg) m_2(\beta_{\min} + \upsilon) \\
        & \stackrel{\eqref{e:reduced-normalization},\eqref{e:fp-m2-deriv-lb}}{\le} \bigg(\gamma + \frac{4 r'(\upsilon_1)^2}{\alpha^2\iota^2}\bigg)
        \bigg(m_2(\beta_{\min}) - \frac{\upsilon \iota}
            {4(C_9)^{3/2}} \bigg)\,.
    \ealn
The above application of \eqref{e:fp-m2-deriv-lb} uses $\upsilon \le \sqrt{C_9}$, which is clearly guaranteed by the condition~\eqref{e:fp-g-one-gamma-parameters-order}. 
Recall that $m_2(\beta_{\min})=1/\gamma$.
By taking $\delta,\upsilon_1$ sufficiently small depending on $\upsilon$ (and on $C_9,\alpha,\iota$), we conclude 
	\[\frac{\|\by(\beta_{\min} + \upsilon)\|^2}{N} < 1\,.\]
Since $\beta \mapsto \|\by(\beta)\|^2$ is strictly decreasing on $\beta\ge\hw_1$, we conclude $\beta < \beta_{\min} + \upsilon$.
Thus, 
    \beq
        \label{e:fp-v-gamma-by-ub-case-no-spike}
        V_\gamma(\by)
        \le \sup \bigg\{g_\gamma(\beta) :
        \beta \in (\beta_{\min} - \upsilon, \beta_{\min} + \upsilon)
        \bigg\} + O(\delta)\,.
    \eeq

\item Otherwise, suppose that we have
$|\hx_j m(\beta)/\alpha - 1| < \upsilon_1$ for some $j\in[r']$.
Then,
    \[
        \frac{\hx}{\alpha} m(\beta)
        \ge 
        \frac{\hx_j}{\alpha} m(\beta)
        > 1-\upsilon_1\,.
    \]
    This implies $\hx > 0$ and
    \[\beta < m^{-1}\bigg(
        \frac{(1-\upsilon_1) \alpha}{\hx}
    \bigg) = m^{-1}\bigg(\frac{\alpha}{\hx}\bigg) + O(\upsilon_1)
    \le \beta_\star + O(\upsilon_1)\,,
    \,,\] 
where $\beta_\star \ge m^{-1}(\alpha/\hx)$ from the definition \eqref{e:fp-betast}, and we have used Lemma~\ref{l:fp-scalar-fn-lipschitz}\eqref{i:fp-m-inv-lipschitz} to estimate
    \[
        \bigg|m^{-1}\bigg(
            \frac{(1-\upsilon_1) \alpha}{\hx}
        \bigg) - m^{-1}\bigg(\frac{\alpha}{\hx}\bigg) \bigg|
        \le \frac{C_6\hx}{\alpha} 
        \bigg|\frac{1}{1-\upsilon_1} - 1\bigg| 
        \stackrel{\eqref{e:reduced-normalization}}{\le} 
       \frac{3 C_6 \upsilon_1}{\alpha\iota}\,.
    \] 
We therefore conclude in this case that
\beq
        \label{e:fp-v-gamma-by-ub-case-with-spike}
        V_\gamma(\by)
        \le \sup\bigg\{ g_\gamma(\beta)
        :\beta \in \Big(
            \beta_{\min} - \upsilon, \beta_\star + O(\upsilon_1)
    \Big)\bigg\}
         + O(\delta)\,.
    \eeq
We also recall from \eqref{e:fp-g-one-gamma-parameters-order} that $\upsilon_1\ll\upsilon$. 
\end{itemize}
Finally, we will show that 
\beq
    \label{e:fp-v-gamma-sup}
    \sup\bigg\{ g_\gamma(\beta) : \beta \in [\beta_{\min} - \upsilon, \beta_\star + \upsilon]
    \bigg\}
    \le g(\gamma) + O(\upsilon)\,.
\eeq
Since $g_\gamma$ is convex, it suffices to show the bound at the endpoints $\beta \in \{\beta_{\min} - \upsilon, \beta_\star + \upsilon\}$.
Recall that $\beta_{\min} - \upsilon \ge \hw_1 + 1/C_{11}$, so for all $\beta \ge \beta_{\min}-\upsilon$,
\[
    -m'(\beta) = m_2(\beta) \le \frac{1}{(\beta - \upsilon - \hw_1)^2} \le (C_{11})^2\,.
\]
It follows that
\begin{align*}
    g_\gamma(\beta_{\min} - \upsilon)
    &= \beta_{\min} - \upsilon + \gamma m(\beta_{\min} - \upsilon)\\
    &\le g_\gamma(\beta_{\min}) + C_9( C_{11})^2 \upsilon
    \le g_\gamma(\beta_\star) + C_9 (C_{11})^2 \upsilon\,.
\end{align*}
On the other hand, recalling that $m$ is decreasing, we also have
\[
    g_\gamma(\beta_{\star} + \upsilon)
    =\beta_{\star} + \upsilon + \gamma m(\beta_{\star} + \upsilon)
    \le g_\gamma(\beta_\star) + \upsilon\,,
\]
Recalling $g_\gamma(\beta_\star) = g(\gamma)$ proves \eqref{e:fp-v-gamma-sup}.
Finally, the supremum in \eqref{e:fp-v-gamma-sup} upper bounds the suprema in \eqref{e:fp-v-gamma-by-ub-case-no-spike} and \eqref{e:fp-v-gamma-by-ub-case-with-spike} (for $\upsilon_1$ sufficiently small depending on $\upsilon$), so we conclude
\[
    V_\gamma(\by) \le g(\gamma) + O(\upsilon) + O(\delta)\,.
\]
Taking $\upsilon, \delta$ sufficiently small depending on $\epsilon'$ completes the proof.
\end{proof}
\end{ppn}

\begin{proof}[\hypertarget{proof:p.fp-g}{Proof of Proposition~\ref{p:fp-g}}]
By increasing $C_9$ if necessary (depending only on the parameters in Remark~\ref{r:fp-rpr}), we can ensure, using Lemma~\ref{l:fp-scalar-fn-lipschitz}\eqref{i:fp-g-sq-lipschitz}, that
\beq
    \label{e:fp-g-C2-choice}
    \frac{4}{(C_9)^{1/2}} + \bigg|g\bigg(\frac1{C_9}\bigg) - g(0)\bigg|
        \le\frac{\epsilon'}{20}\,.
\eeq
Using Lemma~\ref{l:fp-scalar-fn-lipschitz}\eqref{i:fp-g-sq-lipschitz} again, we can find a sequence of $T=O(1)$ numbers $C_9^{-1} = \gamma_1 < \cdots < \gamma_T = C_9$ such that for any $\gamma \in [1/C_9,C_9]$, there exists $t\in [T]$ such that
\beq
    \label{e:fp-g-net}
    4|\gamma^{1/2} - (\gamma_t)^{1/2}| + |g(\gamma) - g(\gamma_t)| 
    \le \frac{\epsilon'}{20}\,.
\eeq
Since $g$ is strictly increasing by Lemma~\ref{l:fp-g2}\eqref{i:fp-g-diff}, \eqref{e:fp-g-C2-choice} in fact implies that for all $\gamma \in [0,C_9]$, there exists $t\in [T]$ such that \eqref{e:fp-g-net} holds.
(In particular, for $\gamma \in [0,1/C_9)$, we can take $\gamma_t = \gamma_1 = 1/C_9$.) Finally, consider the event that $\|\dbg\| \le 2N^{1/2}$, which occurs with probability at least $1-e^{-cN}$ where $c$ is an absolute constant.  On this event, whenever \eqref{e:fp-g-net} is satisfied, we have 
    \baln
    &\bigg|
        \bigg(\sup_{\|\by\|^2=N} V_\gamma(\by) - g(\gamma)\bigg) -
        \bigg(\sup_{\|\by\|^2=N} V_{\gamma_t}(\by) - g(\gamma_t)\bigg)
    \bigg| \\
    &\le 2|\gamma^{1/2} - (\gamma_t)^{1/2}| 
    \sup_{\|\by\|^2=N} 
    \frac{|(\dbg,\by)|}{N} 
    + |g(\gamma) - g(\gamma_t)| \\
    &\le 4|\gamma^{1/2} - \gamma_t^{1/2}| + |g(\gamma) - g(\gamma_t)|
    \le \frac{\epsilon'}{20}\,.
\ealn
Now consider the event that the conclusion of Proposition~\ref{p:fp-g-one-gamma} holds for all $\gamma_t$, $t\in [T]$. This event occurs with probability at least $1-e^{-cN}$, where $c$ depends on $C_9$ and the parameters in Remark~\ref{r:fp-rpr}. On this event, we conclude that for all $\gamma \in [0,C_9]$,
\[
    \sup\Big\{ V_\gamma(\by) : \|\by\|^2=N\Big\}
    \le g(\gamma) + \frac{\epsilon'}{10}\,,
\]
which proves the claim.
\end{proof}

In light of Proposition~\ref{p:fp-g}, we define 
\begin{align}
    \label{e:fp.F.rviii}
    F^{\rviii} 
    &= \sup_{\btheta_+ \in \bbR^{B_+}}
    \bigg\{
        F^{\rviii,+}(\btheta_+)
    \bigg\} 
    + \alpha \sum_{j\in J_-} \frac{\lambda_j \hv_j}{\alpha - \hv_jm(z_\star)}
    - m(z_\star) \gamma_-\,, \\
    \label{e:fp.F.rviii.plus}
    F^{\rviii,+}(\btheta_+)
    &= g\bigg(
        \frac{\|\btheta_+\|^2}{M\alpha}
        + \gamma_-
    \bigg)
    +\frac{2(\hbg_+,\btheta_+)}{M}
    - \frac{((\hbD_+)^{-1} \btheta_+, \btheta_+)}{M}\,.
	\end{align}
We then have the following:

\begin{ppn}\label{p:fp-vi-to-vii}
Let $F^{\rvii}$ and $F^{\rviii}$ be as defined by \eqref{e:fp.F.vii} and \eqref{e:fp.F.rviii} respectively. There exists a positive constant $c$, depending on the parameters in Remark~\ref{r:fp-rpr}, such that we have
    \[
        \P(F^{\rvii} \ge t) \le
        \P\bigg(F^{\rviii} \ge t +
            \frac{\epsilon'}{10}\bigg) + e^{-cN}
    \]
for all $t\in\R$.
That is, $F^{\rvii}$ is approximately stochastically dominated by
$F^{\rviii}$.

\begin{proof}
Recall from \eqref{e:fp.F.vii} that $F^{\rvii}$ involves a supremum over $\|\by\|^2=N$ and a supremum over $\btheta_+\in\bTheta_+$, and we may exchange the order of these to write
    \begin{align*}
    F^{\rvii}
    &=\sup_{\btheta_+\in\bTheta_+}\bigg\{
    \sup_{\|\by\|^2=N}
    \bigg[
     \frac{(\dbD \by, \by)}{N}+\frac{(\dbX \by, \by)}{N}
            +\frac{2\|\btheta_+\|}{(M\alpha)^{1/2}} \frac{(\dbg_+,\by)}{N}
            + \frac{2(\gamma_-)^{1/2}(\dbg_-,\by)}{N}\bigg]\\
            &\qquad+\frac{2(\hbg_+,\btheta_+)}{M}
    -\frac{(\btheta_+,(\hbD_+)^{-1}\btheta_+)}{M}
    +\alpha \sum_{j\in J_-} \frac{\lambda_j \hv_j}{\alpha
        -\hv_j m(z_\star) }
        -m(z_\star)\gamma_-
        \bigg\}\,.
    \end{align*}
We will use Proposition~\ref{p:fp-g} to show that in the first line above, the supremum over $\|\by\|^2=N$ can be upper bounded, up to small additive error, by the $g$ term in \eqref{e:fp.F.rviii.plus}; this will imply the claimed result. 
To this end, consider any $\gamma_+ \ge 0$, and note $(\gamma_+)^{1/2} \dbg_+ + (\gamma_-)^{1/2} \dbg_-$ is equidistributed as $(\gamma_+ + \gamma_-)^{1/2} \dbg$, where $\dbg$ is a standard gaussian in $\R^N$. By Lemma~\ref{l:fp-uniform-bound},
\beq
    \label{e:fp-vi-to-vii-gamma-minus}
    \gamma_- 
    \stackrel{\eqref{e:fp.gamma.minus}}{\le} \gamma_\star
    \le C_7
    = C_7(\alpha,\iota)\,.
\eeq
We now apply Proposition~\ref{p:fp-g} with constant $C_9 = (C_8 / \alpha)^2 + C_7$, with $C_8=C_8(\alpha,\iota)$ as in Lemma~\ref{l:fp-theta-plus-bdd} and \eqref{e:fp.restricted.Theta}.  This implies that with probability at least $1-e^{-cN}$, simultaneously for all $\gamma_+ \ge 0$ such that $\gamma_+ + \gamma_- \in [0,C_9]$, we have
    \beq
        \label{e:fp-v-to-vii}
        \sup_{\|\by\|^2=N} \bigg\{
            \frac{(\dbD \by, \by)}{N}
            + \frac{(\dbX \by, \by)}{N}
            + \frac{2(\gamma_+)^{1/2}(\dbg_+,\by)}{N}
            + \frac{2(\gamma_-)^{1/2}(\dbg_-,\by)}{N}
        \bigg\}
        \le g(\gamma_+ + \gamma_-)
        + \frac{\epsilon'}{10}
        \,.
    \eeq
In particular, by \eqref{e:fp-vi-to-vii-gamma-minus}, the estimate \eqref{e:fp-v-to-vii} holds for all $\gamma_+ \in [0,(C_8 / \alpha)^2]$. Finally, let 
    \[
    \gamma_+ 
    = \frac{\|\btheta_+\|^2}{M\alpha}
    \le \bigg(\frac{C_8}{\alpha}\bigg)^2\,,
    \]
where the last inequality is by the restriction of $\btheta_+$ to $\bTheta_+$, as in \eqref{e:fp.restricted.Theta}. We conclude that with probability at least $1-e^{-cN}$, we have
    \[
        \frac{(\dbD \by, \by)}{N}
        + \frac{(\dbX \by, \by)}{N}
        +\frac{2\|\btheta_+\|}{(M\alpha)^{1/2}} \frac{(\dbg_+,\by)}{N}
        + \frac{2(\gamma_-)^{1/2}(\dbg_-,\by)}{N}
        \le g\bigg(
            \frac{\|\btheta_+\|^2}{M\alpha}
            + \gamma_-
        \bigg) +\frac{\epsilon'}{10}\,.
    \]
    for all $\|\by\|^2=N$ and all $\btheta_+\in\bTheta_+$. This shows that with probability at least $1-e^{-cN}$, we have
    \[
        F^{\rvii}\le
        \sup_{\btheta_+ \in \bTheta_+}
        \bigg\{
            F^{\rviii,+}(\btheta_+)
        \bigg\} 
        + \alpha \sum_{j\in J_-} \frac{\lambda_j \hv_j}{\alpha - \hv_jm(z_\star)}
        - m(z_\star) \gamma_-
        + \frac{\epsilon'}{10}\,.
    \] 
Finally, we note that relaxing the supremum to $\btheta_+ \in \bbR^{B_+}$ as in \eqref{e:fp.F.rviii} only makes the value of $F^{\rviii}$ larger, so the conclusion follows.
\end{proof}
\end{ppn}

\begin{ppn}
    \label{p:fp-vii-estimate}
Let $F^{\rviii}$ be defined by \eqref{e:fp.F.rviii}, and let $V(z_\star)$   be defined by  \eqref{e:freeprob.V} and \eqref{e:fp.z.star}. Then, under both cases \eqref{i:fp-case2} and \eqref{i:fp-case3} of Remark~\ref{r:fp.cases}, we have
    \[
        \P\bigg(F^{\rviii} \ge V(z_\star) + \frac{\epsilon'}{10}\bigg) \le e^{-cN}
    \]
where $c$ is a positive constant depending on all the parameters in Remark~\ref{r:fp-rpr}.
\end{ppn}

In what follows we first present the
\hyperlink{proof:p.fp-vii-estimate-Jplus-empty}{proof of Proposition~\ref{p:fp-vii-estimate} for $\hv_1 \le \delta_2$,} where $\delta_2$ is a small positive constant depending on all the parameters in Remark~\ref{r:fp-rpr}: explicitly, we take 
\beq\label{e:freebrob.delta2}
    \delta_2 = \epsilon' \min\bigg\{
        \frac{\alpha}{5(C_6 + 4\alpha^{1/2})^2},
        \frac{1}{20(\alpha C_7)^2}
    \bigg\} \,.
\eeq
In this case, we will see that either the optimization over $\btheta_+$ is trivial (if $\hv_1 < 0$), or the optimal value of $\btheta_+$ is close to $\bzero$ and the optimization over $\btheta_+$ can be bounded crudely.
The \hyperlink{proof:p.fp-vii-estimate}{proof of Proposition~\ref{p:fp-vii-estimate} for $\hv_1>\delta_2$} is presented later.

\begin{proof}[\hypertarget{proof:p.fp-vii-estimate-Jplus-empty}{Proof of Proposition~\ref{p:fp-vii-estimate} for $\hv_1 \le \delta_2$}]
    We will first show that with probability $1-e^{-cN}$,
    \beq\label{e:fp-Jplus-empty-goal}
    \sup_{\btheta_+ \in \bbR^{B_+}} F^{\rviii,+}(\btheta_+)
        \le g(\gamma_\star) + 
        \frac{\epsilon'}{20}\,.
    \eeq
If $\hv_1 < 0$, then $J_+ = \emptyset$, and the definition \eqref{e:fp.F.rviii.plus} of $F^{\rviii,+}$ simplifies to 
    \[
    F^{\rviii,+}(\btheta_+) 
    =g(\gamma_-)
    \stackrel{\eqref{e:fp.gamma.minus}}{\le}
    g(\gamma_\star)\,,
    \]
where we recall that $g$ is strictly increasing by Lemma~\ref{l:fp-g2}\eqref{i:fp-g-diff}. Thus \eqref{e:fp-Jplus-empty-goal} immediately follows. Therefore, assume $\hv_1 > 0$. We further assume that $\|\hbg_+\| \le 2M^{1/2}$, which occurs with probability $1-e^{-cN}$. Recall from Lemma~\ref{l:fp-scalar-fn-lipschitz}\eqref{i:fp-g-sq-lipschitz} that the function $\gamma \mapsto g(\gamma^2)$ is $C_6$-Lipschitz. Using the assumption $\hv_1\le\delta_2$, with probability $1-e^{-cN}$ we can bound 
    \baln
       &\sup_{\btheta_+ \in \bbR^{B_+}} 
        F^{\rviii,+}(\btheta_+)
        \le g(\gamma_-) +
        \sup_{\btheta_+ \in \bbR^{B_+}}
        \bigg\{  C_6 \bigg[
        \bigg(\frac{\|\btheta_+\|^2}{M\alpha} + \gamma_-\bigg)^{1/2}
            - (\gamma_-)^{1/2}
        \bigg]
        + \frac{4\|\btheta_+\|}{M^{1/2}}
        - \frac{\|\btheta_+\|^2}{\delta_2 M} \bigg\} \\
        &\qquad\le g(\gamma_\star) 
        + \sup_{\btheta_+ \in \bbR^{B_+}}\bigg\{ 
        \frac{(C_6+4\alpha^{1/2})
            \|\btheta_+\|}{(M\alpha)^{1/2}}
        - \frac{\|\btheta_+\|^2}{\delta_2 M} \bigg\}\\
        &\qquad\le g(\gamma_\star) + \frac{(C_6+4\alpha^{1/2})^2 \delta_2}{4\alpha}
        \stackrel{\eqref{e:freebrob.delta2}}{\le} g(\gamma_\star) + \frac{\epsilon'}{20}\,.
    \ealn
Thus we have proved \eqref{e:fp-Jplus-empty-goal} whenever $\hv_1\le\delta_2$. By Lemma~\ref{l:fp-vii-estimate-prep}, $z_\star = \beta_\star(\gamma_\star)$, and we now turn our attention to the main assertion of the proposition.
    It follows from the above that
    \baln
        F^{\rviii}
        &\le g(\gamma_\star) 
        + \alpha \sum_{j\in J_-} \frac{\lambda_j \hv_j}{\alpha - \hv_j m(z_\star)}
        - m(z_\star) \gamma_- + \frac{\epsilon'}{20}\\
        &= z_\star
        + \alpha \sum_{j\in J_-} \frac{\lambda_j \hv_j}{\alpha - \hv_j m(z_\star)}
        + m(z_\star) (\gamma_\star - \gamma_-)
        + \frac{\epsilon'}{20} \\
     &\le V(z_\star) 
        + m(z_\star) (\gamma_\star - \gamma_-)
        + \frac{\epsilon'}{20}\,.
    \ealn
    If $\hv_1 < 0$, then $\gamma_\star = \gamma_-$ and the result follows.
    Otherwise assume $\hv_1 > 0$. 
    Since $z_\star \ge z_0 = m^{-1}(\alpha / \hv_1)$ and $m$ is strictly decreasing,
    \[
        m(z_\star) \le m(z_0) = \frac{\alpha}{\hv_1}\,.
    \]
Combining with the estimate from Lemma~\ref{l:fp-uniform-bound} gives
    \[
        \gamma_\star - \gamma_- 
        = \alpha \sum_{j\in J_+} \lambda_j \bigg(\frac{\hv_j}
            {\alpha -\hv_j m(z_\star)}\bigg)^2
        \le \alpha \bigg(\frac{\hv_1}{\alpha - \hv_1 m(z_\star)}\bigg)^2
        \le \alpha (\hv_1 C_7)^2\,.
    \]
Combining the above estimates gives
    \[
        m(z_\star) (\gamma_\star - \gamma_-)
        \le \frac{\alpha}{\hv_1} \cdot \alpha (\hv_1 C_7)^2
        \le (\alpha C_7)^2 \delta_2
        \stackrel{\eqref{e:freebrob.delta2}}{\le} \frac{\epsilon'}{20}\,.
    \]
    The result follows by substituting into the preceding bound on $F^{\rviii}$.\end{proof}

We next prove several lemmas in preparation for the \hyperlink{proof:p.fp-vii-estimate}{proof of Proposition~\ref{p:fp-vii-estimate} for $\hv_1>\delta_2$}.

\begin{lem}
    \label{l:fp-vii-estimate-prep-OTHER}
    Suppose $J_+ \neq \emptyset$.
    Then, the supremum in \eqref{e:fp.F.rviii} is attained at a point of the form
    \[
        \btheta_+ = \hbA_+(\gamma) \hbg_+\,,
    \]
where we define
    \begin{align}
        \label{e:fp-gamma-btheta-plus}
        \gamma
        &\equiv 
        \gamma(\btheta_+)
        \equiv \frac{\|\btheta_+\|^2}{M\alpha} + \gamma_-\,,\\
        \label{e:fp-viii-hbA}
        \hbA_+(\gamma)
        &\equiv \bigg((\hbD_+)^{-1} -
        \frac{m(\beta_\star(\gamma))}{\alpha} \bI_{B_+} \bigg)^{-1}\,.
    \end{align}
    Furthermore, this value of $\gamma$ satisfies
    \beq
        \label{e:fp-m-betast-gamma-bd}
        m(\beta_\star(\gamma)) \le \frac{\alpha}{\hv_1}\,.
    \eeq

\begin{proof}
Recall that Lemma~\ref{l:fp-g2}\eqref{i:fp-g-growth} shows that $g(\gamma) = \Theta(\gamma^{1/2})$ as $\gamma \to \infty$. It follows that $F^{\rviii,+}(\btheta_+) \to -\infty$ as $\|\btheta_+\| \to +\infty$, so the supremum of $F^{\rviii,+}(\btheta_+)$ is attained at a local maximizer $\btheta_+$.
    We will show that this local maximizer must in fact be a stationary point.
    By Lemma~\ref{l:fp-g2}\eqref{i:fp-g-diff}, $g$ is differentiable on $(0,+\infty)$. Therefore, the only possible non-stationary local maximum $\btheta_+$ must satisfy $\gamma = 0$, which implies $\gamma_-=0$ and $\btheta_+=\bzero$. However, in this case let us consider a small $\bv \in \bbR^{B_+}$ with $(\bv,\hbg_+) = 0$. For such $\bv$, we have
    \[
        F^{\rviii,+}(\bv) - F^{\rviii,+}(\bzero) 
        = g\bigg(\frac{\|\bv\|^2}{M\alpha}\bigg) - g(0)
        - \frac{((\hbD_+)^{-1} \bv,\bv)}{M}\,.
    \]
Lemma~\ref{l:fp-g2}\eqref{i:fp-g-diff} characterizes the behavior of $g(\gamma)$ as $\gamma\downarrow0$, and implies that the above is positive for sufficiently small $\|\bv\|$.
    Thus $\btheta_+=\bzero$ is not a local maximizer, and the only possible local maximizers are stationary points $\btheta_+$.
Therefore, consider any stationary local maximizer $\btheta_+ \in \bbR^{B_+}$ of $F^{\rviii,+}(\btheta_+)$.
    Stationarity implies
    \[
        \btheta_+ = \bigg((\hbD_+)^{-1} -
        \frac{g'(\gamma)}{\alpha} \bI_{B_+} \bigg)^{-1} \hbg_+
    \stackrel{\eqref{e:fp-viii-hbA}}{=}
    \hbA_+(\gamma) \hbg_+\,.
    \]
    In the last equality, we used Lemma~\ref{l:fp-g2}\eqref{i:fp-g-diff} with the above observation that $\gamma>0$. This proves \eqref{e:fp-gamma-btheta-plus}.
    Similarly, for any $\bv$ orthogonal to $\btheta_+$
    (or, if $\btheta_+=\bzero$, for any $\bv$), we have 
    \[
        \frac{d^2}{dt^2}
        F^{\rviii,+}(\btheta_+ + t\bv)
        \bigg|_{t=0}
        = -\frac{2}{M} \bv^\st \bigg(
        (\hbD_+)^{-1} -
        \frac{m(\beta_\star(\gamma))}{\alpha} \bI_{B_+}
        \bigg) \bv\le0\,,
    \]
where the last inequality is because 
$\btheta_+$ is a maximizer. Since $(\hbD_+)^{-1}$ has smallest eigenvalue $1/\hv_1$ with multiplicity at least $\iota N$ (because $(\lambda,\mu) \in \Lambda$), this implies the bound \eqref{e:fp-m-betast-gamma-bd}.
\end{proof}
\end{lem}

Recall that we let $\gamma_+(z)$ denote the contribution to
\eqref{e:fp.gamma.function} from indices $j\in J_+$, so that
    \[\gamma_+(z)
        \equiv\alpha \sum_{j\in J_+} \lambda_j
        \bigg(
        \frac{ \hv_j}{\alpha - \hv_j m(z)}
        \bigg)^2\,.\]
We then have the following:

\begin{lem}\label{l:fp-vartheta}
Suppose $J_+ \neq \emptyset$, and define
    \begin{align}
        \label{e:fp-Delta-plus}
        \Delta_+(\vartheta) 
        &\equiv \vartheta + \gamma_+(z_\star) - \gamma_+(\beta_\star(\gamma_\star + \vartheta))
        \,,\\
        \label{e:fp-Delta-minus}
        \Delta_-(\vartheta) 
        &\equiv \vartheta + \gamma_+(\beta_\star(\gamma_\star - \vartheta)) - \gamma_+(z_\star)\,.
    \end{align}
These functions are continuous and strictly increasing on the intervals $[0,+\infty)$ and $[0,\vartheta_{-,\max})$ respectively, where in the latter case we define
    \[
        \vartheta_{-,\max} \equiv \sup\lt\{
            \vartheta \ge 0 : \beta_\star(\gamma_\star - \vartheta) > z_0
        \rt\}\,.
    \]
    Furthermore, $\Delta_+(0) = \Delta_-(0) = 0$ and 
    \beq
        \label{e:fp-vartheta-limit}
        \lim_{\vartheta \uparrow +\infty} \Delta_+(\vartheta)
        = \lim_{\vartheta \uparrow \vartheta_{-,\max}} \Delta_-(\vartheta)
        = +\infty\,.
    \eeq
    In particular, for any $\delta > 0$, there exist (unique) solutions $\vartheta_+ \in (0,+\infty)$, $\vartheta_- \in (0,\vartheta_{-,\max})$ to the equations
    \beq
        \label{e:fp-vartheta-plus-minus}
        \Delta_+(\vartheta_+) = \Delta_-(\vartheta_-) = \delta\,.
    \eeq
\begin{proof}
    Since $J_+ \neq \emptyset$, we have $\hv_1 > 0$ and $z_0 = m^{-1}(\alpha / \hv_1)$, and thus case \eqref{i:fp-case2} or \eqref{i:fp-case3} of Remark~\ref{r:fp.cases} holds. In either case, we have $z_\star > z_0$. Note that $\gamma_+$ is strictly decreasing for $z > z_0$. Since $\beta_\star(\gamma_\star) = z_\star$, we have $\Delta_+(0) = \Delta_-(0) = 0$.
Recall that $\beta_\star$ is nondecreasing by Lemma~\ref{l:fp-g2}\eqref{i:fp-beta-incr}; it follows that $\beta_\star(\gamma_\star - \vartheta) > z_0$ for all $\vartheta \in [0,\vartheta_{-,\max})$. Moreover, $\beta_\star$ is continuous, so $\Delta_+$ and $\Delta_-$ are continuous and strictly increasing on their respective domains. It remains to verify \eqref{e:fp-vartheta-limit}.
Since $\Delta_+(\vartheta)\ge\vartheta$, it is clear the limit involving $\Delta_+$ is $+\infty$.
Similarly, since $\Delta_-(\vartheta)\ge\vartheta$, 
 the limit involving $\Delta_-$ is also clearly $+\infty$ in the case $\vartheta_{-,\max} = +\infty$.
    Otherwise, we must have $\beta_\star(\gamma_\star - \vartheta) \downarrow z_0$ as $\vartheta \uparrow \vartheta_{-,\max}$, so $\gamma_+(\beta_\star(\gamma_\star - \vartheta)) \uparrow +\infty$.
    This proves \eqref{e:fp-vartheta-limit}.
\end{proof}
\end{lem}

\begin{lem}
    \label{l:fp-m-stability}
Recall from \eqref{e:freeprob.def.z0} the definition of $z_0$.
    Suppose $\hv_1 > \delta_2$, and define
    \[
        \delta_3 = \frac{\iota \delta_2}{2\alpha}\,.
    \]
    Then, we have $z_0 - \delta_3 > \hw_1$ and $m_2(z_0 - \delta_3) \le 1/(\delta_3)^2$.
\begin{proof}
    Let $t = \hw_1 + \mu_1 \hv_1 / \alpha$. Then,
    \[
        m(t) \ge \frac{\mu_1}{t-\hw_1} = \frac{\alpha}{\hv_1} = m(z_0)\,.
    \]
Since $m$ is strictly decreasing, we conclude $t \le z_0$. Therefore
    \[
        z_0 - \delta_3 \ge t - \delta_3 
        > \hw_1 + \frac{\iota \delta_2}{\alpha} - \delta_3
        = \hw_1 + \delta_3\,,
    \]
where the last equality is by the definition of $\delta_3$. This proves the first assertion. Then, since $m_2$ is strictly decreasing, we conclude
    \[
        m_2(z_0 - \delta_3)
        \le m_2(t - \delta_3)
        < m_2(\hw_1 + \delta_3)
        \le \frac1{(\delta_3)^2}\,.
    \]
This proves the second assertion, so we are done.
\end{proof}
\end{lem}

\begin{proof}[\hypertarget{proof:p.fp-vii-estimate}{Proof of Proposition~\ref{p:fp-vii-estimate} for $\hv>\delta_2$}] Since $\hv_1 > \delta_2$, we have $J_+\neq\emptyset$.
Let $\btheta_+$ be a maximizer for \eqref{e:fp.F.rviii}, which exists by Lemma~\ref{l:fp-vii-estimate-prep-OTHER}. Let $\gamma = \gamma(\btheta_+)$ be defined by \eqref{e:fp-gamma-btheta-plus}: we will argue that $\gamma$ must be close to $\gamma_\star=\gamma(z_\star)$, as defined by
    \eqref{e:fp.gamma.function} and \eqref{e:fp.z.star}.
    Let $\delta>0$ be a small constant we will set below, which will depend on only the parameters in Remark~\ref{r:fp-rpr}.
    Let $\vartheta_+,\vartheta_-$ be the solutions to \eqref{e:fp-vartheta-plus-minus}, as guaranteed by Lemma~\ref{l:fp-vartheta}. Recall from the proof of Lemma~\ref{l:fp-vartheta} that the functions $\Delta_+(\vartheta)$ and $\Delta_-(\vartheta)$ are both lower bounded by $\vartheta$; this implies that $\vartheta_+,\vartheta_- \le\delta$. However, we emphasize that (unlike $\delta$) $\vartheta_+,\vartheta_-$ may \textbf{not} be uniformly lower bounded depending on the parameters of Remark~\ref{r:fp-rpr}.
    The constant $c$ in the probabilities $1-e^{-cN}$ below will be nonetheless be uniform in the parameters of Remark~\ref{r:fp-rpr}.

    Since $\tilde\gamma \mapsto \gamma_+(\beta_\star(\tilde\gamma))$ is nonincreasing on $[\gamma_\star - \vartheta_-, \gamma_\star + \vartheta_+]$, for all such $\tilde \gamma$ we have
    \beq
        \label{e:gamma-plus-bound-in-nbd}
        \gamma_+(\beta_\star(\tilde\gamma))
        \le \gamma_+(\beta_\star(\gamma_\star - \vartheta_-))
        = \gamma_+(z_\star) + \Delta_-(\vartheta_-) - \vartheta_-
        \le \gamma_\star + \delta
        \le C_7 + \delta\,,
    \eeq
where the last estimate is by Lemma~\ref{l:fp-uniform-bound}.
    Then, for $\hbA_+$ defined in \eqref{e:fp-viii-hbA}, (recalling $\lambda_1 \ge \iota$ by the assumption $(\lambda,\mu) \in \Lambda$) 
    \baln
        &\|\hbA_+(\gamma_\star + \vartheta_+)\|_{\op} 
        \le \|\hbA_+(\gamma_\star - \vartheta_-)\|_{\op} 
        = \frac{\alpha \hv_1}{\alpha - \hv_1 m(\beta_\star(\gamma_\star - \vartheta_-))} \\
        &\qquad\le 
        \bigg(
        \frac{\alpha \gamma_+(\beta_\star(\gamma_\star - \vartheta_-))}{\lambda_1}
        \bigg)^{1/2}
        \stackrel{
            \eqref{e:gamma-plus-bound-in-nbd}
        }{\le} \bigg(
        \frac{\alpha (C_7+\delta)}{\iota}\bigg)^{1/2}\,.
    \ealn
By Lemma~\ref{l:fp-gaussian-matrix-prod} (with $[\alpha (C_7+\delta)/\iota]^{1/2}$ in place of $C$), it holds with probability at least $1-e^{-cN}$ that 
    \baln
        \frac{\|\hbA_+(\gamma_\star
            + \vartheta_+ ) \hbg_+\|^2}{M\alpha}
        &\le 
        \alpha \sum_{j\in J_+}
        \lambda_j\bigg(
        \frac{ \hv_j}{\alpha -
            \hv_jm(\beta_\star(\gamma_\star + \vartheta_+))}\bigg)^2
        +\delta = \gamma_+(\beta_\star(\gamma_\star+ \vartheta_+))  +\delta \\
    &=  
    \gamma_+(z_\star) + \vartheta_+ - \Delta_+(\vartheta_+) + \delta 
        = \gamma_\star - \gamma_-
        + \vartheta_+ \,.
    \ealn
Note that $\gamma \mapsto \|\hbA_+(\gamma) \hbg_+\|^2$ is nonincreasing for $\gamma$ satisfying \eqref{e:fp-m-betast-gamma-bd}. Therefore,  for $\gamma > \gamma_\star + \vartheta_+$,
    \[
        \gamma - \gamma_-
        = \frac{\|\btheta_+\|^2}{M\alpha}
            =\frac{\|\hbA_+(\gamma) \hbg_+\|^2}{M\alpha}
        \le \frac{\|\hbA_+(\gamma_\star+\vartheta_+) \hbg_+\|^2}{M\alpha}
         \leq 
        \gamma_\star - \gamma_-
        + \vartheta_+ \,,
        \,,
    \]
which yields a contradiction to the assumption $\gamma > \gamma_\star + \vartheta_+$. Similarly, with probability at least $1-e^{-cN}$, another application of Lemma~\ref{l:fp-gaussian-matrix-prod} gives 
    \[\frac{\|\hbA_+(\gamma_\star- \vartheta_- ) \hbg_+\|^2}{M\alpha}
    \ge \gamma_+(\beta_\star(\gamma_\star- \vartheta_- )) - \delta
= \gamma_+(z_+) 
    -\vartheta_-
    + \Delta_-(\vartheta_-) 
    -\delta 
    =\gamma_\star - \gamma_-  - \vartheta_-\,.\]
Therefore,  for $\gamma < \gamma_\star -  \vartheta_-$,     we have
    \[
        \gamma - \gamma_-
        = \frac{\|\btheta_+\|^2}{M\alpha}
        =\frac{\|\hbA_+(\gamma) \hbg_+\|^2}{M\alpha}
        \ge \frac{\|\hbA_+(\gamma_\star- \vartheta_-) \hbg_+\|^2}{M\alpha}
        \geq  \gamma_\star - \gamma_-
- \vartheta_-\,,
    \]
    which again yields a
    contradiction to the assumption $\gamma < \gamma_\star -  \vartheta_-$. Therefore 
        \beq\label{e:freeprob.located.gamma}
    \gamma \in [\gamma_\star - \vartheta_-, \gamma_\star + \vartheta_+]\,.\eeq   
    Having located $\gamma$, we return to the calculation of $F^{\rviii}$, which amounts to substituting the expression for $\btheta_+$ (from Lemma~\ref{l:fp-vii-estimate-prep-OTHER}) into \eqref{e:fp.F.rviii.plus} and evaluating. Applying Lemma~\ref{l:fp-gaussian-matrix-prod} gives
        \begin{align*}
    &\frac{2(\hbg_+,\btheta_+)}{M}
        - \frac{((\hbD_+)^{-1} \btheta_+, \btheta_+)}{M}\\
            &=\frac{2(\hbg_+,\btheta_+)}{M}
        - \frac{( [(\hbD_+)^{-1}
                - m(\beta_\star(\gamma)) \bI /\alpha ]\btheta_+, \btheta_+)}{M}
    -m(\beta_\star(\gamma)) \frac{\|\btheta_+\|^2}{M\alpha} \\
    &=\frac{(\hbg_+,\btheta_+)}{M}
        -m(\beta_\star(\gamma))
        (\gamma - \gamma_-)
    \approx_\delta
    \alpha\sum_{j\in J_+}
    \frac{\lambda_j\hv_j}{\alpha-\hv_j m(\beta_\star(\gamma))}
        -m(\beta_\star(\gamma))
        (\gamma - \gamma_-)\,.
    \end{align*}
Substituting this into the expression \eqref{e:fp.F.rviii} for  $F^{\rviii}$ gives
    \begin{align*}
    F^{\rviii}
    &\approx_\delta
    g(\gamma) +
    \alpha\sum_{j\in J_+}
    \frac{\lambda_j\hv_j}{\alpha-\hv_j m(\beta_\star(\gamma))}
        -m(\beta_\star(\gamma))
        (\gamma - \gamma_-) 
    + \alpha \sum_{j\in J_-} \frac{\lambda_j \hv_j}{\alpha - \hv_jm(z_\star) }
        - m(z_\star) \gamma_- \\
        &=\beta_\star(\gamma)
        +\Big( m(\beta_\star(\gamma))-m(z_\star)\Big) \gamma_-
        +\alpha\sum_{j\in J_+}
    \frac{\lambda_j\hv_j}{\alpha-\hv_j m(\beta_\star(\gamma))}
        +\alpha\sum_{j\in J_-}
    \frac{\lambda_j\hv_j}{\alpha-\hv_j m(z_\star)} \\
    &= V(z_\star)
    + \Big(
        \beta_\star(\gamma) - z_\star
    \Big)
    + \Big( m(\beta_\star(\gamma))-m(z_\star)\Big)  \gamma_-\\
        &\qquad\qquad\qquad\qquad
    \qquad\qquad\qquad\qquad
     + \alpha 
        \sum_{j\in J_+} \bigg\{
        \frac{\lambda_j\hv_j}{\alpha-\hv_j m(\beta_\star(\gamma))}
        -\frac{\lambda_j\hv_j}{\alpha-\hv_j m(z_\star)}
        \bigg\}
        \,.
        \end{align*}
We next bound the error terms above.
    Recall that $\vartheta_+,\vartheta_- \le \delta$, so $|\gamma - \gamma_\star| \le \delta$.
    Since $\gamma \mapsto \beta_\star(\gamma^2)$ is $C_6$-Lipschitz by Lemma~\ref{l:fp-scalar-fn-lipschitz}\eqref{i:fp-betast-sq-lipschitz}, we have
    \[
        \Big|
        \beta_\star(\gamma) - z_\star
        \Big|
        = \Big|\beta_\star(\gamma) - \beta_\star(\gamma_\star)\Big|
        \le C_6\Big|\gamma^{1/2} - (\gamma_\star)^{1/2}\Big| 
        \le C_6 \delta^{1/2}\,.
    \]
Let $\delta_3$ be as defined in Lemma~\ref{l:fp-m-stability}. We will fix $\delta$ later, but it will be chosen small enough so that we have 
$C_6 \delta^{1/2} \le \delta_3$.
    Since $m' = -m_2<0$ and $m_2$ is strictly decreasing, we obtain
    \[
        \Big|m(\beta_\star(\gamma)) - m(z_\star)\Big|
        \le 
        m_2(z_\star - C_6\delta^{1/2})
        \Big|\beta_\star(\gamma) - z_\star\Big| 
        \le 
        m_2(z_0 - \delta_3)
        C_6 \delta^{1/2} 
        \le 
    \frac{ C_6 \delta^{1/2}}{(\delta_3)^2}\,,
    \]
where the last estimate is by Lemma~\ref{l:fp-m-stability}.
    Finally, we bound
    \baln
        &\bigg|
            \sum_{j\in J_+} \bigg\{
            \frac{\lambda_j\hv_j}{\alpha-\hv_j m(\beta_\star(\gamma))}
            - \frac{\lambda_j\hv_j}{\alpha-\hv_j m(z_\star)}\bigg\}
        \bigg| 
        = \bigg|
            \int_{m(z_\star)}^{m(\beta_\star(\gamma))}
            \sum_{j\in J_+} \bigg(\frac{\lambda_j\hv_j}{\alpha-\hv_j t}\bigg)^2 dt
        \bigg| \\
        &\qquad
        \stackrel{\eqref{e:freeprob.located.gamma}}{\le} \Big|m(\beta_\star(\gamma)) - m(z_\star)\Big|
        \sup
      \bigg\{
      \gamma_+(\beta_\star(\tilde \gamma))
      :\tilde\gamma \in [\gamma_\star - \vartheta_-, \gamma_\star + \vartheta_+]
      \bigg\} \\
     &\qquad\stackrel{\eqref{e:gamma-plus-bound-in-nbd}}{\le}
        \Big|m(\beta_\star(\gamma)) - m(z_\star)\Big|(C_7 + \delta)
        \le
        \frac{C_6(C_7 + \delta) \delta^{1/2}}{(\delta_3)^2}\,.
    \ealn
    Altogether this implies
    $F^{\rviii} \le V(z_\star) + O(\delta^{1/2})$. Finally,
    setting $\delta$ sufficiently small (depending on only the parameters in Remark~\ref{r:fp-rpr}) concludes the proof.
\end{proof}

\begin{proof}[\hypertarget{proof:p.fp-v-estimate}{Proof of Proposition~\ref{p:fp-v-estimate}}] The
\hyperlink{proof:p.fp-v-estimate.case.one}{proof of Proposition~\ref{p:fp-v-estimate} under case \eqref{i:fp-case1}
of Remark~\ref{r:fp.cases}} appeared earlier. Under cases \eqref{i:fp-case2} and \eqref{i:fp-case3} of Remark~\ref{r:fp.cases}, the result
follows from
    Proposition~\ref{p:fp-Fdiam-to-V} combined with  Propositions~\ref{p:fp-v-to-vi}, \ref{p:fp-vi-to-vii}, and \ref{p:fp-vii-estimate}.
\end{proof}
\fi

\vfill

\pagebreak\footnotesize
\bibliographystyle{alpha}
\bibliography{bib}

\end{document}